%% file: english.tex
\documentclass[12pt,diff_eng]{article} 
\usepackage[english]{babel}
\usepackage[cp1251]{inputenc}
\usepackage{amssymb,amscd,amsmath}     
\usepackage{diff_eng}

\usepackage{graphicx}
\usepackage{listings}
\usepackage{hyperref}
\usepackage{caption}
\usepackage{float}
\usepackage{caption}
\usepackage{subcaption}
\usepackage{letltxmacro}
\usepackage{xcolor}
\usepackage{mdframed}
\usepackage{sectsty}

\hypersetup{
  raiselinks=true,
  linktocpage,
  linktoc=page
}

\subsubsectionfont{\large}

\definecolor{light-gray}{gray}{0.95}

\AtBeginDocument{
  \LetLtxMacro{\oldref}{\ref}
  \renewcommand{\ref}[1]{\mbox{\oldref{#1}}}
}

\PassOptionsToPackage{linktocpage}{hyperref}

\begin{document}

  \input{diff_eng.tex}
  \input{chapters/theory.tex}

\input{chapters/software.tex}

\input{chapters/bibliography.tex}

\end{document}

%% file: DIFF_ENG.TEX
\thispagestyle{empty}
\large
\bigskip

%% file: chapters/theory.tex
$$
{}
$$

\begin{center}
\huge{\bf Implementation of Strong Numerical Methods of Orders 
0.5, 1.0, 1.5, 2.0, 2.5, and 3.0 for It\^o SDEs with
Multidimensional Non-Commutative Noise Based on the Unified 
Taylor--It\^o and Taylor--Stratonovich Expansions and Multiple Fourier--Legendre Series}
\end{center}

\vspace{1cm}

\begin{center}
\LARGE
Mikhail D. Kuznetsov${}^{1}$ and Dmitriy F. Kuznetsov${}^{2}$
\end{center}

\vspace{1cm}

\begin{center}
\Large
${}^{1}$Faculty of Computer Technologies and Informatics, St. Petersburg Electrotechnical University,
Saint-Petersburg, Russia\\
${}^{2}$Institute of Applied Mathematics and Mechanics,
Peter the Great St. Petersburg Polytechnic University\\
e-mail: sde\_\hspace{0.5mm}kuznetsov@inbox.ru
\end{center}

\vspace{2cm}

\noindent
{\bf Abstract.} The article is devoted to the 
implementation of strong numerical methods with convergence orders
$0.5,$ $1.0,$ $1.5,$ $2.0,$ $2.5,$ and $3.0$  
for It\^{o} stochastic differential equations 
with multidimensional non-commutative noise 
based on multiple Fourier--Legendre
series and unified Taylor--It\^{o} and 
Taylor--Stratonovich expansions.
Algorithms for the implementation of these methods 
are constructed and a package of programs in the Python programming 
language is presented.
An important part of this software package
concerning the mean-square approximation of 
iterated It\^{o} and Stratonovich stochastic integrals
of multiplicities 1 to 6
with respect to components of the multidimensional Wiener process
is based on the method of 
generalized multiple Fourier series.
More precisely, we used multiple Fourier--Legendre series
converging in the sense of norm in Hilbert space 
for the mean-square approximation of iterated It\^{o} and 
Stratonovich stochastic integrals.

\vspace{5mm}

\renewcommand{\baselinestretch}{1.3}{\normalsize{\tableofcontents}}

\newpage

\section{Introduction}

As known, It\^{o} stochastic differential equations (SDEs) have appeared 
in the theory of random processes relatively recently \cite{1} (1951).
Nevertheless, to date, a large number of mathematical models for dynamical 
systems of different physical nature 
under the influence of random perturbations have been built on the basis of such equations
\cite{2}-\cite{16}. Among them we note mathematical models 
in stochastic financial mathematics \cite{5}-\cite{7}, \cite{10}-\cite{12}, 
geophysics \cite{2}, \cite{4}, genetics \cite{13}, hydrology \cite{2},
epidemiology \cite{9}, chemical kinetics \cite{2}, \cite{9}, biology \cite{8}, \cite{15}, 
seismology \cite{2}, electrodynamics \cite{16} and many other fields \cite{2}, \cite{9}, \cite{14}.
In addition, It\^{o} SDEs arise when solving a number of mathematical problems, 
such as filtration \cite{2}, \cite{3}, \cite{17}-\cite{21}, 
stochastic control \cite{2}, \cite{17}, stochastic stability \cite{2}, 
parameter estimation of stochastic systems \cite{2}, \cite{3}, \cite{22}.

Exact solutions of It\^{o} SDEs 
are known in rare cases. 
For this reason, it becomes necessary to construct numerical 
methods for It\^{o} SDEs.
Moreover, the problem of numerical solution of It\^{o} SDEs 
often occurs even in cases when the exact solution of It\^{o} SDE is known.
This means that in some cases, knowing the exact solution of the It\^{o} 
SDE does not allow us to simulate it numerically in a simple way.

This article is devoted to the 
implementation of high-order strong numerical methods 
for systems of It\^{o} SDEs
with multidimensional non-commutative noise. More precisely, 
we consider strong numerical methods 
with convergence orders $1.0,$ $1.5,$ $2.0,$ $2.5,$ and $3.0$.
The article also considers the Euler method, which under suitable conditions \cite{2}
has the order $0.5$ of strong convergence.
To construct the mentioned numerical methods in this article, we use the so-called
unified Taylor--It\^{o} and
Taylor--Stratonovich expansions \cite{24}-\cite{26a}, \cite{58}.
The important components of these expansions are the iterated It\^{o} and Stratonovich 
stochastic integrals, which are functionals of a complex structure with 
respect to the components of a multidimensional Wiener process.

It should be noted that it is impossible to construct a numerical 
method for It\^{o} SDE in a general case 
(multidimensional non-commutative noise)
that includes only increments of the multidimensional Wiener 
processes, but has a higher order of convergence (in the mean-square sense) than the Euler method
(simplest numerical method for It\^{o} SDEs).
This result is known as the "Clark--Cameron paradox" \cite{23} (1980) and well explains the 
need to use high-order numerical methods for It\^{o} SDEs, since the accuracy of 
the Euler method is insufficient for solving a number of 
practical problems related to It\^{o} SDEs \cite{2}.
According to the "Clark--Cameron paradox" \cite{23}, 
avoidance of the problem of mean-square approximation of the mentioned 
iterated stochastic integrals is impossible in the general case when constructing 
high-order strong numerical methods for It\^{o} SDEs.

The problem of mean-square approximation of iterated It\^{o} and 
Stratonovich stochastic integrals in the context of the numerical integration 
of It\^{o} SDEs was considered in a number of works 
\cite{2}, \cite{3}, \cite{7}, \cite{8}, \cite{27}-\cite{39}.

It should be explained why the results of these works are insufficient for 
constructing effective procedures for the implementation of strong numerical 
methods of order $1.5$ and higher for It\^{o} SDEs.

There exists an approach to the 
mean-square approximation of
iterated stochastic integrals based on 
integral sums \cite{27}, \cite{34}, \cite{35}. 
Note that one of the variants of this method is based on reducing 
the problem of mean-square approximation of iterated stochastic 
integrals to the numerical integration of systems of linear It\^{o} SDEs by the Euler method
\cite{39}. However, this approach \cite{27}, \cite{34}, \cite{35}, \cite{39} 
implies the partitioning of the interval 
of integration for iterated stochastic integrals.
It should be noted that the length of this interval
is an integration step for numerical methods for It\^{o} SDEs, which is 
already a fairly small value even without additional partitioning.
Computational experiments show that the 
numerical modeling of iterated stochastic integrals by the method
of integral sums \cite{27}, \cite{34}, \cite{35}, \cite{39} 
leads to unacceptably high computational cost and 
accumulation of computation errors \cite{40}.

More efficient approach of the mean-square approximation of iterated stochastic integrals 
is based on the 
expansion of the so-called Brownian bridge process into the  
trigonometric Fourier series with random terms (version of the so-called Karhunen--Lo\`{e}ve
expansion) \cite{2}, \cite{3}, \cite{7}, \cite{27}, \cite{28}, \cite{33}, \cite{34},
\cite{37}, \cite{38}.
However, in \cite{27}, \cite{33}, \cite{34}, \cite{38}, 
this approach was used to approximate only iterated stochastic integrals 
of multiplicities $1$ and $2$, which makes it possible to implement numerical method 
with order $1.0$ of strong convergence
for It\^{o} SDEs (Milstein method \cite{27}).
In papers \cite{2}, \cite{3}, \cite{7}, \cite{28}, the approximation 
of iterated stochastic integrals of 
multiplicities $1$ to $3$ was considered by the above approach, 
which makes it possible to implement
numerical methods with orders $1.0$ and $1.5$ of strong convergence
for It\^{o} SDEs.
However, formulas concerning integrals of multiplicity $3$ turned out to be too complicated 
and did not find wide application in practice.
Moreover, these formulas (for iterated stochastic 
integrals of multiplicity $3$) were obtained without strict theoretical justification 
and exclude the possibility of effective estimation of the 
mean-square error of approximation (see discussion in \cite{26} (Sect. 2.6.2, 6.2) for details).

It should be noted that in papers \cite{29}, \cite{30}, a similar approach was used to approximate 
iterated stochastic integrals of multiplicities 1 to 3 based on the series expansion of the 
Wiener process using trigonometric functions and Haar functions.
In \cite{Rybakov1} orthonormal expansions of functions in terms of 
Walsh series were used to represent the iterated stochastic integrals.

Note that the iterated stochastic integrals under consideration 
are the random variables with
unknown density functions. The only exception is 
the iterated It\^{o} stochastic integral with multiplicity 2 \cite{31}.
However, the knowledge of density function
of the mentioned stochastic integral 
gives no simple way of its approximation \cite{31}.

In this work, we use a more efficient method
of the mean-square approximation of iterated It\^{o} and Stratonovich stochastic 
integrals than the methods considered above.
This method (the so-called method of generalized multiple Fourier series) 
is based on the theory constructed in Chapters 1, 2, and 5 of monographs \cite{26}, \cite{26a},
\cite{58}.
The method of generalized multiple Fourier series made it possible in 
this work to successfully 
implement the procedures for the mean-square approximation 
of iterated It\^{o} and Stratonovich stochastic integrals of multiplicities 1 to 6.
In this case, we use multiple Fourier--Legendre series, that is, we have chosen 
Legendre polynomials as a basis system of functions for approximating iterated stochastic 
integrals. 
It is important to note that the Legendre polynomials
were first applied in the context of this problem in \cite{41} (1997), 
while in the works of other authors Legendre polynomials were not 
considered as a system of basis functions for approximating 
iterated stochastic integrals (an exception is work \cite{36}).
As shown in \cite{42}, the Legendre polynomials are optimal 
for the implementation of strong numerical methods with 
convergence order 1.5 and higher for It\^{o} SDEs with 
non-commutative noise.

In this article, to build the SDE-MATH software package in the Python programming language, 
we use a database with 270,000 exactly calculated Fourier--Legendre coefficients 
to approximate iterated It\^{o} and Stratonovich stochastic integrals 
of multiplicities $1$ to $6$.
It should be noted that the procedures for the mean-square 
approximation of iterated stochastic integrals of 
multiplicities $4$, $5$, and $6$ constructed in this work have no analogues in the literature.
At the same time, we propose a much more convenient procedure for 
the mean-square approximation of iterated stochastic integrals of multiplicity $3$ than in works 
\cite{2}, \cite{3}, \cite{7}, \cite{28}.
This procedure provides an accurate calculation of the mean-square 
error of approximation of the mentioned stochastic integrals.

Another important feature of the presented software package is the use of 
unified Taylor--It\^{o} and Taylor--Stratonovich expansions 
\cite{24}-\cite{26a}, \cite{58}
for constructing strong numerical methods with convergence orders 
$1.5,$ $2.0$, $2.5$, and $3.0$ for It\^{o} SDEs.
Unified Taylor--It\^{o} and Taylor---Stratonovich expansions make it possible (in contrast with its
classical analogues \cite{2}) to use the minimal sets of 
iterated It\^{o} and Stratonovich stochastic integrals.
This property well explains the motive for using the mentioned unified expansions.

The results of this work on the approximation of iterated stochastic integrals 
can be used to numerically solve various types of SDEs.
For example, for semilinear SPDEs with multiplicative trace class noise
\cite{26}, \cite{26a}, \cite{58} (Chapter 7), \cite{42a}, \cite{42aa}.
This is due to the fact that iterated stochastic integrals 
are a universal tool for constructing high-order strong numerical methods 
for various types of SDEs.
In recent years, the mentioned numerical methods have been constructed for
SDEs with jumps \cite{7},
SPDEs with multiplicative trace class noise
\cite{43}-\cite{45},
McKean SDEs \cite{46}, 
SDEs with switchings
\cite{47}, mean-field SDEs \cite{48}, 
It\^{o}--Volterra stochastic integral equations \cite{45}, etc.

There are many publications in which codes of programs 
in various programming languages are given for 
the numerical solution of SDEs 
\cite{3}, \cite{9}, \cite{14}, \cite{48a}-\cite{56}.
Among them, we note the software described in 
\cite{3}, \cite{49}, \cite{51}, \cite{55}.
Some of the mentioned works \cite{3}, \cite{49}, \cite{51}, \cite{52}, \cite{55} 
are based on the results of monograph \cite{2} 
on the approximation of iterated stochastic integrals
(see above discusson on the disadvantages
of approach \cite{2}).
Other publications \cite{9}, \cite{14}, \cite{48a}, \cite{50} do not use 
the modeling of iterated stochastic 
integrals for the case of multidimensional non-commutative noise at all.

In this article, we develop software for the numerical 
integration of It\^{o} SDEs based on theoretical results 
and MATLAB codes from monographs \cite{53}, \cite{56} for modeling iterated stochastic integrals
of multiplicities 1 to 6 (the case of multidimensional non-commutative noise).
In addition, we provide software (as a part of the SDE-MATH software package) for the 
numerical integration of linear stationary systems of It\^{o} SDEs based on the results 
of article \cite{57} and MATLAB codes from monographs \cite{53}, \cite{56}.

In Sect. 7 we discuss possible directions for the development of the SDE-MATH software package.
In particular, the parallelization of computations, the implementation of methods of the Runge--Kutta type 
\cite{2}, \cite{7}, \cite{40}, \cite{56} and multistep numerical methods for It\^o SDEs \cite{2}, \cite{7}, 
\cite{40}, \cite{56}, the development of a part of the software package for solving filtering problem and 
stochastic optimal control problem \cite{2}, as well as improvement of the graphical user interface.

\newpage

\section{Theoretical Results Underlying the SDE-MATH Software Package}

\subsection{Strong Numerical Methods with Convergence Orders
$0.5,$ $1.0,$ $1.5,$ $2.0,$ $2.5,$ and $3.0$  
for It\^{o} SDEs Based on the Unified Taylor--It\^{o} Expansion}

Let $(\Omega,$ ${\rm F},$ ${\sf P})$ be a complete probability space and let 
$\{{\rm F}_t, t\in[0,T]\}$ be a nondecreasing right-continuous 
family of $\sigma$-algebras of ${\rm F}.$
Let ${\bf w}_t$ be a standard $m$-dimensional Wiener stochastic 
process with independent components ${\bf w}_{t}^{(i)}$ $(i=1,\ldots,m)$, which is
${\rm F}_t$-measurable for any $t\in[0, T].$ 
Consider
an It\^{o} SDE
in the integral form
\begin{equation}
\label{1.5.2}
{\bf x}_t={\bf x}_0+\int\limits_0^t {\bf a}({\bf x}_{\tau},\tau)d\tau+
\sum\limits_{i=1}^m
\int\limits_0^t B_i({\bf x}_{\tau},\tau)d{\bf w}_{\tau}^{(i)},\ \ \
{\bf x}_0={\bf x}(0,\omega),
\end{equation}
where ${\bf x}_t\in\mathbb{R}^n$ is a strong solution of 
the It\^{o} SDE (\ref{1.5.2}),
the nonrandom functions ${\bf a}({\bf x},t): 
\mathbb{R}^n\times[0, T]\to\mathbb{R}^n$,
$B({\bf x},t): \mathbb{R}^n\times[0, T]\to\mathbb{R}^{n\times m}$
guarantee the existence and uniqueness up to stochastic equivalence 
of a strong solution
of (\ref{1.5.2}) \cite{1982}, the second integral on 
the right-hand side of (\ref{1.5.2}) is 
interpreted as an It\^{o} stochastic integral, $B_i({\bf x},t)$
is the $i$th colomn of the matrix function $B({\bf x},t),$
${\bf x}_0$ is an $n$-dimensional and ${\rm F}_0$-measurable random variable, 
${\sf M}\bigl\{\left|{\bf x}_0\right|^2\bigr\}<\infty$ 
(${\sf M}$ is an expectation operator).
We assume that
${\bf x}_0$ and ${\bf w}_t-{\bf w}_0$ are independent when $t>0.$

It is well known that one of the effective approaches 
to the numerical integration of 
It\^{o} SDEs is an approach based on the Taylor--It\^{o} and 
Taylor--Stratonovich expansions
\cite{2}, \cite{7}, \cite{40}. 
The essential feature of such 
expansions is the so-called iterated
It\^{o} and Stratonovich stochastic integrals, which have the 
form
\begin{equation}
\label{ito}
J[\psi^{(k)}]_{T,t}=\int\limits_t^T\psi_k(t_k) \ldots \int\limits_t^{t_{2}}
\psi_1(t_1) d{\bf w}_{t_1}^{(i_1)}\ldots
d{\bf w}_{t_k}^{(i_k)},
\end{equation}
\begin{equation}
\label{str}
J^{*}[\psi^{(k)}]_{T,t}=
\int\limits_t^{*T}\psi_k(t_k) \ldots \int\limits_t^{*t_{2}}
\psi_1(t_1) d{\bf w}_{t_1}^{(i_1)}\ldots
d{\bf w}_{t_k}^{(i_k)},
\end{equation}
where every $\psi_l(\tau)$ $(l=1,\ldots,k)$ is a
nonrandom function 
on $[t,T],$ ${\bf w}_{\tau}^{(i)}$ ($i=1,\ldots,m$) are independent 
standard Wiener processes and
${\bf w}_{\tau}^{(0)}\stackrel{\sf def}{=}\tau,$
$$
\int\limits\ \hbox{and}\ \int\limits^{*}
$$ 
denote It\^{o} and 
Stratonovich stochastic integrals,
respectively; $i_1,\ldots,i_k = 0, 1,\ldots,m$ (see definitions in \cite{2}).

Note that $\psi_l(\tau)\equiv 1$ $(l=1,\ldots,k)$ and
$i_1,\ldots,i_k = 0, 1,\ldots,m$ in  
the classical Taylor--It\^{o} and Taylor--Stratonovich
expansions
\cite{2}. At the same time 
$\psi_l(\tau)\equiv (t-\tau)^{q_l}$ ($l=1,\ldots,k$; 
$q_1,\ldots,q_k=0, 1, 2,\ldots $) and $i_1,\ldots,i_k = 1,\ldots,m$ in
the unified Taylor--It\^{o} and Taylor--Stratonovich
expansions
\cite{24}, \cite{25}
(also see \cite{26}, \cite{26a}, \cite{58}, Chapter 4).

Let $C^{2,1}({\mathbb{R}}^n\times [0,T])$ 
be the space of continuous functions $R({\bf x}, t): 
{\mathbb{R}}^n\times [0, T] \to {\mathbb{R}}^1$ with the 
following property: these functions are twice
continuously differentiable in ${\bf x}$ and have one continuous 
derivative in $t$. Let us consider the following differential operators on
the space $C^{2,1}({\mathbb{R}}^n\times [0,T])$ 
\begin{equation}
\label{2.3}
L= {\partial \over \partial t}
+ \sum^ {n} _ {i=1} {\bf a}^{(i)} ({\bf x},  t) 
{\partial  \over  \partial  {\bf  x}^{(i)}} 
+ {1\over 2} \sum^ {m} _ {j=1} \sum^ {n} _ {l,i=1}
B^ {(lj)} ({\bf x}, t) B^ {(ij)} ({\bf x}, t) {\partial
^{2} \over\partial{\bf x}^{(l)}\partial{\bf x}^{(i)}},
\end{equation}
\begin{equation}
\label{2.4}
G^ {(i)} _ {0} = \sum^ {n} _ {j=1} B^ {(ji)} ({\bf x}, t)
{\partial  \over \partial {\bf x} ^ {(j)}},\ \ \
i=1,\ldots,m,
\end{equation}
where ${\bf a}^{(i)} ({\bf x},  t)$ is the $i$th component of  
the vector function ${\bf a}({\bf x},  t)$ and $B^ {(ij)} ({\bf x}, t)$
is the $ij$th component of the matrix function $B({\bf x}, t)$.

Consider the following sequence of differential 
operators
$$
G_p^{(i)}=\frac{1}{p}\left(
G_{p-1}^{(i)}L-LG_{p-1}^{(i)}\right),\ \ \
p=1, 2,\ldots,\ \ \ i=1,\ldots,m,
$$
where $L$ and $G_0^{(i)},$ $i=1,\ldots,m$
are defined by 
(\ref{2.3}), (\ref{2.4}).

For the further consideration, we need to introduce
the following set of iterated It\^{o} stochastic 
integrals
\begin{equation}
\label{qqq1x}
I_{(l_1\ldots l_k)s,t}^{(i_1\ldots i_k)}=
\int\limits_t^s
(t-t_{k})^{l_{k}}\ldots 
\int\limits_t^{t_2}
(t-t _ {1}) ^ {l_ {1}} d
{\bf w} ^ {(i_ {1})} _ {t_ {1}} \ldots 
d {\bf w} _ {t_ {k}} ^ {(i_ {k})},
\end{equation}
where $l_1,\ldots,l_k=0,1,\ldots$ and $i_1,\ldots,i_k=1,\ldots,m.$

Assume that 
$R({\bf x},t)$, 
${\bf a}({\bf x},t),$ and $B_i({\bf x},t),$ $i=1,\ldots,m$
are enough smooth functions with respect to the 
variables ${\bf x}$ and $t$. Then 
for all $s, t\in[0,T]$ such that $s>t$
we can write the following
unified Taylor--It\^{o} expansion
\cite{24} (also see \cite{26}, \cite{26a}, \cite{58}, Chapter 4)
$$
R({\bf x}_s,s)=
$$

\vspace{-4mm}
$$
=R({\bf x}_t,t)+\sum_{q=1}^r\sum_{(k,j,l_1,\ldots,l_k)\in
{\rm D}_q}
\frac{(s-t)^j}{j!} \sum\limits_{i_1,\ldots,i_k=1}^m
G_{l_1}^{(i_1)}\ldots G_{l_k}^{(i_k)}
L^j R({\bf x}_t,t)\
I^{(i_1\ldots i_k)}_{(l_1\ldots l_k)s,t}+
$$

\begin{equation}
\label{15.001}
+\left(H_{r+1}\right)_{s,t}\ \ \ \hbox{w.\ p.\ 1},
\end{equation}

\vspace{1mm}
\noindent
where 
$$
L^j R({\bf x},t)\stackrel{\sf def}{=}
\begin{cases}\underbrace{L\ldots L}_j
R({\bf x},t)\ &\hbox{for}\ j\ge 1\cr\cr
R({\bf x},t)\ &\hbox{for}\ j=0
\end{cases},
$$

\vspace{-2mm}
\begin{equation}
\label{asas1}
{\rm D}_{q}=\left\{
(k, j, l_ {1},\ldots, l_ {k}):\  k + 2\left(j +
\sum\limits_{p=1}^k l_p\right)= q;\ k, j, l_{1},\ldots, 
l_{k} =0,1,\ldots\right\},
\end{equation}

\noindent
and $\left(H_{r+1}\right)_{s,t}$ is the 
remainder term in integral form
\cite{26}, \cite{26a}, \cite{58}.

Consider the partition 
$\{\tau_p\}_{p=0}^N$ of the interval
$[0,T]$ such that
\begin{equation}
\label{asas2}
0=\tau_0<\tau_1<\ldots<\tau_N=T,\ \ \
\Delta_N=
\max\limits_{0\le j\le N-1}\left|\tau_{j+1}-\tau_j\right|.
\end{equation}

Let ${\bf y}_{\tau_j}\stackrel{\sf def}{=}
{\bf y}_{j},$\ $j=0, 1,\ldots,N$ be a time discrete approximation
of the process ${\bf x}_t,$ $t\in[0,T],$ which is a solution of the It\^{o}
SDE (\ref{1.5.2}). 

{\bf Definiton 1}\ \cite{2}.\
{\it We will say that a time discrete approximation 
${\bf y}_{j},$\ $j=0, 1,\ldots,N,$
corresponding to the maximal step of discretization $\Delta_N,$
converges strongly with order
$\gamma>0$ at time moment 
$T$ to the process ${\bf x}_t,$ $t\in[0,T]$,
if there exists a constant $C>0,$ which does not depend on 
$\Delta_N,$ and a $\delta>0$ such that 
${\sf M}\{|{\bf x}_T-{\bf y}_T|\}\le
C(\Delta_N)^{\gamma}$
for each $\Delta_N\in(0, \delta).$}

From (\ref{15.001}) for $s=\tau_{p+1}$ and
$t=\tau_p$ we obtain the following representation of
explicit one-step strong numerical scheme 
for the It\^{o} SDE (\ref{1.5.2})
$$
{\bf y}_{p+1}={\bf y}_{p}+
\sum_{q=1}^r\sum_{(k,j,l_1,\ldots,l_k)\in{\rm D}_q}
\frac{(\tau_{p+1}-\tau_p)^j}{j!} \sum\limits_{i_1,\ldots,i_k=1}^m
G_{l_1}^{(i_1)}\ldots G_{l_k}^{(i_k)}
L^j\hspace{0.4mm}{\bf y}_{p}\
{\hat I}^{(i_1\ldots i_k)}_{(l_1\ldots l_k)\tau_{p+1},\tau_p}+
$$
\begin{equation}
\label{15.002}
+{\bf 1}_{\{r=2d-1,
d\in \mathbb{N}\}}\frac{(\tau_{p+1}-\tau_p)^{(r+1)/2}}{\left(
(r+1)/2\right)!}L^{(r+1)/2}\hspace{0.4mm}{\bf y}_{p},
\end{equation}

\vspace{1mm}
\noindent
where
$\hat I^{(i_1\ldots i_k)}_{(l_1\ldots l_k)\tau_{p+1},\tau_p}$ 
is an approximation of the iterated It\^{o}
stochastic integral 
(\ref{qqq1x}) and ${\bf 1}_A$ is the indicator of the set $A.$
Note that we understand the equality (\ref{15.002}) componentwise
with respect to the components ${\bf y}_p^{(i)}$ of the column
${\bf y}_p.$
Also for simplicity we put 
$\tau_p=p\Delta$,
$\Delta=T/N,$ $p=0,1,\ldots,N.$

Under the appropriate conditions \cite{2}
the numerical scheme (\ref{15.002}) has strong order $r/2$
($r\in\mathbb{N}$) of convergence.

Below we consider particular cases of the numerical scheme
(\ref{15.002}) for $r=1, 2,3,4,5,$ and $6,$ i.e. 
explicit one-step strong numerical schemes 
with convergence orders $0.5, 1.0, 1.5, 2.0, 2.5,$ and $3.0$ 
for the It\^{o} SDE (\ref{1.5.2}) \cite{26}, \cite{26a}, \cite{1982a}, \cite{1982b}.
At that for simplicity 
we will write ${\bf a},$ $L {\bf a},$ 
$B_i,$ $G_0^{(i)}B_{j}$ etc.
instead of ${\bf a}({\bf y}_p,\tau_p),$ 
$L {\bf a}({\bf y}_p,\tau_p),$ 
$B_i({\bf y}_p,\tau_p),$ 
$G_0^{(i)}B_{j}({\bf y}_p,\tau_p)$ etc. correspondingly.
Moreover, the operators $L$ and $G_0^{(i)},$ $i=1,\ldots,m$
are defined by 
(\ref{2.3}), (\ref{2.4}).

\vspace{6mm}

\centerline{\bf Scheme with strong order 0.5 (Euler scheme)}
\begin{equation}
\label{sm1}
{\bf y}_{p+1}={\bf y}_{p}+\sum_{i_{1}=1}^{m}B_{i_{1}}
\hat I_{(0)\tau_{p+1},\tau_p}^{(i_{1})}+\Delta{\bf a}.
\end{equation}

\vspace{5mm}

\centerline{\bf Scheme with strong order 1.0 (Milstein scheme)}
\begin{equation}
\label{al1}
{\bf y}_{p+1}={\bf y}_{p}+\sum_{i_{1}=1}^{m}B_{i_{1}}
\hat I_{(0)\tau_{p+1},\tau_p}^{(i_{1})}+\Delta{\bf a}
+\sum_{i_{1},i_{2}=1}^{m}G_0^{(i_{1})}
B_{i_{2}}\hat I_{(00)\tau_{p+1},\tau_p}^{(i_{1}i_{2})}.
\end{equation}

\vspace{5mm}

\centerline{\bf Scheme with strong order 1.5}
$$
{\bf y}_{p+1}={\bf y}_{p}+\sum_{i_{1}=1}^{m}B_{i_{1}}
\hat I_{(0)\tau_{p+1},\tau_p}^{(i_{1})}+\Delta{\bf a}
+\sum_{i_{1},i_{2}=1}^{m}G_0^{(i_{1})}
B_{i_{2}}\hat I_{(00)\tau_{p+1},\tau_p}^{(i_{1}i_{2})}+
$$
$$
+
\sum_{i_{1}=1}^{m}\left[G_0^{(i_{1})}{\bf a}\left(
\Delta \hat I_{(0)\tau_{p+1},\tau_p}^{(i_{1})}+
\hat I_{(1)\tau_{p+1},\tau_p}^{(i_{1})}\right)
- LB_{i_{1}}\hat I_{(1)\tau_{p+1},\tau_p}^{(i_{1})}\right]+
$$
\begin{equation}
\label{al2}
+\sum_{i_{1},i_{2},i_{3}=1}^{m} G_0^{(i_{1})}G_0^{(i_{2})}
B_{i_{3}}\hat I_{(000)\tau_{p+1},\tau_p}^{(i_{1}i_{2}i_{3})}+
\frac{\Delta^2}{2}L{\bf a}.
\end{equation}

\newpage
\noindent
\centerline{\bf Scheme with strong order 2.0}

$$
{\bf y}_{p+1}={\bf y}_{p}+\sum_{i_{1}=1}^{m}B_{i_{1}}
\hat I_{(0)\tau_{p+1},\tau_p}^{(i_{1})}+\Delta{\bf a}
+\sum_{i_{1},i_{2}=1}^{m}G_0^{(i_{1})}
B_{i_{2}}\hat I_{(00)\tau_{p+1},\tau_p}^{(i_{1}i_{2})}+
$$
$$
+
\sum_{i_{1}=1}^{m}\left[G_0^{(i_{1})}{\bf a}\left(
\Delta \hat I_{(0)\tau_{p+1},\tau_p}^{(i_{1})}+
\hat I_{(1)\tau_{p+1},\tau_p}^{(i_{1})}\right)
- LB_{i_{1}}\hat I_{(1)\tau_{p+1},\tau_p}^{(i_{1})}\right]+
$$
$$
+\sum_{i_{1},i_{2},i_{3}=1}^{m} G_0^{(i_{1})}G_0^{(i_{2})}
B_{i_{3}}\hat I_{(000)\tau_{p+1},\tau_p}^{(i_{1}i_{2}i_{3})}+
\frac{\Delta^2}{2} L{\bf a}+
$$
$$
+\sum_{i_{1},i_{2}=1}^{m}
\left[G_0^{(i_{1})} LB_{i_{2}}\left(
\hat I_{(10)\tau_{p+1},\tau_p}^{(i_{1}i_{2})}-
\hat I_{(01)\tau_{p+1},\tau_p}^{(i_{1}i_{2})}
\right)
-LG_0^{(i_{1})}
B_{i_{2}}\hat I_{(10)\tau_{p+1},\tau_p}^{(i_{1}i_{2})}
+\right.
$$
$$
\left.+G_0^{(i_{1})}G_0^{(i_{2})}{\bf a}\left(
\hat I_{(01)\tau_{p+1},\tau_p}
^{(i_{1}i_{2})}+\Delta \hat I_{(00)\tau_{p+1},\tau_p}^{(i_{1}i_{2})}
\right)\right]+
$$
\begin{equation}
\label{al3}
+
\sum_{i_{1},i_{2},i_{3},i_{4}=1}^{m}G_0^{(i_{1})}G_0^{(i_{2})}G_0^{(i_{3})}
B_{i_{4}}\hat I_{(0000)\tau_{p+1},\tau_p}^{(i_{1}i_{2}i_{3}i_{4})}.
\end{equation}

\vspace{6mm}

\centerline{\bf Scheme with strong order 2.5}

$$
{\bf y}_{p+1}={\bf y}_{p}+\sum_{i_{1}=1}^{m}B_{i_{1}}
\hat I_{(0)\tau_{p+1},\tau_p}^{(i_{1})}+\Delta{\bf a}
+\sum_{i_{1},i_{2}=1}^{m}G_0^{(i_{1})}
B_{i_{2}}\hat I_{(00)\tau_{p+1},\tau_p}^{(i_{1}i_{2})}+
$$
$$
+
\sum_{i_{1}=1}^{m}\left[G_0^{(i_{1})}{\bf a}\left(
\Delta \hat I_{(0)\tau_{p+1},\tau_p}^{(i_{1})}+
\hat I_{(1)\tau_{p+1},\tau_p}^{(i_{1})}\right)
- LB_{i_{1}}\hat I_{(1)\tau_{p+1},\tau_p}^{(i_{1})}\right]+
$$
$$
+\sum_{i_{1},i_{2},i_{3}=1}^{m} G_0^{(i_{1})}G_0^{(i_{2})}
B_{i_{3}}\hat I_{(000)\tau_{p+1},\tau_p}^{(i_{1}i_{2}i_{3})}+
\frac{\Delta^2}{2} L{\bf a}+
$$
$$
+\sum_{i_{1},i_{2}=1}^{m}
\left[G_0^{(i_{1})} LB_{i_{2}}\left(
\hat I_{(10)\tau_{p+1},\tau_p}^{(i_{1}i_{2})}-
\hat I_{(01)\tau_{p+1},\tau_p}^{(i_{1}i_{2})}
\right)
- LG_0^{(i_{1})}
B_{i_{2}}\hat I_{(10)\tau_{p+1},\tau_p}^{(i_{1}i_{2})}
+\right.
$$
$$
\left.+G_0^{(i_{1})}G_0^{(i_{2})}{\bf a}\left(
\hat I_{(01)\tau_{p+1},\tau_p}
^{(i_{1}i_{2})}+\Delta \hat I_{(00)\tau_{p+1},\tau_p}^{(i_{1}i_{2})}
\right)\right]+
$$

\vspace{-1mm}
$$
+
\sum_{i_{1},i_{2},i_{3},i_{4}=1}^{m}G_0^{(i_{1})}G_0^{(i_{2})}G_0^{(i_{3})}
B_{i_{4}}\hat I_{(0000)\tau_{p+1},\tau_p}^{(i_{1}i_{2}i_{3}i_{4})}+
$$
$$
+\sum_{i_{1}=1}^{m}\Biggl[G_0^{(i_{1})} L{\bf a}\left(\frac{1}{2}
\hat I_{(2)\tau_{p+1},\tau_p}
^{(i_{1})}+\Delta \hat I_{(1)\tau_{p+1},\tau_p}^{(i_{1})}+
\frac{\Delta^2}{2}\hat I_{(0)\tau_{p+1},\tau_p}^{(i_{1})}\right)\Biggr.+
$$
$$
+\frac{1}{2} L LB_{i_{1}}\hat I_{(2)\tau_{p+1},\tau_p}^{(i_{1})}-
LG_0^{(i_{1})}{\bf a}\Biggl.
\left(\hat I_{(2)\tau_{p+1},\tau_p}^{(i_{1})}+
\Delta \hat I_{(1)\tau_{p+1},\tau_p}^{(i_{1})}\right)\Biggr]+
$$
$$
+
\sum_{i_{1},i_{2},i_{3}=1}^m\left[
G_0^{(i_{1})} LG_0^{(i_{2})}B_{i_{3}}
\left(\hat I_{(100)\tau_{p+1},\tau_p}
^{(i_{1}i_{2}i_{3})}-\hat I_{(010)\tau_{p+1},\tau_p}
^{(i_{1}i_{2}i_{3})}\right)
\right.+
$$
$$
+G_0^{(i_{1})}G_0^{(i_{2})} LB_{i_{3}}\left(
\hat I_{(010)\tau_{p+1},\tau_p}^{(i_{1}i_{2}i_{3})}-
\hat I_{(001)\tau_{p+1},\tau_p}^{(i_{1}i_{2}i_{3})}\right)+
$$

\vspace{-1mm}
$$
+G_0^{(i_{1})}G_0^{(i_{2})}G_0^{(i_{3})} {\bf a}
\left(\Delta \hat I_{(000)\tau_{p+1},\tau_p}^{(i_{1}i_{2}i_{3})}+
\hat I_{(001)\tau_{p+1},\tau_p}^{(i_{1}i_{2}i_{3})}\right)-
$$

\vspace{-1mm}
$$
\left.- LG_0^{(i_{1})}G_0^{(i_{2})}B_{i_{3}}
\hat I_{(100)\tau_{p+1},\tau_p}^{(i_{1}i_{2}i_{3})}\right]+
$$
$$
+\sum_{i_{1},i_{2},i_{3},i_{4},i_{5}=1}^m
G_0^{(i_{1})}G_0^{(i_{2})}G_0^{(i_{3})}G_0^{(i_{4})}B_{i_{5}}
\hat I_{(00000)\tau_{p+1},\tau_p}^{(i_{1}i_{2}i_{3}i_{4}i_{5})}+
$$
\begin{equation}
\label{al4}
+
\frac{\Delta^3}{6}LL{\bf a}.
\end{equation}

\vspace{4mm}

\centerline{\bf Scheme with strong order 3.0}
$$
{\bf y}_{p+1}={\bf y}_{p}+\sum_{i_{1}=1}^{m}B_{i_{1}}
\hat I_{(0)\tau_{p+1},\tau_p}^{(i_{1})}+\Delta{\bf a}
+\sum_{i_{1},i_{2}=1}^{m}G_0^{(i_{1})}
B_{i_{2}}\hat I_{(00)\tau_{p+1},\tau_p}^{(i_{1}i_{2})}+
$$
$$
+
\sum_{i_{1}=1}^{m}\left[G_0^{(i_{1})}{\bf a}\left(
\Delta \hat I_{(0)\tau_{p+1},\tau_p}^{(i_{1})}+
\hat I_{(1)\tau_{p+1},\tau_p}^{(i_{1})}\right)
- LB_{i_{1}}\hat I_{(1)\tau_{p+1},\tau_p}^{(i_{1})}\right]+
$$
$$
+\sum_{i_{1},i_{2},i_{3}=1}^{m} G_0^{(i_{1})}G_0^{(i_{2})}
B_{i_{3}}\hat I_{(000)\tau_{p+1},\tau_p}^{(i_{1}i_{2}i_{3})}+
\frac{\Delta^2}{2} L{\bf a}+
$$
$$
+\sum_{i_{1},i_{2}=1}^{m}
\left[G_0^{(i_{1})} LB_{i_{2}}\left(
\hat I_{(10)\tau_{p+1},\tau_p}^{(i_{1}i_{2})}-
\hat I_{(01)\tau_{p+1},\tau_p}^{(i_{1}i_{2})}
\right)
- LG_0^{(i_{1})}
B_{i_{2}}\hat I_{(10)\tau_{p+1},\tau_p}^{(i_{1}i_{2})}
+\right.
$$
$$
\left.+G_0^{(i_{1})}G_0^{(i_{2})}{\bf a}\left(
\hat I_{(01)\tau_{p+1},\tau_p}
^{(i_{1}i_{2})}+\Delta \hat I_{(00)\tau_{p+1},\tau_p}^{(i_{1}i_{2})}
\right)\right]+
$$
\begin{equation}
\label{al5}
+
\sum_{i_{1},i_{2},i_{3},i_{4}=1}^{m}G_0^{(i_{1})}G_0^{(i_{2})}G_0^{(i_{3})}
B_{i_{4}}\hat I_{(0000)\tau_{p+1},\tau_p}^{(i_{1}i_{2}i_{3}i_{4})}+
{\bf q}_{p+1,p}+{\bf r}_{p+1,p},
\end{equation}

\noindent
where
$$
{\bf q}_{p+1,p}=
\sum_{i_{1}=1}^{m}\Biggl[G_0^{(i_{1})} L{\bf a}\left(\frac{1}{2}
\hat I_{(2)\tau_{p+1},\tau_p}
^{(i_{1})}+\Delta \hat I_{(1)\tau_{p+1},\tau_p}^{(i_{1})}+
\frac{\Delta^2}{2}\hat I_{(0)\tau_{p+1},\tau_p}^{(i_{1})}\right)\Biggr.+
$$
$$
+\frac{1}{2} L LB_{i_{1}}\hat I_{(2)\tau_{p+1},\tau_p}^{(i_{1})}-
LG_0^{(i_{1})}{\bf a}\Biggl.
\left(\hat I_{(2)\tau_{p+1},\tau_p}^{(i_{1})}+
\Delta \hat I_{(1)\tau_{p+1},\tau_p}^{(i_{1})}\right)\Biggr]+
$$
$$
+
\sum_{i_{1},i_{2},i_{3}=1}^m\left[
G_0^{(i_{1})} LG_0^{(i_{2})}B_{i_{3}}
\left(\hat I_{(100)\tau_{p+1},\tau_p}
^{(i_{1}i_{2}i_{3})}-\hat I_{(010)\tau_{p+1},\tau_p}
^{(i_{1}i_{2}i_{3})}\right)
\right.+
$$
$$
+G_0^{(i_{1})}G_0^{(i_{2})} LB_{i_{3}}\left(
\hat I_{(010)\tau_{p+1},\tau_p}^{(i_{1}i_{2}i_{3})}-
\hat I_{(001)\tau_{p+1},\tau_p}^{(i_{1}i_{2}i_{3})}\right)+
$$

\vspace{-1mm}
$$
+G_0^{(i_{1})}G_0^{(i_{2})}G_0^{(i_{3})} {\bf a}
\left(\Delta \hat I_{(000)\tau_{p+1},\tau_p}^{(i_{1}i_{2}i_{3})}+
\hat I_{(001)\tau_{p+1},\tau_p}^{(i_{1}i_{2}i_{3})}\right)-
$$

\vspace{-1mm}
$$
\left.- LG_0^{(i_{1})}G_0^{(i_{2})}B_{i_{3}}
\hat I_{(100)\tau_{p+1},\tau_p}^{(i_{1}i_{2}i_{3})}\right]+
$$
$$
+\sum_{i_{1},i_{2},i_{3},i_{4},i_{5}=1}^m
G_0^{(i_{1})}G_0^{(i_{2})}G_0^{(i_{3})}G_0^{(i_{4})}B_{i_{5}}
\hat I_{(00000)\tau_{p+1},\tau_p}^{(i_{1}i_{2}i_{3}i_{4}i_{5})}+
$$
$$
+
\frac{\Delta^3}{6}LL {\bf a},
$$
\noindent
and
$$
{\bf r}_{p+1,p}=\sum_{i_{1},i_{2}=1}^{m}
\Biggl[G_0^{(i_{1})}G_0^{(i_{2})} L {\bf a}\Biggl(
\frac{1}{2}\hat I_{(02)\tau_{p+1},\tau_p}^{(i_{1}i_{2})}
+
\Delta \hat I_{(01)\tau_{p+1},\tau_p}^{(i_{1}i_{2})}
+
\frac{\Delta^2}{2}
\hat I_{(00)\tau_{p+1},\tau_p}^{(i_{1}i_{2})}\Biggr)+\Biggr.
$$
$$
+
\frac{1}{2} L LG_0^{(i_{1})}B_{i_{2}}
\hat I_{(20)\tau_{p+1},\tau_p}^{(i_{1}i_{2})}+
$$

\vspace{-3mm}
$$
+G_0^{(i_{1})} LG_0^{(i_{2})} {\bf a}\left(
\hat I_{(11)\tau_{p+1},\tau_p}
^{(i_{1}i_{2})}-\hat I_{(02)\tau_{p+1},\tau_p}^{(i_{1}i_{2})}+
\Delta\left(\hat I_{(10)\tau_{p+1},\tau_p}
^{(i_{1}i_{2})}-\hat I_{(01)\tau_{p+1},\tau_p}^{(i_{1}i_{2})}
\right)\right)+
$$

\vspace{-1mm}
$$
+ LG_0^{(i_{1})} LB_{i_2}\left(
\hat I_{(11)\tau_{p+1},\tau_p}
^{(i_{1}i_{2})}-\hat I_{(20)\tau_{p+1},\tau_p}^{(i_{1}i_{2})}\right)+
$$
$$
+G_0^{(i_{1})} L LB_{i_2}\Biggl(
\frac{1}{2}\hat I_{(02)\tau_{p+1},\tau_p}^{(i_{1}i_{2})}+
\frac{1}{2}\hat I_{(20)\tau_{p+1},\tau_p}^{(i_{1}i_{2})}-
\hat I_{(11)\tau_{p+1},\tau_p}^{(i_{1}i_{2})}\Biggr)-
$$
$$
\Biggl.- LG_0^{(i_{1})}G_0^{(i_{2})}{\bf a}\left(
\Delta \hat I_{(10)\tau_{p+1},\tau_p}
^{(i_{1}i_{2})}+\hat I_{(11)\tau_{p+1},\tau_p}^{(i_{1}i_{2})}\right)
\Biggr]+
$$
$$
+
\sum_{i_{1},i_2,i_3,i_{4}=1}^m\Biggl[
G_0^{(i_{1})}G_0^{(i_{2})}G_0^{(i_{3})}G_0^{(i_{4})}{\bf a}
\left(\Delta \hat I_{(0000)\tau_{p+1},\tau_p}
^{(i_1i_{2}i_{3}i_{4})}+\hat I_{(0001)\tau_{p+1},\tau_p}
^{(i_1i_{2}i_{3}i_{4})}\right)
+\Biggr.
$$
$$
+G_0^{(i_{1})}G_0^{(i_{2})} LG_0^{(i_{3})}B_{i_4}
\left(\hat I_{(0100)\tau_{p+1},\tau_p}
^{(i_1i_{2}i_{3}i_{4})}-\hat I_{(0010)\tau_{p+1},\tau_p}
^{(i_1i_{2}i_{3}i_{4})}\right)-
$$

\vspace{-1mm}
$$
- LG_0^{(i_{1})}G_0^{(i_{2})}G_0^{(i_{3})}B_{i_4}
\hat I_{(1000)\tau_{p+1},\tau_p}
^{(i_1i_{2}i_{3}i_{4})}+
$$

\vspace{-1mm}
$$
+G_0^{(i_{1})} LG_0^{(i_{2})}G_0^{(i_{3})}B_{i_4}
\left(\hat I_{(1000)\tau_{p+1},\tau_p}
^{(i_1i_{2}i_{3}i_{4})}-\hat I_{(0100)\tau_{p+1},\tau_p}
^{(i_1i_{2}i_{3}i_{4})}\right)+
$$
$$
\Biggl.+G_0^{(i_{1})}G_0^{(i_{2})}G_0^{(i_{3})}LB_{i_4}
\left(\hat I_{(0010)\tau_{p+1},\tau_p}
^{(i_1i_{2}i_{3}i_{4})}-\hat I_{(0001)\tau_{p+1},\tau_p}
^{(i_1i_{2}i_{3}i_{4})}\right)\Biggr]+
$$
$$
+\sum_{i_{1},i_2,i_3,i_4,i_5,i_{6}=1}^m
G_0^{(i_{1})}G_0^{(i_{2})}
G_0^{(i_{3})}G_0^{(i_{4})}G_0^{(i_{5})}B_{i_{6}}
\hat I_{(000000)\tau_{p+1},\tau_p}^{(i_1i_{2}i_{3}i_{4}i_{5}i_{6})}.
$$

\vspace{2mm}

Under the suitable conditions \cite{2}
the numerical schemes (\ref{al1})--(\ref{al5}) 
have strong orders 1.0, 1.5, 2.0, 2.5, and 3.0 of convergence
correspondingly.
Among these conditions we consider only the condition
for approximations of iterated It\^{o} stochastic 
integrals from 
(\ref{al1})--(\ref{al5}) \cite{2} (also see \cite{40})
\begin{equation}
\label{uslov}
{\sf M}\left\{\Biggl(I_{(l_{1}\ldots l_{k})\tau_{p+1},\tau_p}
^{(i_{1}\ldots i_{k})} 
-\hat I_{(l_{1}\ldots l_{k})\tau_{p+1},\tau_p}^{(i_{1}\ldots i_{k})}
\Biggr)^2\right\}\le C\Delta^{r+1},
\end{equation}
where constant $C$ is independent of $\Delta$ and
$r/2$ are the strong convergence orders for the numerical schemes
(\ref{al1})--(\ref{al5}), i.e. $r/2=1.0, 1.5,$ $2.0, 2.5,$ and $3.0.$

Note that
the numerical schemes (\ref{al1})--(\ref{al5})
are unrealizable in practice without 
procedures for the numerical simulation 
of iterated It\^{o} stochastic integrals
from (\ref{15.002}).
In Sect.~2.3
we give a brief overview of the effective method
of the mean-square approximation of
iterated It\^{o} and Stratonovich stochastic integrals
of arbitrary multiplicity $k$ ($k\in\mathbb{N}$).

\subsection{Strong Numerical Methods with Convergence Orders
$1.0,$ $1.5,$ $2.0,$ $2.5,$ and $3.0$  
for It\^{o} SDEs Based on the Unified Taylor--Stra\-to\-no\-vich Expansion}

Let us consider the following differential operator on
the space $C^{2,1}({\mathbb{R}}^n\times [0,T])$ 

\newpage
\noindent
\begin{equation}
\label{2.4a}
{\bar L}=L-
\frac{1}{2}\sum^m_{i=1}G_0^{(i)}G_0^{(i)},
\end{equation}
where operators $L$ and $G_0^{(i)},$ $i=1,\ldots,m$
are defined by 
(\ref{2.3}), (\ref{2.4}).

Define the following sequence of differential 
operators
\begin{equation}
\label{a9x}
\bar G_p^{(i)}\stackrel{\sf def}{=}\frac{1}{p}\left(
\bar G_{p-1}^{(i)}\bar L-\bar L\bar G_{p-1}^{(i)}\right),\ \ \
p=1, 2,\ldots,\ \ \ i=1,\ldots,m,
\end{equation}
where $\bar G_0^{(i)}\stackrel{\sf def}{=}G_0^{(i)},$
$i=1,\ldots,m.$
The operators $\bar L$ and $G_0^{(i)},$ $i=1,\ldots,m$
are defined by 
(\ref{2.4a}) and (\ref{2.4}) correspondingly.

For the further consideration, we need to introduce
the following set of iterated Stratonovich stochastic 
integrals
\begin{equation}
\label{qqq1xx}
I_{(l_1\ldots l_k)s,t}^{*(i_1\ldots i_k)}=
\int\limits_t^{*s}
(t-t_{k})^{l_{k}}\ldots 
\int\limits_t^{* t_2}
(t-t _ {1}) ^ {l_ {1}} d
{\bf w} ^ {(i_ {1})} _ {t_ {1}} \ldots 
d {\bf w} _ {t_ {k}} ^ {(i_ {k})},
\end{equation}
where $l_1,\ldots,l_k=0,1,\ldots$ and $i_1,\ldots,i_k=1,\ldots,m.$

Assume that 
$R({\bf x},t)$, 
${\bf a}({\bf x},t),$ and $B_i({\bf x},t),$ $i=1,\ldots,m$
are enough smooth functions with respect to the 
variables ${\bf x}$ and $t$. Then 
for all $s, t\in[0,T]$ such that $s>t$
we can write the following
unified Taylor--Stratonovich expansion
\cite{25} (also see \cite{26}, \cite{26a}, \cite{58}, Chapter 4) 
$$
R({\bf x}_s,s)=
$$

\vspace{-6mm}
$$
=R({\bf x}_t,t)+\sum_{q=1}^r\sum_{(k,j,l_1,\ldots,l_k)\in
{\rm D}_q}
\frac{(s-t)^j}{j!} \sum\limits_{i_1,\ldots,i_k=1}^m
\bar G_{l_1}^{(i_1)}\ldots \bar G_{l_k}^{(i_k)}
\bar L^j R({\bf x}_t,t)\
I^{*(i_1\ldots i_k)}_{(l_1\ldots l_k)s,t}+
$$

\vspace{-2mm}
\begin{equation}
\label{15.001x}
+\left(\bar H_{r+1}\right)_{s,t}\ \ \ \hbox{w.\ p.\ 1},
\end{equation}

\vspace{2mm}
\noindent
where 
$$
{\bar L}^j R({\bf x},t)\stackrel{\sf def}{=}
\begin{cases}\underbrace{\bar L\ldots \bar L}_j
R({\bf x},t)\ &\hbox{for}\ j\ge 1\cr\cr
R({\bf x},t)\ &\hbox{for}\ j=0
\end{cases},
$$

\vspace{2mm}
\noindent
the set ${\rm D}_{q}$ is defined by the equality (\ref{asas1})
and $\left(\bar H_{r+1}\right)_{s,t}$ is
the remainder term in integral form 
\cite{25} (also see \cite{26}, \cite{26a}, \cite{58}, Chapter 4).

Consider the partition (\ref{asas2})
of the interval
$[0,T].$ 
From (\ref{15.001x}) for $s=\tau_{p+1}$ and
$t=\tau_p$ we obtain the following representation of
explicit one-step strong numerical scheme 
for the It\^{o} SDE (\ref{1.5.2})
$$
{\bf y}_{p+1}={\bf y}_{p}+
\sum_{q=1}^r\sum_{(k,j,l_1,\ldots,l_k)\in{\rm D}_q}
\frac{(\tau_{p+1}-\tau_p)^j}{j!} \sum\limits_{i_1,\ldots,i_k=1}^m
\bar G_{l_1}^{(i_1)}\ldots \bar G_{l_k}^{(i_k)}
\bar L^j\hspace{0.4mm}{\bf y}_{p}\
{\hat I}^{*(i_1\ldots i_k)}_{(l_1\ldots l_k)\tau_{p+1},\tau_p}+
$$
\begin{equation}
\label{15.002x}
+{\bf 1}_{\{r=2d-1,
d\in N\}}\frac{(\tau_{p+1}-\tau_p)^{(r+1)/2}}{\left(
(r+1)/2\right)!}L^{(r+1)/2}\hspace{0.4mm}{\bf y}_{p},
\end{equation}

\vspace{2mm}
\noindent
where
$\hat I^{*(i_1\ldots i_k)}_{(l_1\ldots l_k)\tau_{p+1},\tau_p}$ 
is an approximation of the iterated Stratonovich
stochastic integral (\ref{qqq1xx}) and ${\bf 1}_A$ is the indicator of the set $A.$
Note that we understand the equality (\ref{15.002x}) componentwise
with respect to the components ${\bf y}_p^{(i)}$ of the column
${\bf y}_p.$
Also for simplicity we put 
$\tau_p=p\Delta$,
$\Delta=T/N,$ $p=0,1,\ldots,N.$

Under the appropriate conditions \cite{2}
the numerical scheme (\ref{15.002x}) has strong order $r/2$
($r\in\mathbb{N}$) of convergence.

Denote
$$
\bar{\bf a}({\bf x},t)={\bf a}({\bf x},t)-
\frac{1}{2}\sum\limits_{j=1}^m G_0^{(j)}B_j({\bf x},t),
$$

\noindent
where $B_j({\bf x},t)$ is the $j$th column of the matrix
function $B({\bf x},t).$

It is not difficult to show that (see (\ref{2.4a}))
\begin{equation}
\label{2.4xxx}
{\bar L}=
\frac{\partial }{\partial t}+
\sum\limits_{i=1}^n \bar {\bf a}^{(i)}({\bf x},t)
\frac{\partial }{\partial {\bf x}^{(i)}},
\end{equation}

\noindent
where $\bar {\bf a}^{(i)}({\bf x},t)$ is the $i$th component of the vector
function $\bar {\bf a}({\bf x},t).$

Below we consider particular cases of the numerical scheme
(\ref{15.002x}) for $r=2,3,4,5,$ and $6,$ i.e. 
explicit one-step strong numerical schemes 
with convergence orders $1.0, 1.5, 2.0, 2.5,$ and $3.0$
for the It\^{o} SDE (\ref{1.5.2}) \cite{26}, \cite{26a}, \cite{1982c}, \cite{58}.
At that for simplicity 
we will write $\bar{\bf a},$ $\bar L\bar {\bf a},$ $L{\bf a},$
$B_i,$ $G_0^{(i)}B_{j}$ etc.
instead of $\bar{\bf a}({\bf y}_p,\tau_p),$ 
$\bar L \bar {\bf a}({\bf y}_p,\tau_p),$ 
$L{\bf a}({\bf y}_p,\tau_p),$  $B_i({\bf y}_p,\tau_p),$ 
$G_0^{(i)}B_{j}({\bf y}_p,\tau_p)$ etc. correspondingly.

\vspace{6mm}

\centerline{\bf Scheme with strong order 1.0}
\begin{equation}
\label{al1x}
{\bf y}_{p+1}={\bf y}_{p}+\sum_{i_{1}=1}^{m}B_{i_{1}}
\hat I_{(0)\tau_{p+1},\tau_p}^{*(i_{1})}+\Delta\bar{\bf a}
+\sum_{i_{1},i_{2}=1}^{m}G_0^{(i_{1})}
B_{i_{2}}\hat I_{(00)\tau_{p+1},\tau_p}^{*(i_{1}i_{2})}.
\end{equation}

\newpage
\noindent
\centerline{\bf Scheme with strong order 1.5}

$$
{\bf y}_{p+1}={\bf y}_{p}+\sum_{i_{1}=1}^{m}B_{i_{1}}
\hat I_{(0)\tau_{p+1},\tau_p}^{*(i_{1})}+\Delta\bar{\bf a}
+\sum_{i_{1},i_{2}=1}^{m}G_0^{(i_{1})}
B_{i_{2}}\hat I_{(00)\tau_{p+1},\tau_p}^{*(i_{1}i_{2})}+
$$
$$
+
\sum_{i_{1}=1}^{m}\left[G_0^{(i_{1})}\bar{\bf a}\left(
\Delta \hat I_{(0)\tau_{p+1},\tau_p}^{*(i_{1})}+
\hat I_{(1)\tau_{p+1},\tau_p}^{*(i_{1})}\right)
-\bar LB_{i_{1}}\hat I_{(1)\tau_{p+1},\tau_p}^{*(i_{1})}\right]+
$$
\begin{equation}
\label{al2x}
+\sum_{i_{1},i_{2},i_{3}=1}^{m} G_0^{(i_{1})}G_0^{(i_{2})}
B_{i_{3}}\hat I_{(000)\tau_{p+1},\tau_p}^{*(i_{1}i_{2}i_{3})}+
\frac{\Delta^2}{2}L{\bf a}.
\end{equation}

\vspace{6mm}

\centerline{\bf Scheme with strong order 2.0}

$$
{\bf y}_{p+1}={\bf y}_{p}+\sum_{i_{1}=1}^{m}B_{i_{1}}
\hat I_{(0)\tau_{p+1},\tau_p}^{*(i_{1})}+\Delta\bar{\bf a}
+\sum_{i_{1},i_{2}=1}^{m}G_0^{(i_{1})}
B_{i_{2}}\hat I_{(00)\tau_{p+1},\tau_p}^{*(i_{1}i_{2})}+
$$
$$
+
\sum_{i_{1}=1}^{m}\left[G_0^{(i_{1})}\bar{\bf a}\left(
\Delta \hat I_{(0)\tau_{p+1},\tau_p}^{*(i_{1})}+
\hat I_{(1)\tau_{p+1},\tau_p}^{*(i_{1})}\right)
-\bar LB_{i_{1}}\hat I_{(1)\tau_{p+1},\tau_p}^{*(i_{1})}\right]+
$$
$$
+\sum_{i_{1},i_{2},i_{3}=1}^{m} G_0^{(i_{1})}G_0^{(i_{2})}
B_{i_{3}}\hat I_{(000)\tau_{p+1},\tau_p}^{*(i_{1}i_{2}i_{3})}+
\frac{\Delta^2}{2}\bar L\bar{\bf a}+
$$
$$
+\sum_{i_{1},i_{2}=1}^{m}
\left[G_0^{(i_{1})}\bar LB_{i_{2}}\left(
\hat I_{(10)\tau_{p+1},\tau_p}^{*(i_{1}i_{2})}-
\hat I_{(01)\tau_{p+1},\tau_p}^{*(i_{1}i_{2})}
\right)
-\bar LG_0^{(i_{1})}
B_{i_{2}}\hat I_{(10)\tau_{p+1},\tau_p}^{*(i_{1}i_{2})}
+\right.
$$
$$
\left.+G_0^{(i_{1})}G_0^{(i_{2})}\bar{\bf a}\left(
\hat I_{(01)\tau_{p+1},\tau_p}
^{*(i_{1}i_{2})}+\Delta \hat I_{(00)\tau_{p+1},\tau_p}^{*(i_{1}i_{2})}
\right)\right]+
$$
\begin{equation}
\label{al3x}
+
\sum_{i_{1},i_{2},i_{3},i_{4}=1}^{m}G_0^{(i_{1})}G_0^{(i_{2})}G_0^{(i_{3})}
B_{i_{4}}\hat I_{(0000)\tau_{p+1},\tau_p}^{*(i_{1}i_{2}i_{3}i_{4})}.
\end{equation}

\vspace{6mm}

\centerline{\bf Scheme with strong order 2.5}

$$
{\bf y}_{p+1}={\bf y}_{p}+\sum_{i_{1}=1}^{m}B_{i_{1}}
\hat I_{(0)\tau_{p+1},\tau_p}^{*(i_{1})}+\Delta\bar{\bf a}
+\sum_{i_{1},i_{2}=1}^{m}G_0^{(i_{1})}
B_{i_{2}}\hat I_{(00)\tau_{p+1},\tau_p}^{*(i_{1}i_{2})}+
$$
$$
+
\sum_{i_{1}=1}^{m}\left[G_0^{(i_{1})}\bar{\bf a}\left(
\Delta \hat I_{(0)\tau_{p+1},\tau_p}^{*(i_{1})}+
\hat I_{(1)\tau_{p+1},\tau_p}^{*(i_{1})}\right)
-\bar LB_{i_{1}}\hat I_{(1)\tau_{p+1},\tau_p}^{*(i_{1})}\right]+
$$
$$
+\sum_{i_{1},i_{2},i_{3}=1}^{m} G_0^{(i_{1})}G_0^{(i_{2})}
B_{i_{3}}\hat I_{(000)\tau_{p+1},\tau_p}^{*(i_{1}i_{2}i_{3})}+
\frac{\Delta^2}{2}\bar L\bar{\bf a}+
$$
$$
+\sum_{i_{1},i_{2}=1}^{m}
\left[G_0^{(i_{1})}\bar LB_{i_{2}}\left(
\hat I_{(10)\tau_{p+1},\tau_p}^{*(i_{1}i_{2})}-
\hat I_{(01)\tau_{p+1},\tau_p}^{*(i_{1}i_{2})}
\right)
-\bar LG_0^{(i_{1})}
B_{i_{2}}\hat I_{(10)\tau_{p+1},\tau_p}^{*(i_{1}i_{2})}
+\right.
$$
$$
\left.+G_0^{(i_{1})}G_0^{(i_{2})}\bar{\bf a}\left(
\hat I_{(01)\tau_{p+1},\tau_p}
^{*(i_{1}i_{2})}+\Delta \hat I_{(00)\tau_{p+1},\tau_p}^{*(i_{1}i_{2})}
\right)\right]+
$$
$$
+
\sum_{i_{1},i_{2},i_{3},i_{4}=1}^{m}G_0^{(i_{1})}G_0^{(i_{2})}G_0^{(i_{3})}
B_{i_{4}}\hat I_{(0000)\tau_{p+1},\tau_p}^{*(i_{1}i_{2}i_{3}i_{4})}+
$$
$$
+\sum_{i_{1}=1}^{m}\Biggl[G_0^{(i_{1})}\bar L\bar{\bf a}\left(\frac{1}{2}
\hat I_{(2)\tau_{p+1},\tau_p}
^{*(i_{1})}+\Delta \hat I_{(1)\tau_{p+1},\tau_p}^{*(i_{1})}+
\frac{\Delta^2}{2}\hat I_{(0)\tau_{p+1},\tau_p}^{*(i_{1})}\right)\Biggr.+
$$
$$
+\frac{1}{2}\bar L\bar LB_{i_{1}}\hat I_{(2)\tau_{p+1},\tau_p}^{*(i_{1})}-
\bar LG_0^{(i_{1})}\bar{\bf a}\Biggl.
\left(\hat I_{(2)\tau_{p+1},\tau_p}^{*(i_{1})}+
\Delta \hat I_{(1)\tau_{p+1},\tau_p}^{*(i_{1})}\right)\Biggr]+
$$
$$
+
\sum_{i_{1},i_{2},i_{3}=1}^m\left[
G_0^{(i_{1})}\bar LG_0^{(i_{2})}B_{i_{3}}
\left(\hat I_{(100)\tau_{p+1},\tau_p}
^{*(i_{1}i_{2}i_{3})}-\hat I_{(010)\tau_{p+1},\tau_p}
^{*(i_{1}i_{2}i_{3})}\right)
\right.+
$$
$$
+G_0^{(i_{1})}G_0^{(i_{2})}\bar LB_{i_{3}}\left(
\hat I_{(010)\tau_{p+1},\tau_p}^{*(i_{1}i_{2}i_{3})}-
\hat I_{(001)\tau_{p+1},\tau_p}^{*(i_{1}i_{2}i_{3})}\right)+
$$

\vspace{-1mm}
$$
+G_0^{(i_{1})}G_0^{(i_{2})}G_0^{(i_{3})}\bar {\bf a}
\left(\Delta \hat I_{(000)\tau_{p+1},\tau_p}^{*(i_{1}i_{2}i_{3})}+
\hat I_{(001)\tau_{p+1},\tau_p}^{*(i_{1}i_{2}i_{3})}\right)-
$$

\vspace{-1mm}
$$
\left.-\bar LG_0^{(i_{1})}G_0^{(i_{2})}B_{i_{3}}
\hat I_{(100)\tau_{p+1},\tau_p}^{*(i_{1}i_{2}i_{3})}\right]+
$$
$$
+\sum_{i_{1},i_{2},i_{3},i_{4},i_{5}=1}^m
G_0^{(i_{1})}G_0^{(i_{2})}G_0^{(i_{3})}G_0^{(i_{4})}B_{i_{5}}
\hat I_{(00000)\tau_{p+1},\tau_p}^{*(i_{1}i_{2}i_{3}i_{4}i_{5})}+
$$
\begin{equation}
\label{al4x}
+
\frac{\Delta^3}{6}LL{\bf a}.
\end{equation}

\vspace{6mm}

\centerline{\bf Scheme with strong order 3.0}

\vspace{-2mm}
$$
{\bf y}_{p+1}={\bf y}_{p}+\sum_{i_{1}=1}^{m}B_{i_{1}}
\hat I_{(0)\tau_{p+1},\tau_p}^{*(i_{1})}+\Delta\bar{\bf a}
+\sum_{i_{1},i_{2}=1}^{m}G_0^{(i_{1})}
B_{i_{2}}\hat I_{(00)\tau_{p+1},\tau_p}^{*(i_{1}i_{2})}+
$$
$$
+
\sum_{i_{1}=1}^{m}\left[G_0^{(i_{1})}\bar{\bf a}\left(
\Delta \hat I_{(0)\tau_{p+1},\tau_p}^{*(i_{1})}+
\hat I_{(1)\tau_{p+1},\tau_p}^{*(i_{1})}\right)
-\bar LB_{i_{1}}\hat I_{(1)\tau_{p+1},\tau_p}^{*(i_{1})}\right]+
$$
$$
+\sum_{i_{1},i_{2},i_{3}=1}^{m} G_0^{(i_{1})}G_0^{(i_{2})}
B_{i_{3}}\hat I_{(000)\tau_{p+1},\tau_p}^{*(i_{1}i_{2}i_{3})}+
\frac{\Delta^2}{2}\bar L\bar{\bf a}+
$$
$$
+\sum_{i_{1},i_{2}=1}^{m}
\left[G_0^{(i_{1})}\bar LB_{i_{2}}\left(
\hat I_{(10)\tau_{p+1},\tau_p}^{*(i_{1}i_{2})}-
\hat I_{(01)\tau_{p+1},\tau_p}^{*(i_{1}i_{2})}
\right)
-\bar LG_0^{(i_{1})}
B_{i_{2}}\hat I_{(10)\tau_{p+1},\tau_p}^{*(i_{1}i_{2})}
+\right.
$$
$$
\left.+G_0^{(i_{1})}G_0^{(i_{2})}\bar{\bf a}\left(
\hat I_{(01)\tau_{p+1},\tau_p}
^{*(i_{1}i_{2})}+\Delta \hat I_{(00)\tau_{p+1},\tau_p}^{*(i_{1}i_{2})}
\right)\right]+
$$
\begin{equation}
\label{al5x}
+
\sum_{i_{1},i_{2},i_{3},i_{4}=1}^{m}G_0^{(i_{1})}G_0^{(i_{2})}G_0^{(i_{3})}
B_{i_{4}}\hat I_{(0000)\tau_{p+1},\tau_p}^{*(i_{1}i_{2}i_{3}i_{4})}+
{\bf q}_{p+1,p}+{\bf r}_{p+1,p},
\end{equation}

\noindent
where
$$
{\bf q}_{p+1,p}=
\sum_{i_{1}=1}^{m}\Biggl[G_0^{(i_{1})}\bar L\bar{\bf a}\left(\frac{1}{2}
\hat I_{(2)\tau_{p+1},\tau_p}
^{*(i_{1})}+\Delta \hat I_{(1)\tau_{p+1},\tau_p}^{*(i_{1})}+
\frac{\Delta^2}{2}\hat I_{(0)\tau_{p+1},\tau_p}^{*(i_{1})}\right)\Biggr.+
$$
$$
+\frac{1}{2}\bar L\bar LB_{i_{1}}\hat I_{(2)\tau_{p+1},\tau_p}^{*(i_{1})}-
\bar LG_0^{(i_{1})}\bar{\bf a}\Biggl.
\left(\hat I_{(2)\tau_{p+1},\tau_p}^{*(i_{1})}+
\Delta \hat I_{(1)\tau_{p+1},\tau_p}^{*(i_{1})}\right)\Biggr]+
$$
$$
+
\sum_{i_{1},i_{2},i_{3}=1}^m\left[
G_0^{(i_{1})}\bar LG_0^{(i_{2})}B_{i_{3}}
\left(\hat I_{(100)\tau_{p+1},\tau_p}
^{*(i_{1}i_{2}i_{3})}-\hat I_{(010)\tau_{p+1},\tau_p}
^{*(i_{1}i_{2}i_{3})}\right)
\right.+
$$
$$
+G_0^{(i_{1})}G_0^{(i_{2})}\bar LB_{i_{3}}\left(
\hat I_{(010)\tau_{p+1},\tau_p}^{*(i_{1}i_{2}i_{3})}-
\hat I_{(001)\tau_{p+1},\tau_p}^{*(i_{1}i_{2}i_{3})}\right)+
$$

\vspace{-1mm}
$$
+G_0^{(i_{1})}G_0^{(i_{2})}G_0^{(i_{3})}\bar {\bf a}
\left(\Delta \hat I_{(000)\tau_{p+1},\tau_p}^{*(i_{1}i_{2}i_{3})}+
\hat I_{(001)\tau_{p+1},\tau_p}^{*(i_{1}i_{2}i_{3})}\right)-
$$

\vspace{-1mm}
$$
\left.-\bar LG_0^{(i_{1})}G_0^{(i_{2})}B_{i_{3}}
\hat I_{(100)\tau_{p+1},\tau_p}^{*(i_{1}i_{2}i_{3})}\right]+
$$
$$
+\sum_{i_{1},i_{2},i_{3},i_{4},i_{5}=1}^m
G_0^{(i_{1})}G_0^{(i_{2})}G_0^{(i_{3})}G_0^{(i_{4})}B_{i_{5}}
\hat I_{(00000)\tau_{p+1},\tau_p}^{*(i_{1}i_{2}i_{3}i_{4}i_{5})}+
$$
$$
+
\frac{\Delta^3}{6}\bar L\bar L\bar {\bf a},
$$

\noindent
and
$$
{\bf r}_{p+1,p}=\sum_{i_{1},i_{2}=1}^{m}
\Biggl[G_0^{(i_{1})}G_0^{(i_{2})}\bar L\bar {\bf a}\Biggl(
\frac{1}{2}\hat I_{(02)\tau_{p+1},\tau_p}^{*(i_{1}i_{2})}
+
\Delta \hat I_{(01)\tau_{p+1},\tau_p}^{*(i_{1}i_{2})}
+
\frac{\Delta^2}{2}
\hat I_{(00)\tau_{p+1},\tau_p}^{*(i_{1}i_{2})}\Biggr)+\Biggr.
$$
$$
+
\frac{1}{2}\bar L\bar LG_0^{(i_{1})}B_{i_{2}}
\hat I_{(20)\tau_{p+1},\tau_p}^{*(i_{1}i_{2})}+
$$
$$
+G_0^{(i_{1})}\bar LG_0^{(i_{2})}\bar {\bf a}\left(
\hat I_{(11)\tau_{p+1},\tau_p}
^{*(i_{1}i_{2})}-\hat I_{(02)\tau_{p+1},\tau_p}^{*(i_{1}i_{2})}+
\Delta\left(\hat I_{(10)\tau_{p+1},\tau_p}
^{*(i_{1}i_{2})}-\hat I_{(01)\tau_{p+1},\tau_p}^{*(i_{1}i_{2})}
\right)\right)+
$$

\vspace{-1mm}
$$
+\bar LG_0^{(i_{1})}\bar LB_{i_2}\left(
\hat I_{(11)\tau_{p+1},\tau_p}
^{*(i_{1}i_{2})}-\hat I_{(20)\tau_{p+1},\tau_p}^{*(i_{1}i_{2})}\right)+
$$
$$
+G_0^{(i_{1})}\bar L\bar LB_{i_2}\Biggl(
\frac{1}{2}\hat I_{(02)\tau_{p+1},\tau_p}^{*(i_{1}i_{2})}+
\frac{1}{2}\hat I_{(20)\tau_{p+1},\tau_p}^{*(i_{1}i_{2})}-
\hat I_{(11)\tau_{p+1},\tau_p}^{*(i_{1}i_{2})}\Biggr)-
$$
$$
\Biggl.-\bar LG_0^{(i_{1})}G_0^{(i_{2})}\bar{\bf a}\left(
\Delta \hat I_{(10)\tau_{p+1},\tau_p}
^{*(i_{1}i_{2})}+\hat I_{(11)\tau_{p+1},\tau_p}^{*(i_{1}i_{2})}\right)
\Biggr]+
$$
$$
+
\sum_{i_{1},i_2,i_3,i_{4}=1}^m\Biggl[
G_0^{(i_{1})}G_0^{(i_{2})}G_0^{(i_{3})}G_0^{(i_{4})}\bar{\bf a}
\left(\Delta \hat I_{(0000)\tau_{p+1},\tau_p}
^{*(i_1i_{2}i_{3}i_{4})}+\hat I_{(0001)\tau_{p+1},\tau_p}
^{*(i_1i_{2}i_{3}i_{4})}\right)
+\Biggr.
$$
$$
+G_0^{(i_{1})}G_0^{(i_{2})}\bar LG_0^{(i_{3})}B_{i_4}
\left(\hat I_{(0100)\tau_{p+1},\tau_p}
^{*(i_1i_{2}i_{3}i_{4})}-\hat I_{(0010)\tau_{p+1},\tau_p}
^{*(i_1i_{2}i_{3}i_{4})}\right)-
$$

\vspace{-1mm}
$$
-\bar LG_0^{(i_{1})}G_0^{(i_{2})}G_0^{(i_{3})}B_{i_4}
\hat I_{(1000)\tau_{p+1},\tau_p}
^{*(i_1i_{2}i_{3}i_{4})}+
$$

\vspace{-1mm}
$$
+G_0^{(i_{1})}\bar LG_0^{(i_{2})}G_0^{(i_{3})}B_{i_4}
\left(\hat I_{(1000)\tau_{p+1},\tau_p}
^{*(i_1i_{2}i_{3}i_{4})}-\hat I_{(0100)\tau_{p+1},\tau_p}
^{*(i_1i_{2}i_{3}i_{4})}\right)+
$$
$$
\Biggl.+G_0^{(i_{1})}G_0^{(i_{2})}G_0^{(i_{3})}\bar LB_{i_4}
\left(\hat I_{(0010)\tau_{p+1},\tau_p}
^{*(i_1i_{2}i_{3}i_{4})}-\hat I_{(0001)\tau_{p+1},\tau_p}
^{*(i_1i_{2}i_{3}i_{4})}\right)\Biggr]+
$$
$$
+\sum_{i_{1},i_2,i_3,i_4,i_5,i_{6}=1}^m
G_0^{(i_{1})}G_0^{(i_{2})}
G_0^{(i_{3})}G_0^{(i_{4})}G_0^{(i_{5})}B_{i_{6}}
\hat I_{(000000)\tau_{p+1},\tau_p}^{*(i_1i_{2}i_{3}i_{4}i_{5}i_{6})}.
$$

\vspace{2mm}

Under the suitable conditions \cite{2}
the numerical schemes (\ref{al1x})--(\ref{al5x}) 
have strong orders 1.0, 1.5, 2.0, 2.5, and 3.0 of convergence
correspondingly.
Among these conditions we consider only the condition
for approximations of iterated Stratonovich stochastic 
integrals from (\ref{al1x})--(\ref{al5x}) \cite{2} (also see \cite{40})
\begin{equation}
\label{uslov1}
{\sf M}\left\{\Biggl(I_{(l_{1}\ldots l_{k})\tau_{p+1},\tau_p}
^{*(i_{1}\ldots i_{k})} 
-\hat I_{(l_{1}\ldots l_{k})\tau_{p+1},\tau_p}^{*(i_{1}\ldots i_{k})}
\Biggr)^2\right\}\le C\Delta^{r+1},
\end{equation}
where constant $C$ is independent of $\Delta$ and
$r/2$ are the strong convergence orders for the numerical schemes
(\ref{al1x})--(\ref{al5x}), i.e. $r/2=1.0, 1.5,$ $2.0, 2.5,$ and $3.0.$

Note that
the numerical schemes (\ref{al1x})--(\ref{al5x})
are unrealizable in practice without 
procedures for the numerical simulation 
of iterated Stratonovich stochastic integrals
from (\ref{15.002x}).
The next section is devoted 
to the effective method
of the mean-square approximation of
iterated It\^{o} and Stratonovich stochastic integrals
of arbitrary multiplicity $k$ ($k\in\mathbb{N}$).

\subsection{Method of Expansion and Approximation of
Iterated It\^{o} and Stra\-to\-no\-vich Stochastic Integrals
Based on Generalized Multiple Fourier Series}

Let us consider the effective approach to expansion of iterated It\^{o}
stochastic integrals \cite{40} (2006) (also see
\cite{26}, \cite{26a}, 
\cite{42}-\cite{42aa}, \cite{53}, \cite{56}, \cite{58}, \cite{59}).
This method is reffered to as the 
method of generalized
multiple Fourier series.

Suppose that every $\psi_l(\tau)$ $(l=1,\ldots,k)$ is a 
nonrandom function from the space $L_2([t, T])$.
Define the following function on the hypercube $[t, T]^k$
\begin{equation}
\label{ppp}
K(t_1,\ldots,t_k)=
\begin{cases}
\psi_1(t_1)\ldots \psi_k(t_k)\ \ \hbox{for}\ \ t_1<\ldots<t_k\cr\cr
0\ \ \hbox{otherwise}
\end{cases},\ \ t_1,\ldots,t_k\in[t, T],
\end{equation}
where $k\ge 2$ and 
$K(t_1)\equiv\psi_1(t_1)$ for $t_1\in[t, T].$

Suppose that $\{\phi_j(x)\}_{j=0}^{\infty}$
is a complete orthonormal system of functions in the space
$L_2([t, T])$.

The function $K(t_1,\ldots,t_k)$ (the so-called factorized Volterra--type kernel) belongs to the space
$L_2([t, T]^k).$
At this situation it is well known that the generalized 
multiple Fourier series 
of $K(t_1,\ldots,t_k)\in L_2([t, T]^k)$ is converging 
to $K(t_1,\ldots,t_k)$ in the hypercube $[t, T]^k$ in 
the mean-square sense, i.e.
$$
\hbox{\vtop{\offinterlineskip\halign{
\hfil#\hfil\cr
{\rm lim}\cr
$\stackrel{}{{}_{p_1,\ldots,p_k\to \infty}}$\cr
}} }\Biggl\Vert
K(t_1,\ldots,t_k)-
\sum_{j_1=0}^{p_1}\ldots \sum_{j_k=0}^{p_k}
C_{j_k\ldots j_1}\prod_{l=1}^{k} \phi_{j_l}(t_l)\Biggr\Vert_{L_2([t,T]^k)}=0,
$$
where
\begin{equation}
\label{ppppa}
C_{j_k\ldots j_1}=\int\limits_{[t,T]^k}
K(t_1,\ldots,t_k)\prod_{l=1}^{k}\phi_{j_l}(t_l)dt_1\ldots dt_k
\end{equation}
is the Fourier coefficient and
$$
\left\Vert f\right\Vert_{L_2([t,T]^k)}=\left(\int\limits_{[t,T]^k}
f^2(t_1,\ldots,t_k)dt_1\ldots dt_k\right)^{1/2}.
$$

Consider the partition $\{\tau_j\}_{j=0}^N$ of the interval
$[t,T]$ such that
\begin{equation}
\label{1111}
t=\tau_0<\ldots <\tau_N=T,\ \ \
\Delta_N=
\hbox{\vtop{\offinterlineskip\halign{
\hfil#\hfil\cr
{\rm max}\cr
$\stackrel{}{{}_{0\le j\le N-1}}$\cr
}} }\Delta\tau_j\to 0\ \ \hbox{if}\ \ N\to \infty,\ \ \
\Delta\tau_j=\tau_{j+1}-\tau_j.
\end{equation}

{\bf Theorem 1} \cite{40} (2006) (also see
\cite{26}, \cite{26a}, \cite{42}-\cite{42aa}, \cite{53}, \cite{56}, \cite{58}, \cite{59}).\
{\it Suppose that
every $\psi_l(\tau)$ $(l=1,\ldots, k)$ is a continuous nonrandom function on 
$[t, T]$ and
$\{\phi_j(x)\}_{j=0}^{\infty}$ is a complete orthonormal system  
of continuous functions in the space $L_2([t,T]).$ Then
$$
J[\psi^{(k)}]_{T,t}\  =\ 
\hbox{\vtop{\offinterlineskip\halign{
\hfil#\hfil\cr
{\rm l.i.m.}\cr
$\stackrel{}{{}_{p_1,\ldots,p_k\to \infty}}$\cr
}} }\sum_{j_1=0}^{p_1}\ldots\sum_{j_k=0}^{p_k}
C_{j_k\ldots j_1}\Biggl(
\prod_{l=1}^k\zeta_{j_l}^{(i_l)}\ -
\Biggr.
$$
\begin{equation}
\label{tyyy}
-\ \Biggl.
\hbox{\vtop{\offinterlineskip\halign{
\hfil#\hfil\cr
{\rm l.i.m.}\cr
$\stackrel{}{{}_{N\to \infty}}$\cr
}} }\sum_{(l_1,\ldots,l_k)\in {\rm G}_k}
\phi_{j_{1}}(\tau_{l_1})
\Delta{\bf w}_{\tau_{l_1}}^{(i_1)}\ldots
\phi_{j_{k}}(\tau_{l_k})
\Delta{\bf w}_{\tau_{l_k}}^{(i_k)}\Biggr),
\end{equation}

\vspace{4mm}
\noindent
where $J[\psi^{(k)}]_{T,t}$ is defined by {\rm (\ref{ito}),}
$$
{\rm G}_k={\rm H}_k\backslash{\rm L}_k,\ \ \
{\rm H}_k=\{(l_1,\ldots,l_k):\ l_1,\ldots,l_k=0, 1,\ldots,N-1\},
$$
$$
{\rm L}_k=\{(l_1,\ldots,l_k):\ l_1,\ldots,l_k=0, 1,\ldots,N-1;\
l_g\ne l_r\ (g\ne r);\ g, r=1,\ldots,k\},
$$

\noindent
${\rm l.i.m.}$ is a limit in the mean-square sense$,$
$i_1,\ldots,i_k=0,1,\ldots,m,$ 
\begin{equation}
\label{rr23}
\zeta_{j}^{(i)}=
\int\limits_t^T \phi_{j}(s) d{\bf w}_s^{(i)}
\end{equation} 
are independent standard Gaussian random variables
for various
$i$ or $j$ {\rm(}in the case when $i\ne 0${\rm),}
$C_{j_k\ldots j_1}$ is the Fourier coefficient {\rm(\ref{ppppa}),}
$\Delta{\bf w}_{\tau_{j}}^{(i)}=
{\bf w}_{\tau_{j+1}}^{(i)}-{\bf w}_{\tau_{j}}^{(i)}$
$(i=0, 1,\ldots,m),$
$\left\{\tau_{j}\right\}_{j=0}^{N}$ is a partition of
the interval $[t, T],$ which satisfies the condition {\rm (\ref{1111})}.
}

Note that a number of modifications and 
generalizations of Theorem 1 can be found in \cite{26}, \cite{26a}, \cite{58}
(Chapter~1).

Consider transformed particular cases of (\ref{tyyy}) for 
$k=1,\ldots,6$ 
\cite{26}, \cite{26a}, \cite{53}, \cite{56}, \cite{58}, \cite{59}
\begin{equation}
\label{a1}
J[\psi^{(1)}]_{T,t}
=\hbox{\vtop{\offinterlineskip\halign{
\hfil#\hfil\cr
{\rm l.i.m.}\cr
$\stackrel{}{{}_{p_1\to \infty}}$\cr
}} }\sum_{j_1=0}^{p_1}
C_{j_1}\zeta_{j_1}^{(i_1)},
\end{equation}
\begin{equation}
\label{leto5001}
J[\psi^{(2)}]_{T,t}
=\hbox{\vtop{\offinterlineskip\halign{
\hfil#\hfil\cr
{\rm l.i.m.}\cr
$\stackrel{}{{}_{p_1,p_2\to \infty}}$\cr
}} }\sum_{j_1=0}^{p_1}\sum_{j_2=0}^{p_2}
C_{j_2j_1}\Biggl(\zeta_{j_1}^{(i_1)}\zeta_{j_2}^{(i_2)}
-{\bf 1}_{\{i_1=i_2\ne 0\}}
{\bf 1}_{\{j_1=j_2\}}\Biggr),
\end{equation}

$$
J[\psi^{(3)}]_{T,t}=
\hbox{\vtop{\offinterlineskip\halign{
\hfil#\hfil\cr
{\rm l.i.m.}\cr
$\stackrel{}{{}_{p_1,\ldots,p_3\to \infty}}$\cr
}} }\sum_{j_1=0}^{p_1}\sum_{j_2=0}^{p_2}\sum_{j_3=0}^{p_3}
C_{j_3j_2j_1}\Biggl(
\zeta_{j_1}^{(i_1)}\zeta_{j_2}^{(i_2)}\zeta_{j_3}^{(i_3)}
-\Biggr.
$$
\begin{equation}
\label{leto5002}
\Biggl.-{\bf 1}_{\{i_1=i_2\ne 0\}}
{\bf 1}_{\{j_1=j_2\}}
\zeta_{j_3}^{(i_3)}
-{\bf 1}_{\{i_2=i_3\ne 0\}}
{\bf 1}_{\{j_2=j_3\}}
\zeta_{j_1}^{(i_1)}-
{\bf 1}_{\{i_1=i_3\ne 0\}}
{\bf 1}_{\{j_1=j_3\}}
\zeta_{j_2}^{(i_2)}\Biggr),
\end{equation}

$$
J[\psi^{(4)}]_{T,t}
=
\hbox{\vtop{\offinterlineskip\halign{
\hfil#\hfil\cr
{\rm l.i.m.}\cr
$\stackrel{}{{}_{p_1,\ldots,p_4\to \infty}}$\cr
}} }\sum_{j_1=0}^{p_1}\ldots\sum_{j_4=0}^{p_4}
C_{j_4\ldots j_1}\Biggl(
\prod_{l=1}^4\zeta_{j_l}^{(i_l)}
\Biggr.
-
$$
$$
-
{\bf 1}_{\{i_1=i_2\ne 0\}}
{\bf 1}_{\{j_1=j_2\}}
\zeta_{j_3}^{(i_3)}
\zeta_{j_4}^{(i_4)}
-
{\bf 1}_{\{i_1=i_3\ne 0\}}
{\bf 1}_{\{j_1=j_3\}}
\zeta_{j_2}^{(i_2)}
\zeta_{j_4}^{(i_4)}-
$$
$$
-
{\bf 1}_{\{i_1=i_4\ne 0\}}
{\bf 1}_{\{j_1=j_4\}}
\zeta_{j_2}^{(i_2)}
\zeta_{j_3}^{(i_3)}
-
{\bf 1}_{\{i_2=i_3\ne 0\}}
{\bf 1}_{\{j_2=j_3\}}
\zeta_{j_1}^{(i_1)}
\zeta_{j_4}^{(i_4)}-
$$
$$
-
{\bf 1}_{\{i_2=i_4\ne 0\}}
{\bf 1}_{\{j_2=j_4\}}
\zeta_{j_1}^{(i_1)}
\zeta_{j_3}^{(i_3)}
-
{\bf 1}_{\{i_3=i_4\ne 0\}}
{\bf 1}_{\{j_3=j_4\}}
\zeta_{j_1}^{(i_1)}
\zeta_{j_2}^{(i_2)}+
$$
$$
+
{\bf 1}_{\{i_1=i_2\ne 0\}}
{\bf 1}_{\{j_1=j_2\}}
{\bf 1}_{\{i_3=i_4\ne 0\}}
{\bf 1}_{\{j_3=j_4\}}
+
{\bf 1}_{\{i_1=i_3\ne 0\}}
{\bf 1}_{\{j_1=j_3\}}
{\bf 1}_{\{i_2=i_4\ne 0\}}
{\bf 1}_{\{j_2=j_4\}}+
$$
\begin{equation}
\label{leto5003}
+\Biggl.
{\bf 1}_{\{i_1=i_4\ne 0\}}
{\bf 1}_{\{j_1=j_4\}}
{\bf 1}_{\{i_2=i_3\ne 0\}}
{\bf 1}_{\{j_2=j_3\}}\Biggr),
\end{equation}

\vspace{2mm}
$$
J[\psi^{(5)}]_{T,t}
=\hbox{\vtop{\offinterlineskip\halign{
\hfil#\hfil\cr
{\rm l.i.m.}\cr
$\stackrel{}{{}_{p_1,\ldots,p_5\to \infty}}$\cr
}} }\sum_{j_1=0}^{p_1}\ldots\sum_{j_5=0}^{p_5}
C_{j_5\ldots j_1}\Biggl(
\prod_{l=1}^5\zeta_{j_l}^{(i_l)}
-\Biggr.
$$
$$
-
{\bf 1}_{\{i_1=i_2\ne 0\}}
{\bf 1}_{\{j_1=j_2\}}
\zeta_{j_3}^{(i_3)}
\zeta_{j_4}^{(i_4)}
\zeta_{j_5}^{(i_5)}-
{\bf 1}_{\{i_1=i_3\ne 0\}}
{\bf 1}_{\{j_1=j_3\}}
\zeta_{j_2}^{(i_2)}
\zeta_{j_4}^{(i_4)}
\zeta_{j_5}^{(i_5)}-
$$
$$
-
{\bf 1}_{\{i_1=i_4\ne 0\}}
{\bf 1}_{\{j_1=j_4\}}
\zeta_{j_2}^{(i_2)}
\zeta_{j_3}^{(i_3)}
\zeta_{j_5}^{(i_5)}-
{\bf 1}_{\{i_1=i_5\ne 0\}}
{\bf 1}_{\{j_1=j_5\}}
\zeta_{j_2}^{(i_2)}
\zeta_{j_3}^{(i_3)}
\zeta_{j_4}^{(i_4)}-
$$
$$
-
{\bf 1}_{\{i_2=i_3\ne 0\}}
{\bf 1}_{\{j_2=j_3\}}
\zeta_{j_1}^{(i_1)}
\zeta_{j_4}^{(i_4)}
\zeta_{j_5}^{(i_5)}-
{\bf 1}_{\{i_2=i_4\ne 0\}}
{\bf 1}_{\{j_2=j_4\}}
\zeta_{j_1}^{(i_1)}
\zeta_{j_3}^{(i_3)}
\zeta_{j_5}^{(i_5)}-
$$
$$
-
{\bf 1}_{\{i_2=i_5\ne 0\}}
{\bf 1}_{\{j_2=j_5\}}
\zeta_{j_1}^{(i_1)}
\zeta_{j_3}^{(i_3)}
\zeta_{j_4}^{(i_4)}
-{\bf 1}_{\{i_3=i_4\ne 0\}}
{\bf 1}_{\{j_3=j_4\}}
\zeta_{j_1}^{(i_1)}
\zeta_{j_2}^{(i_2)}
\zeta_{j_5}^{(i_5)}-
$$
$$
-
{\bf 1}_{\{i_3=i_5\ne 0\}}
{\bf 1}_{\{j_3=j_5\}}
\zeta_{j_1}^{(i_1)}
\zeta_{j_2}^{(i_2)}
\zeta_{j_4}^{(i_4)}
-{\bf 1}_{\{i_4=i_5\ne 0\}}
{\bf 1}_{\{j_4=j_5\}}
\zeta_{j_1}^{(i_1)}
\zeta_{j_2}^{(i_2)}
\zeta_{j_3}^{(i_3)}+
$$
$$
+
{\bf 1}_{\{i_1=i_2\ne 0\}}
{\bf 1}_{\{j_1=j_2\}}
{\bf 1}_{\{i_3=i_4\ne 0\}}
{\bf 1}_{\{j_3=j_4\}}\zeta_{j_5}^{(i_5)}+
{\bf 1}_{\{i_1=i_2\ne 0\}}
{\bf 1}_{\{j_1=j_2\}}
{\bf 1}_{\{i_3=i_5\ne 0\}}
{\bf 1}_{\{j_3=j_5\}}\zeta_{j_4}^{(i_4)}+
$$
$$
+
{\bf 1}_{\{i_1=i_2\ne 0\}}
{\bf 1}_{\{j_1=j_2\}}
{\bf 1}_{\{i_4=i_5\ne 0\}}
{\bf 1}_{\{j_4=j_5\}}\zeta_{j_3}^{(i_3)}+
{\bf 1}_{\{i_1=i_3\ne 0\}}
{\bf 1}_{\{j_1=j_3\}}
{\bf 1}_{\{i_2=i_4\ne 0\}}
{\bf 1}_{\{j_2=j_4\}}\zeta_{j_5}^{(i_5)}+
$$
$$
+
{\bf 1}_{\{i_1=i_3\ne 0\}}
{\bf 1}_{\{j_1=j_3\}}
{\bf 1}_{\{i_2=i_5\ne 0\}}
{\bf 1}_{\{j_2=j_5\}}\zeta_{j_4}^{(i_4)}+
{\bf 1}_{\{i_1=i_3\ne 0\}}
{\bf 1}_{\{j_1=j_3\}}
{\bf 1}_{\{i_4=i_5\ne 0\}}
{\bf 1}_{\{j_4=j_5\}}\zeta_{j_2}^{(i_2)}+
$$
$$
+
{\bf 1}_{\{i_1=i_4\ne 0\}}
{\bf 1}_{\{j_1=j_4\}}
{\bf 1}_{\{i_2=i_3\ne 0\}}
{\bf 1}_{\{j_2=j_3\}}\zeta_{j_5}^{(i_5)}+
{\bf 1}_{\{i_1=i_4\ne 0\}}
{\bf 1}_{\{j_1=j_4\}}
{\bf 1}_{\{i_2=i_5\ne 0\}}
{\bf 1}_{\{j_2=j_5\}}\zeta_{j_3}^{(i_3)}+
$$
$$
+
{\bf 1}_{\{i_1=i_4\ne 0\}}
{\bf 1}_{\{j_1=j_4\}}
{\bf 1}_{\{i_3=i_5\ne 0\}}
{\bf 1}_{\{j_3=j_5\}}\zeta_{j_2}^{(i_2)}+
{\bf 1}_{\{i_1=i_5\ne 0\}}
{\bf 1}_{\{j_1=j_5\}}
{\bf 1}_{\{i_2=i_3\ne 0\}}
{\bf 1}_{\{j_2=j_3\}}\zeta_{j_4}^{(i_4)}+
$$
$$
+
{\bf 1}_{\{i_1=i_5\ne 0\}}
{\bf 1}_{\{j_1=j_5\}}
{\bf 1}_{\{i_2=i_4\ne 0\}}
{\bf 1}_{\{j_2=j_4\}}\zeta_{j_3}^{(i_3)}+
{\bf 1}_{\{i_1=i_5\ne 0\}}
{\bf 1}_{\{j_1=j_5\}}
{\bf 1}_{\{i_3=i_4\ne 0\}}
{\bf 1}_{\{j_3=j_4\}}\zeta_{j_2}^{(i_2)}+
$$
$$
+
{\bf 1}_{\{i_2=i_3\ne 0\}}
{\bf 1}_{\{j_2=j_3\}}
{\bf 1}_{\{i_4=i_5\ne 0\}}
{\bf 1}_{\{j_4=j_5\}}\zeta_{j_1}^{(i_1)}+
{\bf 1}_{\{i_2=i_4\ne 0\}}
{\bf 1}_{\{j_2=j_4\}}
{\bf 1}_{\{i_3=i_5\ne 0\}}
{\bf 1}_{\{j_3=j_5\}}\zeta_{j_1}^{(i_1)}+
$$
\begin{equation}
\label{a5}
+\Biggl.
{\bf 1}_{\{i_2=i_5\ne 0\}}
{\bf 1}_{\{j_2=j_5\}}
{\bf 1}_{\{i_3=i_4\ne 0\}}
{\bf 1}_{\{j_3=j_4\}}\zeta_{j_1}^{(i_1)}\Biggr),
\end{equation}

\vspace{2mm}
$$
J[\psi^{(6)}]_{T,t}
=\hbox{\vtop{\offinterlineskip\halign{
\hfil#\hfil\cr
{\rm l.i.m.}\cr
$\stackrel{}{{}_{p_1,\ldots,p_6\to \infty}}$\cr
}} }\sum_{j_1=0}^{p_1}\ldots\sum_{j_6=0}^{p_6}
C_{j_6\ldots j_1}\Biggl(
\prod_{l=1}^6
\zeta_{j_l}^{(i_l)}
-\Biggr.
$$
$$
-
{\bf 1}_{\{i_1=i_6\ne 0\}}
{\bf 1}_{\{j_1=j_6\}}
\zeta_{j_2}^{(i_2)}
\zeta_{j_3}^{(i_3)}
\zeta_{j_4}^{(i_4)}
\zeta_{j_5}^{(i_5)}-
{\bf 1}_{\{i_2=i_6\ne 0\}}
{\bf 1}_{\{j_2=j_6\}}
\zeta_{j_1}^{(i_1)}
\zeta_{j_3}^{(i_3)}
\zeta_{j_4}^{(i_4)}
\zeta_{j_5}^{(i_5)}-
$$
$$
-
{\bf 1}_{\{i_3=i_6\ne 0\}}
{\bf 1}_{\{j_3=j_6\}}
\zeta_{j_1}^{(i_1)}
\zeta_{j_2}^{(i_2)}
\zeta_{j_4}^{(i_4)}
\zeta_{j_5}^{(i_5)}-
{\bf 1}_{\{i_4=i_6\ne 0\}}
{\bf 1}_{\{j_4=j_6\}}
\zeta_{j_1}^{(i_1)}
\zeta_{j_2}^{(i_2)}
\zeta_{j_3}^{(i_3)}
\zeta_{j_5}^{(i_5)}-
$$
$$
-
{\bf 1}_{\{i_5=i_6\ne 0\}}
{\bf 1}_{\{j_5=j_6\}}
\zeta_{j_1}^{(i_1)}
\zeta_{j_2}^{(i_2)}
\zeta_{j_3}^{(i_3)}
\zeta_{j_4}^{(i_4)}-
{\bf 1}_{\{i_1=i_2\ne 0\}}
{\bf 1}_{\{j_1=j_2\}}
\zeta_{j_3}^{(i_3)}
\zeta_{j_4}^{(i_4)}
\zeta_{j_5}^{(i_5)}
\zeta_{j_6}^{(i_6)}-
$$
$$
-
{\bf 1}_{\{i_1=i_3\ne 0\}}
{\bf 1}_{\{j_1=j_3\}}
\zeta_{j_2}^{(i_2)}
\zeta_{j_4}^{(i_4)}
\zeta_{j_5}^{(i_5)}
\zeta_{j_6}^{(i_6)}-
{\bf 1}_{\{i_1=i_4\ne 0\}}
{\bf 1}_{\{j_1=j_4\}}
\zeta_{j_2}^{(i_2)}
\zeta_{j_3}^{(i_3)}
\zeta_{j_5}^{(i_5)}
\zeta_{j_6}^{(i_6)}-
$$
$$
-
{\bf 1}_{\{i_1=i_5\ne 0\}}
{\bf 1}_{\{j_1=j_5\}}
\zeta_{j_2}^{(i_2)}
\zeta_{j_3}^{(i_3)}
\zeta_{j_4}^{(i_4)}
\zeta_{j_6}^{(i_6)}-
{\bf 1}_{\{i_2=i_3\ne 0\}}
{\bf 1}_{\{j_2=j_3\}}
\zeta_{j_1}^{(i_1)}
\zeta_{j_4}^{(i_4)}
\zeta_{j_5}^{(i_5)}
\zeta_{j_6}^{(i_6)}-
$$
$$
-
{\bf 1}_{\{i_2=i_4\ne 0\}}
{\bf 1}_{\{j_2=j_4\}}
\zeta_{j_1}^{(i_1)}
\zeta_{j_3}^{(i_3)}
\zeta_{j_5}^{(i_5)}
\zeta_{j_6}^{(i_6)}-
{\bf 1}_{\{i_2=i_5\ne 0\}}
{\bf 1}_{\{j_2=j_5\}}
\zeta_{j_1}^{(i_1)}
\zeta_{j_3}^{(i_3)}
\zeta_{j_4}^{(i_4)}
\zeta_{j_6}^{(i_6)}-
$$
$$
-
{\bf 1}_{\{i_3=i_4\ne 0\}}
{\bf 1}_{\{j_3=j_4\}}
\zeta_{j_1}^{(i_1)}
\zeta_{j_2}^{(i_2)}
\zeta_{j_5}^{(i_5)}
\zeta_{j_6}^{(i_6)}-
{\bf 1}_{\{i_3=i_5\ne 0\}}
{\bf 1}_{\{j_3=j_5\}}
\zeta_{j_1}^{(i_1)}
\zeta_{j_2}^{(i_2)}
\zeta_{j_4}^{(i_4)}
\zeta_{j_6}^{(i_6)}-
$$
$$
-
{\bf 1}_{\{i_4=i_5\ne 0\}}
{\bf 1}_{\{j_4=j_5\}}
\zeta_{j_1}^{(i_1)}
\zeta_{j_2}^{(i_2)}
\zeta_{j_3}^{(i_3)}
\zeta_{j_6}^{(i_6)}+
$$
$$
+
{\bf 1}_{\{i_1=i_2\ne 0\}}
{\bf 1}_{\{j_1=j_2\}}
{\bf 1}_{\{i_3=i_4\ne 0\}}
{\bf 1}_{\{j_3=j_4\}}
\zeta_{j_5}^{(i_5)}
\zeta_{j_6}^{(i_6)}+
$$
$$
+
{\bf 1}_{\{i_1=i_2\ne 0\}}
{\bf 1}_{\{j_1=j_2\}}
{\bf 1}_{\{i_3=i_5\ne 0\}}
{\bf 1}_{\{j_3=j_5\}}
\zeta_{j_4}^{(i_4)}
\zeta_{j_6}^{(i_6)}+
$$
$$
+
{\bf 1}_{\{i_1=i_2\ne 0\}}
{\bf 1}_{\{j_1=j_2\}}
{\bf 1}_{\{i_4=i_5\ne 0\}}
{\bf 1}_{\{j_4=j_5\}}
\zeta_{j_3}^{(i_3)}
\zeta_{j_6}^{(i_6)}
+
$$
$$
+
{\bf 1}_{\{i_1=i_3\ne 0\}}
{\bf 1}_{\{j_1=j_3\}}
{\bf 1}_{\{i_2=i_4\ne 0\}}
{\bf 1}_{\{j_2=j_4\}}
\zeta_{j_5}^{(i_5)}
\zeta_{j_6}^{(i_6)}+
$$
$$
+
{\bf 1}_{\{i_1=i_3\ne 0\}}
{\bf 1}_{\{j_1=j_3\}}
{\bf 1}_{\{i_2=i_5\ne 0\}}
{\bf 1}_{\{j_2=j_5\}}
\zeta_{j_4}^{(i_4)}
\zeta_{j_6}^{(i_6)}
+
$$
$$
+{\bf 1}_{\{i_1=i_3\ne 0\}}
{\bf 1}_{\{j_1=j_3\}}
{\bf 1}_{\{i_4=i_5\ne 0\}}
{\bf 1}_{\{j_4=j_5\}}
\zeta_{j_2}^{(i_2)}
\zeta_{j_6}^{(i_6)}+
$$
$$
+
{\bf 1}_{\{i_1=i_4\ne 0\}}
{\bf 1}_{\{j_1=j_4\}}
{\bf 1}_{\{i_2=i_3\ne 0\}}
{\bf 1}_{\{j_2=j_3\}}
\zeta_{j_5}^{(i_5)}
\zeta_{j_6}^{(i_6)}
+
$$
$$
+
{\bf 1}_{\{i_1=i_4\ne 0\}}
{\bf 1}_{\{j_1=j_4\}}
{\bf 1}_{\{i_2=i_5\ne 0\}}
{\bf 1}_{\{j_2=j_5\}}
\zeta_{j_3}^{(i_3)}
\zeta_{j_6}^{(i_6)}+
$$
$$
+
{\bf 1}_{\{i_1=i_4\ne 0\}}
{\bf 1}_{\{j_1=j_4\}}
{\bf 1}_{\{i_3=i_5\ne 0\}}
{\bf 1}_{\{j_3=j_5\}}
\zeta_{j_2}^{(i_2)}
\zeta_{j_6}^{(i_6)}
+
$$
$$
+
{\bf 1}_{\{i_1=i_5\ne 0\}}
{\bf 1}_{\{j_1=j_5\}}
{\bf 1}_{\{i_2=i_3\ne 0\}}
{\bf 1}_{\{j_2=j_3\}}
\zeta_{j_4}^{(i_4)}
\zeta_{j_6}^{(i_6)}+
$$
$$
+
{\bf 1}_{\{i_1=i_5\ne 0\}}
{\bf 1}_{\{j_1=j_5\}}
{\bf 1}_{\{i_2=i_4\ne 0\}}
{\bf 1}_{\{j_2=j_4\}}
\zeta_{j_3}^{(i_3)}
\zeta_{j_6}^{(i_6)}
+
$$
$$
+
{\bf 1}_{\{i_1=i_5\ne 0\}}
{\bf 1}_{\{j_1=j_5\}}
{\bf 1}_{\{i_3=i_4\ne 0\}}
{\bf 1}_{\{j_3=j_4\}}
\zeta_{j_2}^{(i_2)}
\zeta_{j_6}^{(i_6)}+
$$
$$
+
{\bf 1}_{\{i_2=i_3\ne 0\}}
{\bf 1}_{\{j_2=j_3\}}
{\bf 1}_{\{i_4=i_5\ne 0\}}
{\bf 1}_{\{j_4=j_5\}}
\zeta_{j_1}^{(i_1)}
\zeta_{j_6}^{(i_6)}
+
$$
$$
+{\bf 1}_{\{i_2=i_4\ne 0\}}
{\bf 1}_{\{j_2=j_4\}}
{\bf 1}_{\{i_3=i_5\ne 0\}}
{\bf 1}_{\{j_3=j_5\}}
\zeta_{j_1}^{(i_1)}
\zeta_{j_6}^{(i_6)}+
$$
$$
+
{\bf 1}_{\{i_2=i_5\ne 0\}}
{\bf 1}_{\{j_2=j_5\}}
{\bf 1}_{\{i_3=i_4\ne 0\}}
{\bf 1}_{\{j_3=j_4\}}
\zeta_{j_1}^{(i_1)}
\zeta_{j_6}^{(i_6)}
+
$$
$$
+{\bf 1}_{\{i_6=i_1\ne 0\}}
{\bf 1}_{\{j_6=j_1\}}
{\bf 1}_{\{i_3=i_4\ne 0\}}
{\bf 1}_{\{j_3=j_4\}}
\zeta_{j_2}^{(i_2)}
\zeta_{j_5}^{(i_5)}+
$$
$$
+
{\bf 1}_{\{i_6=i_1\ne 0\}}
{\bf 1}_{\{j_6=j_1\}}
{\bf 1}_{\{i_3=i_5\ne 0\}}
{\bf 1}_{\{j_3=j_5\}}
\zeta_{j_2}^{(i_2)}
\zeta_{j_4}^{(i_4)}
+
$$
$$
+{\bf 1}_{\{i_6=i_1\ne 0\}}
{\bf 1}_{\{j_6=j_1\}}
{\bf 1}_{\{i_2=i_5\ne 0\}}
{\bf 1}_{\{j_2=j_5\}}
\zeta_{j_3}^{(i_3)}
\zeta_{j_4}^{(i_4)}+
$$
$$
+
{\bf 1}_{\{i_6=i_1\ne 0\}}
{\bf 1}_{\{j_6=j_1\}}
{\bf 1}_{\{i_2=i_4\ne 0\}}
{\bf 1}_{\{j_2=j_4\}}
\zeta_{j_3}^{(i_3)}
\zeta_{j_5}^{(i_5)}
+
$$
$$
+{\bf 1}_{\{i_6=i_1\ne 0\}}
{\bf 1}_{\{j_6=j_1\}}
{\bf 1}_{\{i_4=i_5\ne 0\}}
{\bf 1}_{\{j_4=j_5\}}
\zeta_{j_2}^{(i_2)}
\zeta_{j_3}^{(i_3)}+
$$
$$
+
{\bf 1}_{\{i_6=i_1\ne 0\}}
{\bf 1}_{\{j_6=j_1\}}
{\bf 1}_{\{i_2=i_3\ne 0\}}
{\bf 1}_{\{j_2=j_3\}}
\zeta_{j_4}^{(i_4)}
\zeta_{j_5}^{(i_5)}
+
$$
$$
+{\bf 1}_{\{i_6=i_2\ne 0\}}
{\bf 1}_{\{j_6=j_2\}}
{\bf 1}_{\{i_3=i_5\ne 0\}}
{\bf 1}_{\{j_3=j_5\}}
\zeta_{j_1}^{(i_1)}
\zeta_{j_4}^{(i_4)}+
$$
$$
+
{\bf 1}_{\{i_6=i_2\ne 0\}}
{\bf 1}_{\{j_6=j_2\}}
{\bf 1}_{\{i_4=i_5\ne 0\}}
{\bf 1}_{\{j_4=j_5\}}
\zeta_{j_1}^{(i_1)}
\zeta_{j_3}^{(i_3)}
+
$$
$$
+{\bf 1}_{\{i_6=i_2\ne 0\}}
{\bf 1}_{\{j_6=j_2\}}
{\bf 1}_{\{i_3=i_4\ne 0\}}
{\bf 1}_{\{j_3=j_4\}}
\zeta_{j_1}^{(i_1)}
\zeta_{j_5}^{(i_5)}+
$$
$$
+
{\bf 1}_{\{i_6=i_2\ne 0\}}
{\bf 1}_{\{j_6=j_2\}}
{\bf 1}_{\{i_1=i_5\ne 0\}}
{\bf 1}_{\{j_1=j_5\}}
\zeta_{j_3}^{(i_3)}
\zeta_{j_4}^{(i_4)}
+
$$
$$
+{\bf 1}_{\{i_6=i_2\ne 0\}}
{\bf 1}_{\{j_6=j_2\}}
{\bf 1}_{\{i_1=i_4\ne 0\}}
{\bf 1}_{\{j_1=j_4\}}
\zeta_{j_3}^{(i_3)}
\zeta_{j_5}^{(i_5)}+
$$
$$
+
{\bf 1}_{\{i_6=i_2\ne 0\}}
{\bf 1}_{\{j_6=j_2\}}
{\bf 1}_{\{i_1=i_3\ne 0\}}
{\bf 1}_{\{j_1=j_3\}}
\zeta_{j_4}^{(i_4)}
\zeta_{j_5}^{(i_5)}
+
$$
$$
+{\bf 1}_{\{i_6=i_3\ne 0\}}
{\bf 1}_{\{j_6=j_3\}}
{\bf 1}_{\{i_2=i_5\ne 0\}}
{\bf 1}_{\{j_2=j_5\}}
\zeta_{j_1}^{(i_1)}
\zeta_{j_4}^{(i_4)}+
$$
$$
+
{\bf 1}_{\{i_6=i_3\ne 0\}}
{\bf 1}_{\{j_6=j_3\}}
{\bf 1}_{\{i_4=i_5\ne 0\}}
{\bf 1}_{\{j_4=j_5\}}
\zeta_{j_1}^{(i_1)}
\zeta_{j_2}^{(i_2)}
+
$$
$$
+{\bf 1}_{\{i_6=i_3\ne 0\}}
{\bf 1}_{\{j_6=j_3\}}
{\bf 1}_{\{i_2=i_4\ne 0\}}
{\bf 1}_{\{j_2=j_4\}}
\zeta_{j_1}^{(i_1)}
\zeta_{j_5}^{(i_5)}+
$$
$$
+
{\bf 1}_{\{i_6=i_3\ne 0\}}
{\bf 1}_{\{j_6=j_3\}}
{\bf 1}_{\{i_1=i_5\ne 0\}}
{\bf 1}_{\{j_1=j_5\}}
\zeta_{j_2}^{(i_2)}
\zeta_{j_4}^{(i_4)}
+
$$
$$
+{\bf 1}_{\{i_6=i_3\ne 0\}}
{\bf 1}_{\{j_6=j_3\}}
{\bf 1}_{\{i_1=i_4\ne 0\}}
{\bf 1}_{\{j_1=j_4\}}
\zeta_{j_2}^{(i_2)}
\zeta_{j_5}^{(i_5)}+
$$
$$
+
{\bf 1}_{\{i_6=i_3\ne 0\}}
{\bf 1}_{\{j_6=j_3\}}
{\bf 1}_{\{i_1=i_2\ne 0\}}
{\bf 1}_{\{j_1=j_2\}}
\zeta_{j_4}^{(i_4)}
\zeta_{j_5}^{(i_5)}
+
$$
$$
+{\bf 1}_{\{i_6=i_4\ne 0\}}
{\bf 1}_{\{j_6=j_4\}}
{\bf 1}_{\{i_3=i_5\ne 0\}}
{\bf 1}_{\{j_3=j_5\}}
\zeta_{j_1}^{(i_1)}
\zeta_{j_2}^{(i_2)}+
$$
$$
+
{\bf 1}_{\{i_6=i_4\ne 0\}}
{\bf 1}_{\{j_6=j_4\}}
{\bf 1}_{\{i_2=i_5\ne 0\}}
{\bf 1}_{\{j_2=j_5\}}
\zeta_{j_1}^{(i_1)}
\zeta_{j_3}^{(i_3)}
+
$$
$$
+{\bf 1}_{\{i_6=i_4\ne 0\}}
{\bf 1}_{\{j_6=j_4\}}
{\bf 1}_{\{i_2=i_3\ne 0\}}
{\bf 1}_{\{j_2=j_3\}}
\zeta_{j_1}^{(i_1)}
\zeta_{j_5}^{(i_5)}+
$$
$$
+
{\bf 1}_{\{i_6=i_4\ne 0\}}
{\bf 1}_{\{j_6=j_4\}}
{\bf 1}_{\{i_1=i_5\ne 0\}}
{\bf 1}_{\{j_1=j_5\}}
\zeta_{j_2}^{(i_2)}
\zeta_{j_3}^{(i_3)}
+
$$
$$
+{\bf 1}_{\{i_6=i_4\ne 0\}}
{\bf 1}_{\{j_6=j_4\}}
{\bf 1}_{\{i_1=i_3\ne 0\}}
{\bf 1}_{\{j_1=j_3\}}
\zeta_{j_2}^{(i_2)}
\zeta_{j_5}^{(i_5)}+
$$
$$
+
{\bf 1}_{\{i_6=i_4\ne 0\}}
{\bf 1}_{\{j_6=j_4\}}
{\bf 1}_{\{i_1=i_2\ne 0\}}
{\bf 1}_{\{j_1=j_2\}}
\zeta_{j_3}^{(i_3)}
\zeta_{j_5}^{(i_5)}
+
$$
$$
+{\bf 1}_{\{i_6=i_5\ne 0\}}
{\bf 1}_{\{j_6=j_5\}}
{\bf 1}_{\{i_3=i_4\ne 0\}}
{\bf 1}_{\{j_3=j_4\}}
\zeta_{j_1}^{(i_1)}
\zeta_{j_2}^{(i_2)}+
$$
$$
+
{\bf 1}_{\{i_6=i_5\ne 0\}}
{\bf 1}_{\{j_6=j_5\}}
{\bf 1}_{\{i_2=i_4\ne 0\}}
{\bf 1}_{\{j_2=j_4\}}
\zeta_{j_1}^{(i_1)}
\zeta_{j_3}^{(i_3)}
+
$$
$$
+{\bf 1}_{\{i_6=i_5\ne 0\}}
{\bf 1}_{\{j_6=j_5\}}
{\bf 1}_{\{i_2=i_3\ne 0\}}
{\bf 1}_{\{j_2=j_3\}}
\zeta_{j_1}^{(i_1)}
\zeta_{j_4}^{(i_4)}+
$$
$$
+
{\bf 1}_{\{i_6=i_5\ne 0\}}
{\bf 1}_{\{j_6=j_5\}}
{\bf 1}_{\{i_1=i_4\ne 0\}}
{\bf 1}_{\{j_1=j_4\}}
\zeta_{j_2}^{(i_2)}
\zeta_{j_3}^{(i_3)}
+
$$
$$
+{\bf 1}_{\{i_6=i_5\ne 0\}}
{\bf 1}_{\{j_6=j_5\}}
{\bf 1}_{\{i_1=i_3\ne 0\}}
{\bf 1}_{\{j_1=j_3\}}
\zeta_{j_2}^{(i_2)}
\zeta_{j_4}^{(i_4)}+
$$
$$
+
{\bf 1}_{\{i_6=i_5\ne 0\}}
{\bf 1}_{\{j_6=j_5\}}
{\bf 1}_{\{i_1=i_2\ne 0\}}
{\bf 1}_{\{j_1=j_2\}}
\zeta_{j_3}^{(i_3)}
\zeta_{j_4}^{(i_4)}-
$$
$$
-
{\bf 1}_{\{i_6=i_1\ne 0\}}
{\bf 1}_{\{j_6=j_1\}}
{\bf 1}_{\{i_2=i_5\ne 0\}}
{\bf 1}_{\{j_2=j_5\}}
{\bf 1}_{\{i_3=i_4\ne 0\}}
{\bf 1}_{\{j_3=j_4\}}-
$$
$$
-
{\bf 1}_{\{i_6=i_1\ne 0\}}
{\bf 1}_{\{j_6=j_1\}}
{\bf 1}_{\{i_2=i_4\ne 0\}}
{\bf 1}_{\{j_2=j_4\}}
{\bf 1}_{\{i_3=i_5\ne 0\}}
{\bf 1}_{\{j_3=j_5\}}-
$$
$$
-
{\bf 1}_{\{i_6=i_1\ne 0\}}
{\bf 1}_{\{j_6=j_1\}}
{\bf 1}_{\{i_2=i_3\ne 0\}}
{\bf 1}_{\{j_2=j_3\}}
{\bf 1}_{\{i_4=i_5\ne 0\}}
{\bf 1}_{\{j_4=j_5\}}-
$$
$$
-
{\bf 1}_{\{i_6=i_2\ne 0\}}
{\bf 1}_{\{j_6=j_2\}}
{\bf 1}_{\{i_1=i_5\ne 0\}}
{\bf 1}_{\{j_1=j_5\}}
{\bf 1}_{\{i_3=i_4\ne 0\}}
{\bf 1}_{\{j_3=j_4\}}-
$$
$$
-
{\bf 1}_{\{i_6=i_2\ne 0\}}
{\bf 1}_{\{j_6=j_2\}}
{\bf 1}_{\{i_1=i_4\ne 0\}}
{\bf 1}_{\{j_1=j_4\}}
{\bf 1}_{\{i_3=i_5\ne 0\}}
{\bf 1}_{\{j_3=j_5\}}-
$$
$$
-
{\bf 1}_{\{i_6=i_2\ne 0\}}
{\bf 1}_{\{j_6=j_2\}}
{\bf 1}_{\{i_1=i_3\ne 0\}}
{\bf 1}_{\{j_1=j_3\}}
{\bf 1}_{\{i_4=i_5\ne 0\}}
{\bf 1}_{\{j_4=j_5\}}-
$$
$$
-
{\bf 1}_{\{i_6=i_3\ne 0\}}
{\bf 1}_{\{j_6=j_3\}}
{\bf 1}_{\{i_1=i_5\ne 0\}}
{\bf 1}_{\{j_1=j_5\}}
{\bf 1}_{\{i_2=i_4\ne 0\}}
{\bf 1}_{\{j_2=j_4\}}-
$$
$$
-
{\bf 1}_{\{i_6=i_3\ne 0\}}
{\bf 1}_{\{j_6=j_3\}}
{\bf 1}_{\{i_1=i_4\ne 0\}}
{\bf 1}_{\{j_1=j_4\}}
{\bf 1}_{\{i_2=i_5\ne 0\}}
{\bf 1}_{\{j_2=j_5\}}-
$$
$$
-
{\bf 1}_{\{i_3=i_6\ne 0\}}
{\bf 1}_{\{j_3=j_6\}}
{\bf 1}_{\{i_1=i_2\ne 0\}}
{\bf 1}_{\{j_1=j_2\}}
{\bf 1}_{\{i_4=i_5\ne 0\}}
{\bf 1}_{\{j_4=j_5\}}-
$$
$$
-
{\bf 1}_{\{i_6=i_4\ne 0\}}
{\bf 1}_{\{j_6=j_4\}}
{\bf 1}_{\{i_1=i_5\ne 0\}}
{\bf 1}_{\{j_1=j_5\}}
{\bf 1}_{\{i_2=i_3\ne 0\}}
{\bf 1}_{\{j_2=j_3\}}-
$$
$$
-
{\bf 1}_{\{i_6=i_4\ne 0\}}
{\bf 1}_{\{j_6=j_4\}}
{\bf 1}_{\{i_1=i_3\ne 0\}}
{\bf 1}_{\{j_1=j_3\}}
{\bf 1}_{\{i_2=i_5\ne 0\}}
{\bf 1}_{\{j_2=j_5\}}-
$$
$$
-
{\bf 1}_{\{i_6=i_4\ne 0\}}
{\bf 1}_{\{j_6=j_4\}}
{\bf 1}_{\{i_1=i_2\ne 0\}}
{\bf 1}_{\{j_1=j_2\}}
{\bf 1}_{\{i_3=i_5\ne 0\}}
{\bf 1}_{\{j_3=j_5\}}-
$$
$$
-
{\bf 1}_{\{i_6=i_5\ne 0\}}
{\bf 1}_{\{j_6=j_5\}}
{\bf 1}_{\{i_1=i_4\ne 0\}}
{\bf 1}_{\{j_1=j_4\}}
{\bf 1}_{\{i_2=i_3\ne 0\}}
{\bf 1}_{\{j_2=j_3\}}-
$$
$$
-
{\bf 1}_{\{i_6=i_5\ne 0\}}
{\bf 1}_{\{j_6=j_5\}}
{\bf 1}_{\{i_1=i_2\ne 0\}}
{\bf 1}_{\{j_1=j_2\}}
{\bf 1}_{\{i_3=i_4\ne 0\}}
{\bf 1}_{\{j_3=j_4\}}-
$$
\begin{equation}
\label{a6}
\Biggl.-
{\bf 1}_{\{i_6=i_5\ne 0\}}
{\bf 1}_{\{j_6=j_5\}}
{\bf 1}_{\{i_1=i_3\ne 0\}}
{\bf 1}_{\{j_1=j_3\}}
{\bf 1}_{\{i_2=i_4\ne 0\}}
{\bf 1}_{\{j_2=j_4\}}\Biggr),
\end{equation}

\vspace{2mm}
\noindent
where ${\bf 1}_A$ is the indicator of the set $A$.
 
A detailed discussion of advantages of the method based on Theorem 1 over the approximation 
methods from works \cite{2}, \cite{3}, \cite{7}, \cite{8}, \cite{27}-\cite{35}, \cite{37}-\cite{39} 
can be found in \cite{26} (Sect.~1.1.10) or in \cite{26a}, \cite{58}.

For further consideration, let us 
consider the generalization of formulas (\ref{a1})--(\ref{a6})                 
for the case of an arbitrary multiplicity $k$ $(k\in\mathbb{N})$ of 
the iterated Ito stochastic integral $J[\psi^{(k)}]_{T,t}$ defined by (\ref{ito}).
In order to do this, let us
introduce some notations. 
Consider the unordered
set $\{1, 2, \ldots, k\}$ 
and separate it into two parts:
the first part consists of $r$ unordered 
pairs (sequence order of these pairs is also unimportant) and the 
second one consists of the 
remaining $k-2r$ numbers.
So, we have
\begin{equation}
\label{leto5007}
(\{
\underbrace{\{g_1, g_2\}, \ldots, 
\{g_{2r-1}, g_{2r}\}}_{\small{\hbox{part 1}}}
\},
\{\underbrace{q_1, \ldots, q_{k-2r}}_{\small{\hbox{part 2}}}
\}),
\end{equation}

\noindent
where 
$$
\{g_1, g_2, \ldots, 
g_{2r-1}, g_{2r}, q_1, \ldots, q_{k-2r}\}=\{1, 2, \ldots, k\},
$$

\noindent
braces   
mean an unordered 
set, and pa\-ren\-the\-ses mean an ordered set.

We will say that (\ref{leto5007}) is a partition 
and consider the sum with respect to all possible
partitions
\begin{equation}
\label{leto5008}
\sum_{\stackrel{(\{\{g_1, g_2\}, \ldots, 
\{g_{2r-1}, g_{2r}\}\}, \{q_1, \ldots, q_{k-2r}\})}
{{}_{\{g_1, g_2, \ldots, 
g_{2r-1}, g_{2r}, q_1, \ldots, q_{k-2r}\}=\{1, 2, \ldots, k\}}}}
a_{g_1 g_2, \ldots, 
g_{2r-1} g_{2r}, q_1 \ldots q_{k-2r}}.
\end{equation}

\vspace{2mm}

Below there are several examples of sums in the form (\ref{leto5008})
$$
\sum_{\stackrel{(\{g_1, g_2\})}{{}_{\{g_1, g_2\}=\{1, 2\}}}}
a_{g_1 g_2}=a_{12},
$$

\vspace{2mm}
$$
\sum_{\stackrel{(\{\{g_1, g_2\}, \{g_3, g_4\}\})}
{{}_{\{g_1, g_2, g_3, g_4\}=\{1, 2, 3, 4\}}}}
a_{g_1 g_2, g_3 g_4}=a_{12,34} + a_{13,24} + a_{23,14},
$$

\vspace{2mm}
$$
\sum_{\stackrel{(\{g_1, g_2\}, \{q_1, q_{2}\})}
{{}_{\{g_1, g_2, q_1, q_{2}\}=\{1, 2, 3, 4\}}}}
a_{g_1 g_2, q_1 q_{2}}=
$$
$$
=a_{12,34}+a_{13,24}+a_{14,23}
+a_{23,14}+a_{24,13}+a_{34,12},
$$
$$
\sum_{\stackrel{(\{g_1, g_2\}, \{q_1, q_{2}, q_3\})}
{{}_{\{g_1, g_2, q_1, q_{2}, q_3\}=\{1, 2, 3, 4, 5\}}}}
a_{g_1 g_2, q_1 q_{2}q_3}
=
$$
$$
=a_{12,345}+a_{13,245}+a_{14,235}
+a_{15,234}+a_{23,145}+a_{24,135}+
$$
$$
+a_{25,134}+a_{34,125}+a_{35,124}+a_{45,123},
$$

$$
\sum_{\stackrel{(\{\{g_1, g_2\}, \{g_3, g_{4}\}\}, \{q_1\})}
{{}_{\{g_1, g_2, g_3, g_{4}, q_1\}=\{1, 2, 3, 4, 5\}}}}
a_{g_1 g_2, g_3 g_{4},q_1}
=
$$
$$
=
a_{12,34,5}+a_{13,24,5}+a_{14,23,5}+
a_{12,35,4}+a_{13,25,4}+a_{15,23,4}+
$$
$$
+a_{12,54,3}+a_{15,24,3}+a_{14,25,3}+a_{15,34,2}+a_{13,54,2}+a_{14,53,2}+
$$
$$
+
a_{52,34,1}+a_{53,24,1}+a_{54,23,1}.
$$

\vspace{4mm}

Now we can write (\ref{tyyy}) as 
$$
J[\psi^{(k)}]_{T,t}=
\hbox{\vtop{\offinterlineskip\halign{
\hfil#\hfil\cr
{\rm l.i.m.}\cr
$\stackrel{}{{}_{p_1,\ldots,p_k\to \infty}}$\cr
}} }
\sum\limits_{j_1=0}^{p_1}\ldots
\sum\limits_{j_k=0}^{p_k}
C_{j_k\ldots j_1}\Biggl(
\prod_{l=1}^k\zeta_{j_l}^{(i_l)}+\sum\limits_{r=1}^{[k/2]}
(-1)^r \times
\Biggr.
$$
\begin{equation}
\label{leto6000hh}
\times
\sum_{\stackrel{(\{\{g_1, g_2\}, \ldots, 
\{g_{2r-1}, g_{2r}\}\}, \{q_1, \ldots, q_{k-2r}\})}
{{}_{\{g_1, g_2, \ldots, 
g_{2r-1}, g_{2r}, q_1, \ldots, q_{k-2r}\}=\{1, 2, \ldots, k\}}}}
\prod\limits_{s=1}^r
{\bf 1}_{\{i_{g_{{}_{2s-1}}}=~i_{g_{{}_{2s}}}\ne 0\}}
\Biggl.{\bf 1}_{\{j_{g_{{}_{2s-1}}}=~j_{g_{{}_{2s}}}\}}
\prod_{l=1}^{k-2r}\zeta_{j_{q_l}}^{(i_{q_l})}\Biggr),
\end{equation}

\vspace{1mm}
\noindent
where $[x]$ is an integer part of a real number $x,$
$\prod\limits_{\emptyset}
\stackrel{\sf def}{=}1,$ $\sum\limits_{\emptyset}
\stackrel{\sf def}{=}0;$
another notations are the same as in Theorem 1.

In particular, from (\ref{leto6000hh}) for $k=5$ we obtain
$$
J[\psi^{(5)}]_{T,t}=
\hbox{\vtop{\offinterlineskip\halign{
\hfil#\hfil\cr
{\rm l.i.m.}\cr
$\stackrel{}{{}_{p_1,\ldots,p_5\to \infty}}$\cr
}} }\sum_{j_1=0}^{p_1}\ldots\sum_{j_5=0}^{p_5}
C_{j_5\ldots j_1}\Biggl(
\prod_{l=1}^5\zeta_{j_l}^{(i_l)}-\Biggr.
$$
$$
-
\sum\limits_{\stackrel{(\{g_1, g_2\}, \{q_1, q_{2}, q_3\})}
{{}_{\{g_1, g_2, q_{1}, q_{2}, q_3\}=\{1, 2, 3, 4, 5\}}}}
{\bf 1}_{\{i_{g_{{}_{1}}}=~i_{g_{{}_{2}}}\ne 0\}}
{\bf 1}_{\{j_{g_{{}_{1}}}=~j_{g_{{}_{2}}}\}}
\prod_{l=1}^{3}\zeta_{j_{q_l}}^{(i_{q_l})}+
$$
$$
+
\sum_{\stackrel{(\{\{g_1, g_2\}, 
\{g_{3}, g_{4}\}\}, \{q_1\})}
{{}_{\{g_1, g_2, g_{3}, g_{4}, q_1\}=\{1, 2, 3, 4, 5\}}}}
{\bf 1}_{\{i_{g_{{}_{1}}}=~i_{g_{{}_{2}}}\ne 0\}}
{\bf 1}_{\{j_{g_{{}_{1}}}=~j_{g_{{}_{2}}}\}}
\Biggl.{\bf 1}_{\{i_{g_{{}_{3}}}=~i_{g_{{}_{4}}}\ne 0\}}
{\bf 1}_{\{j_{g_{{}_{3}}}=~j_{g_{{}_{4}}}\}}
\zeta_{j_{q_1}}^{(i_{q_1})}\Biggr).
$$

\vspace{1mm}
\noindent
The last equality obviously agrees with
(\ref{a5}).

Let us consider the generalization of Theorem 1 for the case
of an arbitrary complete orthonormal systems  
of functions in the space $L_2([t,T])$ 
and $\psi_1(\tau),\ldots,\psi_k(\tau)\in L_2([t, T]).$

{\bf Theorem~2}\ \cite{58} (Sect.~1.11, 1.14), \cite{add1011} (Sect.~15, 18),
\cite{diffjournal}.\
{\it Suppose that
$\psi_1(\tau),$ $\ldots,$ $\psi_k(\tau)\in L_2([t, T])$ and
$\{\phi_j(x)\}_{j=0}^{\infty}$ is an arbitrary complete orthonormal system  
of functions in the space $L_2([t,T]).$
Then the following expansion
$$
J[\psi^{(k)}]_{T,t}=
\hbox{\vtop{\offinterlineskip\halign{
\hfil#\hfil\cr
{\rm l.i.m.}\cr
$\stackrel{}{{}_{p_1,\ldots,p_k\to \infty}}$\cr
}} }
\sum\limits_{j_1=0}^{p_1}\ldots
\sum\limits_{j_k=0}^{p_k}
C_{j_k\ldots j_1}\Biggl(
\prod_{l=1}^k\zeta_{j_l}^{(i_l)}+\sum\limits_{r=1}^{[k/2]}
(-1)^r \times
\Biggr.
$$
\begin{equation}
\label{leto6000}
\times
\sum_{\stackrel{(\{\{g_1, g_2\}, \ldots, 
\{g_{2r-1}, g_{2r}\}\}, \{q_1, \ldots, q_{k-2r}\})}
{{}_{\{g_1, g_2, \ldots, 
g_{2r-1}, g_{2r}, q_1, \ldots, q_{k-2r}\}=\{1, 2, \ldots, k\}}}}
\prod\limits_{s=1}^r
{\bf 1}_{\{i_{g_{{}_{2s-1}}}=~i_{g_{{}_{2s}}}\ne 0\}}
\Biggl.{\bf 1}_{\{j_{g_{{}_{2s-1}}}=~j_{g_{{}_{2s}}}\}}
\prod_{l=1}^{k-2r}\zeta_{j_{q_l}}^{(i_{q_l})}\Biggr)
\end{equation}

\vspace{1mm}
\noindent
con\-verg\-ing in the mean-square sense is valid,
where $[x]$ is an integer part of a real number $x,$
$\prod\limits_{\emptyset}
\stackrel{\sf def}{=}1,$ $\sum\limits_{\emptyset}
\stackrel{\sf def}{=}0;$
another notations are the same as in Theorem~{\rm 1}.}

Note that an analogue of Theorem 2 (based on the explicit use
of Hermite polynomials) was considered 
in \cite{Rybakov1000}. 
We use another notations and proof
\cite{58} (Sect.~1.11, 1.14), \cite{add1011} (Sect.~15, 18), \cite{diffjournal}
in comparison with \cite{Rybakov1000}.

As it turned out, Theorems 1, 2 can be adapted
for the iterated  Stratonovich 
stochastic integrals
(\ref{str}) of multiplicities 1 to 6
\cite{26}, \cite{26a}, \cite{58}-\cite{59a} (also see bibliography therein).
Let as collect some old results in the following theorem.

{\bf Theorem 3} \cite{26}, \cite{26a}, \cite{58}-\cite{59a}. {\it Suppose that 
$\{\phi_j(x)\}_{j=0}^{\infty}$ is a complete orthonormal system of 
Legendre polynomials or trigonometric functions in the space $L_2([t, T]).$
At the same time $\psi_2(\tau)$ is a continuously differentiable 
function on $[t, T]$ and $\psi_1(\tau),$ $\psi_3(\tau)$ are twice
continuously differentiable functions on $[t, T]$. Then 
\begin{equation}
\label{a}
J^{*}[\psi^{(2)}]_{T,t}=
\hbox{\vtop{\offinterlineskip\halign{
\hfil#\hfil\cr
{\rm l.i.m.}\cr
$\stackrel{}{{}_{p_1,p_2\to \infty}}$\cr
}} }\sum_{j_1=0}^{p_1}\sum_{j_2=0}^{p_2}
C_{j_2j_1}\zeta_{j_1}^{(i_1)}\zeta_{j_2}^{(i_2)}\ \ (i_1,i_2=1,\ldots,m),
\end{equation}
\begin{equation}
\label{feto19000ab}
J^{*}[\psi^{(3)}]_{T,t}=
\hbox{\vtop{\offinterlineskip\halign{
\hfil#\hfil\cr
{\rm l.i.m.}\cr
$\stackrel{}{{}_{p_1,p_2,p_3\to \infty}}$\cr
}} }\sum_{j_1=0}^{p_1}\sum_{j_2=0}^{p_2}\sum_{j_3=0}^{p_3}
C_{j_3 j_2 j_1}\zeta_{j_1}^{(i_1)}\zeta_{j_2}^{(i_2)}\zeta_{j_3}^{(i_3)}\ \
(i_1,i_2,i_3=0, 1,\ldots,m),
\end{equation}
\begin{equation}
\label{feto19000a}
J^{*}[\psi^{(3)}]_{T,t}=
\hbox{\vtop{\offinterlineskip\halign{
\hfil#\hfil\cr
{\rm l.i.m.}\cr
$\stackrel{}{{}_{p\to \infty}}$\cr
}} }
\sum\limits_{j_1,j_2,j_3=0}^{p}
C_{j_3 j_2 j_1}\zeta_{j_1}^{(i_1)}\zeta_{j_2}^{(i_2)}\zeta_{j_3}^{(i_3)}\ \
(i_1,i_2,i_3=1,\ldots,m),
\end{equation}
\begin{equation}
\label{uu}
J^{*}[\psi^{(4)}]_{T,t}=
\hbox{\vtop{\offinterlineskip\halign{
\hfil#\hfil\cr
{\rm l.i.m.}\cr
$\stackrel{}{{}_{p\to \infty}}$\cr
}} }
\sum\limits_{j_1, \ldots, j_4=0}^{p}
C_{j_4 j_3 j_2 j_1}\zeta_{j_1}^{(i_1)}
\zeta_{j_2}^{(i_2)}\zeta_{j_3}^{(i_3)}\zeta_{j_4}^{(i_4)}\ \
(i_1,\ldots,i_4=0, 1,\ldots,m),
\end{equation}
where $J^{*}[\psi^{(k)}]_{T,t}$ is defined by {\rm (\ref{str})}, and
$\psi_l(\tau)\equiv 1$ $(l=1,\ldots,4)$ in {\rm (\ref{feto19000ab})}, 
{\rm (\ref{uu});} another notations are the same as in Theorems {\rm 1, 2.}
}

In 2022, a new approach to the expansion and mean-square 
approximation of iterated Stratonovich stochastic integrals has been obtained
\cite{58} (Sect.~2.10--2.16), \cite{add1014} (Sect.~13--19), 
\cite{add1019} (Sect.~5--11), \cite{add1020} (Sect.~7--13),
\cite{new-art-1-xxy} (Sect.~4--9).

Let us formulate four theorems that were obtained using this approach.

{\bf Theorem 4}\ \cite{58}, \cite{add1014}, \cite{add1019}, \cite{add1020}, \cite{new-art-1-xxy}.\
{\it Suppose 
that $\{\phi_j(x)\}_{j=0}^{\infty}$ is a complete orthonormal system of 
Legendre polynomials or trigonometric functions in the space $L_2([t, T]).$
Furthermore, let $\psi_1(\tau),\psi_2(\tau),\psi_3(\tau)$ are continuously dif\-ferentiable 
nonrandom functions on $[t, T].$ 
Then, for the 
iterated Stra\-to\-no\-vich stochastic integral of third multiplicity
\begin{equation}
\label{hehe111}
J^{*}[\psi^{(3)}]_{T,t}={\int\limits_t^{*}}^T\psi_3(t_3)
{\int\limits_t^{*}}^{t_3}\psi_2(t_2)
{\int\limits_t^{*}}^{t_2}\psi_1(t_1)
d{\bf w}_{t_1}^{(i_1)}
d{\bf w}_{t_2}^{(i_2)}d{\bf w}_{t_3}^{(i_3)}
\end{equation}
the following 
relations
\begin{equation}
\label{fin1}
J^{*}[\psi^{(3)}]_{T,t}
=\hbox{\vtop{\offinterlineskip\halign{
\hfil#\hfil\cr
{\rm l.i.m.}\cr
$\stackrel{}{{}_{p\to \infty}}$\cr
}} }
\sum\limits_{j_1, j_2, j_3=0}^{p}
C_{j_3 j_2 j_1}\zeta_{j_1}^{(i_1)}\zeta_{j_2}^{(i_2)}\zeta_{j_3}^{(i_3)},
\end{equation}
\begin{equation}
\label{fin2}
{\sf M}\left\{\left(
J^{*}[\psi^{(3)}]_{T,t}-
\sum\limits_{j_1, j_2, j_3=0}^{p}
C_{j_3 j_2 j_1}\zeta_{j_1}^{(i_1)}\zeta_{j_2}^{(i_2)}\zeta_{j_3}^{(i_3)}\right)^2\right\}
\le \frac{C}{p}
\end{equation}

\vspace{1mm}
\noindent
are fulfilled, where $i_1, i_2, i_3=0,1,\ldots,m$ in {\rm (\ref{hehe111}), (\ref{fin1})} and 
$i_1, i_2, i_3=1,\ldots,m$ in {\rm (\ref{fin2})},
constant $C$ is independent of $p,$
$$
C_{j_3 j_2 j_1}=\int\limits_t^T\psi_3(t_3)\phi_{j_3}(t_3)
\int\limits_t^{t_3}\psi_2(t_2)\phi_{j_2}(t_2)
\int\limits_t^{t_2}\psi_1(t_1)\phi_{j_1}(t_1)dt_1dt_2dt_3
$$
and
$$
\zeta_{j}^{(i)}=
\int\limits_t^T \phi_{j}(\tau) d{\bf w}_{\tau}^{(i)}
$$ 
are independent standard Gaussian random variables for various 
$i$ or $j$ {\rm (}in the case when $i\ne 0${\rm );} 
another notations are the same as in Theorems~{\rm 1, 2}.}

{\bf Theorem 5}\ \cite{58}, \cite{add1014}, \cite{add1019}, \cite{add1020}, \cite{new-art-1-xxy}.\ {\it Let
$\{\phi_j(x)\}_{j=0}^{\infty}$ be a complete orthonormal system of 
Legendre polynomials or trigonometric functions in the space $L_2([t, T]).$
Furthermore, let $\psi_1(\tau), \ldots,\psi_4(\tau)$ be continuously dif\-ferentiable 
nonrandom functions on $[t, T].$ 
Then, for the 
iterated Stra\-to\-no\-vich stochastic integral of fourth multiplicity
\begin{equation}
\label{fin0}
J^{*}[\psi^{(4)}]_{T,t}={\int\limits_t^{*}}^T\psi_4(t_4)
{\int\limits_t^{*}}^{t_4}\psi_3(t_3)
{\int\limits_t^{*}}^{t_3}\psi_2(t_2)
{\int\limits_t^{*}}^{t_2}\psi_1(t_1)
d{\bf w}_{t_1}^{(i_1)}
d{\bf w}_{t_2}^{(i_2)}d{\bf w}_{t_3}^{(i_3)}d{\bf w}_{t_4}^{(i_4)}
\end{equation}
the following 
relations
\begin{equation}
\label{fin3}
J^{*}[\psi^{(4)}]_{T,t}
=\hbox{\vtop{\offinterlineskip\halign{
\hfil#\hfil\cr
{\rm l.i.m.}\cr
$\stackrel{}{{}_{p\to \infty}}$\cr
}} }
\sum\limits_{j_1, j_2, j_3,j_4=0}^{p}
C_{j_4j_3 j_2 j_1}\zeta_{j_1}^{(i_1)}\zeta_{j_2}^{(i_2)}\zeta_{j_3}^{(i_3)}\zeta_{j_4}^{(i_4)},
\end{equation}
\begin{equation}
\label{fin4}
{\sf M}\left\{\left(
J^{*}[\psi^{(4)}]_{T,t}-
\sum\limits_{j_1, j_2, j_3, j_4=0}^{p}
C_{j_4 j_3 j_2 j_1}\zeta_{j_1}^{(i_1)}\zeta_{j_2}^{(i_2)}\zeta_{j_3}^{(i_3)}
\zeta_{j_4}^{(i_4)}
\right)^2\right\}
\le \frac{C}{p^{1-\varepsilon}}
\end{equation}

\vspace{1mm}
\noindent
are fulfilled, where $i_1, \ldots , i_4=0,1,\ldots,m$ in {\rm (\ref{fin0}),} {\rm (\ref{fin3})} 
and $i_1, \ldots, i_4=1,\ldots,m$ in {\rm (\ref{fin4}),}
constant $C$ does not depend on $p;$
$\varepsilon$ is an arbitrary
small positive real number 
for the case of complete orthonormal system of 
Legendre polynomials in the space $L_2([t, T])$
and $\varepsilon=0$ for the case of
complete orthonormal system of 
trigonometric functions in the space $L_2([t, T]),$
$$
C_{j_4 j_3 j_2 j_1}=
$$
$$
=
\int\limits_t^T\psi_4(t_4)\phi_{j_4}(t_4)
\int\limits_t^{t_4}\psi_3(t_3)\phi_{j_3}(t_3)
\int\limits_t^{t_3}\psi_2(t_2)\phi_{j_2}(t_2)
\int\limits_t^{t_2}\psi_1(t_1)\phi_{j_1}(t_1)dt_1dt_2dt_3dt_4;
$$
another notations are the same as in Theorem~{\rm 4}.}

{\bf Theorem 6}\ \cite{58}, \cite{add1014}, \cite{add1019}, \cite{add1020}, \cite{new-art-1-xxy}.\
{\it Assume 
that $\{\phi_j(x)\}_{j=0}^{\infty}$ is a complete orthonormal system of 
Legendre polynomials or trigonometric functions in the space $L_2([t, T])$
and $\psi_1(\tau), \ldots,$ $\psi_5(\tau)$ are continuously dif\-ferentiable 
nonrandom functions on $[t, T].$ 
Then, for the 
iterated Stra\-to\-no\-vich stochastic integral of fifth multiplicity
$$
J^{*}[\psi^{(5)}]_{T,t}=
$$
\begin{equation}
\label{fin7}
={\int\limits_t^{*}}^T\hspace{-1mm}\psi_5(t_5)
{\int\limits_t^{*}}^{t_5}\hspace{-1mm}\psi_4(t_4)
{\int\limits_t^{*}}^{t_4}\hspace{-1mm}\psi_3(t_3)
{\int\limits_t^{*}}^{t_3}\hspace{-1mm}\psi_2(t_2)
{\int\limits_t^{*}}^{t_2}\hspace{-1mm}\psi_1(t_1)
d{\bf w}_{t_1}^{(i_1)}
d{\bf w}_{t_2}^{(i_2)}d{\bf w}_{t_3}^{(i_3)}d{\bf w}_{t_4}^{(i_4)}
d{\bf w}_{t_5}^{(i_5)}
\end{equation}
the following 
relations
\begin{equation}
\label{fin8}
J^{*}[\psi^{(5)}]_{T,t}
=\hbox{\vtop{\offinterlineskip\halign{
\hfil#\hfil\cr
{\rm l.i.m.}\cr
$\stackrel{}{{}_{p\to \infty}}$\cr
}} }
\sum\limits_{j_1, j_2, j_3, j_4, j_5=0}^{p}
C_{j_5 j_4 j_3 j_2 j_1}\zeta_{j_1}^{(i_1)}\zeta_{j_2}^{(i_2)}
\zeta_{j_3}^{(i_3)}\zeta_{j_4}^{(i_4)}
\zeta_{j_5}^{(i_5)},
\end{equation}
$$
{\sf M}\left\{\left(
J^{*}[\psi^{(5)}]_{T,t}-
\sum\limits_{j_1, j_2, j_3, j_4, j_5=0}^{p}
C_{j_5 j_4 j_3 j_2 j_1}\zeta_{j_1}^{(i_1)}\zeta_{j_2}^{(i_2)}
\zeta_{j_3}^{(i_3)}\zeta_{j_4}^{(i_4)}
\zeta_{j_5}^{(i_5)}
\right)^2\right\}
\le 
$$
\begin{equation}
\label{fin9}
\le\frac{C}{p^{1-\varepsilon}}
\end{equation}

\vspace{1mm}
\noindent
are fulfilled, where $i_1, \ldots , i_5=0,1,\ldots,m$ in {\rm (\ref{fin7}),} {\rm (\ref{fin8})} 
and $i_1, \ldots, i_5=1,\ldots,m$ in {\rm (\ref{fin9}),}
constant $C$ is independent of $p;$
$\varepsilon$ is an arbitrary
small positive real number 
for the case of complete orthonormal system of 
Legendre polynomials in the space $L_2([t, T])$
and $\varepsilon=0$ for the case of
complete orthonormal system of 
trigonometric functions in the space $L_2([t, T]),$
$$
C_{j_5 j_4 j_3 j_2 j_1}=
\int\limits_t^T\psi_5(t_5)\phi_{j_5}(t_5)
\int\limits_t^{t_5}\psi_4(t_4)\phi_{j_4}(t_4)
\int\limits_t^{t_4}\psi_3(t_3)\phi_{j_3}(t_3)\times
$$
$$
\times
\int\limits_t^{t_3}\psi_2(t_2)\phi_{j_2}(t_2)
\int\limits_t^{t_2}\psi_1(t_1)\phi_{j_1}(t_1)
dt_1dt_2dt_3dt_4dt_5;
$$
another notations are the same as in Theorems~{\rm 4, 5}.}

{\bf Theorem 7}\ \cite{58}, \cite{add1014}, \cite{add1019}, \cite{add1020}.\
{\it Suppose that 
$\{\phi_j(x)\}_{j=0}^{\infty}$ is a complete orthonormal system of 
Legendre polynomials or trigonometric functions in the space $L_2([t, T]).$
Then, for the 
iterated Stratonovich stochastic integral of sixth multiplicity
$$
I_{T,t}^{*(i_1 i_2 i_3 i_4 i_5 i_6)}=
$$
$$
={\int\limits_t^{*}}^T
{\int\limits_t^{*}}^{t_6}
{\int\limits_t^{*}}^{t_5}
{\int\limits_t^{*}}^{t_4}
{\int\limits_t^{*}}^{t_3}
{\int\limits_t^{*}}^{t_2}
d{\bf w}_{t_1}^{(i_1)}
d{\bf w}_{t_2}^{(i_2)}
d{\bf w}_{t_3}^{(i_3)}
d{\bf w}_{t_4}^{(i_4)}
d{\bf w}_{t_5}^{(i_5)}
d{\bf w}_{t_6}^{(i_6)}
$$

\noindent
the following 
expansion 
$$
I_{T,t}^{*(i_1 i_2 i_3 i_4 i_5 i_6)}=
$$
$$
=
\hbox{\vtop{\offinterlineskip\halign{
\hfil#\hfil\cr
{\rm l.i.m.}\cr
$\stackrel{}{{}_{p\to \infty}}$\cr
}} }
\sum\limits_{j_1, j_2, j_3, j_4, j_5, j_6=0}^{p}
C_{j_6 j_5 j_4 j_3 j_2 j_1}\zeta_{j_1}^{(i_1)}
\zeta_{j_2}^{(i_2)}\zeta_{j_3}^{(i_3)}\zeta_{j_4}^{(i_4)}\zeta_{j_5}^{(i_5)}
\zeta_{j_6}^{(i_6)}
$$

\noindent
that converges in the mean-square sense is valid, where
$i_1, \ldots, i_6=0, 1,\ldots,m,$
$$
C_{j_6 j_5 j_4 j_3 j_2 j_1}=
$$
$$
=
\int\limits_t^T\phi_{j_6}(t_6)
\int\limits_t^{t_6}\phi_{j_5}(t_5)
\int\limits_t^{t_5}\phi_{j_4}(t_4)
\int\limits_t^{t_4}\phi_{j_3}(t_3)
\int\limits_t^{t_3}\phi_{j_2}(t_2)
\int\limits_t^{t_2}\phi_{j_1}(t_1)dt_1 dt_2 dt_3 dt_4 dt_5 dt_6;
$$

\noindent
another notations are the same as in Theorems~{\rm 4--6}.}

Consider the following Hypothesis on expansion of the iterated 
Stratonovich stochastic integrals (\ref{str}) of 
arbitrary multiplicity $k$ ($k\in\mathbb{N}$).

{\bf Hypothesis 1} \cite{26}, \cite{26a}, \cite{58}, \cite{59}.
{\it Assume that
$\{\phi_j(x)\}_{j=0}^{\infty}$ is a complete orthonormal
system of Legendre polynomials or trigonometric functions
in the space $L_2([t, T])$. Moreover,
every $\psi_l(\tau)$ $(l=1,\ldots,k)$ is 
an enough smooth nonrandom function
on $[t,T].$
Then, for the iterated Stratonovich stochastic integral {\rm (\ref{str})}
of multiplicity $k$
the following 
expansion 
\begin{equation}
\label{feto1900otitarrr}
J^{*}[\psi^{(k)}]_{T,t}=
\hbox{\vtop{\offinterlineskip\halign{
\hfil#\hfil\cr
{\rm l.i.m.}\cr
$\stackrel{}{{}_{p\to \infty}}$\cr
}} }
\sum\limits_{j_1,\ldots j_k=0}^{p}
C_{j_k \ldots j_1}\prod\limits_{l=1}^k \zeta_{j_l}^{(i_l)}
\end{equation}
that converges in the mean-square sense is valid, where 
notations are the same as in Theorems {\rm 1, 2}.
}

Hypothesis 1 allows to approximate the iterated
Stratonovich stochastic integral 
$J^{*}[\psi^{(k)}]_{T,t}$
by the sum
\begin{equation}
\label{otit567r}
J^{*}[\psi^{(k)}]_{T,t}^p=
\sum\limits_{j_1,\ldots j_k=0}^{p}
C_{j_k \ldots j_1}\prod\limits_{l=1}^k
\zeta_{j_l}^{(i_l)},
\end{equation}
where
$$
\lim_{p\to\infty}{\sf M}\left\{\Biggl(
J^{*}[\psi^{(k)}]_{T,t}-
J^{*}[\psi^{(k)}]_{T,t}^p\Biggl)^2\right\}=0.
$$

Note that Hypothesis~1 is proved in \cite{58} (Sect.~2.10)
under the condition of convergence of trace series
(also see \cite{add1014}, 
\cite{add1019}, \cite{add1020}).

Recently (in 2024-2025), the above approach has been generalized
to the case of an arbitrary complete orthonormal
system of functions in the space $L_2([t, T])$
\cite{58} (Sect.~2.1.4, 2.18--2.37), 
\cite{add1014} (Sect.~22--40), 
\cite{add1019} (Sect.~14--31), \cite{add1020} (Sect.~18--36),
\cite{2024xx}, \cite{2024xxx}.
In particular, the following five theorems were proved in 
\cite{58}, \cite{add1014}, 
\cite{add1019}, \cite{add1020}, \cite{2024xx}, \cite{2024xxx}.

{\bf Theorem~A}\ \cite{58}, \cite{add1014}, 
\cite{add1019}, \cite{add1020}, \cite{2024xx}.\ 
{\it Suppose that 
$\{\phi_j(x)\}_{j=0}^{\infty}$ is an arbitrary complete orthonormal system of 
functions in the space $L_2([t, T]).$
Moreover$,$ $\psi_1(\tau), \psi_2(\tau)$ are continuous 
functions on $[t, T].$ 
Then$,$ 
for the iterated Stra\-to\-novich stochastic integral of second
multiplicity
$$
J^{*}[\psi^{(2)}]_{T,t}={\int\limits_t^{*}}^T\psi_2(t_2)
{\int\limits_t^{*}}^{t_2}\psi_1(t_1)d{\bf w}_{t_1}^{(i_1)}
d{\bf w}_{t_2}^{(i_2)}\ \ \ (i_1, i_2=1,\ldots,m)
$$
the following expansion  
$$
J^{*}[\psi^{(2)}]_{T,t}=\hbox{\vtop{\offinterlineskip\halign{
\hfil#\hfil\cr
{\rm l.i.m.}\cr
$\stackrel{}{{}_{p_1,p_2\to \infty}}$\cr
}} }\sum_{j_1=0}^{p_1}\sum_{j_2=0}^{p_2}
C_{j_2j_1}\zeta_{j_1}^{(i_1)}\zeta_{j_2}^{(i_2)}
$$
that converges in the mean-square
sence is valid$,$ where the notations are the same as in Theorems {\rm 1--3.}
}

{\bf Theorem~B}\ \cite{58}, \cite{add1014}, 
\cite{add1019}, \cite{add1020}, \cite{2024xxx}.\ 
{\it Suppose that
$\{\phi_j(x)\}_{j=0}^{\infty}$ is an arbitrary complete orthonormal system of 
functions in the space $L_2([t,T]).$
Then$,$ for the iterated Stra\-to\-no\-vich stochastic integral
of third multiplicity 
$$
I_{(l_1l_2l_3)T,t}^{*(i_1i_2i_3)}={\int\limits_t^{*}}^T (t-t_3)^{l_3}
{\int\limits_t^{*}}^{t_3}(t-t_2)^{l_2}
{\int\limits_t^{*}}^{t_2}(t-t_1)^{l_1}
d{\bf w}_{t_1}^{(i_1)}
d{\bf w}_{t_2}^{(i_2)}d{\bf w}_{t_3}^{(i_3)}
$$
the following expansion 
$$
I_{(l_1l_2l_3)T,t}^{*(i_1i_2i_3)}=
\hbox{\vtop{\offinterlineskip\halign{
\hfil#\hfil\cr
{\rm l.i.m.}\cr
$\stackrel{}{{}_{p\to \infty}}$\cr
}} }\sum_{j_1,j_2,j_3=0}^{p}
C_{j_3 j_2 j_1}\zeta_{j_1}^{(i_1)}\zeta_{j_2}^{(i_2)}\zeta_{j_3}^{(i_3)}
$$
that converges in the mean-square sense is valid, where 
$i_1,i_2,i_3=0,1,\ldots,m;$ $l_1,l_2,l_3=0,1,2,\ldots,$
$$
C_{j_3 j_2 j_1}=\int\limits_t^T
(t-t_3)^{l_3}\phi_{j_3}(t_3)\int\limits_t^{t_3}
(t-t_2)^{l_2}
\phi_{j_2}(t_2)
\int\limits_t^{t_2}
(t-t_1)^{l_1}\phi_{j_1}(t_1)dt_1dt_2dt_3
$$
and
$$
\zeta_{j}^{(i)}=
\int\limits_t^T \phi_{j}(\tau) d{\bf w}_{\tau}^{(i)}
$$ 
are independent standard Gaussian random variables for various 
$i$ or $j$ {\rm (}in the case when $i\ne 0${\rm ),}
${\bf w}_{\tau}^{(0)}=\tau.$}

{\bf Theorem~C}\ \cite{58}, \cite{add1014}, 
\cite{add1019}, \cite{add1020}, \cite{2024xxx}.\ 
{\it Suppose that
$\{\phi_j(x)\}_{j=0}^{\infty}$ is an arbitrary complete orthonormal system of 
functions in the space $L_2([t,T]).$
Then$,$ for the iterated Stra\-to\-no\-vich stochastic integral
of fourth multiplicity 
$$
I_{(l_1l_2l_3 l_4)T,t}^{*(i_1i_2i_3 i_4)}=
{\int\limits_t^{*}}^T (t-t_4)^{l_4}{\int\limits_t^{*}}^{t_4} (t-t_3)^{l_3}
{\int\limits_t^{*}}^{t_3}(t-t_2)^{l_2}
{\int\limits_t^{*}}^{t_2}(t-t_1)^{l_1}\times
$$
$$
\times
d{\bf w}_{t_1}^{(i_1)}
d{\bf w}_{t_2}^{(i_2)}d{\bf w}_{t_3}^{(i_3)}d{\bf w}_{t_4}^{(i_4)}
$$
the following expansion 
$$
I_{(l_1l_2l_3 l_4)T,t}^{*(i_1i_2i_3i_4)}=
\hbox{\vtop{\offinterlineskip\halign{
\hfil#\hfil\cr
{\rm l.i.m.}\cr
$\stackrel{}{{}_{p\to \infty}}$\cr
}} }\sum_{j_1,j_2,j_3,j_4=0}^{p}
C_{j_4 j_3 j_2 j_1}\zeta_{j_1}^{(i_1)}\zeta_{j_2}^{(i_2)}\zeta_{j_3}^{(i_3)}\zeta_{j_4}^{(i_4)}
$$
that converges in the mean-square sense is valid, where 
$i_1,i_2,i_3,i_4=0,1,\ldots,m;$ $l_1,l_2,l_3,l_4=0,1,2,\ldots,$
$$
C_{j_4 j_3 j_2 j_1}=\int\limits_t^T
(t-t_4)^{l_4}\phi_{j_4}(t_4)\int\limits_t^{t_4}
(t-t_3)^{l_3}\phi_{j_3}(t_3)\int\limits_t^{t_3}
(t-t_2)^{l_2}
\phi_{j_2}(t_2)
\int\limits_t^{t_2}
(t-t_1)^{l_1}\phi_{j_1}(t_1)\times
$$
$$
\times
dt_1dt_2dt_3dt_4;
$$
another notations are the same as in Theorem~{\rm B}}.

{\bf Theorem~D}\ \cite{58}, \cite{add1014}, 
\cite{add1019}, \cite{add1020}, \cite{2024xxx}.\ 
{\it Suppose that
$\{\phi_j(x)\}_{j=0}^{\infty}$ is an arbitrary complete orthonormal system of 
functions in the space $L_2([t,T]).$
Then$,$ for the iterated Stra\-to\-no\-vich stochastic integral
of fifth multiplicity 
$$
I_{T,t}^{*(i_1 i_2 i_3 i_4 i_5)}=
{\int\limits_t^{*}}^T
{\int\limits_t^{*}}^{t_5}
{\int\limits_t^{*}}^{t_4}
{\int\limits_t^{*}}^{t_3}
{\int\limits_t^{*}}^{t_2}
d{\bf w}_{t_1}^{(i_1)}
d{\bf w}_{t_2}^{(i_2)}d{\bf w}_{t_3}^{(i_3)}d{\bf w}_{t_4}^{(i_4)}
d{\bf w}_{t_5}^{(i_5)}
$$
the following 
expansion 
$$
I_{T,t}^{*(i_1 i_2 i_3 i_4 i_5)}=
\hbox{\vtop{\offinterlineskip\halign{
\hfil#\hfil\cr
{\rm l.i.m.}\cr
$\stackrel{}{{}_{p\to \infty}}$\cr
}} }
\sum\limits_{j_1, j_2, j_3, j_4, j_5=0}^{p}
C_{j_5 j_4 j_3 j_2 j_1}\zeta_{j_1}^{(i_1)}\zeta_{j_2}^{(i_2)}
\zeta_{j_3}^{(i_3)}\zeta_{j_4}^{(i_4)}
\zeta_{j_5}^{(i_5)}
$$
that converges in the mean-square sense is valid, where 
$i_1,\ldots,i_5=0, 1,\ldots,m,$
$$
C_{j_5 j_4 j_3 j_2 j_1}=
\int\limits_t^T\phi_{j_5}(t_5)
\int\limits_t^{t_5}\phi_{j_4}(t_4)
\int\limits_t^{t_4}\phi_{j_3}(t_3)
\int\limits_t^{t_3}\phi_{j_2}(t_2)
\int\limits_t^{t_2}\phi_{j_1}(t_1)dt_1dt_2dt_3dt_4dt_5;
$$
another notations are the same as in Theorem~{\rm B}}.

{\bf Theorem~E}\ \cite{58}, \cite{add1014}, 
\cite{add1019}, \cite{add1020}.\ 
{\it Suppose that
$\{\phi_j(x)\}_{j=0}^{\infty}$ is an arbitrary complete orthonormal system of 
functions in the space $L_2([t,T]).$
Then$,$ for the iterated Stra\-to\-no\-vich stochastic integral
of sixth multiplicity 
$$
I_{T,t}^{*(i_1 i_2 i_3 i_4 i_5 i_6)}
={\int\limits_t^{*}}^T
{\int\limits_t^{*}}^{t_6}
{\int\limits_t^{*}}^{t_5}
{\int\limits_t^{*}}^{t_4}
{\int\limits_t^{*}}^{t_3}
{\int\limits_t^{*}}^{t_2}
d{\bf w}_{t_1}^{(i_1)}
d{\bf w}_{t_2}^{(i_2)}
d{\bf w}_{t_3}^{(i_3)}
d{\bf w}_{t_4}^{(i_4)}
d{\bf w}_{t_5}^{(i_5)}
d{\bf w}_{t_6}^{(i_6)}
$$

\noindent
the following 
expansion 
$$
I_{T,t}^{*(i_1 i_2 i_3 i_4 i_5 i_6)}=
\hbox{\vtop{\offinterlineskip\halign{
\hfil#\hfil\cr
{\rm l.i.m.}\cr
$\stackrel{}{{}_{p\to \infty}}$\cr
}} }
\sum\limits_{j_1, j_2, j_3, j_4, j_5, j_6=0}^{p}
C_{j_6 j_5 j_4 j_3 j_2 j_1}\zeta_{j_1}^{(i_1)}
\zeta_{j_2}^{(i_2)}\zeta_{j_3}^{(i_3)}\zeta_{j_4}^{(i_4)}\zeta_{j_5}^{(i_5)}
\zeta_{j_6}^{(i_6)}
$$

\noindent
that converges in the mean-square sense is valid, where
$i_1, \ldots, i_6=0, 1,\ldots,m;$
another notations are the same as in Theorem~{\rm 7}.}

Note that Theorem~7 is generalized \cite{58} to the case $k=7$ and $k=8.$

In connection with Theorems~A--E, Hypothesis~1 was revisited 
in \cite{58}, \cite{add1014}, 
\cite{add1019}, \cite{add1020}.
Namely, the following hypothesis was formulated.

{\bf Hypothesis 2}\ \cite{58} (Sect.~2.28), \cite{add1014}, 
\cite{add1019}, \cite{add1020}.\
{\it Assume that
$\{\phi_j(x)\}_{j=0}^{\infty}$ is an arbitrary complete orthonormal system of 
functions in the space $L_2([t,T])$
and $\psi_1(\tau),\ldots,\psi_k(\tau)$ are 
continuous functions on
$[t,T].$
Then, for the iterated Stratonovich stochastic integral {\rm (\ref{str})}
of multiplicity $k$
the following 
expansion 
$$
J^{*}[\psi^{(k)}]_{T,t}=
\hbox{\vtop{\offinterlineskip\halign{
\hfil#\hfil\cr
{\rm l.i.m.}\cr
$\stackrel{}{{}_{p\to \infty}}$\cr
}} }
\sum\limits_{j_1,\ldots j_k=0}^{p}
C_{j_k \ldots j_1}\prod\limits_{l=1}^k \zeta_{j_l}^{(i_l)}
$$
that converges in the mean-square sense is valid, where 
notations are the same as in Theorems {\rm 1, 2}.
}

Note that Hypothesis~2 has been proved in 
\cite{58} (Theorems~59, 61), \cite{add1014}, 
\cite{add1019}, \cite{add1020} but under one additional
condition.

Suppose that $J[\psi^{(k)}]_{T,t}^{p}$ is the approximation 
of (\ref{ito}), which is
the expression on the right-hand side of (\ref{leto6000}) before passing to the limit
for the case $p_1=\ldots =p_k=p$, i.e.
$$
J[\psi^{(k)}]_{T,t}^{p}=
\sum\limits_{j_1,\ldots, j_k=0}^{p}
C_{j_k\ldots j_1}\Biggl(
\prod_{l=1}^k\zeta_{j_l}^{(i_l)}+\sum\limits_{r=1}^{[k/2]}
(-1)^r \times
\Biggr.
$$
\begin{equation}
\label{letusdenote}
\times
\sum_{\stackrel{(\{\{g_1, g_2\}, \ldots, 
\{g_{2r-1}, g_{2r}\}\}, \{q_1, \ldots, q_{k-2r}\})}
{{}_{\{g_1, g_2, \ldots, 
g_{2r-1}, g_{2r}, q_1, \ldots, q_{k-2r}\}=\{1, 2, \ldots, k\}}}}
\prod\limits_{s=1}^r
{\bf 1}_{\{i_{g_{{}_{2s-1}}}=~i_{g_{{}_{2s}}}\ne 0\}}
\Biggl.{\bf 1}_{\{j_{g_{{}_{2s-1}}}=~j_{g_{{}_{2s}}}\}}
\prod_{l=1}^{k-2r}\zeta_{j_{q_l}}^{(i_{q_l})}\Biggr),
\end{equation}

\vspace{1mm}
\noindent
where $[x]$ is an integer part of a real number $x;$
another notations are the same as in Theorems~{\rm 1, 2}.

Let us denote
$$
{\sf M}\left\{\left(J[\psi^{(k)}]_{T,t}-
J[\psi^{(k)}]_{T,t}^{p}\right)^2\right\}\stackrel{{\rm def}}
{=}E_k^{p},
$$
$$
\left\Vert K\right\Vert^2_{L_2([t,T]^k)}=\int\limits_{[t,T]^k}
K^2(t_1,\ldots,t_k)dt_1\ldots dt_k\stackrel{{\rm def}}{=}I_k.
$$

For the futher consideration,
we need the following useful estimate
\cite{26}, \cite{26a}, \cite{58} 
\begin{equation}
\label{star00011}
E_k^p
\le k!\left(I_k
-\sum_{j_1,\ldots,j_k=0}^{p}C_{j_k\ldots j_1}^2\right),
\end{equation}

\vspace{1mm}
\noindent
where $i_1,\ldots,i_k=1,\ldots,m$ for $T-t\in (0,\infty)$ and 
$i_1,\ldots,i_k=0, 1,\ldots,m$ for $T-t\in (0, 1);$
another notations are the same as in Theorems 1, 2.

The value $E_k^{p}$
can be calculated exactly.

{\bf Theorem 8} \cite{58} (Sect.~1.12), \cite{add1012} (Sect.~6).
{\it Suppose that $\{\phi_j(x)\}_{j=0}^{\infty}$ 
is an arbitrary complete orthonormal system  
of functions in the space $L_2([t,T])$ and
$\psi_1(\tau),\ldots,\psi_k(\tau)\in L_2([t, T]),$  $i_1,\ldots, i_k=1,\ldots,m$.
Then
$$
E_k^p=I_k- 
$$
\begin{equation}
\label{tttr11}
-\sum_{j_1,\ldots, j_k=0}^{p}
C_{j_k\ldots j_1}
{\sf M}\left\{J[\psi^{(k)}]_{T,t}
\sum\limits_{(j_1,\ldots,j_k)}
\int\limits_t^T \phi_{j_k}(t_k)
\ldots
\int\limits_t^{t_{2}}\phi_{j_{1}}(t_{1})
d{\bf w}_{t_1}^{(i_1)}\ldots
d{\bf w}_{t_k}^{(i_k)}\right\},
\end{equation}
where $i_1,\ldots,i_k = 1,\ldots,m;$
the expression 
$$
\sum\limits_{(j_1,\ldots,j_k)}
$$ 

\noindent
means the sum with respect to all
possible permutations 
$(j_1,\ldots,j_k)$. At the same time if 
$j_r$ swapped with $j_q$ in the permutation $(j_1,\ldots,j_k),$
then $i_r$ swapped with $i_q$ in the permutation
$(i_1,\ldots,i_k);$
another notations are the same as in Theorems {\rm 1, 2.}
}

Note that 
$$
{\sf M}\left\{J[\psi^{(k)}]_{T,t}
\int\limits_t^T \phi_{j_k}(t_k)
\ldots
\int\limits_t^{t_{2}}\phi_{j_{1}}(t_{1})
d{\bf w}_{t_1}^{(i_1)}\ldots
d{\bf w}_{t_k}^{(i_k)}\right\}=C_{j_k\ldots j_1}.
$$

Then from Theorem 8
we obtain
\begin{equation}
\label{formula0}
E_k^p= I_k- \sum_{j_1,\ldots,j_k=0}^{p}
C_{j_k\ldots j_1}^2\ \ \ (i_1,\ldots,i_k\ \hbox{are pairwise different}),
\end{equation}
\begin{equation}
\label{formula0xx}
E_k^p= I_k - \sum_{j_1,\ldots,j_k=0}^{p}
C_{j_k\ldots j_1}\left(\sum\limits_{(j_1,\ldots,j_k)}
C_{j_k\ldots j_1}\right)\ \ \ (i_1=\ldots=i_k).
\end{equation}

Consider some examples of the application of Theorem 8
$(i_1,\ldots ,i_5=1,\ldots,m)$:
\begin{equation}
\label{formula0yy}
E_2^p
=I_2
-\sum_{j_1,j_2=0}^p
C_{j_2j_1}^2-
\sum_{j_1,j_2=0}^p
C_{j_2j_1}C_{j_1j_2}\ \ \ (i_1=i_2),
\end{equation}
\begin{equation}
\label{formula1}
E_3^p=I_3
-\sum_{j_3,j_2,j_1=0}^p C_{j_3j_2j_1}^2-
\sum_{j_3,j_2,j_1=0}^p C_{j_3j_1j_2}C_{j_3j_2j_1}\ \ \ (i_1=i_2\ne i_3),
\end{equation}

\vspace{-3mm}
\begin{equation}
\label{formula2}
E_3^p=I_3-
\sum_{j_3,j_2,j_1=0}^p C_{j_3j_2j_1}^2-
\sum_{j_3,j_2,j_1=0}^p C_{j_2j_3j_1}C_{j_3j_2j_1}\ \ \ (i_1\ne i_2=i_3),
\end{equation}

\vspace{-3mm}
\begin{equation}
\label{formula3}
E_3^p=I_3
-\sum_{j_3,j_2,j_1=0}^p C_{j_3j_2j_1}^2-
\sum_{j_3,j_2,j_1=0}^p C_{j_3j_2j_1}C_{j_1j_2j_3}\ \ \ (i_1=i_3\ne i_2),
\end{equation}

\vspace{-3mm}
\begin{equation}
\label{formula4}
E_4^p=I_4 - \sum_{j_1,\ldots,j_4=0}^{p}
C_{j_4 \ldots j_1}\left(\sum\limits_{(j_1,j_2)}
C_{j_4 \ldots j_1}\right)\ \ \ (i_1=i_2\ne i_3, i_4;\ i_3\ne i_4),
\end{equation}

\vspace{-3mm}
\begin{equation}
\label{formula5}
E_4^p=I_4 - \sum_{j_1,\ldots,j_4=0}^{p}
C_{j_4 \ldots j_1}\left(\sum\limits_{(j_1,j_3)}
C_{j_4 \ldots j_1}\right)\ \ \ (i_1=i_3\ne i_2, i_4;\ i_2\ne i_4),
\end{equation}

\vspace{-3mm}
\begin{equation}
\label{formula6}
E_4^p=I_4 - \sum_{j_1,\ldots,j_4=0}^{p}
C_{j_4 \ldots j_1}\left(\sum\limits_{(j_2,j_3)}
C_{j_4 \ldots j_1}\right)\ \ \ (i_2=i_3\ne i_1, i_4;\ i_1\ne i_4),
\end{equation}

\vspace{-3mm}
\begin{equation}
\label{formula7}
E_4^p=I_4 - \sum_{j_1,\ldots,j_4=0}^{p}
C_{j_4 \ldots j_1}\left(\sum\limits_{(j_1,j_4)}
C_{j_4 \ldots j_1}\right)\ \ \ (i_1=i_4\ne i_2, i_3;\ i_2\ne i_3),
\end{equation}

\vspace{-3mm}
\begin{equation}
\label{mark1}
E_4^p=I_4 - \sum_{j_1,\ldots,j_4=0}^{p}
C_{j_4 \ldots j_1}\left(\sum\limits_{(j_1,j_4)}\left(\sum\limits_{(j_2,j_3)}
C_{j_4 \ldots j_1}\right)\right)\ \ \ (i_1=i_4\ne i_2=i_3),
\end{equation}

\vspace{-3mm}
\begin{equation}
\label{mark2}
E_4^p=I_4 - \sum_{j_1,\ldots,j_4=0}^{p}
C_{j_4 \ldots j_1}\left(\sum\limits_{(j_1,j_2,j_3)}
C_{j_4 \ldots j_1}\right)\ \ \ (i_1=i_2=i_3\ne i_4),
\end{equation}

\vspace{-3mm}
\begin{equation}
\label{mark3}
E_4^p = I_4 -
 \sum_{j_1,\ldots,j_4=0}^{p}
C_{j_4\ldots j_1}\left(\sum\limits_{(j_2,j_3,j_4)}
C_{j_4\ldots j_1}\right)\ \ \ (i_2=i_3=i_4\ne i_1),
\end{equation}
\begin{equation}
\label{mark4}
E_4^p = I_4 -
 \sum_{j_1,\ldots,j_4=0}^{p}
C_{j_4\ldots j_1}\left(\sum\limits_{(j_1,j_2,j_4)}
C_{j_4\ldots j_1}\right),\ \ \ (i_1=i_2=i_4\ne i_3),
\end{equation}

\vspace{-3mm}
\begin{equation}
\label{mark5}
E_4^p = I_4 -
 \sum_{j_1,\ldots,j_4=0}^{p}
C_{j_4\ldots j_1}\left(\sum\limits_{(j_1,j_3,j_4)}
C_{j_4\ldots j_1}\right)\ \ \ (i_1=i_3=i_4\ne i_2),
\end{equation}

\vspace{-3mm}
\begin{equation}
\label{formula8}
E_5^p=I_5 - \sum_{j_1,\ldots,j_5=0}^{p}
C_{j_5 \ldots j_1}\left(\sum\limits_{(j_1,j_2)}
C_{j_5 \ldots j_1}\right),
\end{equation}

\noindent
where $i_1=i_2\ne i_3, i_4, i_5$ and
$i_3, i_4, i_5$ are pairwise different,

\vspace{-1mm}
\begin{equation}
\label{formula9}
E_5^p=I_5 - \sum_{j_1,\ldots,j_5=0}^{p}
C_{j_5 \ldots j_1}\left(\sum\limits_{(j_2,j_3)}
C_{j_5 \ldots j_1}\right),
\end{equation}

\noindent
where $i_2=i_3\ne i_1, i_4, i_5$ and
$i_1, i_4, i_5$ are pairwise different,

\vspace{-1mm}
\begin{equation}
\label{mark9}
E^p_5 = I_5 - \sum_{j_1,\ldots,j_5=0}^{p}
C_{j_5\ldots j_1}\left(\sum\limits_{(j_1,j_3)}
C_{j_5\ldots j_1}\right),
\end{equation}

\noindent
where $i_1=i_3\ne i_2,i_4,i_5$ and $i_2,i_4,i_5$ are pairwise different,

\vspace{-1mm}
\begin{equation}
\label{mark10}
E^p_5 = I_5 - \sum_{j_1,\ldots,j_5=0}^{p}
C_{j_5\ldots j_1}\left(\sum\limits_{(j_1,j_4)}
C_{j_5\ldots j_1}\right),
\end{equation}

\noindent
where $i_1=i_4\ne i_2,i_3,i_5$ and $i_2,i_3,i_5$ are pairwise different,

\vspace{-1mm}
\begin{equation}
\label{mark11}
E^p_5 = I_5 - \sum_{j_1,\ldots,j_5=0}^{p}
C_{j_5\ldots j_1}\left(\sum\limits_{(j_1,j_5)}
C_{j_5\ldots j_1}\right),
\end{equation}

\noindent
where $i_1=i_5\ne i_2,i_3,i_4$ and $i_2,i_3,i_4$  are pairwise different.

\vspace{-1mm}
\begin{equation}
\label{mark12}
E^p_5 = I_5 - \sum_{j_1,\ldots,j_5=0}^{p}
C_{j_5\ldots j_1}\left(\sum\limits_{(j_2,j_4)}
C_{j_5\ldots j_1}\right),
\end{equation}

\noindent
where $i_2=i_4\ne i_1,i_3,i_5$ and $i_1,i_3,i_5$ are pairwise different.

\newpage
\noindent
\begin{equation}
\label{formula10}
E_5^p=I_5 - \sum_{j_1,\ldots,j_5=0}^{p}
C_{j_5 \ldots j_1}\left(\sum\limits_{(j_4,j_5)}
C_{j_5 \ldots j_1}\right),
\end{equation}
where $i_4=i_5\ne i_1, i_2, i_3$ and
$i_1, i_2, i_3$ are pairwise different,
\begin{equation}
\label{formula11}
E_5^p = I_5 - \sum_{j_1,\ldots,j_5=0}^{p}
C_{j_5\ldots j_1}\left(\sum\limits_{(j_2,j_4)}\left(
\sum\limits_{(j_3,j_5)}
C_{j_5\ldots j_1}\right)\right)\ \ \ (i_1\ne i_2=i_4\ne i_3=i_5\ne i_1).
\end{equation}

\subsection{Approximations of Iterated 
It\^{o} Stochastic 
Integrals from the Numerical Schemes (\ref{sm1})--(\ref{al5})
Using Legendre Polynomials}

This section is devoted to 
approximation of the iterated It\^{o} 
stochastic integrals (\ref{qqq1x})
of multiplicities 1 to 6 based on 
Theorems 1, 2. At that we will use
multiple Fourier--Legendre series for 
approximation of the mentioned stochastic integrals.

The numerical schemes (\ref{sm1})--(\ref{al5})
contain the following set (see (\ref{qqq1x}))
of iterated It\^{o} 
stochastic integrals 
\begin{equation}
\label{sm10}
I_{(0)T,t}^{(i_1)},\ \ \ I_{(1)T,t}^{(i_1)},\ \ \ I_{(2)T,t}^{(i_1)},\ \ \ 
I_{(00)T,t}^{(i_1 i_2)},\ \ \ I_{(10)T,t}^{(i_1 i_2)},\ \ \ 
I_{(01)T,t}^{(i_1 i_2)},\ \ \ I_{(000)T,t}^{(i_1 i_2 i_3)},\ \ \ 
I_{(0000)T,t}^{(i_1 i_2 i_3 i_4)},\ \ \
\end{equation}
\begin{equation}
\label{sm11}
I_{(00000)T,t}^{(i_1 i_2 i_3 i_4 i_5)},\ \ \ 
I_{(02)T,t}^{(i_1 i_2)},\ \ \ I_{(20)T,t}^{(i_1 i_2)},\ \ \ 
I_{(11)T,t}^{(i_1 i_2)},\ \ \ 
I_{(100)T,t}^{(i_1 i_2 i_3)},\ \ \ I_{(010)T,t}^{(i_1 i_2 i_3)},\ \ \ 
I_{(001)T,t}^{(i_1 i_2 i_3)},
\end{equation}
\begin{equation}
\label{sm12}
I_{(0001)T,t}^{(i_1 i_2 i_3 i_4)},\ \ \ 
I_{(0010)T,t}^{(i_1 i_2 i_3 i_4)},\ \ \ 
I_{(0100)T,t}^{(i_1 i_2 i_3 i_4)},\ \ \ I_{(1000)T,t}^{(i_1 i_2 i_3 i_4)},\ \ 
\ I_{(000000)T,t}^{(i_1 i_2 i_3 i_4 i_5 i_6)}.
\end{equation}

Let us consider 
the complete orthonormal system of Legendre polynomials in the 
space $L_2([t,T])$ 
\begin{equation}
\label{4009}
\phi_j(x)=\sqrt{\frac{2j+1}{T-t}}P_j\left(\left(
x-\frac{T+t}{2}\right)\frac{2}{T-t}\right),\ \ \ j=0, 1, 2,\ldots,
\end{equation}
where $P_j(x)$ is the Legendre polynomial
$$
P_j(x)=\frac{1}{2^j j!} \frac{d^j}{dx^j}\left(x^2-1\right)^j.
$$

Using Theorems 1, 2 and 
well known properties of the Legendre polynomials,
we obtain the following formulas for numerical 
modeling of the stochastic integrals
(\ref{sm10})--(\ref{sm12})
\cite{26}, \cite{26a}, \cite{40}-\cite{42aa},
\cite{58}, \cite{59},
\cite{60}-\cite{62}
\begin{equation}
\label{desy1}
I_{(0)T,t}^{(i_1)}=\sqrt{T-t}\zeta_0^{(i_1)},
\end{equation}
\begin{equation}
\label{desy2}
I_{(1)T,t}^{(i_1)}=-\frac{(T-t)^{3/2}}{2}\left(\zeta_0^{(i_1)}+
\frac{1}{\sqrt{3}}\zeta_1^{(i_1)}\right),
\end{equation}
\begin{equation}
\label{desy3}
I_{(2)T,t}^{(i_1)}=\frac{(T-t)^{5/2}}{3}\left(\zeta_0^{(i_1)}+
\frac{\sqrt{3}}{2}\zeta_1^{(i_1)}+
\frac{1}{2\sqrt{5}}\zeta_2^{(i_1)}\right),
\end{equation}
\begin{equation}
\label{nach1}
I_{(00)T,t}^{(i_1 i_2)q}=
\frac{T-t}{2}\left(\zeta_0^{(i_1)}\zeta_0^{(i_2)}+\sum_{i=1}^{q}
\frac{1}{\sqrt{4i^2-1}}\left(
\zeta_{i-1}^{(i_1)}\zeta_{i}^{(i_2)}-
\zeta_i^{(i_1)}\zeta_{i-1}^{(i_2)}\right) - {\bf 1}_{\{i_1=i_2\}}\right),
\end{equation}
$$
I_{(000)T,t}^{(i_1i_2i_3)q_1}
=
\sum_{j_1,j_2,j_3=0}^{q_1}
C_{j_3j_2j_1}^{000}
\Biggl(
\zeta_{j_1}^{(i_1)}\zeta_{j_2}^{(i_2)}\zeta_{j_3}^{(i_3)}
-{\bf 1}_{\{i_1=i_2\}}
{\bf 1}_{\{j_1=j_2\}}
\zeta_{j_3}^{(i_3)}-
\Biggr.
$$
\begin{equation}
\label{desy4}
\Biggl.
-{\bf 1}_{\{i_2=i_3\}}
{\bf 1}_{\{j_2=j_3\}}
\zeta_{j_1}^{(i_1)}-
{\bf 1}_{\{i_1=i_3\}}
{\bf 1}_{\{j_1=j_3\}}
\zeta_{j_2}^{(i_2)}\Biggr),
\end{equation}

\vspace{2mm}
\begin{equation}
\label{ogo1}
I_{(10)T,t}^{(i_1 i_2)q_2}=
\sum_{j_1,j_2=0}^{q_2}
C_{j_2j_1}^{10}\Biggl(\zeta_{j_1}^{(i_1)}\zeta_{j_2}^{(i_2)}
-{\bf 1}_{\{i_1=i_2\}}
{\bf 1}_{\{j_1=j_2\}}\Biggr),
\end{equation}
\begin{equation}
\label{ogo2}
I_{(01)T,t}^{(i_1 i_2)q_2}=
\sum_{j_1,j_2=0}^{q_2}
C_{j_2j_1}^{01}\Biggl(\zeta_{j_1}^{(i_1)}\zeta_{j_2}^{(i_2)}
-{\bf 1}_{\{i_1=i_2\}}
{\bf 1}_{\{j_1=j_2\}}\Biggr),
\end{equation}

\vspace{2mm}

$$
I_{(0000)T,t}^{(i_1 i_2 i_3 i_4)q_3}
=
\sum_{j_1,j_2,j_3,j_4=0}^{q_3}
C_{j_4 j_3 j_2 j_1}^{0000}\Biggl(
\zeta_{j_1}^{(i_1)}\zeta_{j_2}^{(i_2)}\zeta_{j_3}^{(i_3)}\zeta_{j_4}^{(i_4)}
-\Biggr.
$$
$$
-
{\bf 1}_{\{i_1=i_2\}}
{\bf 1}_{\{j_1=j_2\}}
\zeta_{j_3}^{(i_3)}
\zeta_{j_4}^{(i_4)}
-
{\bf 1}_{\{i_1=i_3\}}
{\bf 1}_{\{j_1=j_3\}}
\zeta_{j_2}^{(i_2)}
\zeta_{j_4}^{(i_4)}-
$$
$$
-
{\bf 1}_{\{i_1=i_4\}}
{\bf 1}_{\{j_1=j_4\}}
\zeta_{j_2}^{(i_2)}
\zeta_{j_3}^{(i_3)}
-
{\bf 1}_{\{i_2=i_3\}}
{\bf 1}_{\{j_2=j_3\}}
\zeta_{j_1}^{(i_1)}
\zeta_{j_4}^{(i_4)}-
$$
$$
-
{\bf 1}_{\{i_2=i_4\}}
{\bf 1}_{\{j_2=j_4\}}
\zeta_{j_1}^{(i_1)}
\zeta_{j_3}^{(i_3)}
-
{\bf 1}_{\{i_3=i_4\}}
{\bf 1}_{\{j_3=j_4\}}
\zeta_{j_1}^{(i_1)}
\zeta_{j_2}^{(i_2)}+
$$
$$
+
{\bf 1}_{\{i_1=i_2\}}
{\bf 1}_{\{j_1=j_2\}}
{\bf 1}_{\{i_3=i_4\}}
{\bf 1}_{\{j_3=j_4\}}+
{\bf 1}_{\{i_1=i_3\}}
{\bf 1}_{\{j_1=j_3\}}
{\bf 1}_{\{i_2=i_4\}}
{\bf 1}_{\{j_2=j_4\}}+
$$
\begin{equation}
\label{desy5}
+\Biggl.
{\bf 1}_{\{i_1=i_4\}}
{\bf 1}_{\{j_1=j_4\}}
{\bf 1}_{\{i_2=i_3\}}
{\bf 1}_{\{j_2=j_3\}}\Biggr),
\end{equation}

\vspace{2mm}

$$
I_{(00000)T,t}^{(i_1 i_2 i_3 i_4 i_5)q_4}
=
\sum_{j_1,j_2,j_3,j_4,j_5=0}^{q_4}
C_{j_5 j_4 j_3 j_2 j_1}^{00000}\Biggl(
\prod_{l=1}^5\zeta_{j_l}^{(i_l)}
-\Biggr.
$$
$$
-
{\bf 1}_{\{i_1=i_2\}}
{\bf 1}_{\{j_1=j_2\}}
\zeta_{j_3}^{(i_3)}
\zeta_{j_4}^{(i_4)}
\zeta_{j_5}^{(i_5)}-
{\bf 1}_{\{i_1=i_3\}}
{\bf 1}_{\{j_1=j_3\}}
\zeta_{j_2}^{(i_2)}
\zeta_{j_4}^{(i_4)}
\zeta_{j_5}^{(i_5)}-
$$
$$
-
{\bf 1}_{\{i_1=i_4\}}
{\bf 1}_{\{j_1=j_4\}}
\zeta_{j_2}^{(i_2)}
\zeta_{j_3}^{(i_3)}
\zeta_{j_5}^{(i_5)}-
{\bf 1}_{\{i_1=i_5\}}
{\bf 1}_{\{j_1=j_5\}}
\zeta_{j_2}^{(i_2)}
\zeta_{j_3}^{(i_3)}
\zeta_{j_4}^{(i_4)}-
$$
$$
-
{\bf 1}_{\{i_2=i_3\}}
{\bf 1}_{\{j_2=j_3\}}
\zeta_{j_1}^{(i_1)}
\zeta_{j_4}^{(i_4)}
\zeta_{j_5}^{(i_5)}-
{\bf 1}_{\{i_2=i_4\}}
{\bf 1}_{\{j_2=j_4\}}
\zeta_{j_1}^{(i_1)}
\zeta_{j_3}^{(i_3)}
\zeta_{j_5}^{(i_5)}-
$$
$$
-
{\bf 1}_{\{i_2=i_5\}}
{\bf 1}_{\{j_2=j_5\}}
\zeta_{j_1}^{(i_1)}
\zeta_{j_3}^{(i_3)}
\zeta_{j_4}^{(i_4)}
-{\bf 1}_{\{i_3=i_4\}}
{\bf 1}_{\{j_3=j_4\}}
\zeta_{j_1}^{(i_1)}
\zeta_{j_2}^{(i_2)}
\zeta_{j_5}^{(i_5)}-
$$
$$
-
{\bf 1}_{\{i_3=i_5\}}
{\bf 1}_{\{j_3=j_5\}}
\zeta_{j_1}^{(i_1)}
\zeta_{j_2}^{(i_2)}
\zeta_{j_4}^{(i_4)}
-{\bf 1}_{\{i_4=i_5\}}
{\bf 1}_{\{j_4=j_5\}}
\zeta_{j_1}^{(i_1)}
\zeta_{j_2}^{(i_2)}
\zeta_{j_3}^{(i_3)}+
$$
$$
+
{\bf 1}_{\{i_1=i_2\}}
{\bf 1}_{\{j_1=j_2\}}
{\bf 1}_{\{i_3=i_4\}}
{\bf 1}_{\{j_3=j_4\}}\zeta_{j_5}^{(i_5)}+
{\bf 1}_{\{i_1=i_2\}}
{\bf 1}_{\{j_1=j_2\}}
{\bf 1}_{\{i_3=i_5\}}
{\bf 1}_{\{j_3=j_5\}}\zeta_{j_4}^{(i_4)}+
$$
$$
+
{\bf 1}_{\{i_1=i_2\}}
{\bf 1}_{\{j_1=j_2\}}
{\bf 1}_{\{i_4=i_5\}}
{\bf 1}_{\{j_4=j_5\}}\zeta_{j_3}^{(i_3)}+
{\bf 1}_{\{i_1=i_3\}}
{\bf 1}_{\{j_1=j_3\}}
{\bf 1}_{\{i_2=i_4\}}
{\bf 1}_{\{j_2=j_4\}}\zeta_{j_5}^{(i_5)}+
$$
$$
+
{\bf 1}_{\{i_1=i_3\}}
{\bf 1}_{\{j_1=j_3\}}
{\bf 1}_{\{i_2=i_5\}}
{\bf 1}_{\{j_2=j_5\}}\zeta_{j_4}^{(i_4)}+
{\bf 1}_{\{i_1=i_3\}}
{\bf 1}_{\{j_1=j_3\}}
{\bf 1}_{\{i_4=i_5\}}
{\bf 1}_{\{j_4=j_5\}}\zeta_{j_2}^{(i_2)}+
$$
$$
+
{\bf 1}_{\{i_1=i_4\}}
{\bf 1}_{\{j_1=j_4\}}
{\bf 1}_{\{i_2=i_3\}}
{\bf 1}_{\{j_2=j_3\}}\zeta_{j_5}^{(i_5)}+
{\bf 1}_{\{i_1=i_4\}}
{\bf 1}_{\{j_1=j_4\}}
{\bf 1}_{\{i_2=i_5\}}
{\bf 1}_{\{j_2=j_5\}}\zeta_{j_3}^{(i_3)}+
$$
$$
+
{\bf 1}_{\{i_1=i_4\}}
{\bf 1}_{\{j_1=j_4\}}
{\bf 1}_{\{i_3=i_5\}}
{\bf 1}_{\{j_3=j_5\}}\zeta_{j_2}^{(i_2)}+
{\bf 1}_{\{i_1=i_5\}}
{\bf 1}_{\{j_1=j_5\}}
{\bf 1}_{\{i_2=i_3\}}
{\bf 1}_{\{j_2=j_3\}}\zeta_{j_4}^{(i_4)}+
$$
$$
+
{\bf 1}_{\{i_1=i_5\}}
{\bf 1}_{\{j_1=j_5\}}
{\bf 1}_{\{i_2=i_4\}}
{\bf 1}_{\{j_2=j_4\}}\zeta_{j_3}^{(i_3)}+
{\bf 1}_{\{i_1=i_5\}}
{\bf 1}_{\{j_1=j_5\}}
{\bf 1}_{\{i_3=i_4\}}
{\bf 1}_{\{j_3=j_4\}}\zeta_{j_2}^{(i_2)}+
$$
$$
+
{\bf 1}_{\{i_2=i_3\}}
{\bf 1}_{\{j_2=j_3\}}
{\bf 1}_{\{i_4=i_5\}}
{\bf 1}_{\{j_4=j_5\}}\zeta_{j_1}^{(i_1)}+
{\bf 1}_{\{i_2=i_4\}}
{\bf 1}_{\{j_2=j_4\}}
{\bf 1}_{\{i_3=i_5\}}
{\bf 1}_{\{j_3=j_5\}}\zeta_{j_1}^{(i_1)}+
$$
\begin{equation}
\label{desy6}
+\Biggl.
{\bf 1}_{\{i_2=i_5\}}
{\bf 1}_{\{j_2=j_5\}}
{\bf 1}_{\{i_3=i_4\}}
{\bf 1}_{\{j_3=j_4\}}\zeta_{j_1}^{(i_1)}\Biggr),
\end{equation}

\vspace{4mm}

\begin{equation}
\label{ogo3}
I_{(20)T,t}^{(i_1 i_2)q_5}=
\sum_{j_1,j_2=0}^{q_5}
C_{j_2j_1}^{20}\Biggl(\zeta_{j_1}^{(i_1)}\zeta_{j_2}^{(i_2)}
-{\bf 1}_{\{i_1=i_2\}}
{\bf 1}_{\{j_1=j_2\}}\Biggr),
\end{equation}

\vspace{2mm}
\begin{equation}
\label{ogo4}
I_{(11)T,t}^{(i_1 i_2)q_6}=
\sum_{j_1,j_2=0}^{q_6}
C_{j_2j_1}^{11}\Biggl(\zeta_{j_1}^{(i_1)}\zeta_{j_2}^{(i_2)}
-{\bf 1}_{\{i_1=i_2\}}
{\bf 1}_{\{j_1=j_2\}}\Biggr),
\end{equation}

\vspace{2mm}
\begin{equation}
\label{ogo5}
I_{(02)T,t}^{(i_1 i_2)q_{7}}=
\sum_{j_1,j_2=0}^{q_{7}}
C_{j_2j_1}^{02}\Biggl(\zeta_{j_1}^{(i_1)}\zeta_{j_2}^{(i_2)}
-{\bf 1}_{\{i_1=i_2\}}
{\bf 1}_{\{j_1=j_2\}}\Biggr),
\end{equation}

\vspace{5mm}
$$
I_{(001)T,t}^{(i_1i_2i_3)q_8}
=
\sum_{j_1,j_2,j_3=0}^{q_8}
C_{j_3j_2j_1}^{001}\Biggl(
\zeta_{j_1}^{(i_1)}\zeta_{j_2}^{(i_2)}\zeta_{j_3}^{(i_3)}
-{\bf 1}_{\{i_1=i_2\}}
{\bf 1}_{\{j_1=j_2\}}
\zeta_{j_3}^{(i_3)}-
\Biggr.
$$
\begin{equation}
\label{desy7}
\Biggl.
-{\bf 1}_{\{i_2=i_3\}}
{\bf 1}_{\{j_2=j_3\}}
\zeta_{j_1}^{(i_1)}-
{\bf 1}_{\{i_1=i_3\}}
{\bf 1}_{\{j_1=j_3\}}
\zeta_{j_2}^{(i_2)}\Biggr),
\end{equation}

\newpage
\noindent
$$
I_{(010)T,t}^{(i_1i_2i_3)q_9}
=
\sum_{j_1,j_2,j_3=0}^{q_9}
C_{j_3j_2j_1}^{010}\Biggl(
\zeta_{j_1}^{(i_1)}\zeta_{j_2}^{(i_2)}\zeta_{j_3}^{(i_3)}
-{\bf 1}_{\{i_1=i_2\}}
{\bf 1}_{\{j_1=j_2\}}
\zeta_{j_3}^{(i_3)}-
\Biggr.
$$
\begin{equation}
\label{desy8}
\Biggl.
-{\bf 1}_{\{i_2=i_3\}}
{\bf 1}_{\{j_2=j_3\}}
\zeta_{j_1}^{(i_1)}-
{\bf 1}_{\{i_1=i_3\}}
{\bf 1}_{\{j_1=j_3\}}
\zeta_{j_2}^{(i_2)}\Biggr),
\end{equation}

\vspace{4mm}
$$
I_{(100)T,t}^{(i_1i_2i_3)q_{10}}
=
\sum_{j_1,j_2,j_3=0}^{q_{10}}
C_{j_3j_2j_1}^{100}\Biggl(
\zeta_{j_1}^{(i_1)}\zeta_{j_2}^{(i_2)}\zeta_{j_3}^{(i_3)}
-{\bf 1}_{\{i_1=i_2\}}
{\bf 1}_{\{j_1=j_2\}}
\zeta_{j_3}^{(i_3)}-
\Biggr.
$$
\begin{equation}
\label{desy9}
\Biggl.
-{\bf 1}_{\{i_2=i_3\}}
{\bf 1}_{\{j_2=j_3\}}
\zeta_{j_1}^{(i_1)}-
{\bf 1}_{\{i_1=i_3\}}
{\bf 1}_{\{j_1=j_3\}}
\zeta_{j_2}^{(i_2)}\Biggr),
\end{equation}

\vspace{4mm}
$$
I_{(0001)T,t}^{(i_1 i_2 i_3 i_4)q_{11}}
=
\sum_{j_1,j_2,j_3,j_4=0}^{q_{11}}
C_{j_4 j_3 j_2 j_1}^{0001}\Biggl(
\zeta_{j_1}^{(i_1)}\zeta_{j_2}^{(i_2)}\zeta_{j_3}^{(i_3)}\zeta_{j_4}^{(i_4)}
-\Biggr.
$$
$$
-
{\bf 1}_{\{i_1=i_2\}}
{\bf 1}_{\{j_1=j_2\}}
\zeta_{j_3}^{(i_3)}
\zeta_{j_4}^{(i_4)}
-
{\bf 1}_{\{i_1=i_3\}}
{\bf 1}_{\{j_1=j_3\}}
\zeta_{j_2}^{(i_2)}
\zeta_{j_4}^{(i_4)}-
$$
$$
-
{\bf 1}_{\{i_1=i_4\}}
{\bf 1}_{\{j_1=j_4\}}
\zeta_{j_2}^{(i_2)}
\zeta_{j_3}^{(i_3)}
-
{\bf 1}_{\{i_2=i_3\}}
{\bf 1}_{\{j_2=j_3\}}
\zeta_{j_1}^{(i_1)}
\zeta_{j_4}^{(i_4)}-
$$
$$
-
{\bf 1}_{\{i_2=i_4\}}
{\bf 1}_{\{j_2=j_4\}}
\zeta_{j_1}^{(i_1)}
\zeta_{j_3}^{(i_3)}
-
{\bf 1}_{\{i_3=i_4\}}
{\bf 1}_{\{j_3=j_4\}}
\zeta_{j_1}^{(i_1)}
\zeta_{j_2}^{(i_2)}+
$$
$$
+
{\bf 1}_{\{i_1=i_2\}}
{\bf 1}_{\{j_1=j_2\}}
{\bf 1}_{\{i_3=i_4\}}
{\bf 1}_{\{j_3=j_4\}}+
{\bf 1}_{\{i_1=i_3\}}
{\bf 1}_{\{j_1=j_3\}}
{\bf 1}_{\{i_2=i_4\}}
{\bf 1}_{\{j_2=j_4\}}+
$$
\begin{equation}
\label{desy10}
+\Biggl.
{\bf 1}_{\{i_1=i_4\}}
{\bf 1}_{\{j_1=j_4\}}
{\bf 1}_{\{i_2=i_3\}}
{\bf 1}_{\{j_2=j_3\}}\Biggr),
\end{equation}

\vspace{4mm}

$$
I_{(0010)T,t}^{(i_1 i_2 i_3 i_4)q_{12}}
\sum_{j_1,j_2,j_3,j_4=0}^{q_{12}}
C_{j_4 j_3 j_2 j_1}^{0010}\Biggl(
\zeta_{j_1}^{(i_1)}\zeta_{j_2}^{(i_2)}\zeta_{j_3}^{(i_3)}\zeta_{j_4}^{(i_4)}
-\Biggr.
$$
$$
-
{\bf 1}_{\{i_1=i_2\}}
{\bf 1}_{\{j_1=j_2\}}
\zeta_{j_3}^{(i_3)}
\zeta_{j_4}^{(i_4)}
-
{\bf 1}_{\{i_1=i_3\}}
{\bf 1}_{\{j_1=j_3\}}
\zeta_{j_2}^{(i_2)}
\zeta_{j_4}^{(i_4)}-
$$
$$
-
{\bf 1}_{\{i_1=i_4\}}
{\bf 1}_{\{j_1=j_4\}}
\zeta_{j_2}^{(i_2)}
\zeta_{j_3}^{(i_3)}
-
{\bf 1}_{\{i_2=i_3\}}
{\bf 1}_{\{j_2=j_3\}}
\zeta_{j_1}^{(i_1)}
\zeta_{j_4}^{(i_4)}-
$$
$$
-
{\bf 1}_{\{i_2=i_4\}}
{\bf 1}_{\{j_2=j_4\}}
\zeta_{j_1}^{(i_1)}
\zeta_{j_3}^{(i_3)}
-
{\bf 1}_{\{i_3=i_4\}}
{\bf 1}_{\{j_3=j_4\}}
\zeta_{j_1}^{(i_1)}
\zeta_{j_2}^{(i_2)}+
$$
$$
+
{\bf 1}_{\{i_1=i_2\}}
{\bf 1}_{\{j_1=j_2\}}
{\bf 1}_{\{i_3=i_4\}}
{\bf 1}_{\{j_3=j_4\}}+
{\bf 1}_{\{i_1=i_3\}}
{\bf 1}_{\{j_1=j_3\}}
{\bf 1}_{\{i_2=i_4\}}
{\bf 1}_{\{j_2=j_4\}}+
$$
\begin{equation}
\label{desy11}
+\Biggl.
{\bf 1}_{\{i_1=i_4\}}
{\bf 1}_{\{j_1=j_4\}}
{\bf 1}_{\{i_2=i_3\}}
{\bf 1}_{\{j_2=j_3\}}\Biggr),
\end{equation}

\newpage
\noindent
$$
I_{(0100)T,t}^{(i_1 i_2 i_3 i_4)q_{13}}
=
\sum_{j_1,j_2,j_3,j_4=0}^{q_{13}}
C_{j_4 j_3 j_2 j_1}^{0100}\Biggl(
\zeta_{j_1}^{(i_1)}\zeta_{j_2}^{(i_2)}\zeta_{j_3}^{(i_3)}\zeta_{j_4}^{(i_4)}
-\Biggr.
$$
$$
-
{\bf 1}_{\{i_1=i_2\}}
{\bf 1}_{\{j_1=j_2\}}
\zeta_{j_3}^{(i_3)}
\zeta_{j_4}^{(i_4)}
-
{\bf 1}_{\{i_1=i_3\}}
{\bf 1}_{\{j_1=j_3\}}
\zeta_{j_2}^{(i_2)}
\zeta_{j_4}^{(i_4)}-
$$
$$
-
{\bf 1}_{\{i_1=i_4\}}
{\bf 1}_{\{j_1=j_4\}}
\zeta_{j_2}^{(i_2)}
\zeta_{j_3}^{(i_3)}
-
{\bf 1}_{\{i_2=i_3\}}
{\bf 1}_{\{j_2=j_3\}}
\zeta_{j_1}^{(i_1)}
\zeta_{j_4}^{(i_4)}-
$$
$$
-
{\bf 1}_{\{i_2=i_4\}}
{\bf 1}_{\{j_2=j_4\}}
\zeta_{j_1}^{(i_1)}
\zeta_{j_3}^{(i_3)}
-
{\bf 1}_{\{i_3=i_4\}}
{\bf 1}_{\{j_3=j_4\}}
\zeta_{j_1}^{(i_1)}
\zeta_{j_2}^{(i_2)}+
$$
$$
+
{\bf 1}_{\{i_1=i_2\}}
{\bf 1}_{\{j_1=j_2\}}
{\bf 1}_{\{i_3=i_4\}}
{\bf 1}_{\{j_3=j_4\}}+
{\bf 1}_{\{i_1=i_3\}}
{\bf 1}_{\{j_1=j_3\}}
{\bf 1}_{\{i_2=i_4\}}
{\bf 1}_{\{j_2=j_4\}}+
$$
\begin{equation}
\label{desy12}
+\Biggl.
{\bf 1}_{\{i_1=i_4\}}
{\bf 1}_{\{j_1=j_4\}}
{\bf 1}_{\{i_2=i_3\}}
{\bf 1}_{\{j_2=j_3\}}\Biggr),
\end{equation}

\vspace{3mm}

$$
I_{(1000)T,t}^{(i_1 i_2 i_3 i_4)q_{14}}
=
\sum_{j_1,j_2,j_3,j_4=0}^{q_{14}}
C_{j_4 j_3 j_2 j_1}^{1000}\Biggl(
\zeta_{j_1}^{(i_1)}\zeta_{j_2}^{(i_2)}\zeta_{j_3}^{(i_3)}\zeta_{j_4}^{(i_4)}
-\Biggr.
$$
$$
-
{\bf 1}_{\{i_1=i_2\}}
{\bf 1}_{\{j_1=j_2\}}
\zeta_{j_3}^{(i_3)}
\zeta_{j_4}^{(i_4)}
-
{\bf 1}_{\{i_1=i_3\}}
{\bf 1}_{\{j_1=j_3\}}
\zeta_{j_2}^{(i_2)}
\zeta_{j_4}^{(i_4)}-
$$
$$
-
{\bf 1}_{\{i_1=i_4\}}
{\bf 1}_{\{j_1=j_4\}}
\zeta_{j_2}^{(i_2)}
\zeta_{j_3}^{(i_3)}
-
{\bf 1}_{\{i_2=i_3\}}
{\bf 1}_{\{j_2=j_3\}}
\zeta_{j_1}^{(i_1)}
\zeta_{j_4}^{(i_4)}-
$$
$$
-
{\bf 1}_{\{i_2=i_4\}}
{\bf 1}_{\{j_2=j_4\}}
\zeta_{j_1}^{(i_1)}
\zeta_{j_3}^{(i_3)}
-
{\bf 1}_{\{i_3=i_4\}}
{\bf 1}_{\{j_3=j_4\}}
\zeta_{j_1}^{(i_1)}
\zeta_{j_2}^{(i_2)}+
$$
$$
+
{\bf 1}_{\{i_1=i_2\}}
{\bf 1}_{\{j_1=j_2\}}
{\bf 1}_{\{i_3=i_4\}}
{\bf 1}_{\{j_3=j_4\}}+
{\bf 1}_{\{i_1=i_3\}}
{\bf 1}_{\{j_1=j_3\}}
{\bf 1}_{\{i_2=i_4\}}
{\bf 1}_{\{j_2=j_4\}}+
$$
\begin{equation}
\label{desy13}
+\Biggl.
{\bf 1}_{\{i_1=i_4\}}
{\bf 1}_{\{j_1=j_4\}}
{\bf 1}_{\{i_2=i_3\}}
{\bf 1}_{\{j_2=j_3\}}\Biggr),
\end{equation}

\vspace{6mm}

$$
I_{(000000)T,t}^{(i_1 i_2 i_3 i_4 i_5 i_6)q_{15}}
=\sum_{j_1,j_2,j_3,j_4,j_5,j_6=0}^{q_{15}}
C_{j_6 j_5 j_4 j_3 j_2 j_1}^{000000}\Biggl(
\prod_{l=1}^6
\zeta_{j_l}^{(i_l)}
-\Biggr.
$$
$$
-
{\bf 1}_{\{j_1=j_6\}}
{\bf 1}_{\{i_1=i_6\}}
\zeta_{j_2}^{(i_2)}
\zeta_{j_3}^{(i_3)}
\zeta_{j_4}^{(i_4)}
\zeta_{j_5}^{(i_5)}-
{\bf 1}_{\{j_2=j_6\}}
{\bf 1}_{\{i_2=i_6\}}
\zeta_{j_1}^{(i_1)}
\zeta_{j_3}^{(i_3)}
\zeta_{j_4}^{(i_4)}
\zeta_{j_5}^{(i_5)}-
$$
$$
-
{\bf 1}_{\{j_3=j_6\}}
{\bf 1}_{\{i_3=i_6\}}
\zeta_{j_1}^{(i_1)}
\zeta_{j_2}^{(i_2)}
\zeta_{j_4}^{(i_4)}
\zeta_{j_5}^{(i_5)}-
{\bf 1}_{\{j_4=j_6\}}
{\bf 1}_{\{i_4=i_6\}}
\zeta_{j_1}^{(i_1)}
\zeta_{j_2}^{(i_2)}
\zeta_{j_3}^{(i_3)}
\zeta_{j_5}^{(i_5)}-
$$
$$
-
{\bf 1}_{\{j_5=j_6\}}
{\bf 1}_{\{i_5=i_6\}}
\zeta_{j_1}^{(i_1)}
\zeta_{j_2}^{(i_2)}
\zeta_{j_3}^{(i_3)}
\zeta_{j_4}^{(i_4)}-
{\bf 1}_{\{j_1=j_2\}}
{\bf 1}_{\{i_1=i_2\}}
\zeta_{j_3}^{(i_3)}
\zeta_{j_4}^{(i_4)}
\zeta_{j_5}^{(i_5)}
\zeta_{j_6}^{(i_6)}-
$$
$$
-
{\bf 1}_{\{j_1=j_3\}}
{\bf 1}_{\{i_1=i_3\}}
\zeta_{j_2}^{(i_2)}
\zeta_{j_4}^{(i_4)}
\zeta_{j_5}^{(i_5)}
\zeta_{j_6}^{(i_6)}-
{\bf 1}_{\{j_1=j_4\}}
{\bf 1}_{\{i_1=i_4\}}
\zeta_{j_2}^{(i_2)}
\zeta_{j_3}^{(i_3)}
\zeta_{j_5}^{(i_5)}
\zeta_{j_6}^{(i_6)}-
$$
$$
-
{\bf 1}_{\{j_1=j_5\}}
{\bf 1}_{\{i_1=i_5\}}
\zeta_{j_2}^{(i_2)}
\zeta_{j_3}^{(i_3)}
\zeta_{j_4}^{(i_4)}
\zeta_{j_6}^{(i_6)}-
{\bf 1}_{\{j_2=j_3\}}
{\bf 1}_{\{i_2=i_3\}}
\zeta_{j_1}^{(i_1)}
\zeta_{j_4}^{(i_4)}
\zeta_{j_5}^{(i_5)}
\zeta_{j_6}^{(i_6)}-
$$
$$
-
{\bf 1}_{\{j_2=j_4\}}
{\bf 1}_{\{i_2=i_4\}}
\zeta_{j_1}^{(i_1)}
\zeta_{j_3}^{(i_3)}
\zeta_{j_5}^{(i_5)}
\zeta_{j_6}^{(i_6)}-
{\bf 1}_{\{j_2=j_5\}}
{\bf 1}_{\{i_2=i_5\}}
\zeta_{j_1}^{(i_1)}
\zeta_{j_3}^{(i_3)}
\zeta_{j_4}^{(i_4)}
\zeta_{j_6}^{(i_6)}-
$$
$$
-
{\bf 1}_{\{j_3=j_4\}}
{\bf 1}_{\{i_3=i_4\}}
\zeta_{j_1}^{(i_1)}
\zeta_{j_2}^{(i_2)}
\zeta_{j_5}^{(i_5)}
\zeta_{j_6}^{(i_6)}-
{\bf 1}_{\{j_3=j_5\}}
{\bf 1}_{\{i_3=i_5\}}
\zeta_{j_1}^{(i_1)}
\zeta_{j_2}^{(i_2)}
\zeta_{j_4}^{(i_4)}
\zeta_{j_6}^{(i_6)}-
$$
$$
-
{\bf 1}_{\{j_4=j_5\}}
{\bf 1}_{\{i_4=i_5\}}
\zeta_{j_1}^{(i_1)}
\zeta_{j_2}^{(i_2)}
\zeta_{j_3}^{(i_3)}
\zeta_{j_6}^{(i_6)}+
$$
$$
+
{\bf 1}_{\{j_1=j_2\}}
{\bf 1}_{\{i_1=i_2\}}
{\bf 1}_{\{j_3=j_4\}}
{\bf 1}_{\{i_3=i_4\}}
\zeta_{j_5}^{(i_5)}
\zeta_{j_6}^{(i_6)}
+
{\bf 1}_{\{j_1=j_2\}}
{\bf 1}_{\{i_1=i_2\}}
{\bf 1}_{\{j_3=j_5\}}
{\bf 1}_{\{i_3=i_5\}}
\zeta_{j_4}^{(i_4)}
\zeta_{j_6}^{(i_6)}+
$$
$$
+
{\bf 1}_{\{j_1=j_2\}}
{\bf 1}_{\{i_1=i_2\}}
{\bf 1}_{\{j_4=j_5\}}
{\bf 1}_{\{i_4=i_5\}}
\zeta_{j_3}^{(i_3)}
\zeta_{j_6}^{(i_6)}
+
{\bf 1}_{\{j_1=j_3\}}
{\bf 1}_{\{i_1=i_3\}}
{\bf 1}_{\{j_2=j_4\}}
{\bf 1}_{\{i_2=i_4\}}
\zeta_{j_5}^{(i_5)}
\zeta_{j_6}^{(i_6)}+
$$
$$
+
{\bf 1}_{\{j_1=j_3\}}
{\bf 1}_{\{i_1=i_3\}}
{\bf 1}_{\{j_2=j_5\}}
{\bf 1}_{\{i_2=i_5\}}
\zeta_{j_4}^{(i_4)}
\zeta_{j_6}^{(i_6)}
+
{\bf 1}_{\{j_1=j_3\}}
{\bf 1}_{\{i_1=i_3\}}
{\bf 1}_{\{j_4=j_5\}}
{\bf 1}_{\{i_4=i_5\}}
\zeta_{j_2}^{(i_2)}
\zeta_{j_6}^{(i_6)}+
$$
$$
+
{\bf 1}_{\{j_1=j_4\}}
{\bf 1}_{\{i_1=i_4\}}
{\bf 1}_{\{j_2=j_3\}}
{\bf 1}_{\{i_2=i_3\}}
\zeta_{j_5}^{(i_5)}
\zeta_{j_6}^{(i_6)}
+
{\bf 1}_{\{j_1=j_4\}}
{\bf 1}_{\{i_1=i_4\}}
{\bf 1}_{\{j_2=j_5\}}
{\bf 1}_{\{i_2=i_5\}}
\zeta_{j_3}^{(i_3)}
\zeta_{j_6}^{(i_6)}+
$$
$$
+
{\bf 1}_{\{j_1=j_4\}}
{\bf 1}_{\{i_1=i_4\}}
{\bf 1}_{\{j_3=j_5\}}
{\bf 1}_{\{i_3=i_5\}}
\zeta_{j_2}^{(i_2)}
\zeta_{j_6}^{(i_6)}
+
{\bf 1}_{\{j_1=j_5\}}
{\bf 1}_{\{i_1=i_5\}}
{\bf 1}_{\{j_2=j_3\}}
{\bf 1}_{\{i_2=i_3\}}
\zeta_{j_4}^{(i_4)}
\zeta_{j_6}^{(i_6)}+
$$
$$
+
{\bf 1}_{\{j_1=j_5\}}
{\bf 1}_{\{i_1=i_5\}}
{\bf 1}_{\{j_2=j_4\}}
{\bf 1}_{\{i_2=i_4\}}
\zeta_{j_3}^{(i_3)}
\zeta_{j_6}^{(i_6)}
+
{\bf 1}_{\{j_1=j_5\}}
{\bf 1}_{\{i_1=i_5\}}
{\bf 1}_{\{j_3=j_4\}}
{\bf 1}_{\{i_3=i_4\}}
\zeta_{j_2}^{(i_2)}
\zeta_{j_6}^{(i_6)}+
$$
$$
+
{\bf 1}_{\{j_2=j_3\}}
{\bf 1}_{\{i_2=i_3\}}
{\bf 1}_{\{j_4=j_5\}}
{\bf 1}_{\{i_4=i_5\}}
\zeta_{j_1}^{(i_1)}
\zeta_{j_6}^{(i_6)}
+
{\bf 1}_{\{j_2=j_4\}}
{\bf 1}_{\{i_2=i_4\}}
{\bf 1}_{\{j_3=j_5\}}
{\bf 1}_{\{i_3=i_5\}}
\zeta_{j_1}^{(i_1)}
\zeta_{j_6}^{(i_6)}+
$$
$$
+
{\bf 1}_{\{j_2=j_5\}}
{\bf 1}_{\{i_2=i_5\}}
{\bf 1}_{\{j_3=j_4\}}
{\bf 1}_{\{i_3=i_4\}}
\zeta_{j_1}^{(i_1)}
\zeta_{j_6}^{(i_6)}
+
{\bf 1}_{\{j_6=j_1\}}
{\bf 1}_{\{i_6=i_1\}}
{\bf 1}_{\{j_3=j_4\}}
{\bf 1}_{\{i_3=i_4\}}
\zeta_{j_2}^{(i_2)}
\zeta_{j_5}^{(i_5)}+
$$
$$
+
{\bf 1}_{\{j_6=j_1\}}
{\bf 1}_{\{i_6=i_1\}}
{\bf 1}_{\{j_3=j_5\}}
{\bf 1}_{\{i_3=i_5\}}
\zeta_{j_2}^{(i_2)}
\zeta_{j_4}^{(i_4)}
+
{\bf 1}_{\{j_6=j_1\}}
{\bf 1}_{\{i_6=i_1\}}
{\bf 1}_{\{j_2=j_5\}}
{\bf 1}_{\{i_2=i_5\}}
\zeta_{j_3}^{(i_3)}
\zeta_{j_4}^{(i_4)}+
$$
$$
+
{\bf 1}_{\{j_6=j_1\}}
{\bf 1}_{\{i_6=i_1\}}
{\bf 1}_{\{j_2=j_4\}}
{\bf 1}_{\{i_2=i_4\}}
\zeta_{j_3}^{(i_3)}
\zeta_{j_5}^{(i_5)}
+
{\bf 1}_{\{j_6=j_1\}}
{\bf 1}_{\{i_6=i_1\}}
{\bf 1}_{\{j_4=j_5\}}
{\bf 1}_{\{i_4=i_5\}}
\zeta_{j_2}^{(i_2)}
\zeta_{j_3}^{(i_3)}+
$$
$$
+
{\bf 1}_{\{j_6=j_1\}}
{\bf 1}_{\{i_6=i_1\}}
{\bf 1}_{\{j_2=j_3\}}
{\bf 1}_{\{i_2=i_3\}}
\zeta_{j_4}^{(i_4)}
\zeta_{j_5}^{(i_5)}
+
{\bf 1}_{\{j_6=j_2\}}
{\bf 1}_{\{i_6=i_2\}}
{\bf 1}_{\{j_3=j_5\}}
{\bf 1}_{\{i_3=i_5\}}
\zeta_{j_1}^{(i_1)}
\zeta_{j_4}^{(i_4)}+
$$
$$
+
{\bf 1}_{\{j_6=j_2\}}
{\bf 1}_{\{i_6=i_2\}}
{\bf 1}_{\{j_4=j_5\}}
{\bf 1}_{\{i_4=i_5\}}
\zeta_{j_1}^{(i_1)}
\zeta_{j_3}^{(i_3)}
+
{\bf 1}_{\{j_6=j_2\}}
{\bf 1}_{\{i_6=i_2\}}
{\bf 1}_{\{j_3=j_4\}}
{\bf 1}_{\{i_3=i_4\}}
\zeta_{j_1}^{(i_1)}
\zeta_{j_5}^{(i_5)}+
$$
$$
+
{\bf 1}_{\{j_6=j_2\}}
{\bf 1}_{\{i_6=i_2\}}
{\bf 1}_{\{j_1=j_5\}}
{\bf 1}_{\{i_1=i_5\}}
\zeta_{j_3}^{(i_3)}
\zeta_{j_4}^{(i_4)}
+
{\bf 1}_{\{j_6=j_2\}}
{\bf 1}_{\{i_6=i_2\}}
{\bf 1}_{\{j_1=j_4\}}
{\bf 1}_{\{i_1=i_4\}}
\zeta_{j_3}^{(i_3)}
\zeta_{j_5}^{(i_5)}+
$$
$$
+
{\bf 1}_{\{j_6=j_2\}}
{\bf 1}_{\{i_6=i_2\}}
{\bf 1}_{\{j_1=j_3\}}
{\bf 1}_{\{i_1=i_3\}}
\zeta_{j_4}^{(i_4)}
\zeta_{j_5}^{(i_5)}
+
{\bf 1}_{\{j_6=j_3\}}
{\bf 1}_{\{i_6=i_3\}}
{\bf 1}_{\{j_2=j_5\}}
{\bf 1}_{\{i_2=i_5\}}
\zeta_{j_1}^{(i_1)}
\zeta_{j_4}^{(i_4)}+
$$
$$
+
{\bf 1}_{\{j_6=j_3\}}
{\bf 1}_{\{i_6=i_3\}}
{\bf 1}_{\{j_4=j_5\}}
{\bf 1}_{\{i_4=i_5\}}
\zeta_{j_1}^{(i_1)}
\zeta_{j_2}^{(i_2)}
+
{\bf 1}_{\{j_6=j_3\}}
{\bf 1}_{\{i_6=i_3\}}
{\bf 1}_{\{j_2=j_4\}}
{\bf 1}_{\{i_2=i_4\}}
\zeta_{j_1}^{(i_1)}
\zeta_{j_5}^{(i_5)}+
$$
$$
+
{\bf 1}_{\{j_6=j_3\}}
{\bf 1}_{\{i_6=i_3\}}
{\bf 1}_{\{j_1=j_5\}}
{\bf 1}_{\{i_1=i_5\}}
\zeta_{j_2}^{(i_2)}
\zeta_{j_4}^{(i_4)}
+
{\bf 1}_{\{j_6=j_3\}}
{\bf 1}_{\{i_6=i_3\}}
{\bf 1}_{\{j_1=j_4\}}
{\bf 1}_{\{i_1=i_4\}}
\zeta_{j_2}^{(i_2)}
\zeta_{j_5}^{(i_5)}+
$$
$$
+
{\bf 1}_{\{j_6=j_3\}}
{\bf 1}_{\{i_6=i_3\}}
{\bf 1}_{\{j_1=j_2\}}
{\bf 1}_{\{i_1=i_2\}}
\zeta_{j_4}^{(i_4)}
\zeta_{j_5}^{(i_5)}
+
{\bf 1}_{\{j_6=j_4\}}
{\bf 1}_{\{i_6=i_4\}}
{\bf 1}_{\{j_3=j_5\}}
{\bf 1}_{\{i_3=i_5\}}
\zeta_{j_1}^{(i_1)}
\zeta_{j_2}^{(i_2)}+
$$
$$
+
{\bf 1}_{\{j_6=j_4\}}
{\bf 1}_{\{i_6=i_4\}}
{\bf 1}_{\{j_2=j_5\}}
{\bf 1}_{\{i_2=i_5\}}
\zeta_{j_1}^{(i_1)}
\zeta_{j_3}^{(i_3)}
+
{\bf 1}_{\{j_6=j_4\}}
{\bf 1}_{\{i_6=i_4\}}
{\bf 1}_{\{j_2=j_3\}}
{\bf 1}_{\{i_2=i_3\}}
\zeta_{j_1}^{(i_1)}
\zeta_{j_5}^{(i_5)}+
$$
$$
+
{\bf 1}_{\{j_6=j_4\}}
{\bf 1}_{\{i_6=i_4\}}
{\bf 1}_{\{j_1=j_5\}}
{\bf 1}_{\{i_1=i_5\}}
\zeta_{j_2}^{(i_2)}
\zeta_{j_3}^{(i_3)}
+
{\bf 1}_{\{j_6=j_4\}}
{\bf 1}_{\{i_6=i_4\}}
{\bf 1}_{\{j_1=j_3\}}
{\bf 1}_{\{i_1=i_3\}}
\zeta_{j_2}^{(i_2)}
\zeta_{j_5}^{(i_5)}+
$$
$$
+
{\bf 1}_{\{j_6=j_4\}}
{\bf 1}_{\{i_6=i_4\}}
{\bf 1}_{\{j_1=j_2\}}
{\bf 1}_{\{i_1=i_2\}}
\zeta_{j_3}^{(i_3)}
\zeta_{j_5}^{(i_5)}
+
{\bf 1}_{\{j_6=j_5\}}
{\bf 1}_{\{i_6=i_5\}}
{\bf 1}_{\{j_3=j_4\}}
{\bf 1}_{\{i_3=i_4\}}
\zeta_{j_1}^{(i_1)}
\zeta_{j_2}^{(i_2)}+
$$
$$
+
{\bf 1}_{\{j_6=j_5\}}
{\bf 1}_{\{i_6=i_5\}}
{\bf 1}_{\{j_2=j_4\}}
{\bf 1}_{\{i_2=i_4\}}
\zeta_{j_1}^{(i_1)}
\zeta_{j_3}^{(i_3)}
+
{\bf 1}_{\{j_6=j_5\}}
{\bf 1}_{\{i_6=i_5\}}
{\bf 1}_{\{j_2=j_3\}}
{\bf 1}_{\{i_2=i_3\}}
\zeta_{j_1}^{(i_1)}
\zeta_{j_4}^{(i_4)}+
$$
$$
+
{\bf 1}_{\{j_6=j_5\}}
{\bf 1}_{\{i_6=i_5\}}
{\bf 1}_{\{j_1=j_4\}}
{\bf 1}_{\{i_1=i_4\}}
\zeta_{j_2}^{(i_2)}
\zeta_{j_3}^{(i_3)}
+
{\bf 1}_{\{j_6=j_5\}}
{\bf 1}_{\{i_6=i_5\}}
{\bf 1}_{\{j_1=j_3\}}
{\bf 1}_{\{i_1=i_3\}}
\zeta_{j_2}^{(i_2)}
\zeta_{j_4}^{(i_4)}+
$$
$$
+
{\bf 1}_{\{j_6=j_5\}}
{\bf 1}_{\{i_6=i_5\}}
{\bf 1}_{\{j_1=j_2\}}
{\bf 1}_{\{i_1=i_2\}}
\zeta_{j_3}^{(i_3)}
\zeta_{j_4}^{(i_4)}-
$$
$$
-
{\bf 1}_{\{j_6=j_1\}}
{\bf 1}_{\{i_6=i_1\}}
{\bf 1}_{\{j_2=j_5\}}
{\bf 1}_{\{i_2=i_5\}}
{\bf 1}_{\{j_3=j_4\}}
{\bf 1}_{\{i_3=i_4\}}-
$$
$$
-
{\bf 1}_{\{j_6=j_1\}}
{\bf 1}_{\{i_6=i_1\}}
{\bf 1}_{\{j_2=j_4\}}
{\bf 1}_{\{i_2=i_4\}}
{\bf 1}_{\{j_3=j_5\}}
{\bf 1}_{\{i_3=i_5\}}-
$$
$$
-
{\bf 1}_{\{j_6=j_1\}}
{\bf 1}_{\{i_6=i_1\}}
{\bf 1}_{\{j_2=j_3\}}
{\bf 1}_{\{i_2=i_3\}}
{\bf 1}_{\{j_4=j_5\}}
{\bf 1}_{\{i_4=i_5\}}-
$$
$$
-               
{\bf 1}_{\{j_6=j_2\}}
{\bf 1}_{\{i_6=i_2\}}
{\bf 1}_{\{j_1=j_5\}}
{\bf 1}_{\{i_1=i_5\}}
{\bf 1}_{\{j_3=j_4\}}
{\bf 1}_{\{i_3=i_4\}}-
$$
$$
-
{\bf 1}_{\{j_6=j_2\}}
{\bf 1}_{\{i_6=i_2\}}
{\bf 1}_{\{j_1=j_4\}}
{\bf 1}_{\{i_1=i_4\}}
{\bf 1}_{\{j_3=j_5\}}
{\bf 1}_{\{i_3=i_5\}}-
$$
$$
-
{\bf 1}_{\{j_6=j_2\}}
{\bf 1}_{\{i_6=i_2\}}
{\bf 1}_{\{j_1=j_3\}}
{\bf 1}_{\{i_1=i_3\}}
{\bf 1}_{\{j_4=j_5\}}
{\bf 1}_{\{i_4=i_5\}}-
$$
$$
-
{\bf 1}_{\{j_6=j_3\}}
{\bf 1}_{\{i_6=i_3\}}
{\bf 1}_{\{j_1=j_5\}}
{\bf 1}_{\{i_1=i_5\}}
{\bf 1}_{\{j_2=j_4\}}
{\bf 1}_{\{i_2=i_4\}}-
$$
$$
-
{\bf 1}_{\{j_6=j_3\}}
{\bf 1}_{\{i_6=i_3\}}
{\bf 1}_{\{j_1=j_4\}}
{\bf 1}_{\{i_1=i_4\}}
{\bf 1}_{\{j_2=j_5\}}
{\bf 1}_{\{i_2=i_5\}}-
$$
$$
-
{\bf 1}_{\{j_3=j_6\}}
{\bf 1}_{\{i_3=i_6\}}
{\bf 1}_{\{j_1=j_2\}}
{\bf 1}_{\{i_1=i_2\}}
{\bf 1}_{\{j_4=j_5\}}
{\bf 1}_{\{i_4=i_5\}}-
$$
$$
-
{\bf 1}_{\{j_6=j_4\}}
{\bf 1}_{\{i_6=i_4\}}
{\bf 1}_{\{j_1=j_5\}}
{\bf 1}_{\{i_1=i_5\}}
{\bf 1}_{\{j_2=j_3\}}
{\bf 1}_{\{i_2=i_3\}}-
$$
$$
-
{\bf 1}_{\{j_6=j_4\}}
{\bf 1}_{\{i_6=i_4\}}
{\bf 1}_{\{j_1=j_3\}}
{\bf 1}_{\{i_1=i_3\}}
{\bf 1}_{\{j_2=j_5\}}
{\bf 1}_{\{i_2=i_5\}}-
$$
$$
-
{\bf 1}_{\{j_6=j_4\}}
{\bf 1}_{\{i_6=i_4\}}
{\bf 1}_{\{j_1=j_2\}}
{\bf 1}_{\{i_1=i_2\}}
{\bf 1}_{\{j_3=j_5\}}
{\bf 1}_{\{i_3=i_5\}}-
$$
$$
-
{\bf 1}_{\{j_6=j_5\}}
{\bf 1}_{\{i_6=i_5\}}
{\bf 1}_{\{j_1=j_4\}}
{\bf 1}_{\{i_1=i_4\}}
{\bf 1}_{\{j_2=j_3\}}
{\bf 1}_{\{i_2=i_3\}}-
$$
$$
-
{\bf 1}_{\{j_6=j_5\}}
{\bf 1}_{\{i_6=i_5\}}
{\bf 1}_{\{j_1=j_2\}}
{\bf 1}_{\{i_1=i_2\}}
{\bf 1}_{\{j_3=j_4\}}
{\bf 1}_{\{i_3=i_4\}}-
$$
\begin{equation}
\label{desy14}
\Biggl.-
{\bf 1}_{\{j_6=j_5\}}
{\bf 1}_{\{i_6=i_5\}}
{\bf 1}_{\{j_1=j_3\}}
{\bf 1}_{\{i_1=i_3\}}
{\bf 1}_{\{j_2=j_4\}}
{\bf 1}_{\{i_2=i_4\}}\Biggr),
\end{equation}

\vspace{3mm}
\noindent
where ${\bf 1}_A$ is the indicator of the set $A$ and

\begin{equation}
\label{ma1}
C_{j_3j_2j_1}^{000}
=\frac{\sqrt{(2j_1+1)(2j_2+1)(2j_3+1)}}{8}(T-t)^{3/2}\bar
C_{j_3j_2j_1}^{000},
\end{equation}

\begin{equation}
\label{ma2}
C_{j_2j_1}^{01}
=\frac{\sqrt{(2j_1+1)(2j_2+1)}}{8}(T-t)^{2}\bar
C_{j_2j_1}^{01},
\end{equation}

\begin{equation}
\label{ma3}
C_{j_2j_1}^{10}
=\frac{\sqrt{(2j_1+1)(2j_2+1)}}{8}(T-t)^{2}\bar
C_{j_2j_1}^{10},
\end{equation}

\begin{equation}
\label{ma4}
C_{j_4j_3j_2j_1}^{0000}=
\frac{\sqrt{(2j_1+1)(2j_2+1)(2j_3+1)(2j_4+1)}}{16}(T-t)^{2}\bar
C_{j_4j_3j_2j_1}^{0000},
\end{equation}

\begin{equation}
\label{ma3a1}
C_{j_2j_1}^{02}=
\frac{\sqrt{(2j_1+1)(2j_2+1)}}{16}(T-t)^{3}\bar
C_{j_2j_1}^{02},
\end{equation}

\newpage
\noindent
\begin{equation}
\label{ma3a2}
C_{j_2j_1}^{20}=
\frac{\sqrt{(2j_1+1)(2j_2+1)}}{16}(T-t)^{3}\bar
C_{j_2j_1}^{20},
\end{equation}

\begin{equation}
\label{ma3a3}
C_{j_2j_1}^{11}=
\frac{\sqrt{(2j_1+1)(2j_2+1)}}{16}(T-t)^{3}\bar
C_{j_2j_1}^{11}, 
\end{equation}

\begin{equation}
\label{ma5}
C_{j_3j_2j_1}^{001}
=\frac{\sqrt{(2j_1+1)(2j_2+1)(2j_3+1)}}{16}(T-t)^{5/2}\bar
C_{j_3j_2j_1}^{001},
\end{equation}

\begin{equation}
\label{ma6}
C_{j_3j_2j_1}^{010}
=\frac{\sqrt{(2j_1+1)(2j_2+1)(2j_3+1)}}{16}(T-t)^{5/2}\bar
C_{j_3j_2j_1}^{010},
\end{equation}

\begin{equation}
\label{ma7}
C_{j_3j_2j_1}^{100}
=\frac{\sqrt{(2j_1+1)(2j_2+1)(2j_3+1)}}{16}(T-t)^{5/2}\bar
C_{j_3j_2j_1}^{100},
\end{equation}

\begin{equation}
\label{ma8}
C_{j_5j_4 j_3 j_2 j_1}^{00000}=
\frac{\sqrt{(2j_1+1)(2j_2+1)(2j_3+1)(2j_4+1)(2j_5+1)}}{32}(T-t)^{5/2}\bar
C_{j_5j_4 j_3 j_2 j_1}^{00000},
\end{equation}

\vspace{-4mm}
\begin{equation}
\label{ma9}
C_{j_4j_3j_2j_1}^{0001}
=\frac{\sqrt{(2j_1+1)(2j_2+1)(2j_3+1)(2j_4+1)}}{32}(T-t)^{3}\bar
C_{j_4j_3j_2j_1}^{0001},
\end{equation}

\begin{equation}
\label{ma10}
C_{j_3j_2j_1}^{0010}
=\frac{\sqrt{(2j_1+1)(2j_2+1)(2j_3+1)(2j_4+1)}}{32}(T-t)^{3}\bar
C_{j_4j_3j_2j_1}^{0010},
\end{equation}

\begin{equation}
\label{ma11}
C_{j_4j_3j_2j_1}^{0100}=
\frac{\sqrt{(2j_1+1)(2j_2+1)(2j_3+1)(2j_4+1)}}{32}(T-t)^{3}\bar
C_{j_3j_2j_1}^{0100},
\end{equation}

\begin{equation}
\label{ma12}
C_{j_4j_3j_2j_1}^{1000}
=\frac{\sqrt{(2j_1+1)(2j_2+1)(2j_3+1)(2j_4+1)}}{32}(T-t)^{3}\bar
C_{j_4j_3j_2j_1}^{1000},
\end{equation}

\vspace{1mm}
$$
C_{j_6j_5j_4 j_3 j_2 j_1}^{000000}
=
$$

\vspace{-4mm}
\begin{equation}
\label{ma13}
=\frac{\sqrt{(2j_1+1)(2j_2+1)(2j_3+1)
(2j_4+1)(2j_5+1)(2j_6+1)}}{64}(T-t)^{3}\bar
C_{j_6j_5j_4 j_3 j_2 j_1}^{000000},
\end{equation}

\vspace{2mm}
\noindent
where
\begin{equation}
\label{ma14}
\bar C_{j_3j_2j_1}^{000}=\int\limits_{-1}^{1}P_{j_3}(z)
\int\limits_{-1}^{z}P_{j_2}(y)
\int\limits_{-1}^{y}
P_{j_1}(x)dx dy dz,
\end{equation}
\begin{equation}
\label{ma15}
\bar C_{j_2j_1}^{01}=-
\int\limits_{-1}^{1}(1+y)P_{j_2}(y)
\int\limits_{-1}^{y}
P_{j_1}(x)dx dy,
\end{equation}
\begin{equation}
\label{ma16}
\bar C_{j_2j_1}^{10}=-
\int\limits_{-1}^{1}P_{j_2}(y)
\int\limits_{-1}^{y}
(1+x)P_{j_1}(x)dx dy,
\end{equation}

\vspace{-2mm}
\begin{equation}
\label{ma17}
\bar C_{j_4j_3j_2j_1}^{0000}=\int\limits_{-1}^{1}P_{j_4}(u)
\int\limits_{-1}^{u}P_{j_3}(z)
\int\limits_{-1}^{z}P_{j_2}(y)
\int\limits_{-1}^{y}
P_{j_1}(x)dx dy dz du,
\end{equation}

\vspace{-2mm}
\begin{equation}
\label{ma17a1}
\bar C_{j_2j_1}^{02}=
\int\limits_{-1}^{1}P_{j_2}(y)(y+1)^2
\int\limits_{-1}^{y}
P_{j_1}(x)dx dy,
\end{equation}

\vspace{-2mm}
\begin{equation}
\label{ma17a2}
\bar C_{j_2j_1}^{20}=
\int\limits_{-1}^{1}P_{j_2}(y)
\int\limits_{-1}^{y}
P_{j_1}(x)(x+1)^2 dx dy,
\end{equation}

\vspace{-2mm}
\begin{equation}
\label{ma17a3}
\bar C_{j_2j_1}^{11}=
\int\limits_{-1}^{1}P_{j_2}(y)(y+1)
\int\limits_{-1}^{y}
P_{j_1}(x)(x+1)dx dy,
\end{equation}

\vspace{-2mm}
\begin{equation}
\label{ma18}
\bar C_{j_3j_2j_1}^{001}=-
\int\limits_{-1}^{1}P_{j_3}(z)(z+1)
\int\limits_{-1}^{z}P_{j_2}(y)
\int\limits_{-1}^{y}
P_{j_1}(x)dx dy dz,
\end{equation}

\vspace{-2mm}
\begin{equation}
\label{ma19}
\bar C_{j_3j_2j_1}^{010}=-
\int\limits_{-1}^{1}P_{j_3}(z)
\int\limits_{-1}^{z}P_{j_2}(y)(y+1)
\int\limits_{-1}^{y}
P_{j_1}(x)dx dy dz,
\end{equation}

\vspace{-2mm}
\begin{equation}
\label{ma20}
\bar C_{j_3j_2j_1}^{100}=-
\int\limits_{-1}^{1}P_{j_3}(z)
\int\limits_{-1}^{z}P_{j_2}(y)
\int\limits_{-1}^{y}
P_{j_1}(x)(x+1)dx dy dz,
\end{equation}

\vspace{-2mm}
\begin{equation}
\label{ma21}
\bar C_{j_5j_4 j_3 j_2 j_1}^{00000}=
\int\limits_{-1}^{1}P_{j_5}(v)
\int\limits_{-1}^{v}P_{j_4}(u)
\int\limits_{-1}^{u}P_{j_3}(z)
\int\limits_{-1}^{z}P_{j_2}(y)
\int\limits_{-1}^{y}
P_{j_1}(x)dx dy dz du dv,
\end{equation}

\newpage
\noindent
\begin{equation}
\label{ma22}
\bar C_{j_4j_3j_2j_1}^{1000}=-
\int\limits_{-1}^{1}P_{j_4}(u)
\int\limits_{-1}^{u}P_{j_3}(z)
\int\limits_{-1}^{z}P_{j_2}(y)
\int\limits_{-1}^{y}
P_{j_1}(x)(x+1)dx dy dz du,
\end{equation}

\vspace{-2mm}
\begin{equation}
\label{ma23}
\bar C_{j_4j_3j_2j_1}^{0100}=-
\int\limits_{-1}^{1}P_{j_4}(u)
\int\limits_{-1}^{u}P_{j_3}(z)
\int\limits_{-1}^{z}P_{j_2}(y)(y+1)
\int\limits_{-1}^{y}
P_{j_1}(x)dx dy dz du,
\end{equation}

\vspace{-2mm}
\begin{equation}
\label{ma24}
\bar C_{j_4j_3j_2j_1}^{0010}=-
\int\limits_{-1}^{1}P_{j_4}(u)
\int\limits_{-1}^{u}P_{j_3}(z)(z+1)
\int\limits_{-1}^{z}P_{j_2}(y)
\int\limits_{-1}^{y}
P_{j_1}(x)dx dy dz du,
\end{equation}

\vspace{-2mm}
\begin{equation}
\label{ma25}
\bar C_{j_4j_3j_2j_1}^{0001}=-
\int\limits_{-1}^{1}P_{j_4}(u)(u+1)
\int\limits_{-1}^{u}P_{j_3}(z)
\int\limits_{-1}^{z}P_{j_2}(y)
\int\limits_{-1}^{y}
P_{j_1}(x)dx dy dz du,
\end{equation}

$$
\bar C_{j_6j_5j_4 j_3 j_2 j_1}^{000000}=
$$
\begin{equation}
\label{ma26}
=
\int\limits_{-1}^{1}P_{j_6}(w)
\int\limits_{-1}^{w}P_{j_5}(v)
\int\limits_{-1}^{v}P_{j_4}(u)
\int\limits_{-1}^{u}P_{j_3}(z)
\int\limits_{-1}^{z}P_{j_2}(y)
\int\limits_{-1}^{y}
P_{j_1}(x)dx dy dz du dv dw;
\end{equation}

\noindent
another notations are the same as in Theorems 1, 2.

\subsection{Optimization of Approximations of Iterated 
It\^{o} Stochastic 
Integrals from the Numerical Schemes (\ref{al1})--(\ref{al5})}

This section is devoted to the optimization
of approximations of iterated 
It\^{o} stochastic 
integrals from the numerical schemes (\ref{al1})--(\ref{al5}).
More precisely, we discuss how to minimize the numbers 
$q, q_1, q_2,\ldots, q_{15}$ from Sect.~2.4.

Suppose that $\varepsilon >0$ is the mean-square 
accuracy of approximation of the iterated 
It\^{o} stochastic integrals (\ref{qqq1x}), i.e.
$$
E^{(i_1\ldots i_k)p}_{(l_1\ldots l_k)T,t}\stackrel{\sf def}{=}{\sf M}\left\{\left(
I_{(l_1\ldots l_k)T,t}^{(i_1\ldots i_k)}-
I_{(l_1\ldots l_k)T,t}^{(i_1\ldots i_k)p}\right)^2\right\}
\le \varepsilon,
$$
where $I_{(l_1\ldots l_k)T,t}^{(i_1\ldots i_k)p},$ $p\in\mathbb{N}$
is the approximation of the iterated
It\^{o} stochastic integral $I_{(l_1\ldots l_k)T,t}^{(i_1\ldots i_k)}$.
Then from (\ref{nach1}) and (\ref{star00011}) we obtain the following 
conditions for choosing the numbers $q, q_1, q_2,\ldots,$ $q_{15}$
for approximations of the iterated It\^{o} stochastic integrals
(\ref{sm10})--(\ref{sm12}) \cite{26}, \cite{26a}, \cite{58}
\begin{equation}
\label{nach2}
E^{(i_1i_2)q}_{(00)T,t}
=\frac{(T-t)^2}{2}\left(\frac{1}{2}-\sum_{i=1}^q
\frac{1}{4i^2-1}\right)\le \varepsilon\ \ \ (i_1\ne i_2),
\end{equation}
\begin{equation}
\label{nach3}
E^{(i_1i_2i_3)q_1}_{(000)T,t}
\le 6\left(\frac{(T-t)^{3}}{6}-\sum_{j_1,j_2,j_3=0}^{q_1}
\left(C_{j_3j_2j_1}^{000}\right)^2\right)\le \varepsilon,
\end{equation}
\begin{equation}
\label{nach4}
E^{(i_1i_2)q_2}_{(01)T,t}\le 2\left(\frac{(T-t)^{4}}{4}-\sum_{j_1,j_2=0}^{q_2}
\left(C_{j_2j_1}^{01}\right)^2\right)\le \varepsilon,
\end{equation}
\begin{equation}
\label{nach5}
E^{(i_1i_2)q_2}_{(10)T,t}\le 2\left(\frac{(T-t)^{4}}{12}-\sum_{j_1,j_2=0}^{q_2}
\left(C_{j_2j_1}^{10}\right)^2\right)\le \varepsilon,
\end{equation}
\begin{equation}
\label{nach6}
E^{(i_1\ldots i_4)q_3}_{(0000)T,t}\le 24\Biggl(\frac{(T-t)^4}{24}-\sum_{j_1,j_2,j_3,j_4=0}^{q_3}
\left(C_{j_4j_3j_2j_1}^{0000}\right)^2\Biggr)\le \varepsilon,
\end{equation}
\begin{equation}
\label{nach7}
~~~~~~~~~E^{(i_1\ldots i_5)q_4}_{(00000)T,t}\le
120\left(\frac{(T-t)^{5}}{120}-\sum_{j_1,j_2,j_3,j_4,j_5=0}^{q_4}
\left(C_{j_5j_4j_3j_2j_1}^{00000}\right)^2\right)\le \varepsilon,
\end{equation}
\begin{equation}
\label{nach8}
E^{(i_1i_2)q_5}_{(20)T,t}\le
2\Biggl(\frac{(T-t)^6}{30}-\sum_{j_2,j_1=0}^{q_5}
\left(C_{j_2j_1}^{20}\right)^2\Biggr)\le \varepsilon,
\end{equation}
\begin{equation}
\label{nach9}
E^{(i_1i_2)q_6}_{(11)T,t}\le 2\Biggl(\frac{(T-t)^6}{18}-\sum_{j_2,j_1=0}^{q_6}
\left(C_{j_2j_1}^{11}\right)^2\Biggr)\le \varepsilon,
\end{equation}
\begin{equation}
\label{nach10}
E^{(i_1i_2)q_7}_{(02)T,t}\le 2\Biggl(\frac{(T-t)^6}{6}-\sum_{j_2,j_1=0}^{q_{7}}
\left(C_{j_2j_1}^{02}\right)^2\Biggr)\le \varepsilon,
\end{equation}
\begin{equation}
\label{nach11}
E^{(i_1i_2i_3)q_8}_{(001)T,t}\le 6\Biggl(\frac{(T-t)^5}{10}-\sum_{j_1,j_2,j_3=0}^{q_8}
\left(C_{j_3j_2j_1}^{001}\right)^2\Biggr)\le \varepsilon,
\end{equation}
\begin{equation}
\label{nach12}
E^{(i_1i_2i_3)q_9}_{(010)T,t}\le
6\Biggl(\frac{(T-t)^{5}}{20}-\sum_{j_1,j_2,j_3=0}^{q_9}
\left(C_{j_3j_2j_1}^{010}\right)^2\Biggr)\le \varepsilon,
\end{equation}
\begin{equation}
\label{nach13}
E^{(i_1i_2i_3)q_{10}}_{(100)T,t}\le
6\Biggl(\frac{(T-t)^{5}}{60}-\sum_{j_1,j_2,j_3=0}^{q_{10}}
\left(C_{j_3j_2j_1}^{100}\right)^2\Biggr)\le \varepsilon,
\end{equation}
\begin{equation}
\label{nach14}
E^{(i_1\ldots i_4)q_{11}}_{(0001)T,t}\le
24\Biggl(\frac{(T-t)^6}{36}-\sum_{j_1,j_2,j_3, j_4=0}^{q_{11}}
\left(C_{j_4j_3j_2j_1}^{0001}\right)^2\Biggr)\le \varepsilon,
\end{equation}
\begin{equation}
\label{nach15}
E^{(i_1\ldots i_4)q_{12}}_{(0010)T,t}\le
24\Biggl(\frac{(T-t)^6}{60}-\sum_{j_1,j_2,j_3, j_4=0}^{q_{12}}
\left(C_{j_4j_3j_2j_1}^{0010}\right)^2\Biggr)\le \varepsilon,
\end{equation}
\begin{equation}
\label{nach16}
E^{(i_1\ldots i_4)q_{13}}_{(0100)T,t}\le
24\Biggl(\frac{(T-t)^{6}}{120}-\sum_{j_1,j_2,j_3, j_4=0}^{q_{13}}
\left(C_{j_4j_3j_2j_1}^{0100}\right)^2\Biggr)\le \varepsilon,
\end{equation}
\begin{equation}
\label{nach17}
E^{(i_1\ldots i_4)q_{14}}_{(1000)T,t}\le
24\Biggl(\frac{(T-t)^{6}}{360}-\sum_{j_1,j_2,j_3, j_4=0}^{q_{14}}
\left(C_{j_4j_3j_2j_1}^{1000}\right)^2\Biggr)\le \varepsilon,
\end{equation}
\begin{equation}
\label{nach18}
~~~~~~~E^{(i_1\ldots i_6)q_{15}}_{(000000)T,t}\le
720\left(\frac{(T-t)^{6}}{720}-\sum_{j_1,j_2,j_3,j_4,j_5,j_6=0}^{q_{15}}
\left(C_{j_6 j_5 j_4 j_3 j_2 j_1}\right)^2\right)\le \varepsilon.
\end{equation}

\vspace{1mm}

Taking into account (\ref{uslov}) and (\ref{ma1})--(\ref{ma26}),
(\ref{nach2})--(\ref{nach18}),
we obtain the following conditions for choosing the numbers
$q, q_1, q_2,\ldots, q_{15}$ for the numerical schemes
(\ref{al1})--(\ref{al5}) (constant $C$ is independent of $T-t$ (see below)).

\vspace{6mm}

\centerline{\bf Milstein scheme (\ref{al1})}
$$
\frac{1}{2}\left(\frac{1}{2}-\sum_{i=1}^q
\frac{1}{4i^2-1}\right)\le C (T-t).
$$

\vspace{6mm}

\centerline{\bf Strong Taylor--It\^{o} scheme (\ref{al2})
with convergence order 1.5}
$$
\frac{1}{2}\left(\frac{1}{2}-\sum_{i=1}^q
\frac{1}{4i^2-1}\right)\le C (T-t)^2,
$$
\begin{equation}
\label{fas1}
~~~~~~~6\left(\frac{1}{6}-\frac{1}{64}\sum_{j_1,j_2,j_3=0}^{q_1}
(2j_1+1)(2j_2+1)(2j_3+1)
\left(\bar C_{j_3j_2j_1}^{000}\right)^2\right)\le C(T-t).
\end{equation}

\vspace{6mm}

\centerline{\bf Strong Taylor--It\^{o} scheme (\ref{al3})
with convergence order 2.0}
$$
\frac{1}{2}\left(\frac{1}{2}-\sum_{i=1}^q
\frac{1}{4i^2-1}\right)\le C (T-t)^3,
$$
\begin{equation}
\label{fas2}
~~~~~6\left(\frac{1}{6}-\frac{1}{64}\sum_{j_1,j_2,j_3=0}^{q_1}
(2j_1+1)(2j_2+1)(2j_3+1)
\left(\bar C_{j_3j_2j_1}^{000}\right)^2\right)\le C(T-t)^2,
\end{equation}
\begin{equation}
\label{fas3}
2\left(\frac{1}{4}-\frac{1}{64}\sum_{j_1,j_2=0}^{q_2}
(2j_1+1)(2j_2+1)
\left(\bar C_{j_2j_1}^{01}\right)^2\right)\le C(T-t),
\end{equation}
\begin{equation}
\label{fas4}
2\left(\frac{1}{12}-\frac{1}{64}\sum_{j_1,j_2=0}^{q_2}
(2j_1+1)(2j_2+1)
\left(\bar C_{j_2j_1}^{10}\right)^2\right)\le C(T-t),
\end{equation}
$$
\hspace{-7mm}24\left(\frac{1}{24}-\frac{1}{256}\sum_{j_1,\ldots,j_4=0}^{q_3}
(2j_1+1)(2j_2+1)(2j_3+1)(2j_4+1)
\left(\bar C_{j_4\ldots j_1}^{0000}\right)^2\right)\le 
$$
\begin{equation}
\label{fas5}
\le C(T-t).
\end{equation}

\vspace{6mm}

\centerline{\bf Strong Taylor--It\^{o} scheme (\ref{al4})
with convergence order 2.5}
$$
\frac{1}{2}\left(\frac{1}{2}-\sum_{i=1}^q
\frac{1}{4i^2-1}\right)\le C (T-t)^4,
$$
\begin{equation}
\label{fas6}
~~~6\left(\frac{1}{6}-\frac{1}{64}\sum_{j_1,j_2,j_3=0}^{q_1}
(2j_1+1)(2j_2+1)(2j_3+1)
\left(\bar C_{j_3j_2j_1}^{000}\right)^2\right)\le C(T-t)^3,
\end{equation}
\begin{equation}
\label{fas7}
~~~~\hspace{-11mm}
2\left(\frac{1}{4}-\frac{1}{64}\sum_{j_1,j_2=0}^{q_2}
(2j_1+1)(2j_2+1)
\left(\bar C_{j_2j_1}^{01}\right)^2\right)\le C(T-t)^2,
\end{equation}
\begin{equation}
\label{fas8}
~~~~~~2\left(\frac{1}{12}-\frac{1}{64}\sum_{j_1,j_2=0}^{q_2}
(2j_1+1)(2j_2+1)
\left(\bar C_{j_2j_1}^{10}\right)^2\right)\le C(T-t)^2,
\end{equation}

\vspace{-1mm}
$$
\hspace{-11mm}24\left(\frac{1}{24}-\frac{1}{256}\sum_{j_1,\ldots,j_4=0}^{q_3}
(2j_1+1)(2j_2+1)(2j_3+1)(2j_4+1)
\left(\bar C_{j_4\ldots j_1}^{0000}\right)^2\right)\le 
$$
\begin{equation}
\label{fas9}
\le C(T-t)^2,
\end{equation}
$$
6\left(\frac{1}{10}-\frac{1}{256}\sum_{j_1,j_2,j_3=0}^{q_8}
(2j_1+1)(2j_2+1)(2j_3+1)
\left(\bar C_{j_3j_2j_1}^{001}\right)^2\right)\le 
$$
\begin{equation}
\label{fas10}
\le C(T-t),
\end{equation}
$$
6\left(\frac{1}{20}-\frac{1}{256}\sum_{j_1,j_2,j_3=0}^{q_9}
(2j_1+1)(2j_2+1)(2j_3+1)
\left(\bar C_{j_3j_2j_1}^{010}\right)^2\right)\le 
$$
\begin{equation}
\label{fas11}
\le C(T-t),
\end{equation}
$$
6\left(\frac{1}{60}-\frac{1}{256}\sum_{j_1,j_2,j_3=0}^{q_{10}}
(2j_1+1)(2j_2+1)(2j_3+1)
\left(\bar C_{j_3j_2j_1}^{100}\right)^2\right)\le 
$$
\begin{equation}
\label{fas12}
\le C(T-t),
\end{equation}

\vspace{-5mm}
$$
120\Biggl(\frac{1}{120}-\frac{1}{32^2}\sum_{j_1,\ldots,j_5=0}^{q_{4}}
(2j_1+1)(2j_2+1)(2j_3+1)(2j_4+1)(2j_5+1)\times\Biggr.
$$
\begin{equation}
\label{fas13}
~~~~~~~\Biggl.\times\left(\bar C_{j_5\ldots j_1}^{00000}\right)^2\Biggr)\le C(T-t).
\end{equation}

\vspace{5mm}

\centerline{\bf Strong Taylor--It\^{o} scheme (\ref{al5})
with convergence order 3.0}
$$
\frac{1}{2}\left(\frac{1}{2}-\sum_{i=1}^q
\frac{1}{4i^2-1}\right)\le C (T-t)^5,
$$
\begin{equation}
\label{fas14}
~~~6\left(\frac{1}{6}-\frac{1}{64}\sum_{j_1,j_2,j_3=0}^{q_1}
(2j_1+1)(2j_2+1)(2j_3+1)
\left(\bar C_{j_3j_2j_1}^{000}\right)^2\right)\le C(T-t)^4,
\end{equation}
\begin{equation}
\label{fas15}
\hspace{-11mm}
2\left(\frac{1}{4}-\frac{1}{64}\sum_{j_1,j_2=0}^{q_2}
(2j_1+1)(2j_2+1)
\left(\bar C_{j_2j_1}^{01}\right)^2\right)\le C(T-t)^3,
\end{equation}
\begin{equation}
\label{fas16}
~~2\left(\frac{1}{12}-\frac{1}{64}\sum_{j_1,j_2=0}^{q_2}
(2j_1+1)(2j_2+1)
\left(\bar C_{j_2j_1}^{10}\right)^2\right)\le C(T-t)^3,
\end{equation}

\vspace{-2mm}
$$
\hspace{-11mm}24\left(\frac{1}{24}-\frac{1}{256}\sum_{j_1,\ldots,j_4=0}^{q_3}
(2j_1+1)(2j_2+1)(2j_3+1)(2j_4+1)
\left(\bar C_{j_4\ldots j_1}^{0000}\right)^2\right)\le 
$$
\begin{equation}
\label{fas17}
\le C(T-t)^3,
\end{equation}
$$
6\left(\frac{1}{10}-\frac{1}{256}\sum_{j_1,j_2,j_3=0}^{q_8}
(2j_1+1)(2j_2+1)(2j_3+1)
\left(\bar C_{j_3j_2j_1}^{001}\right)^2\right)\le 
$$
\begin{equation}
\label{fas18}
\le C(T-t)^2,
\end{equation}
$$
6\left(\frac{1}{20}-\frac{1}{256}\sum_{j_1,j_2,j_3=0}^{q_9}
(2j_1+1)(2j_2+1)(2j_3+1)
\left(\bar C_{j_3j_2j_1}^{010}\right)^2\right)\le 
$$
\begin{equation}
\label{fas19}
\le C(T-t)^2,
\end{equation}
$$
6\left(\frac{1}{60}-\frac{1}{256}\sum_{j_1,j_2,j_3=0}^{q_{10}}
(2j_1+1)(2j_2+1)(2j_3+1)
\left(\bar C_{j_3j_2j_1}^{100}\right)^2\right)\le 
$$
\begin{equation}
\label{fas20}
\le C(T-t)^2,
\end{equation}

\vspace{-5mm}
$$
120\Biggl(\frac{1}{120}-\frac{1}{32^2}\sum_{j_1,\ldots,j_5=0}^{q_{4}}
(2j_1+1)(2j_2+1)(2j_3+1)(2j_4+1)(2j_5+1)\times\Biggr.
$$
\begin{equation}
\label{fas21}
~~~~~~~\Biggl.\times
\left(\bar C_{j_5\ldots j_1}^{00000}\right)^2\Biggr)\le C(T-t)^2,
\end{equation}

\vspace{-2mm}
\begin{equation}
\label{fas22}
~~~~~~~~~2\left(\frac{1}{30}-\frac{1}{256}\sum_{j_1,j_2=0}^{q_5}
(2j_1+1)(2j_2+1)
\left(\bar C_{j_2j_1}^{20}\right)^2\right)\le C(T-t),
\end{equation}
\begin{equation}
\label{fas23}
~~~~~~~~~2\left(\frac{1}{18}-\frac{1}{256}\sum_{j_1,j_2=0}^{q_6}
(2j_1+1)(2j_2+1)
\left(\bar C_{j_2j_1}^{11}\right)^2\right)\le C(T-t),
\end{equation}
\begin{equation}
\label{fas24}
~~~~~~~~~2\left(\frac{1}{6}-\frac{1}{256}\sum_{j_1,j_2=0}^{q_7}
(2j_1+1)(2j_2+1)
\left(\bar C_{j_2j_1}^{02}\right)^2\right)\le C(T-t),
\end{equation}
$$
24\left(\frac{1}{36}-\frac{1}{32^2}\sum_{j_1,\ldots,j_4=0}^{q_{11}}
(2j_1+1)(2j_2+1)(2j_3+1)(2j_4+1)
\left(\bar C_{j_4\ldots j_1}^{0001}\right)^2\right)\le 
$$
\begin{equation}
\label{fas25}
\le C(T-t),
\end{equation}
$$
24\left(\frac{1}{60}-\frac{1}{32^2}\sum_{j_1,\ldots,j_4=0}^{q_{12}}
(2j_1+1)(2j_2+1)(2j_3+1)(2j_4+1)
\left(\bar C_{j_4\ldots j_1}^{0010}\right)^2\right)\le 
$$
\begin{equation}
\label{fas26}
\le C(T-t),
\end{equation}
$$
24\left(\frac{1}{120}-\frac{1}{32^2}\sum_{j_1,\ldots,j_4=0}^{q_{13}}
(2j_1+1)(2j_2+1)(2j_3+1)(2j_4+1)
\left(\bar C_{j_4\ldots j_1}^{0100}\right)^2\right)\le 
$$
\begin{equation}
\label{fas27}
\le C(T-t),
\end{equation}
$$
24\left(\frac{1}{360}-\frac{1}{32^2}\sum_{j_1,\ldots,j_4=0}^{q_{14}}
(2j_1+1)(2j_2+1)(2j_3+1)(2j_4+1)
\left(\bar C_{j_4\ldots j_1}^{1000}\right)^2\right)\le 
$$
\begin{equation}
\label{fas28}
\le C(T-t),
\end{equation}
$$
720\Biggl(\frac{1}{720}-\frac{1}{64^2}\sum_{j_1,\ldots,j_6=0}^{q_{15}}
(2j_1+1)(2j_2+1)(2j_3+1)(2j_4+1)(2j_5+1)(2j_6+1)\times\Biggr.
$$
\begin{equation}
\label{fas29}
\Biggl.
~~~~~~\times
\left(\bar C_{j_6\ldots j_1}^{000000}\right)^2\Biggr)\le C(T-t).
\end{equation}

\vspace{1mm}

Taking into account Theorem 8 and the results of Listings 
\ref{lst:q_hypothesis} and \ref{lst:q_hypothesis2} (see Sect.~5)
we decided to exclude the multiplier factors $k!$
from the left-hand sides 
of (\ref{fas1}), 
(\ref{fas2})--(\ref{fas5}), (\ref{fas6})--(\ref{fas13}),
(\ref{fas14})--(\ref{fas29}). The detailed 
numerical confirmation of the mentioned possibility
can be found in \cite{63}.
This means that we will use the following conditions 
for choosing the numbers
$q, q_1, q_2,\ldots, q_{15}$ for the numerical schemes
(\ref{al1})--(\ref{al5}) (constant $C$ is independent of $T-t$ (see below)).

\vspace{5mm}

\centerline{\bf Milstein scheme (\ref{al1})}
\begin{equation}
\label{ress1}
\frac{1}{2}\left(\frac{1}{2}-\sum_{i=1}^q
\frac{1}{4i^2-1}\right)\le C (T-t).
\end{equation}

\vspace{5mm}

\centerline{\bf Strong Taylor--It\^{o} scheme (\ref{al2})
with convergence order 1.5}
\begin{equation}
\label{ress2}
\frac{1}{2}\left(\frac{1}{2}-\sum_{i=1}^q
\frac{1}{4i^2-1}\right)\le C (T-t)^2,
\end{equation}
\begin{equation}
\label{fas1a}
\frac{1}{6}-\frac{1}{64}\sum_{j_1,j_2,j_3=0}^{q_1}
(2j_1+1)(2j_2+1)(2j_3+1)
\left(\bar C_{j_3j_2j_1}^{000}\right)^2\le C(T-t).
\end{equation}

\vspace{5mm}

\centerline{\bf Strong Taylor--It\^{o} scheme (\ref{al3})
with convergence order 2.0}

\begin{equation}
\label{ress3}
\frac{1}{2}\left(\frac{1}{2}-\sum_{i=1}^q
\frac{1}{4i^2-1}\right)\le C (T-t)^3,
\end{equation}
\begin{equation}
\label{fas2a}
\frac{1}{6}-\frac{1}{64}\sum_{j_1,j_2,j_3=0}^{q_1}
(2j_1+1)(2j_2+1)(2j_3+1)
\left(\bar C_{j_3j_2j_1}^{000}\right)^2\le C(T-t)^2,
\end{equation}
\begin{equation}
\label{fas3a}
\frac{1}{4}-\frac{1}{64}\sum_{j_1,j_2=0}^{q_2}
(2j_1+1)(2j_2+1)
\left(\bar C_{j_2j_1}^{01}\right)^2\le C(T-t),
\end{equation}
\begin{equation}
\label{fas4a}
\frac{1}{12}-\frac{1}{64}\sum_{j_1,j_2=0}^{q_2}
(2j_1+1)(2j_2+1)
\left(\bar C_{j_2j_1}^{10}\right)^2\le C(T-t),
\end{equation}
\begin{equation}
\label{fas5a}
\frac{1}{24}-\frac{1}{256}\sum_{j_1,\ldots,j_4=0}^{q_3}
(2j_1+1)(2j_2+1)(2j_3+1)(2j_4+1)
\left(\bar C_{j_4\ldots j_1}^{0000}\right)^2\le C(T-t).
\end{equation}

\vspace{5mm}

\centerline{\bf Strong Taylor--It\^{o} scheme (\ref{al4})
with convergence order 2.5}
\begin{equation}
\label{ress4}
\frac{1}{2}\left(\frac{1}{2}-\sum_{i=1}^q
\frac{1}{4i^2-1}\right)\le C (T-t)^4,
\end{equation}
\begin{equation}
\label{fas6a}
\frac{1}{6}-\frac{1}{64}\sum_{j_1,j_2,j_3=0}^{q_1}
(2j_1+1)(2j_2+1)(2j_3+1)
\left(\bar C_{j_3j_2j_1}^{000}\right)^2\le C(T-t)^3,
\end{equation}
\begin{equation}
\label{fas7a}
\frac{1}{4}-\frac{1}{64}\sum_{j_1,j_2=0}^{q_2}
(2j_1+1)(2j_2+1)
\left(\bar C_{j_2j_1}^{01}\right)^2\le C(T-t)^2,
\end{equation}
\begin{equation}
\label{fas8a}
\frac{1}{12}-\frac{1}{64}\sum_{j_1,j_2=0}^{q_2}
(2j_1+1)(2j_2+1)
\left(\bar C_{j_2j_1}^{10}\right)^2\le C(T-t)^2,
\end{equation}
\begin{equation}
\label{fas9a}
\frac{1}{24}-\frac{1}{256}\sum_{j_1,\ldots,j_4=0}^{q_3}
(2j_1+1)(2j_2+1)(2j_3+1)(2j_4+1)
\left(\bar C_{j_4\ldots j_1}^{0000}\right)^2\le C(T-t)^2,
\end{equation}
\begin{equation}
\label{fas10a}
\frac{1}{10}-\frac{1}{256}\sum_{j_1,j_2,j_3=0}^{q_8}
(2j_1+1)(2j_2+1)(2j_3+1)
\left(\bar C_{j_3j_2j_1}^{001}\right)^2\le C(T-t),
\end{equation}
\begin{equation}
\label{fas11a}
\frac{1}{20}-\frac{1}{256}\sum_{j_1,j_2,j_3=0}^{q_9}
(2j_1+1)(2j_2+1)(2j_3+1)
\left(\bar C_{j_3j_2j_1}^{010}\right)^2\le C(T-t),
\end{equation}
\begin{equation}
\label{fas12a}
\frac{1}{60}-\frac{1}{256}\sum_{j_1,j_2,j_3=0}^{q_{10}}
(2j_1+1)(2j_2+1)(2j_3+1)
\left(\bar C_{j_3j_2j_1}^{100}\right)^2\le C(T-t),
\end{equation}
$$
\frac{1}{120}-\frac{1}{32^2}\sum_{j_1,\ldots,j_5=0}^{q_{4}}
(2j_1+1)(2j_2+1)(2j_3+1)(2j_4+1)(2j_5+1)
\left(\bar C_{j_5\ldots j_1}^{00000}\right)^2\le 
$$
\begin{equation}
\label{fas13a}
\le C(T-t).
\end{equation}

\vspace{5mm}

\centerline{\bf Strong Taylor--It\^{o} scheme (\ref{al5}) with convergence order 3.0}
\begin{equation}
\label{ress5}
\frac{1}{2}\left(\frac{1}{2}-\sum_{i=1}^q
\frac{1}{4i^2-1}\right)\le C (T-t)^5,
\end{equation}
\begin{equation}
\label{fas14a}
\frac{1}{6}-\frac{1}{64}\sum_{j_1,j_2,j_3=0}^{q_1}
(2j_1+1)(2j_2+1)(2j_3+1)
\left(\bar C_{j_3j_2j_1}^{000}\right)^2\le C(T-t)^4,
\end{equation}
\begin{equation}
\label{fas15a}
\frac{1}{4}-\frac{1}{64}\sum_{j_1,j_2=0}^{q_2}
(2j_1+1)(2j_2+1)
\left(\bar C_{j_2j_1}^{01}\right)^2\le C(T-t)^3,
\end{equation}
\begin{equation}
\label{fas16a}
\frac{1}{12}-\frac{1}{64}\sum_{j_1,j_2=0}^{q_2}
(2j_1+1)(2j_2+1)
\left(\bar C_{j_2j_1}^{10}\right)^2\le C(T-t)^3,
\end{equation}
\begin{equation}
\label{fas17a}
\frac{1}{24}-\frac{1}{256}\sum_{j_1,\ldots,j_4=0}^{q_3}
(2j_1+1)(2j_2+1)(2j_3+1)(2j_4+1)
\left(\bar C_{j_4\ldots j_1}^{0000}\right)^2\le C(T-t)^3,
\end{equation}
\begin{equation}
\label{fas18a}
\frac{1}{10}-\frac{1}{256}\sum_{j_1,j_2,j_3=0}^{q_8}
(2j_1+1)(2j_2+1)(2j_3+1)
\left(\bar C_{j_3j_2j_1}^{001}\right)^2\le C(T-t)^2,
\end{equation}
\begin{equation}
\label{fas19a}
\frac{1}{20}-\frac{1}{256}\sum_{j_1,j_2,j_3=0}^{q_9}
(2j_1+1)(2j_2+1)(2j_3+1)
\left(\bar C_{j_3j_2j_1}^{010}\right)^2\le C(T-t)^2,
\end{equation}
\begin{equation}
\label{fas20a}
\frac{1}{60}-\frac{1}{256}\sum_{j_1,j_2,j_3=0}^{q_{10}}
(2j_1+1)(2j_2+1)(2j_3+1)
\left(\bar C_{j_3j_2j_1}^{100}\right)^2\le C(T-t)^2,
\end{equation}
$$
\frac{1}{120}-\frac{1}{32^2}\sum_{j_1,\ldots,j_5=0}^{q_{4}}
(2j_1+1)(2j_2+1)(2j_3+1)(2j_4+1)(2j_5+1)
\left(\bar C_{j_5\ldots j_1}^{00000}\right)^2\le 
$$
\begin{equation}
\label{fas21a}
\le C(T-t)^2,
\end{equation}
\begin{equation}
\label{fas22a}
\frac{1}{30}-\frac{1}{256}\sum_{j_1,j_2=0}^{q_5}
(2j_1+1)(2j_2+1)
\left(\bar C_{j_2j_1}^{20}\right)^2\le C(T-t),
\end{equation}
\begin{equation}
\label{fas23a}
\frac{1}{18}-\frac{1}{256}\sum_{j_1,j_2=0}^{q_6}
(2j_1+1)(2j_2+1)
\left(\bar C_{j_2j_1}^{11}\right)^2\le C(T-t),
\end{equation}
\begin{equation}
\label{fas24a}
\frac{1}{6}-\frac{1}{256}\sum_{j_1,j_2=0}^{q_7}
(2j_1+1)(2j_2+1)
\left(\bar C_{j_2j_1}^{02}\right)^2\le C(T-t),
\end{equation}
\begin{equation}
\label{fas25a}
\frac{1}{36}-\frac{1}{32^2}\sum_{j_1,\ldots,j_4=0}^{q_{11}}
(2j_1+1)(2j_2+1)(2j_3+1)(2j_4+1)
\left(\bar C_{j_4\ldots j_1}^{0001}\right)^2\le C(T-t),
\end{equation}
\begin{equation}
\label{fas26a}
\frac{1}{60}-\frac{1}{32^2}\sum_{j_1,\ldots,j_4=0}^{q_{12}}
(2j_1+1)(2j_2+1)(2j_3+1)(2j_4+1)
\left(\bar C_{j_4\ldots j_1}^{0010}\right)^2\le C(T-t),
\end{equation}
\begin{equation}
\label{fas27a}
\frac{1}{120}-\frac{1}{32^2}\sum_{j_1,\ldots,j_4=0}^{q_{13}}
(2j_1+1)(2j_2+1)(2j_3+1)(2j_4+1)
\left(\bar C_{j_4\ldots j_1}^{0100}\right)^2\le C(T-t),
\end{equation}
\begin{equation}
\label{fas28a}
\frac{1}{360}-\frac{1}{32^2}\sum_{j_1,\ldots,j_4=0}^{q_{14}}
(2j_1+1)(2j_2+1)(2j_3+1)(2j_4+1)
\left(\bar C_{j_4\ldots j_1}^{1000}\right)^2\le C(T-t),
\end{equation}
$$
\frac{1}{720}-\frac{1}{64^2}\sum_{j_1,\ldots,j_6=0}^{q_{15}}
(2j_1+1)(2j_2+1)(2j_3+1)(2j_4+1)(2j_5+1)(2j_6+1)\times
$$
\begin{equation}
\label{fas29a}
\times
\left(\bar C_{j_6\ldots j_1}^{000000}\right)^2\le C(T-t).
\end{equation}

\subsection{Approximations of Iterated 
Stratonovich Stochastic 
Integrals from the Numerical Schemes (\ref{al1x})--(\ref{al5x})
Using Legendre Polynomials}

This section is devoted to 
approximation of the Stratonovich 
stochastic integrals (\ref{qqq1xx})
of multiplicities 1 to 6 based on 
Theorems 3--7. At that we will use
multiple Fourier--Legendre series for 
approximation of the mentioned stochastic integrals.

The numerical schemes (\ref{al1x})--(\ref{al5x})
contain the following set (see (\ref{qqq1xx}))
of iterated Stratonovich 
stochastic integrals 
\begin{equation}
\label{sm10x}
I_{(0)T,t}^{*(i_1)},\ \ \ I_{(1)T,t}^{*(i_1)},\ \ \ I_{(2)T,t}^{*(i_1)},\ \ \ 
I_{(00)T,t}^{*(i_1 i_2)},\ \ \ I_{(10)T,t}^{*(i_1 i_2)},\ \ \ 
I_{(01)T,t}^{*(i_1 i_2)},\ \ \ I_{(000)T,t}^{*(i_1 i_2 i_3)},\ \ \ 
I_{(0000)T,t}^{*(i_1 i_2 i_3 i_4)},\ \ \
\end{equation}
\begin{equation}
\label{sm11x}
I_{(00000)T,t}^{*(i_1 i_2 i_3 i_4 i_5)},\ \ \ 
I_{(02)T,t}^{*(i_1 i_2)},\ \ \ I_{(20)T,t}^{*(i_1 i_2)},\ \ \ 
I_{(11)T,t}^{*(i_1 i_2)},\ \ \ 
I_{(100)T,t}^{*(i_1 i_2 i_3)},\ \ \ I_{(010)T,t}^{*(i_1 i_2 i_3)},\ \ \ 
I_{(001)T,t}^{*(i_1 i_2 i_3)},
\end{equation}
\begin{equation}
\label{sm12x}
I_{(0001)T,t}^{*(i_1 i_2 i_3 i_4)},\ \ \ 
I_{(0010)T,t}^{*(i_1 i_2 i_3 i_4)},\ \ \ 
I_{(0100)T,t}^{*(i_1 i_2 i_3 i_4)},\ \ \ I_{(1000)T,t}^{*(i_1 i_2 i_3 i_4)},\ \ 
\ I_{(000000)T,t}^{*(i_1 i_2 i_3 i_4 i_5 i_6)}.
\end{equation}

Using Theorems 3--7 and 
well known properties of the Legendre polynomials,
we obtain the following formulas for numerical 
modeling of the stochastic integrals
(\ref{sm10x})--(\ref{sm12x})
\cite{26}, \cite{26a}, \cite{40}-\cite{42aa}, \cite{58}, \cite{59},
\cite{60}-\cite{62}

\begin{equation}
\label{key1}
I_{(0)T,t}^{*(i_1)}=\sqrt{T-t}\zeta_0^{(i_1)},
\end{equation}
\begin{equation}
\label{key2}
I_{(1)T,t}^{*(i_1)}=-\frac{(T-t)^{3/2}}{2}\left(\zeta_0^{(i_1)}+
\frac{1}{\sqrt{3}}\zeta_1^{(i_1)}\right),
\end{equation}
\begin{equation}
\label{key3}
I_{(2)T,t}^{*(i_1)}=\frac{(T-t)^{5/2}}{3}\left(\zeta_0^{(i_1)}+
\frac{\sqrt{3}}{2}\zeta_1^{(i_1)}+
\frac{1}{2\sqrt{5}}\zeta_2^{(i_1)}\right),
\end{equation}
\begin{equation}
\label{dr1}
~~~~~~~~I_{(00)T,t}^{*(i_1 i_2)q}=
\frac{T-t}{2}\left(\zeta_0^{(i_1)}\zeta_0^{(i_2)}+\sum_{i=1}^{q}
\frac{1}{\sqrt{4i^2-1}}\left(
\zeta_{i-1}^{(i_1)}\zeta_{i}^{(i_2)}-
\zeta_i^{(i_1)}\zeta_{i-1}^{(i_2)}\right)\right),
\end{equation}
\begin{equation}
\label{key4}
I_{(000)T,t}^{*(i_1i_2i_3)q_1}
=
\sum_{j_1,j_2,j_3=0}^{q_1}
C_{j_3j_2j_1}^{000}
\zeta_{j_1}^{(i_1)}\zeta_{j_2}^{(i_2)}\zeta_{j_3}^{(i_3)},
\end{equation}
\begin{equation}
\label{key5}
I_{(10)T,t}^{*(i_1 i_2)q_2}=
\sum_{j_1,j_2=0}^{q_2}
C_{j_2j_1}^{10}\zeta_{j_1}^{(i_1)}\zeta_{j_2}^{(i_2)},
\end{equation}
\begin{equation}
\label{key6}
I_{(01)T,t}^{*(i_1 i_2)q_2}=
\sum_{j_1,j_2=0}^{q_2}
C_{j_2j_1}^{01}\zeta_{j_1}^{(i_1)}\zeta_{j_2}^{(i_2)},
\end{equation}
\begin{equation}
\label{key7}
I_{(0000)T,t}^{(i_1 i_2 i_3 i_4)q_3}
=
\sum_{j_1,j_2,j_3,j_4=0}^{q_3}
C_{j_4 j_3 j_2 j_1}^{0000}
\zeta_{j_1}^{(i_1)}\zeta_{j_2}^{(i_2)}\zeta_{j_3}^{(i_3)}\zeta_{j_4}^{(i_4)},
\end{equation}

\begin{equation}
\label{key8}
I_{(00000)T,t}^{(i_1 i_2 i_3 i_4 i_5)q_4}
=
\sum_{j_1,j_2,j_3,j_4,j_5=0}^{q_4}
C_{j_5 j_4 j_3 j_2 j_1}^{00000}
\zeta_{j_1}^{(i_1)}\zeta_{j_2}^{(i_2)}\zeta_{j_3}^{(i_3)}
\zeta_{j_4}^{(i_4)}\zeta_{j_5}^{(i_5)},
\end{equation}

\begin{equation}
\label{key9}
I_{(20)T,t}^{*(i_1 i_2)q_5}=
\sum_{j_1,j_2=0}^{q_5}
C_{j_2j_1}^{20}\zeta_{j_1}^{(i_1)}\zeta_{j_2}^{(i_2)},
\end{equation}
\begin{equation}
\label{key10}
I_{(11)T,t}^{(i_1 i_2)q_6}=
\sum_{j_1,j_2=0}^{q_6}
C_{j_2j_1}^{11}\zeta_{j_1}^{(i_1)}\zeta_{j_2}^{(i_2)},
\end{equation}
\begin{equation}
\label{key11}
I_{(02)T,t}^{(i_1 i_2)q_{7}}=
\sum_{j_1,j_2=0}^{q_{7}}
C_{j_2j_1}^{02}\zeta_{j_1}^{(i_1)}\zeta_{j_2}^{(i_2)},
\end{equation}
\begin{equation}
\label{key12}
I_{(001)T,t}^{*(i_1i_2i_3)q_8}
=
\sum_{j_1,j_2,j_3=0}^{q_8}
C_{j_3j_2j_1}^{001}
\zeta_{j_1}^{(i_1)}\zeta_{j_2}^{(i_2)}\zeta_{j_3}^{(i_3)},
\end{equation}
\begin{equation}
\label{key13}
I_{(010)T,t}^{(i_1i_2i_3)q_9}
=
\sum_{j_1,j_2,j_3=0}^{q_9}
C_{j_3j_2j_1}^{010}
\zeta_{j_1}^{(i_1)}\zeta_{j_2}^{(i_2)}\zeta_{j_3}^{(i_3)},
\end{equation}
\begin{equation}
\label{key14}
I_{(100)T,t}^{(i_1i_2i_3)q_{10}}
=
\sum_{j_1,j_2,j_3=0}^{q_{10}}
C_{j_3j_2j_1}^{100}
\zeta_{j_1}^{(i_1)}\zeta_{j_2}^{(i_2)}\zeta_{j_3}^{(i_3)},
\end{equation}
\begin{equation}
\label{key15}
I_{(0001)T,t}^{(i_1 i_2 i_3 i_4)q_{11}}
=
\sum_{j_1,j_2,j_3,j_4=0}^{q_{11}}
C_{j_4 j_3 j_2 j_1}^{0001}
\zeta_{j_1}^{(i_1)}\zeta_{j_2}^{(i_2)}\zeta_{j_3}^{(i_3)}\zeta_{j_4}^{(i_4)},
\end{equation}
\begin{equation}
\label{key16}
I_{(0010)T,t}^{(i_1 i_2 i_3 i_4)q_{12}}
=\sum_{j_1,j_2,j_3,j_4=0}^{q_{12}}
C_{j_4 j_3 j_2 j_1}^{0010}
\zeta_{j_1}^{(i_1)}\zeta_{j_2}^{(i_2)}\zeta_{j_3}^{(i_3)}\zeta_{j_4}^{(i_4)},
\end{equation}
\begin{equation}
\label{key17}
I_{(0100)T,t}^{(i_1 i_2 i_3 i_4)q_{13}}
=
\sum_{j_1,j_2,j_3,j_4=0}^{q_{13}}
C_{j_4 j_3 j_2 j_1}^{0100}
\zeta_{j_1}^{(i_1)}\zeta_{j_2}^{(i_2)}\zeta_{j_3}^{(i_3)}\zeta_{j_4}^{(i_4)},
\end{equation}
\begin{equation}
\label{key18}
I_{(1000)T,t}^{(i_1 i_2 i_3 i_4)q_{14}}
=
\sum_{j_1,j_2,j_3,j_4=0}^{q_{14}}
C_{j_4 j_3 j_2 j_1}^{1000}
\zeta_{j_1}^{(i_1)}\zeta_{j_2}^{(i_2)}\zeta_{j_3}^{(i_3)}\zeta_{j_4}^{(i_4)},
\end{equation}
\begin{equation}
\label{key19}
~~~~~~~~I_{(000000)T,t}^{(i_1 i_2 i_3 i_4 i_5 i_6)q_{15}}
=\sum_{j_1,j_2,j_3,j_4,j_5,j_6=0}^{q_{15}}
C_{j_6 j_5 j_4 j_3 j_2 j_1}^{000000}
\zeta_{j_1}^{(i_1)}\zeta_{j_2}^{(i_2)}\zeta_{j_3}^{(i_3)}\zeta_{j_4}^{(i_4)}
\zeta_{j_5}^{(i_5)}\zeta_{j_6}^{(i_6)},
\end{equation}

\vspace{2mm}
\noindent
where ${\bf 1}_A$ is the indicator of the set $A;$
another notations are the same as in Sect.~2.4.

The question on choosing the numbers $q_1, q_2,\ldots,q_{15}$
in (\ref{key4})--(\ref{key19}) turned out to be nontrivial 
\cite{26}, \cite{26a}, \cite{58} (Chapter 5).
The expansions (\ref{key4})--(\ref{key19}) for iterated Stratonovich stochastic integrals
are simpler 
than their analogues (\ref{desy4})--(\ref{desy14}) for
iterated It\^{o} stochastic integrals. However, the calculation of the mean-square
approximation error for iterated Stratonovich
stochastic integrals turns out to be much more difficult than for 
iterated It\^{o} stochastic integrals \cite{26}, \cite{26a}, \cite{58} (Chapter 5).
Below we give some reasoning regarding this problem.

Denote
$$
E^{*(i_1\ldots i_k)p}_{(l_1\ldots l_k)T,t}\stackrel{\sf def}{=}{\sf M}\left\{\left(
I_{(l_1\ldots l_k)T,t}^{*(i_1\ldots i_k)}-
I_{(l_1\ldots l_k)T,t}^{*(i_1\ldots i_k)p}\right)^2\right\},
$$
where $I_{(l_1\ldots l_k)T,t}^{*(i_1\ldots i_k)p},$ $p\in\mathbb{N}$
is the approximation of the iterated
Stratonovich stochastic integral $I_{(l_1\ldots l_k)T,t}^{*(i_1\ldots i_k)}$.

From (\ref{dr1}) for $i_1\ne i_2$ we obtain \cite{26}, \cite{26a}, \cite{58}
$$
E^{*(i_1i_2)q}_{(00)T,t}
=\frac{(T-t)^2}{2}
\sum\limits_{i=q+1}^{\infty}\frac{1}{4i^2-1}
\le \frac{(T-t)^2}{2}\int\limits_{q}^{\infty}
\frac{1}{4x^2-1}dx
=
$$
\begin{equation}
\label{teacxx}
=-\frac{(T-t)^2}{8}{\rm ln}\left|
1-\frac{2}{2q+1}\right|\le C_1\frac{(T-t)^2}{q},
\end{equation}
where constant $C_1$ is independent of $q$.

It is easy to notice that for a sufficiently
small $T-t$ (recall that $T-t\ll 1$ since it is a step of integration
for numerical schemes for It\^{o} SDEs) there 
exists a constant $C_2$ such that
\begin{equation}
\label{teac3xx}
E^{*(i_1\ldots i_k)q}_{(l_1\ldots l_k)T,t}
\le C_2 E^{*(i_1i_2)q}_{(00)T,t}.
\end{equation}

From (\ref{teacxx}) and (\ref{teac3xx}) we finally obtain
\begin{equation}
\label{teac4}
E^{*(i_1\ldots i_k)q}_{(l_1\ldots l_k)T,t}
\le C \frac{(T-t)^2}{q},
\end{equation}
where constant $C$ does not depend on $T-t$.
The same idea can be found in \cite{2} in the framework of 
the method based
on the trigonometric expansion of the
Brownian bridge process. Note that, in contrast to (\ref{teac4}), 
the constant $C$ in Theorems 4--6
does not depend on $q.$

Obviously, we can get more information about the numbers $q_1, q_2,\ldots,q_{15}$ (these
numbers are different for different iterated Stratonovich
stochastic integrals)
using the another approach.
Since 
$$
J^{*}[\psi^{(k)}]_{T,t}=J[\psi^{(k)}]_{T,t}\ \ \ \hbox{w.\ p.\ 1}
$$
for pairwise different $i_1,\ldots,i_k=1,\ldots,m$, 
where $J[\psi^{(k)}]_{T,t},$ $J^{*}[\psi^{(k)}]_{T,t}$
are defined by (\ref{ito}) and (\ref{str}) correspondingly,
then 
for pairwise different 
$i_1,\ldots,i_6=1,\ldots,m$ from (\ref{formula0}) 
we obtain \cite{26}, \cite{26a}, \cite{58}
\begin{equation}
\label{form1}
E^{*(i_1i_2)q}_{(00)T,t}=
\frac{(T-t)^2}{2}\left(\frac{1}{2}-\sum_{i=1}^q
\frac{1}{4i^2-1}\right),
\end{equation}
\begin{equation}
\label{form2}
E^{*(i_1i_2i_3)q_1}_{(000)T,t}=
\frac{(T-t)^{3}}{6}-\sum_{j_3,j_2,j_1=0}^{q_1}
\left(C_{j_3j_2j_1}^{000}\right)^2,
\end{equation}
\begin{equation}
\label{form3}
E^{*(i_1i_2)q_2}_{(01)T,t}=
\frac{(T-t)^{4}}{4}-\sum_{j_1,j_2=0}^{q_2}
\left(C_{j_2j_1}^{01}\right)^2,
\end{equation}
\begin{equation}
\label{form4}
E^{*(i_1i_2)q_2}_{(10)T,t}=
\frac{(T-t)^{4}}{12}-\sum_{j_1,j_2=0}^{q_2}
\left(C_{j_2j_1}^{10}\right)^2,
\end{equation}
\begin{equation}
\label{form5}
E^{*(i_1\ldots i_4)q_3}_{(0000)T,t}=
\frac{(T-t)^{4}}{24}-\sum_{j_1,j_2,j_3,j_4=0}^{q_3}
\left(C_{j_4j_3j_2j_1}^{0000}\right)^2,
\end{equation}
\begin{equation}
\label{form6}
E^{*(i_1\ldots i_5)q_4}_{(00000)T,t}=
\frac{(T-t)^{5}}{120}-\sum_{j_1,j_2,j_3,j_4,j_5=0}^{q_4}
\left(C_{j_5 i_4 i_3 i_2 j_1}^{00000}\right)^2,
\end{equation}
\begin{equation}
\label{form7}
E^{*(i_1i_2)q_5}_{(20)T,t}=
\frac{(T-t)^6}{30}-\sum_{j_2,j_1=0}^{q_5}
\left(C_{j_2j_1}^{20}\right)^2,
\end{equation}
\begin{equation}
\label{form8}
E^{*(i_1i_2)q_6}_{(11)T,t}=
\frac{(T-t)^6}{18}-\sum_{j_2,j_1=0}^{q_6}
\left(C_{j_2j_1}^{11}\right)^2,
\end{equation}
\begin{equation}
\label{form9}
E^{*(i_1i_2)q_7}_{(02)T,t}=
\frac{(T-t)^6}{6}-\sum_{j_2,j_1=0}^{q_7}
\left(C_{j_2j_1}^{02}\right)^2,
\end{equation}
\begin{equation}
\label{form10}
E^{*(i_1i_2i_3)q_8}_{(001)T,t}=
\frac{(T-t)^5}{10}-\sum_{j_1,j_2,j_3=0}^{q_8}
\left(C_{j_3j_2j_1}^{001}\right)^2,
\end{equation}
\begin{equation}
\label{form11}
E^{*(i_1i_2i_3)q_9}_{(010)T,t}=
\frac{(T-t)^{5}}{20}-\sum_{j_1,j_2,j_3=0}^{q_9}
\left(C_{j_3j_2j_1}^{010}\right)^2,
\end{equation}
\begin{equation}
\label{form12}
E^{*(i_1i_2i_3)q_{10}}_{(100)T,t}=
\frac{(T-t)^{5}}{60}-\sum_{j_1,j_2,j_3=0}^{q_{10}}
\left(C_{j_3j_2j_1}^{100}\right)^2,
\end{equation}
\begin{equation}
\label{form13}
E^{*(i_1\ldots i_4)q_{11}}_{(0001)T,t}=
\frac{(T-t)^6}{36}-\sum_{j_1,j_2,j_3, j_4=0}^{q_{11}}
\left(C_{j_4j_3j_2j_1}^{0001}\right)^2,
\end{equation}
\begin{equation}
\label{form14}
E^{*(i_1\ldots i_4)q_{12}}_{(0010)T,t}=
\frac{(T-t)^6}{60}-\sum_{j_1,j_2,j_3, j_4=0}^{q_{12}}
\left(C_{j_4j_3j_2j_1}^{0010}\right)^2,
\end{equation}
\begin{equation}
\label{form15}
E^{*(i_1\ldots i_4)q_{13}}_{(0100)T,t}=
\frac{(T-t)^{6}}{120}-\sum_{j_1,j_2,j_3, j_4=0}^{q_{13}}
\left(C_{j_4j_3j_2j_1}^{0100}\right)^2,
\end{equation}
\begin{equation}
\label{form16}
E^{*(i_1\ldots i_4)q_{14}}_{(1000)T,t}=
\frac{(T-t)^{6}}{360}-\sum_{j_1,j_2,j_3, j_4=0}^{q_{14}}
\left(C_{j_4j_3j_2j_1}^{1000}\right)^2,
\end{equation}
\begin{equation}
\label{form17}
E^{*(i_1\ldots i_6)q_{15}}_{(000000)T,t}=
\frac{(T-t)^{6}}{720}-\sum_{j_1,j_2,j_3,j_4,j_5,j_6=0}^{q_{15}}
\left(C_{j_6 j_5 j_4 j_3 j_2 j_1}^{000000}\right)^2.
\end{equation}

\vspace{2mm}

Taking into account (\ref{form1})--(\ref{form17})
and the results of paper \cite{63},
we use in the SDE-MATH software package the conditions from Table 1
for choosing the numbers
$q, q_1, q_2,\ldots, q_{15}$ for the numerical schemes
(\ref{al1x})--(\ref{al5x}).

\begin{figure}
\begin{center}
\centerline{Table 1. High-order strong Taylor--Stratonovich schemes.}
\vspace{4mm}
\begin{tabular}{|c|c|c|}
\hline
Order of convergence&Scheme&Conditions for choosing the numbers $q, q_1,\ldots,q_{15}$\\
\hline
1.0&(\ref{al1x})&(\ref{ress1})\\
\hline
1.5 &(\ref{al2x})&(\ref{ress2}), (\ref{fas1a})\\
\hline
2.0 &(\ref{al3x})&(\ref{ress3})--(\ref{fas5a})\\
\hline
2.5 &(\ref{al4x})&(\ref{ress4})--(\ref{fas13a})\\
\hline
3.0 &(\ref{al5x})&(\ref{ress5})--(\ref{fas29a})\\
\hline
\end{tabular}
\end{center}
\end{figure}

Note that in the SDE-MATH software package, which is presented in the following sections,
we use the following upper bounds $b$ on the numbers $q_1,\ldots, q_{15}$
$$
b=56\ \ \ \hbox{for}\ \ \ q_1,\ \ \ \ \ \ 
b=15\ \ \ \hbox{for}\ \ \ q_2,\ q_3,\ \ \ \ \ \ b=6\ \ \ \hbox{for}\ \ \ 
q_4,\ q_8,\ q_9,\ q_{10},
$$
$$
b=2\ \ \ \hbox{for}\ \ \ q_5,\ q_6,\ q_7,\ q_{11},\ q_{12},\ q_{13},\ q_{14},\ q_{15}.
$$

This means that for the implementing of the numerical methods
(\ref{al2})--(\ref{al5}) and (\ref{al2x})--(\ref{al5x})
we use in the SDE-MATH software package the following quantities of the 
exactly calculated Fourier--Legendre
coefficients
$$
57^3=185,\hspace{-1mm}193\ \ \ \hbox{for}\ \ \ C_{j_3j_2j_1}^{000},
$$
$$
16^3=4,\hspace{-1mm}096\ \ \ \hbox{for each of}\ \ \ C_{j_2j_1}^{10},\ C_{j_2j_1}^{01},
$$
$$
16^4=65,\hspace{-1mm}536\ \ \ \hbox{for}\ \ \ C_{j_4j_3j_2j_1}^{0000},
$$ 
$$
7^3=343\ \ \ \hbox{for each of}\ \ \ C_{j_3j_2j_1}^{100},\ 
C_{j_3j_2j_1}^{010},\ C_{j_3j_2j_1}^{001},
$$
$$
7^5=16,\hspace{-1mm}807\ \ \ \hbox{for}\ \ \ C_{j_5j_4j_3j_2j_1}^{00000},
$$ 
$$
3^2=9\ \ \ \hbox{for each of}\ \ \ C_{j_2j_1}^{20},\ C_{j_2j_1}^{02},\
C_{j_2j_1}^{11},
$$
$$
3^4=81\ \ \ \hbox{for each of}\ \ \ C_{j_4j_3j_2j_1}^{1000},\ 
C_{j_4j_3j_2j_1}^{0100},\ C_{j_4j_3j_2j_1}^{0010},\ 
C_{j_4j_3j_2j_1}^{0001},
$$
$$
3^6=729 \ \ \ \hbox{for}\ \ \ C_{j_6j_5j_4j_3j_2j_1}^{000000}.
$$

It should be noted that unlike the method based on Theorems 1--7, 8,
existing and well-known approaches to the mean-square approximation 
of iterated stochastic integrals based on the trigonometric basis functions
\cite{2}, \cite{3}, \cite{7}, \cite{27},
\cite{28}, \cite{34}, \cite{37} do not allow choosing theoretically
different numbers $q$ for approximations of different 
iterated stochastic integrals (starting form the multiplicity 2
of stochastic integrals). Moreover, the noted
approaches \cite{2}, \cite{3}, \cite{7}, \cite{27},
\cite{28}, \cite{34}, \cite{37} exclude the possibility
for obtaining of approximate and exact expressions
for the mean-square approximation error similar to the formulas
(\ref{star00011}), (\ref{tttr11}).

\subsection{Numerical Algorithm for Linear Stationary
Systems of It\^{o} SDEs Based on Spectral Decomposition}

Consider the following linear stationary
system of It\^{o} SDEs
\begin{equation}
\label{lin1}
d{\bf x}_t=\left(A{\bf x}_t + B{\bf u}(t)\right)dt+Fd{\bf w}_t,\ \ \ 
{\bf x}_0={\bf x}(0),\ \ \ t\in [0, T],
\end{equation}

\noindent
where ${\bf x}_t\in \mathbb{R}^n$ is a solution of the system (\ref{lin1}),
${\bf u}(t):$ $[0, T]\to\mathbb{R}^k$ is a nonrandom function,
$A\in \mathbb{R}^{n\times n},$
$F \in \mathbb{R}^{n\times m}$, $B\in \mathbb{R}^{n\times k}$, and
${\bf w}_t$ is a standard $m$-dimensional Wiener process 
with independent components ${\bf w}_t^{(i)},$ $i=1,\ldots,m.$
Also we suppose that $n,$ $m,$ $k\ge 1.$
The process ${\bf y}_t=H{\bf x}_t\in \mathbb{R}^{1}$ is interpreted as an output process
of the system (\ref{lin1}), where
$H\in \mathbb{R}^{1\times n}.$

It is well-known that the solution of (\ref{lin1}) has the form \cite{4}
\begin{equation}
\label{lin2}
{\bf x}_{t} =e^{A(t-t_0)}{\bf x}_{t_0} +
\int\limits_{t_0}^{t}  e^{A(t-s)}B{\bf u}(s)ds +
\int\limits_{t_0}^{t} e^{A(t-s)}Fd{\bf w}_{s},\ \ \ 0\le t_0\le t\le T,
\end{equation}
where  
$e^C$ is a matrix exponent
$$
\sum\limits_{j=0}^{\infty}\frac{C^j}{j!}\stackrel{\sf def}{=}e^{C},
$$
$C$ is a square matrix, and $C^0\stackrel{\sf def}{=}I$ is a unity matrix.

Consider the partition $\{\tau_p\}_{p=0}^{N}$ of $[0, T]$ such that
$\tau_p=p\Delta,$ $\Delta>0.$
For simplicity, we will suppose that ${\bf u}(s),$ $s\in [0, T]$
can be approximated by the step function, i.e.
${\bf u}(s)\approx \hat {\bf u}(s),$ $s\in [0, T],$ where $\hat {\bf u}(s)={\bf u}(\tau_p)$ 
for $s\in[\tau_p,\tau_{p+1}),$
$p=0, 1,\ldots,$ $N-1$ (more accurate approximations of ${\bf u}(s)$
are discussed in \cite{57} (also see \cite{53}, \cite{56})).
Substituting $t=\tau_{p+1},$ $t_0=\tau_p,$ and $\hat {\bf u}(s)$ instead
of ${\bf u}(s)$
into (\ref{lin2}), we obtain
\begin{equation}
\label{lin3}
\hat {\bf x}_{p+1}=e^{A\Delta}\hat {\bf x}_p+A^{-1}\bigl(
e^{A\Delta}-I\bigr)B{\bf u}(p\Delta) +\tilde{\bf w}_{p+1}(\Delta),\ \ \ {\bf x}_0={\bf x}(0),
\end{equation}
where $\hat {\bf x}_{p}$ is the approximation of ${\bf x}_{\tau_p}$ and
$$
\int\limits_{0}^{\Delta} e^{A(\Delta-s)}Fd{\bf w}_{s+p\Delta}
\stackrel{\sf def}{=}
\tilde {\bf w}_{p+1}(\Delta).
$$

Also we assume that $\hat {\bf y}_{p}=H\hat {\bf x}_{p}$, 
where $\hat {\bf y}_{p}$ is the approximation of ${\bf y}_{\tau_p}.$
The random column $\tilde {\bf w}_{p+1}(\Delta)$ admits the following representation
\cite{4}
\begin{equation}
\label{lin5}
\tilde {\bf w}_{p+1}(\Delta)=S_{D}(\Delta)
\Lambda_D(\Delta)\bar{\bf w}_{p+1},
\end{equation}
where $\bar{\bf w}_{p}\in \mathbb{R}^n$ is a column of independent
standard Gaussian random variables such that
${\sf M}\left\{\bar{\bf w}_{p}\bar{\bf w}_{q}^{\sf T}\right\}=\cal O$ 
for $p\ne q$, $\cal O$ is a zero matrix of size $n\times n$,
$S_{D}(\Delta)$ is a matrix of orthonormal 
eigenvectors of the matrix
$D_f(\Delta)$ and $\Lambda_D^2(\Delta)$ is a diagonal matrix on the main diagonal of which
are the eigenvalues of the matrix $D_f(\Delta)$, the matrix $D_f(\Delta)$
is defined by
$$
D_f(\Delta)={\sf M}\left\{\tilde {\bf w}_{p+1}(\Delta)
\tilde {\bf w}_{p+1}^{\sf T}(\Delta)
\right\}=\int\limits_0^{\Delta}{\rm exp}(A(\Delta-s))
FF^{\sf T} {\rm exp}(A^{\sf T}(\Delta-s))ds,
$$
where $C^{\sf T}$ is a transposed matrix $C$.
Moreover, $D_f(\Delta)=D_f(t)\bigl\vert_{t=\Delta}\bigr.$, where $D_f(t)$
is a solution of the following Cauchy problem \cite{4}
$$
\frac{d{D}_{f}}{dt}(t) = AD_{f}(t) + D_{f}(t)A^{\sf T} +
FF^{\sf T},\ \ \ 
D_{f}(0)=\cal O.
$$

In the SDE-MATH software package, we implement the 
numerical modeling of the system (\ref{lin1})
by the formulas (\ref{lin3}), (\ref{lin5}).
At that we use Algorithms 2.3--2.6 from \cite{57} (also see \cite{56}, Chapter 11)
for the implemetation of (\ref{lin5}).

%% file: chapters/software.tex
\footnotetext[1]{All programs in Python programming language
from this paper were written by the first author}

\section{The Structure of the SDE-MATH Software Package}

\subsection{Development Tools}

The software package was implemented with
Python programming language. The main
reason to use it is a huge community and
significant amount of helpful libraries
for calculations and mathematics. The
development was performed in free
to use Atom text editor\footnotemark[1].

\subsection{Dependency Libraries}

\vspace{1mm}

In the development of the SDE-MATH software package
such libraries as SymPy, NumPy, PyQt5, and
Matplotlib were involved. All these libraries
and tools are free and open source.

\vspace{1mm}

\begin{itemize}
\item SymPy is a Python library
able to perform symbolic algebra
calculations.
\item NumPy is a library which
specialization is efficient
mathematical calculations. Most
part of this library is written
in C programming language that
guarantees high calculation performance.
\item The database is SQLite3.
This is a tiny database for a
local usage on one machine.
\item Matplotlib library is a piece
of software used to present
obtained results in a best way.
\item PyQt5 is a library used to build 
graphical user interface for the SDE-MATH software package.
\end{itemize}

\vspace{-1mm}
\subsection{Architecture}

\vspace{1mm}

\begin{figure}[H]
    \vspace{7mm}
    \centering
    \includegraphics[width=.9\textwidth]{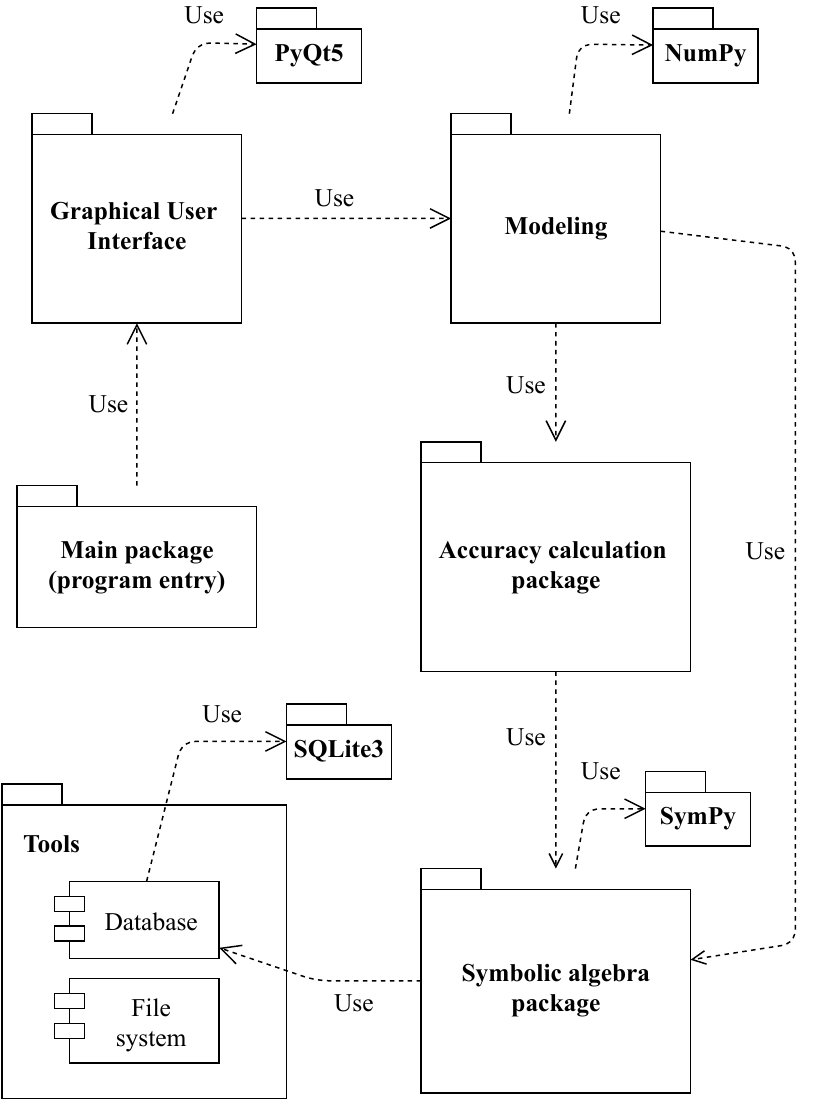}
    \caption{The SDE-MATH software package architecture\label{fig:architecture}}
\end{figure}

Taking into account, that the SDE-MATH
software package is oriented on a
numerical modeling its architecture
is clear. There are two main statements.
The first is that mathematical formulas
are strongly integrated with SymPy
library. By that we mean that they
completely rely on SymPy. And the
second is usage of database to make
some calculations able for caching.
The architecture itself is provided
on Figure \ref{fig:architecture}. Here
all parts of the software package can be seen.

The main package is responsible
for startup, so it decides which
part of the software package must
be started. The software package has 
several modes of operation.
The objectives now are

\vspace{1mm}

\begin{itemize}
\item Run program to calculate
and store the Fourier--Legendre coefficients in
few text files with further loading in database.
\item Run program with graphical user 
interface. This is the main program entry for the SDE-MATH software package.
\end{itemize}

\newpage
\newpage

On the current state of development the main entry package booting up
PyQt5 library with all necessary widgets. More detailed description
of this process will be provided later.

Moving further, the 
modeling package comes up. This
package responsible for all work
referenced to modeling including
initialization of modeling environment,
calculations loops and more. Also, it
depends on accuracy calculation module
deciding which amount of members in
each approximation of iterated stochastic 
integral should be used in modeling
of the It\^{o} SDE (\ref{1.5.2}) solution.

Accuracy calculation module accepts the
order of strong numerical scheme for the It\^o SDE 
(\ref{1.5.2}) and its integration step and then calculates 
necessary amount of members in approximations of 
iterated It\^{o} and Stratonovich stochastic integrals.

Symbolic algebra module is the
construction part which combines
many supplementary differential operators with
strong numerical schemes for the It\^{o} SDE (\ref{1.5.2}). Having
these components combined this module
performs simplification of resulting formula
so the modeling package can do its modeling work.

Tools module provides some functionality
related to bootstrap of runtime environment and
external instruments such as database and
file system.

\vspace{-2mm}

\subsubsection{Integration with SymPy}

Class inheritance tree was extended
to implement strong numerical schemes for It\^{o} SDEs.
While numerical schemes for It\^{o} SDEs were being implemented it
was also necessary to implement 
supplementary subprograms. SymPy is a
Python library able to perform symbolic
algebra calculations. This is a core part
of the project since it is differentiates
input functions, builds and simplifies
strong numerical schemes for It\^{o} SDEs to model
the It\^{o} SDE (\ref{1.5.2}) solution. Without
this part the program package cannot be able
to provide such flexible input of data.

\vspace{-2mm}

\subsubsection{Purpose of NumPy}

NumPy is a library that helps with
calculation optimizations in this project.
The library specialization is efficient
mathematical calculations. The main usage
case is to calculate compiled symbolic
formulas with it. It has integration with
SymPy to replace symbolic functions with
high performance numerical functions.

\subsubsection{Purpose of SQLite Database}

The database was used to store
the precalculated Fourier--Le\-gen\-dre coefficients, so getting
them from there made numerical modeling
much faster, because calculation process
for these Fourier--Legendre coefficients involve high-cost
symbolic operations. The database contains
only one table, and might be thought
redundant, but modeling needs hundreds (or even thousands)
of precalculated coefficients. Obviously,
calculation of them at runtime is terribly
inefficient, but text files also not the
best choice. Text files provide a sequential
access memory and combining different
accuracy values $q_1,\ldots, q_{15}$ it causes sequential
search which extends time to give the result.
That is where database comes up. The random
access allows to get any Fourier--Legendre coefficient or any
quantity of them which makes solution as
flexible as it possible.

The download of precalculated Fourier--Legendre coefficients
is built in supplemental subprograms to
provide fluent calculation pipeline. Having
the precalculated Fourier--Legendre coefficient not found,
subprogram initiates calculation for it
with following store in the database.

\begin{figure}[H]
    \centering
    \includegraphics[width=.9\textwidth]{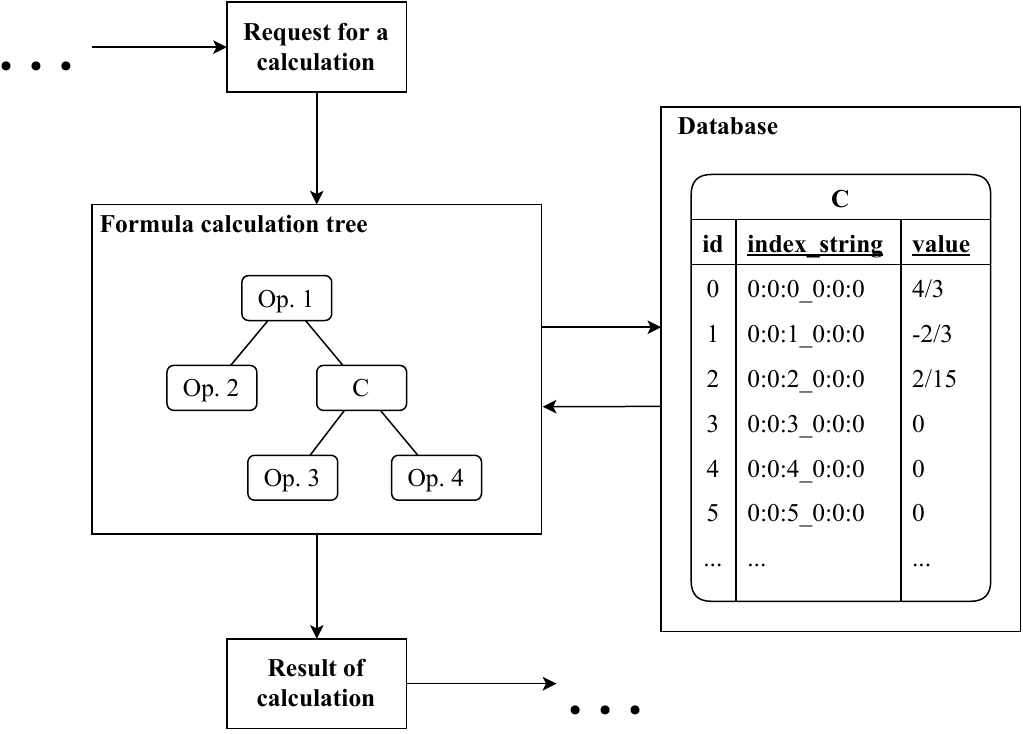}
    \caption{Fourier--Legendre coefficients calculations explanation\label{fig:calculations}}
\end{figure}

It is an interesting note that having
mentioned earlier optimization done, the
calculations performance were increased
in several times. Now the most heavy operation
is symbolic simplifications before modeling.
The actual modeling takes seconds,
for thousands of iterations on $m$ components
of stochastic process, so it is not such
important how long modeling period of time,
as the accuracy that needs to be accomplished.

The scheme of calculation process is
presented on Figure \ref{fig:calculations}.

\subsubsection{Purpose of Matplotlib}

Matplotlib library is a piece of
software used to present obtained
results in a best way. This library
has many features, but feature that
needed in this project is to print
charts with modeling results in a
PyQt5 widget. Thus the data visualization
is integrated in graphical user interface.

\subsection{Implementation Plan}

The implementation of SDE-MATH software
package was performed sequentially.
The components of SDE-MATH software package
were implemented in order of their
necessity for calculation pipeline
completion.

\subsubsection{Calculation of the Fourier--Legendre Coefficients}

The Fourier--Legendre coefficients for the approximations
of iterated It\^{o} and Stratonovich stochastic integrals
were implemented and placed in
Listings \ref{lst:polinomial}--\ref{lst:c000000}. This was the first step since
the Fourier--Legendre coefficients involved almost in every
strong numerical scheme for the It\^o SDE (\ref{1.5.2}).

Also it is important to note that the SDE-MATH software package contains a Python script
intended for generating of Fourier--Legendre coefficients using multiprocessing. This script 
placed in Listing \ref{lst:main_new_c} and already contains tasks that 
were performed to generate about 300,000 Fourier--Legendre coefficients.
Similarly, user can run and calculate additional Fourier--Legendre coefficients if they are needed.
To determine which Fourier--Legendre coefficients will be calculated user must specify pairs of 
starting and ending values of components in lower multi-index and specify upper 
multi-index of the Fourier--Legendre coefficient. For example 
$(((0, 15), (0, 15), (0, 15)), (0, 1, 0))$. This means that program 
calculates the Fourier--Legendre coefficients $C_{j_3j_2j_1}^{010},$ 
where $j_1, j_2, j_3 \in \{0, 1, \ldots, 14\}$.

\subsubsection{Differential Operators $L,$ ${\bar L},$ $G_0^{(i)},$ $i = 1, \ldots, m$}

Moving further, strong numerical schemes for It\^o SDEs
rely on the differential operators (\ref{2.3}), (\ref{2.4}), and (\ref{2.4xxx}). They
were implemented and placed in Listings \ref{lst:l}--\ref{lst:lj}.

\vspace{-1mm}

\subsubsection{Approximations of Iterated Stochastic Integrals}

The next step is implementation of approximations
of iterated It\^{o} and Stra\-to\-no\-vich
stochastic integrals for the numerical
schemes (\ref{al1})--(\ref{al5}), (\ref{al1x})--(\ref{al5x}). 
They are implemented and definition of
their classes are placed in
Listings \ref{lst:i0}--\ref{lst:j000000}.

\vspace{-1mm}

\subsubsection{Strong Numerical Schemes for It\^{o} SDEs}

The strong numerical schemes (\ref{al1})--(\ref{al5}),
(\ref{al1x})--(\ref{al5x}) for It\^{o} SDEs
were implemented. They are placed in Listings 
\ref{lst:euler_dri}--\ref{lst:strong_taylor_stratonovich_3p0}.

\vspace{-1mm}

\subsubsection{Graphical User Interface}

Finally, the graphical user interface was implemented.
The source codes referenced to graphical user interface
are placed in Listings \ref{lst:gui_config}--\ref{lst:gui_linear_step8}.

\vspace{-2mm}

\section{Software Package Graphical User Interface}

For the SDE-MATH software package mentioned above the graphical user interface was developed.
The graphical user interface is important and massive part of SDE-MATH software package because
it allows user to perform modeling experiments without programming skills and understanding of
program package architecture and principles of work.

\vspace{-2mm}

\subsection{Information Model of The Graphical User Interface}

The development of graphical user interface was started from consideration of

\begin{figure}[H]
    \centering
    \includegraphics[width=.9\textwidth]{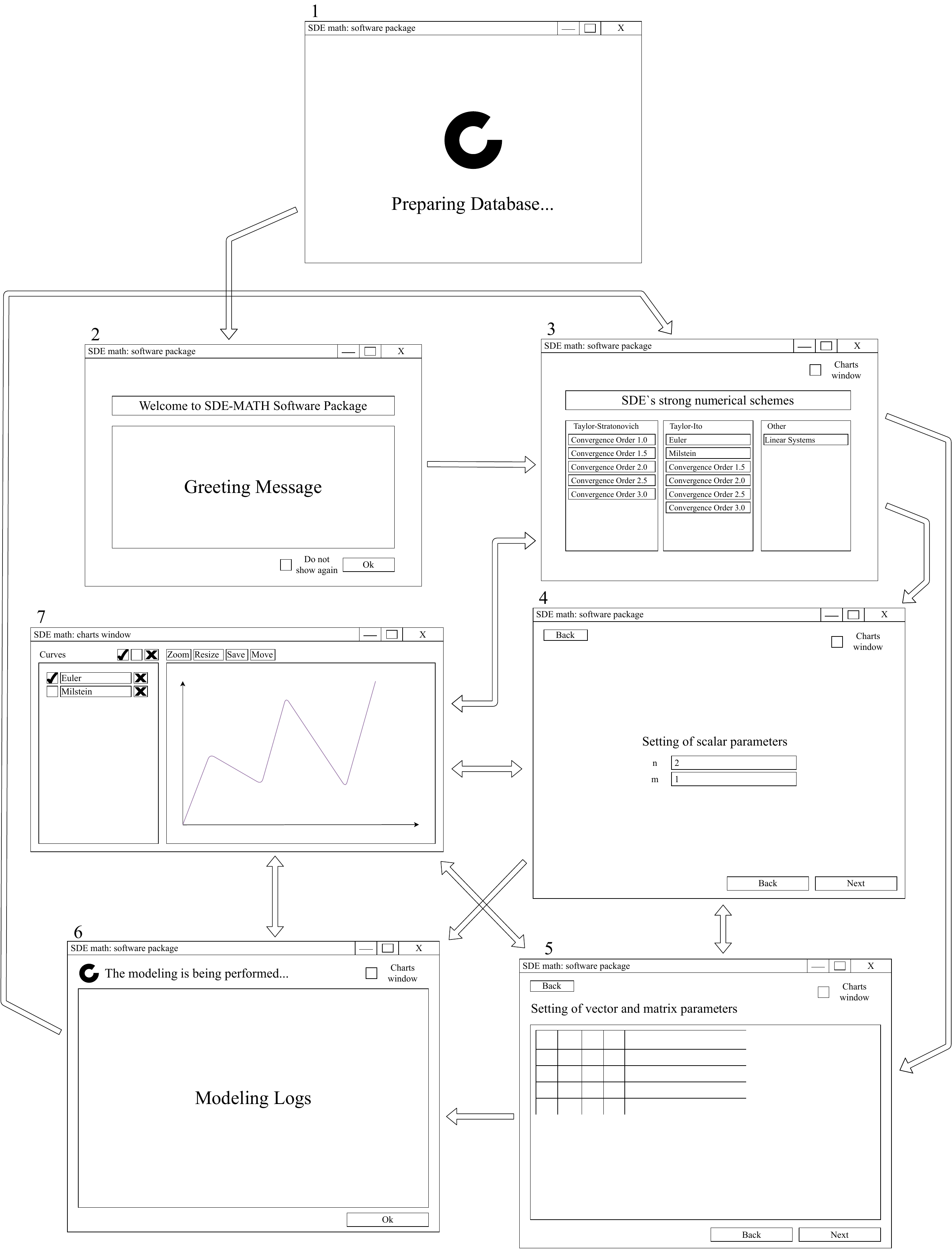}
    \caption{Information model of graphical user interface\label{fig:info_model}}
\end{figure}

\noindent experiments and routines which can be performed with the SDE-MATH software package.
The graphical user interface is aimed on provision of user capabilities to perform 
nonlinear and linear systems of It\^{o} SDEs modeling experiments. The information model which schematically
describes the graphical user interface structure is presented on Figure \ref{fig:info_model}.

\vspace{-3mm}

\subsubsection{Processing Screens}
To represent long duration processes the graphical user interface has two dialogs which can be seen 
on Figure \ref{fig:info_model}, in Windows 1 and 6. The first one represents database preparing 
process on application very first run. During this process the Fourier--Legendre coefficients are being loaded
into the SQLite database. This screen appears also when user calculates new Fourier--Legendre coefficients.
The other screen shows logs during modeling experiment.

\vspace{-3mm}

\subsubsection{Greetings Dialog}
After the SDE-MATH software package has completed the database preparation, it shows greeting dialog which 
represents short information about its purposes. The greeting dialog can be seen on Figure \ref{fig:info_model}
in Window 2.

\vspace{-3mm}

\subsubsection{Main Menu Dialog}
In the main menu of the SDE-MATH software package user can choose one of strong numerical schemes for It\^o SDEs
to perform modeling experiments. The main menu dialog can be seen on Figure \ref{fig:info_model}
in Window 3.

\vspace{-3mm}

\subsubsection{Visualization Tool}
It is important to note that the main SDE-MATH software package window has a checkbox in right upper corner
which do switching on and off charts window. In any time user can call this window or hide it if it
is not needed. The charts window is universal utility for modeling experiments
results visualization. This window has few instruments on it. The left side bar contains all curves
labels, and control elements for hiding, showing, and deleting curves. On the right side of the window 
there are plot which draws the curves. The charts window can be seen on Figure \ref{fig:info_model},
it is Window 3.

\vspace{-3mm}

\subsubsection{Data Input Dialogs}
Since the software package has options to perform linear It\^o SDEs modeling experiments
it is necessary to provide user with input fields for numerical data both scalar and matrix.
On the other side, for nonlinear It\^o SDEs it is necessary to provide symbolic input.
The choice of control elements is conditioned by the above obstacles. On Figure \ref{fig:info_model}, 
and especially in Windows 4 and 5, these input controls can be seen. There are "LineEditWidget" and "TableWidget" 
which are sufficient to provide input abilities. The topic of input data validation is also 
important but to be more accurate referenced to user experience rather than to information model, so
it will be described further.

\vspace{-3mm}

\subsection{The User Experience and Implementation Results}

\begin{figure}[H]
    \centering
    \includegraphics[width=.7\textwidth]{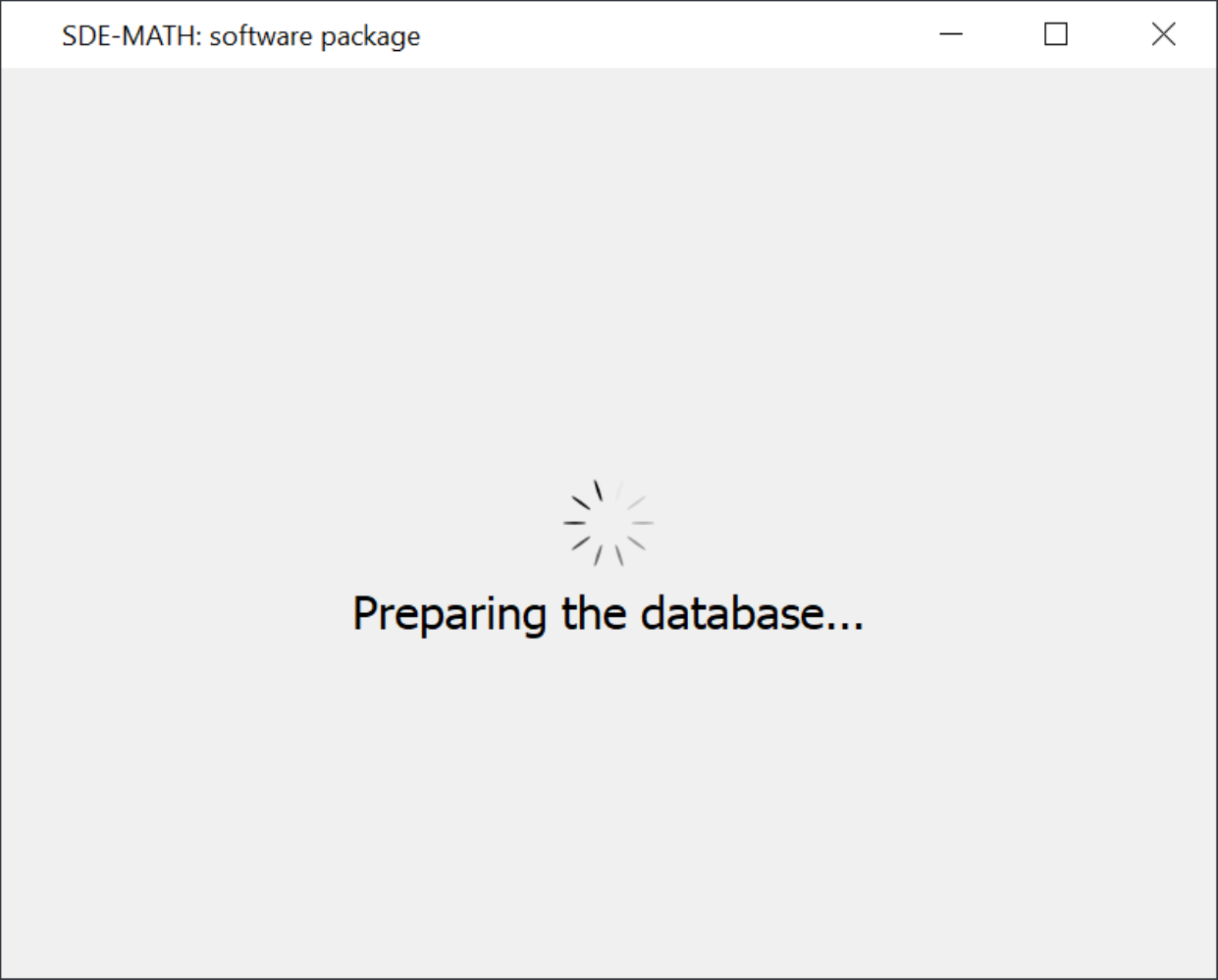}
    \caption{Fourier--Legendre coefficients database preparation screen\label{fig:preparing}}
\end{figure}

\begin{figure}[H]
    \vspace{6mm}
    \centering
    \includegraphics[width=.7\textwidth]{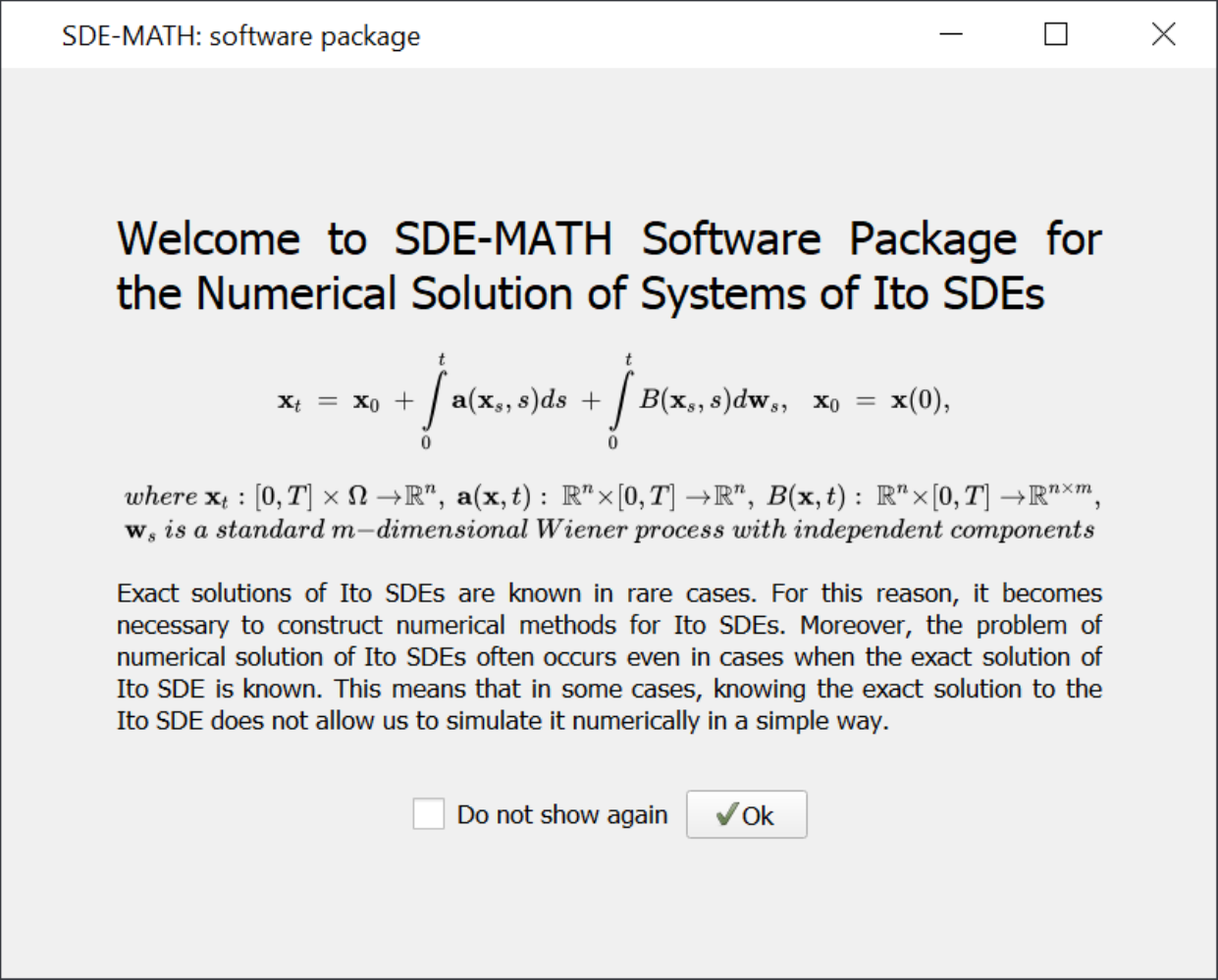}
    \caption{Greetings screen\label{fig:greet}}
\end{figure}

\begin{figure}[H]
    \centering
    \includegraphics[width=.7\textwidth]{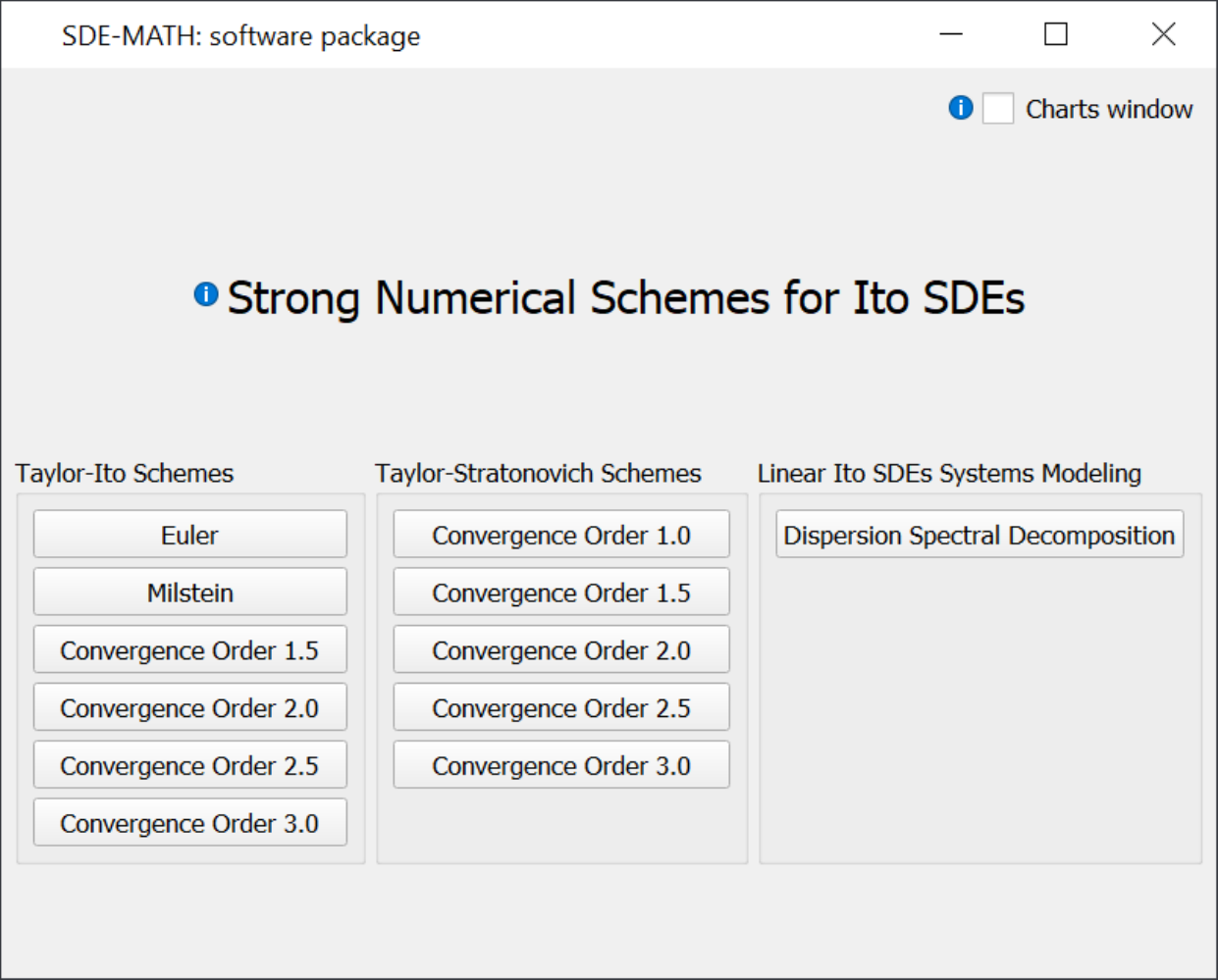}
    \caption{Main menu dialog\label{fig:main_menu}}
\end{figure}

\begin{figure}[H]
    \vspace{5mm}
    \centering
    \includegraphics[width=.7\textwidth]{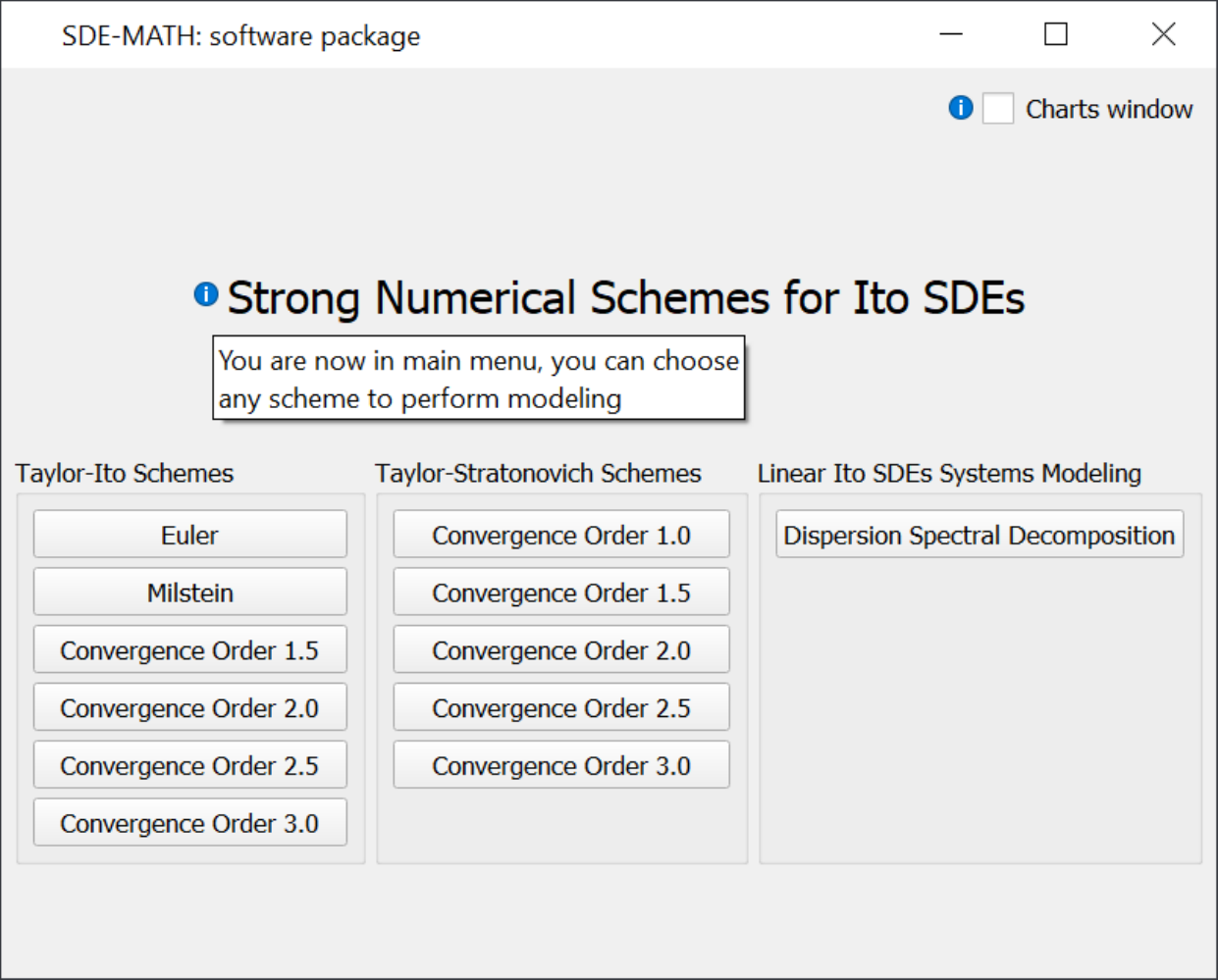}
    \caption{Tooltip\label{fig:tooltip}}
\end{figure}

\begin{figure}[H]
    \centering
    \includegraphics[width=.9\textwidth]{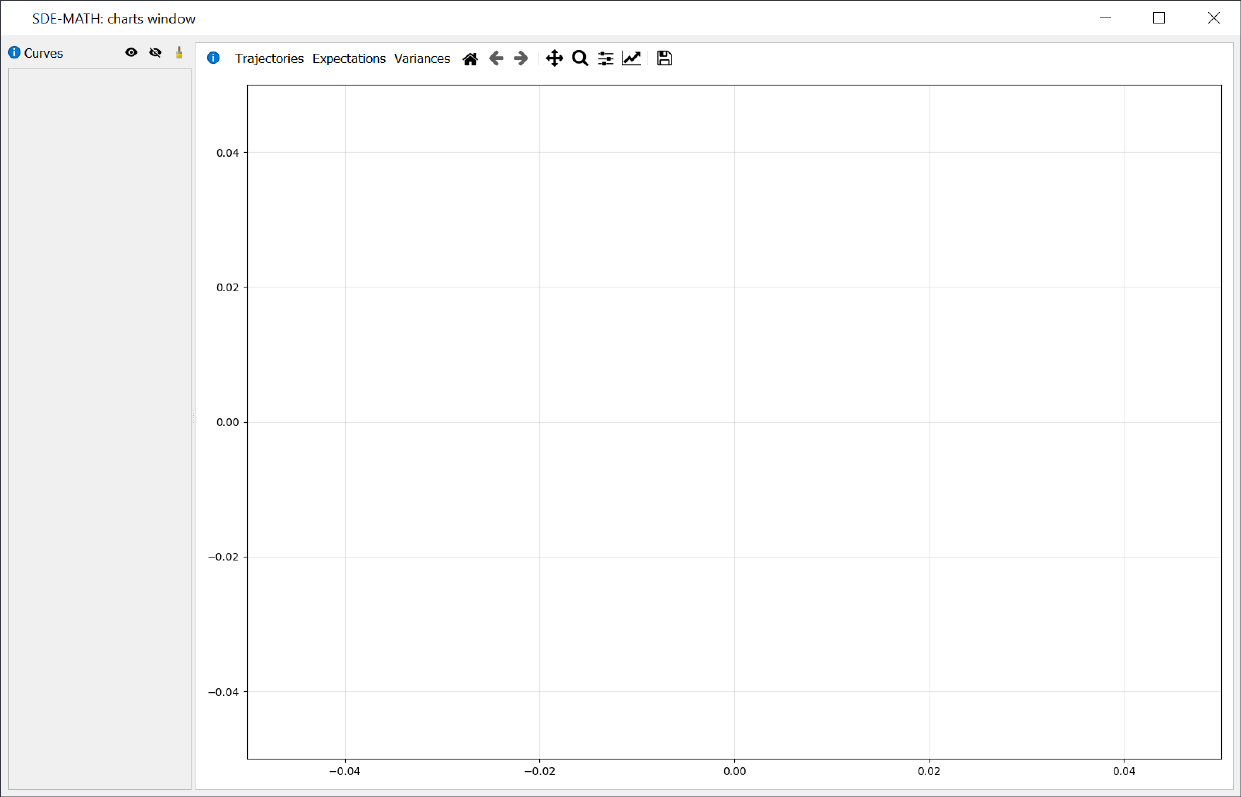}
    \caption{Charts window\label{fig:plots}}
\end{figure}

\begin{figure}[H]
    \vspace{6mm}
    \centering
    \includegraphics[width=.7\textwidth]{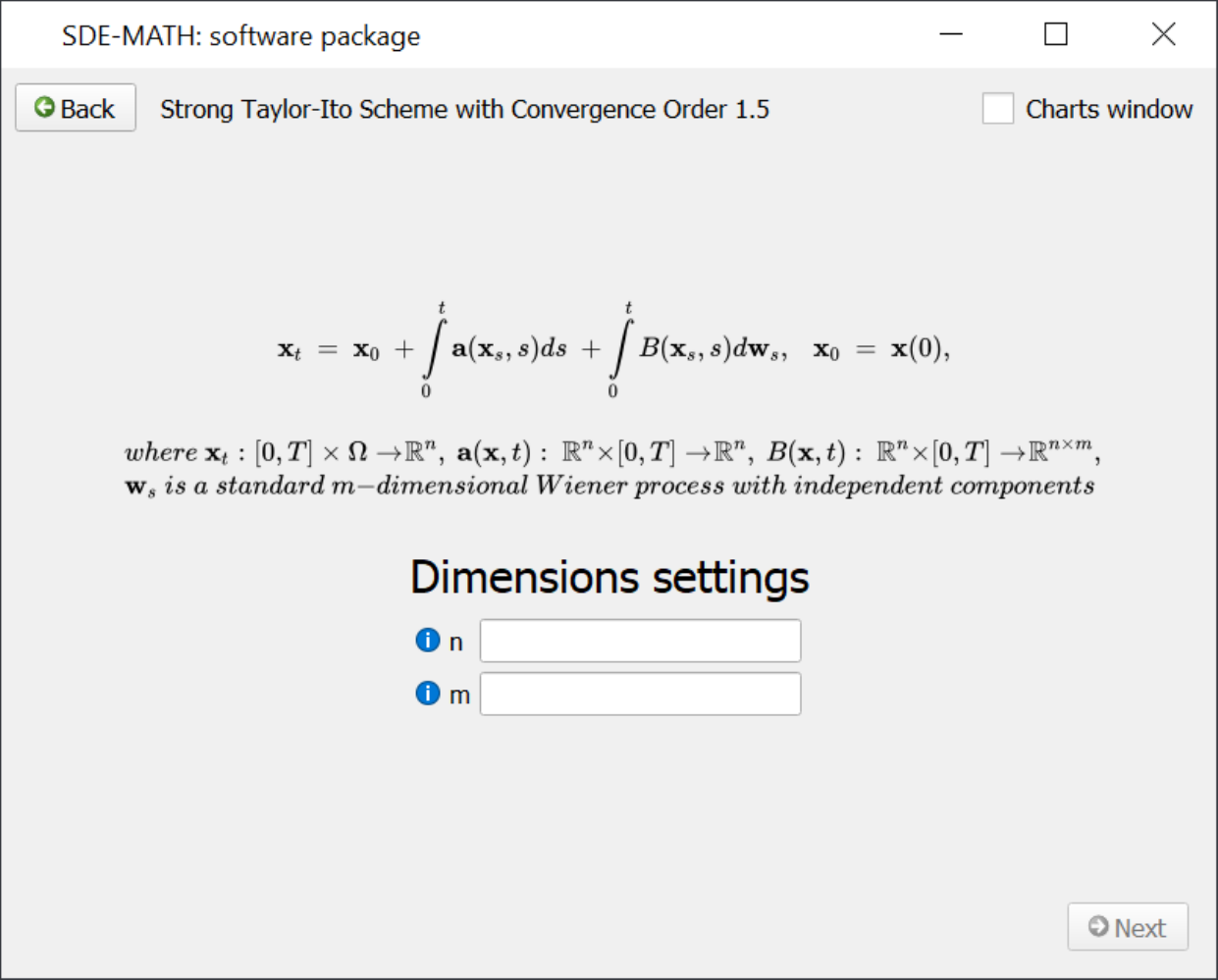}
    \caption{Nonlinear system of It\^o SDEs data input\label{fig:nonlinear_1}}
\end{figure}

\begin{figure}[H]
    \centering
    \includegraphics[width=.7\textwidth]{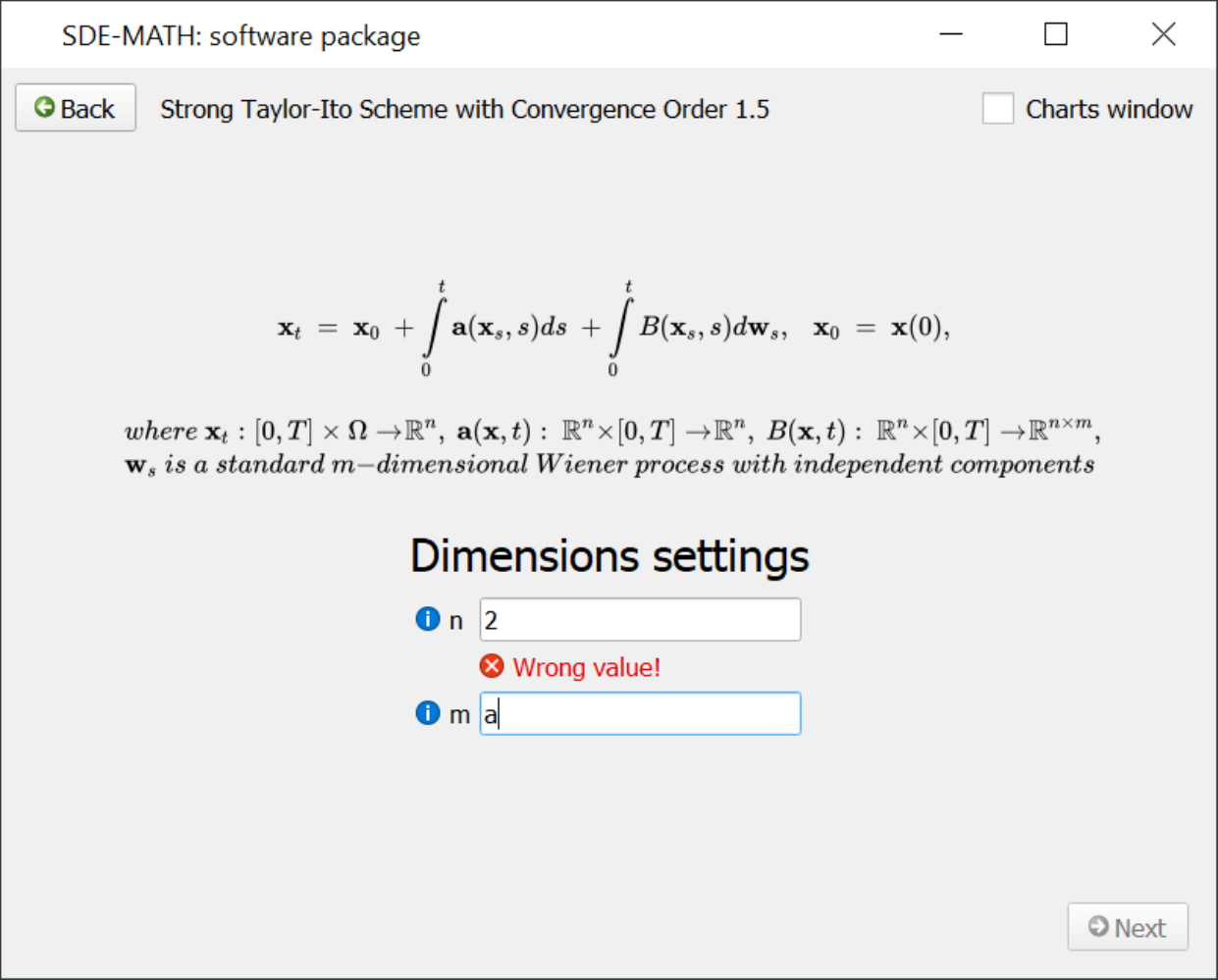}
    \caption{Wrong data input\label{fig:nonlinear_2}}
\end{figure}

\begin{figure}[H]
    \vspace{6mm}
    \centering
    \includegraphics[width=.7\textwidth]{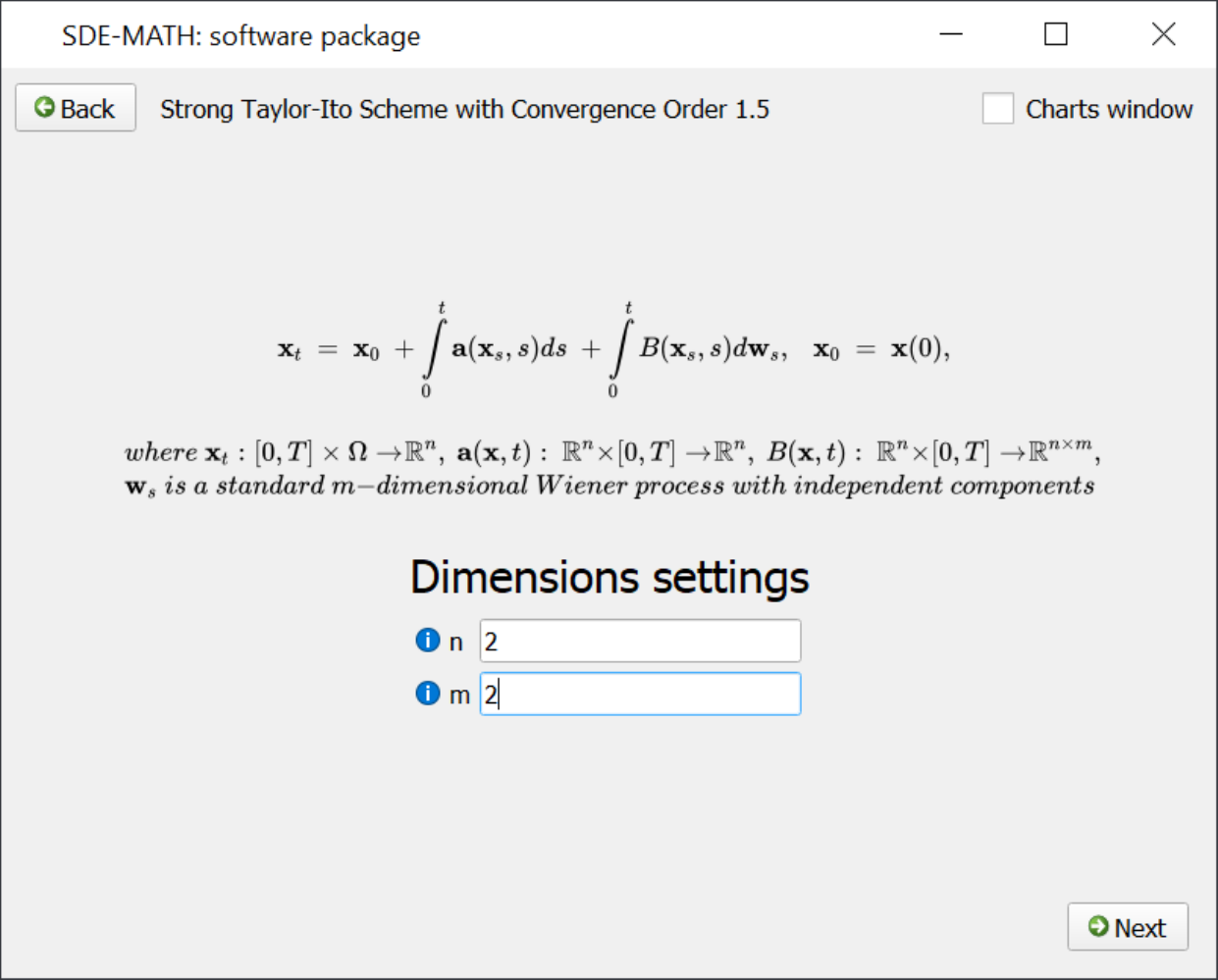}
    \caption{Correct data input\label{fig:nonlinear_3}}
\end{figure}

\begin{figure}[H]
    \centering
    \includegraphics[width=.7\textwidth]{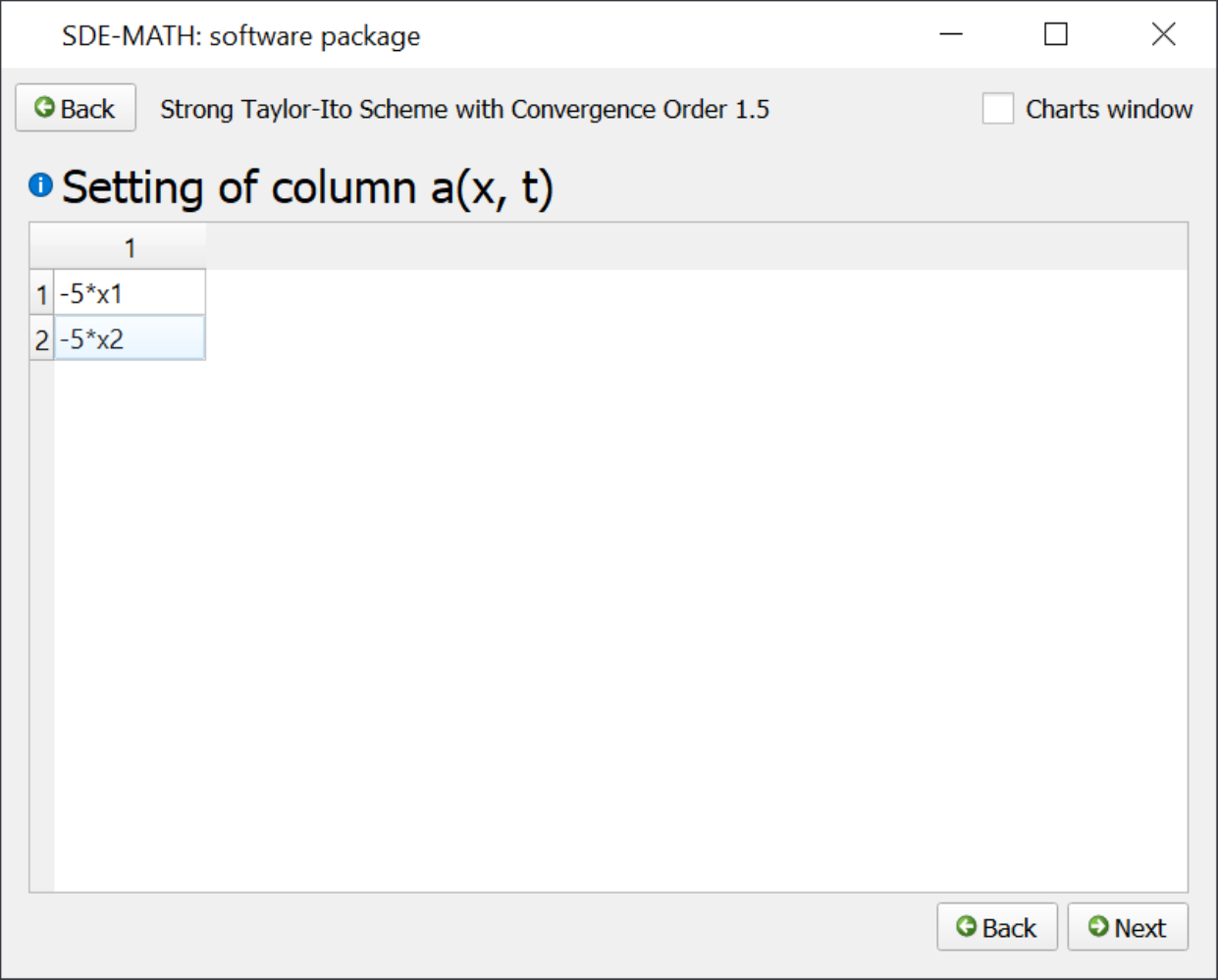}
    \caption{Vector function ${\bf a}({\bf x},t)$ input\label{fig:nonlinear_4}}
\end{figure}

\begin{figure}[H]
    \vspace{6mm}
    \centering
    \includegraphics[width=.7\textwidth]{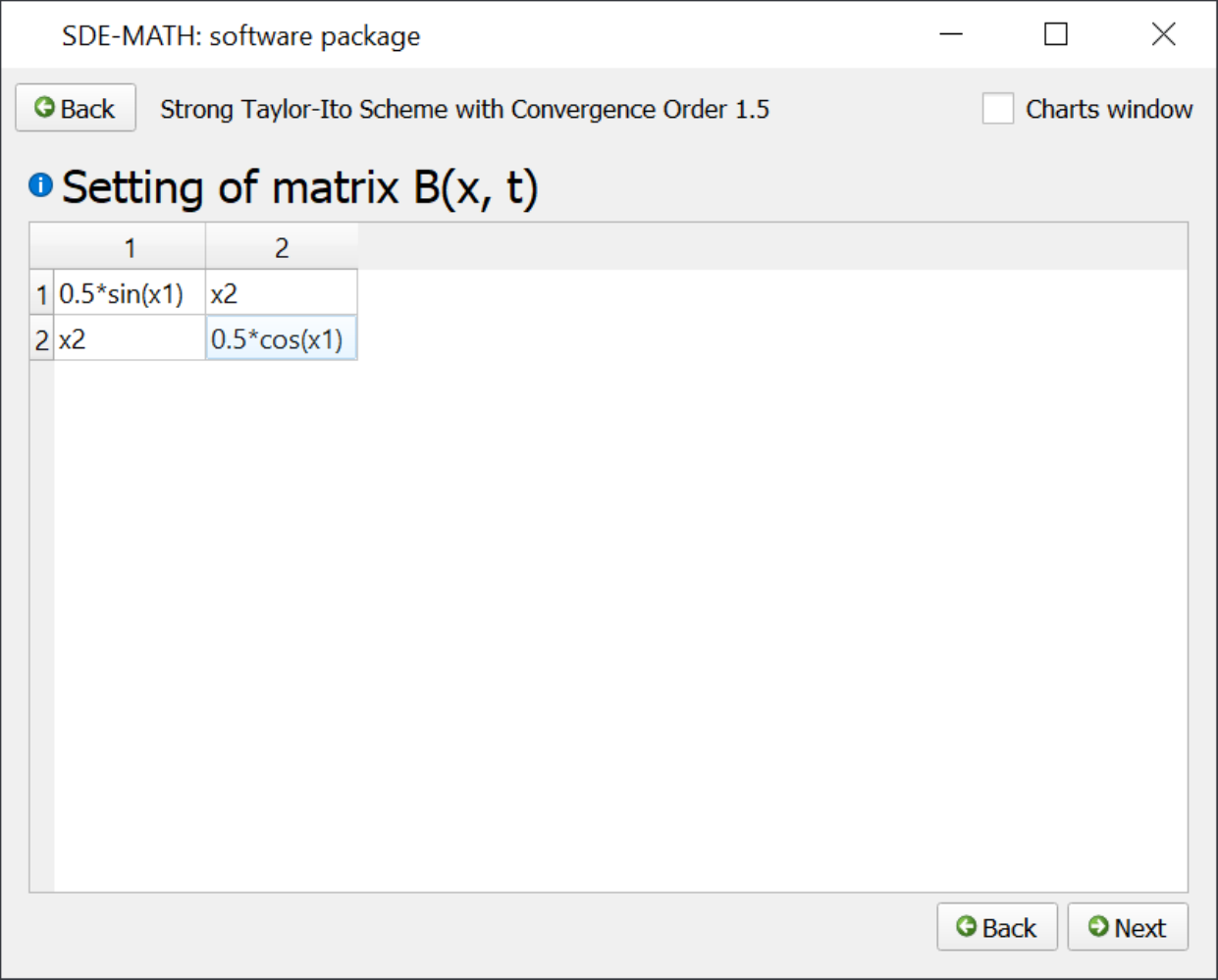}
    \caption{Matrix function $B({\bf x},t)$ input\label{fig:nonlinear_5}}
\end{figure}

\begin{figure}[H]
    \centering
    \includegraphics[width=.7\textwidth]{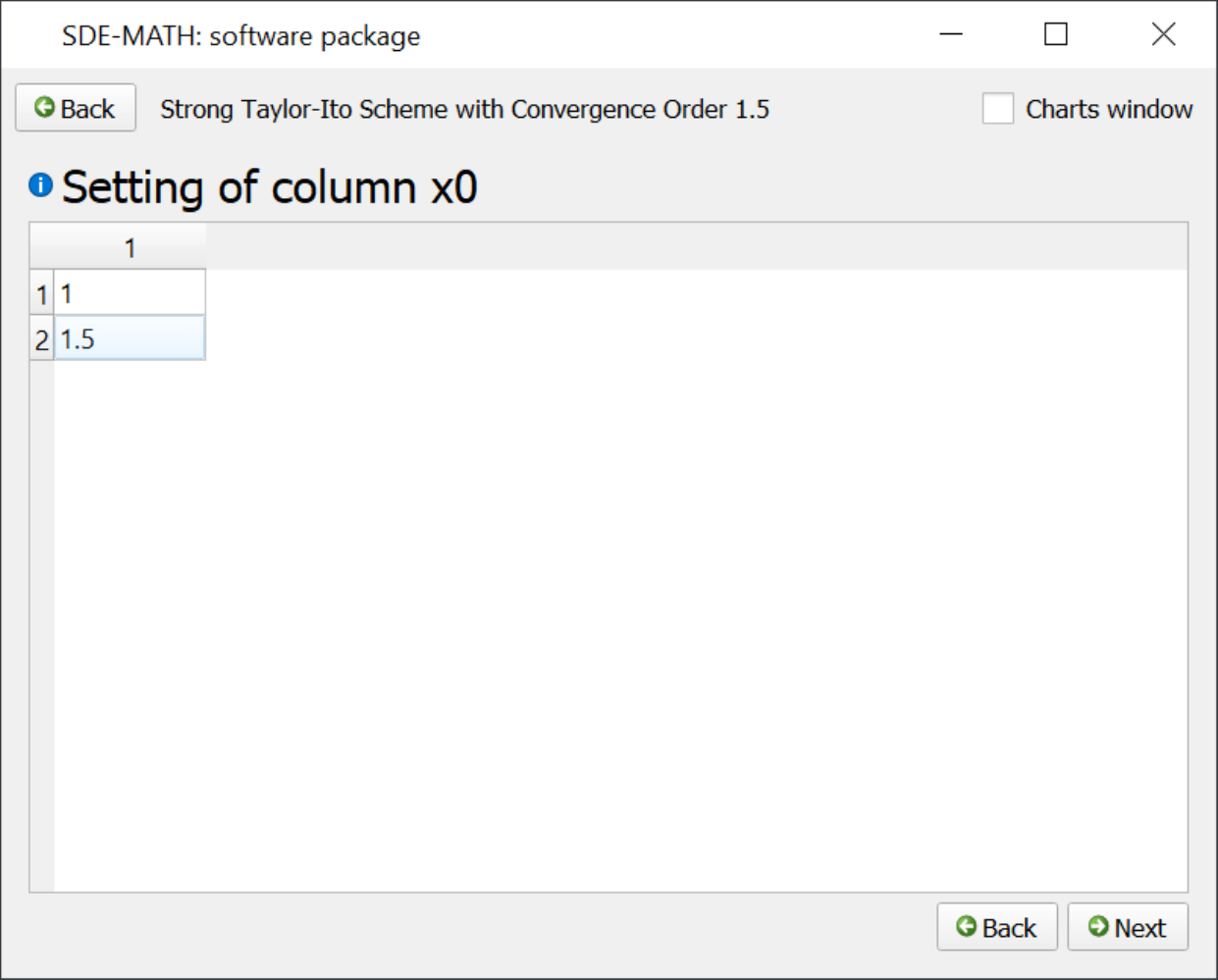}
    \caption{Initial data input\label{fig:nonlinear_6}}
\end{figure}

\begin{figure}[H]
    \vspace{6mm}
    \centering
    \includegraphics[width=.7\textwidth]{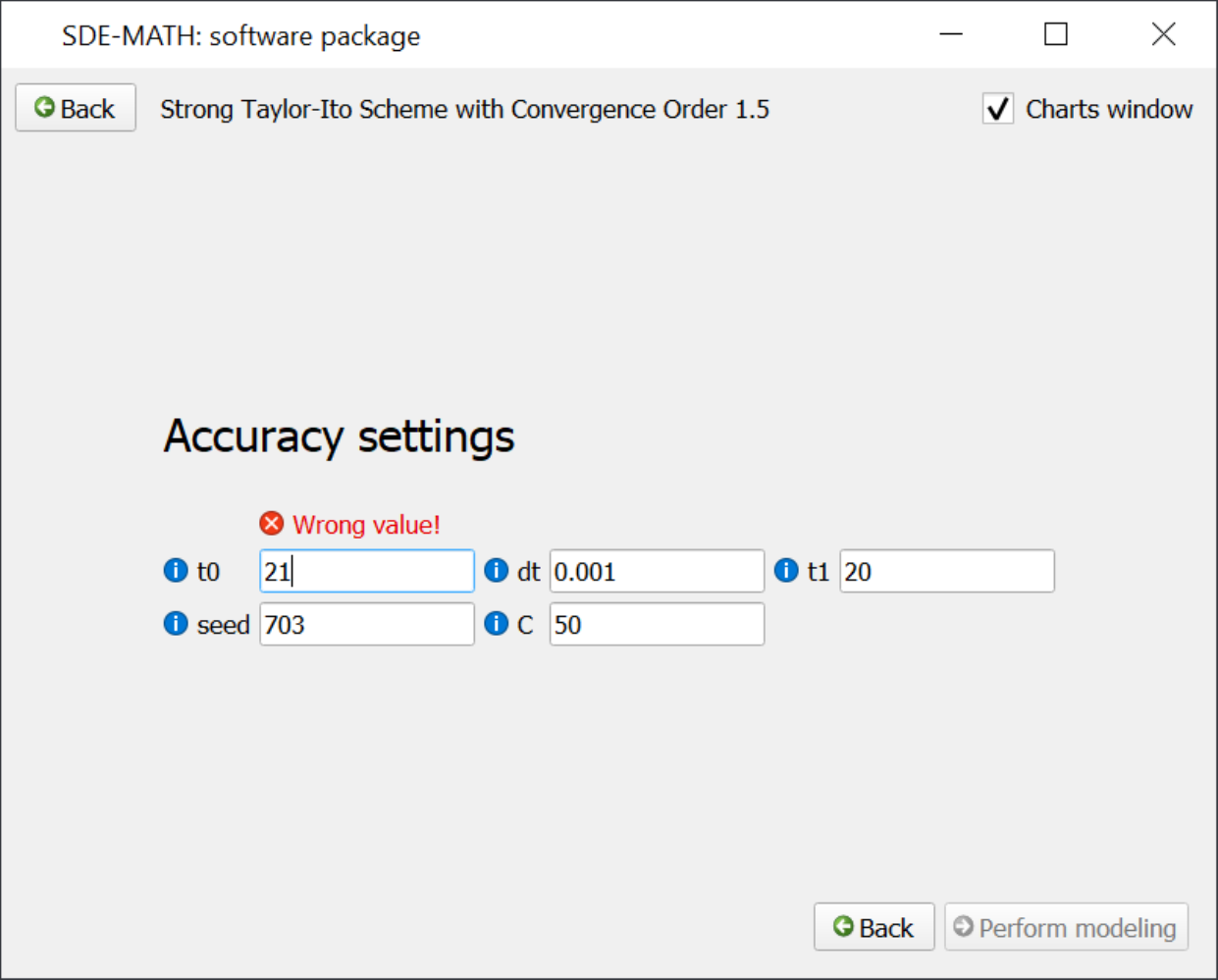}
    \caption{Wrong data input\label{fig:nonlinear_7}}
\end{figure}

\begin{figure}[H]
    \centering
    \includegraphics[width=.7\textwidth]{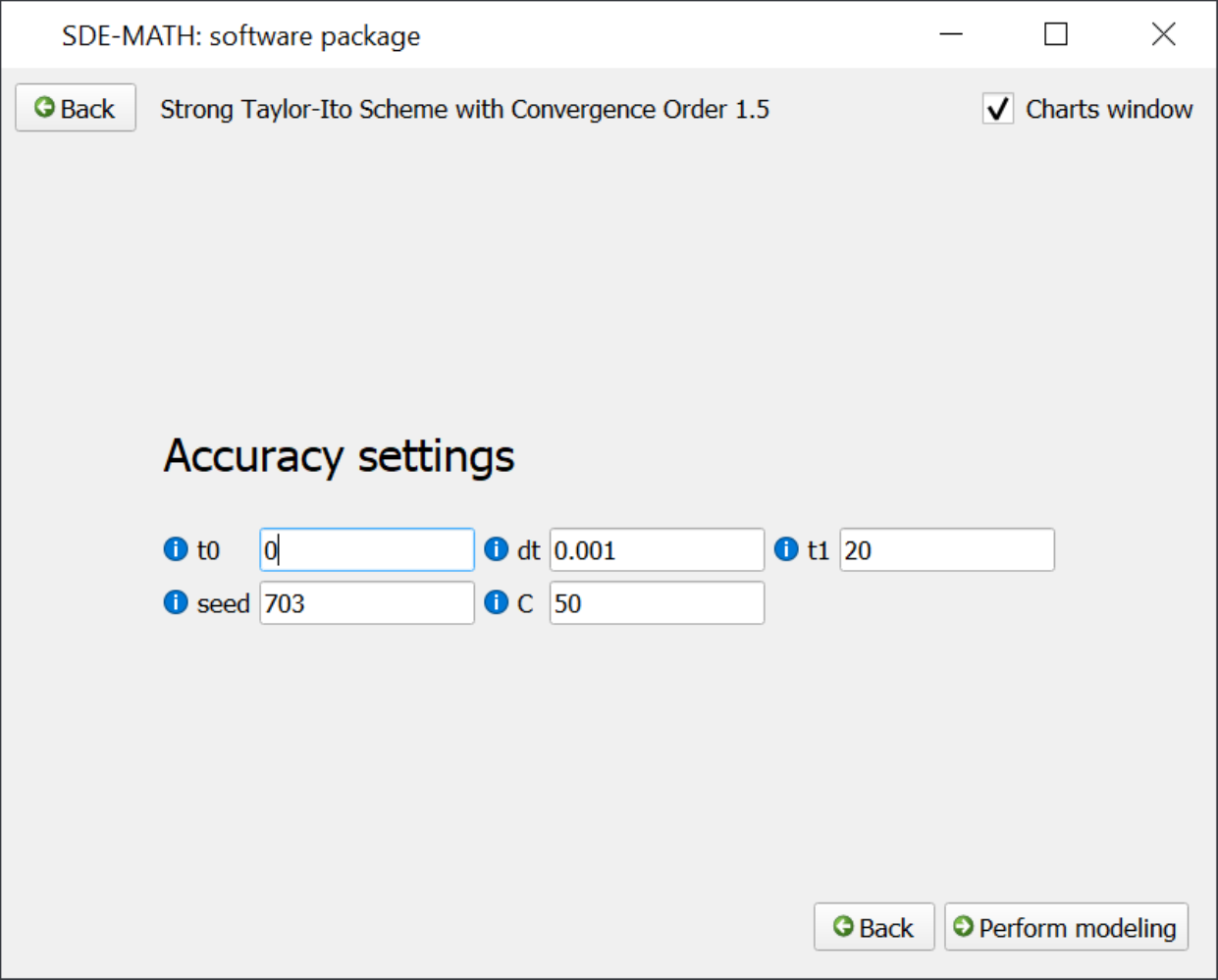}
    \caption{Correct data input\label{fig:nonlinear_8}}
\end{figure}

\begin{figure}[H]
    \vspace{6mm}
    \centering
    \includegraphics[width=.7\textwidth]{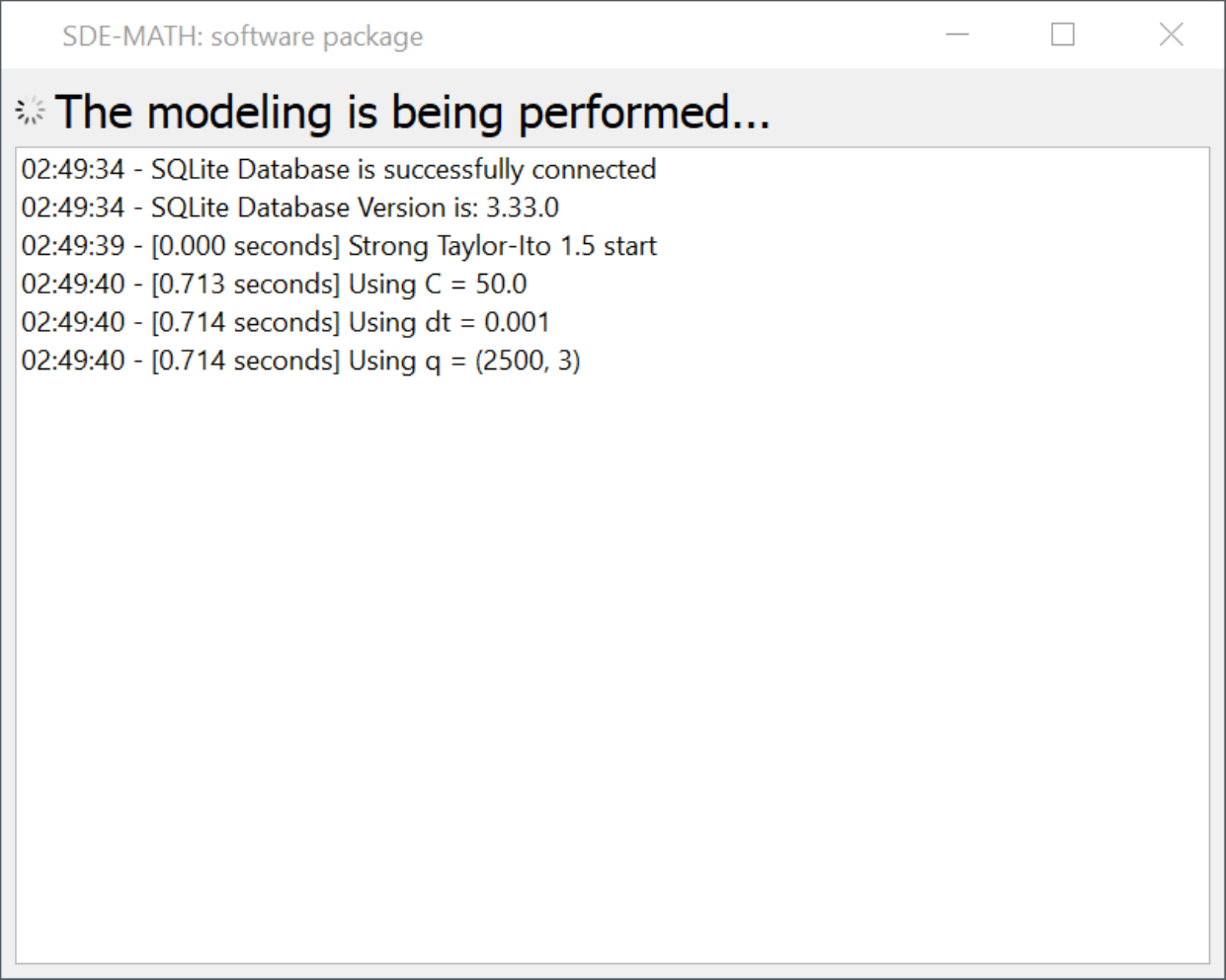}
    \caption{Modeling logs\label{fig:nonlinear_9}}
\end{figure}

\begin{figure}[H]
    \centering
    \includegraphics[width=.7\textwidth]{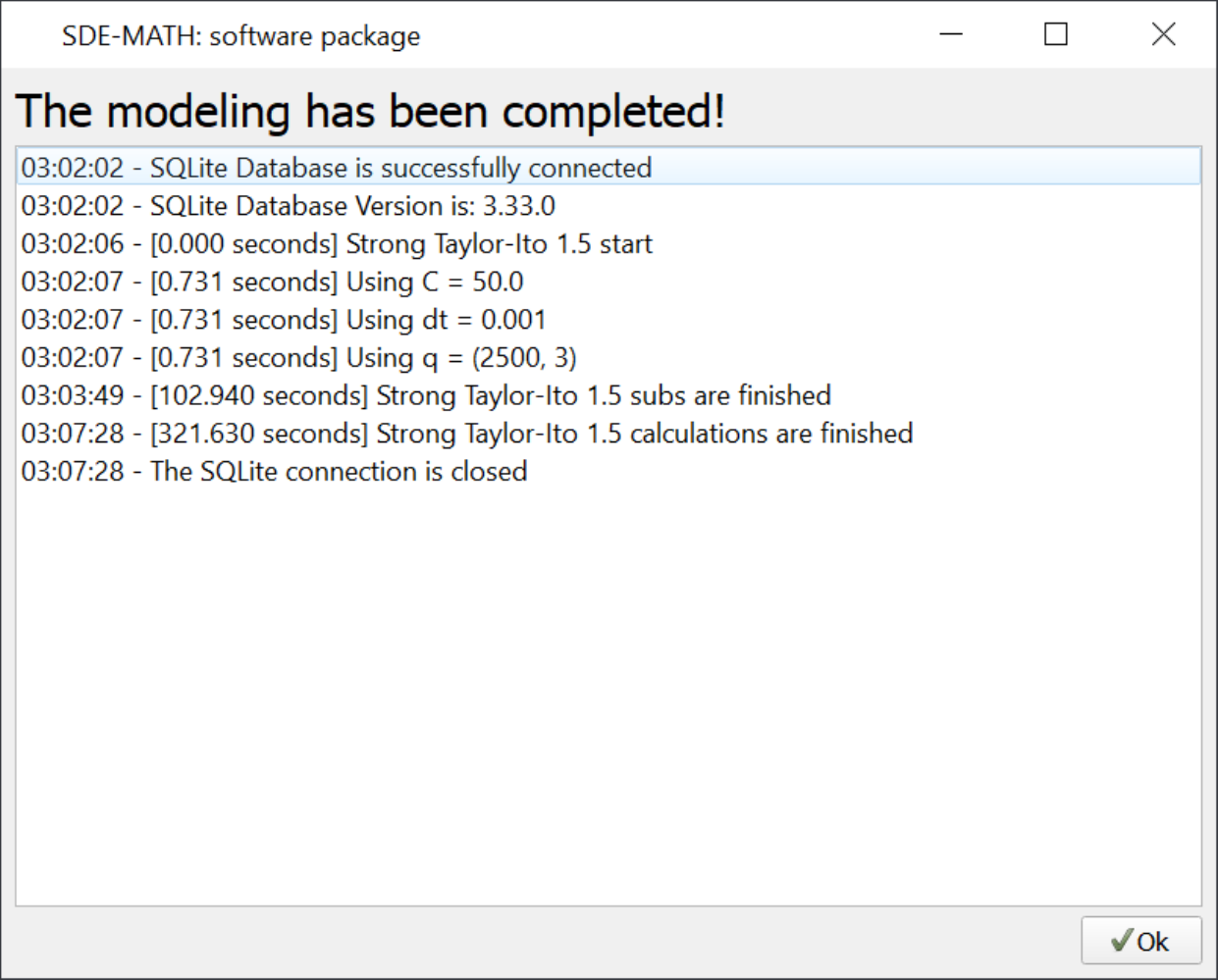}
    \caption{Modeling logs\label{fig:nonlinear_10}}
\end{figure}

\begin{figure}[H]
    \vspace{5mm}
    \centering
    \includegraphics[width=.9\textwidth]{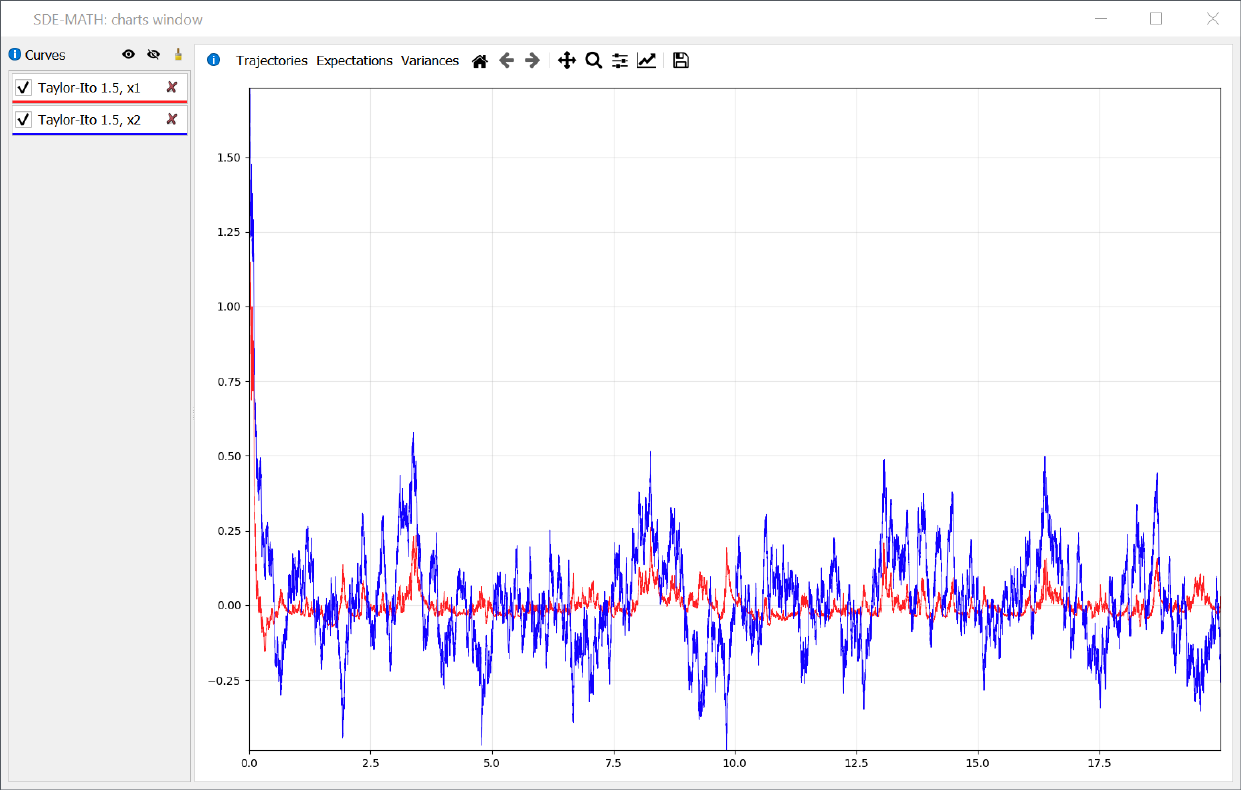}
    \caption{Modeling results\label{fig:nonlinear_11}}
\end{figure}

\begin{figure}[H]
    \centering
    \includegraphics[width=.7\textwidth]{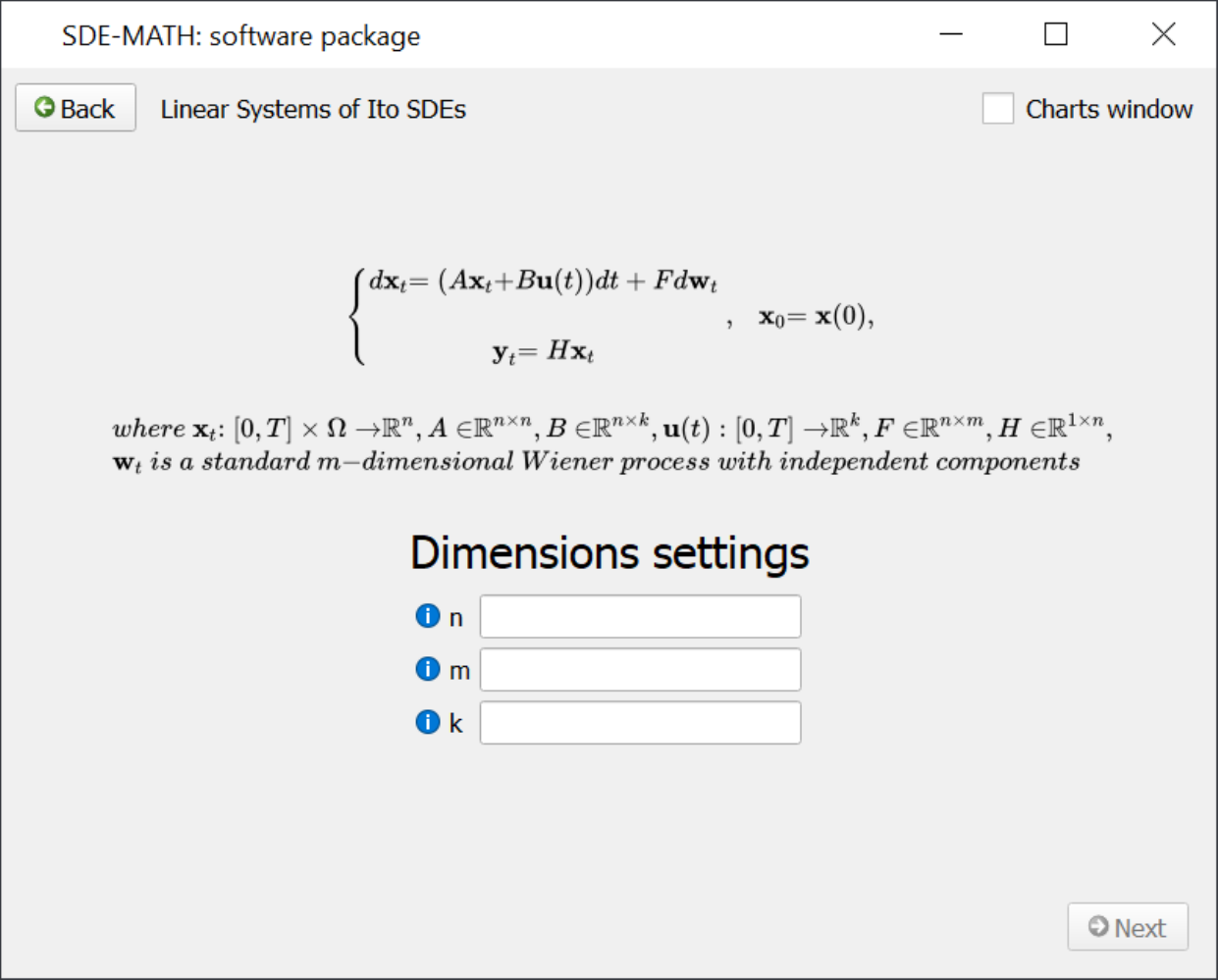}
    \caption{Linear system of It\^o SDEs data input\label{fig:linear_1}}
\end{figure}

\begin{figure}[H]
    \vspace{6mm}
    \centering
    \includegraphics[width=.7\textwidth]{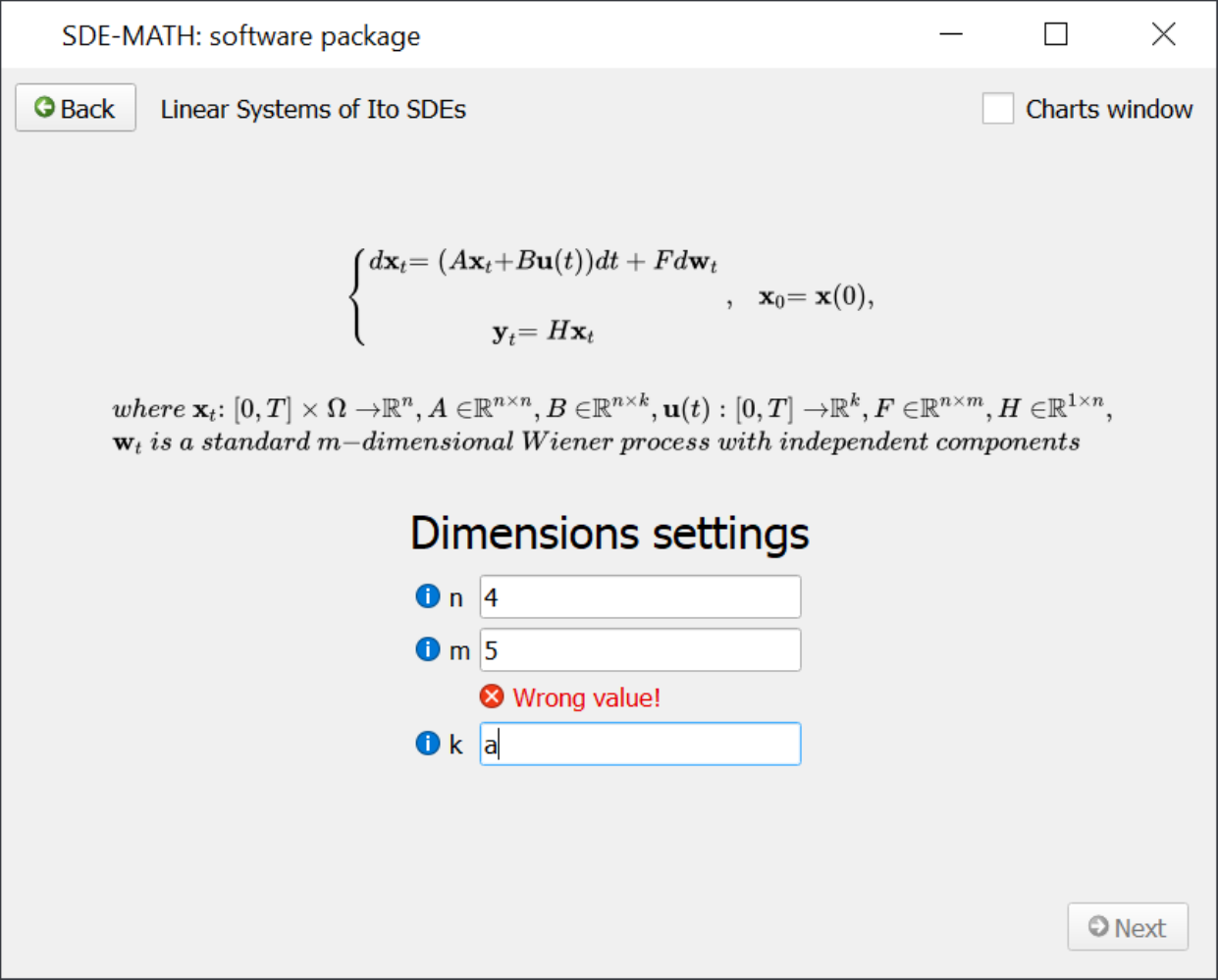}
    \caption{Wrong data input\label{fig:linear_2}}
\end{figure}

\begin{figure}[H]
    \centering
    \includegraphics[width=.7\textwidth]{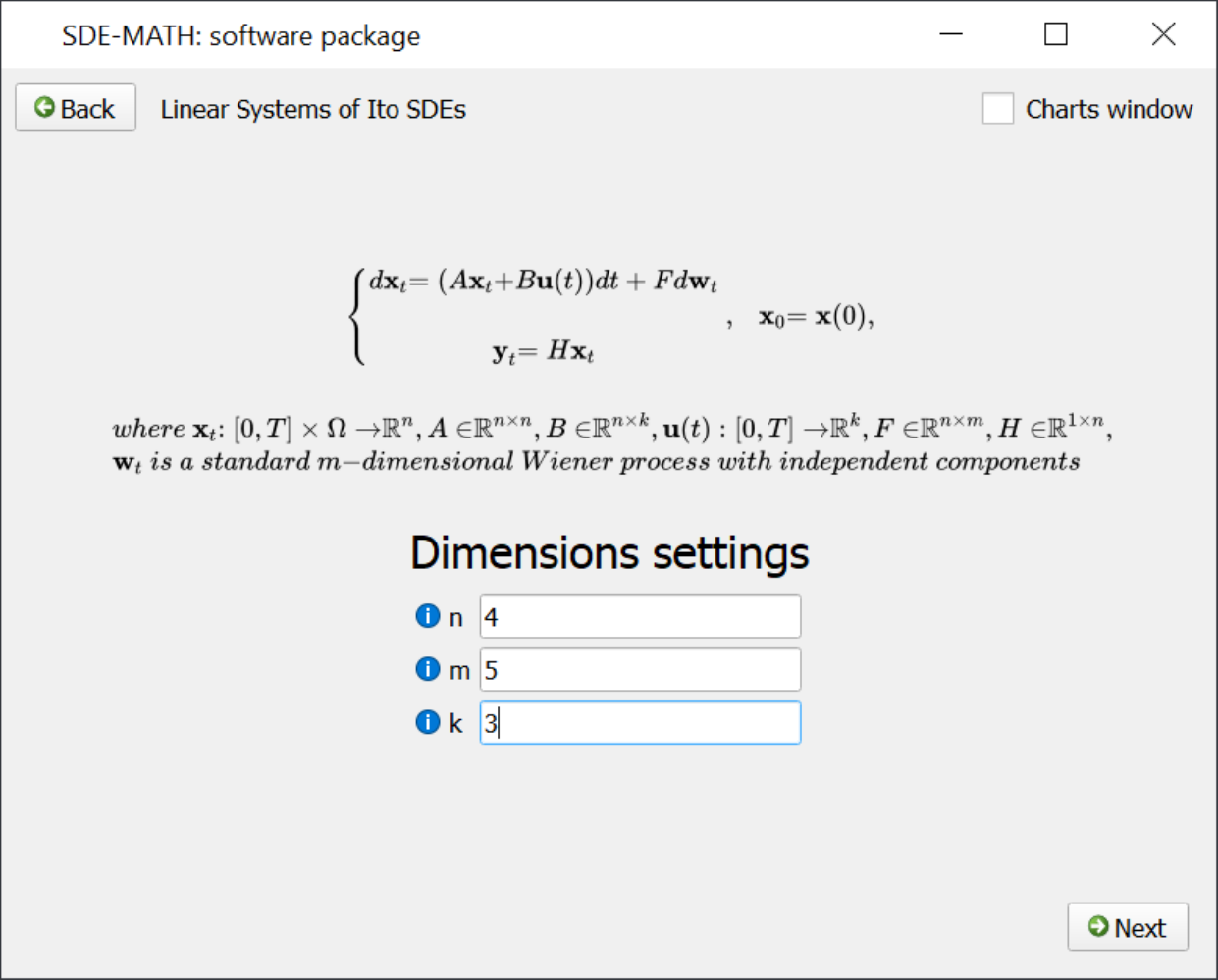}
    \caption{Correct data input\label{fig:linear_3}}
\end{figure}

\begin{figure}[H]
    \vspace{6mm}
    \centering
    \includegraphics[width=.7\textwidth]{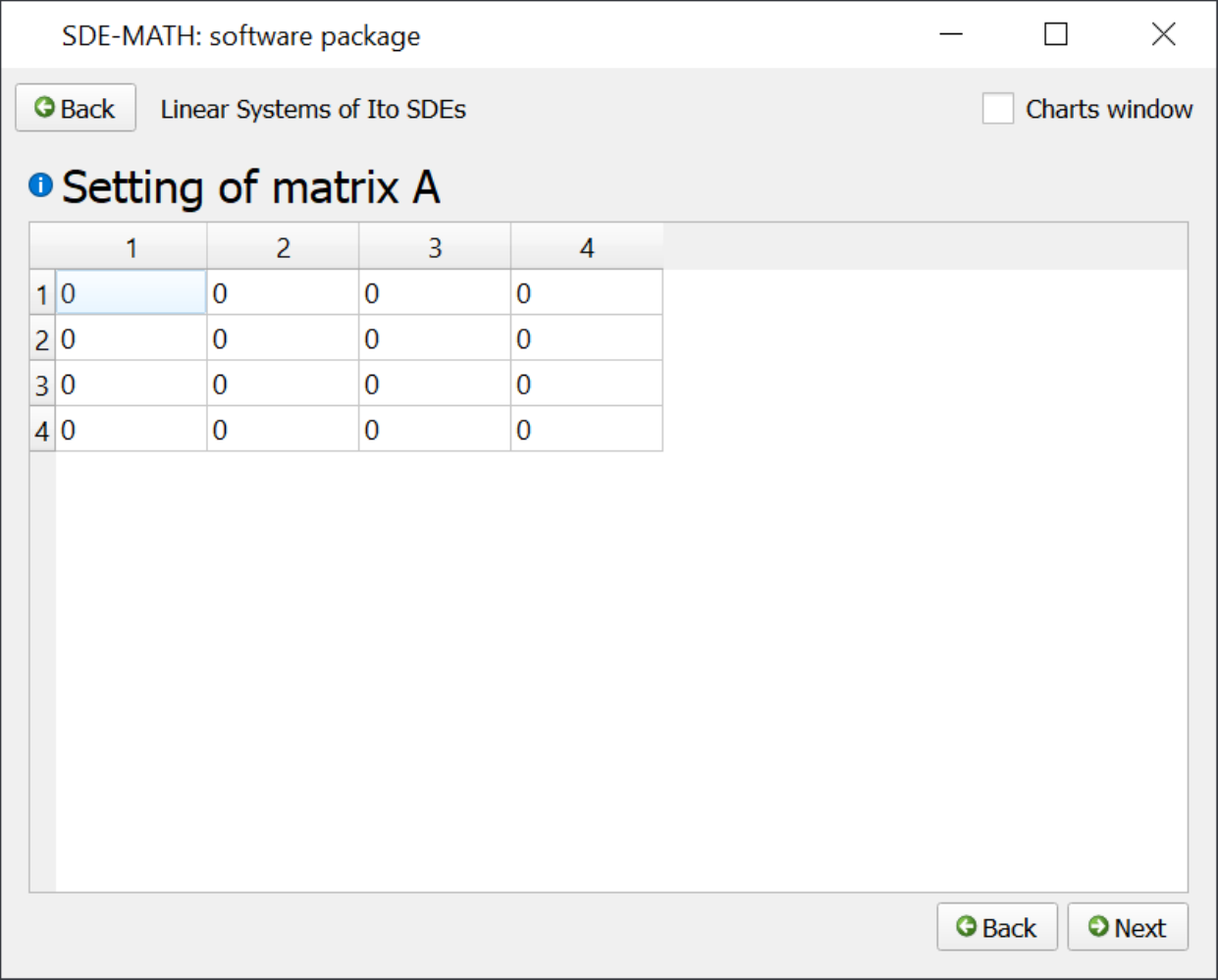}
    \caption{Matrix $A$ input\label{fig:linear_4}}
\end{figure}

\begin{figure}[H]
    \centering
    \includegraphics[width=.7\textwidth]{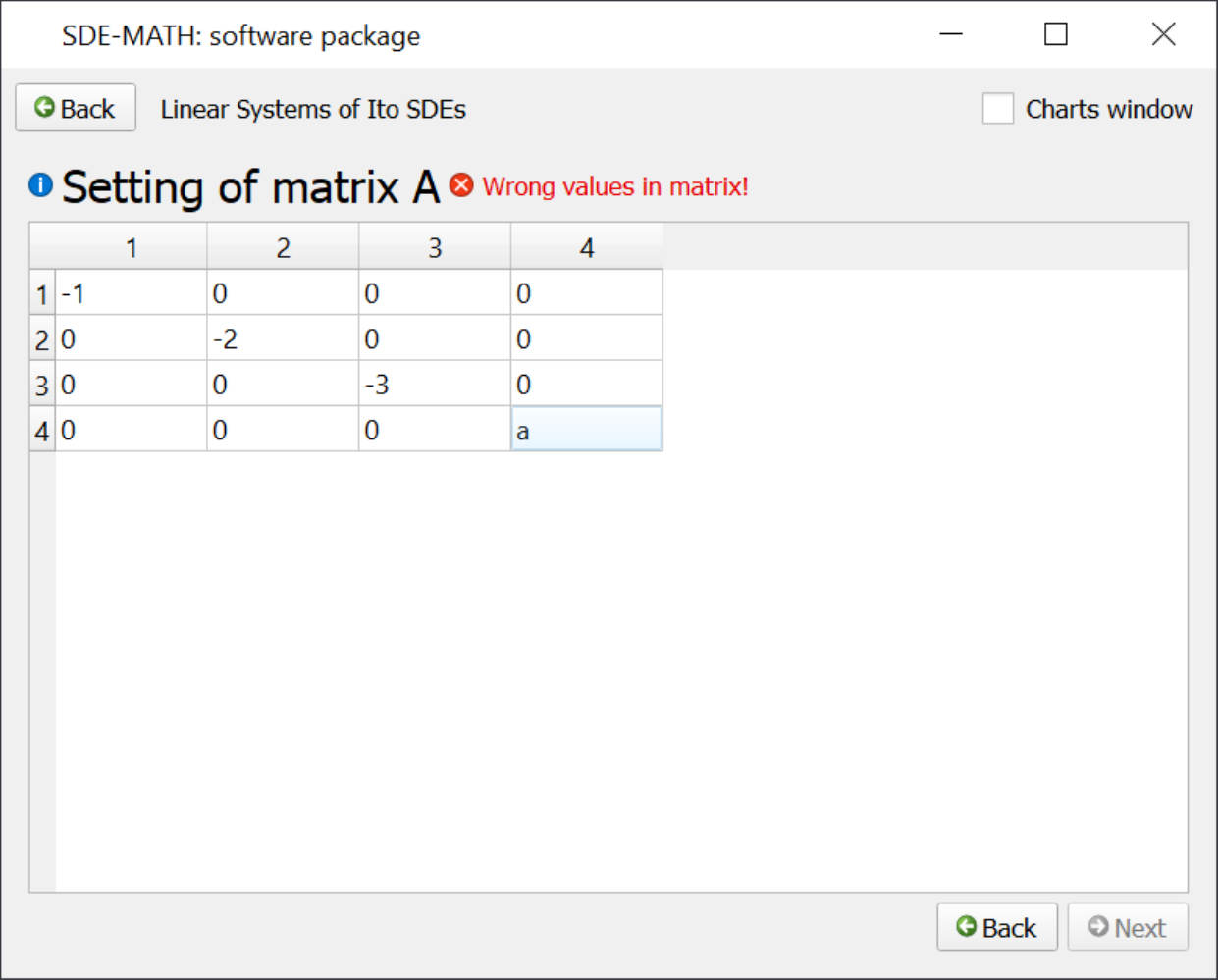}
    \caption{Wrong matrix $A$ input\label{fig:linear_5}}
\end{figure}

\begin{figure}[H]
    \vspace{6mm}
    \centering
    \includegraphics[width=.7\textwidth]{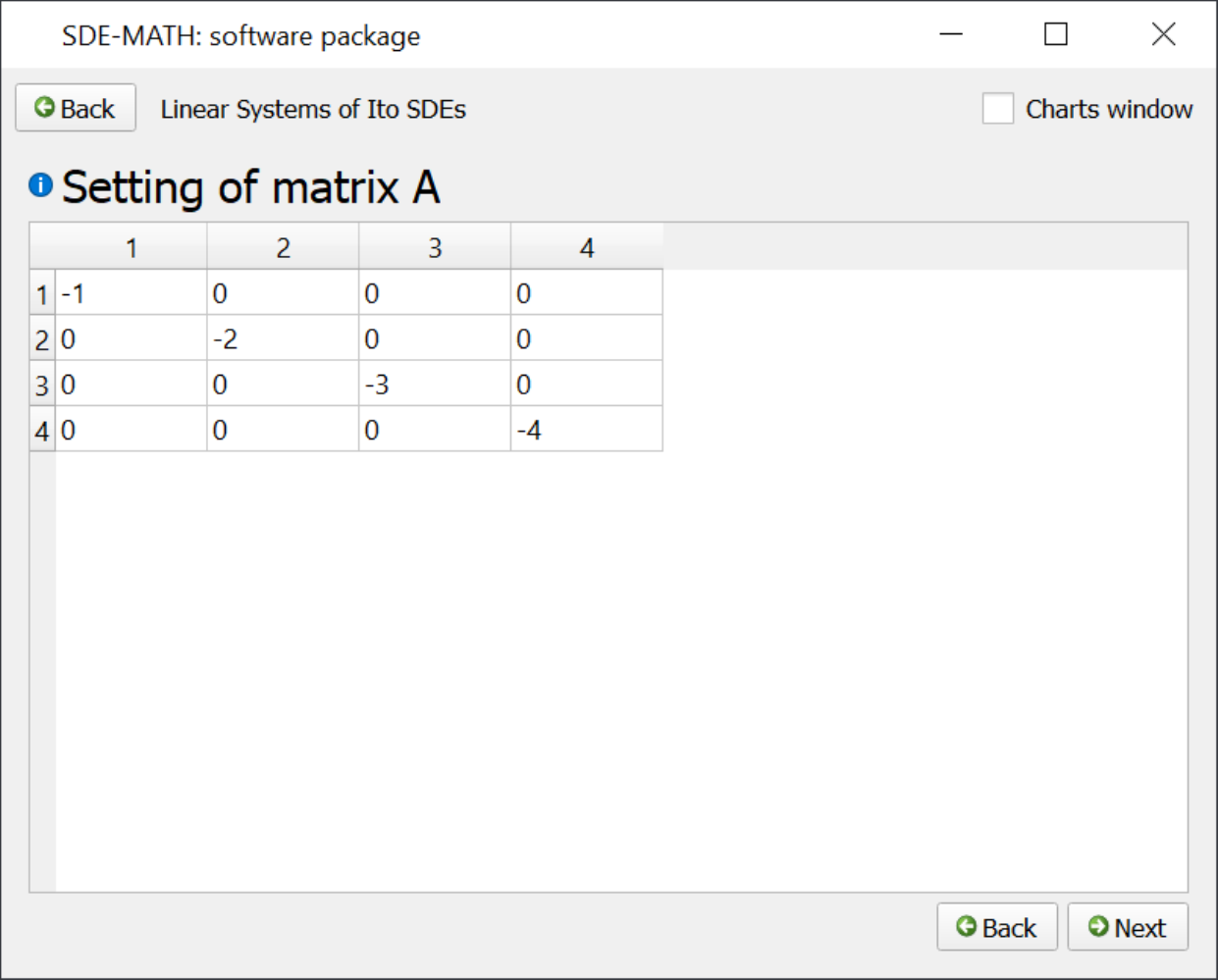}
    \caption{Correct matrix $A$ input\label{fig:linear_6}}
\end{figure}

\begin{figure}[H]
    \centering
    \includegraphics[width=.7\textwidth]{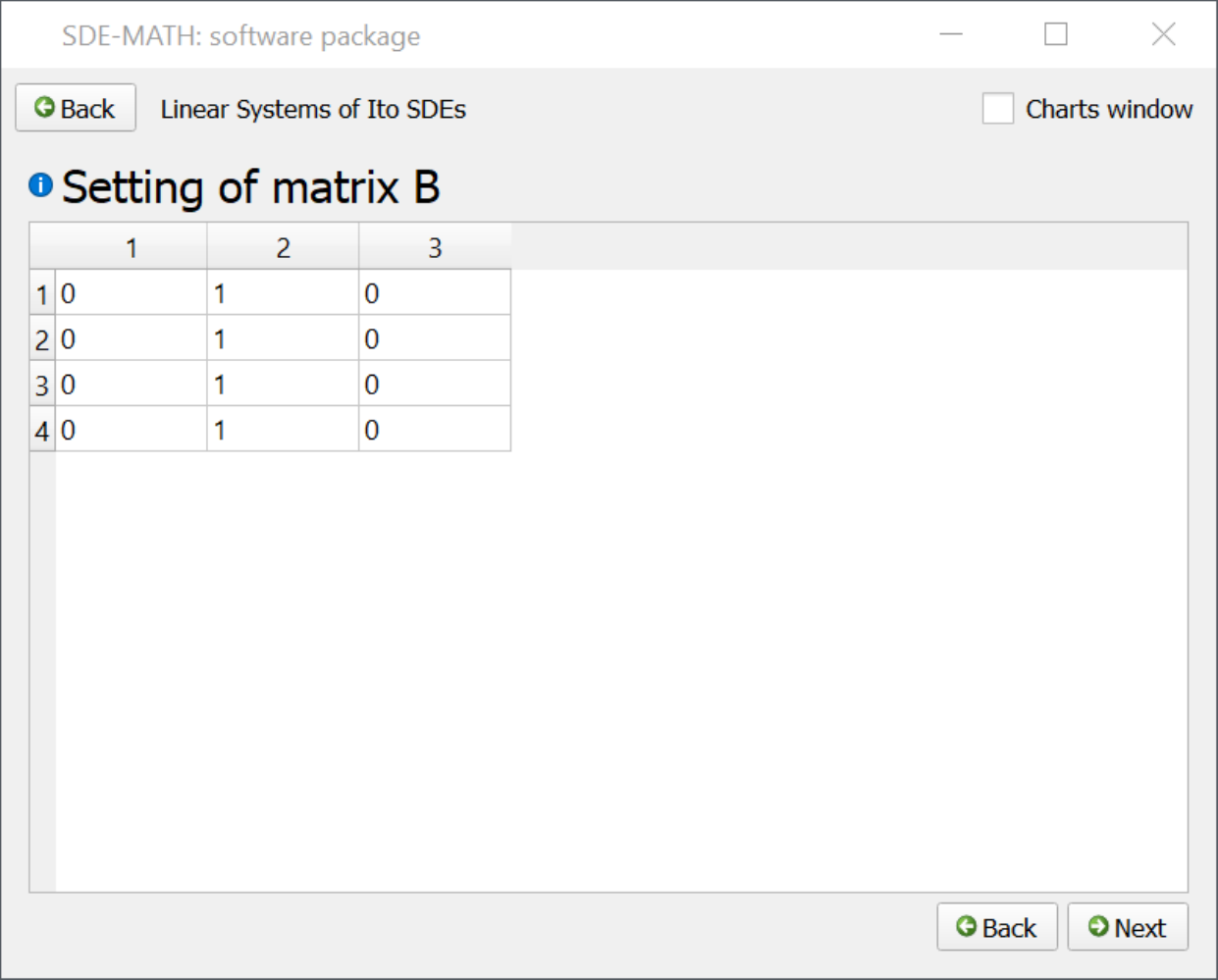}
    \caption{Matrix $B$ input\label{fig:linear_7}}
\end{figure}

\begin{figure}[H]
    \vspace{6mm}
    \centering
    \includegraphics[width=.7\textwidth]{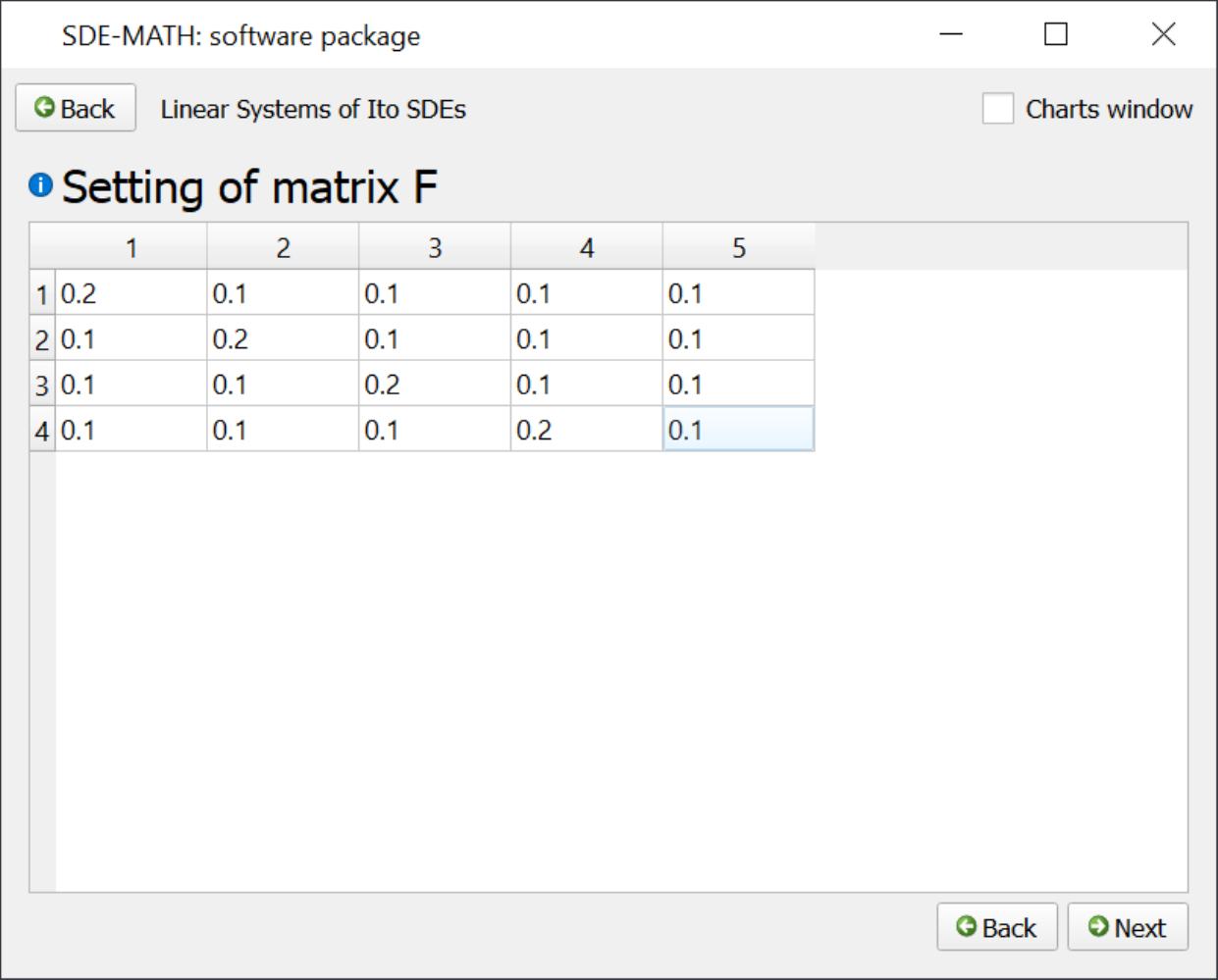}
    \caption{Matrix $F$ input\label{fig:linear_8}}
\end{figure}

\begin{figure}[H]
    \centering
    \includegraphics[width=.7\textwidth]{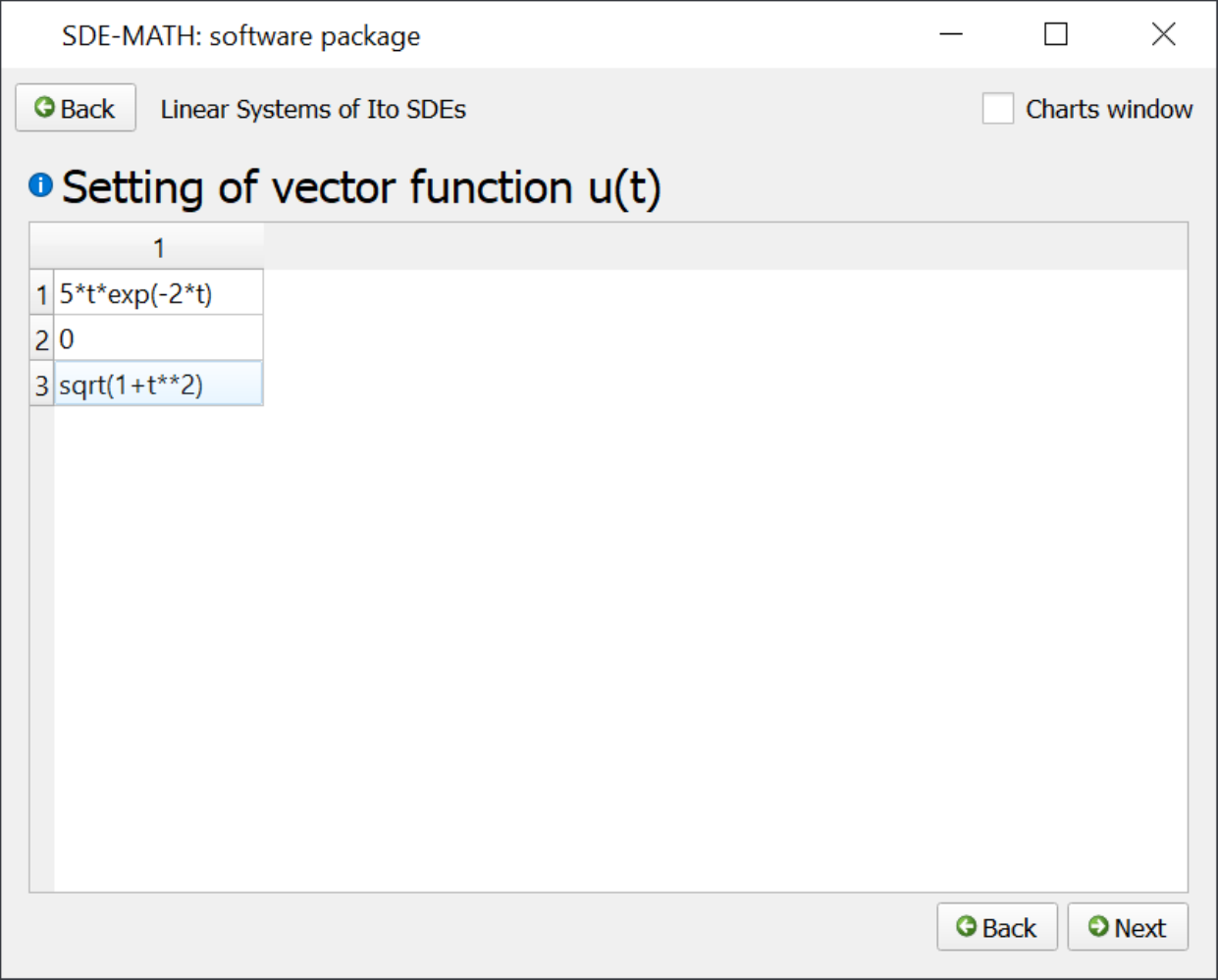}
    \caption{Vector function ${\bf u}(t)$ input\label{fig:linear_9}}
\end{figure}

\begin{figure}[H]
    \vspace{6mm}
    \centering
    \includegraphics[width=.7\textwidth]{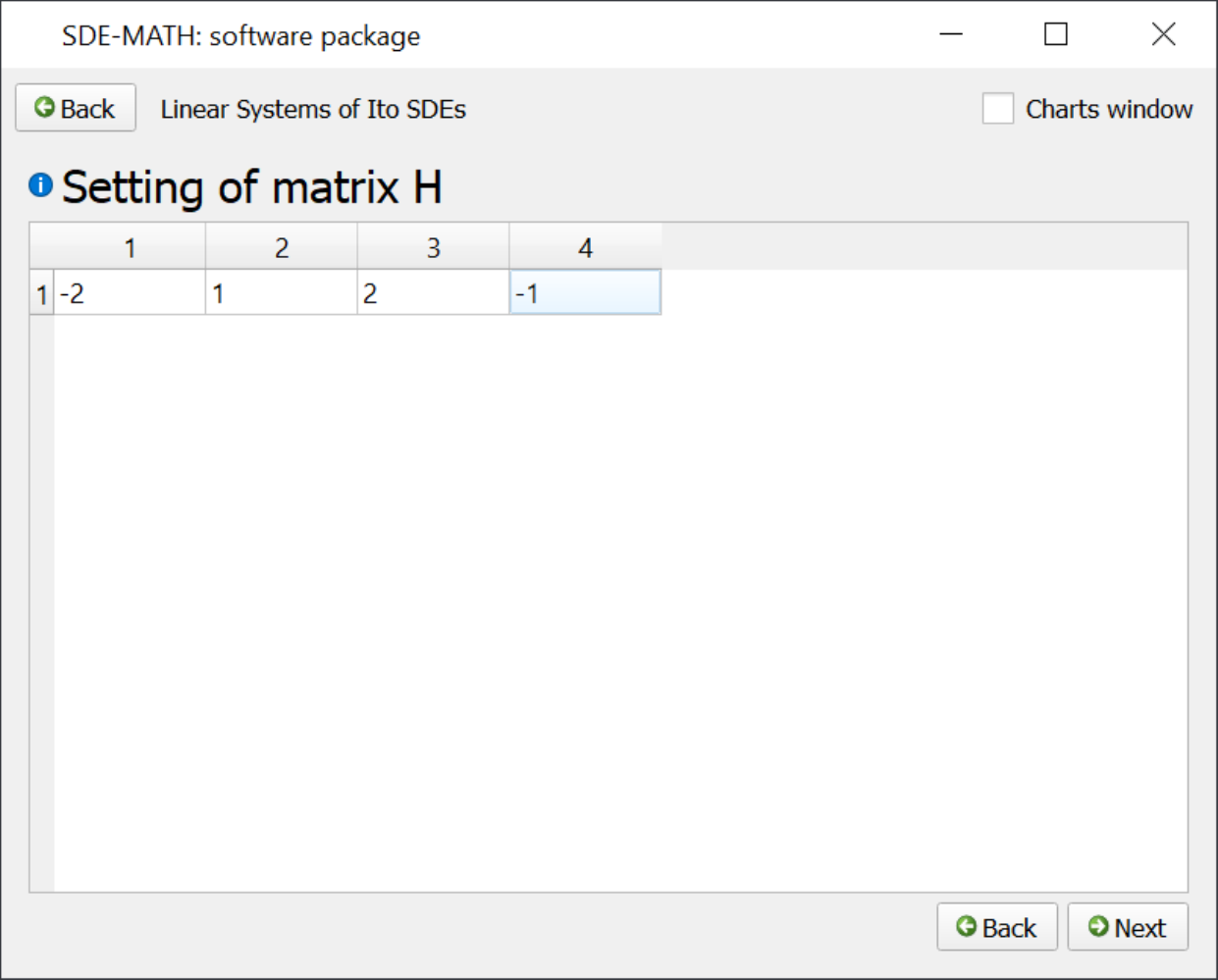}
    \caption{Matrix $H$ input\label{fig:linear_10}}
\end{figure}

\begin{figure}[H]
    \centering
    \includegraphics[width=.7\textwidth]{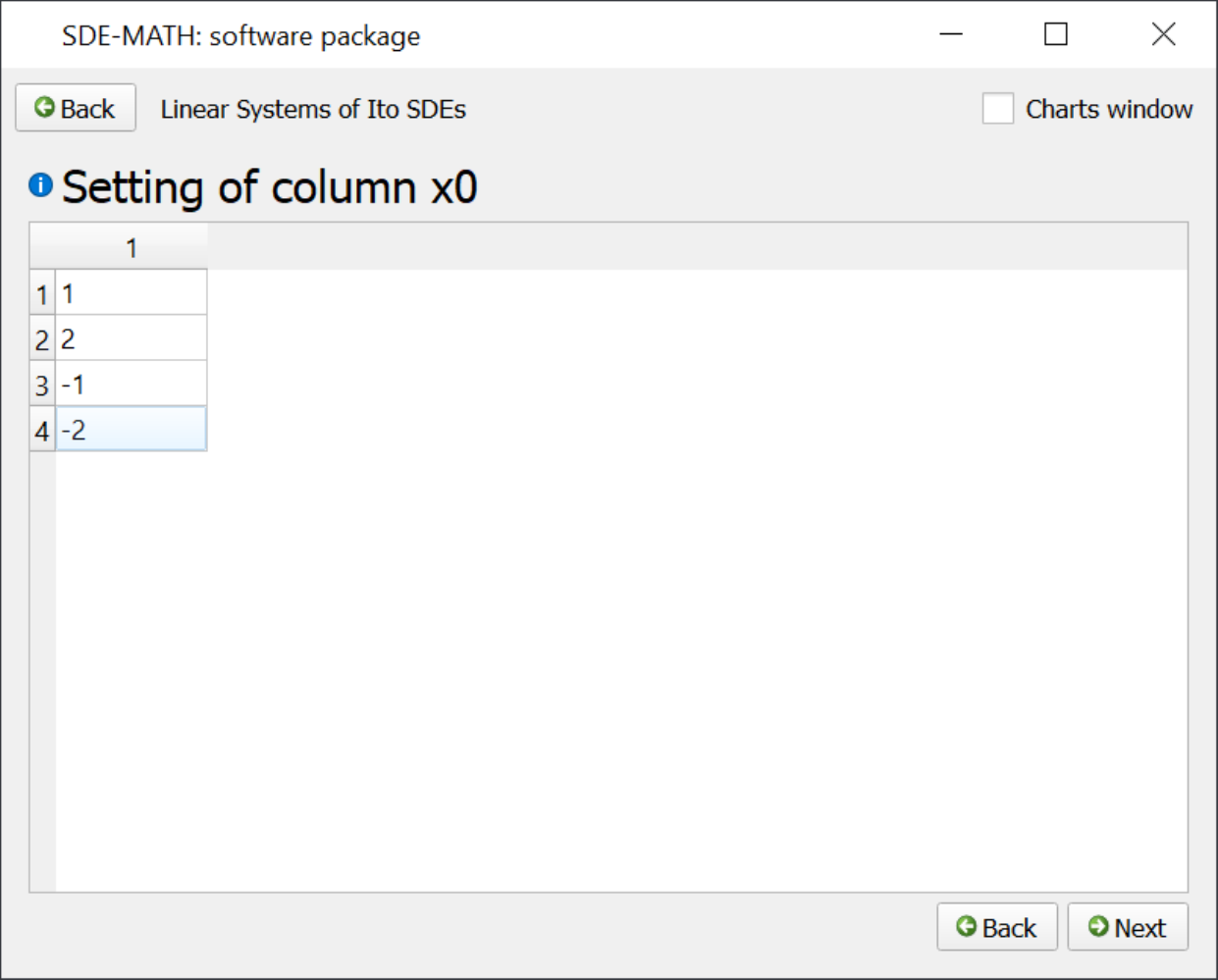}
    \caption{Initial data input\label{fig:linear_11}}
\end{figure}

\begin{figure}[H]
    \vspace{6mm}
    \centering
    \includegraphics[width=.7\textwidth]{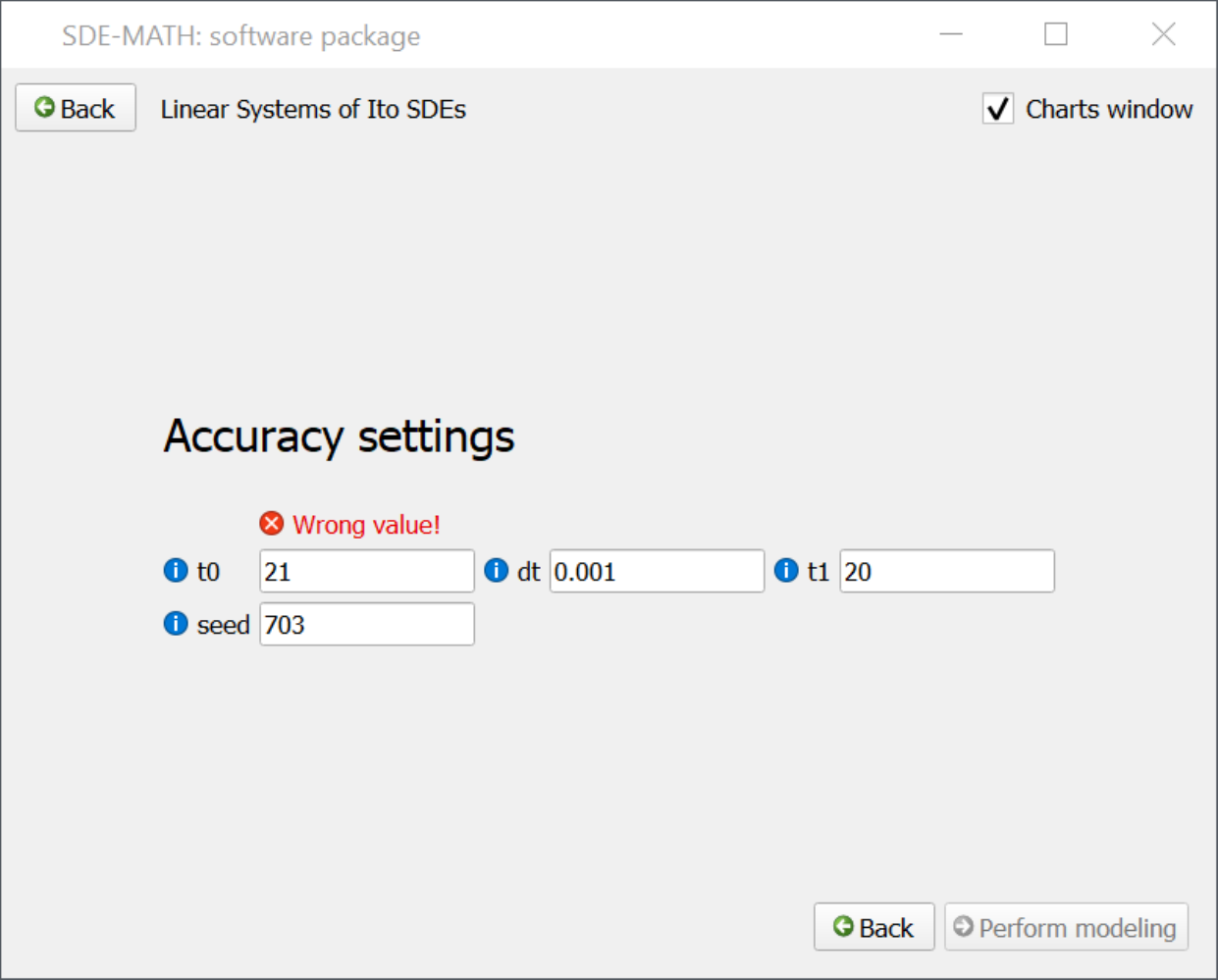}
    \caption{Wrong data input\label{fig:linear_12}}
\end{figure}

\begin{figure}[H]
    \centering
    \includegraphics[width=.7\textwidth]{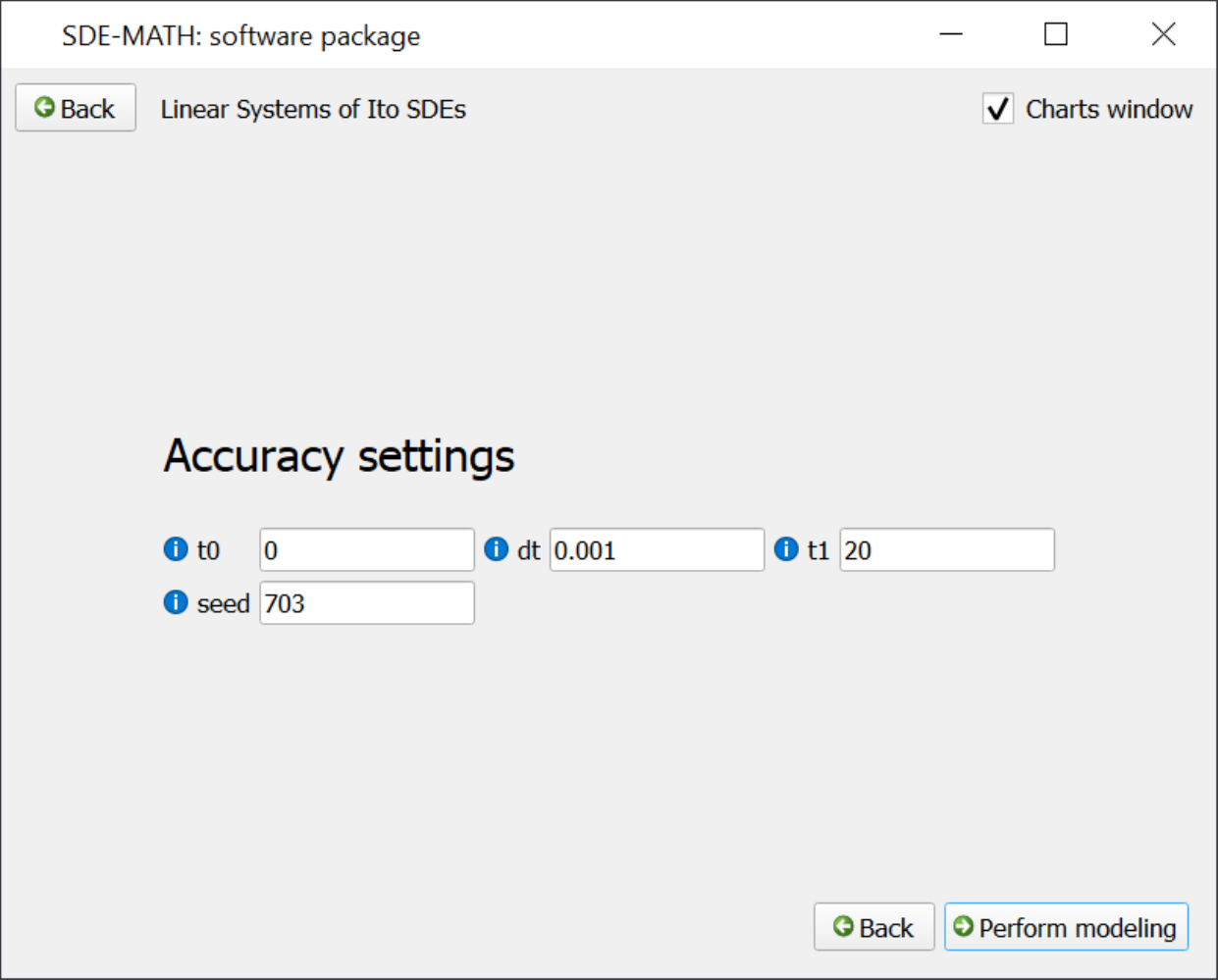}
    \caption{Correct data input\label{fig:linear_13}}
\end{figure}

\begin{figure}[H]
    \vspace{6mm}
    \centering
    \includegraphics[width=.7\textwidth]{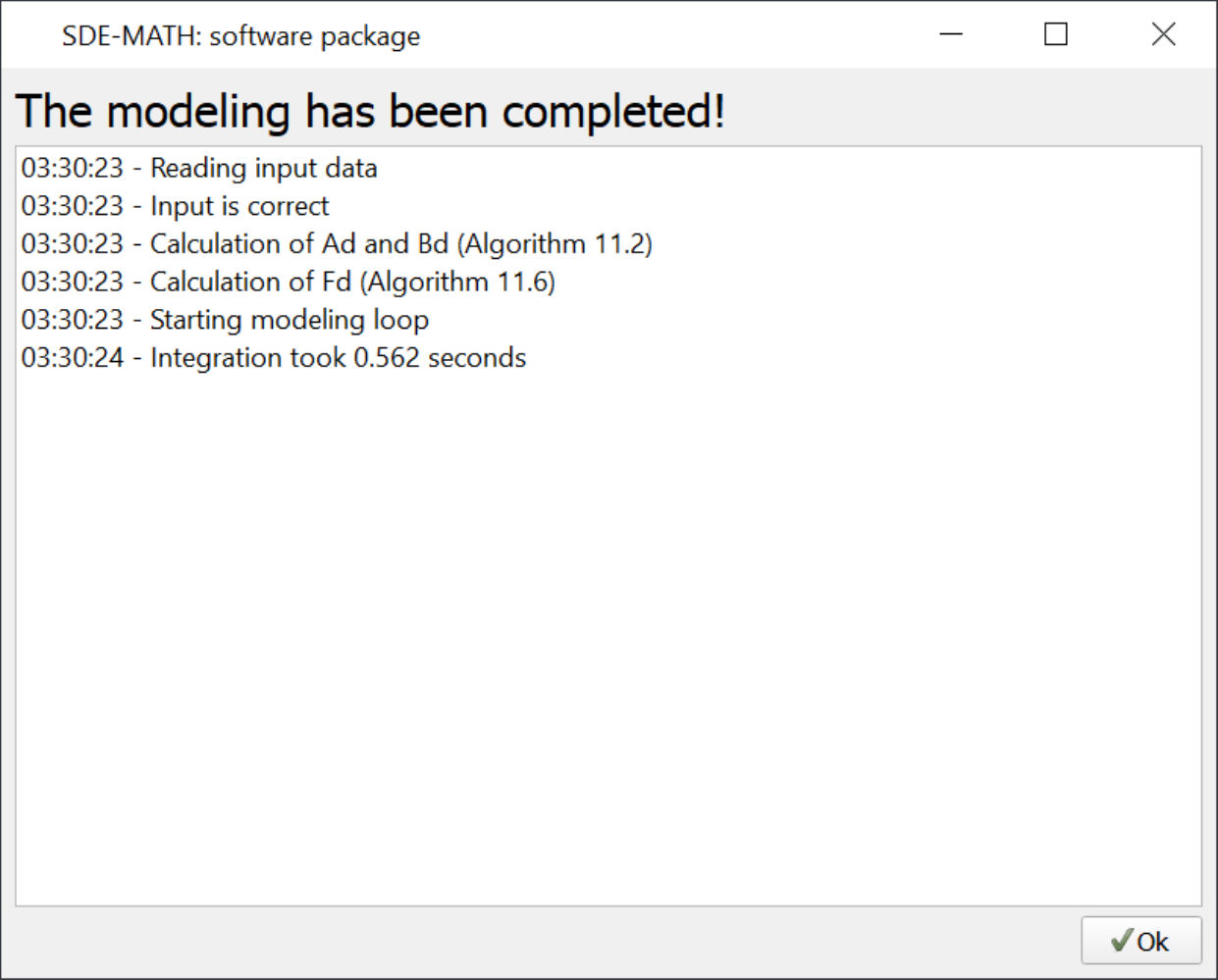}
    \caption{Modeling logs\label{fig:linear_14}}
\end{figure}

\begin{figure}[H]
    \centering
    \includegraphics[width=.9\textwidth]{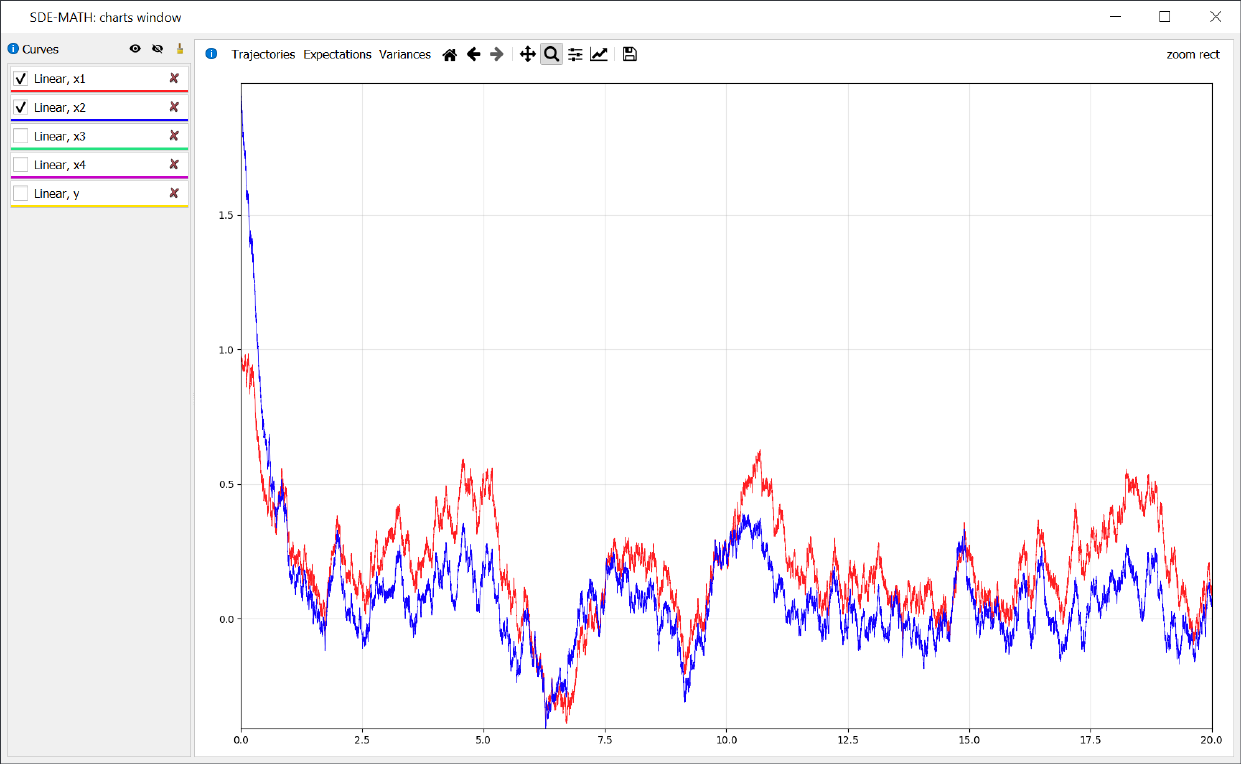}
    \caption{Modeling results (components of solution)\label{fig:linear_15}}
\end{figure}

\begin{figure}[H]
    \vspace{7mm}
    \centering
    \includegraphics[width=.9\textwidth]{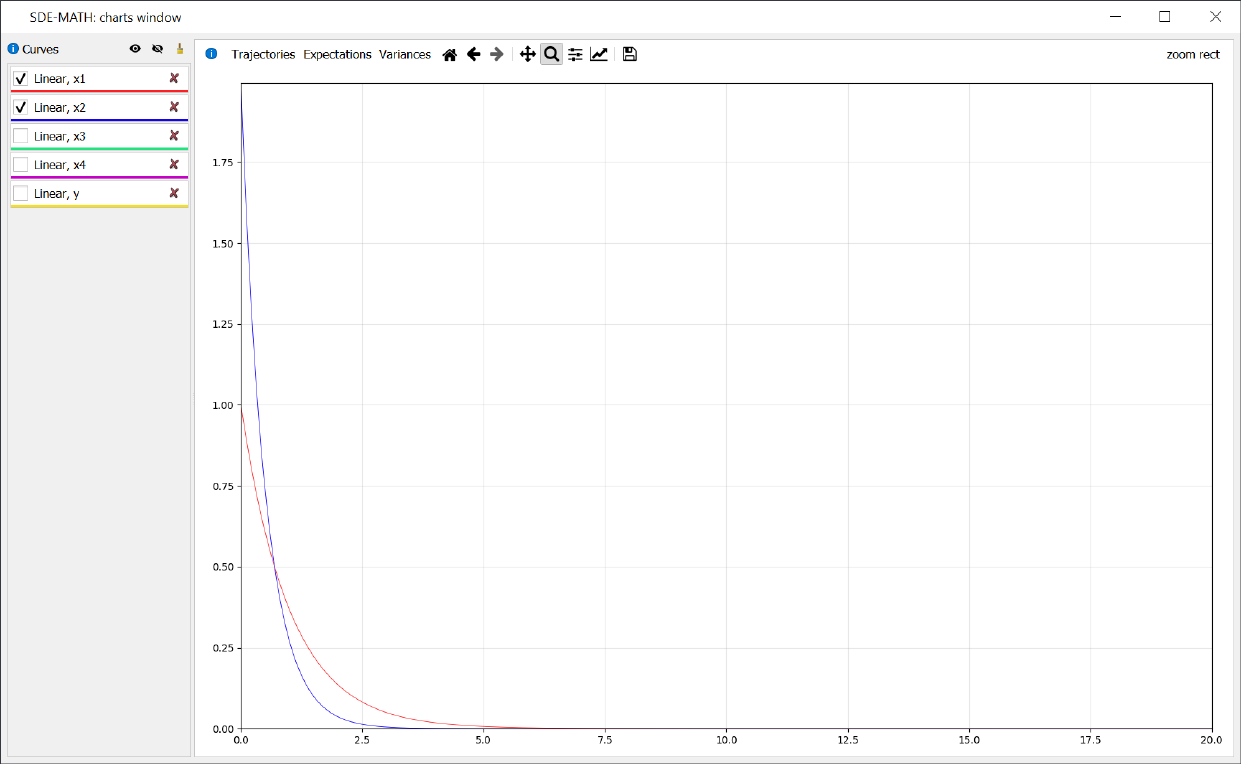}
    \caption{Modeling results (expectations)\label{fig:linear_16}}
\end{figure}

\begin{figure}[H]
    \centering
    \includegraphics[width=.9\textwidth]{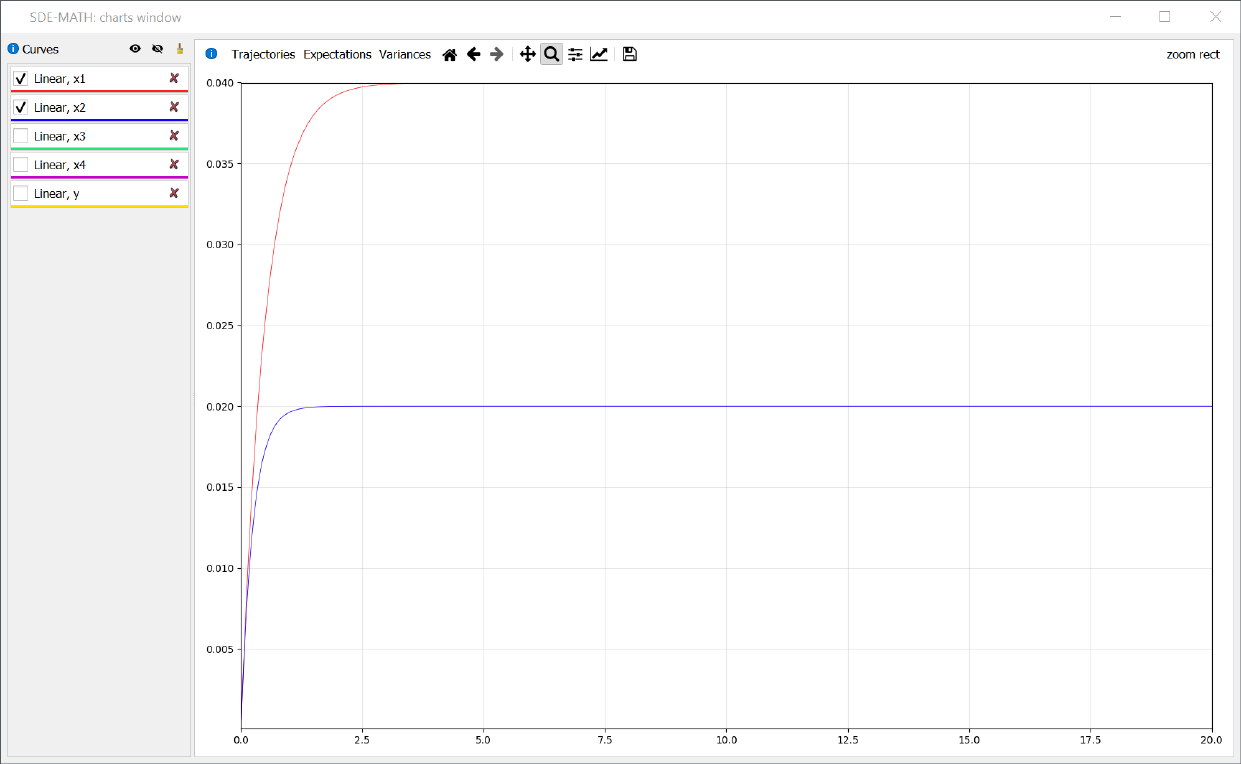}
    \caption{Modeling results (variances)\label{fig:linear_17}}
\end{figure}

The above part represents the structure of software package but not the
dynamics and user experience of it. Let us discuss the SDE-MATH software package
user experience on few examples provided further on Figures \ref{fig:preparing}--\ref{fig:linear_17}.
This examples represent two scenarios of the SDE-MATH software package use.

The database preparation screen is presented on Figure \ref{fig:preparing}. During the database
preparation this screen displays informational message and spinning visualizer of process continuation.

The screen that presented on Figure \ref{fig:greet} appears 
every time when software package runs unless user presses "Ok" button with marked 
checkbox. In such case this message screen will not be shown again.

On Figure \ref{fig:main_menu} the main menu dialog is presented. In this dialog user
can choose any strong numerical scheme for It\^o SDEs to perform modeling experiment.

The tooltip example can be seen on Figure \ref{fig:tooltip}. Such tooltips displayed with 
characteristic icon are placed all over software package interface to help user with explanations.

As noted earlier, the dedicated charts window is 
universal tool for visualization. The specific 
examples of such visualization are presented on Figure 
\ref{fig:plots},\ref{fig:nonlinear_11}, \ref{fig:linear_15}--\ref{fig:linear_17}. 

The initial state of input dialogs for nonlinear and linear It\^o 
SDEs are displayed on Figures \ref{fig:nonlinear_1} and \ref{fig:linear_1}.
At that moment user can start to input the data.

The example of wrong scalar data input is presented on Figures \ref{fig:nonlinear_2}, 
\ref{fig:nonlinear_7}, \ref{fig:linear_2}, and \ref{fig:linear_12}.
When user input wrong data the error message appears and "Next" or "Perform modeling"
button is blocked. The input field is being checked all the user data input process, 
and as soon as wrong character is entered notification pops up.

If scalar data is correct the "Next" button is automatically unblocked.
On Figures \ref{fig:nonlinear_3}, \ref{fig:nonlinear_8}, 
\ref{fig:linear_3}, and \ref{fig:linear_13} the examples of scenario are displayed.

On Figures \ref{fig:nonlinear_4}, \ref{fig:nonlinear_5}, and \ref{fig:linear_9} the example 
of correct matrix data input is presented. In this particular case the input is symbolic. 
Symbolic algebra input errors are much harder to determine so this 
is done on further stages, in modeling runtime.

In the other case when matrix input data are numerical, the validation is performed right after
user has finished input. The examples of incorrect matrix numerical input can be found on Figure
\ref{fig:linear_5}. 

When user finishes input with a success the "Next" or "Perform modeling" button is automatically unblocked.
On Figures \ref{fig:linear_4}, \ref{fig:linear_6}--\ref{fig:linear_8}, 
\ref{fig:linear_10}, and \ref{fig:linear_11} that can be clearly seen. 

The Figures \ref{fig:nonlinear_9},  \ref{fig:nonlinear_10}, and  \ref{fig:linear_14} displays 
sequence of log messages emerged during the modeling process.

After modeling has been done the focus moves to the charts window where obtained modeling results can be seen.
The results of modeling is displayed on Figures \ref{fig:nonlinear_11}, \ref{fig:linear_15}--\ref{fig:linear_17}. 
On Figures \ref{fig:linear_16} and \ref{fig:linear_17}
the expectations and variances of obtained components of solution are displayed.

\section{The Results Obtained Using the SDE-MATH Software Package}

This section represents the results
that were obtained with the SDE-MATH software
package at the current stage of the development.

\vspace{-2mm}

\subsection{The Calculated Fourier--Legendre Coefficients}

When application runs first time it
performs loading of Fourier--Legendre coefficients basic pack in the database from the files. 
Further, in Listings \ref{lst:c_xxx}--\ref{lst:c_xxxxxx} few examples of them can be seen.

\begin{lstlisting}[label={lst:c_xxx}, caption=\bfseries{The Fourier--Legendre coefficients $C^{000}_{j_3j_2j_1}$ examples}]
C_0:0:0 = 4/3
C_0:0:1 = -2/3
C_0:0:2 = 2/15
C_0:0:3 = 0
...
C_0:6:4 = -4/429
C_0:6:5 = 2/143
C_0:6:6 = 2/2145
C_1:0:0 = 2/3
...
C_47:33:44 = 3874457388633368/31334948307735906710660485
C_47:33:45 = 0
C_47:33:46 = 52892292737827468/2224781329849249376456894435
C_47:34:0 = 0
C_47:34:1 = 0
\end{lstlisting}
\begin{lstlisting}[label={lst:c_xxxx}, caption=\bfseries{The Fourier--Legendre coefficients $C^{0000}_{j_4j_3j_2j_1}$ examples}]
C_0:0:0:0 = 2/3
C_0:0:0:1 = -2/5
C_0:0:0:2 = 2/15
C_0:0:1:0 = -2/15
C_0:0:1:1 = 2/15
...
C_1:1:0:1 = -2/35
C_1:1:0:2 = 0
C_1:1:1:0 = 2/105
C_1:1:1:1 = 0
...
C_20:20:20:1 = -2401828/165607444685315115
C_20:20:20:2 = 0
C_20:20:20:3 = -1241929832/77669891557412788935
C_20:20:20:4 = 0
\end{lstlisting}
\begin{lstlisting}[label={lst:c_xxxxx}, caption=\bfseries{The Fourier--Legendre coefficients $C^{00000}_{j_5j_4j_3j_2j_1}$ examples}]
C_0:0:0:0:0 = 4/15
C_0:0:0:0:1 = -8/45
C_0:0:0:1:0 = -4/45
C_0:0:0:1:1 = 8/105
C_0:0:1:0:0 = 0
C_0:0:1:0:1 = 4/315
...
C_1:1:0:1:0 = -4/315
C_1:1:0:1:1 = 4/315
C_1:1:1:0:0 = 2/105
C_1:1:1:0:1 = -8/945
C_1:1:1:1:0 = 2/945
C_1:1:1:1:1 = 0
...
C_20:20:20:20:1 = 0
C_20:20:20:20:2 = -249877207023610010/10969028984480026752856704371
C_20:20:20:20:3 = 0
C_20:20:20:20:4 = -307937246954575016/571102494076952592887484313076115
\end{lstlisting}
\begin{lstlisting}[label={lst:c_xxxxxx}, caption=\bfseries{The Fourier--Legendre coefficients $C^{000000}_{j_6j_5j_4j_3j_2j_1}$ examples}]
C_0:0:0:0:0:0 = 4/45
C_0:0:0:0:0:1 = -4/63
C_0:0:0:0:0:2 = 2/63
C_0:0:0:0:1:0 = -4/105
...
C_2:1:0:1:0:2 = -2/1575
C_2:1:0:1:1:0 = 38/22275
C_2:1:0:1:1:1 = -2/1575
C_2:1:0:1:1:2 = 68/81081
...
C_15:15:15:15:15:15 = 0
C_15:15:15:15:15:16 = -798538765964/243076352242280511713913783475
C_15:15:15:15:15:17 = 0
C_15:15:15:15:15:18 = -59075427603328/17302616709609603697454044769175

\end{lstlisting}

\vspace{-2mm}

\subsection{Accuracy Settings}

From Theorem 7 (see formulas (\ref{formula0})--(\ref{formula11})) it follows that the number $p$ in the 
formula (\ref{letusdenote}) should be chosen individually for various combinations of indices 
$i_1,\ldots,i_k\in\{1,\ldots,m\}.$  As follows from Listing \ref{lst:q_hypothesis} (see below)
and the results of work \cite{63}, these numbers $p$ in the overwhelming majority of 
cases do not exceed the number $p$ from the formula (\ref{formula0}).
Moreover, all the mentioned numbers $p$ are many times less than the number $p$ selected 
using the formula (\ref{star00011}) (due to the presence of the multiplier factor $k!$ on the left-hand side of (\ref{star00011})).

In this work, we have replaced the mentioned numbers $p$ for all possible 
combinations of indices $i_1,\ldots,i_k\in\{1,\ldots,m\}$ 
with the number $p$ according to the formula (\ref{formula0}).
This is possible due to the results of Listing \ref{lst:q_hypothesis2}. This listing shows that the above replacement 
does not lead to noticeable accuracy loss of the
mean square approximation of iterated It\^{o} stochastic integrals (for more details see \cite{63}).

Thus, in this paper we decided
to exclude the multiplier factor $k!$ in the
conditions for choosing the numbers
$q_1,\ldots,q_{15}$ (see (\ref{ress1})--(\ref{fas29a})). Recall that
these numbers are used to construct
the approximations of iterated
It\^o and Stratonovich stochastic integrals
from the numerical schemes
(\ref{al2})--(\ref{al5}),
(\ref{al2x})--(\ref{al5x}).
The test script was
written. The results of its work are presented in Listings 
\ref{lst:q_hypothesis} and \ref{lst:q_hypothesis2}, where

\begin{enumerate}

    \item dt is the integration step;

    \item q1(1,2) means $p$ from (\ref{formula1}),
    q1(2,3) means $p$ from (\ref{formula2}),
    q1(1,3) means $p$ from (\ref{formula3}),
    q1 means $p$ from (\ref{formula0}) for $k=3$;

    \item $C = 1$ (see (\ref{uslov}) and (\ref{uslov1}));

    \item error 1 means the left-hand side of (\ref{fas1a});

    \item error 2 means the left-hand side of (\ref{formula1}) divided by $(T-t)^3$;

    \item error 3 means the left-hand side of (\ref{formula3}) divided by $(T-t)^3$;

    \item error 4 means the left-hand side of (\ref{formula2}) divided by $(T-t)^3$.

\end{enumerate}

The above idea of calculation of the numbers
$q_1,\ldots,q_{15}$
is described in Listing \ref{lst:q}.

\vspace{3mm}

\begin{lstlisting}[label={lst:q_hypothesis}, caption=\bfseries{Accuracy calculation module}]

dt = 0.011
	q1        = 12
	q1 (1, 2) = 6
	q1 (1, 3) = 12
	q1 (2, 3) = 6

dt = 0.008
	q1        = 16
	q1 (1, 2) = 8
	q1 (1, 3) = 16
	q1 (2, 3) = 8

dt = 0.0045
	q1        = 28
	q1 (1, 2) = 14
	q1 (1, 3) = 28
	q1 (2, 3) = 14

dt = 0.0035
	q1        = 36
	q1 (1, 2) = 18
	q1 (1, 3) = 36
	q1 (2, 3) = 18

dt = 0.0027
	q1        = 47
	q1 (1, 2) = 23
	q1 (1, 3) = 47
	q1 (2, 3) = 23

dt = 0.0025
	q1        = 50
	q1 (1, 2) = 25
	q1 (1, 3) = 51
	q1 (2, 3) = 25

Process finished with exit code 0

\end{lstlisting}

\vspace{7mm}

\begin{lstlisting}[label={lst:q_hypothesis2}, caption=\bfseries{Accuracy calculation module}]

dt = 0.011
	error 1 = 0.010153888451696458
	q1        = 12
	error 2 = 0.005076944225848201
	q1 (1, 2) = 12
	error 3 = 0.010307776903394072
	q1 (1, 3) = 12
	error 4 = 0.005076944225848284
	q1 (2, 3) = 12

dt = 0.008
	error 1 = 0.007681193827577537
	q1        = 16
	error 2 = 0.003840596913789046
	q1 (1, 2) = 16
	error 3 = 0.0077866300793989485
	q1 (1, 3) = 16
	error 4 = 0.003840596913789157
	q1 (2, 3) = 16

dt = 0.0045
	error 1 = 0.004432832059862973
	q1        = 28
	error 2 = 0.0022164160299319446
	q1 (1, 2) = 28
	error 3 = 0.004479699207443705
	q1 (1, 3) = 28
	error 4 = 0.002216416029932139
	q1 (2, 3) = 28

dt = 0.0035
	error 1 = 0.0034564405520411956
	q1        = 36
	error 2 = 0.0017282202760207088
	q1 (1, 2) = 36
	error 3 = 0.003488223569838411
	q1 (1, 3) = 36
	error 4 = 0.0017282202760210141
	q1 (2, 3) = 36

dt = 0.0027
	error 1 = 0.0026523659377455377
	q1        = 47
	error 2 = 0.00132618296887127
	q1 (1, 2) = 47
	error 3 = 0.0026731529281250332
	q1 (1, 3) = 47
	error 4 = 0.0013261829688714366
	q1 (2, 3) = 47

dt = 0.0025
	error 1 = 0.002494053620431952
	q1        = 50
	error 2 = 0.0012470268102122428
	q1 (1, 2) = 50
	error 3 = 0.0025128597161119537
	q1 (1, 3) = 50
	error 4 = 0.0012470268102123538
	q1 (2, 3) = 50

Process finished with exit code 0
\end{lstlisting}

\vspace{7mm}

\subsection{Testing Example (Nonlinear System of It\^o SDEs)}

\vspace{1mm}

The input data for testing of the SDE-MATH 
software package correspond to the
autonomous variant of nonlinear system of It\^{o} SDE (\ref{1.5.2})
with multidimensional non-commutative noise. More precisely,
we choose $n=2,$ $m=2$, ${\bf x}_0^{(1)}=1,$ ${\bf x}_0^{(2)}=1.5,$

\newpage

\noindent
$$
    {\bf a}\left({\bf x}^{(1)},{\bf x}^{(2)}\right) = \left(\begin{matrix}
    -5 {\bf x}^{(1)}\\\\
    -5 {\bf x}^{(2)}
    \end{matrix}\right),
$$

$$
    B\left({\bf x}^{(1)},{\bf x}^{(2)}\right) = \left(\begin{matrix}
        0.5 \cdot \sin\left({\bf x}^{(1)}\right) & {\bf x}^{(2)}\\\\
        {\bf x}^{(2)} & 0.5 \cdot \cos\left({\bf x}^{(1)}\right)
    \end{matrix}\right).
$$

\vspace{5mm}

\begin{figure}[H]
    \centering
    \includegraphics[width=.9\textwidth]{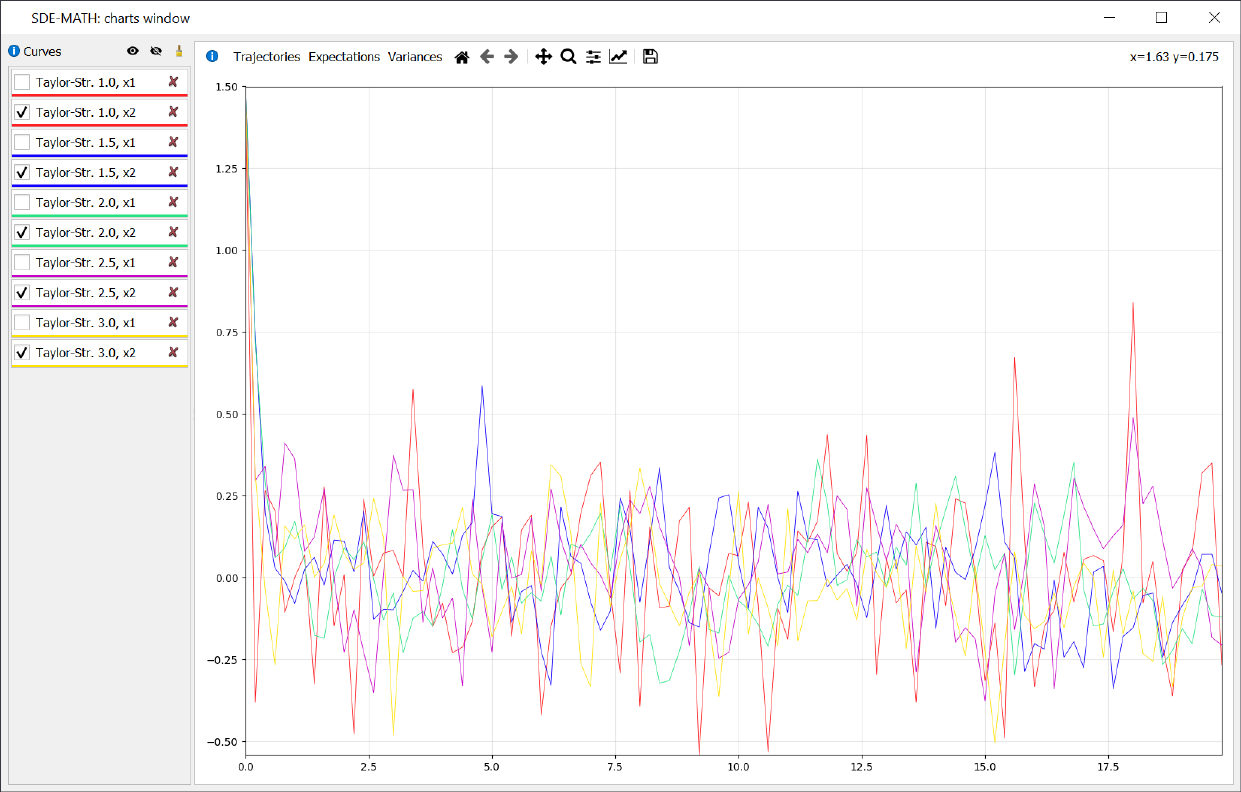}
    \caption{Strong Taylor--Stratonovich schemes of orders 1.0, 1.5, 2.0, 2.5, and 3.0 (${\bf x}_t^{(2)}$ component, $C = 4,$ $dt = 0.2$)\label{fig:straton_3p0_big_7}}
\end{figure}

Figures \ref{fig:ito_1p5_small_logs}--\ref{fig:straton_3p0_big_7}
related to the strong high-order Taylor--It\^{o} and Taylor--Stra\-to\-no\-vich
schemes (\ref{al1})--(\ref{al5}), (\ref{al1x})--(\ref{al5x}) for 
the It\^o SDE (\ref{1.5.2}) represent modeling results.

Test machine specifications are CPU
with maximum core frequency 4.2 GHz and 16GB of RAM.


\vspace{3mm}

\subsection{Visualization and Numerical Results for Nonlinear System of It\^{o} SDEs Obtained via the SDE-MATH Software Package}

\begin{figure}[H]
    \vspace{5mm}
    \hspace*{\fill}
    \begin{subfigure}[b]{.45\textwidth}
        \includegraphics[width=\textwidth]{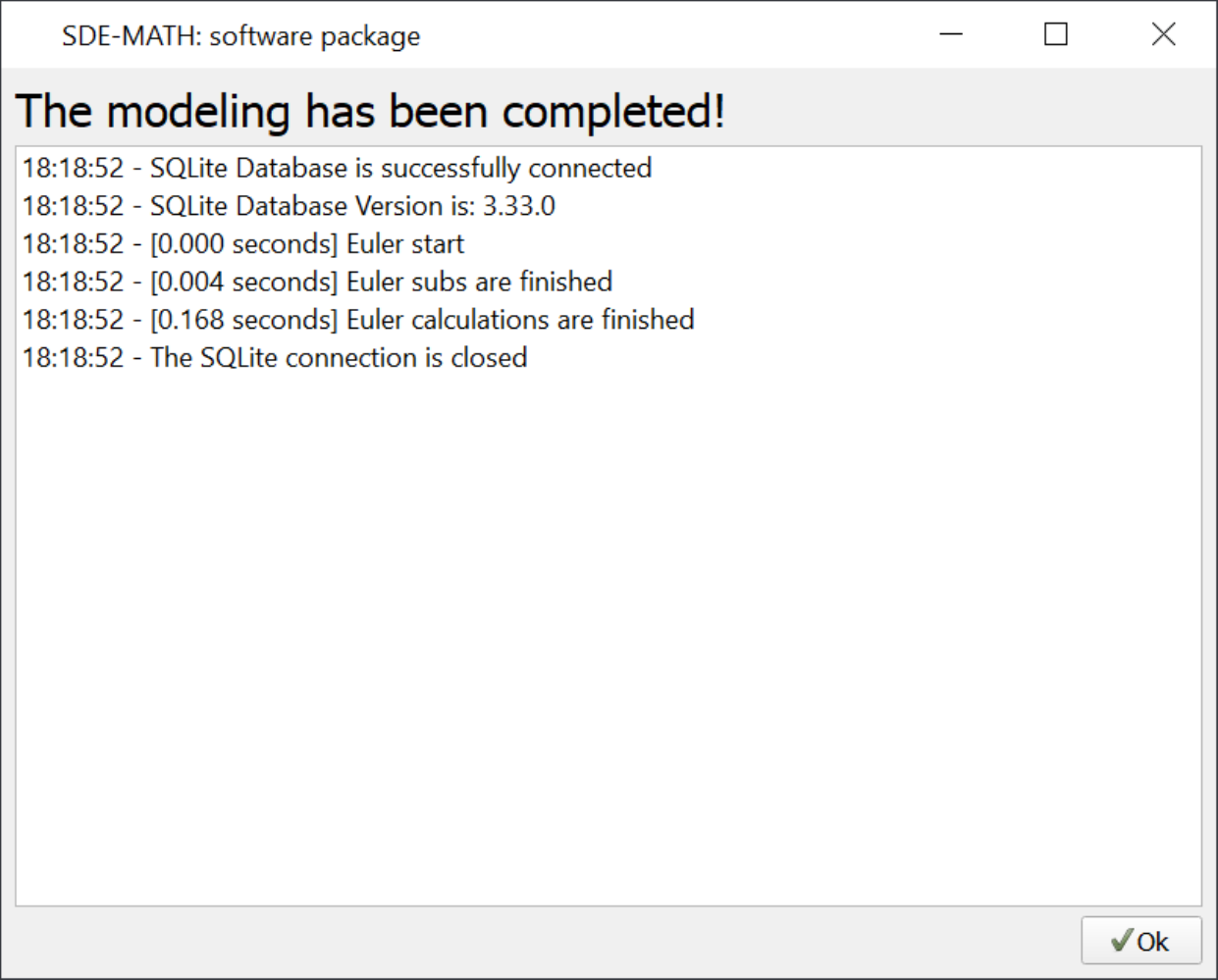}
        \caption*{Euler scheme ($dt = 0.001$)\label{fig:ito_1p5_small_1}}
    \end{subfigure}
    \hfill
    \begin{subfigure}[b]{.45\textwidth}
        \includegraphics[width=\textwidth]{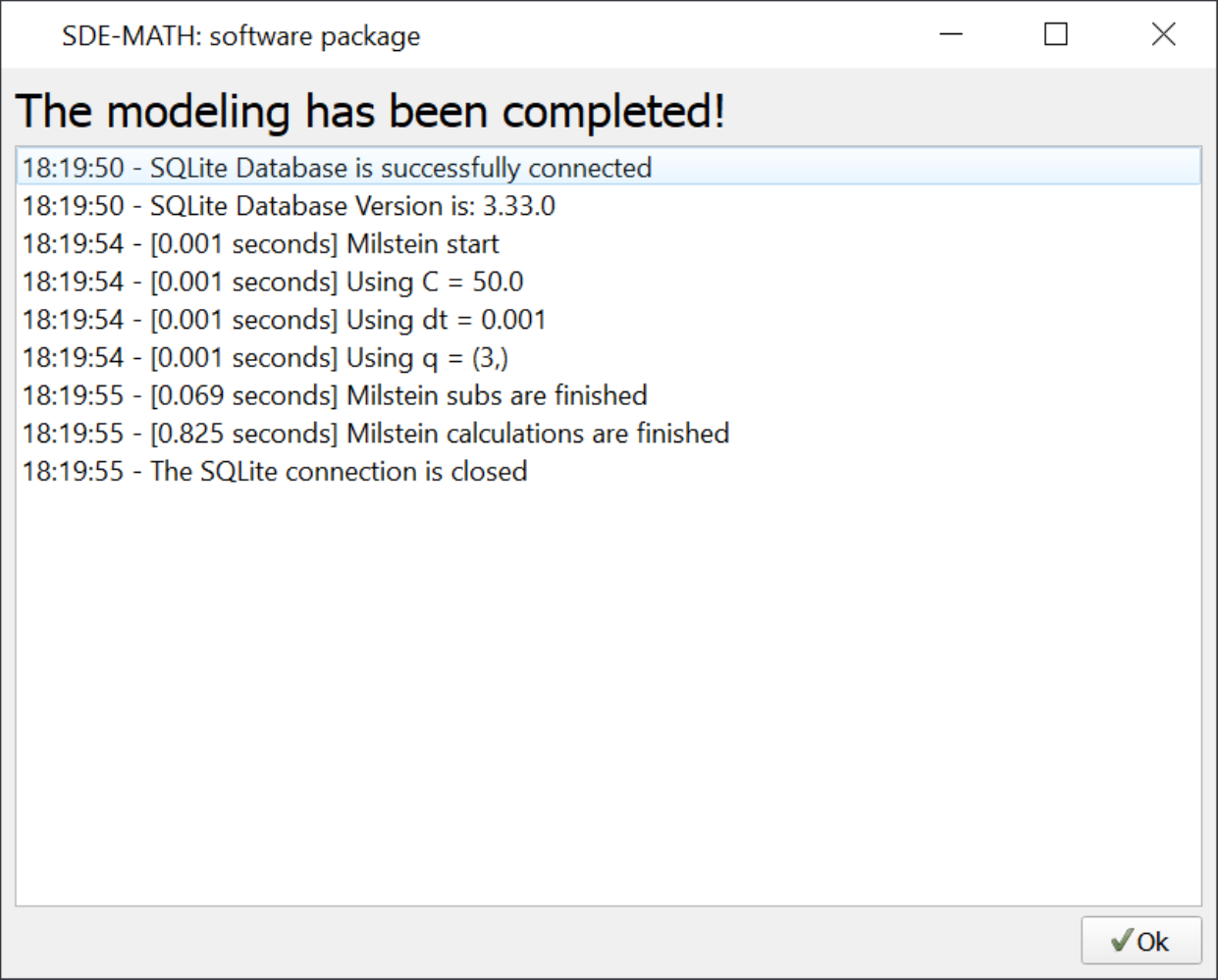}
        \caption*{Milstein scheme ($C = 50,$ $dt = 0.001$)\label{fig:ito_1p5_small_2}}
    \end{subfigure}
    \hspace*{\fill}

    \caption{Modeling logs\label{fig:ito_1p5_small_logs}}
\end{figure}

\begin{figure}[H]
    \vspace{13mm}
    \centering
    \includegraphics[width=.45\textwidth]{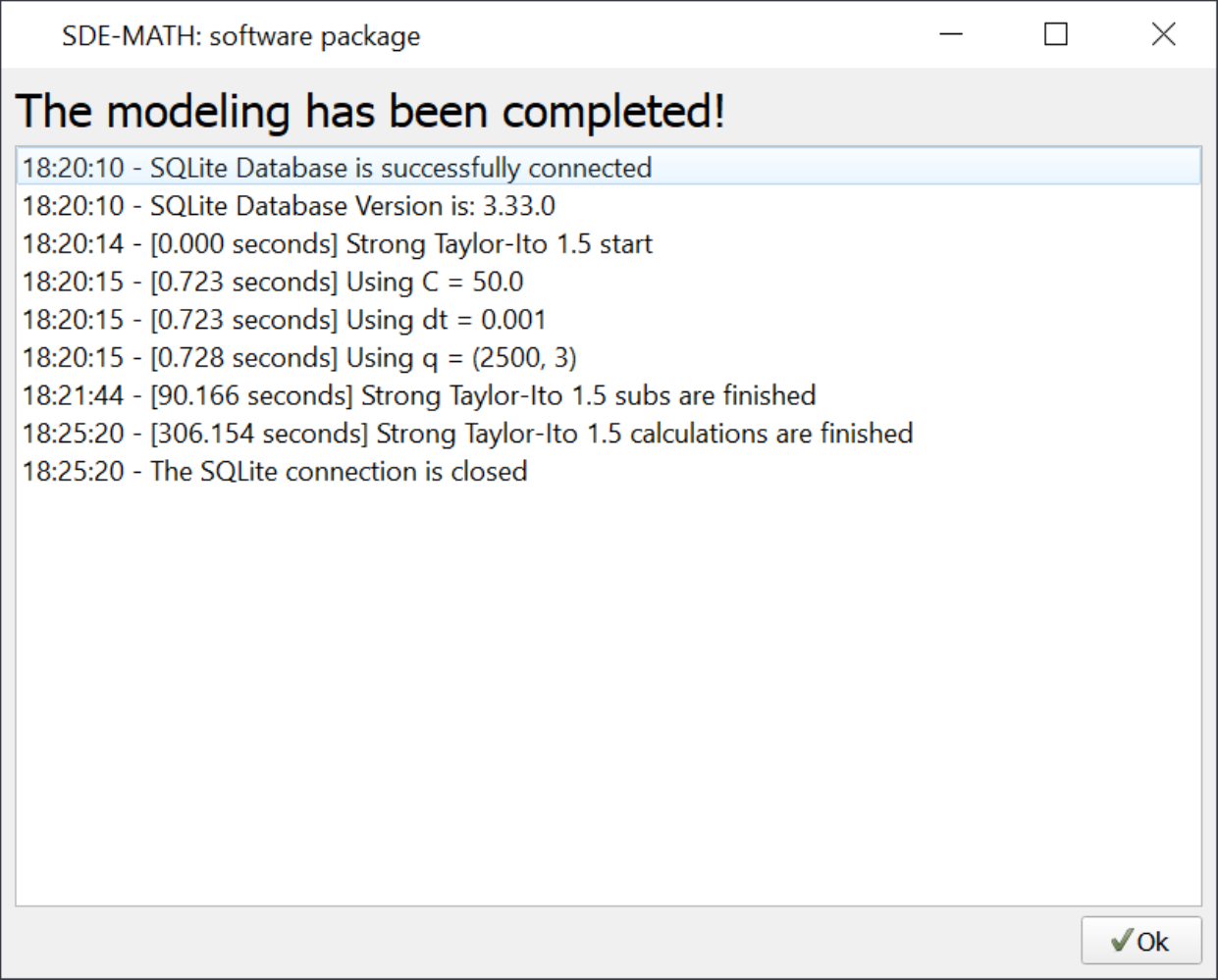}
    \caption{Strong Taylor--It\^o scheme of order 1.5 ($C = 50,$ $dt = 0.001$)\label{fig:ito_1p5_small_3}}
\end{figure}

\begin{figure}[H]
    \vspace{10mm}
    \centering
    \includegraphics[width=.9\textwidth]{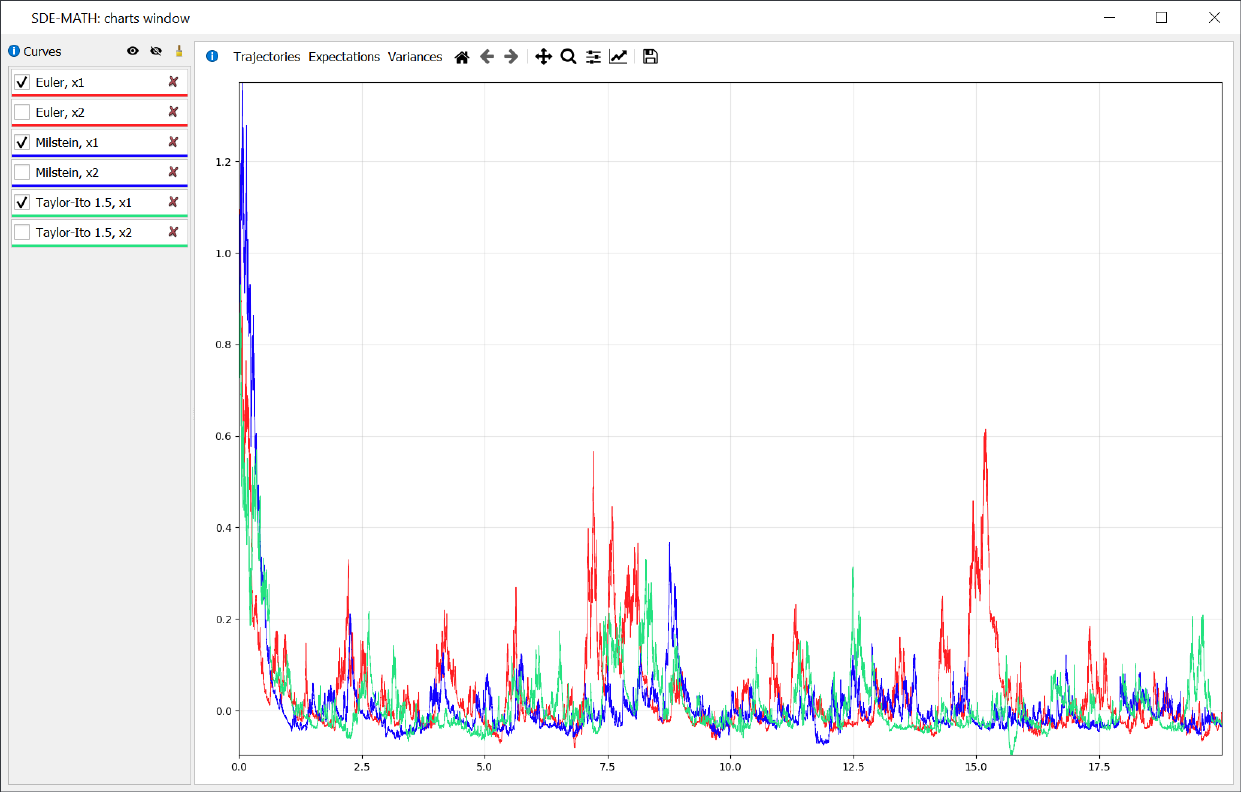}
    \caption{Strong Taylor--It\^o schemes of orders $0.5,$ $1.0,$ and $1.5$ (${\bf x}_t^{(1)}$ component, $C = 50,$ $dt = 0.001$)\label{fig:ito_1p5_small_4}}
\end{figure}

\begin{figure}[H]
    \vspace{10mm}
    \centering
    \includegraphics[width=.9\textwidth]{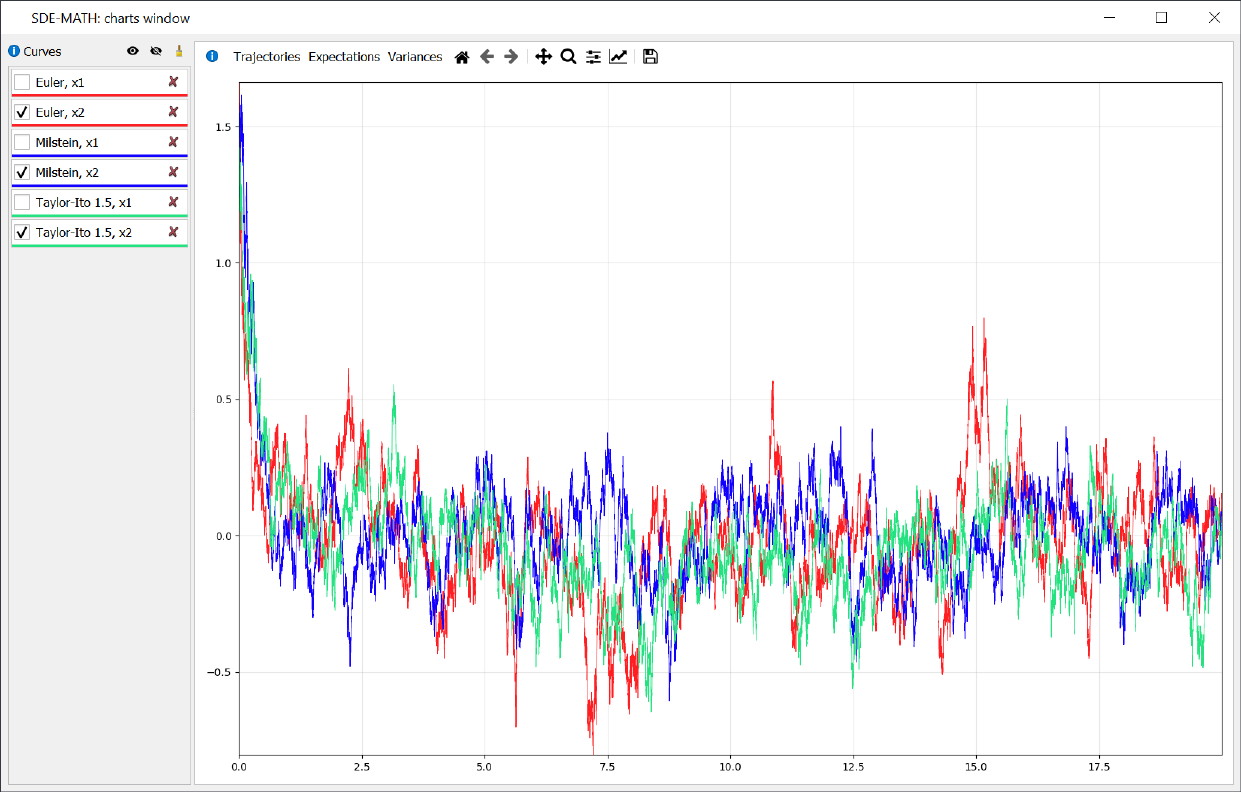}
    \caption{Strong Taylor--It\^o schemes of orders $0.5,$ $1.0,$ and $1.5$ (${\bf x}_t^{(2)}$ component, $C = 50,$ $dt = 0.001$)\label{fig:ito_1p5_small_5}}
\end{figure}

\begin{figure}[H]
    \vspace{8mm}
    \centering
    \hspace*{\fill}
    \begin{subfigure}[b]{.45\textwidth}
        \includegraphics[width=\textwidth]{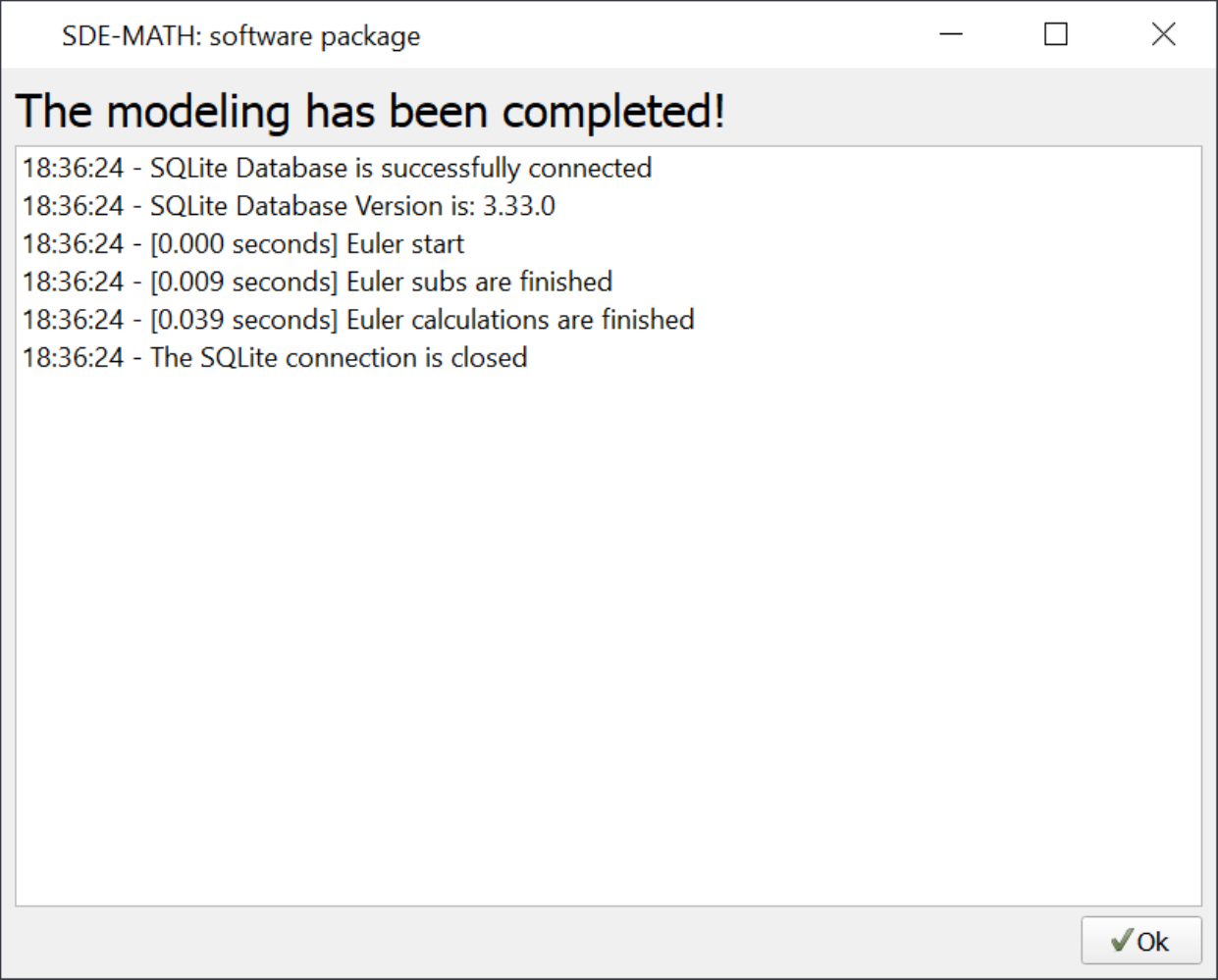}
        \caption*{Euler scheme ($dt = 0.005$)\label{fig:ito_2p0_small_1}}
    \end{subfigure}
    \hfill
    \begin{subfigure}[b]{.45\textwidth}
        \includegraphics[width=\textwidth]{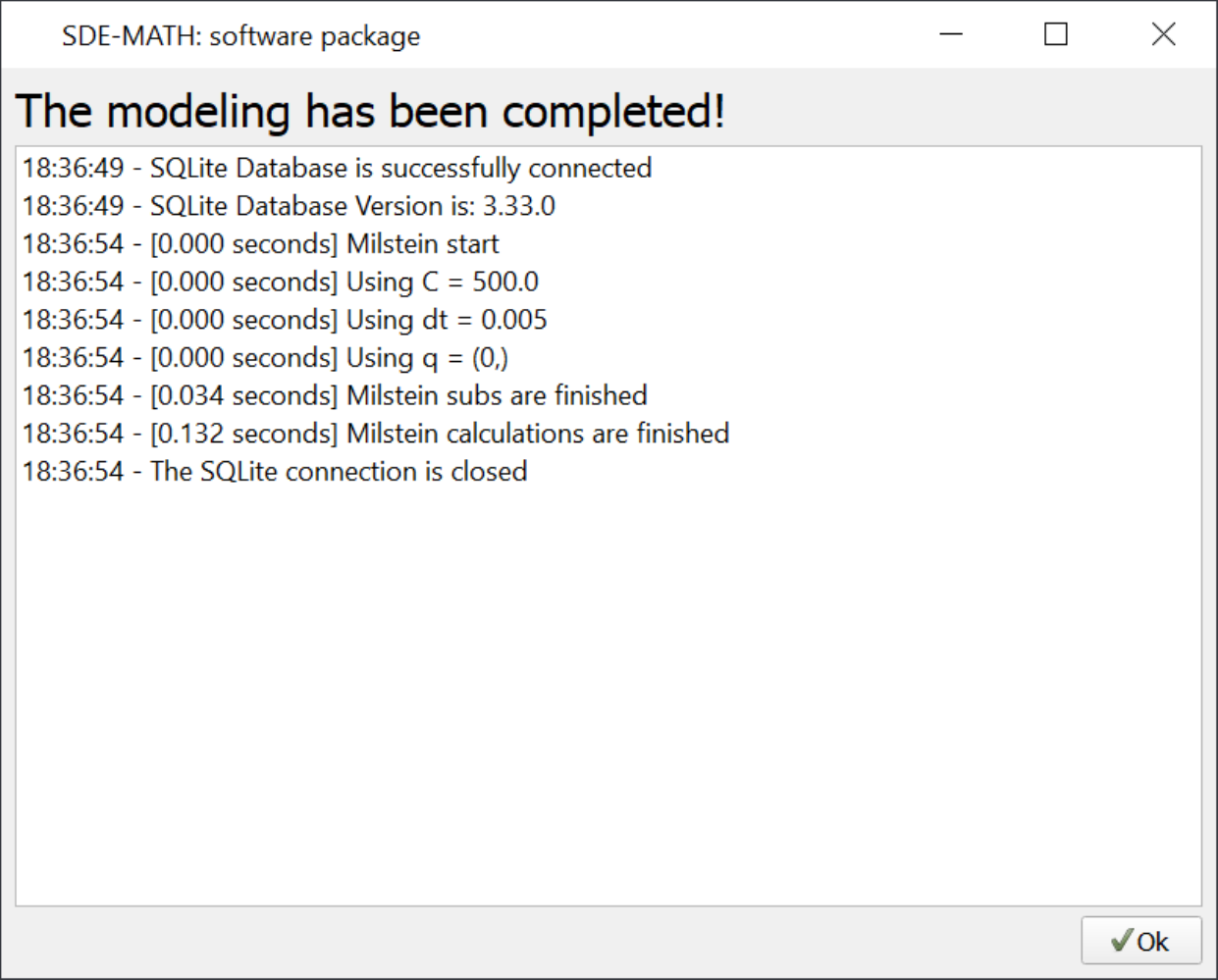}
        \caption*{Milstein scheme ($C = 500,$ $dt = 0.005$)\label{fig:ito_2p0_small_2}}
    \end{subfigure}
    \hspace*{\fill}

    \caption{Modeling logs\label{fig:ito_2p0_small_logs1}}

\end{figure}

\begin{figure}[H]
    \vspace{10mm}
    \hspace*{\fill}
    \begin{subfigure}[b]{.45\textwidth}
        \includegraphics[width=\textwidth]{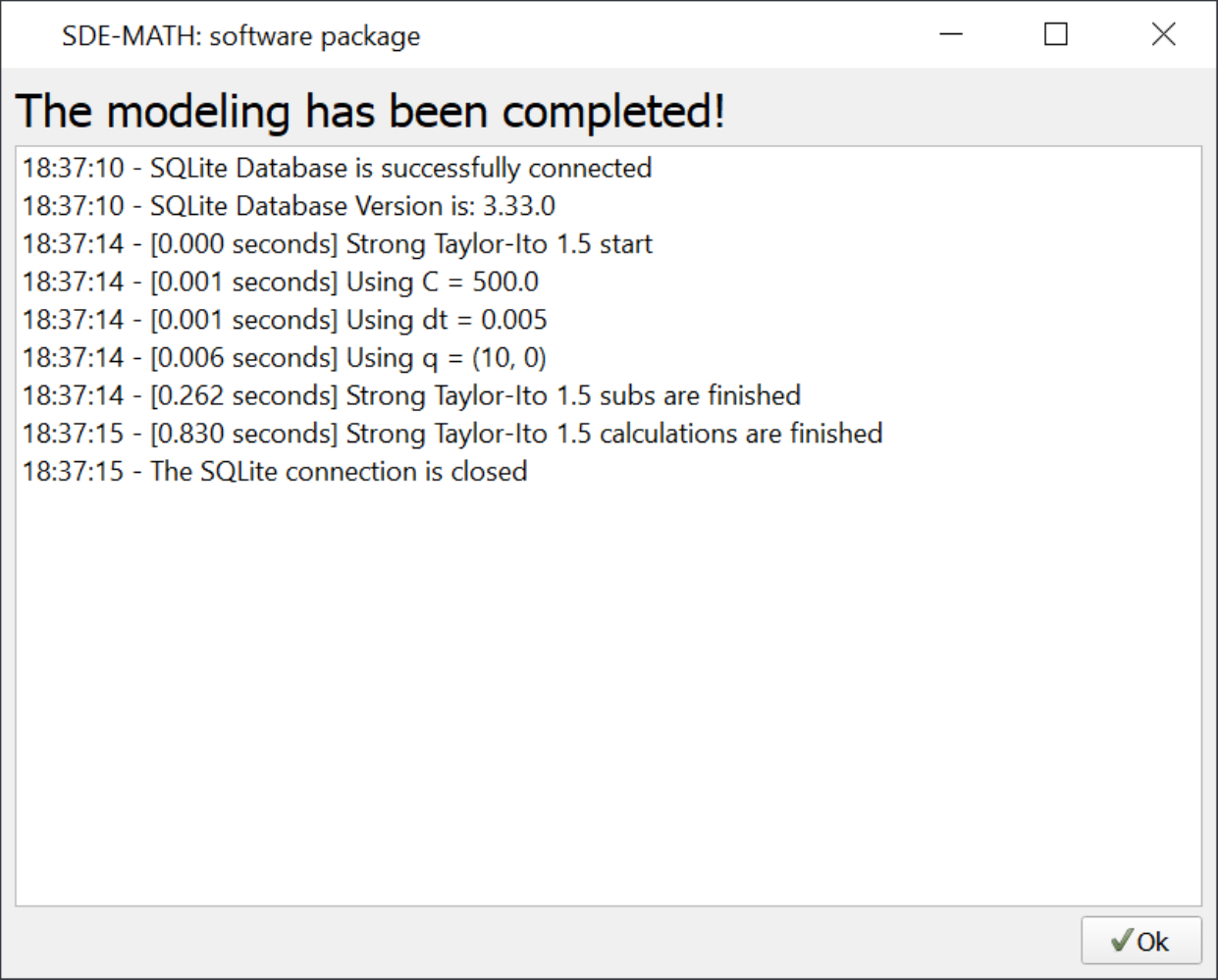}
        \caption*{Strong Taylor--It\^o scheme of order 1.5 ($C = 500,$ $dt = 0.005$)\label{fig:ito_2p0_small_3}}
    \end{subfigure}
    \hfill
    \begin{subfigure}[b]{.45\textwidth}
        \includegraphics[width=\textwidth]{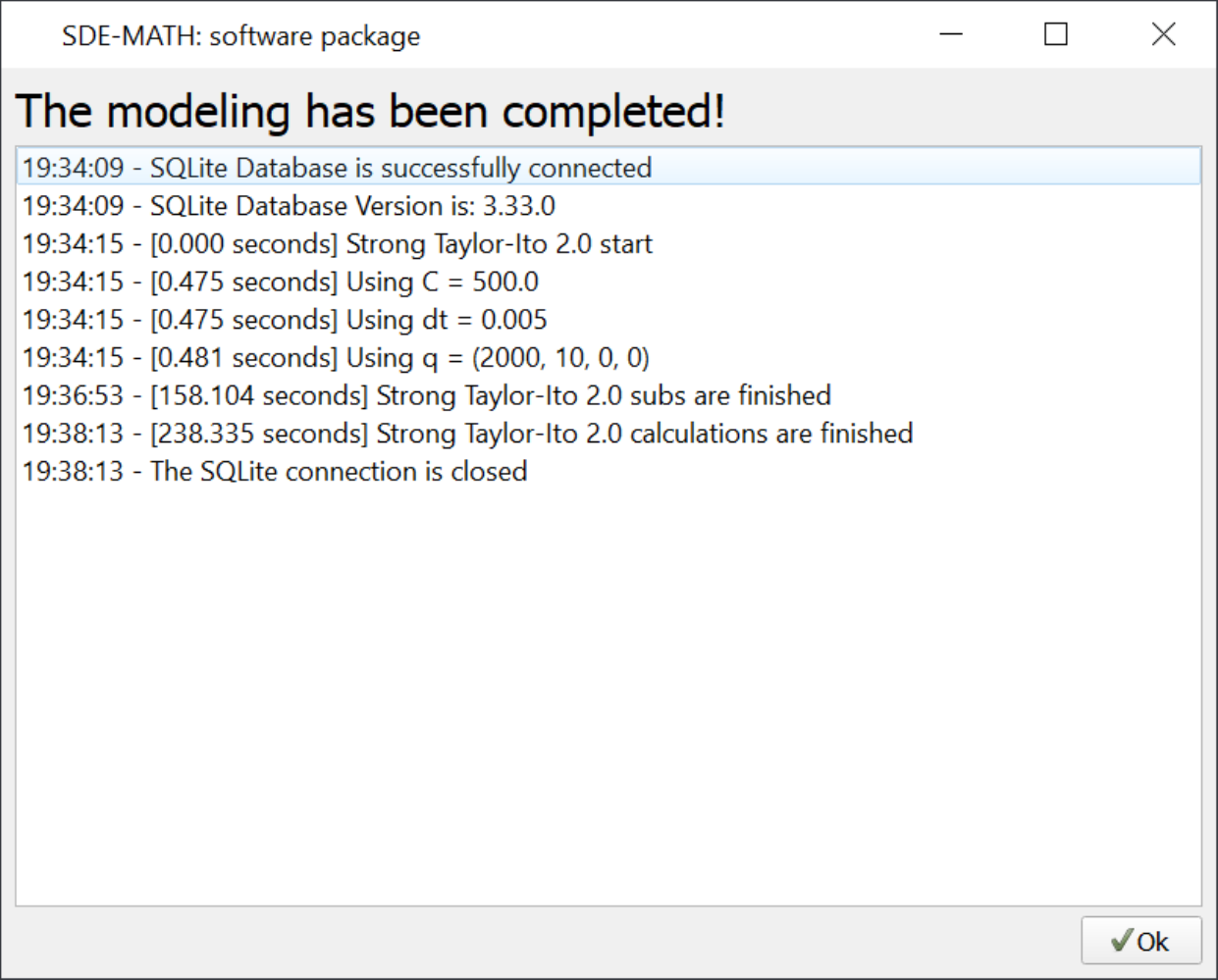}
        \caption*{Strong Taylor--It\^o scheme of order 2.0 ($C = 500,$ $dt = 0.005$)\label{fig:ito_2p0_small_4}}
    \end{subfigure}
    \hspace*{\fill}

    \caption{Modeling logs\label{fig:ito_2p0_small_logs2}}
    
\end{figure}

\begin{figure}[H]
    \vspace{5mm}
    \centering
    \includegraphics[width=.9\textwidth]{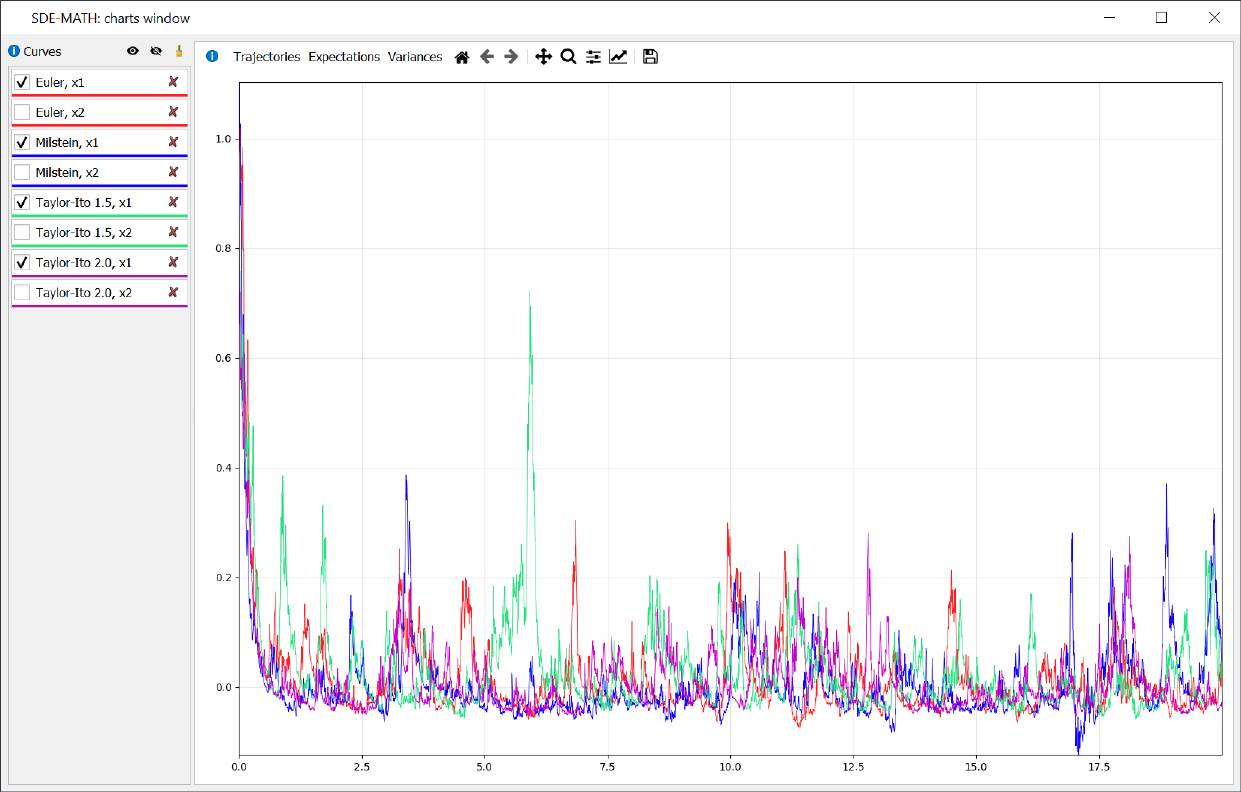}
    \caption{Strong Taylor--It\^o schemes of orders 0.5, 1.0, 1.5, and 2.0 (${\bf x}_t^{(1)}$ component, $C = 500,$ $dt = 0.005$)\label{fig:ito_2p0_small_5}}
\end{figure}

\begin{figure}[H]
    \centering
    \includegraphics[width=.9\textwidth]{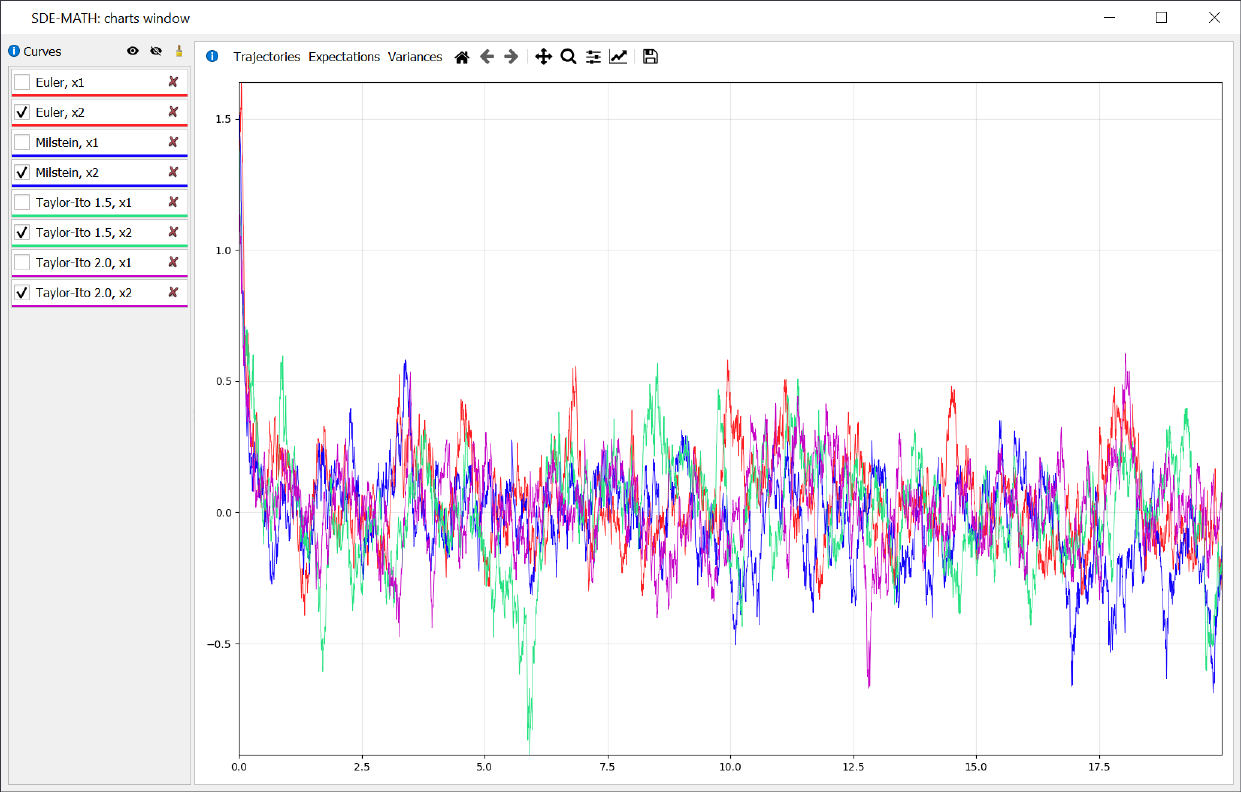}
    \caption{Strong Taylor--It\^o schemes of orders 0.5, 1.0, 1.5, and 2.0 (${\bf x}_t^{(2)}$ component, $C = 500,$ $dt = 0.005$)\label{fig:ito_2p0_small_6}}
\end{figure}

\begin{figure}[H]
    \centering
    \includegraphics[width=.9\textwidth]{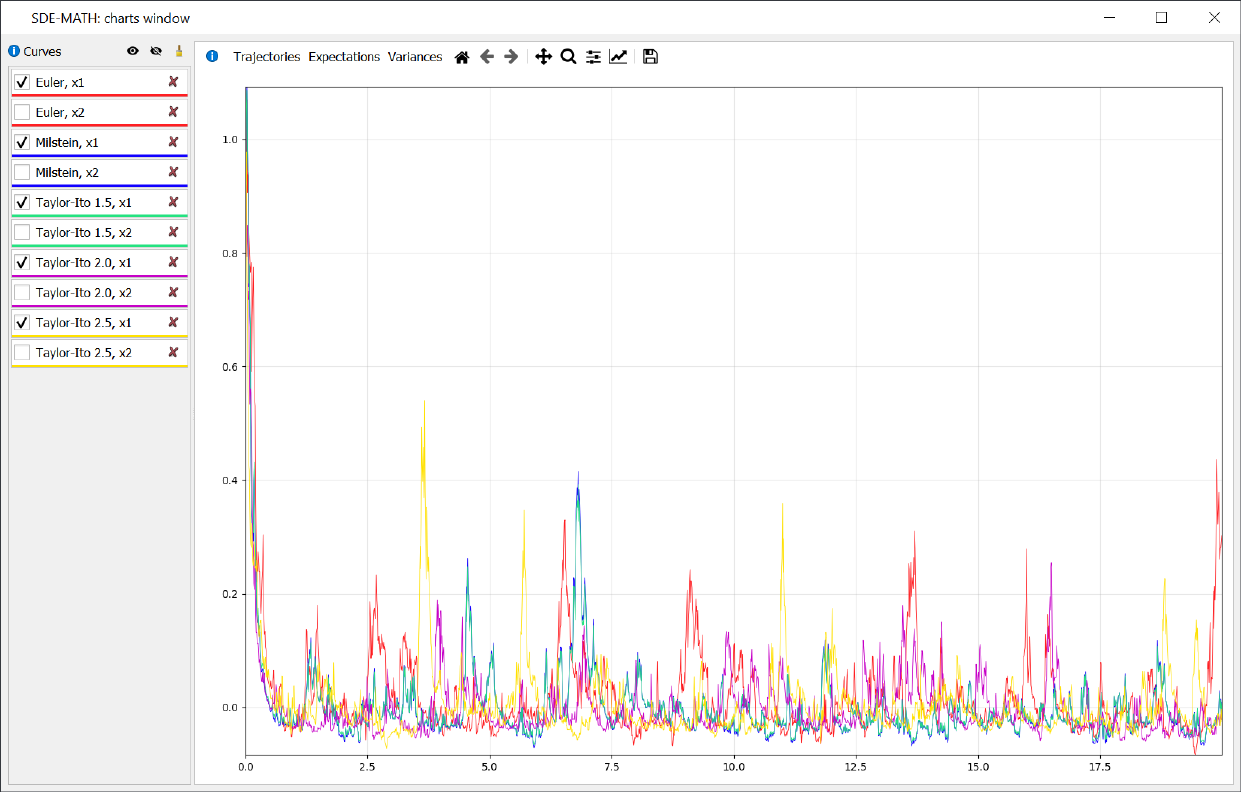}
    \caption{Strong Taylor--It\^o schemes of orders 0.5, 1.0, 1.5, 2.0, and 2.5 (${\bf x}_t^{(1)}$ component, $C = 7500,$ $dt = 0.01$)\label{fig:ito_2p5_small_6}}
\end{figure}

\begin{figure}[H]
    \centering

    \hspace*{\fill}
    \begin{subfigure}[b]{.45\textwidth}
        \includegraphics[width=\textwidth]{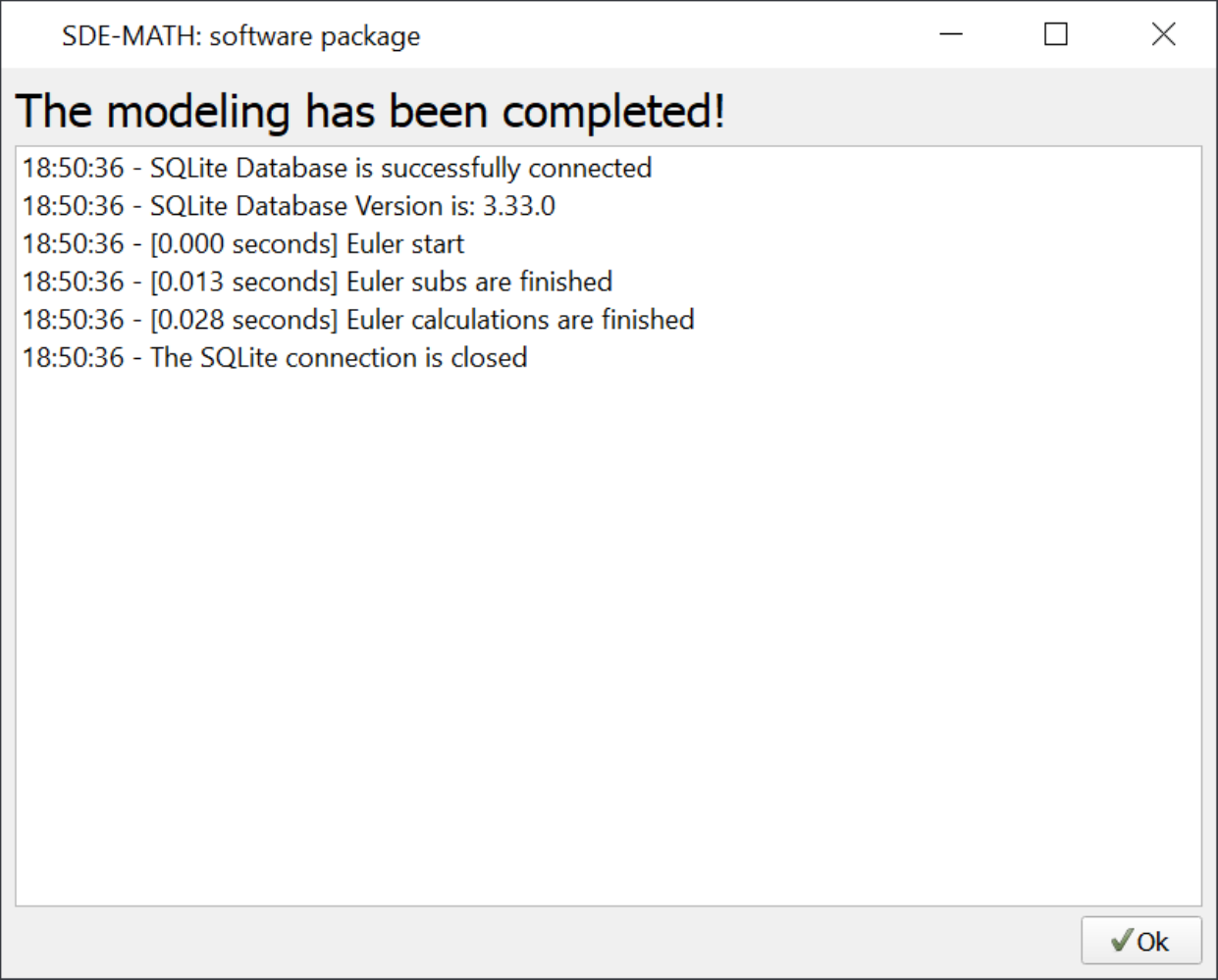}
        \caption*{Euler scheme ($dt = 0.01$)\label{fig:ito_2p5_small_1}}
    \end{subfigure}
    \hfill
    \begin{subfigure}[b]{.45\textwidth}
        \includegraphics[width=\textwidth]{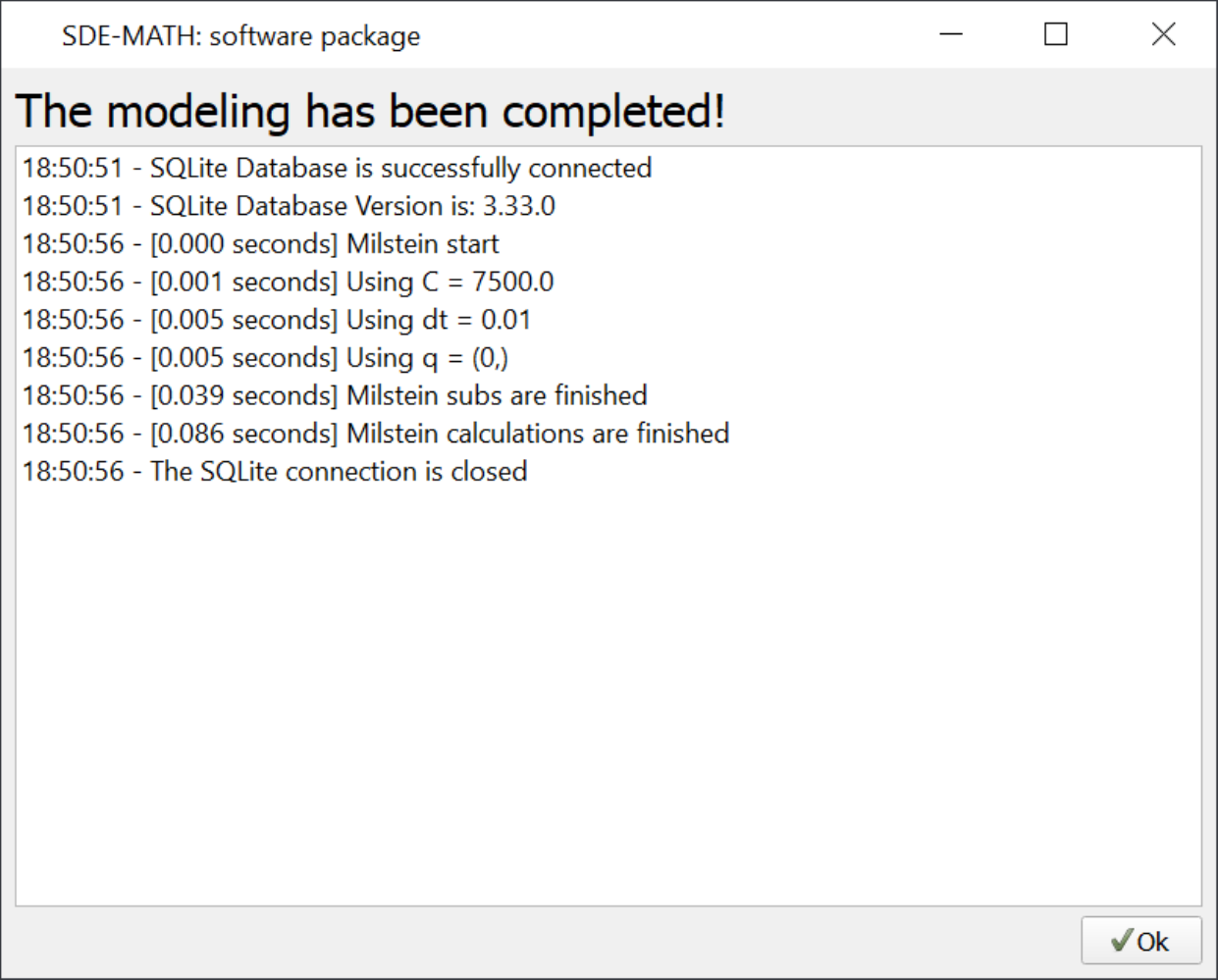}
        \caption*{Milstein scheme ($C = 7500,$ $dt = 0.01$)\label{fig:ito_2p5_small_2}}
    \end{subfigure}
    \hspace*{\fill}

    \vspace{2mm}
    \hspace*{\fill}
    \begin{subfigure}[b]{.45\textwidth}
        \includegraphics[width=\textwidth]{figures/little_ito_3/2.pdf}
        \caption*{Strong Taylor--It\^o scheme of order 1.5 ($C = 7500,$ $dt = 0.01$)\label{fig:ito_2p5_small_3}}
    \end{subfigure}
    \hfill
    \begin{subfigure}[b]{.45\textwidth}
        \includegraphics[width=\textwidth]{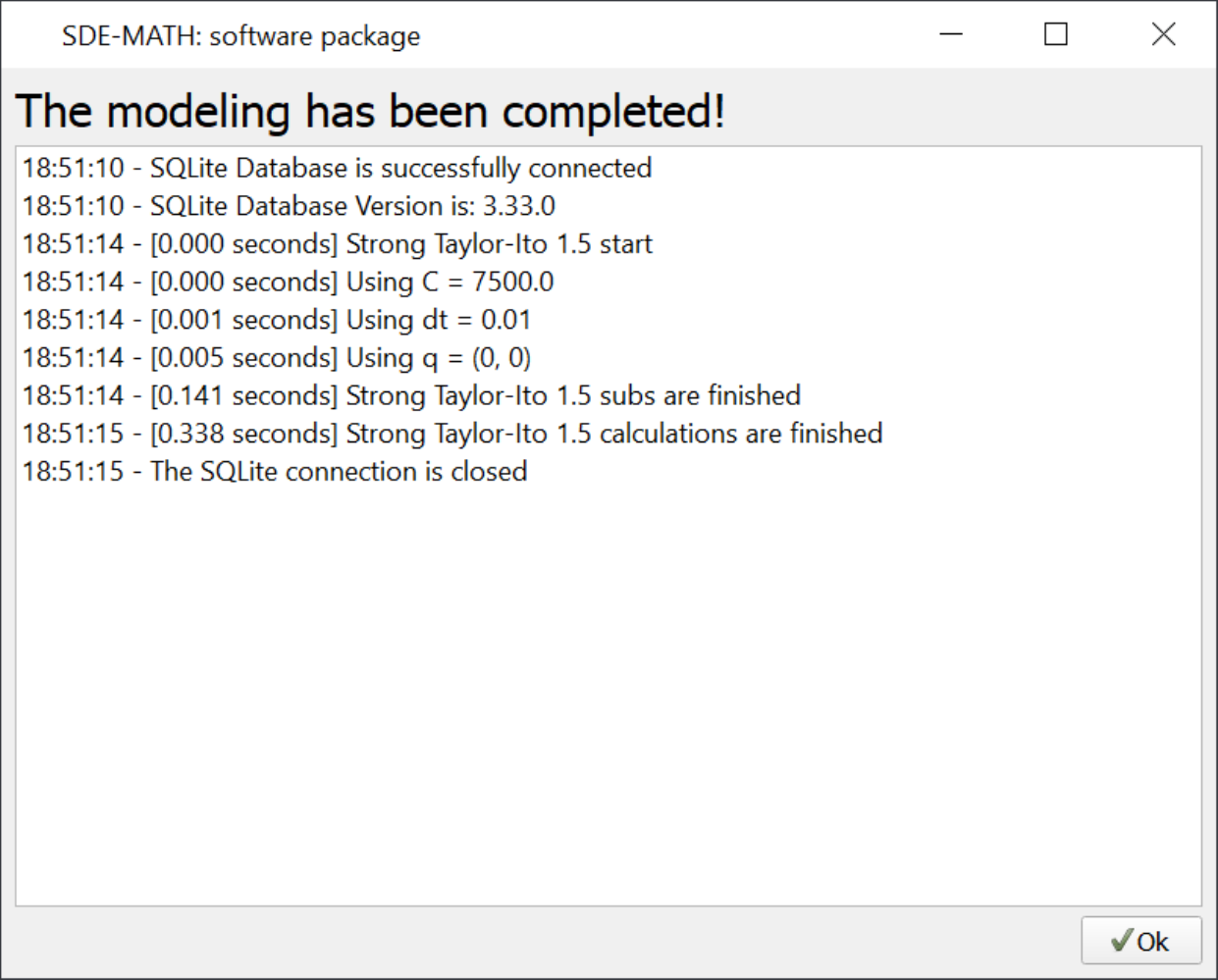}
        \caption*{Strong Taylor--It\^o scheme of order 2.0 ($C = 7500,$ $dt = 0.01$)\label{fig:ito_2p5_small_4}}
    \end{subfigure}
    \hspace*{\fill}

    \vspace{2mm}
    \hspace*{\fill}
    \begin{subfigure}[b]{.45\textwidth}
        \includegraphics[width=\textwidth]{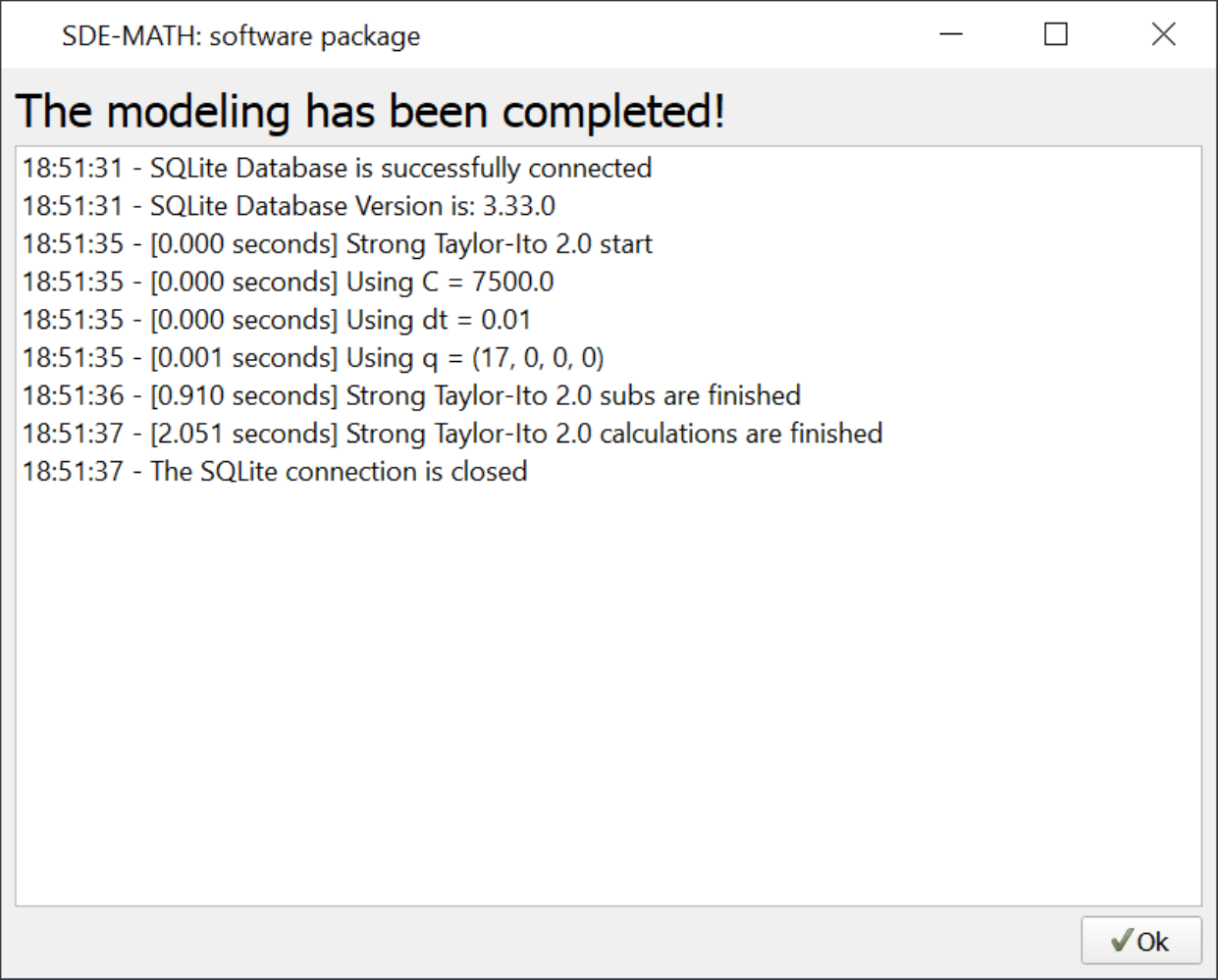}
        \caption*{Strong Taylor--It\^o scheme of order 2.5 ($C = 7500,$ $dt = 0.01$)\label{fig:ito_2p5_small_5}}
    \end{subfigure}
    \hspace*{\fill}

    \caption{Modeling logs\label{fig:ito_2p5_small_logs}}

\end{figure}

\begin{figure}[H]
    \centering
    \includegraphics[width=.9\textwidth]{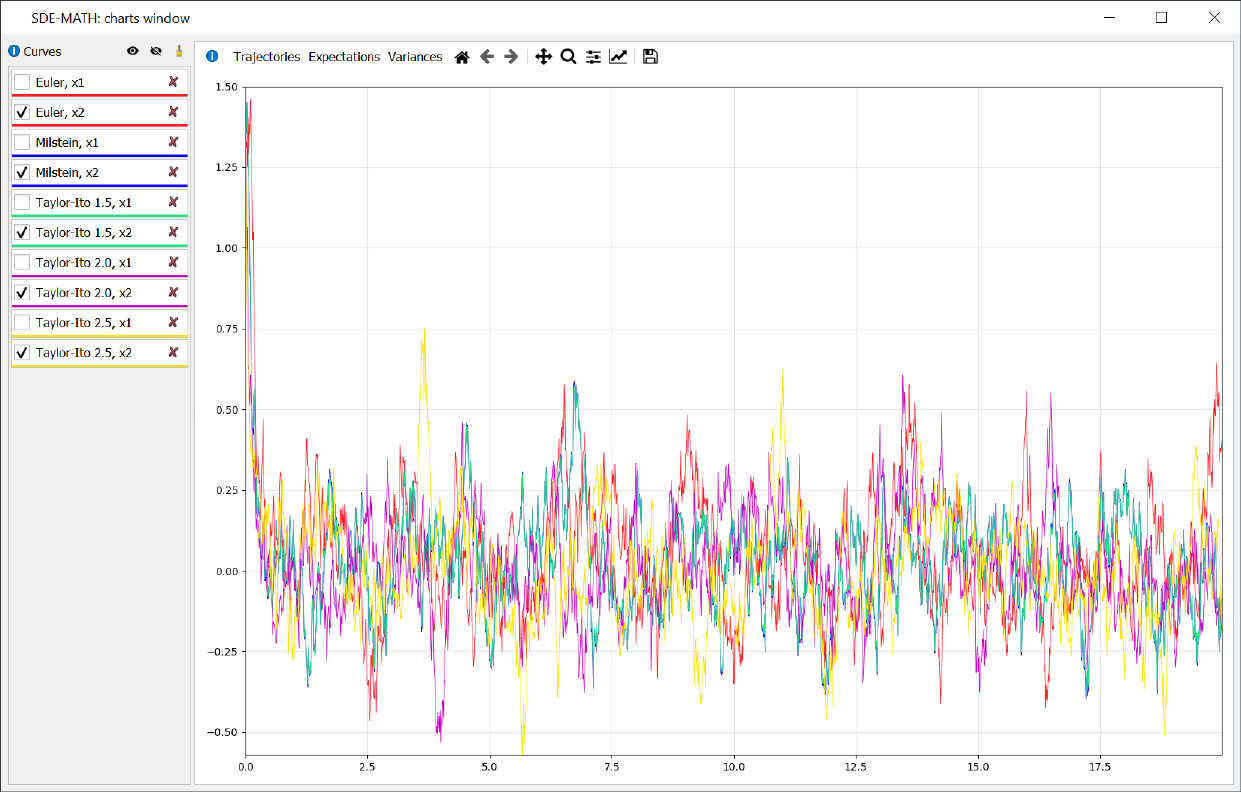}
    \caption{Strong Taylor--It\^o schemes of orders 0.5, 1.0, 1.5, 2.0, and 2.5 (${\bf x}_t^{(2)}$ component, $C = 7500,$ $dt = 0.01$)\label{fig:ito_2p5_small_7}}
\end{figure}

\begin{figure}[H]
    \centering

    \hspace*{\fill}
    \begin{subfigure}[b]{.45\textwidth}
        \includegraphics[width=\textwidth]{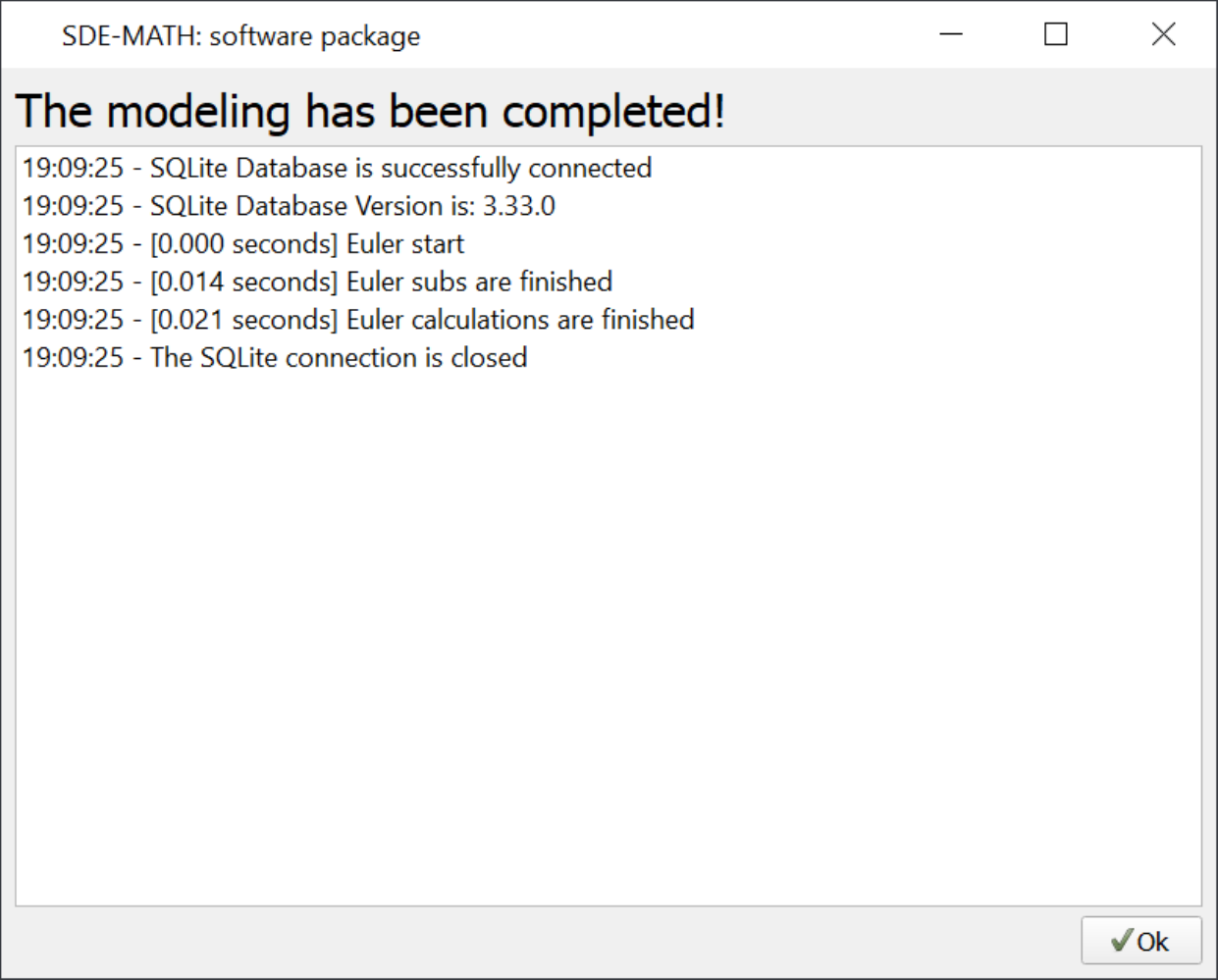}
        \caption*{Euler scheme ($dt = 0.025$)\label{fig:ito_3p0_small_1}}
    \end{subfigure}
    \hfill
    \begin{subfigure}[b]{.45\textwidth}
        \includegraphics[width=\textwidth]{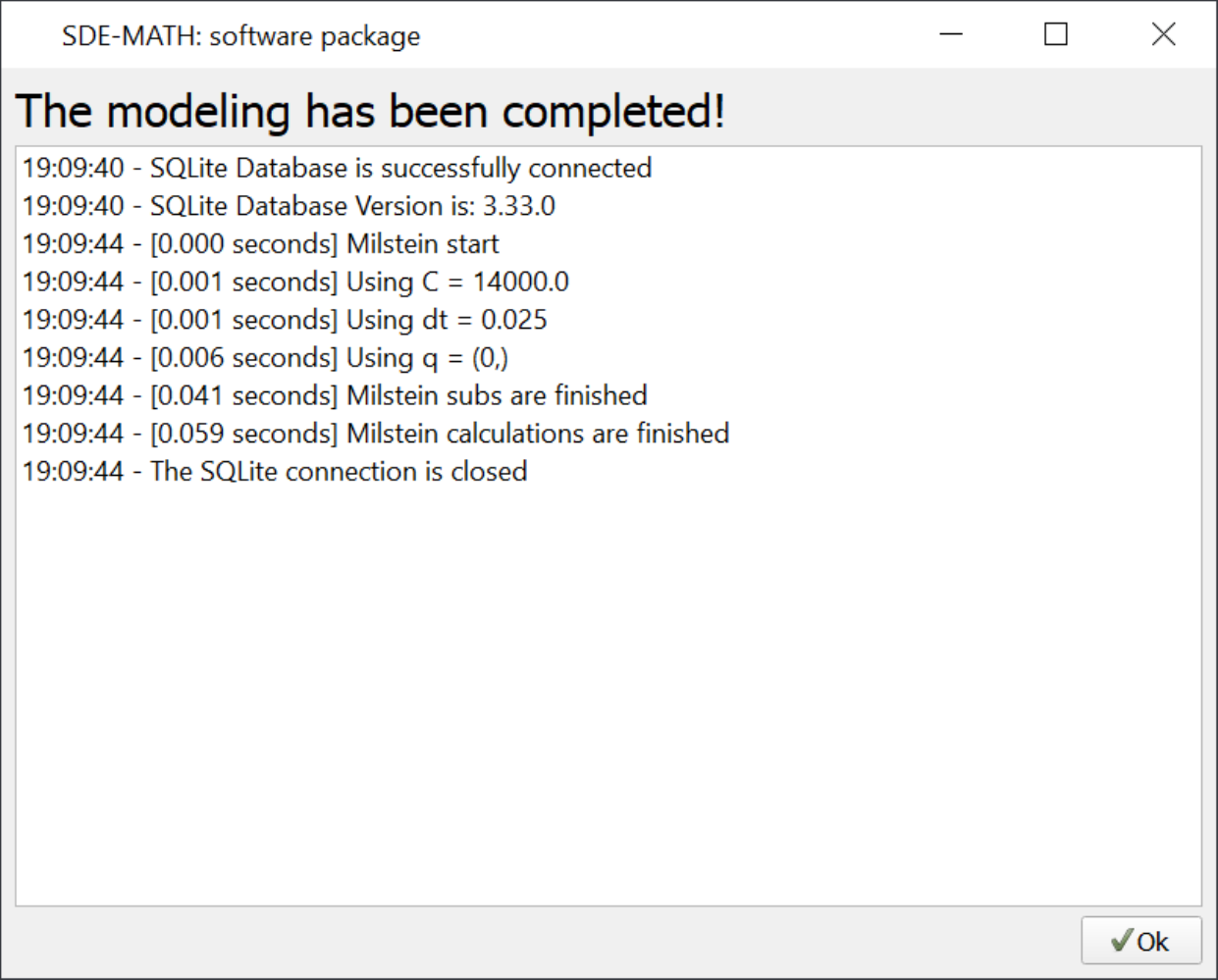}
        \caption*{Milstein scheme ($C = 14000,$ $dt = 0.025$)\label{fig:ito_3p0_small_2}}
    \end{subfigure}
    \hspace*{\fill}

    \vspace{2mm}
    \hspace*{\fill}
    \begin{subfigure}[b]{.45\textwidth}
        \includegraphics[width=\textwidth]{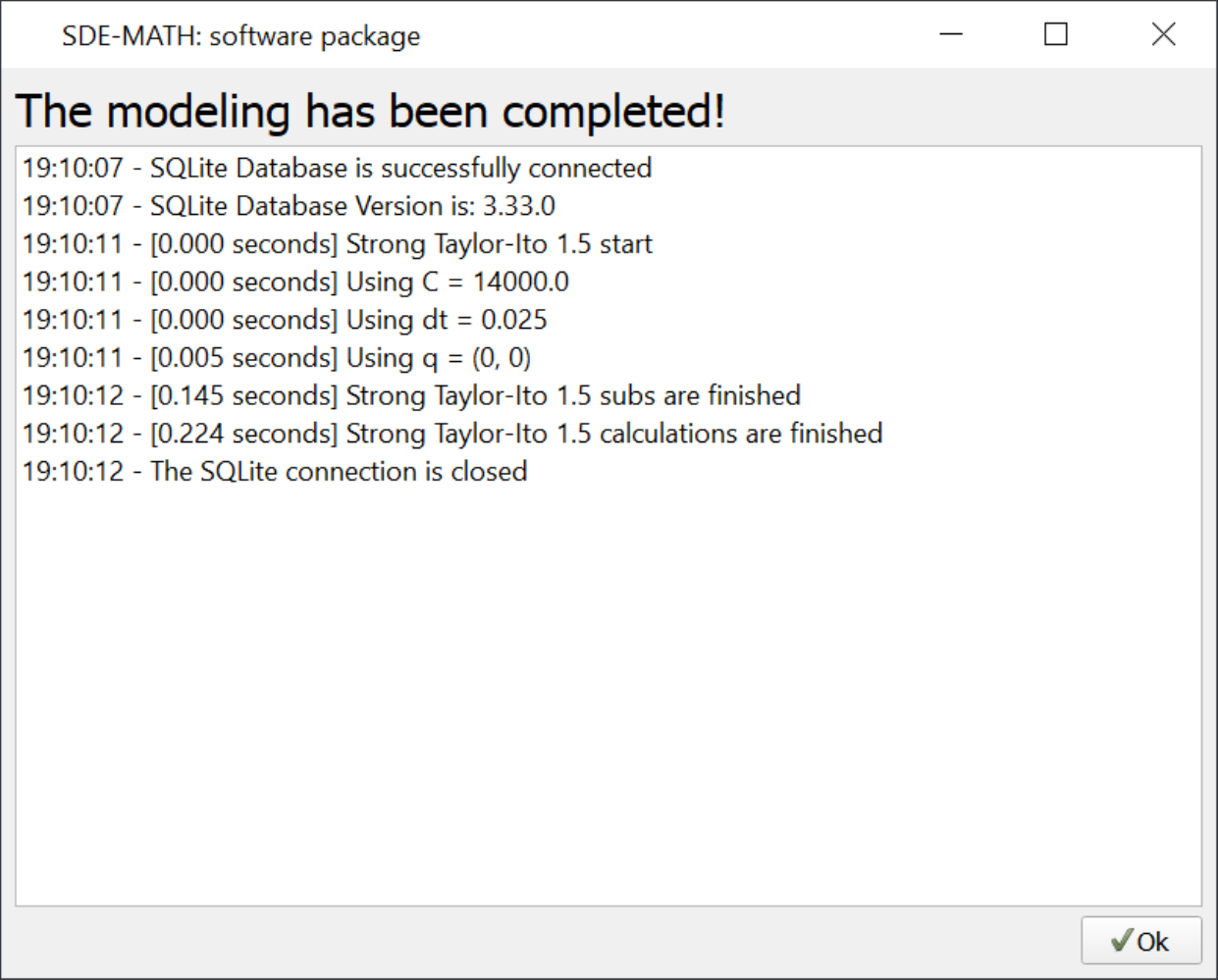}
        \caption*{Strong Taylor--It\^o scheme of order 1.5 ($C = 14000,$ $dt = 0.025$)\label{fig:ito_3p0_small_3}}
    \end{subfigure}
    \hfill
    \begin{subfigure}[b]{.45\textwidth}
        \includegraphics[width=\textwidth]{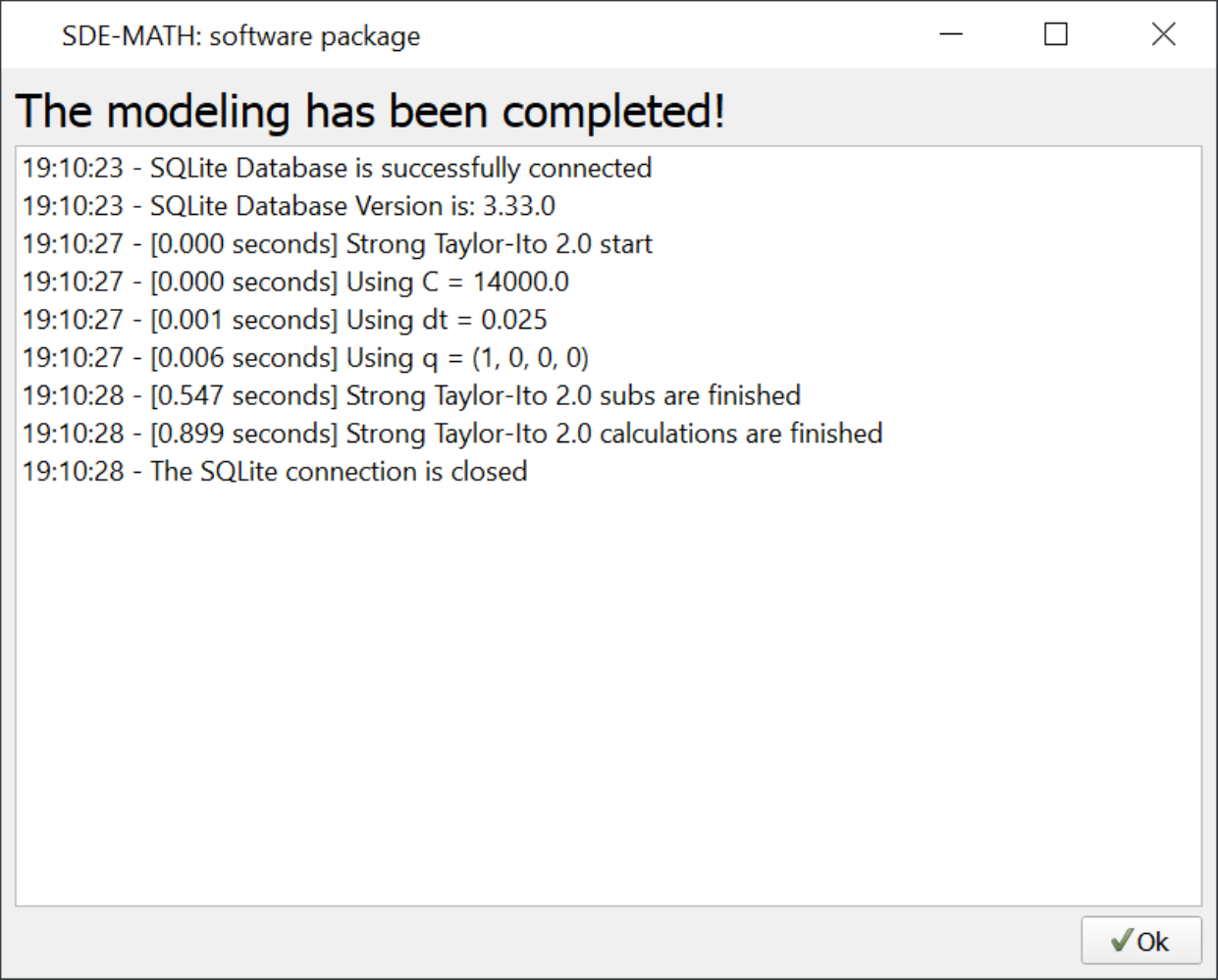}
        \caption*{Strong Taylor--It\^o scheme of order 2.0 ($C = 14000,$ $dt = 0.025$)\label{fig:ito_3p0_small_4}}
    \end{subfigure}
    \hspace*{\fill}

    \vspace{2mm}
    \hspace*{\fill}
    \begin{subfigure}[b]{.45\textwidth}
        \includegraphics[width=\textwidth]{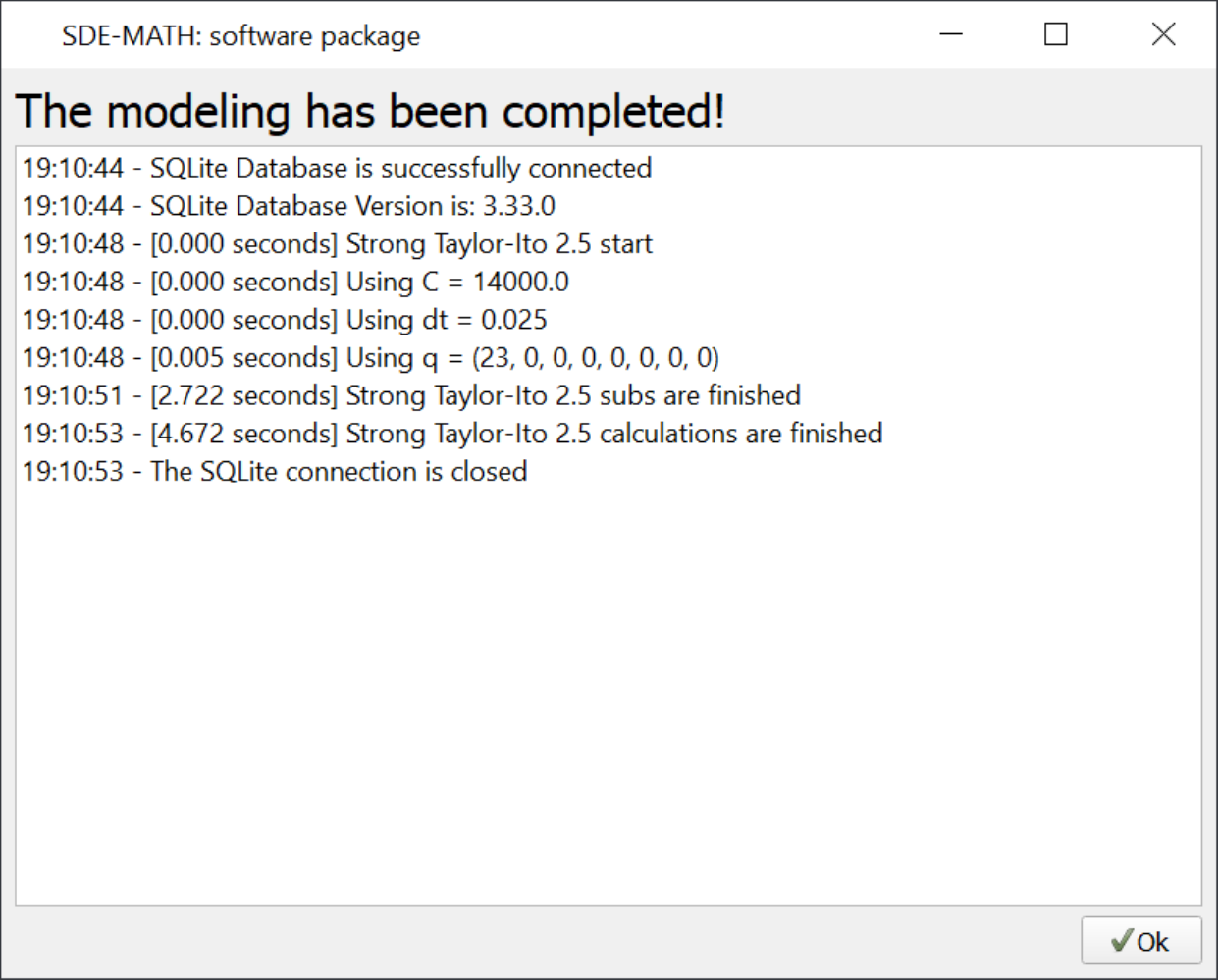}
        \caption*{Strong Taylor--It\^o scheme of order 2.5 ($C = 14000,$ $dt = 0.025$)\label{fig:ito_3p0_small_5}}
    \end{subfigure}
    \hfill
    \begin{subfigure}[b]{.45\textwidth}
        \includegraphics[width=\textwidth]{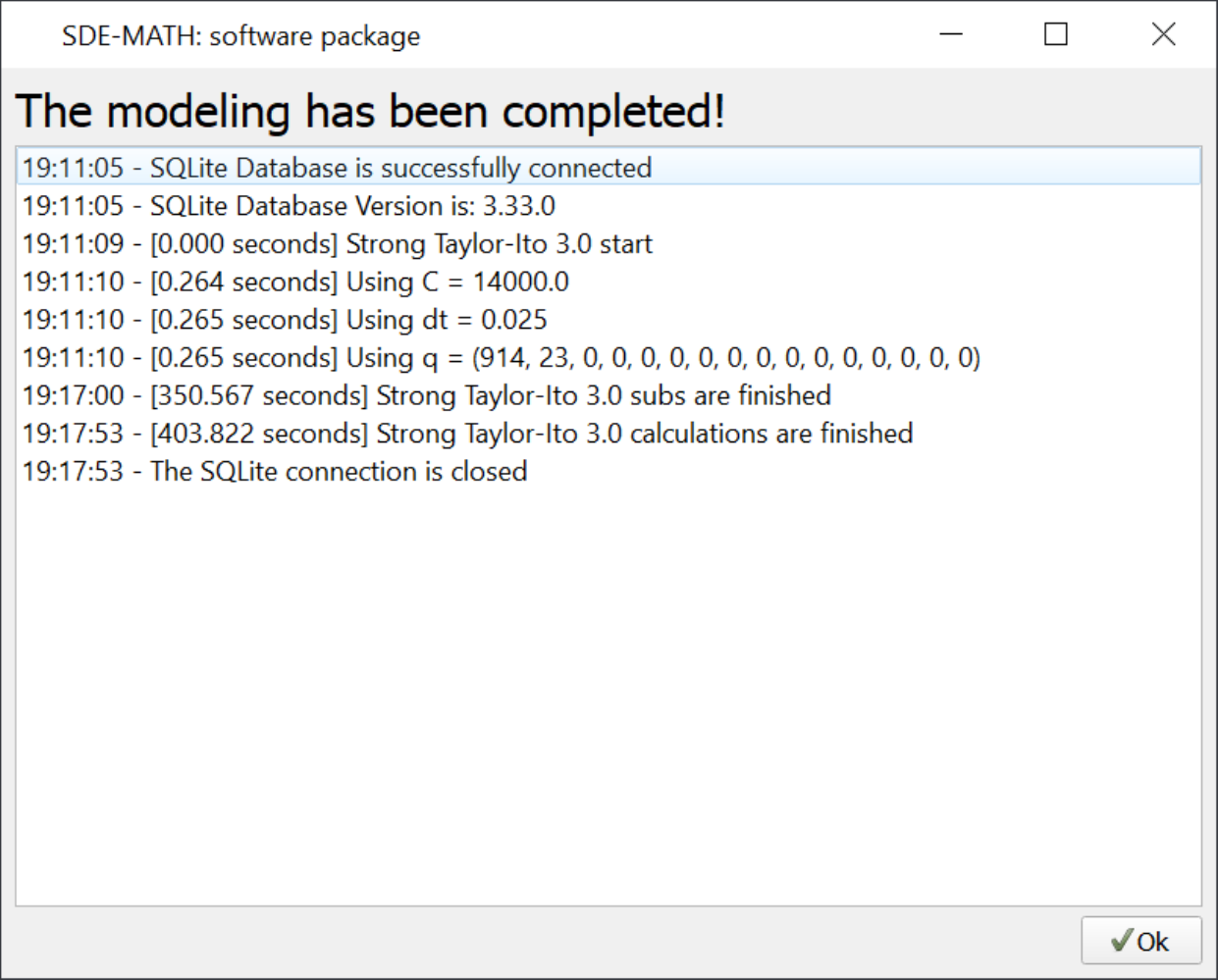}
        \caption*{Strong Taylor--It\^o scheme of order 3.0 ($C = 14000,$ $dt = 0.025$)\label{fig:ito_3p0_small_6}}
    \end{subfigure}
    \hspace*{\fill}

    \caption{Modeling logs\label{fig:ito_3p0_small_logs}}

\end{figure}

\begin{figure}[H]
    \vspace{7mm}
    \centering

    \hspace*{\fill}
    \begin{subfigure}[b]{.45\textwidth}
        \includegraphics[width=\textwidth]{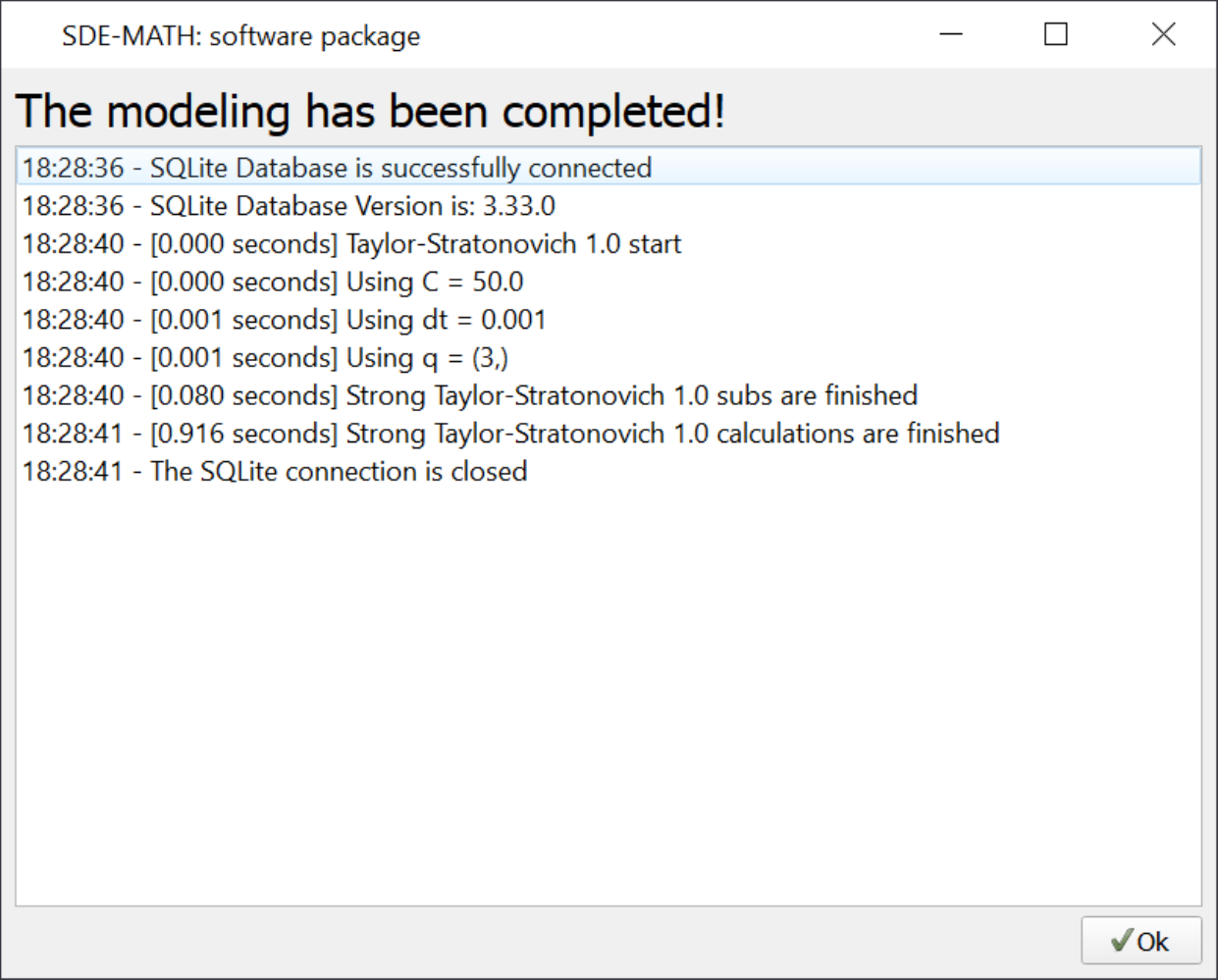}
        \caption*{Strong Taylor--Stratonovich scheme of order 1.0 ($C = 50,$ $dt = 0.001$)\label{fig:straton_1p5_small_1}}
    \end{subfigure}
    \hfill
    \begin{subfigure}[b]{.45\textwidth}
        \includegraphics[width=\textwidth]{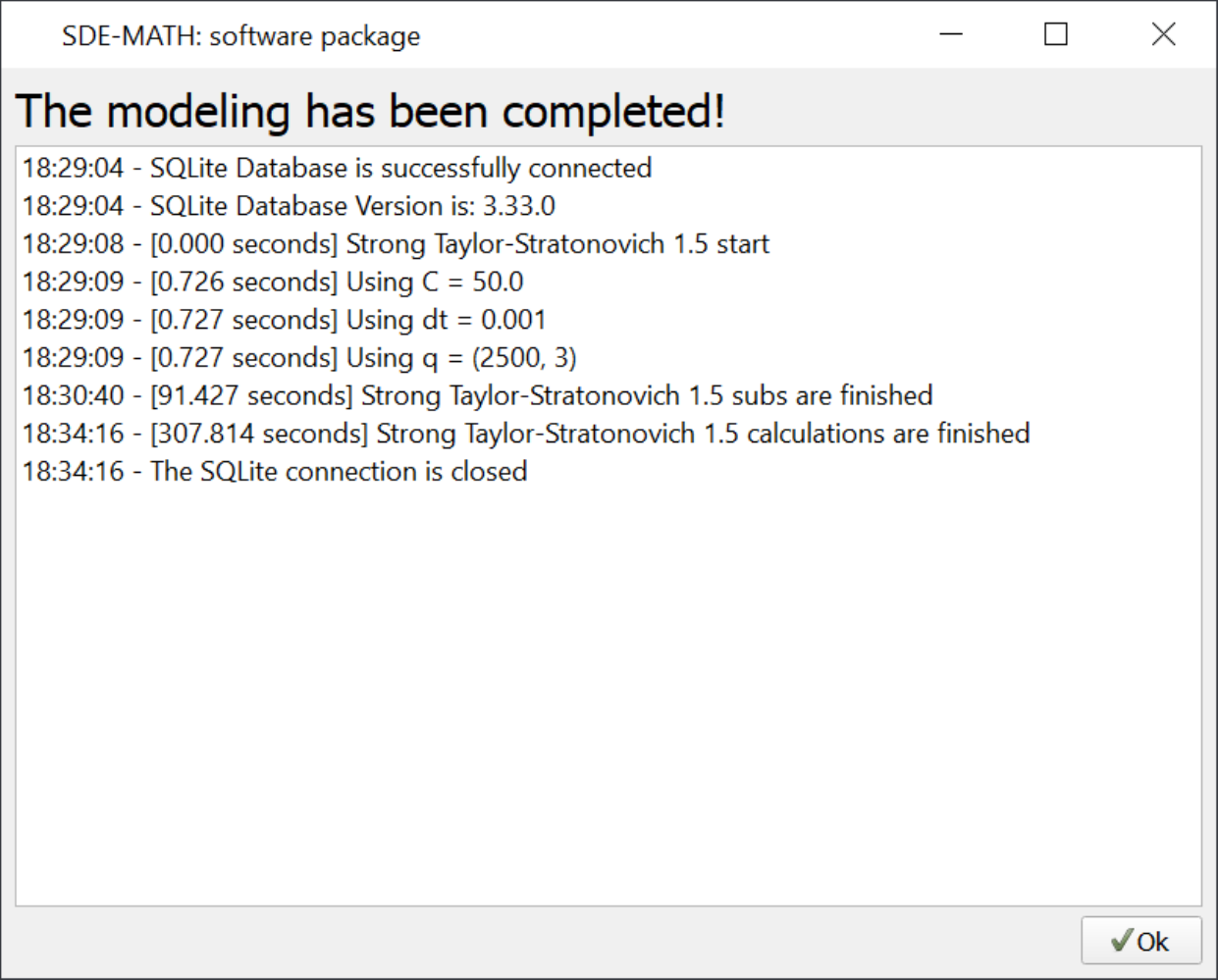}
        \caption*{Strong Taylor--Stratonovich scheme of order 1.5 ($C = 50,$ $dt = 0.001$)\label{fig:straton_1p5_small_2}}
    \end{subfigure}
    \hspace*{\fill}

    \caption{Modeling logs\label{fig:straton_1p5_small_logs}}

\end{figure}

\begin{figure}[H]
    \centering
    \includegraphics[width=.9\textwidth]{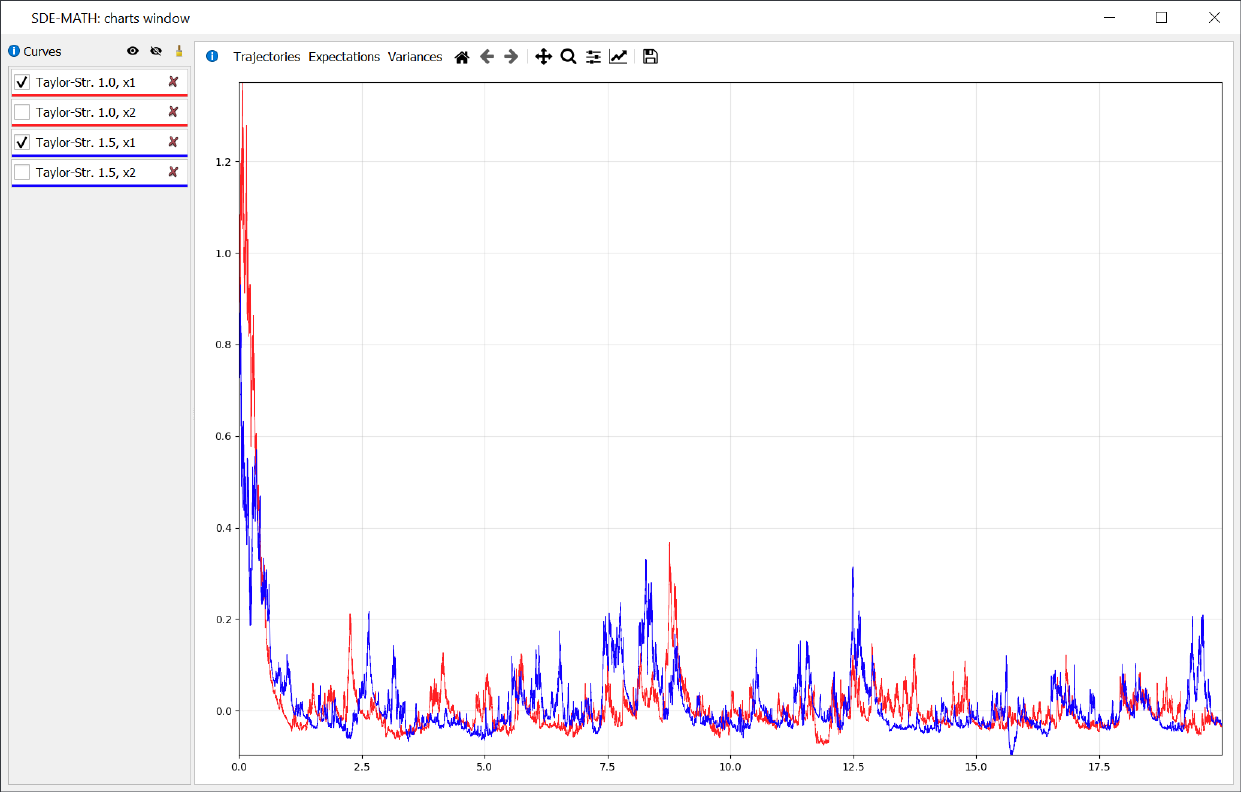}
    \caption{Strong Taylor--Stratonovich schemes of orders 1.0 and 1.5 (${\bf x}_t^{(1)}$ component, $C = 50,$ $dt = 0.001$)\label{fig:straton_1p5_small_3}}
\end{figure}

\begin{figure}[H]
    \centering
    \includegraphics[width=.9\textwidth]{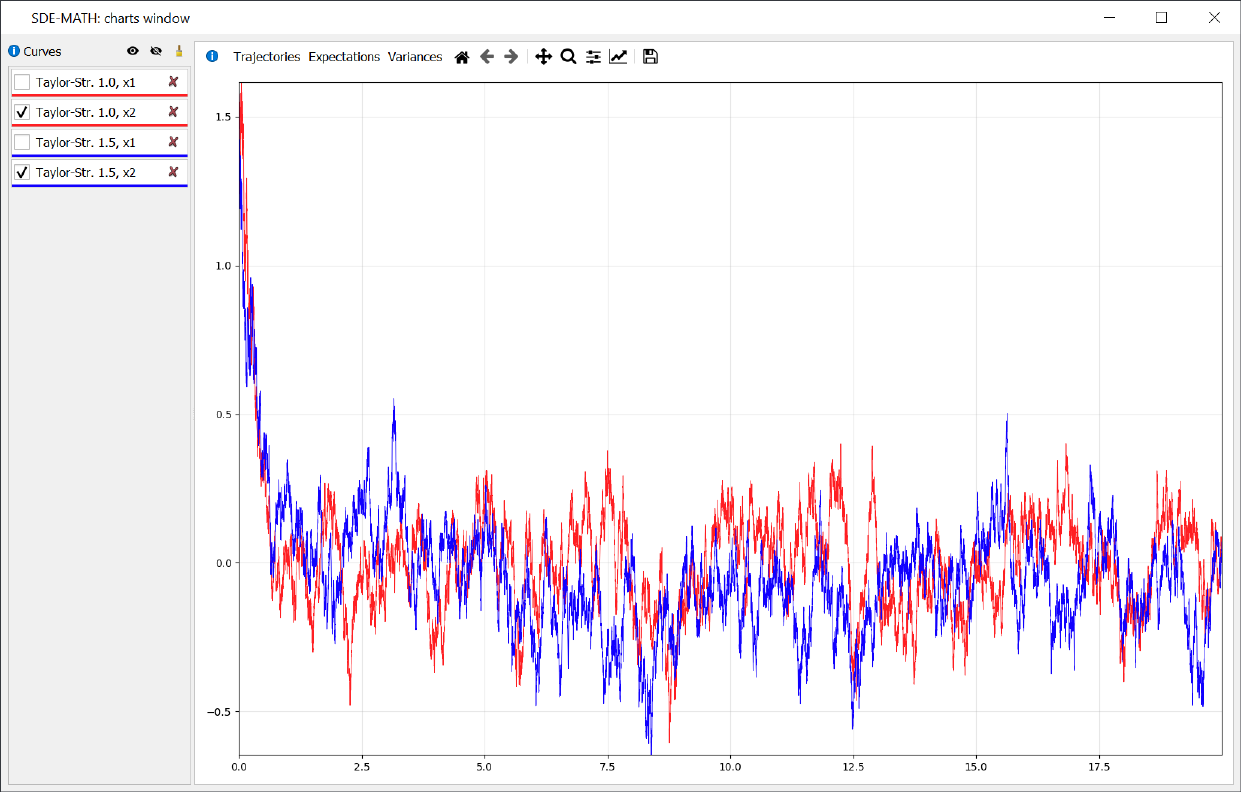}
    \caption{Strong Taylor--Stratonovich schemes of orders 1.0 and 1.5 (${\bf x}_t^{(2)}$ component, $C = 50,$ $dt = 0.001$)\label{fig:straton_1p5_small_4}}
\end{figure}

\begin{figure}[H]
    \vspace{10mm}
    \centering
    \includegraphics[width=.9\textwidth]{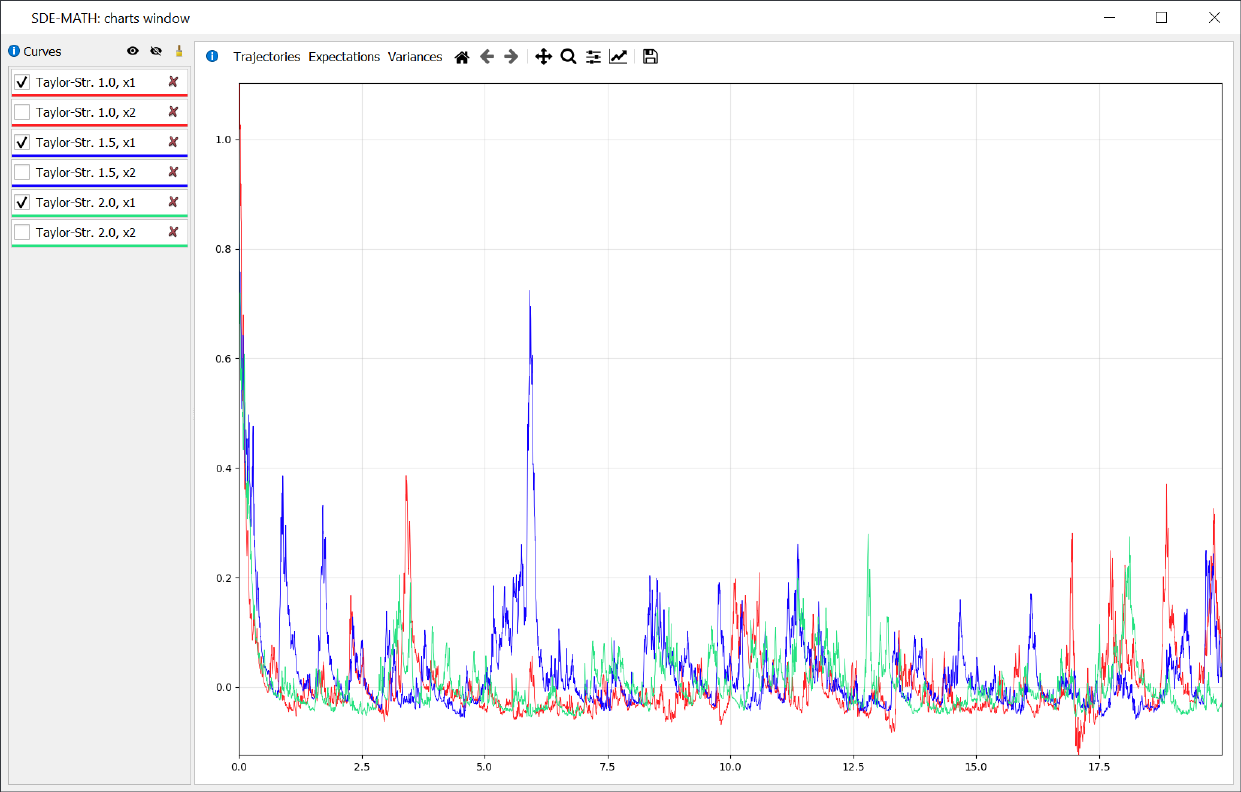}
    \caption{Strong Taylor--Stratonovich schemes of orders 1.0, 1.5, and 2.0 (${\bf x}_t^{(1)}$ component, $C = 500,$ $dt = 0.005$)\label{fig:straton_2p0_small_4}}
\end{figure}

\begin{figure}[H]
    \vspace{7mm}
    \centering

    \hspace*{\fill}
    \begin{subfigure}[b]{.45\textwidth}
        \includegraphics[width=\textwidth]{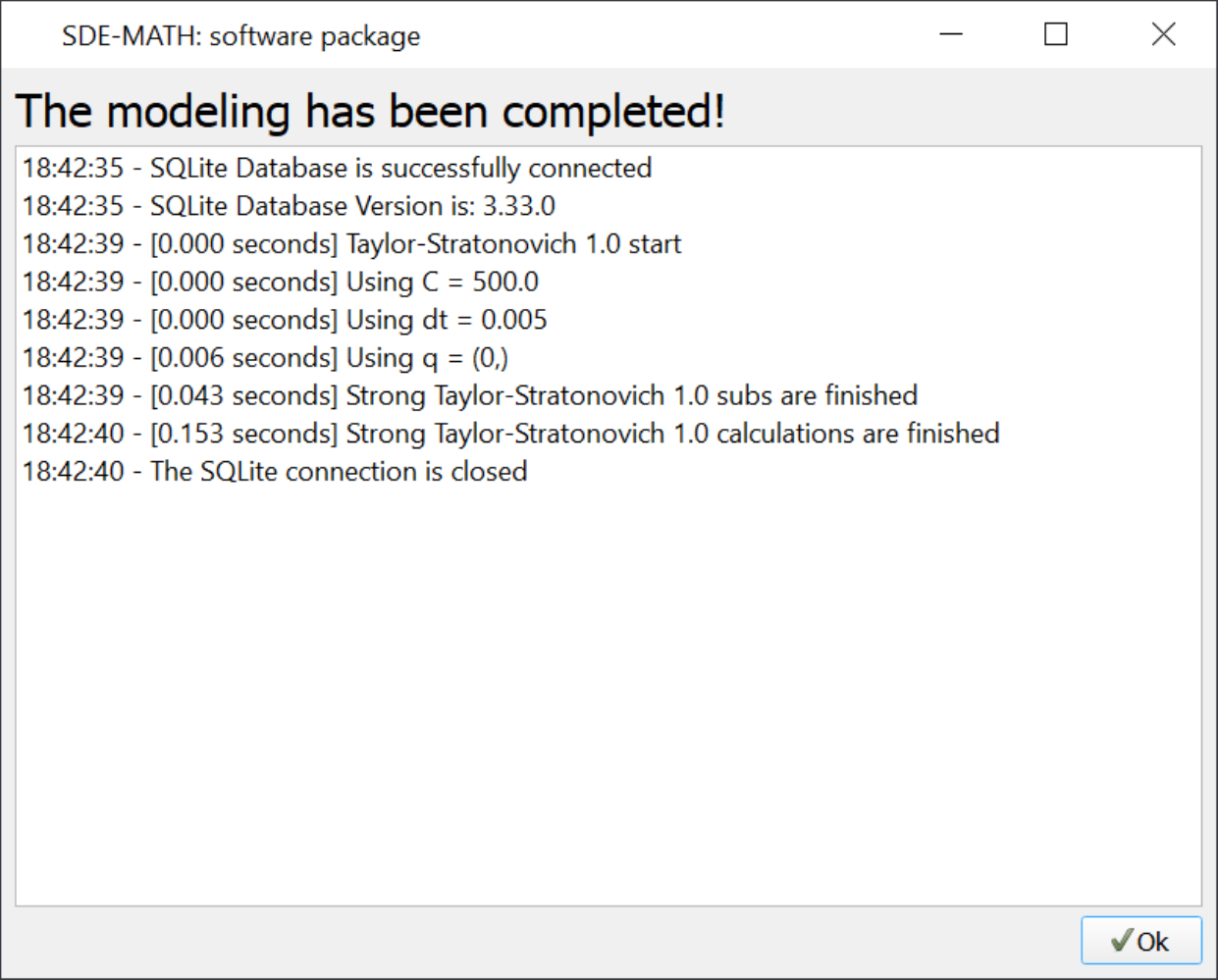}
        \caption*{Strong Taylor--Stratonovich scheme of order 1.0 ($C = 500,$ $dt = 0.005$)\label{fig:straton_2p0_small_1}}
    \end{subfigure}
    \hfill
    \begin{subfigure}[b]{.45\textwidth}
        \includegraphics[width=\textwidth]{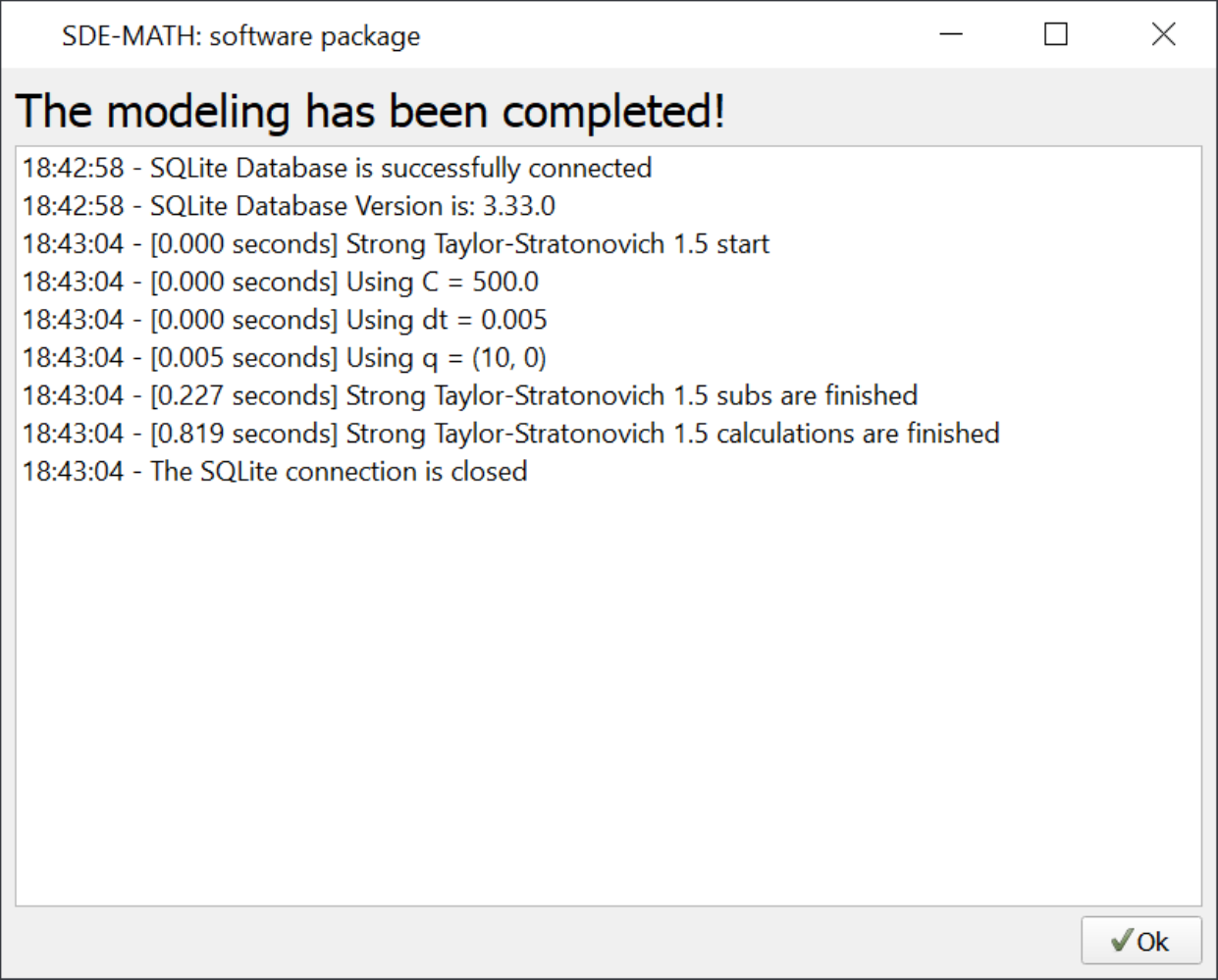}
        \caption*{Strong Taylor--Stratonovich scheme of order 1.5 ($C = 500,$ $dt = 0.005$)\label{fig:straton_2p0_small_2}}
    \end{subfigure}
    \hspace*{\fill}

    \caption{Modeling logs\label{fig:straton_2p0_small_logs}}

\end{figure}

\begin{figure}[H]
    \vspace{13mm}
    \centering
    \includegraphics[width=.45\textwidth]{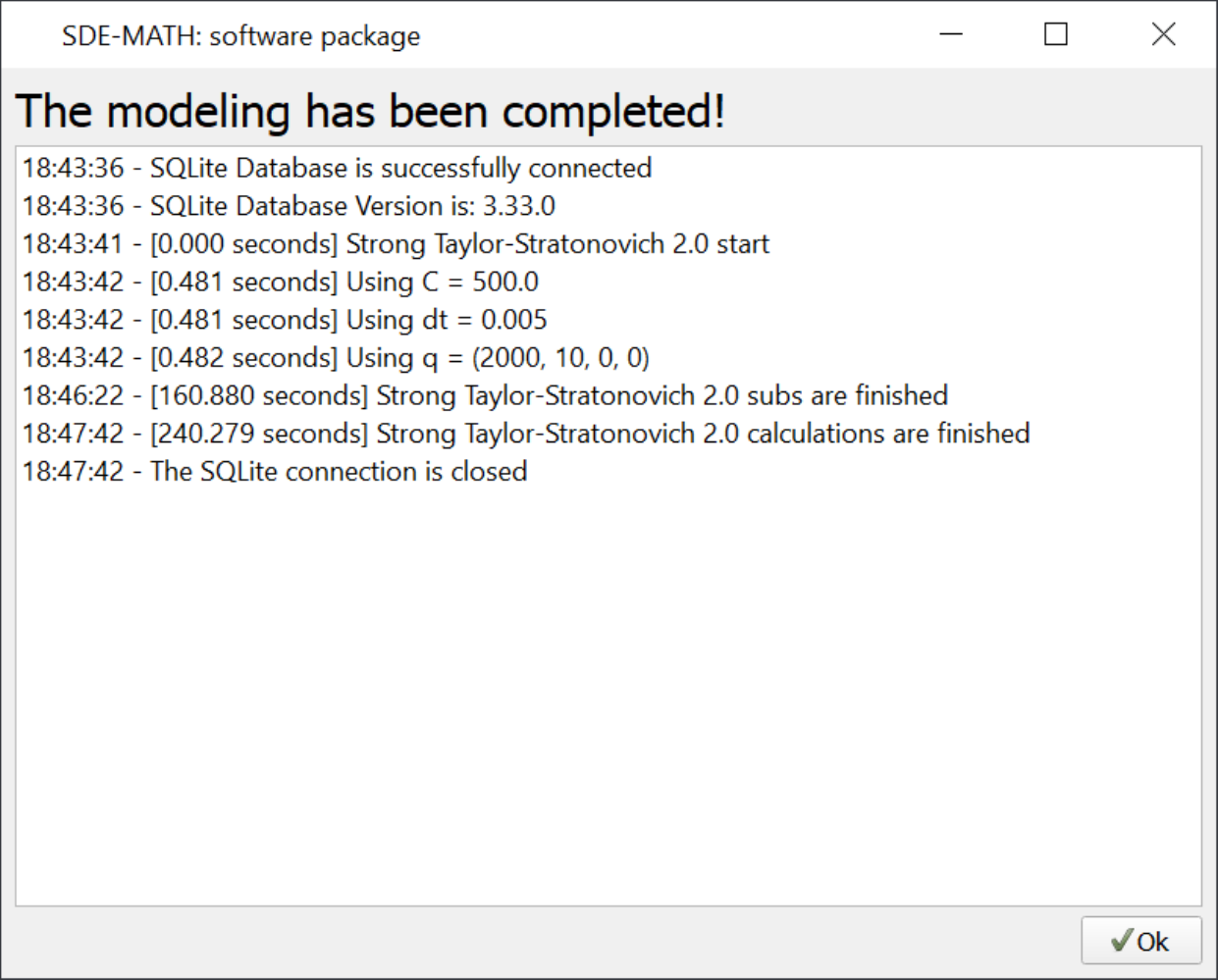}
    \caption{Strong Taylor--Stratonovich scheme of order 2.0 ($C = 500,$ $dt = 0.005$)\label{fig:straton_2p0_small_3}}
\end{figure}

\begin{figure}[H]
    \vspace{10mm}
    \centering
    \includegraphics[width=.9\textwidth]{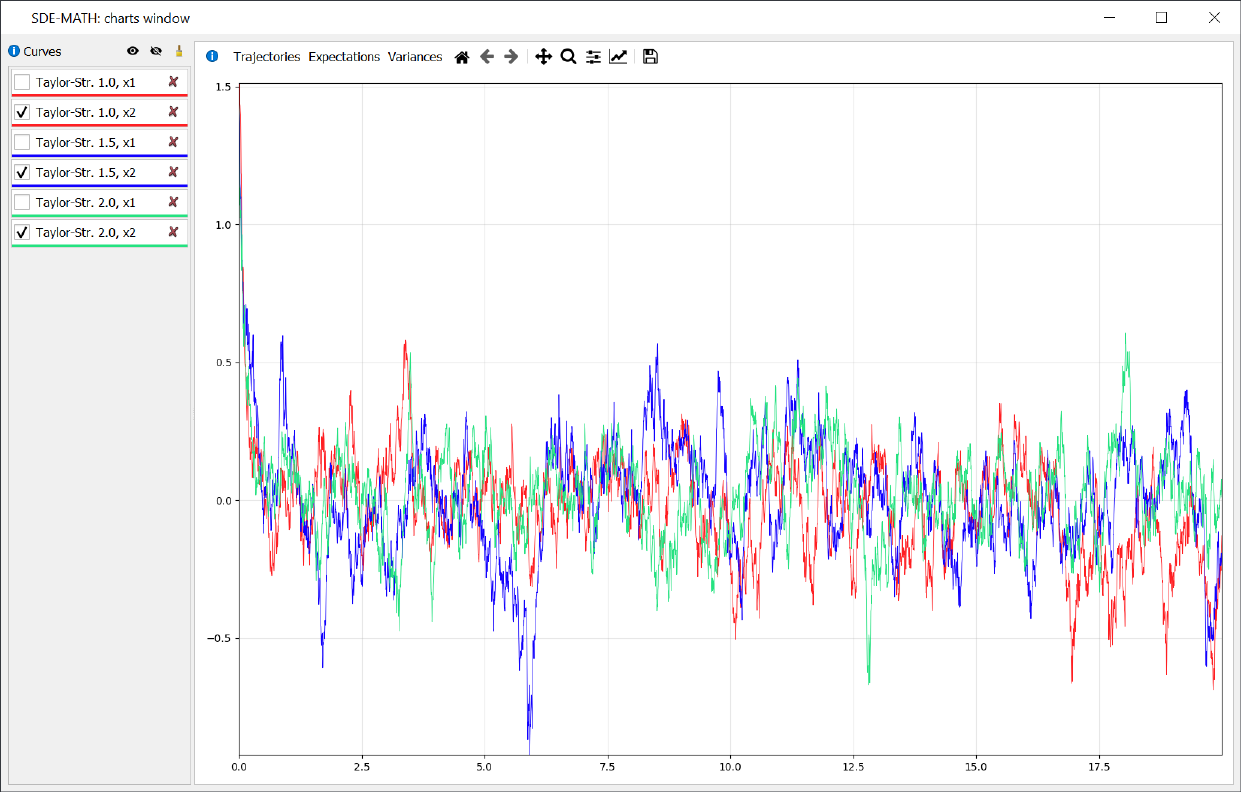}
    \caption{Strong Taylor--Stratonovich schemes of orders 1.0, 1.5, and 2.0 (${\bf x}_t^{(2)}$ component, $C = 500,$ $dt = 0.005$)\label{fig:straton_2p0_small_5}}
\end{figure}

\begin{figure}[H]
    \vspace{10mm}
    \centering
    \includegraphics[width=.9\textwidth]{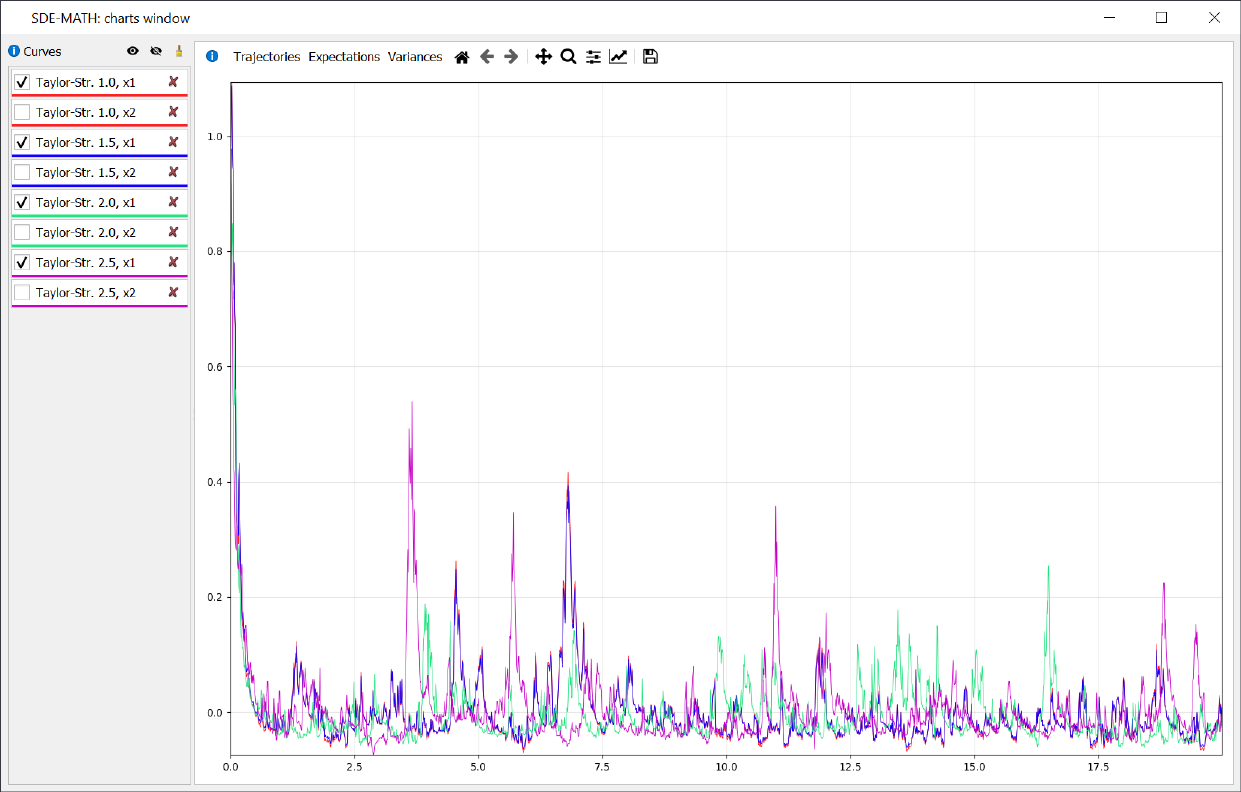}
    \caption{Strong Taylor--Stratonovich schemes of orders 1.0, 1.5, 2.0, and 2.5 (${\bf x}_t^{(1)}$ component, $C = 7500,$ $dt = 0.01$)\label{fig:straton_2p5_small_5}}
\end{figure}

\begin{figure}[H]
    \vspace{7mm}
    \centering

    \hspace*{\fill}
    \begin{subfigure}[b]{.45\textwidth}
        \includegraphics[width=\textwidth]{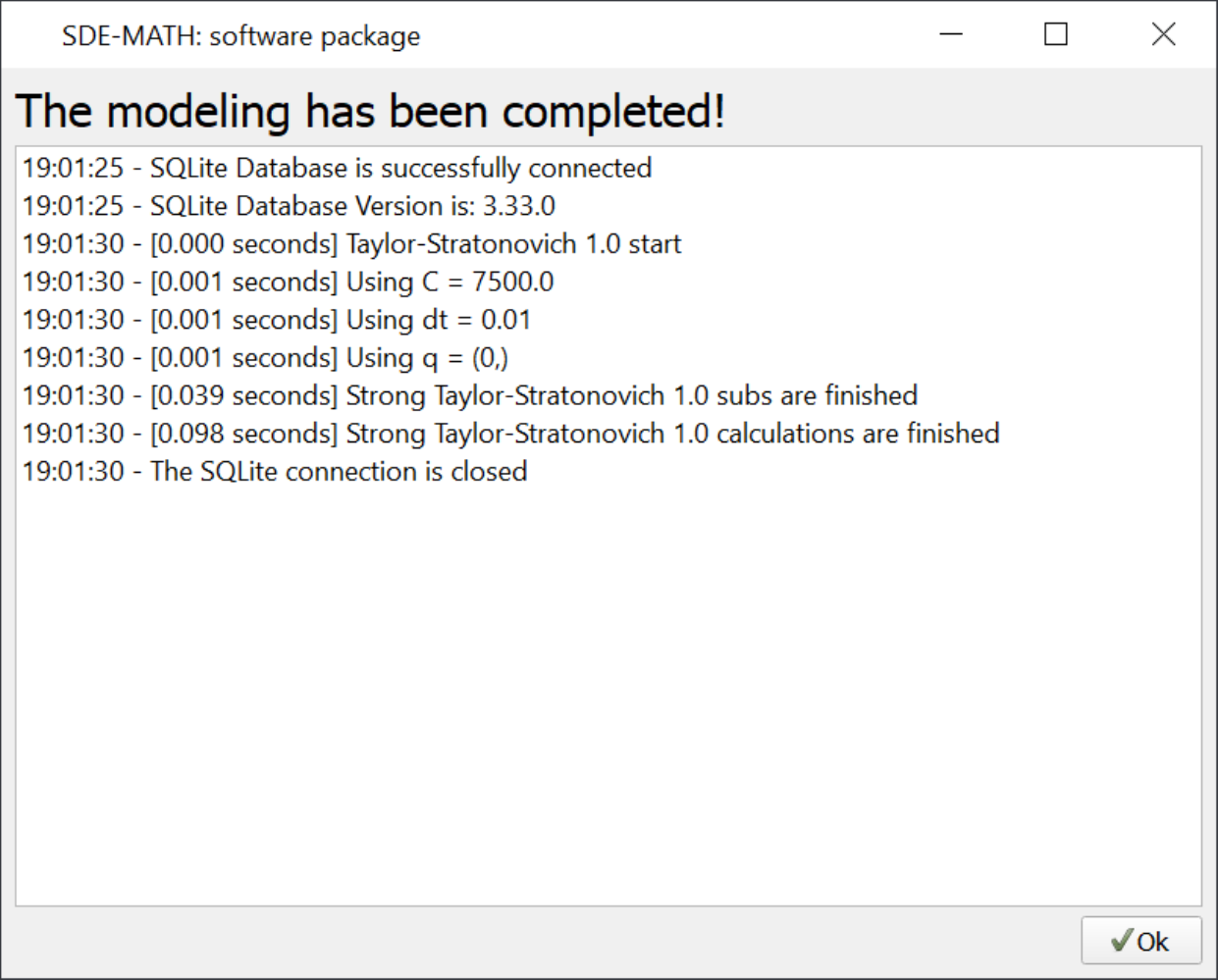}
        \caption*{Strong Taylor--Stratonovich scheme of order 1.0 ($C = 7500,$ $dt = 0.01$)\label{fig:straton_2p5_small_1}}
    \end{subfigure}
    \hfill
    \begin{subfigure}[b]{.45\textwidth}
        \includegraphics[width=\textwidth]{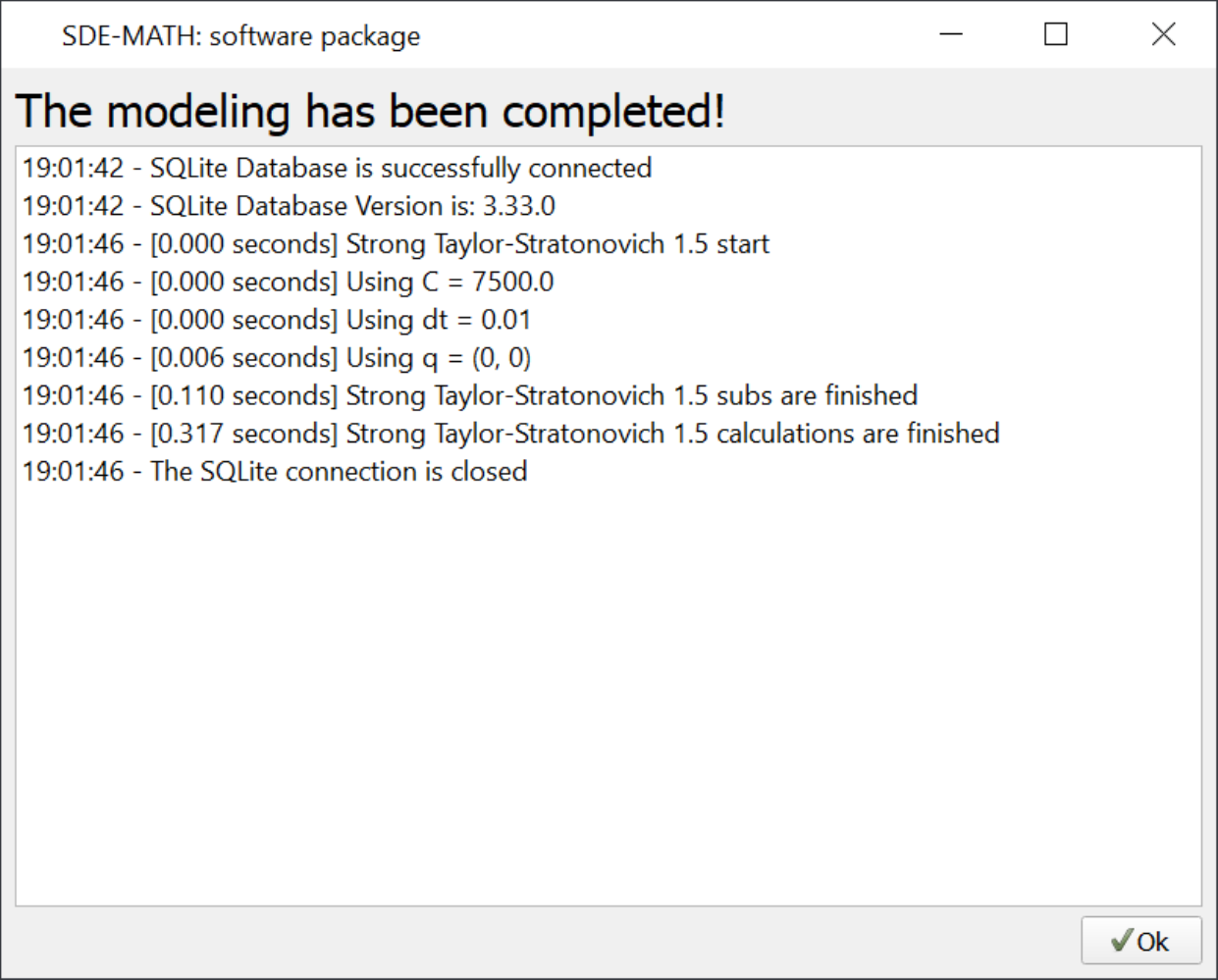}
        \caption*{Strong Taylor--Stratonovich scheme of order 1.5 ($C = 7500,$ $dt = 0.01$)\label{fig:straton_2p5_small_2}}
    \end{subfigure}
    \hspace*{\fill}

    \caption{Modeling logs\label{fig:straton_2p5_small_logs1}}

\end{figure}

\begin{figure}[H]
    \vspace{10mm}
    \centering
    \hspace*{\fill}
    \begin{subfigure}[b]{.45\textwidth}
        \includegraphics[width=\textwidth]{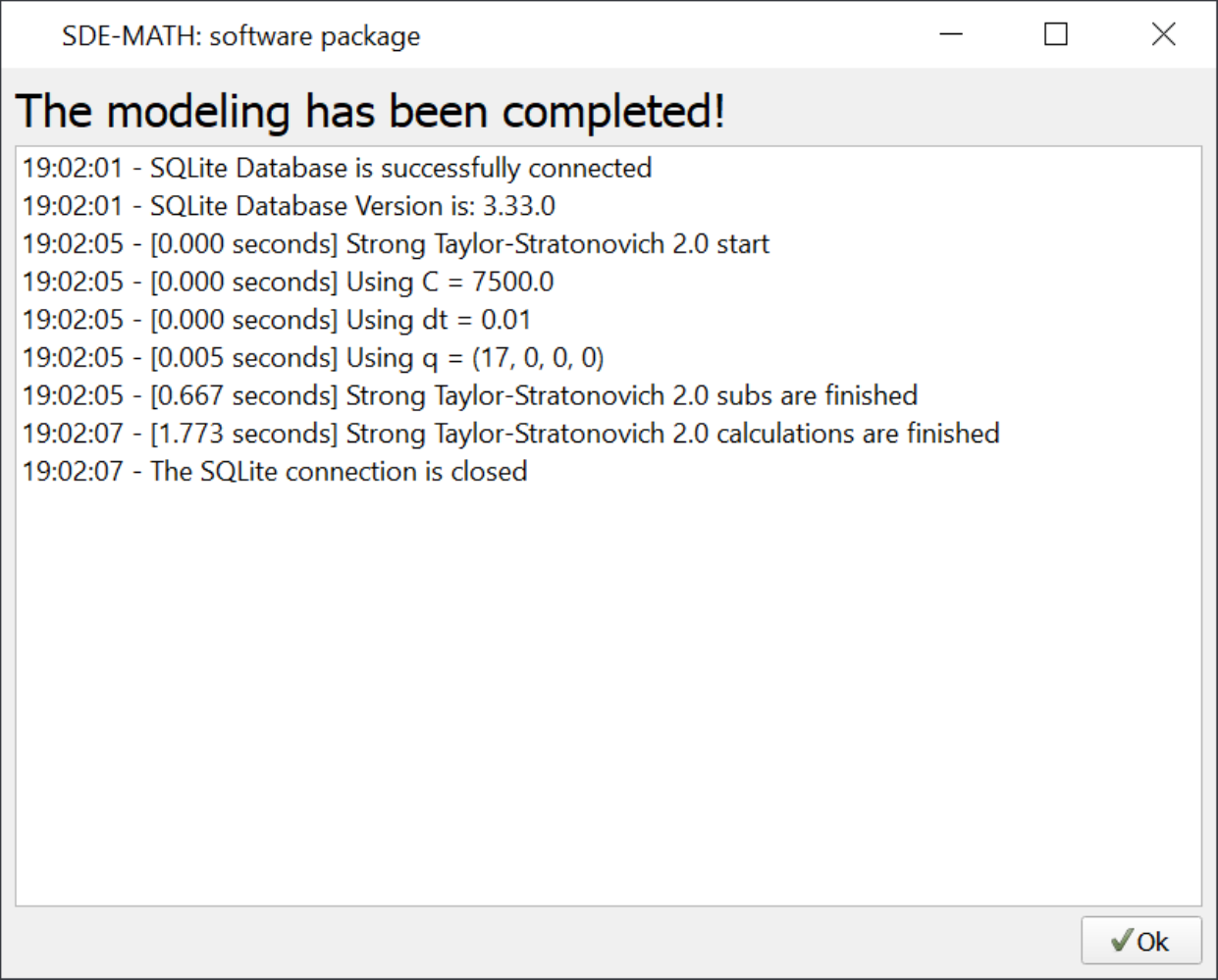}
        \caption*{Strong Taylor--Stratonovich scheme of order 2.0 ($C = 7500,$ $dt = 0.01$)\label{fig:straton_2p5_small_3}}
    \end{subfigure}
    \hfill
    \begin{subfigure}[b]{.45\textwidth}
        \includegraphics[width=\textwidth]{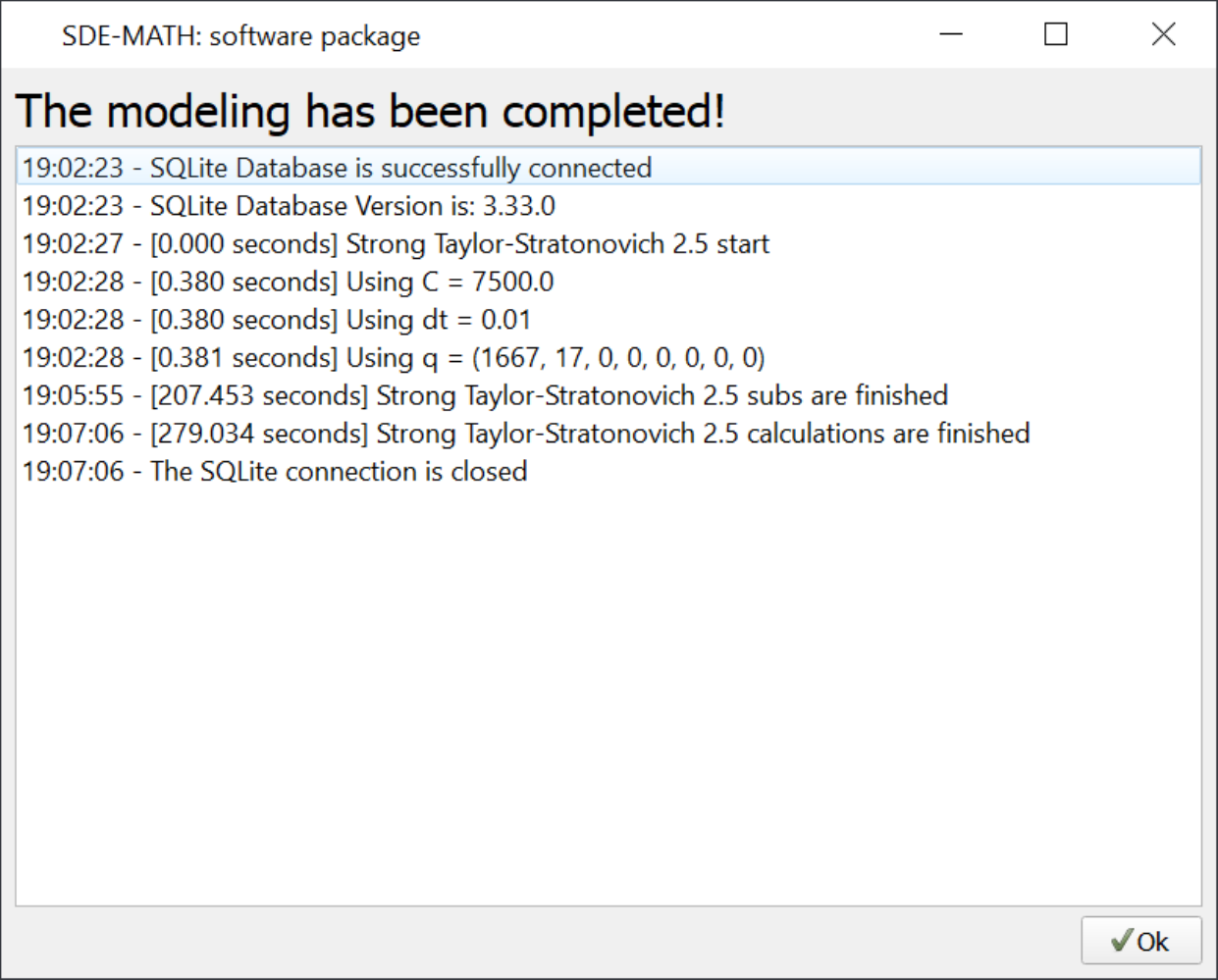}
        \caption*{Strong Taylor--Stratonovich scheme of order 2.5 ($C = 7500,$ $dt = 0.01$)\label{fig:straton_2p5_small_4}}
    \end{subfigure}
    \hspace*{\fill}

    \caption{Modeling logs\label{fig:straton_2p5_small_logs2}}

\end{figure}

\begin{figure}[H]
    \vspace{8mm}
    \centering
    \includegraphics[width=.9\textwidth]{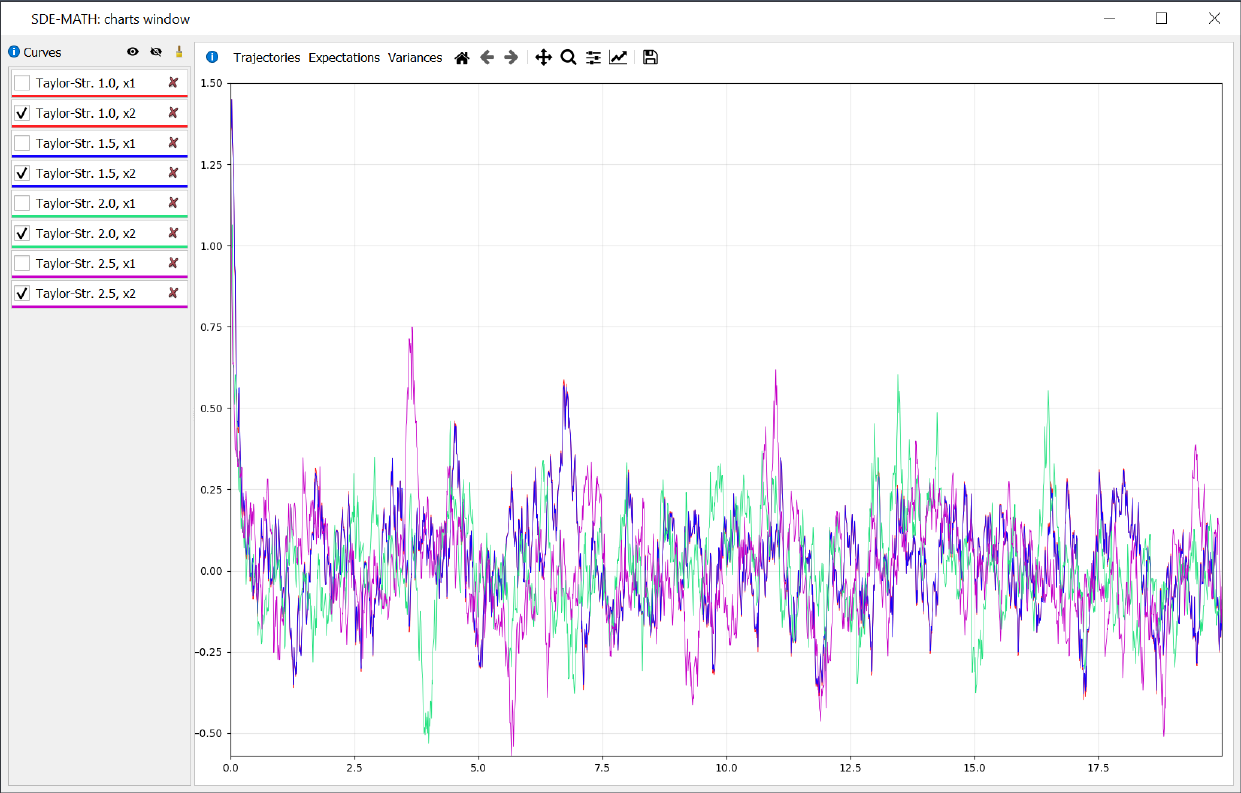}
    \caption{Strong Taylor--Stratonovich schemes of orders 1.0, 1.5, 2.0, and 2.5 (${\bf x}_t^{(2)}$ component, $C = 7500,$ $dt = 0.01$)\label{fig:straton_2p5_small_6}}
\end{figure}

\begin{figure}[H]
\centering

\hspace*{\fill}
\begin{subfigure}[b]{.45\textwidth}
    \includegraphics[width=\textwidth]{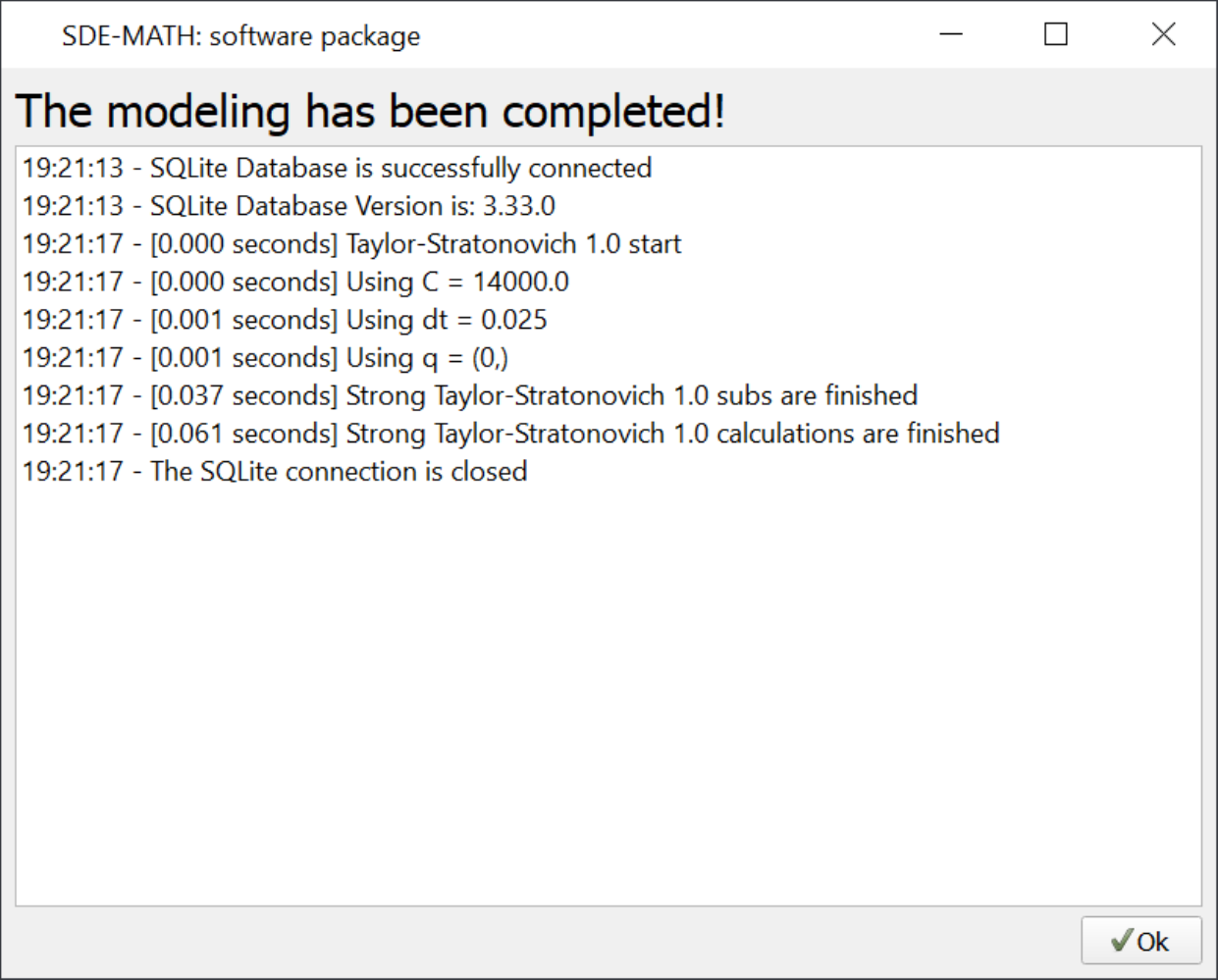}
    \caption*{Strong Taylor--Stratonovich scheme of order 1.0 ($C = 14000,$ $dt = 0.025$)\label{fig:straton_small_3p0_1}}
\end{subfigure}
\hfill
\begin{subfigure}[b]{.45\textwidth}
    \includegraphics[width=\textwidth]{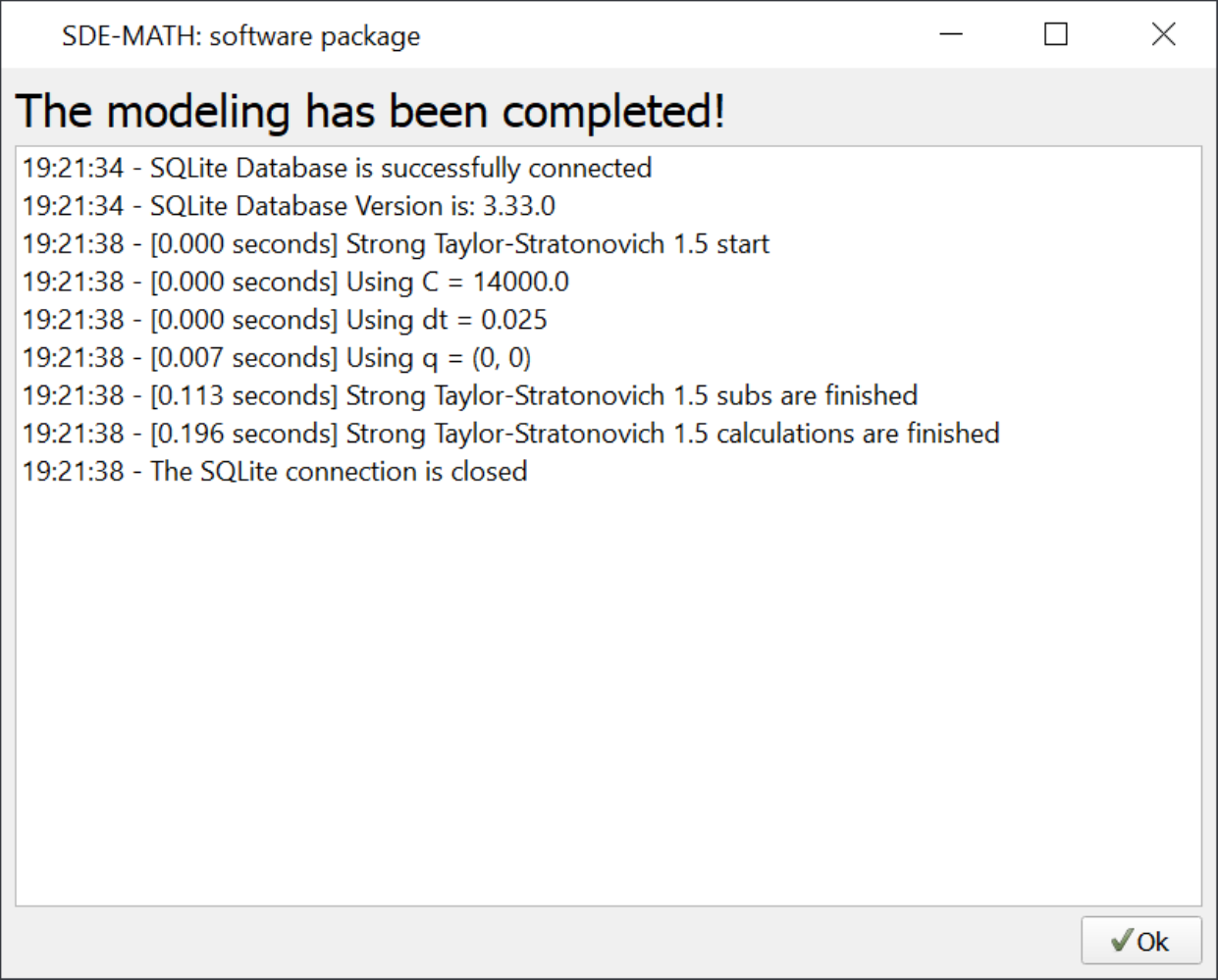}
    \caption*{Strong Taylor--Stratonovich scheme of order 1.5 ($C = 14000,$ $dt = 0.025$)\label{fig:straton_small_3p0_2}}
\end{subfigure}
\hspace*{\fill}

\vspace{2mm}
\hspace*{\fill}
\begin{subfigure}[b]{.45\textwidth}
    \includegraphics[width=\textwidth]{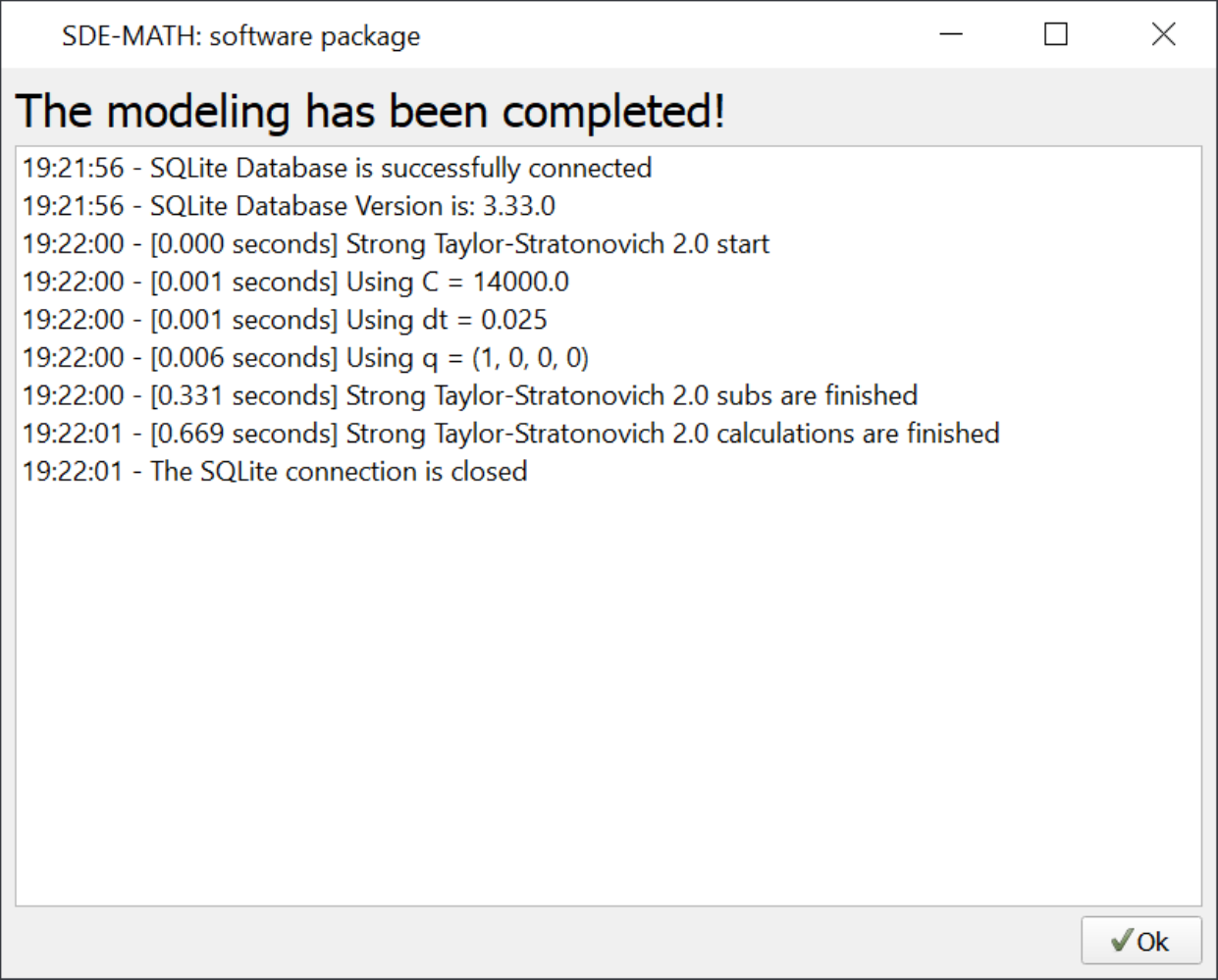}
    \caption*{Strong Taylor--Stratonovich scheme of order 2.0 ($C = 14000,$ $dt = 0.025$)\label{fig:straton_small_3p0_3}}
\end{subfigure}
\hfill
\begin{subfigure}[b]{.45\textwidth}
    \includegraphics[width=\textwidth]{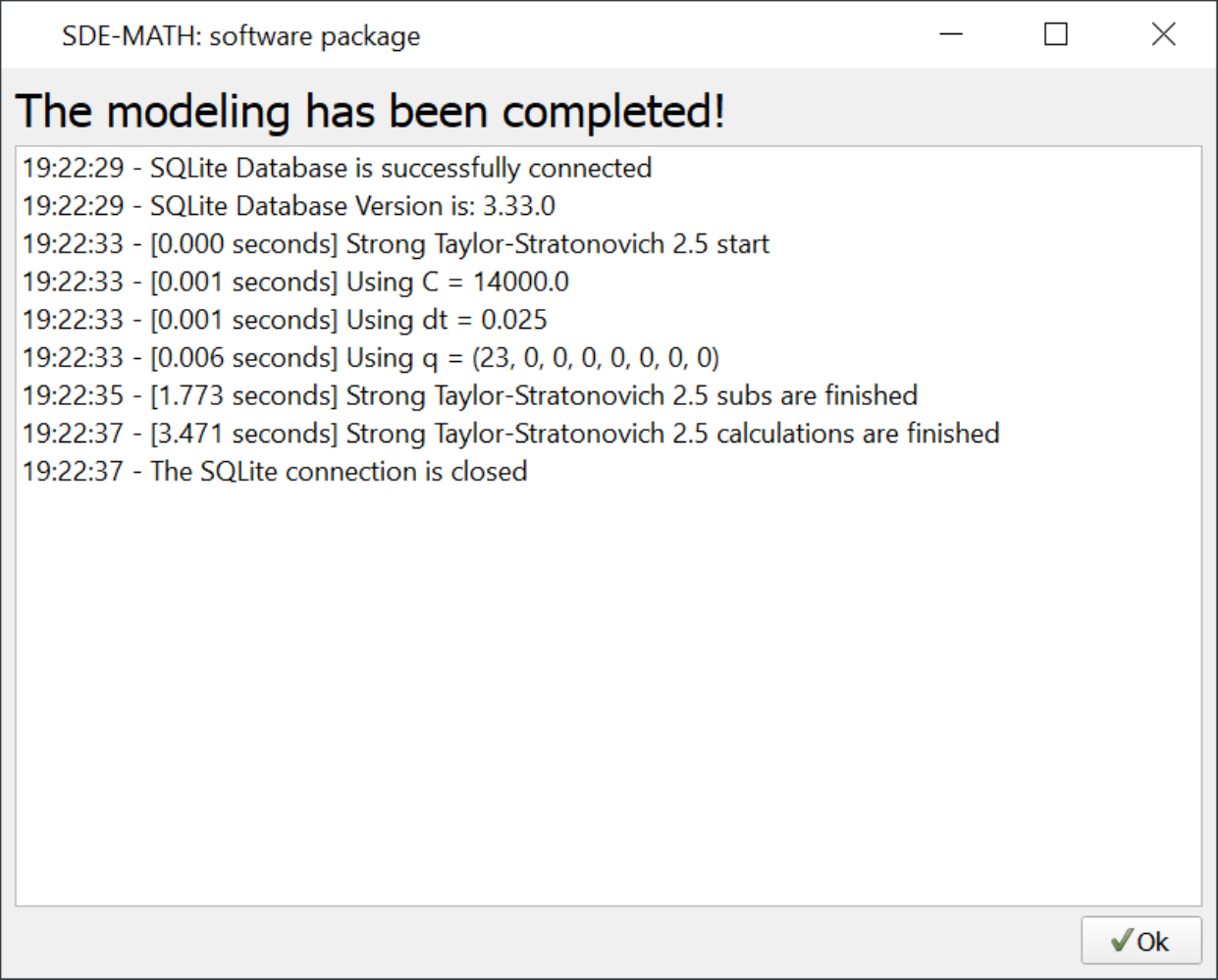}
    \caption*{Strong Taylor--Stratonovich scheme of order 2.5 ($C = 14000,$ $dt = 0.025$)\label{fig:straton_small_3p0_4}}
\end{subfigure}
\hspace*{\fill}

\vspace{2mm}
\hspace*{\fill}
\begin{subfigure}[b]{.45\textwidth}
    \includegraphics[width=\textwidth]{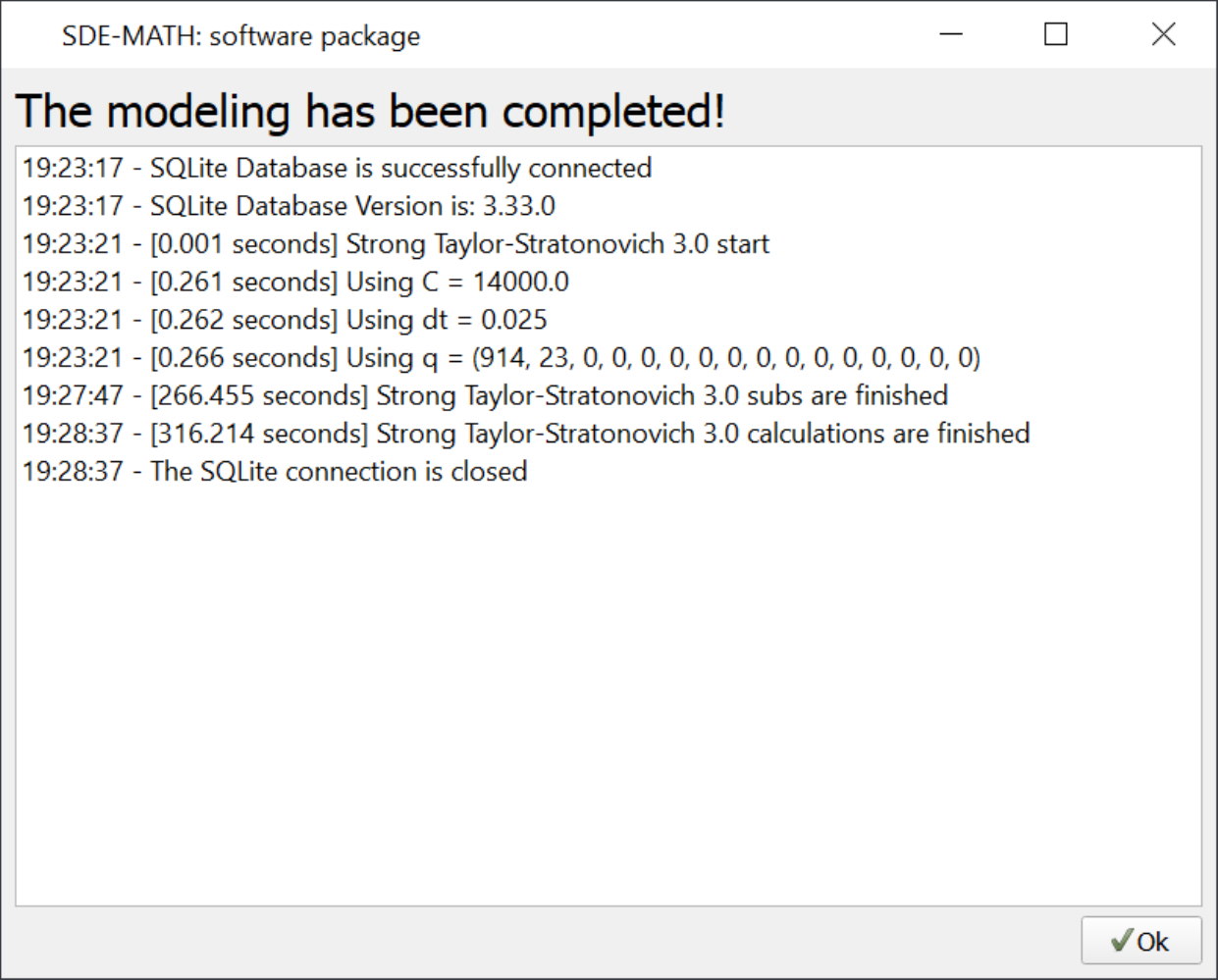}
    \caption*{Strong Taylor--Stratonovich scheme of order 3.0 ($C = 14000,$ $dt = 0.025$)\label{fig:straton_small_3p0_5}}
\end{subfigure}
\hspace*{\fill}

\caption{Modeling logs\label{fig:straton_small_3p0_logs}}

\end{figure}

\begin{figure}[H]
    \centering
    \includegraphics[width=.9\textwidth]{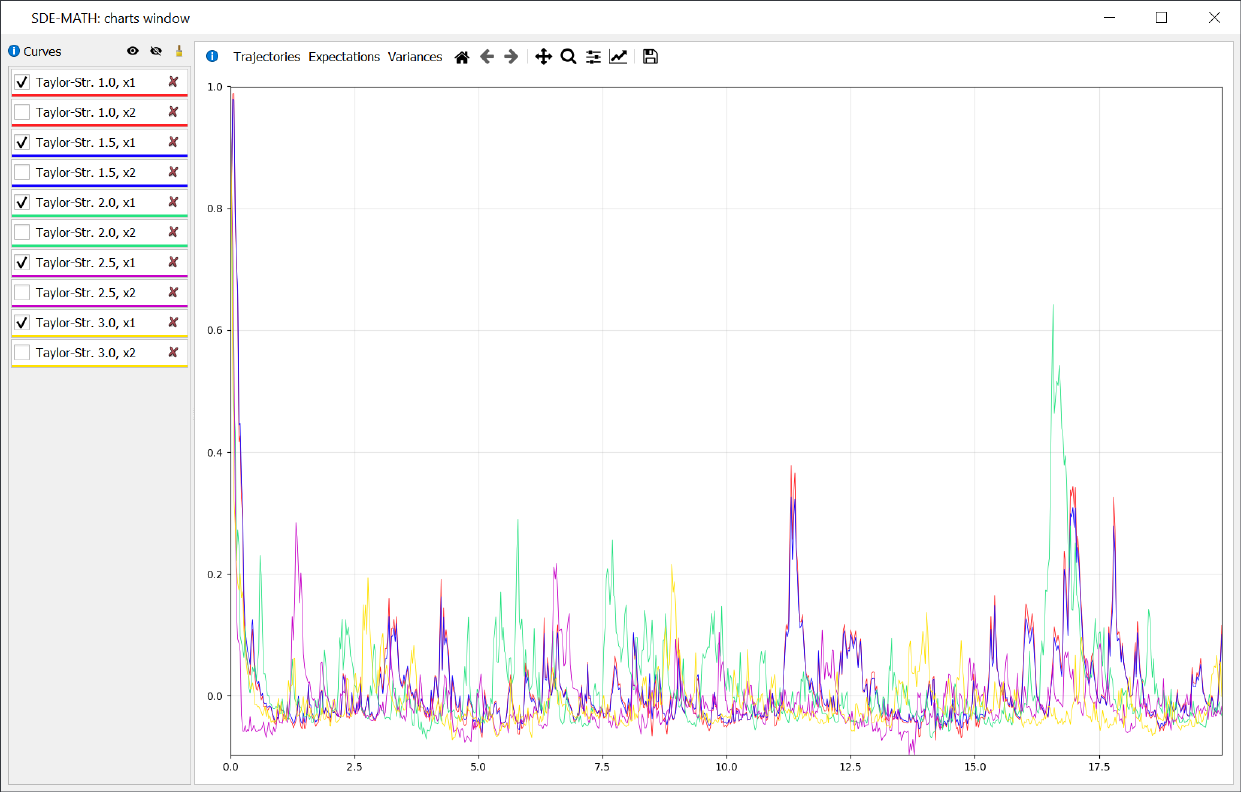}
    \caption{Strong Taylor--Stratonovich schemes of orders 1.0, 1.5, 2.0, 2.5, and 3.0 (${\bf x}_t^{(1)}$ component, $C = 14000,$ $dt = 0.025$)\label{fig:straton_small_3p0_6}}
\end{figure}

\begin{figure}[H]
    \centering
    \includegraphics[width=.9\textwidth]{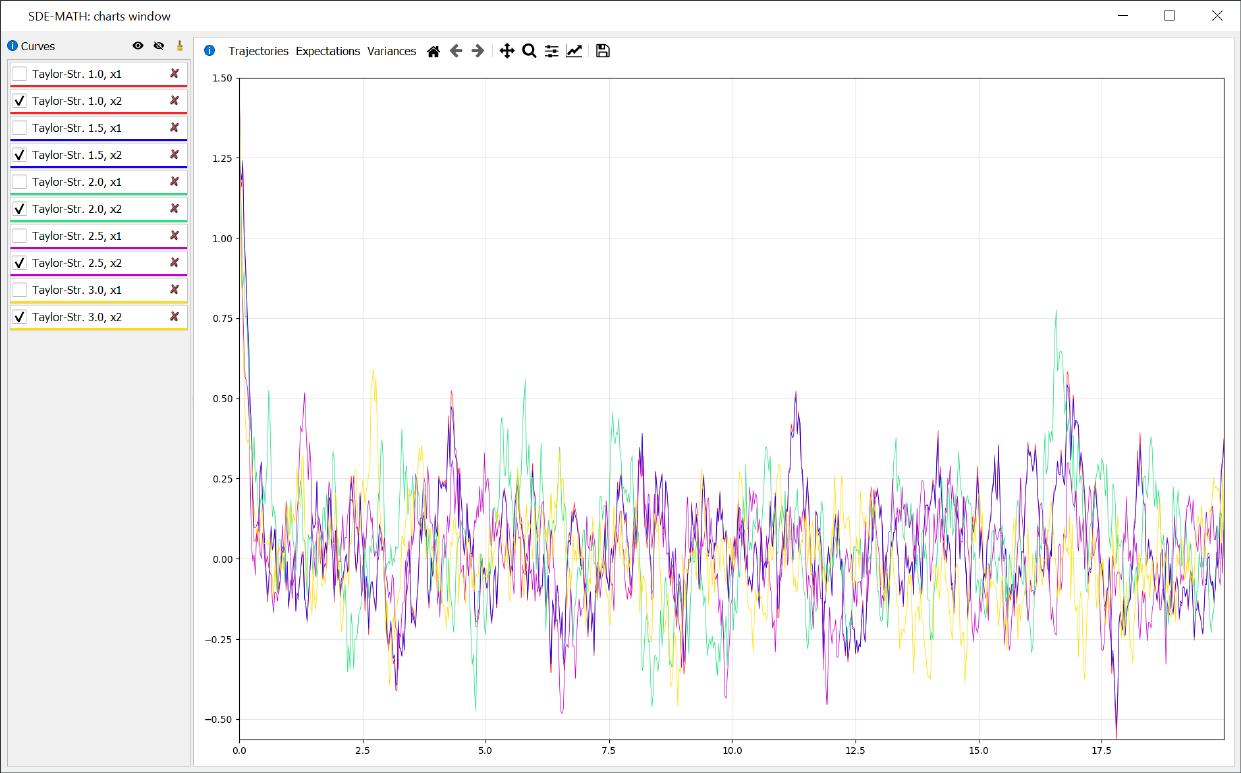}
    \caption{Strong Taylor--Stratonovich schemes of orders 1.0, 1.5, 2.0, 2.5, and 3.0 (${\bf x}_t^{(2)}$ component, $C = 14000,$ $dt = 0.025$)\label{fig:straton_small_3p0_7}}
\end{figure}

\begin{figure}[H]
    \vspace{10mm}
    \centering
    \includegraphics[width=.9\textwidth]{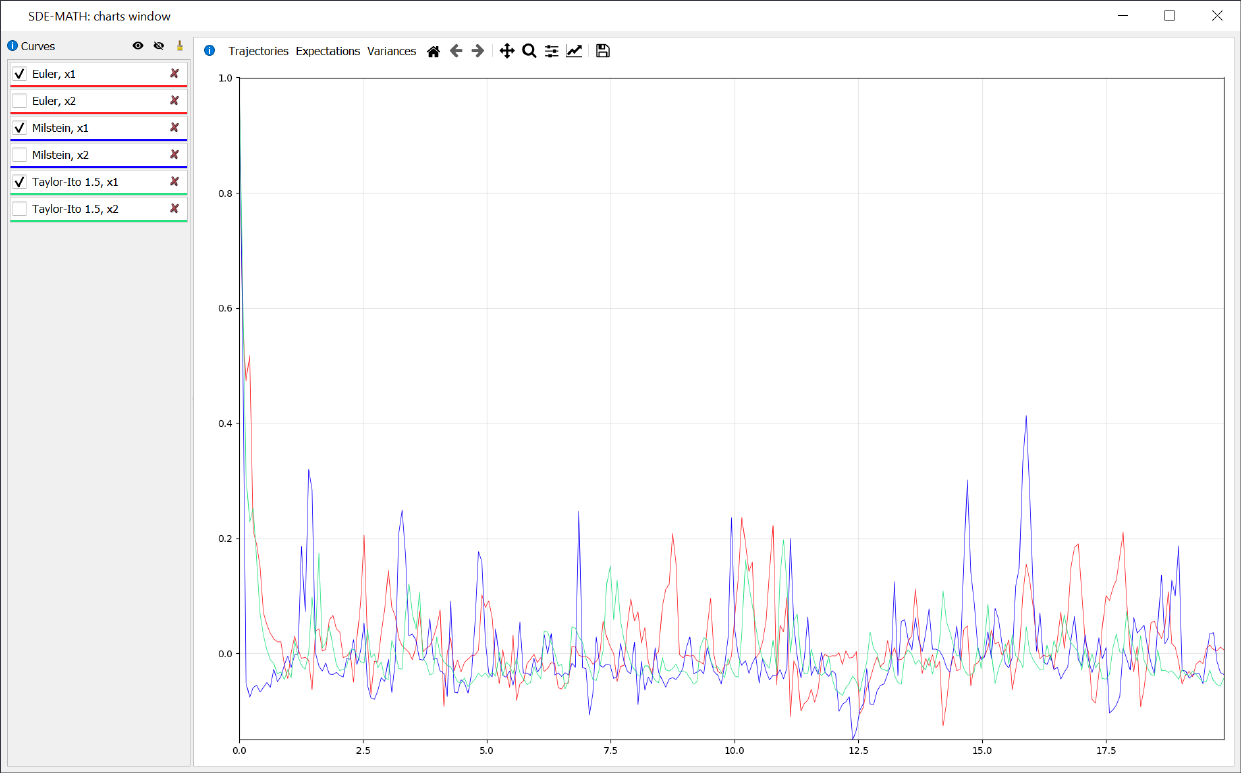}
    \caption{Strong Taylor--It\^o schemes of orders 0.5, 1.0, and 1.5 (${\bf x}_t^{(1)}$ component, $C = 0.1,$ $dt = 0.07$)\label{fig:ito_1p5_big_4}}
\end{figure}

\begin{figure}[H]
    \vspace{7mm}
    \centering

    \hspace*{\fill}
    \begin{subfigure}[b]{.45\textwidth}
        \includegraphics[width=\textwidth]{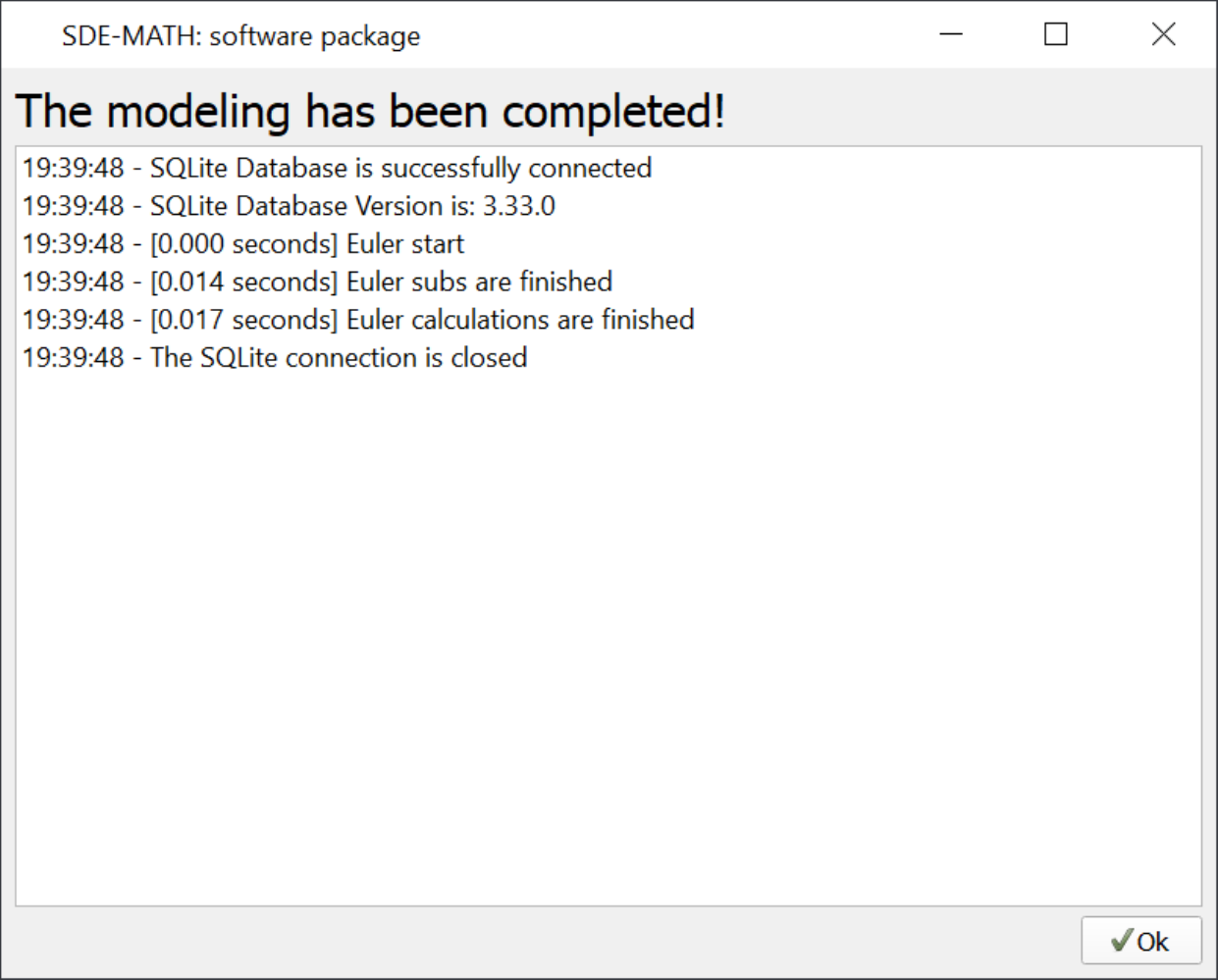}
        \caption*{Euler scheme ($dt = 0.07$)\label{fig:ito_1p5_big_1}}
    \end{subfigure}
    \hfill
    \begin{subfigure}[b]{.45\textwidth}
        \includegraphics[width=\textwidth]{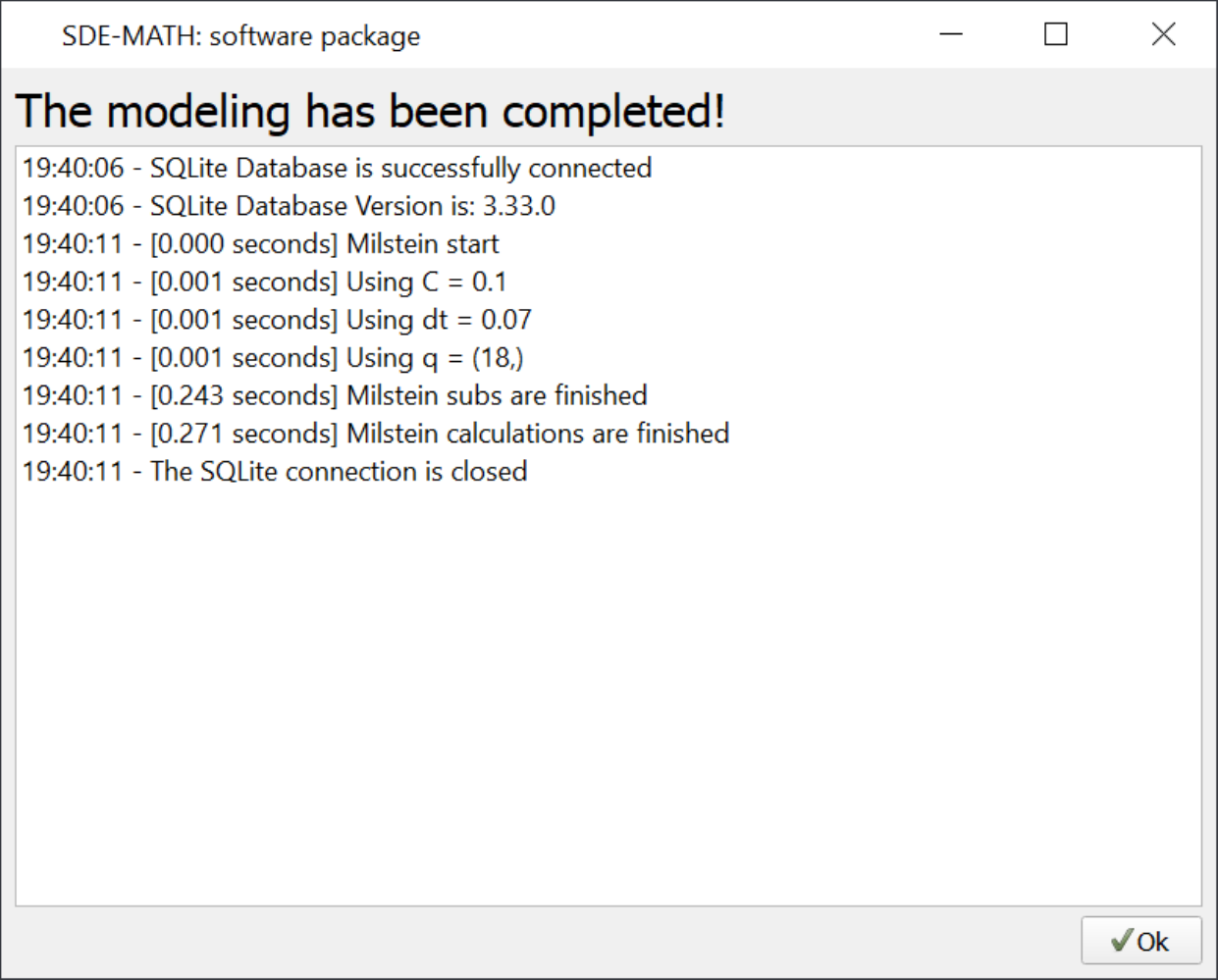}
        \caption*{Milstein scheme ($C = 0.1,$ $dt = 0.07$)\label{fig:ito_1p5_big_2}}
    \end{subfigure}
    \hspace*{\fill}

    \caption{Modeling logs\label{fig:ito_1p5_big_logs1}}

\end{figure}

\begin{figure}[H]
    \vspace{13mm}
    \centering
    \includegraphics[width=.45\textwidth]{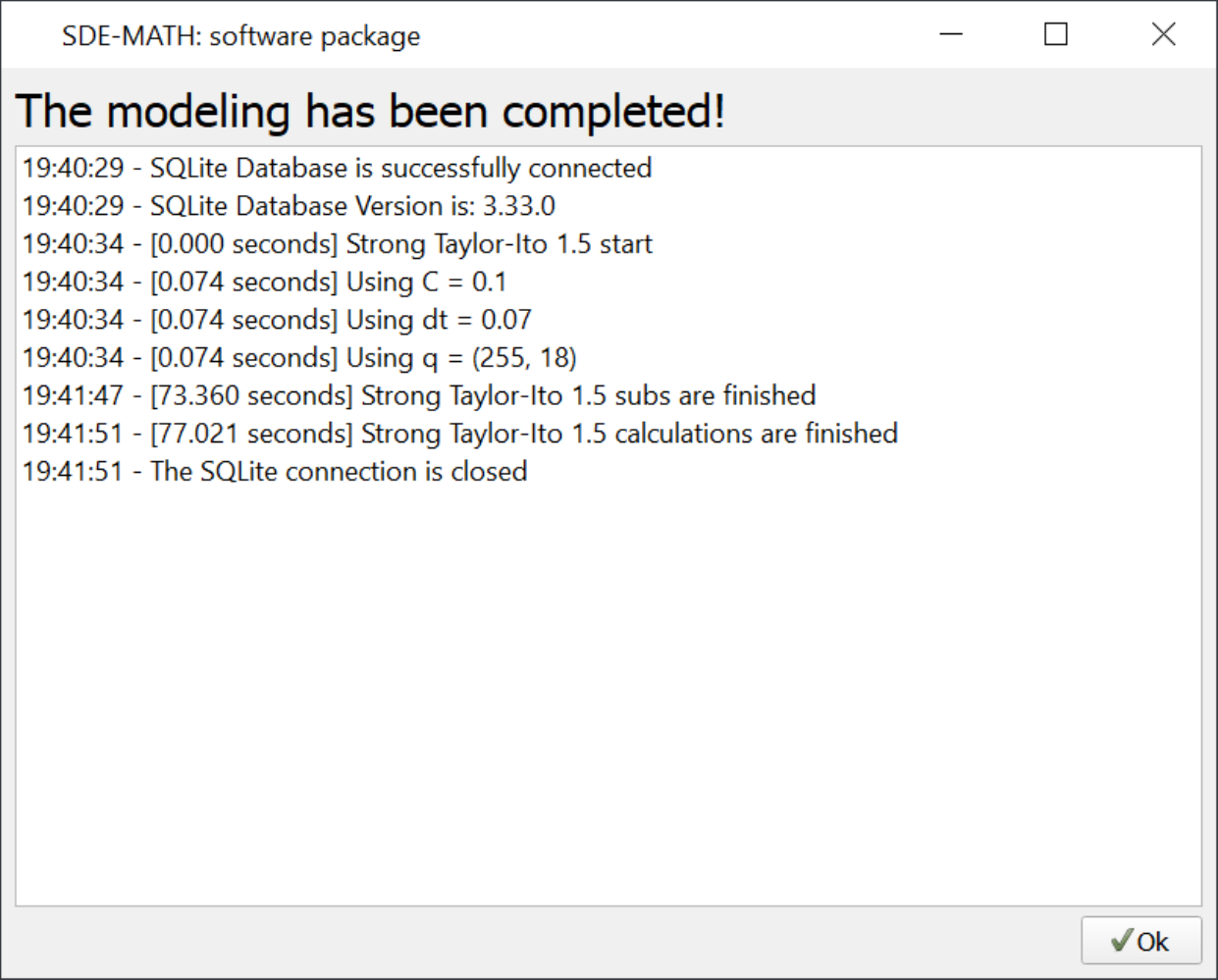}
    \caption{Strong Taylor--It\^o scheme of order 1.5 ($C = 0.1,$ $dt = 0.07$)\label{fig:ito_1p5_big_3}}
\end{figure}

\begin{figure}[H]
    \vspace{10mm}
    \centering
    \includegraphics[width=.9\textwidth]{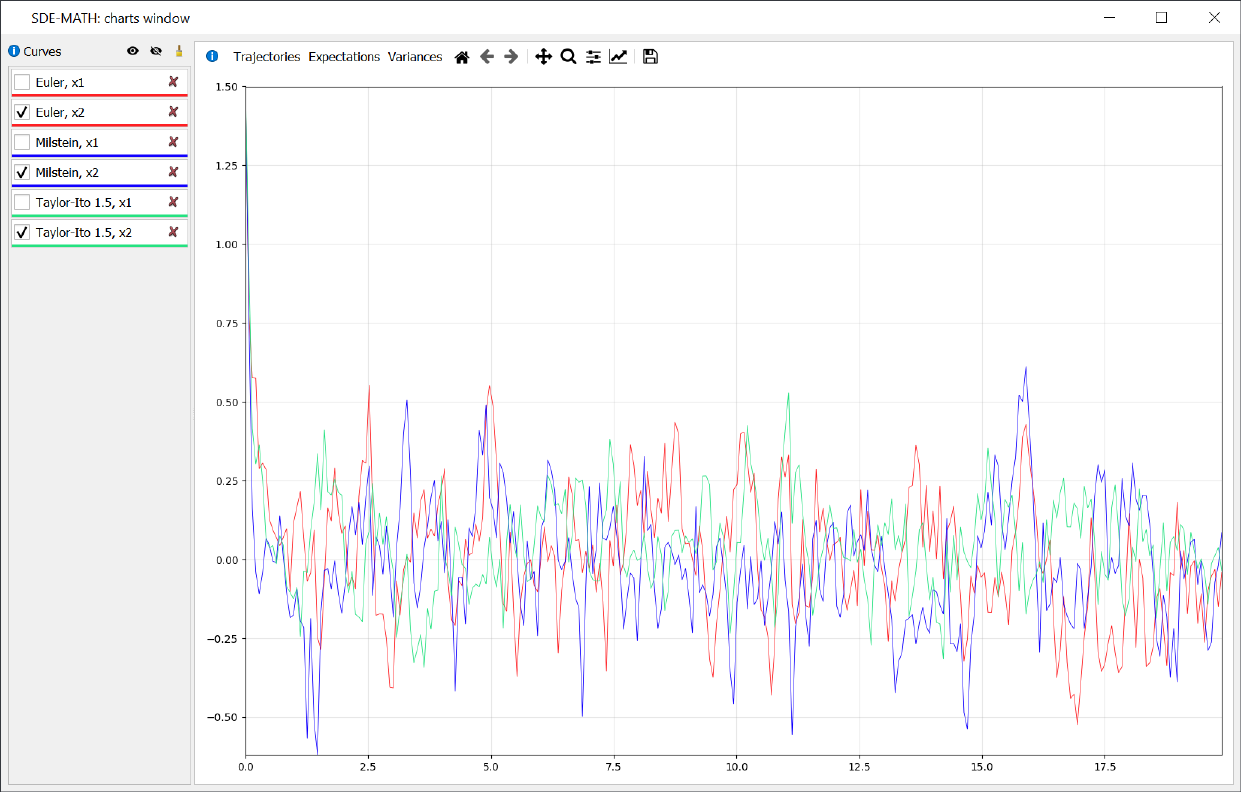}
    \caption{Strong Taylor--It\^o schemes of orders 0.5, 1.0, and 1.5 (${\bf x}_t^{(2)}$ component, $C = 0.1,$ $dt = 0.07$)\label{fig:ito_1p5_big_5}}
\end{figure}

\begin{figure}[H]
    \vspace{10mm}
    \centering
    \includegraphics[width=.9\textwidth]{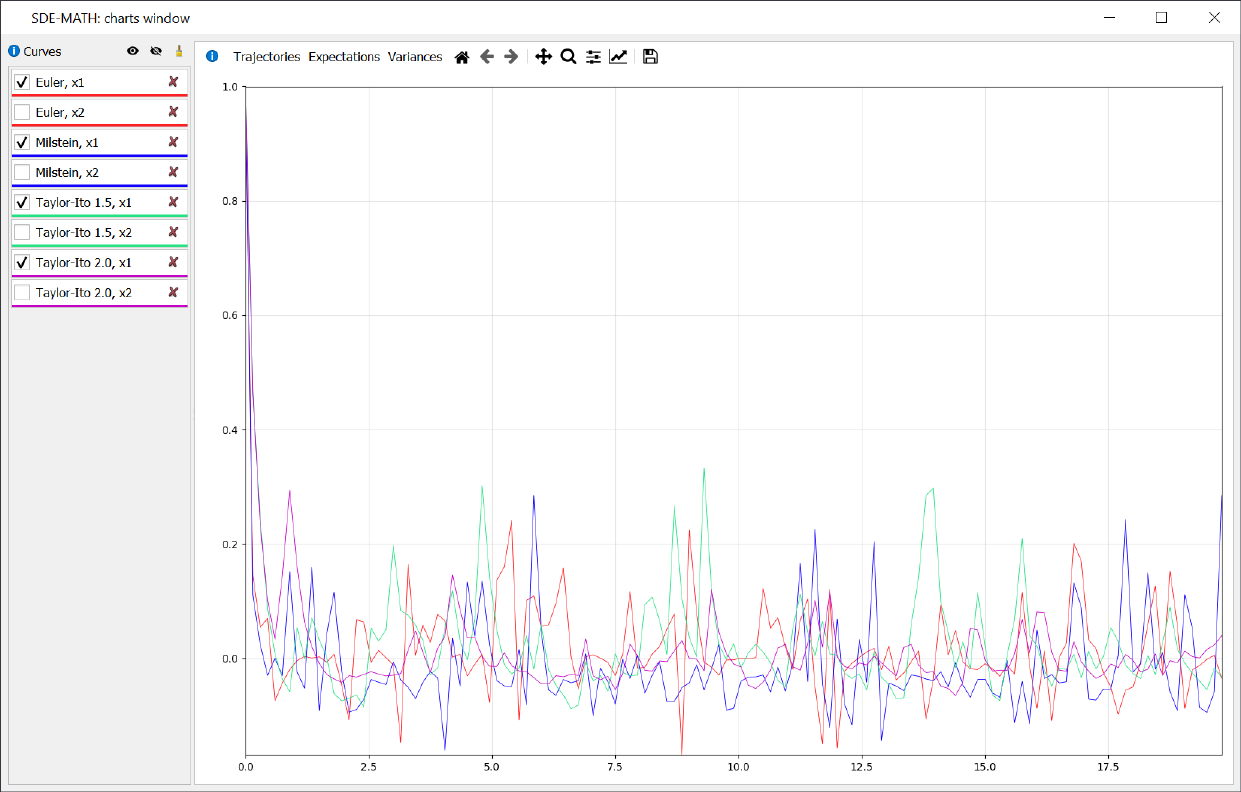}
    \caption{Strong Taylor--It\^o schemes of orders 0.5, 1.0, 1.5, and 2.0 (${\bf x}_t^{(1)}$ component, $C = 0.5,$ $dt = 0.15$)\label{fig:ito_2p0_big_5}}
\end{figure}

\begin{figure}[H]
    \vspace{7mm}
    \centering

    \hspace*{\fill}
    \begin{subfigure}[b]{.45\textwidth}
        \includegraphics[width=\textwidth]{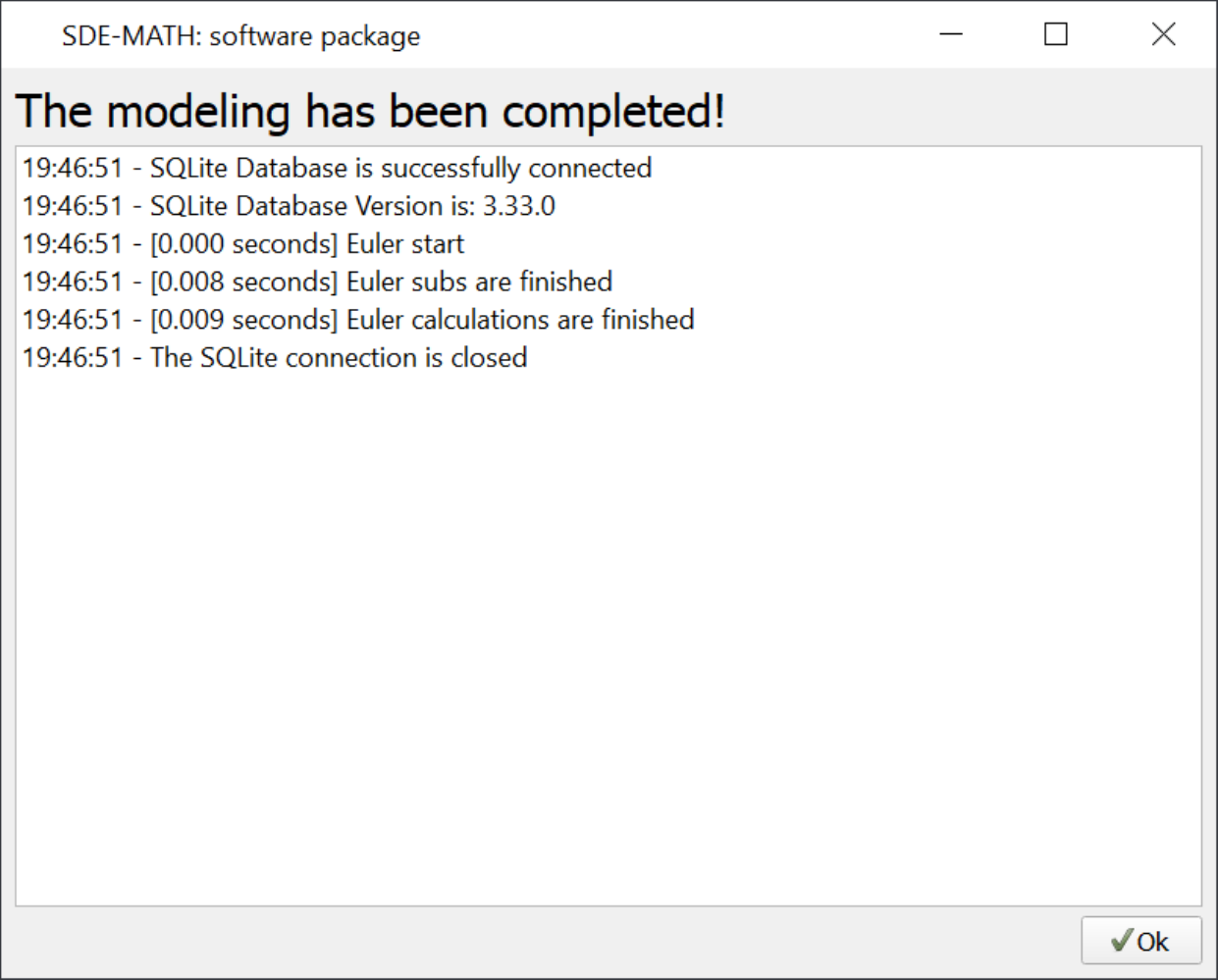}
        \caption*{Euler scheme ($dt = 0.15$)\label{fig:ito_2p0_big_1}}
    \end{subfigure}
    \hfill
    \begin{subfigure}[b]{.45\textwidth}
        \includegraphics[width=\textwidth]{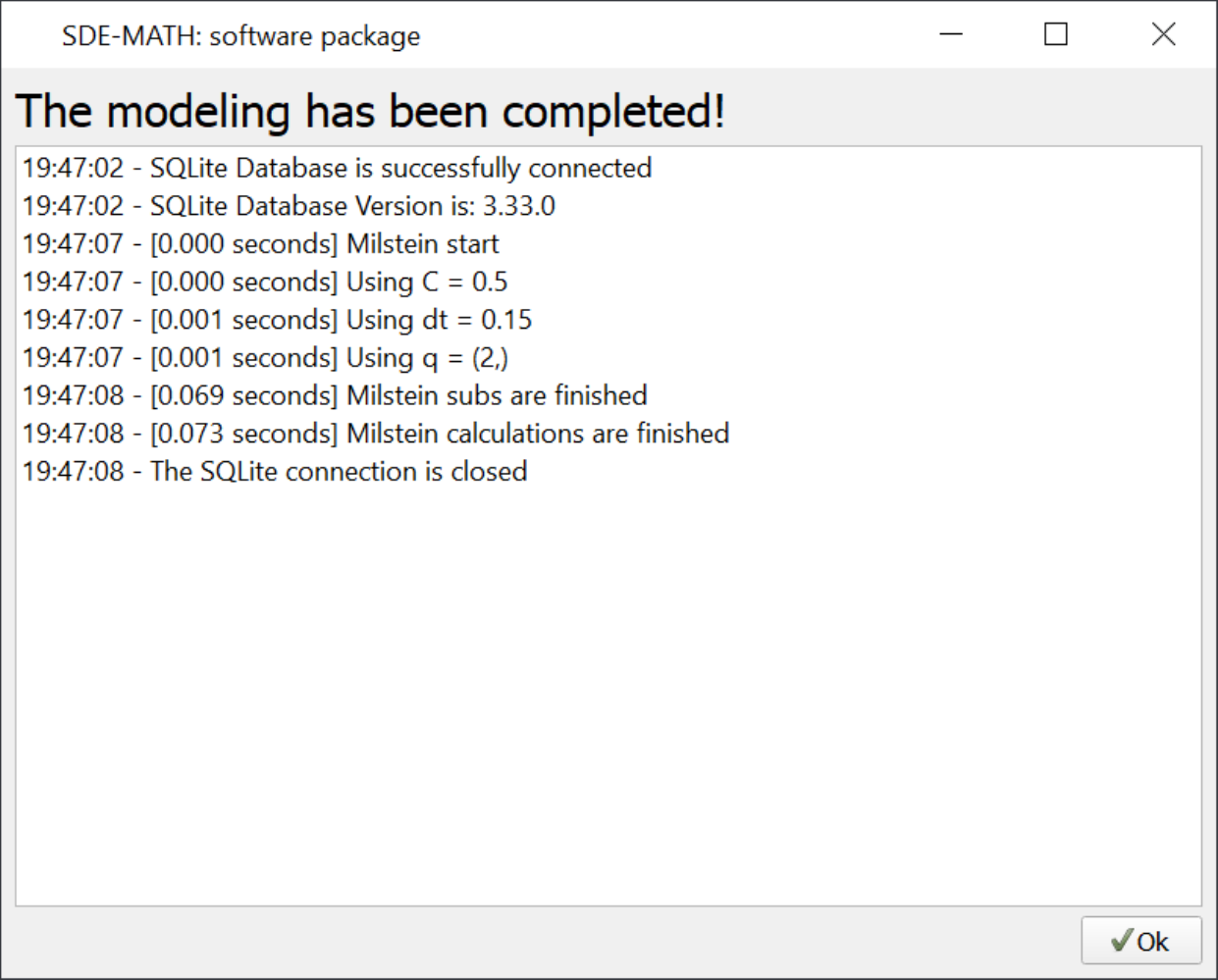}
        \caption*{Milstein scheme ($C = 0.5,$ $dt = 0.15$)\label{fig:ito_2p0_big_2}}
    \end{subfigure}
    \hspace*{\fill}

    \caption{Modeling logs\label{fig:ito_2p0_big_logs1}}

\end{figure}

\begin{figure}[H]
    \vspace{10mm}
    \centering
    \hspace*{\fill}
    \begin{subfigure}[b]{.45\textwidth}
        \includegraphics[width=\textwidth]{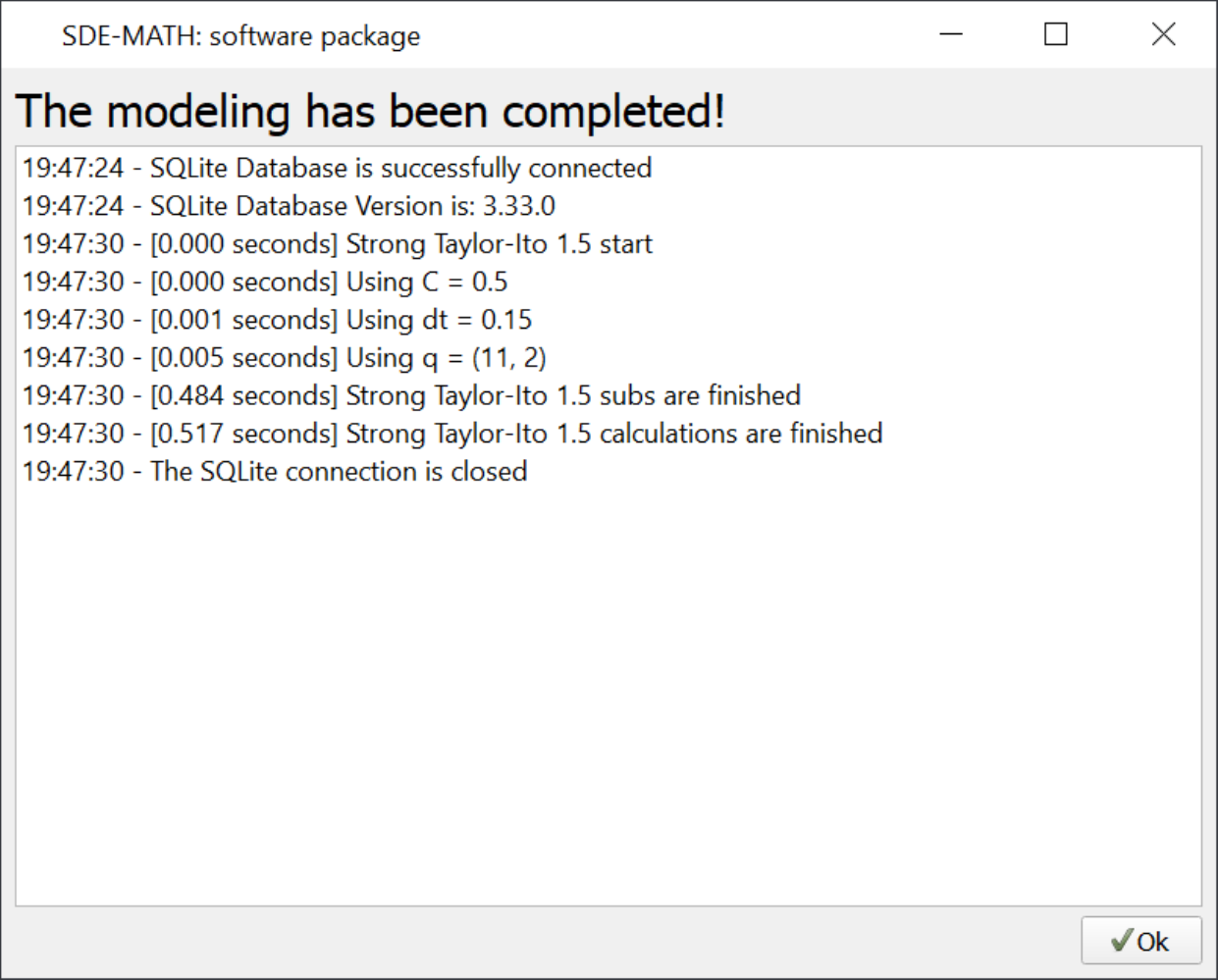}
        \caption*{Strong Taylor--It\^o scheme of order 1.5 ($C = 0.5,$ $dt = 0.15$)\label{fig:ito_2p0_big_3}}
    \end{subfigure}
    \hfill
    \begin{subfigure}[b]{.45\textwidth}
        \includegraphics[width=\textwidth]{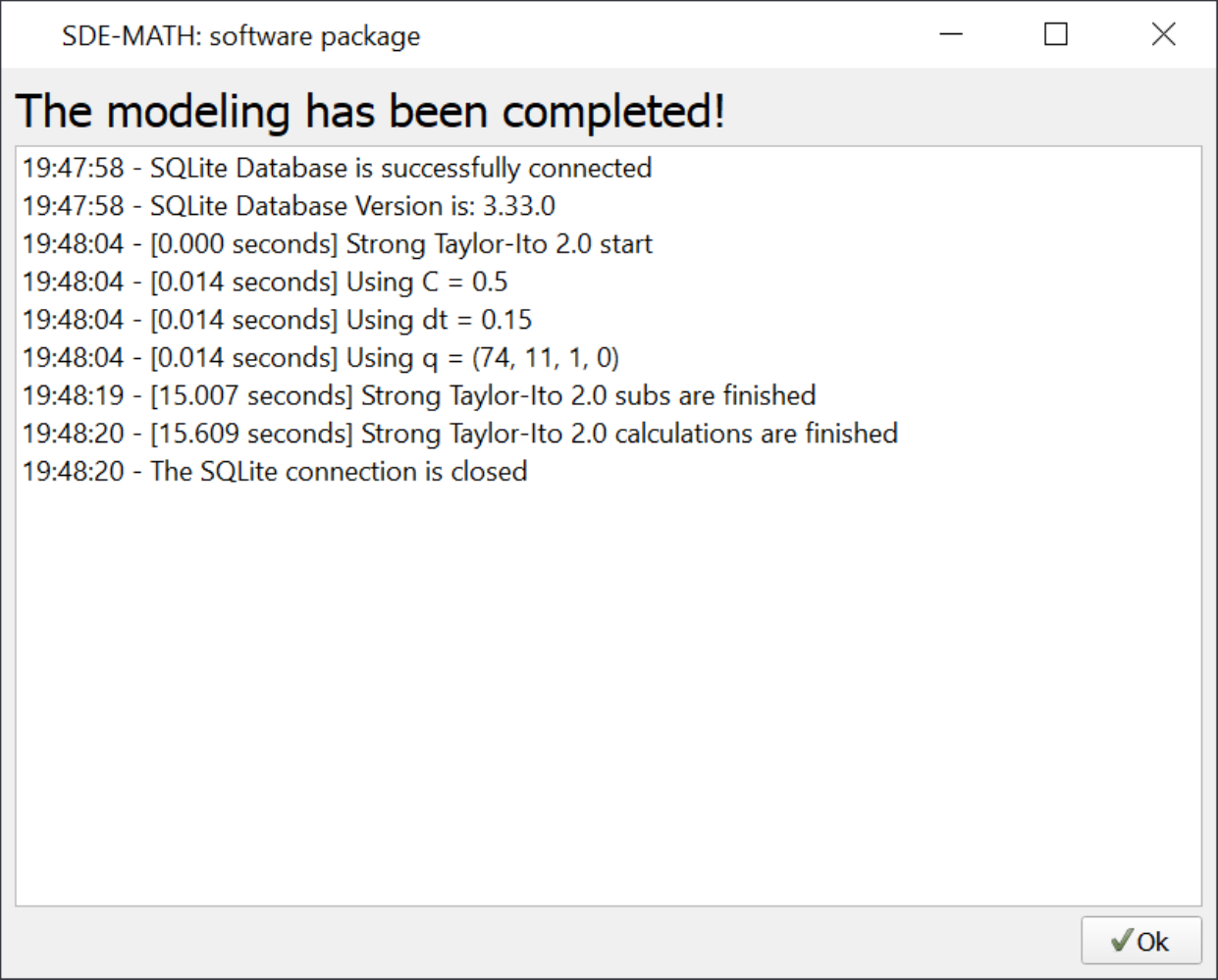}
        \caption*{Strong Taylor--It\^o scheme of order 2.0 ($C = 0.5,$ $dt = 0.15$)\label{fig:ito_2p0_big_4}}
    \end{subfigure}
    \hspace*{\fill}

    \caption{Modeling logs\label{fig:ito_2p0_big_logs2}}

\end{figure}

\begin{figure}[H]
    \vspace{7mm}
    \centering
    \includegraphics[width=.9\textwidth]{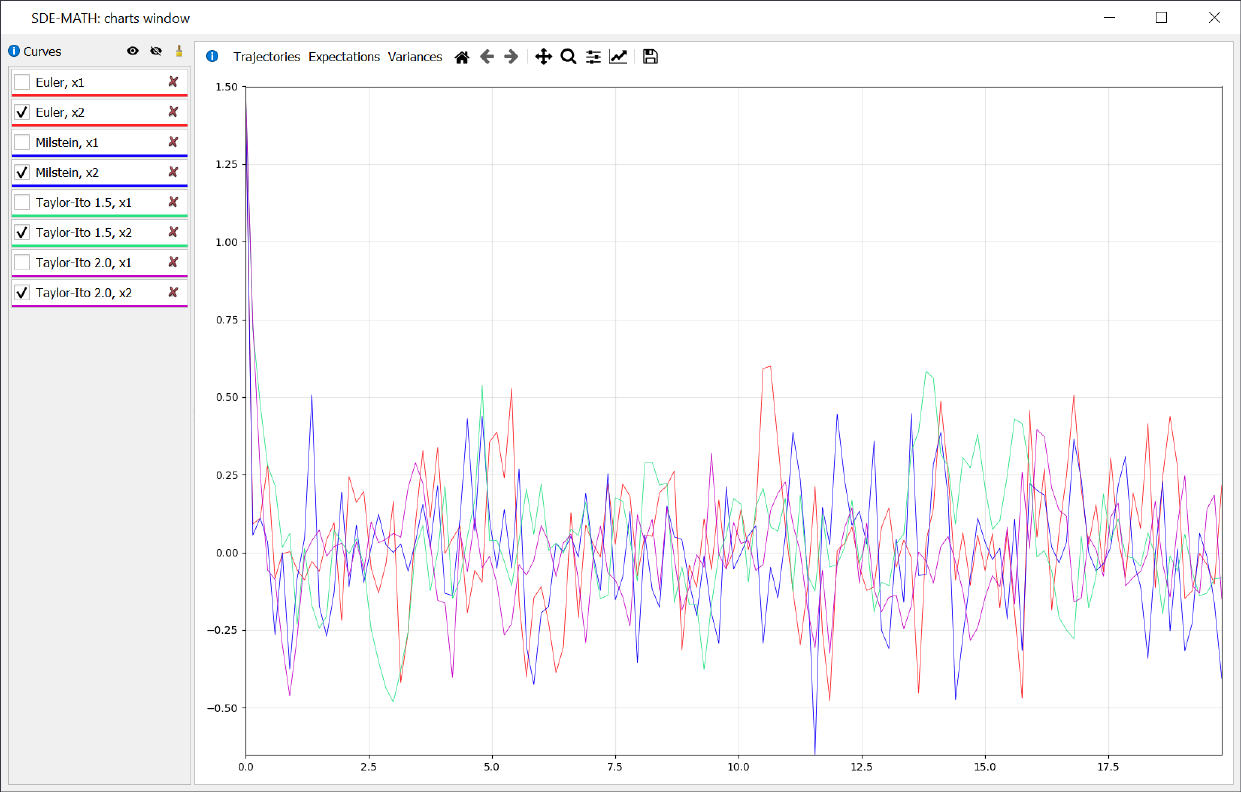}
    \caption{Strong Taylor--It\^o schemes of orders 0.5, 1.0, 1.5, and 2.0 (${\bf x}_t^{(2)}$ component, $C = 0.5,$ $dt = 0.15$)\label{fig:ito_2p0_big_6}}
\end{figure}

\begin{figure}[H]
    \centering

    \hspace*{\fill}
    \begin{subfigure}[b]{.45\textwidth}
        \includegraphics[width=\textwidth]{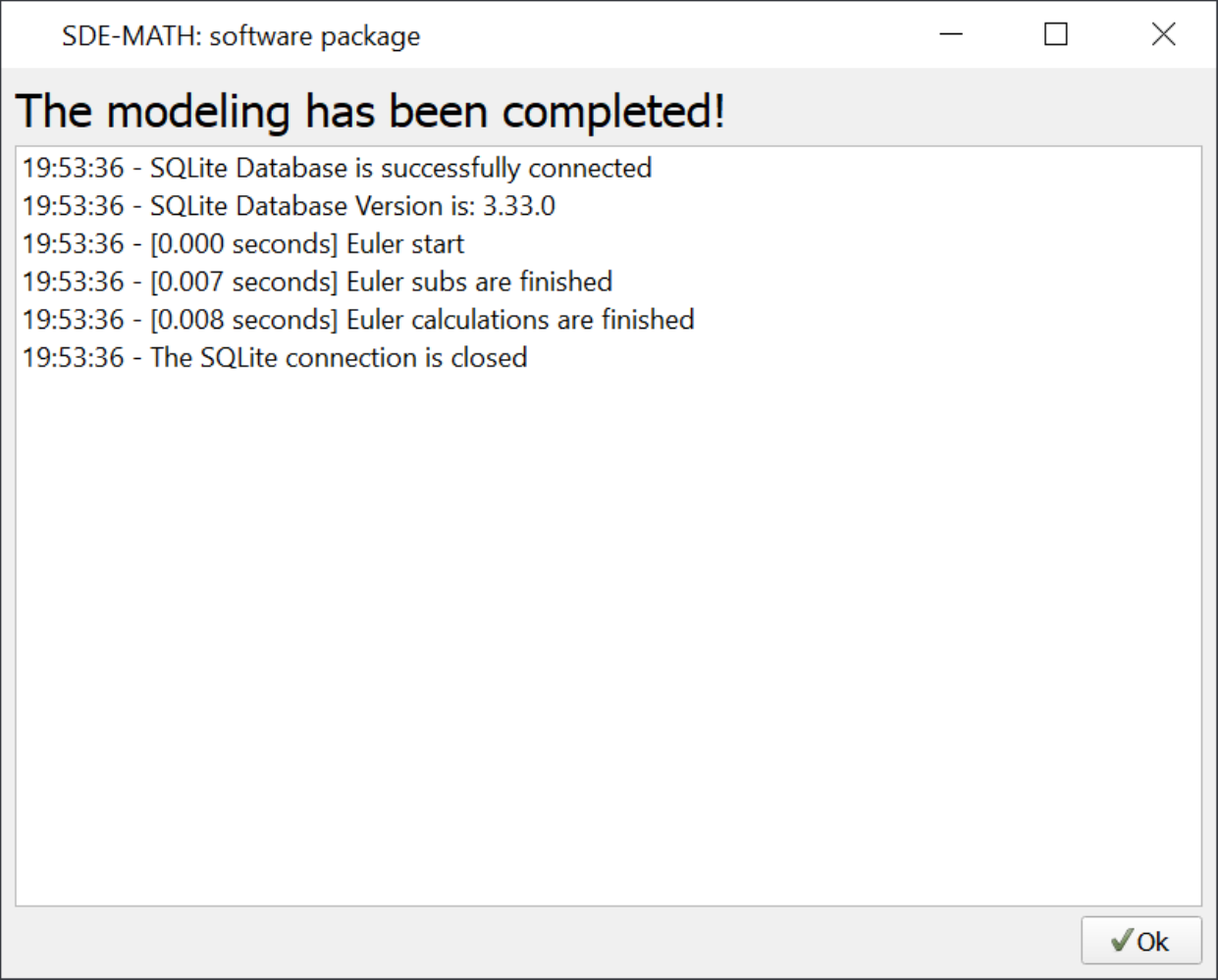}
        \caption*{Euler scheme ($dt = 0.2$)\label{fig:ito_2p5_big_1}}
    \end{subfigure}
    \hfill
    \begin{subfigure}[b]{.45\textwidth}
        \includegraphics[width=\textwidth]{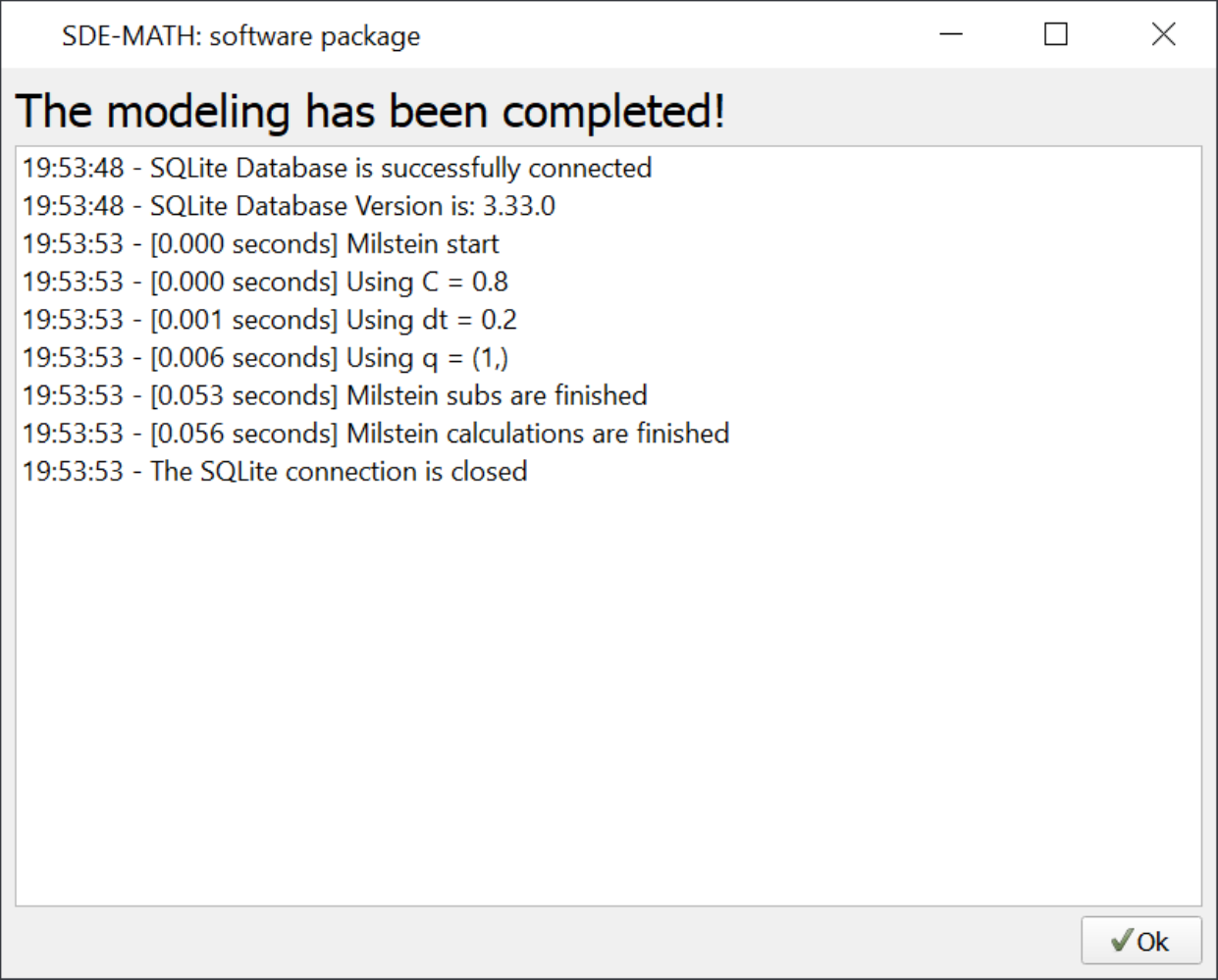}
        \caption*{Milstein scheme ($C = 0.8,$ $dt = 0.2$)\label{fig:ito_2p5_big_2}}
    \end{subfigure}
    \hspace*{\fill}

    \vspace{2mm}
    \hspace*{\fill}
    \begin{subfigure}[b]{.45\textwidth}
        \includegraphics[width=\textwidth]{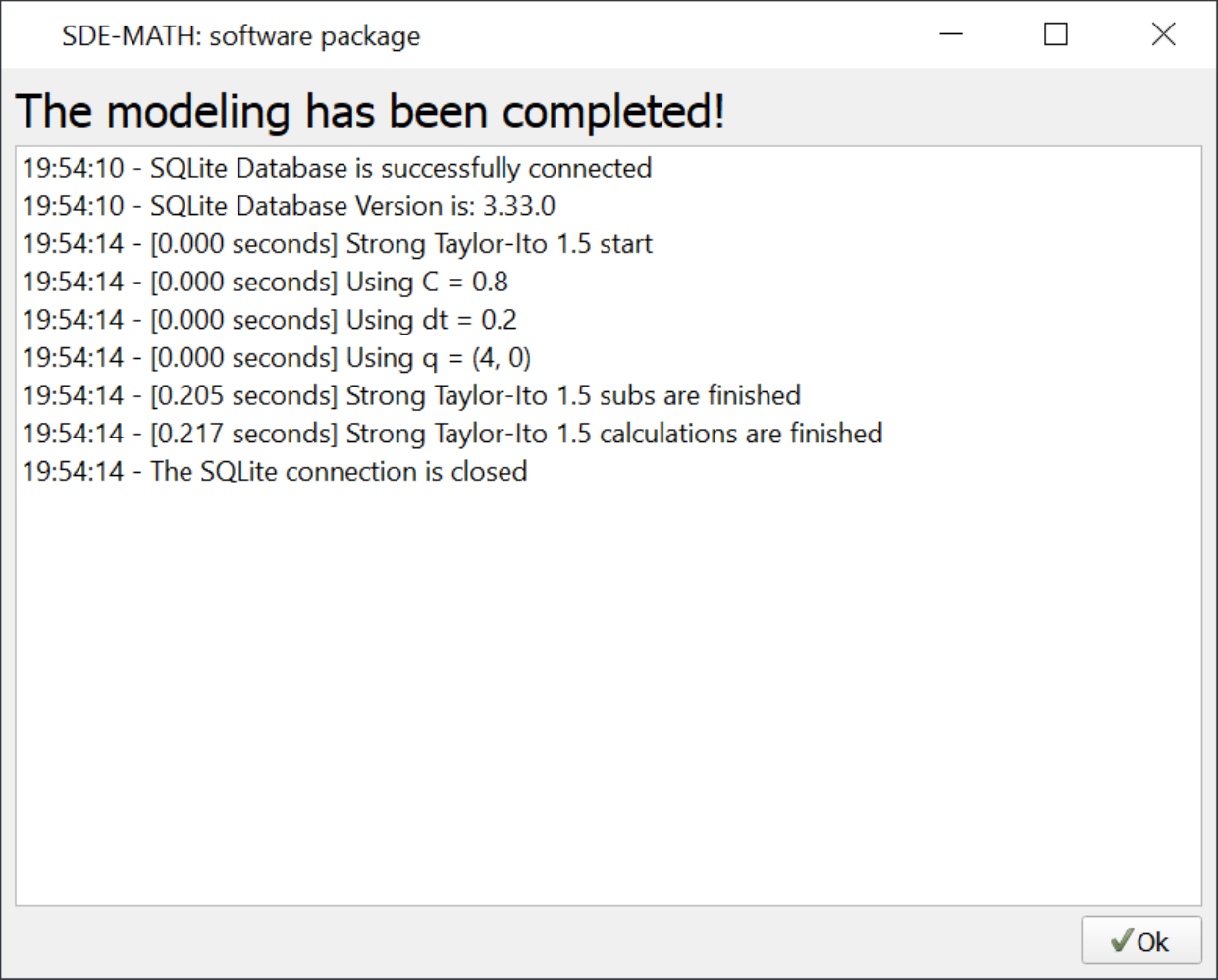}
        \caption*{Strong Taylor--It\^o scheme of order 1.5 ($C = 0.8,$ $dt = 0.2$)\label{fig:ito_2p5_big_3}}
    \end{subfigure}
    \hfill
    \begin{subfigure}[b]{.45\textwidth}
        \includegraphics[width=\textwidth]{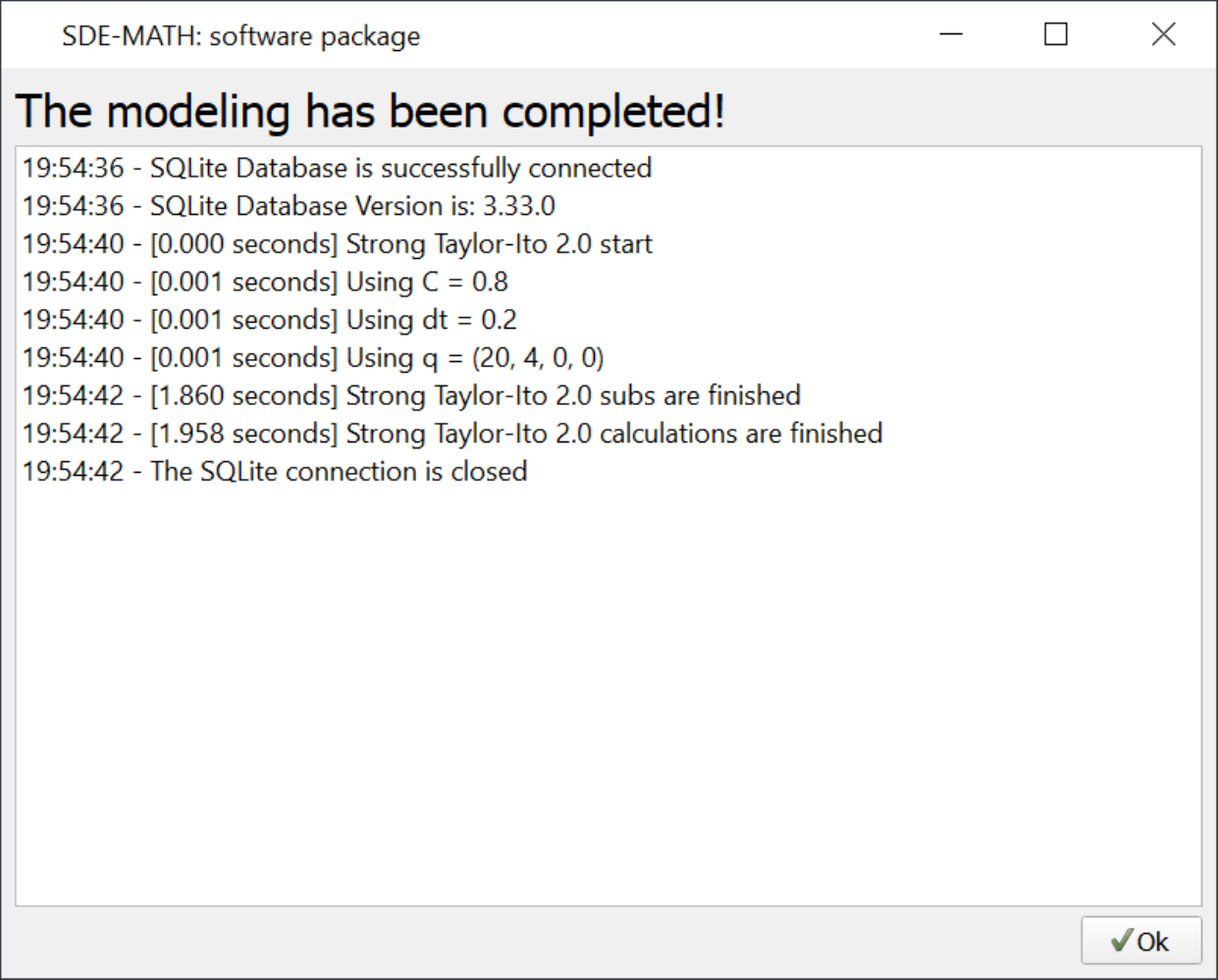}
        \caption*{Strong Taylor--It\^o scheme of order 2.0 ($C = 0.8,$ $dt = 0.2$)\label{fig:ito_2p5_big_4}}
    \end{subfigure}
    \hspace*{\fill}

    \vspace{2mm}
    \hspace*{\fill}
    \begin{subfigure}[b]{.45\textwidth}
        \includegraphics[width=\textwidth]{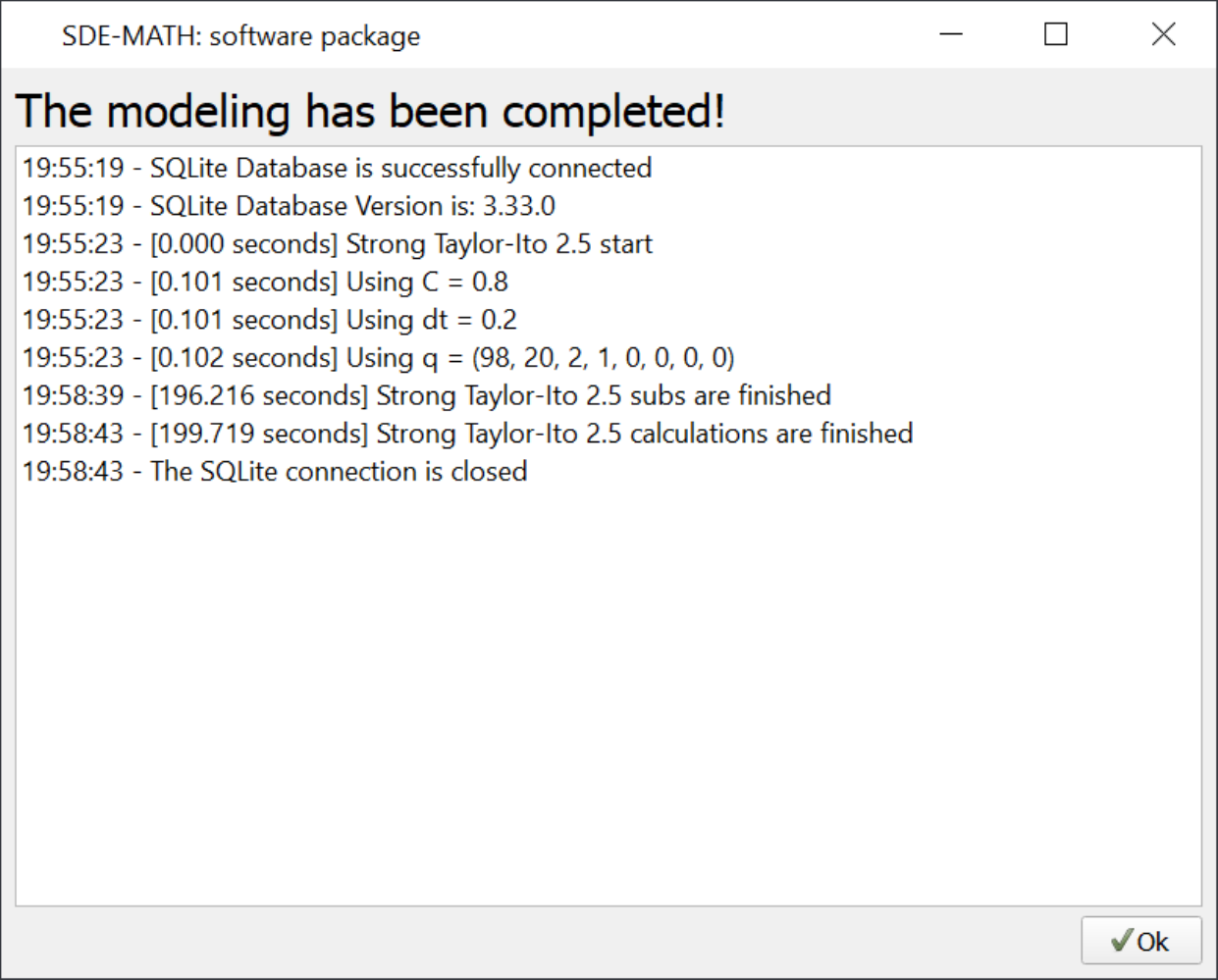}
        \caption*{Strong Taylor--It\^o scheme of order 2.5 ($C = 0.8,$ $dt = 0.2$)\label{fig:ito_2p5_big_5}}
    \end{subfigure}
    \hspace*{\fill}

    \caption{Modeling logs\label{fig:ito_2p5_big_logs}}

\end{figure}

\begin{figure}[H]
    \centering
    \includegraphics[width=.9\textwidth]{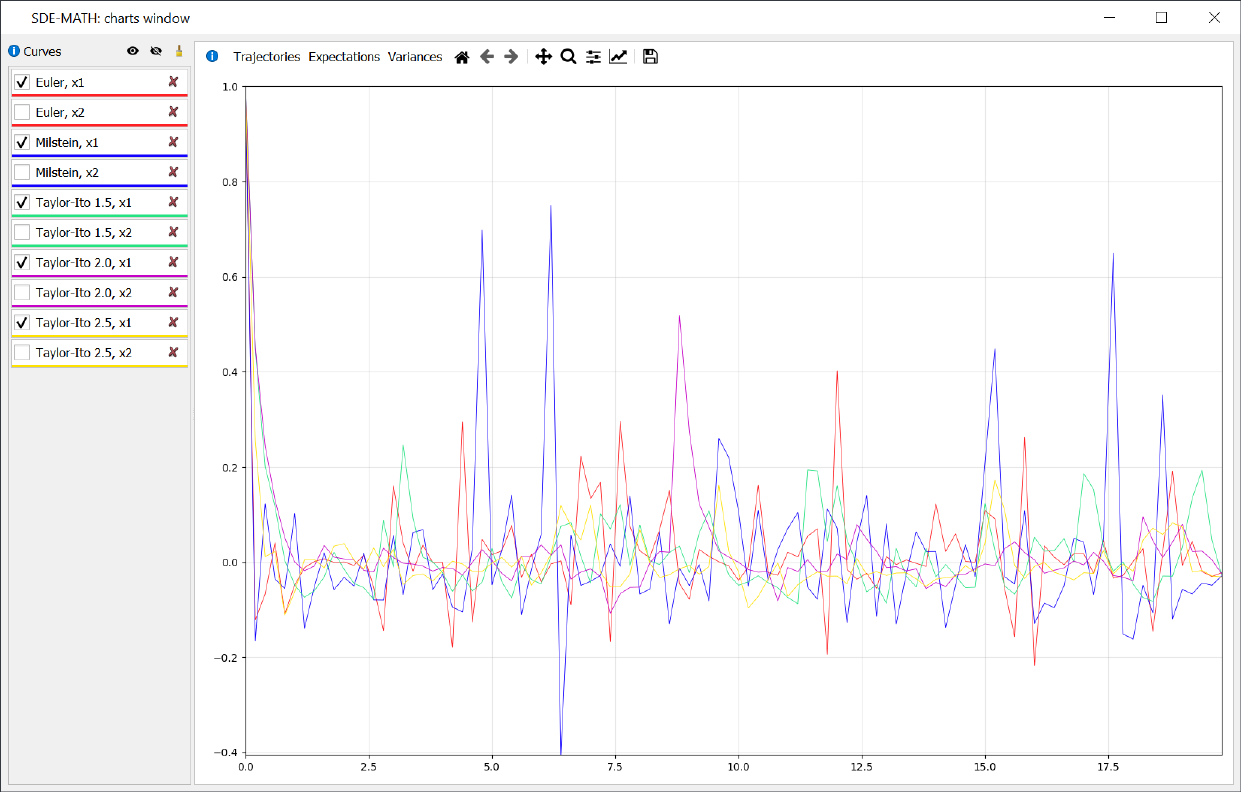}
    \caption{Strong Taylor--It\^o schemes of orders 0.5, 1.0, 1.5, 2.0, and 2.5 (${\bf x}_t^{(1)}$ component, $C = 0.8,$ $dt = 0.2$)\label{fig:ito_2p5_big_6}}
\end{figure}

\begin{figure}[H]
    \centering
    \includegraphics[width=.9\textwidth]{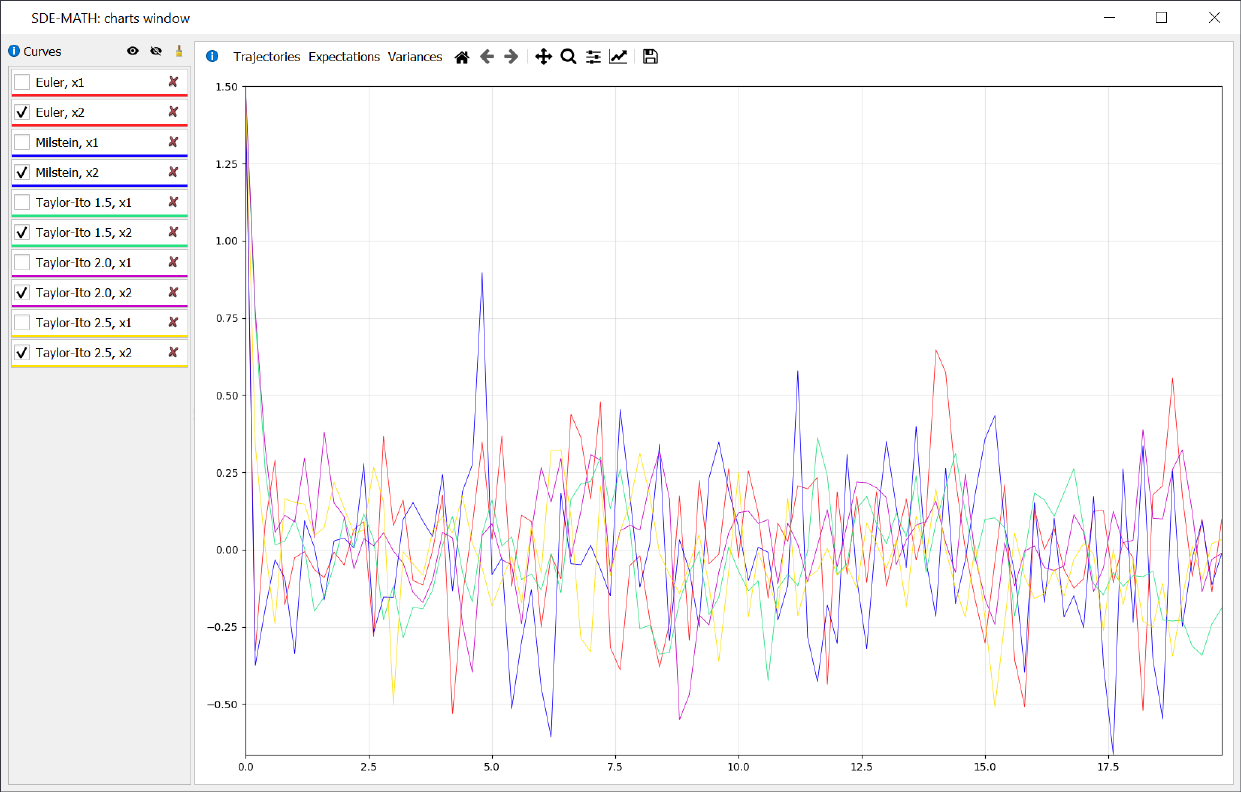}
    \caption{Strong Taylor--It\^o schemes of orders 0.5, 1.0, 1.5, 2.0, and 2.5 (${\bf x}_t^{(2)}$ component, $C = 0.8,$ $dt = 0.2$)\label{fig:ito_2p5_big_7}}
\end{figure}

\begin{figure}[H]
    \centering

    \hspace*{\fill}
    \begin{subfigure}[b]{.45\textwidth}
        \includegraphics[width=\textwidth]{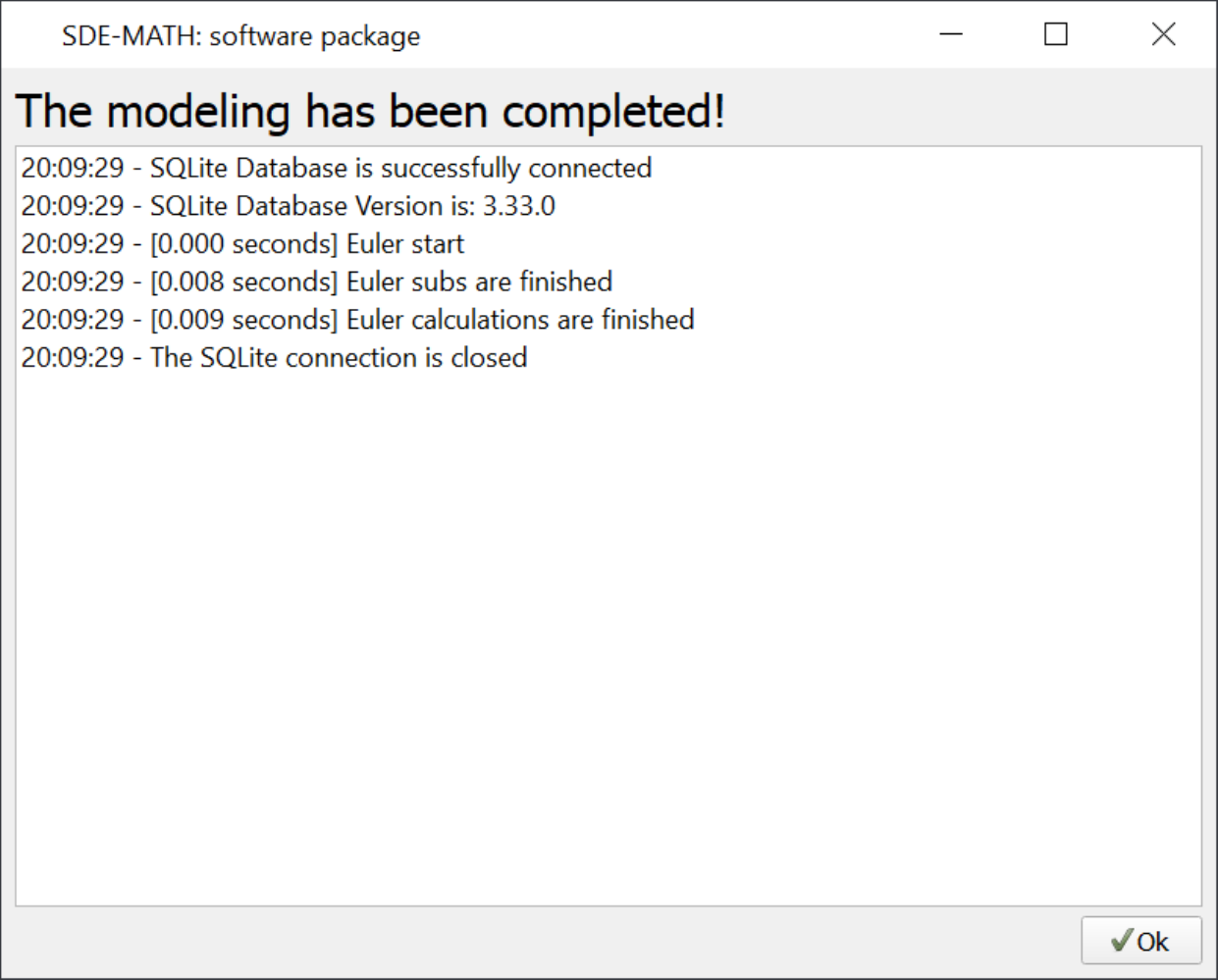}
        \caption*{Euler scheme ($dt = 0.2$)\label{fig:ito_3p0_big_1}}
    \end{subfigure}
    \hfill
    \begin{subfigure}[b]{.45\textwidth}
        \includegraphics[width=\textwidth]{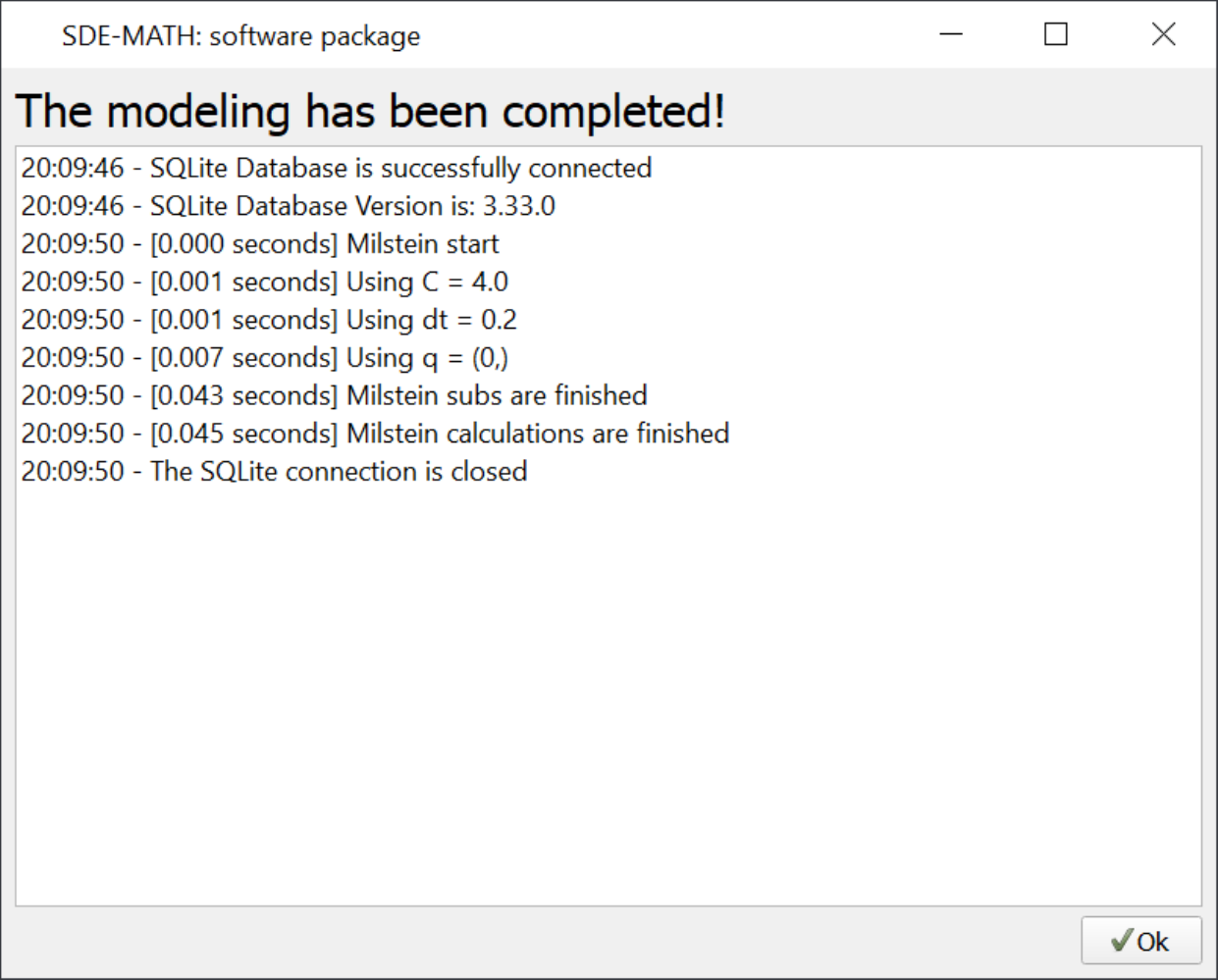}
        \caption*{Milstein scheme ($C = 4,$ $dt = 0.2$)\label{fig:ito_3p0_big_2}}
    \end{subfigure}
    \hspace*{\fill}

    \vspace{2mm}
    \hspace*{\fill}
    \begin{subfigure}[b]{.45\textwidth}
        \includegraphics[width=\textwidth]{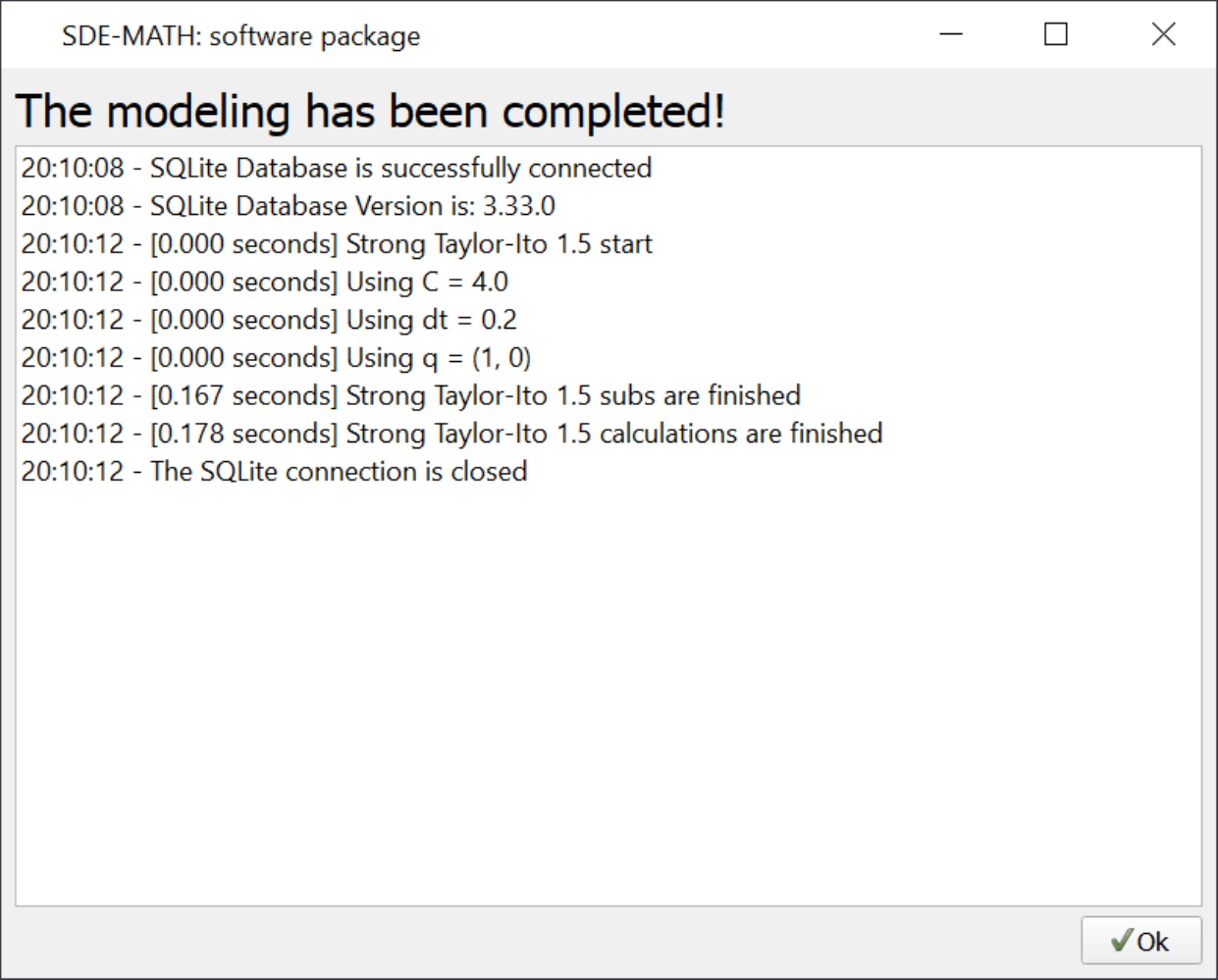}
        \caption*{Strong Taylor--It\^o scheme of order 1.5 ($C = 4,$ $dt = 0.2$)\label{fig:ito_3p0_big_3}}
    \end{subfigure}
    \hfill
    \begin{subfigure}[b]{.45\textwidth}
        \includegraphics[width=\textwidth]{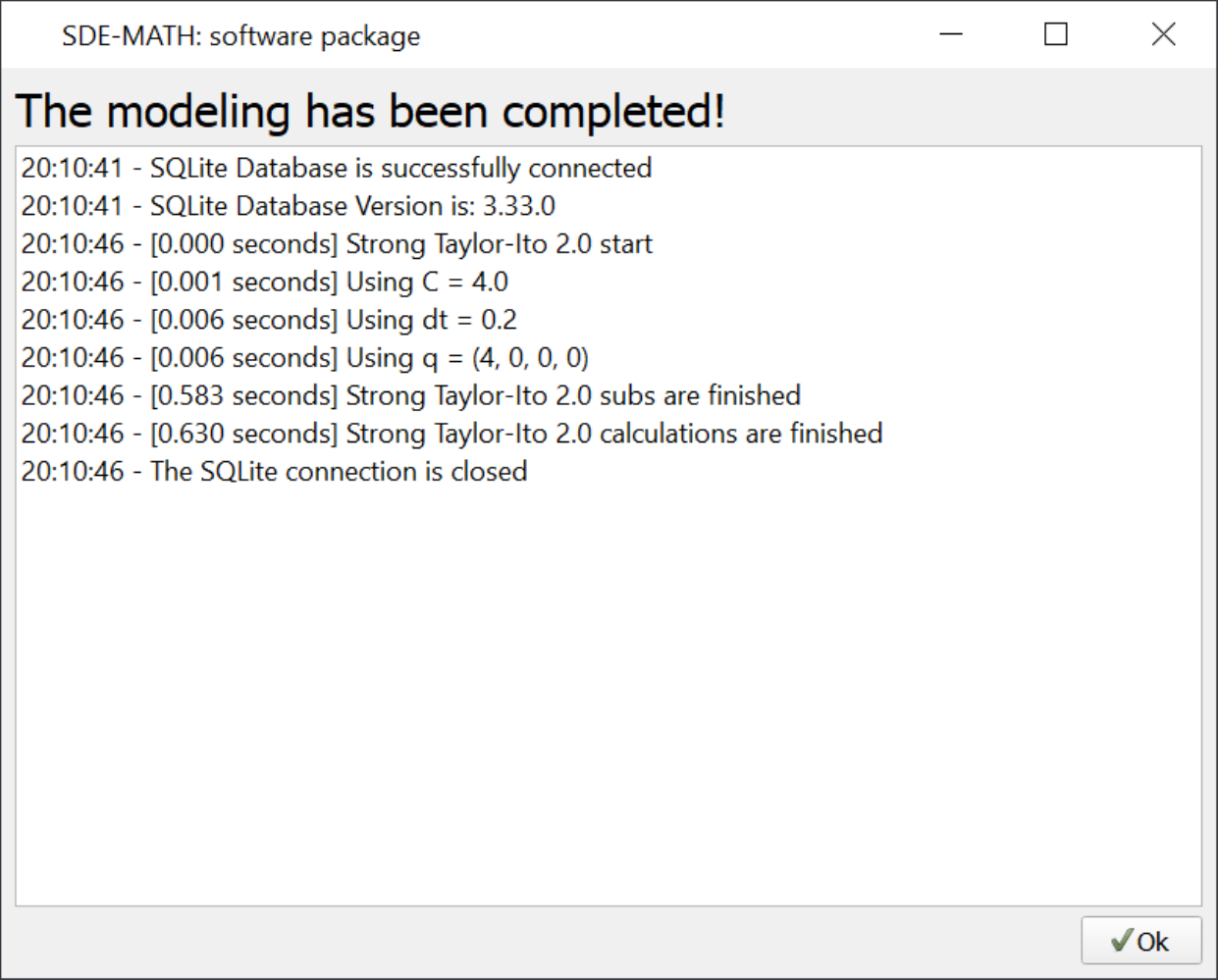}
        \caption*{Strong Taylor--It\^o scheme of order 2.0 ($C = 4,$ $dt = 0.2$)\label{fig:ito_3p0_big_4}}
    \end{subfigure}
    \hspace*{\fill}

    \vspace{2mm}
    \hspace*{\fill}
    \begin{subfigure}[b]{.45\textwidth}
        \includegraphics[width=\textwidth]{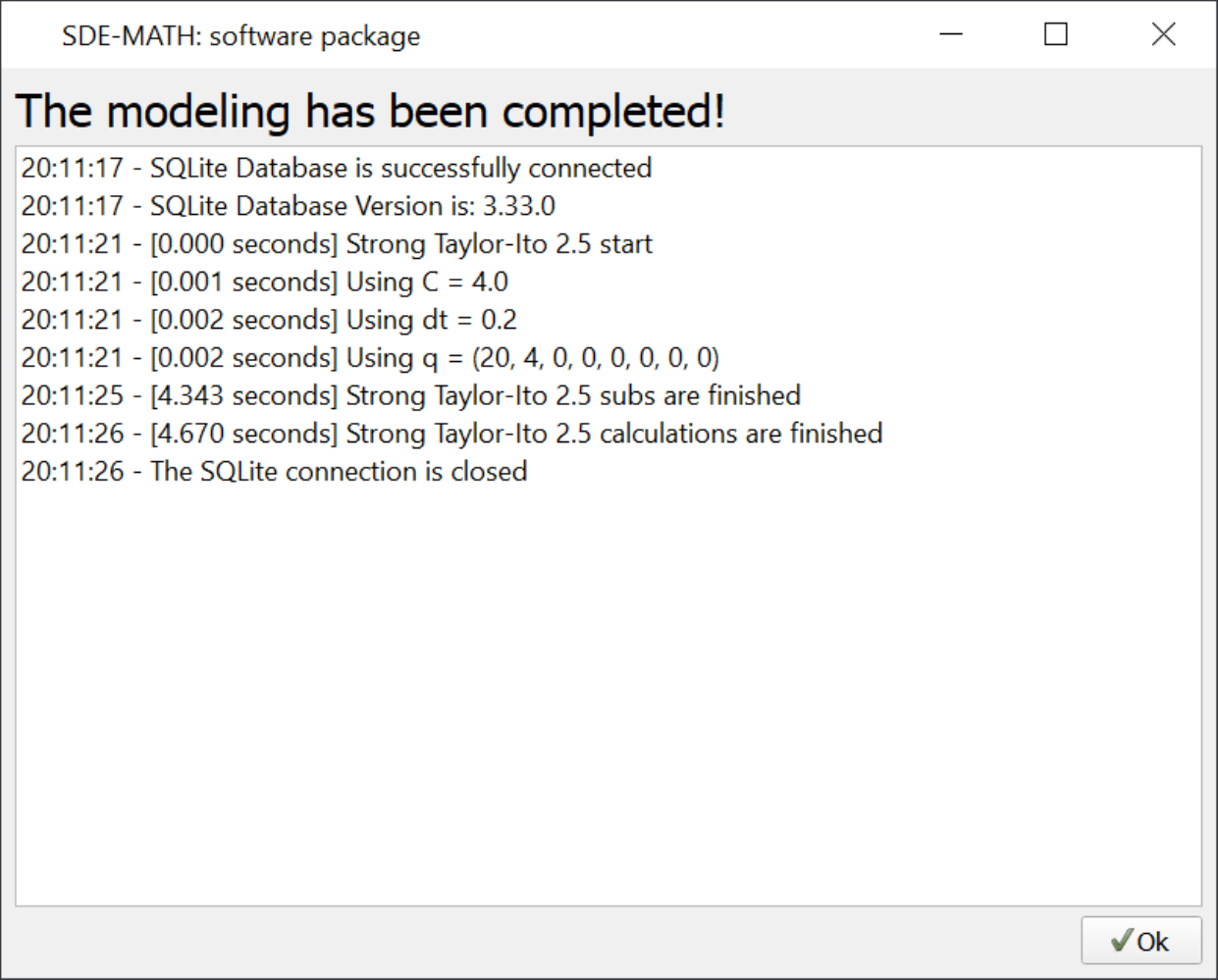}
        \caption*{Strong Taylor--It\^o scheme of order 2.5 ($C = 4,$ $dt = 0.2$)\label{fig:ito_3p0_big_5}}
    \end{subfigure}
    \hfill
    \begin{subfigure}[b]{.45\textwidth}
        \includegraphics[width=\textwidth]{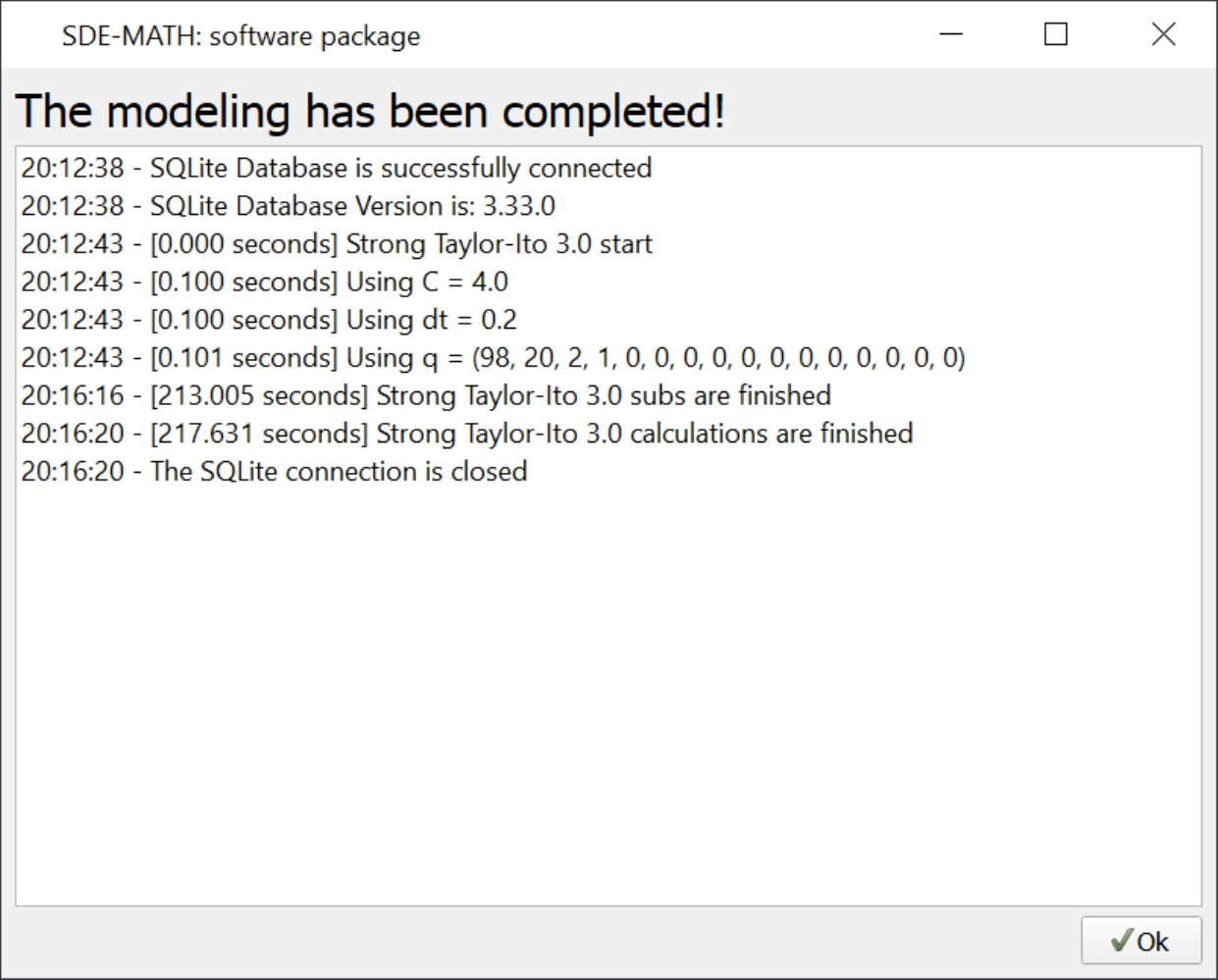}
        \caption*{Strong Taylor--It\^o scheme of order 3.0 ($C = 4,$ $dt = 0.2$)\label{fig:ito_3p0_big_6}}
    \end{subfigure}
    \hspace*{\fill}

    \caption{Modeling logs\label{fig:ito_3p0_big_logs}}

\end{figure}

\begin{figure}[H]
    \centering
    \includegraphics[width=.9\textwidth]{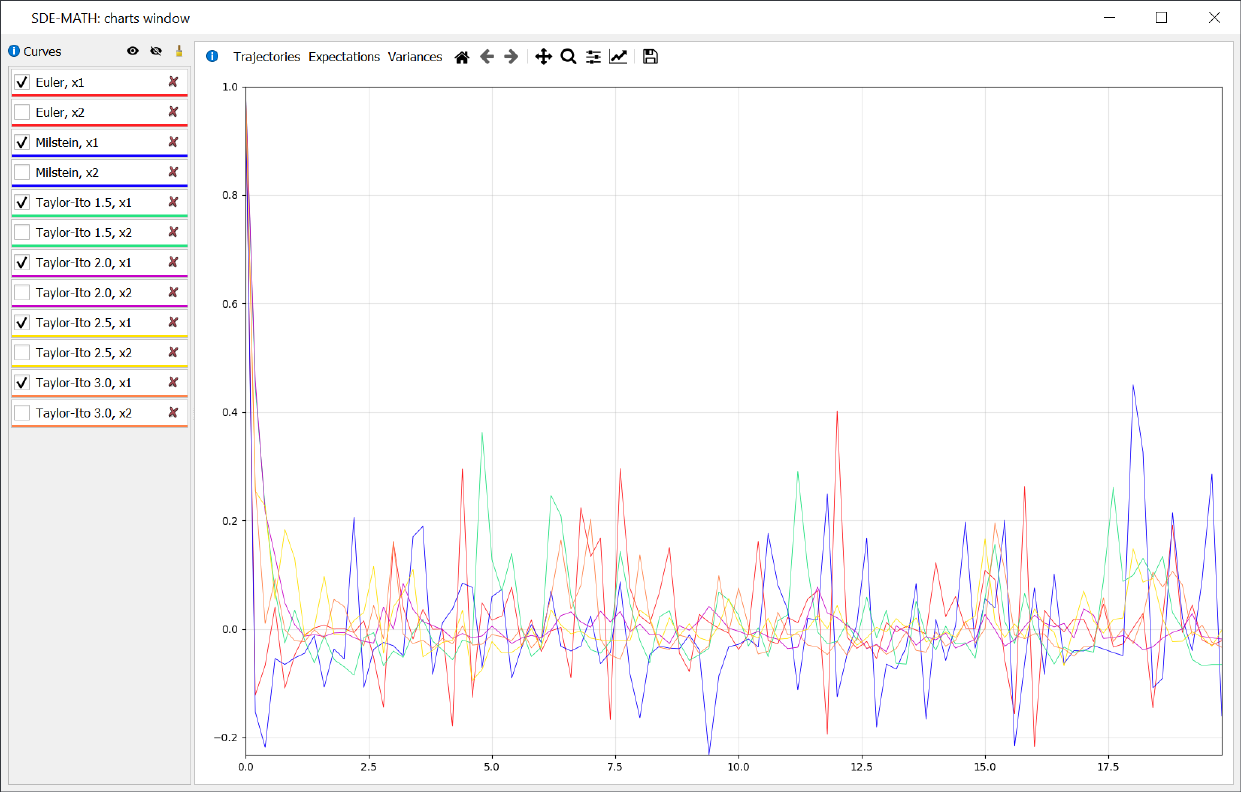}
    \caption{Strong Taylor--It\^o schemes of orders 0.5, 1.0, 1.5, 2.0, 2.5, and 3.0 (${\bf x}_t^{(1)}$ component, $C = 4,$ $dt = 0.2$)\label{fig:ito_3p0_big_7}}
\end{figure}

\begin{figure}[H]
    \centering
    \includegraphics[width=.9\textwidth]{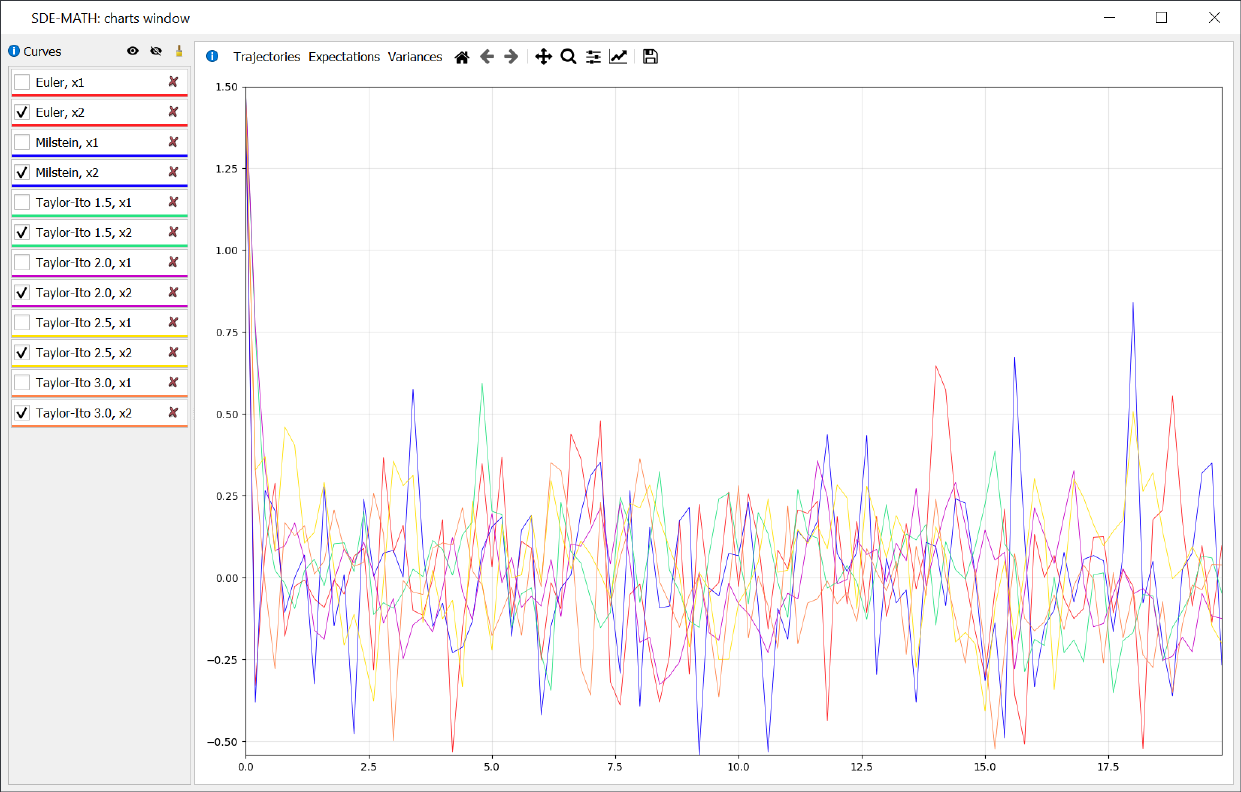}
    \caption{Strong Taylor--It\^o schemes of orders 0.5, 1.0, 1.5, 2.0, 2.5, and 3.0 (${\bf x}_t^{(2)}$ component, $C = 4,$ $dt = 0.2$)\label{fig:ito_3p0_big_8}}
\end{figure}

\begin{figure}[H]
    \vspace{10mm}
    \centering

    \hspace*{\fill}
    \begin{subfigure}[b]{.45\textwidth}
        \includegraphics[width=\textwidth]{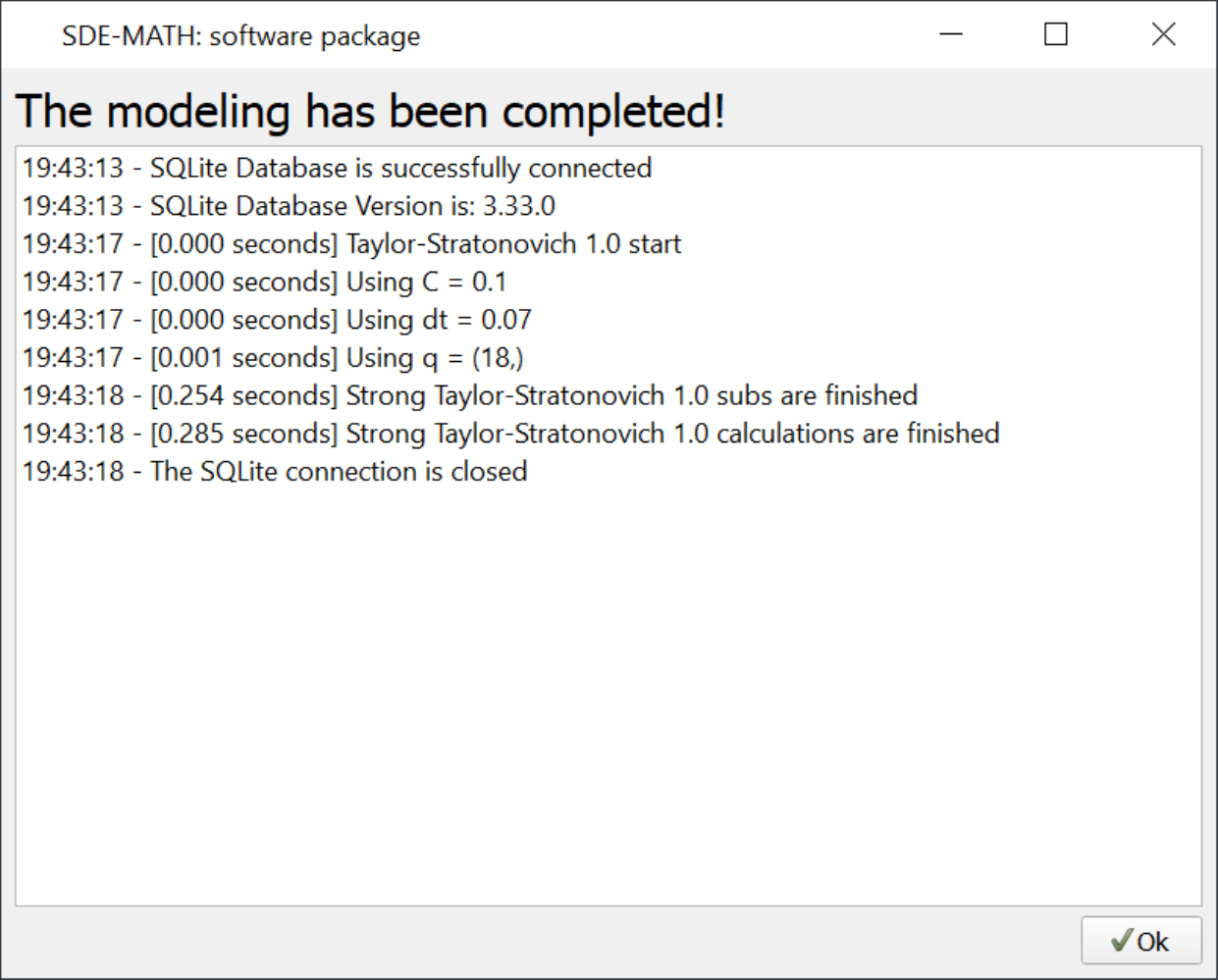}
        \caption*{Strong Taylor--Stratonovich scheme of order 1.0 ($C = 0.1,$ $dt = 0.07$)\label{fig:straton_1p5_big_1}}
    \end{subfigure}
    \hfill
    \begin{subfigure}[b]{.45\textwidth}
        \includegraphics[width=\textwidth]{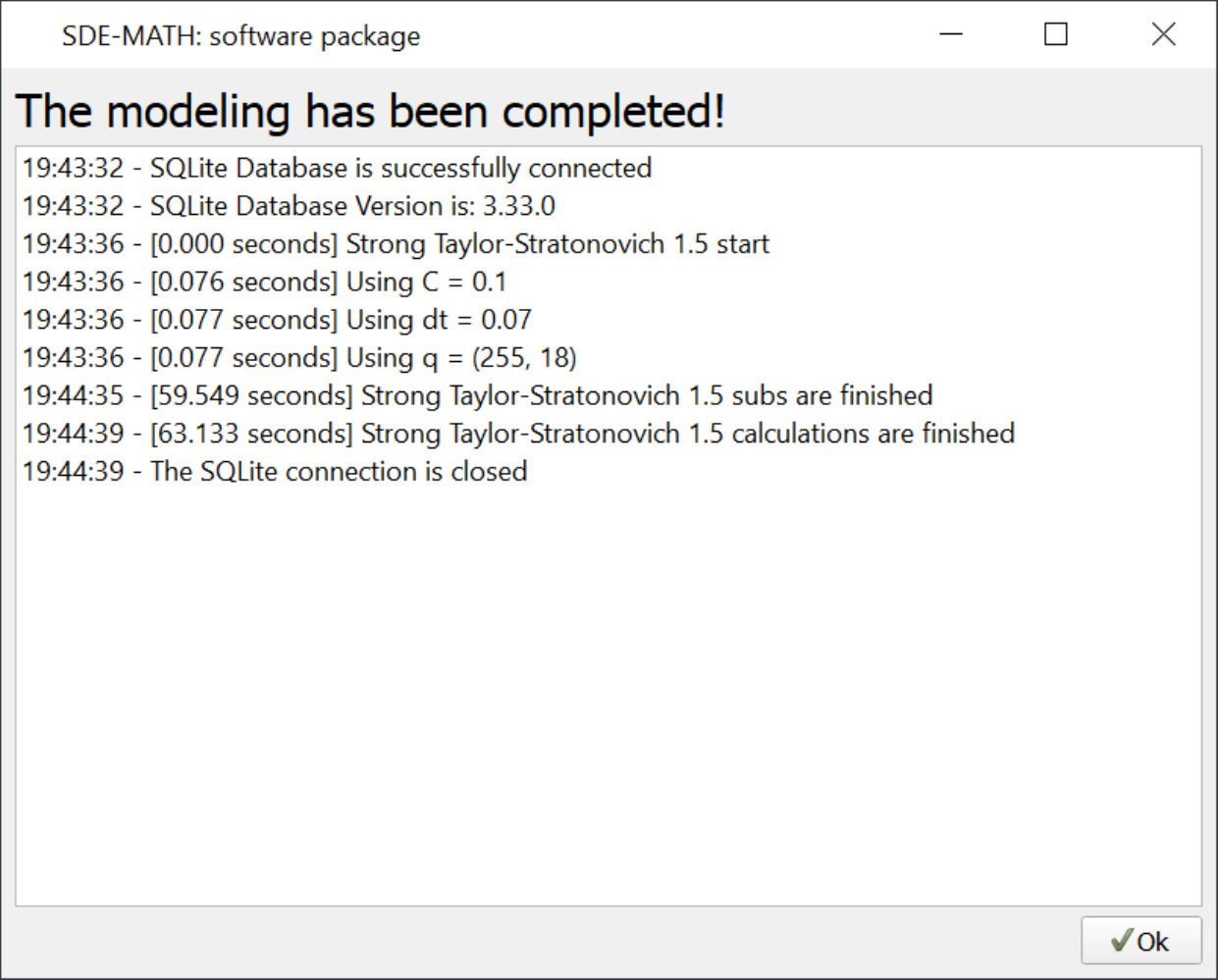}
        \caption*{Strong Taylor--Stratonovich scheme of order 1.5 ($C = 0.1,$ $dt = 0.07$)\label{fig:straton_1p5_big_2}}
    \end{subfigure}
    \hspace*{\fill}

    \caption{Modeling logs\label{fig:straton_1p5_big_logs}}

\end{figure}

\begin{figure}[H]
    \vspace{7mm}
    \centering
    \includegraphics[width=.9\textwidth]{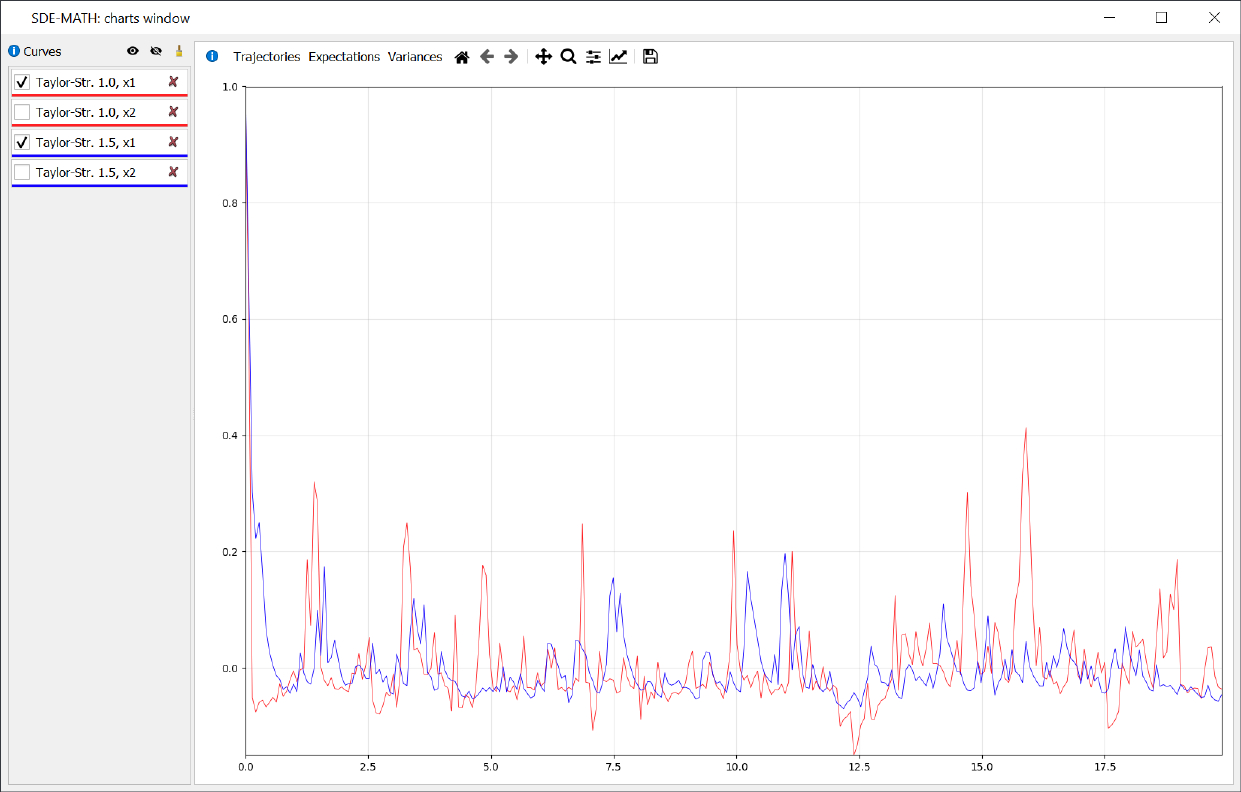}
    \caption{Strong Taylor--Stratonovich schemes of orders 1.0 and 1.5 (${\bf x}_t^{(1)}$ component, $C = 0.1,$ $dt = 0.07$)\label{fig:straton_1p5_big_3}}
\end{figure}

\begin{figure}[H]
    \vspace{10mm}
    \centering
    \includegraphics[width=.9\textwidth]{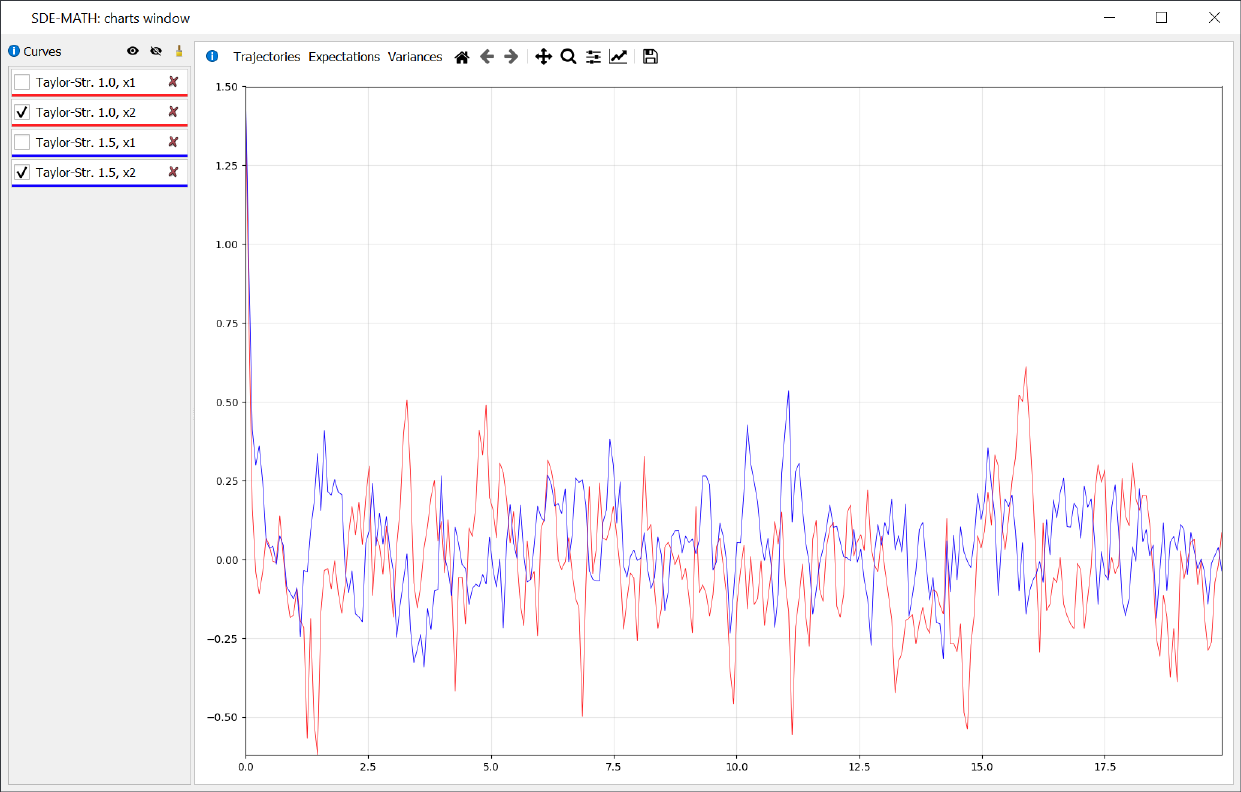}
    \caption{Strong Taylor--Stratonovich schemes of orders 1.0 and 1.5 (${\bf x}_t^{(2)}$ component, $C = 0.1,$ $dt = 0.07$)\label{fig:straton_1p5_big_4}}
\end{figure}

\begin{figure}[H]
    \vspace{7mm}
    \centering

    \hspace*{\fill}
    \begin{subfigure}[b]{.45\textwidth}
        \includegraphics[width=\textwidth]{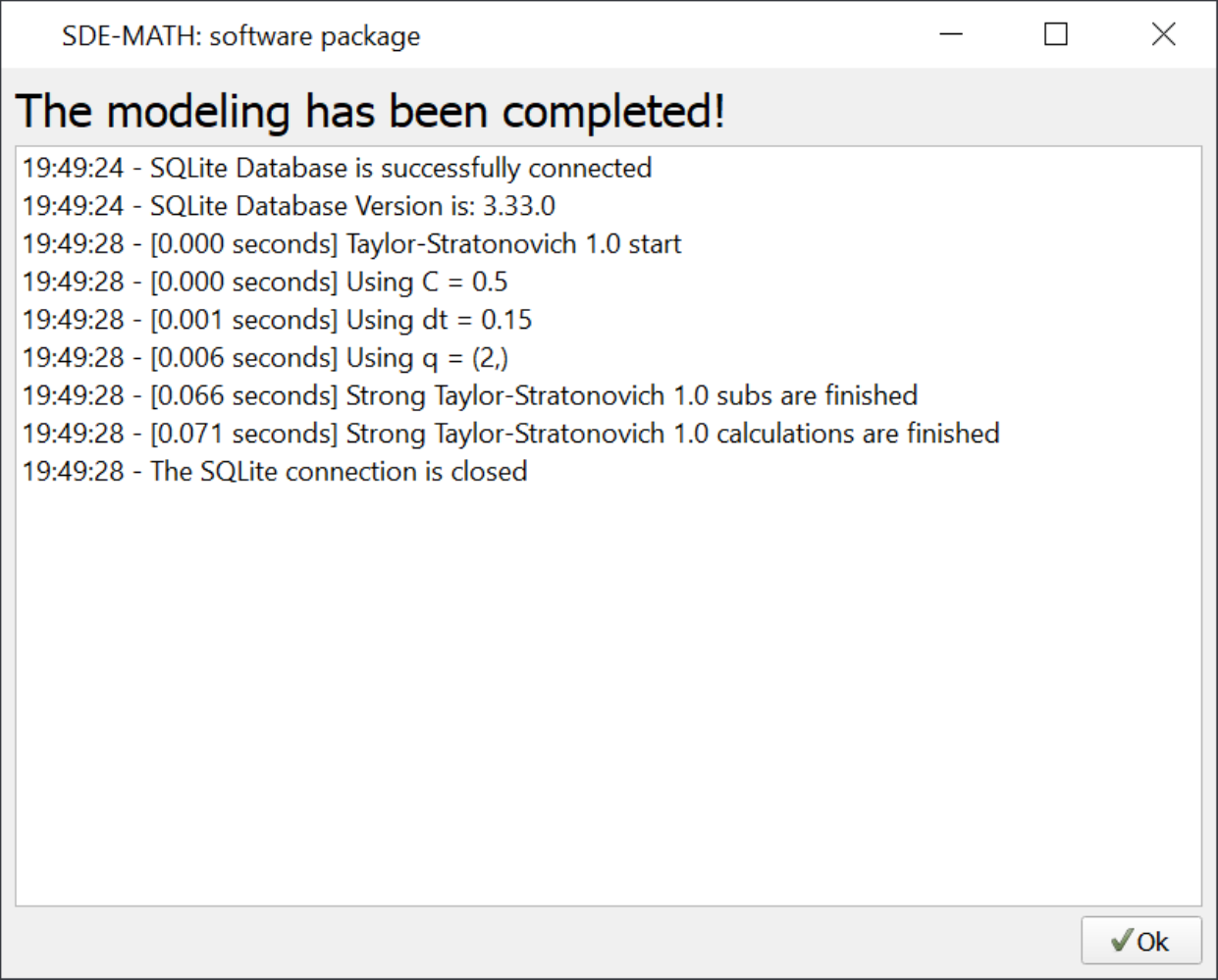}
        \caption*{Strong Taylor--Stratonovich scheme of order 1.0 ($C = 0.5,$ $dt = 0.15$)\label{fig:straton_2p0_big_1}}
    \end{subfigure}
    \hfill
    \begin{subfigure}[b]{.45\textwidth}
        \includegraphics[width=\textwidth]{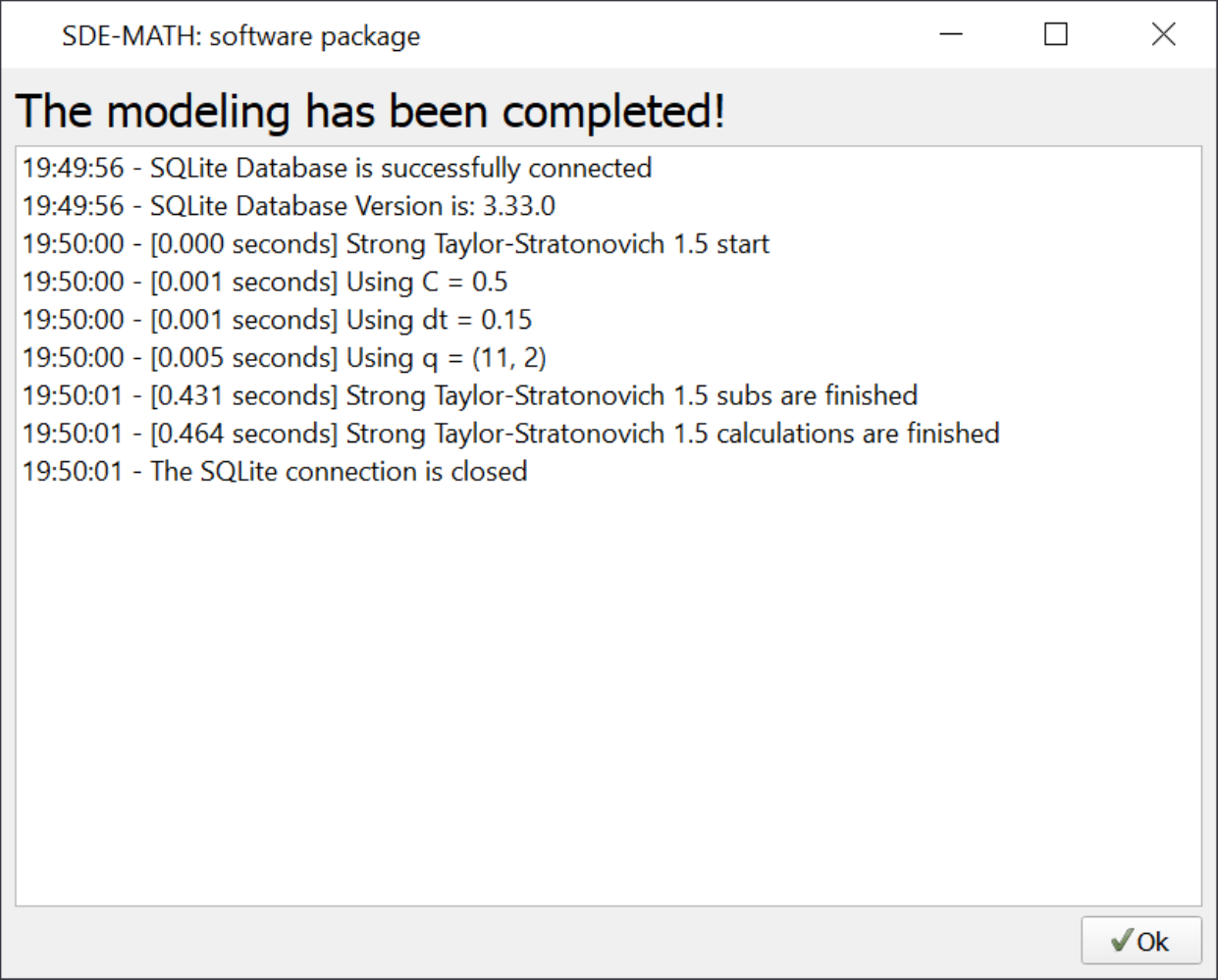}
        \caption*{Strong Taylor--Stratonovich scheme of order 1.5 ($C = 0.5,$ $dt = 0.15$)\label{fig:straton_2p0_big_2}}
    \end{subfigure}
    \hspace*{\fill}
    \caption{Modeling logs\label{fig:straton_2p0_big_logs}}

\end{figure}

\begin{figure}[H]
    \vspace{13mm}
    \centering
    \includegraphics[width=.45\textwidth]{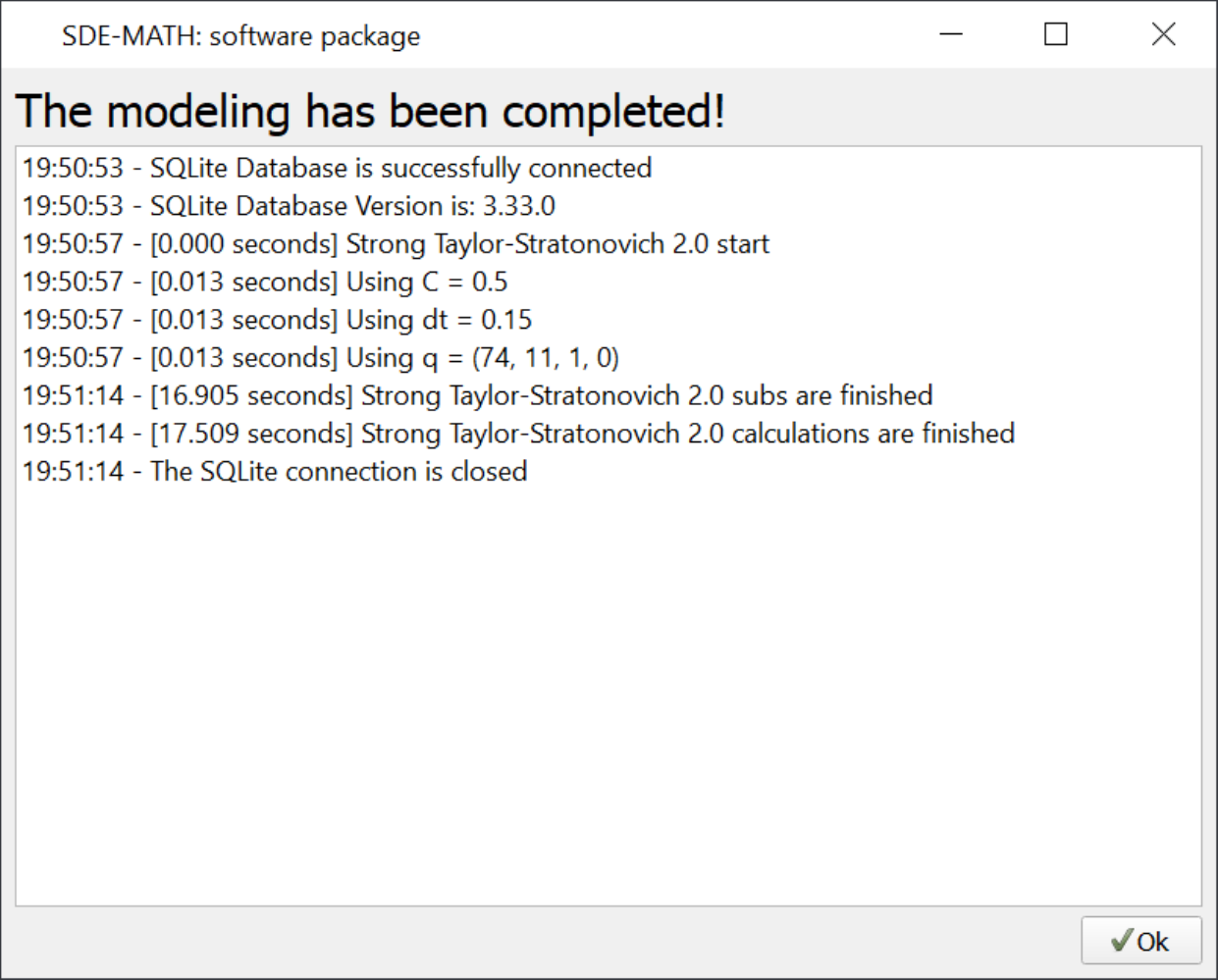}
    \caption{Strong Taylor--Stratonovich scheme of order 2.0 ($C = 0.5,$ $dt = 0.15$)\label{fig:straton_2p0_big_3}}
\end{figure}

\begin{figure}[H]
    \vspace{10mm}
    \centering
    \includegraphics[width=.9\textwidth]{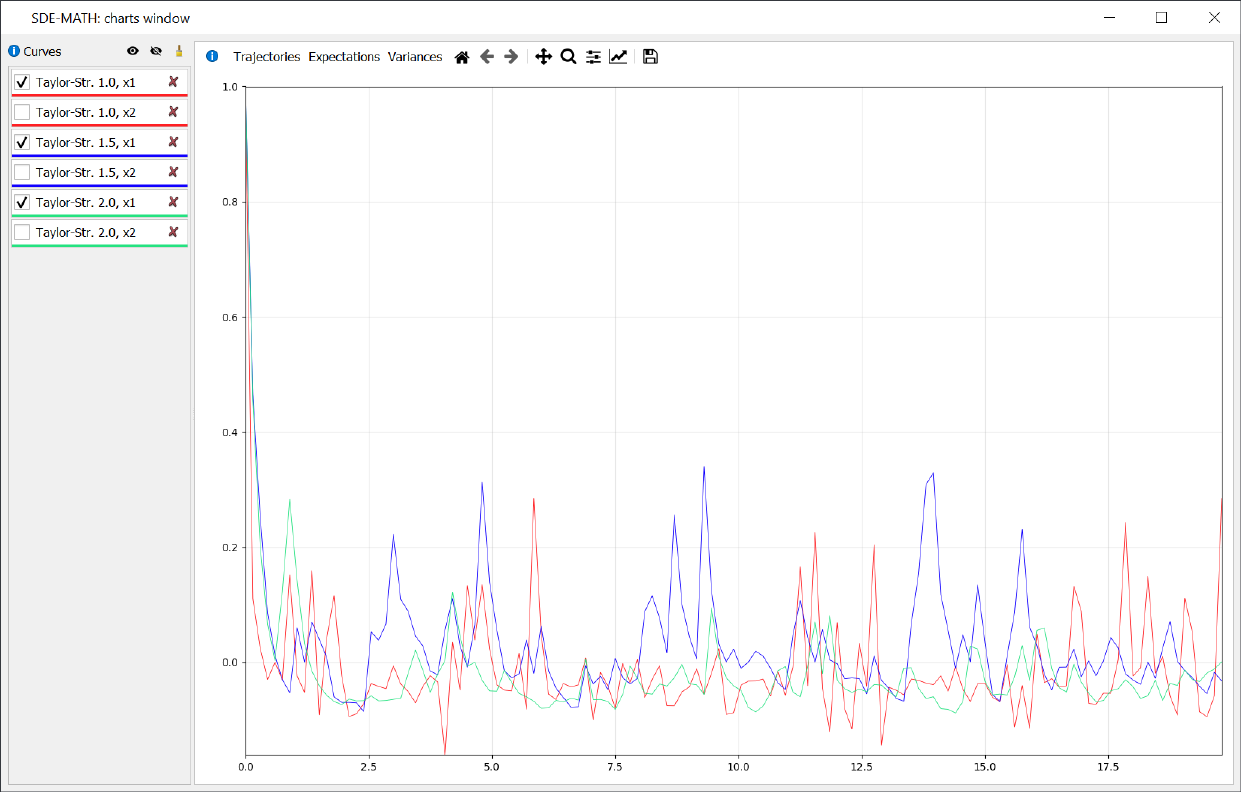}
    \caption{Strong Taylor--Stratonovich schemes of orders 1.0, 1.5, and 2.0 (${\bf x}_t^{(1)}$ component, $C = 0.5,$ $dt = 0.15$)\label{fig:straton_2p0_big_4}}
\end{figure}

\begin{figure}[H]
    \vspace{10mm}
    \centering
    \includegraphics[width=.9\textwidth]{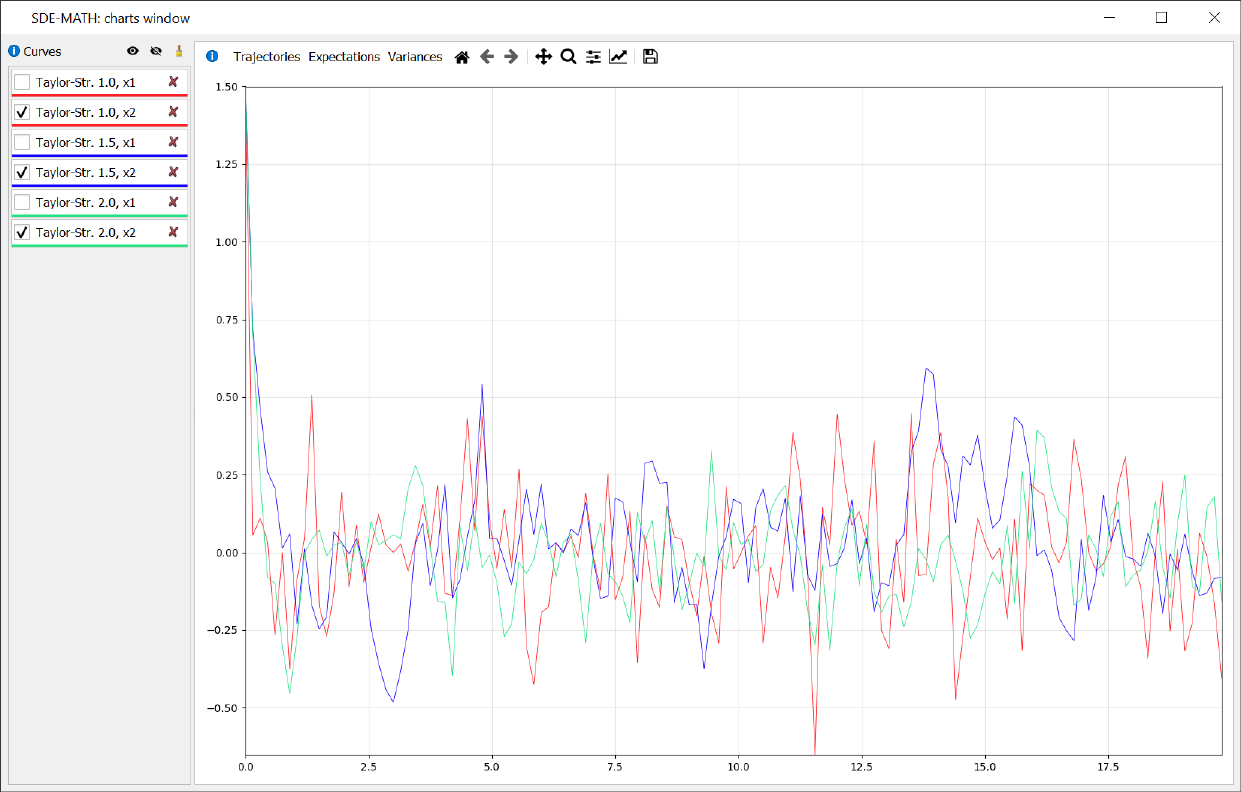}
    \caption{Strong Taylor--Stratonovich schemes of orders 1.0, 1.5, and 2.0 (${\bf x}_t^{(2)}$ component, $C = 0.5,$ $dt = 0.15$)\label{fig:straton_2p0_big_5}}
\end{figure}

\begin{figure}[H]
    \vspace{7mm}
    \centering

    \hspace*{\fill}
    \begin{subfigure}[b]{.45\textwidth}
        \includegraphics[width=\textwidth]{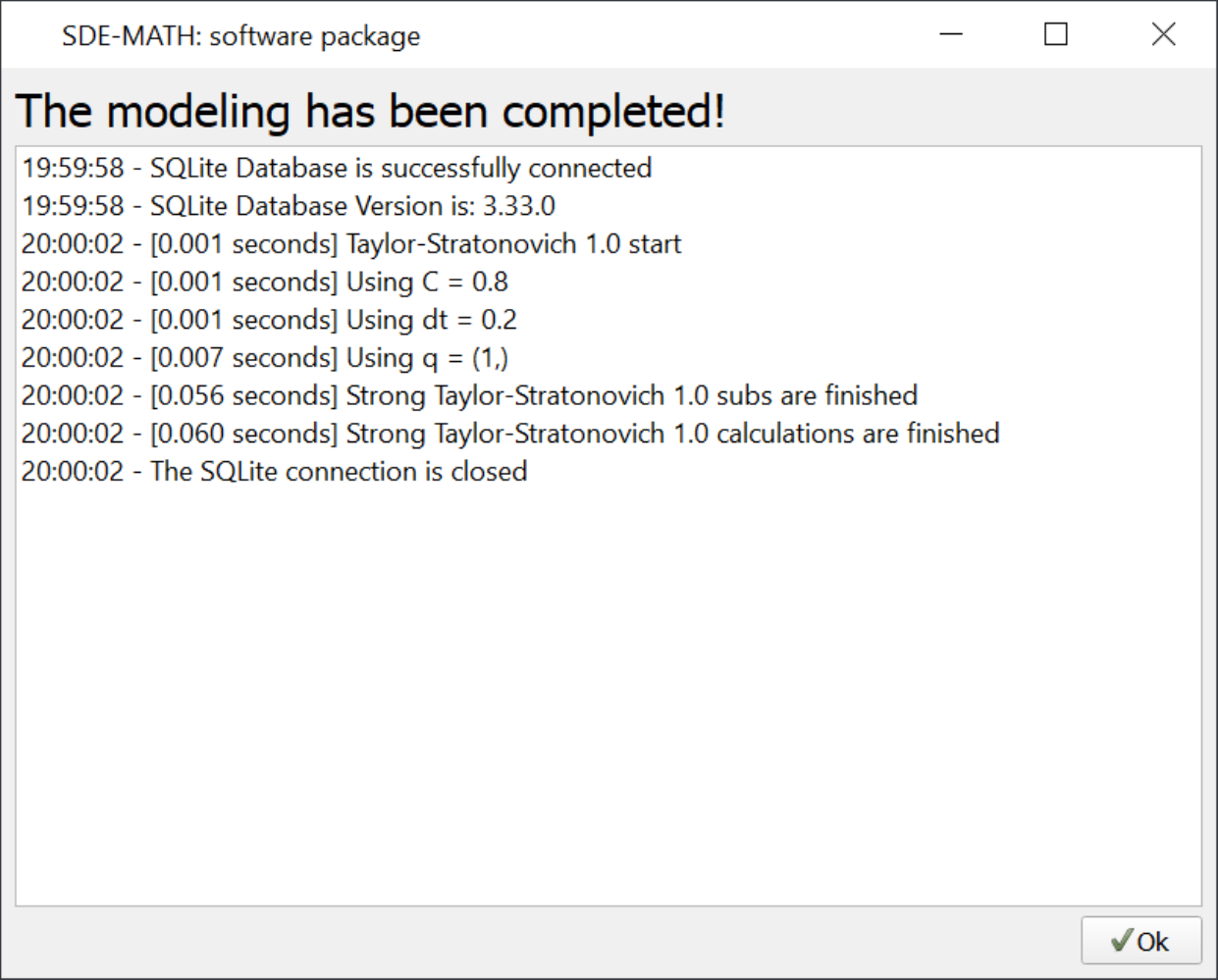}
        \caption*{Strong Taylor--Stratonovich scheme of order 1.0 ($C = 0.8,$ $dt = 0.2$)\label{fig:straton_2p5_big_1}}
    \end{subfigure}
    \hfill
    \begin{subfigure}[b]{.45\textwidth}
        \includegraphics[width=\textwidth]{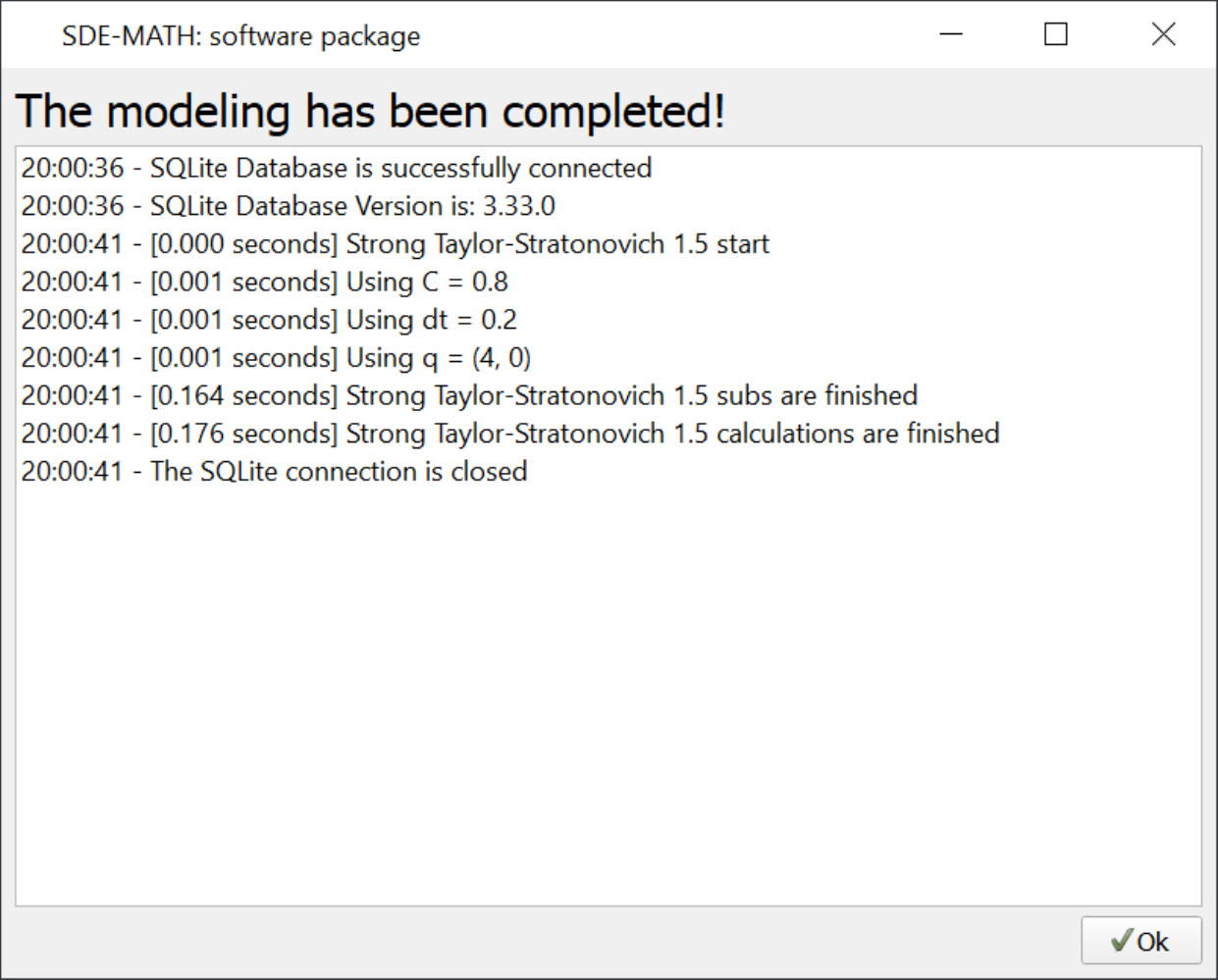}
        \caption*{Strong Taylor--Stratonovich scheme of order 1.5 ($C = 0.8,$ $dt = 0.2$)\label{fig:straton_2p5_big_2}}
    \end{subfigure}
    \hspace*{\fill}

    \caption{Modeling logs\label{fig:straton_2p5_big_logs1}}

\end{figure}

\begin{figure}[H]
    \vspace{10mm}
    \centering
    \vspace{2mm}
    \hspace*{\fill}
    \begin{subfigure}[b]{.45\textwidth}
        \includegraphics[width=\textwidth]{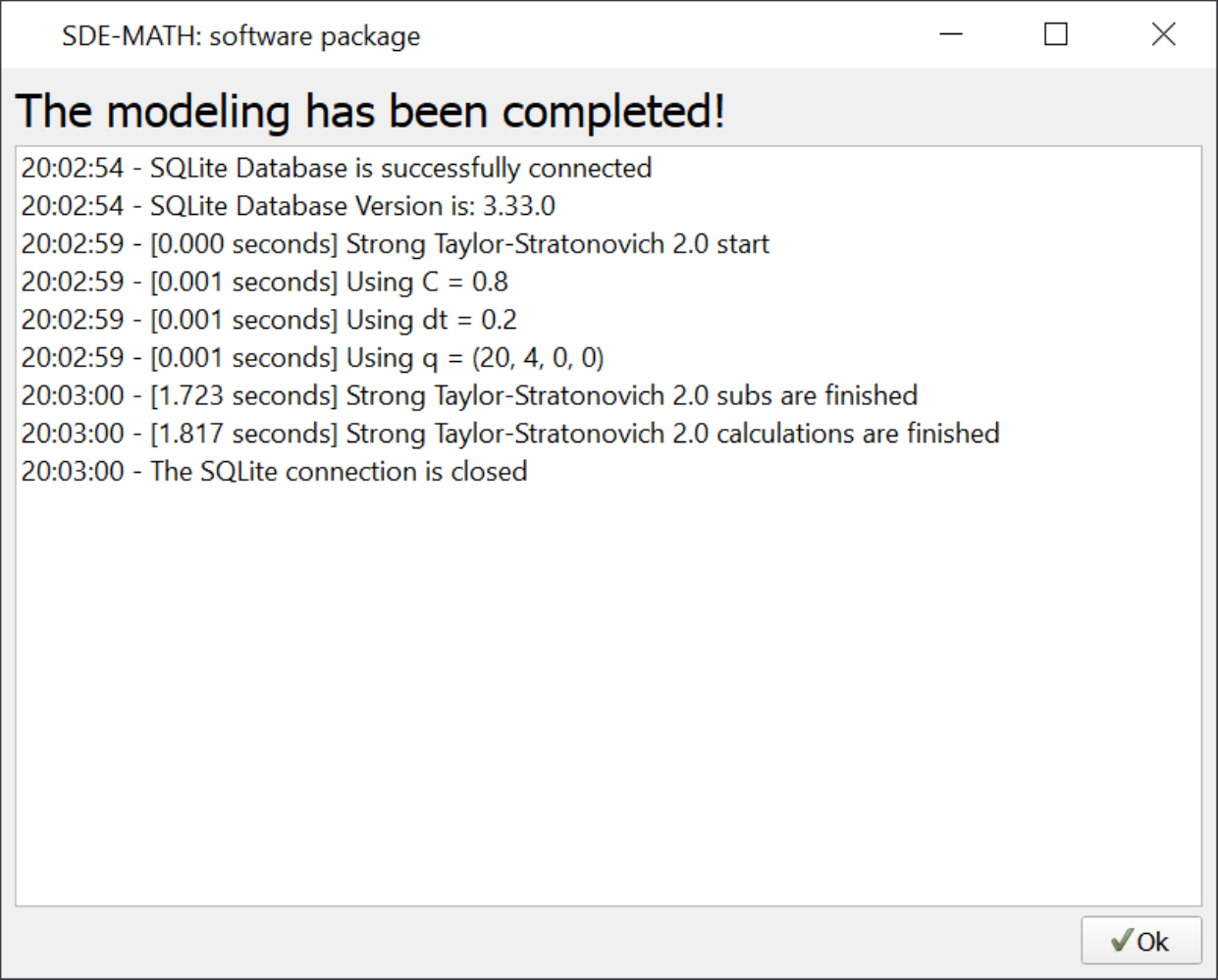}
        \caption*{Strong Taylor--Stratonovich scheme of order 2.0 ($C = 0.8,$ $dt = 0.2$)\label{fig:straton_2p5_big_3}}
    \end{subfigure}
    \hfill
    \begin{subfigure}[b]{.45\textwidth}
        \includegraphics[width=\textwidth]{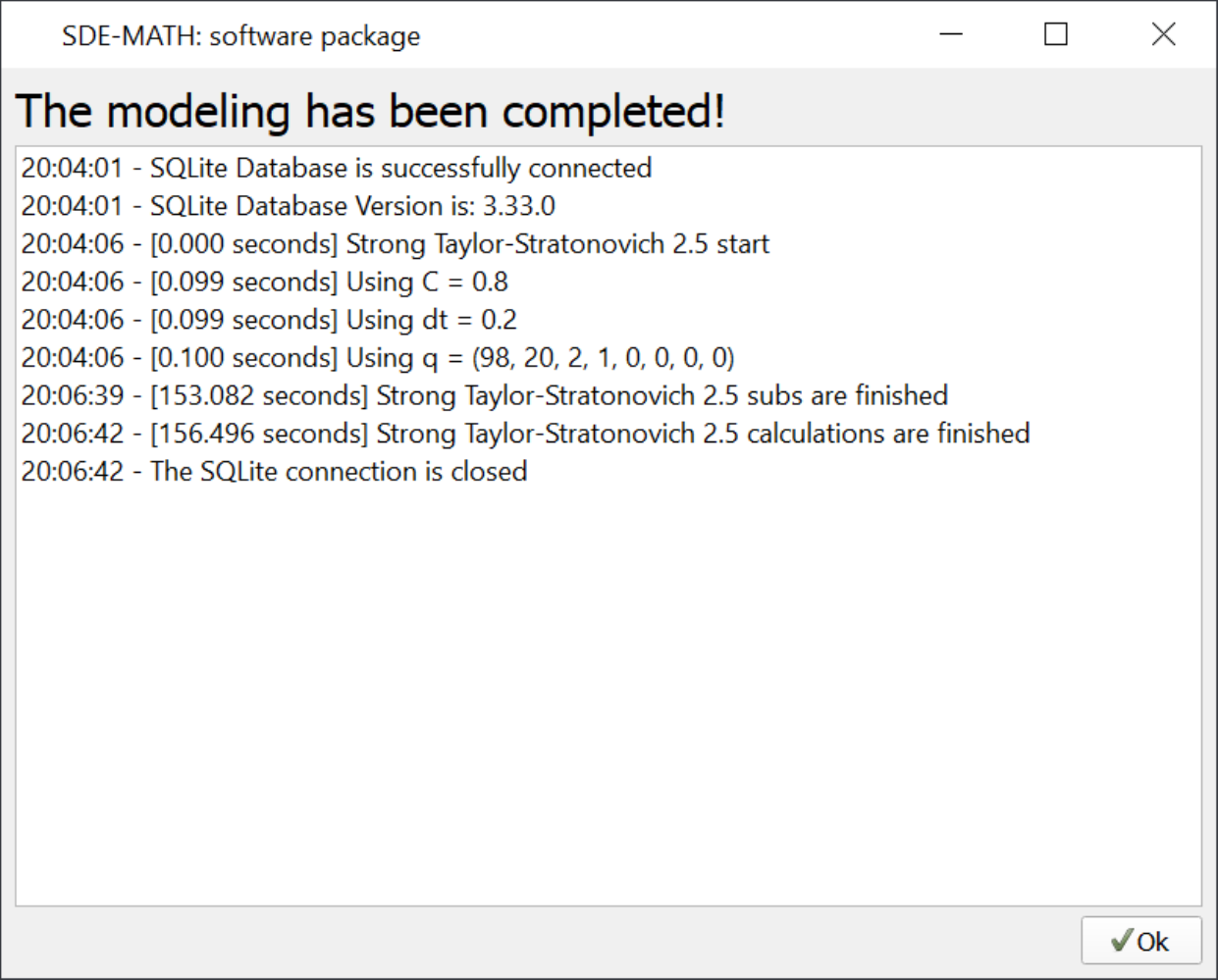}
        \caption*{Strong Taylor--Stratonovich scheme of order 2.5 ($C = 0.8,$ $dt = 0.2$)\label{fig:straton_2p5_big_4}}
    \end{subfigure}
    \hspace*{\fill}

    \caption{Modeling logs\label{fig:straton_2p5_big_logs2}}

\end{figure}

\begin{figure}[H]
    \centering

    \hspace*{\fill}
    \begin{subfigure}[b]{.45\textwidth}
        \includegraphics[width=\textwidth]{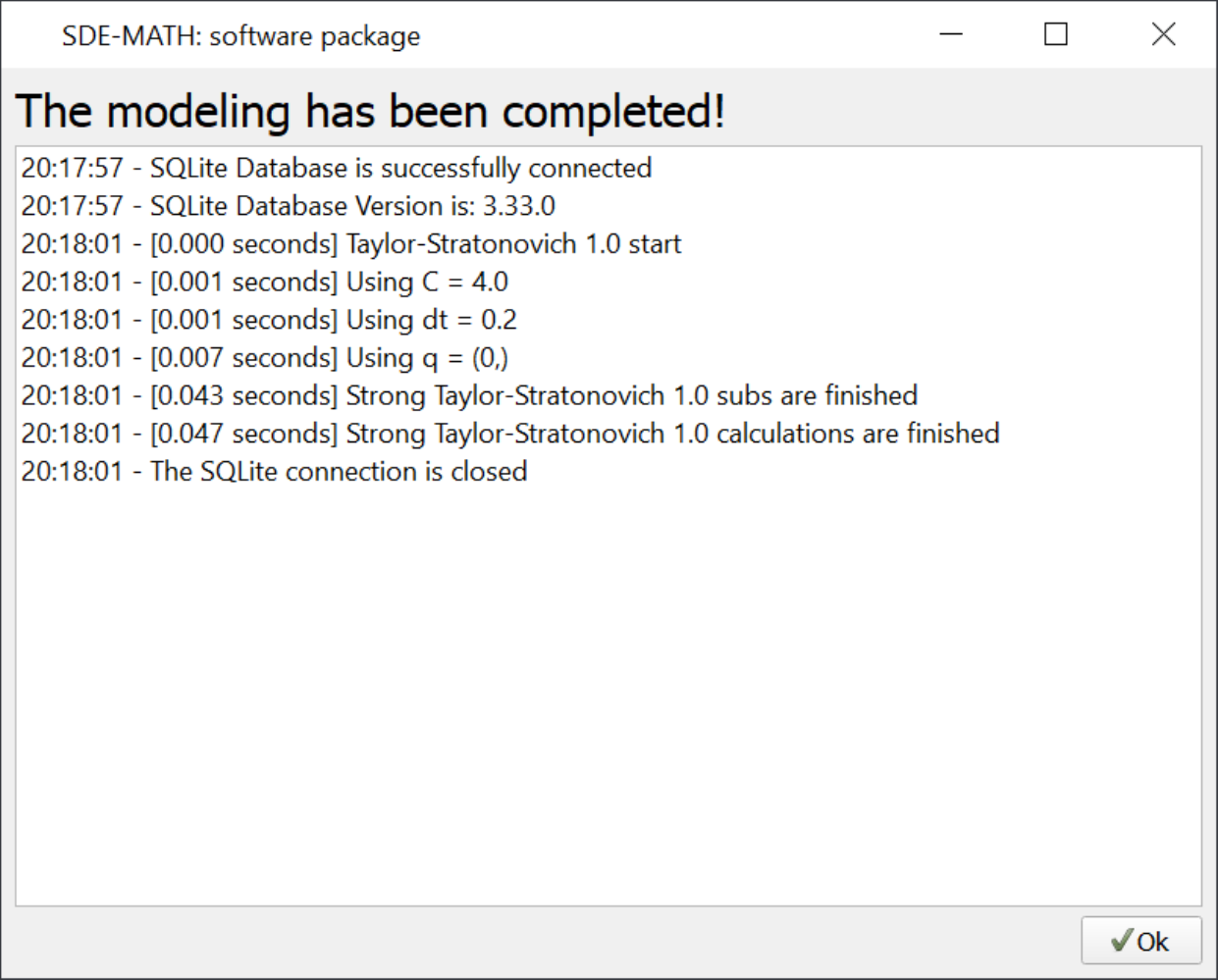}
        \caption*{Strong Taylor--Stratonovich scheme of order 1.0 ($C = 4,$ $dt = 0.2$)\label{fig:straton_3p0_big_1}}
        
    \end{subfigure}
    \hfill
    \begin{subfigure}[b]{.45\textwidth}
        \includegraphics[width=\textwidth]{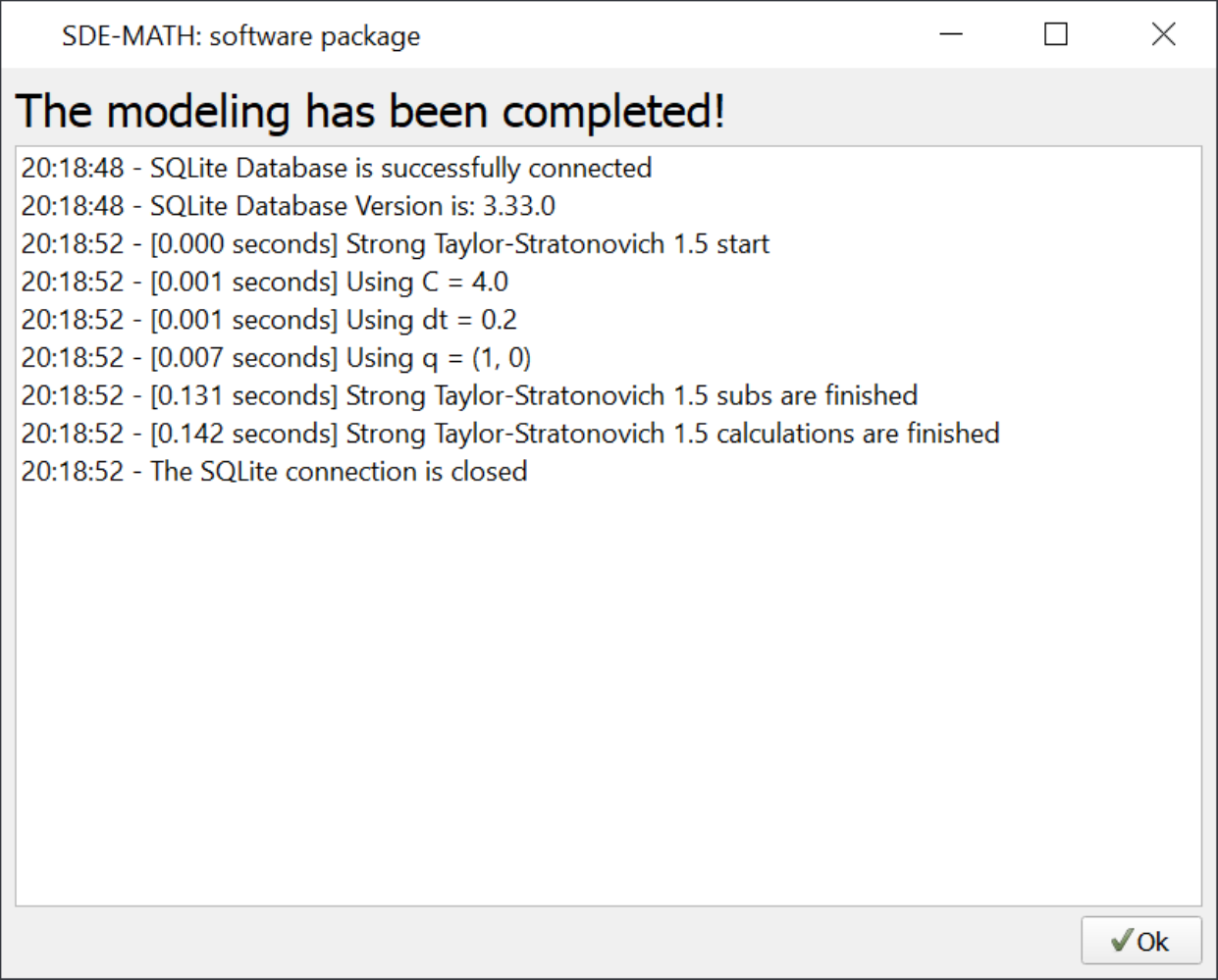}
        \caption*{Strong Taylor--Stratonovich scheme of order 1.5 ($C = 4,$ $dt = 0.2$)\label{fig:straton_3p0_big_2}}
    \end{subfigure}
    \hspace*{\fill}

    \vspace{2mm}
    \hspace*{\fill}
    \begin{subfigure}[b]{.45\textwidth}
        \includegraphics[width=\textwidth]{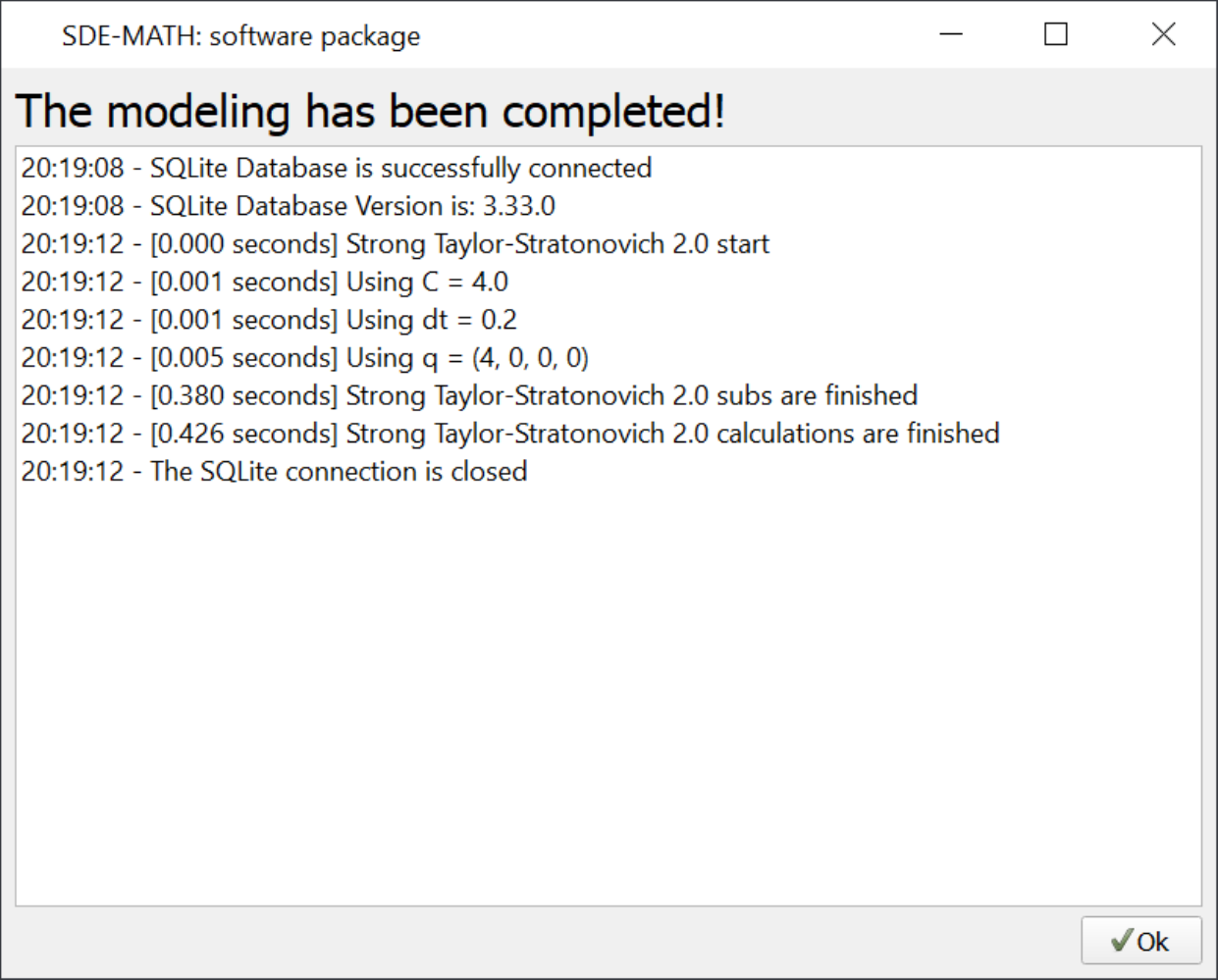}
        \caption*{Strong Taylor--Stratonovich scheme of order 2.0 ($C = 4,$ $dt = 0.2$)\label{fig:straton_3p0_big_3}}
    \end{subfigure}
    \hfill
    \begin{subfigure}[b]{.45\textwidth}
        \includegraphics[width=\textwidth]{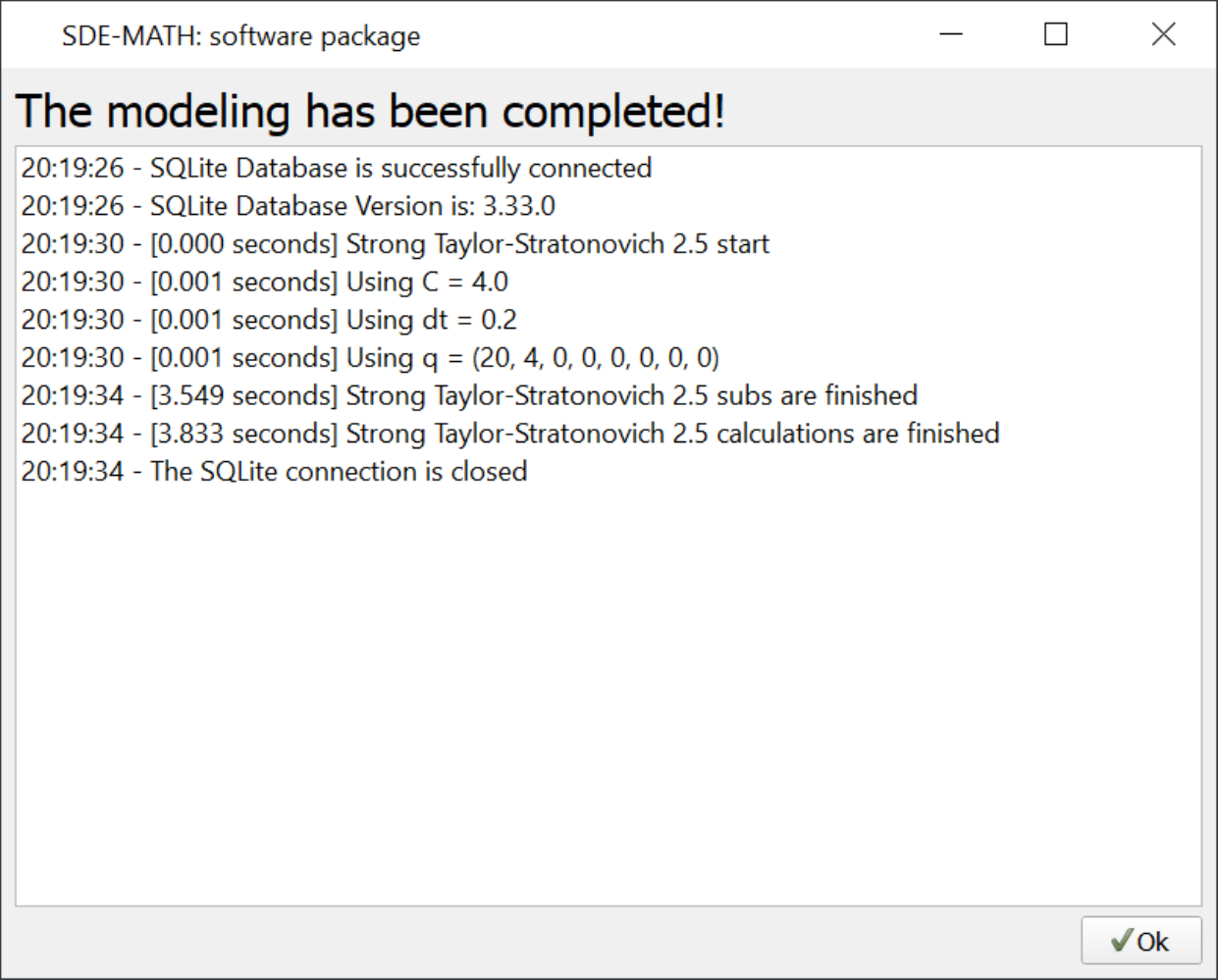}
        \caption*{Strong Taylor--Stratonovich scheme of order 2.5 ($C = 4,$ $dt = 0.2$)\label{fig:straton_3p0_big_4}}
    \end{subfigure}
    \hspace*{\fill}

    \vspace{2mm}
    \hspace*{\fill}
    \begin{subfigure}[b]{.45\textwidth}
        \includegraphics[width=\textwidth]{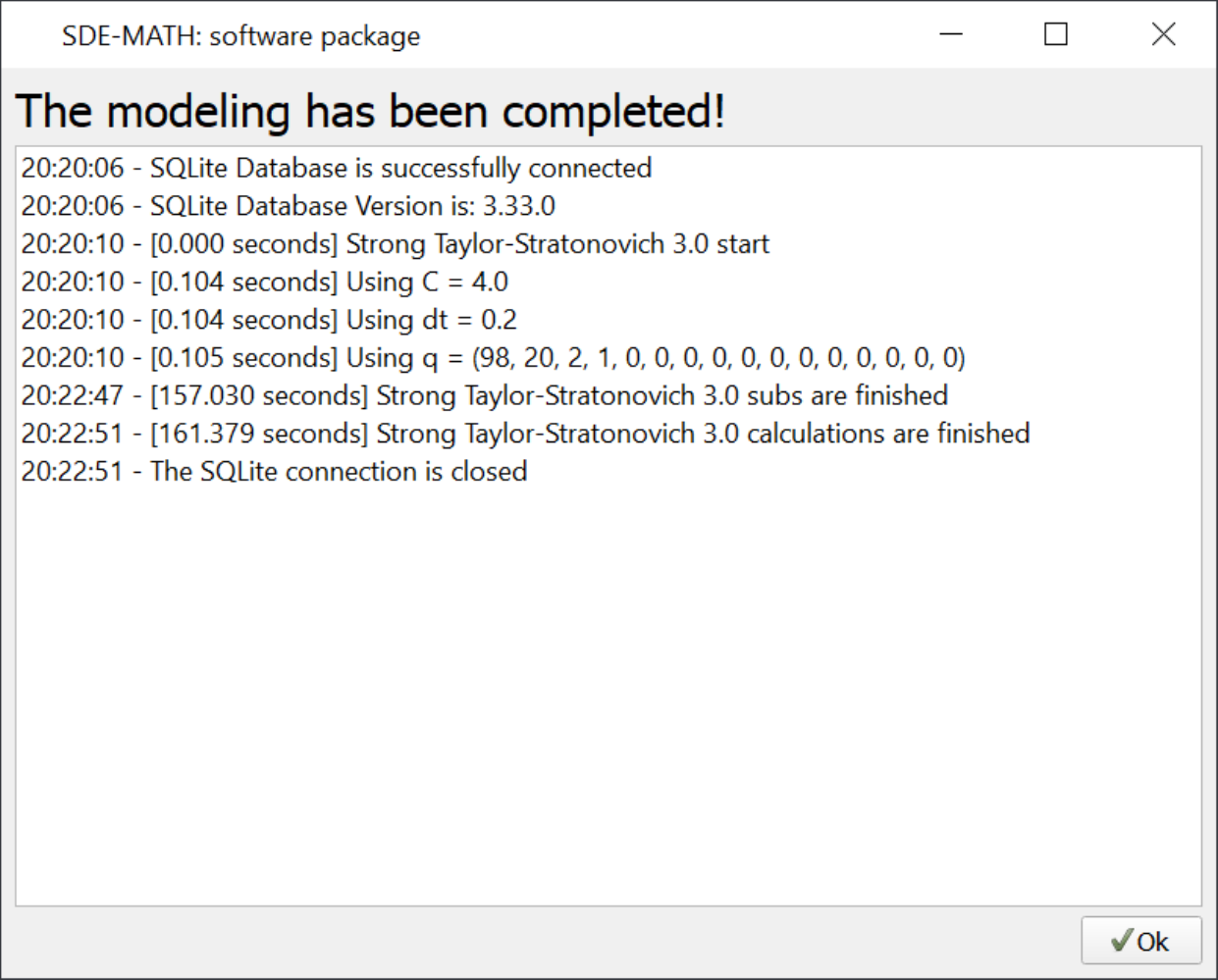}
        \caption*{Strong Taylor--Stratonovich scheme of order 3.0 ($C = 4,$ $dt = 0.2$)\label{fig:straton_3p0_big_5}}
    \end{subfigure}
    \hspace*{\fill}

    \caption{Modeling logs\label{fig:straton_3p0_big_logs}}

\end{figure}

This subsection is fully devoted to modeling logs and results visualization. 
They are presented on Figures \ref{fig:ito_1p5_small_logs}--\ref{fig:straton_3p0_big_logs}

%
%

\vspace{5mm}

\begin{figure}[H]
    \centering
    \includegraphics[width=.9\textwidth]{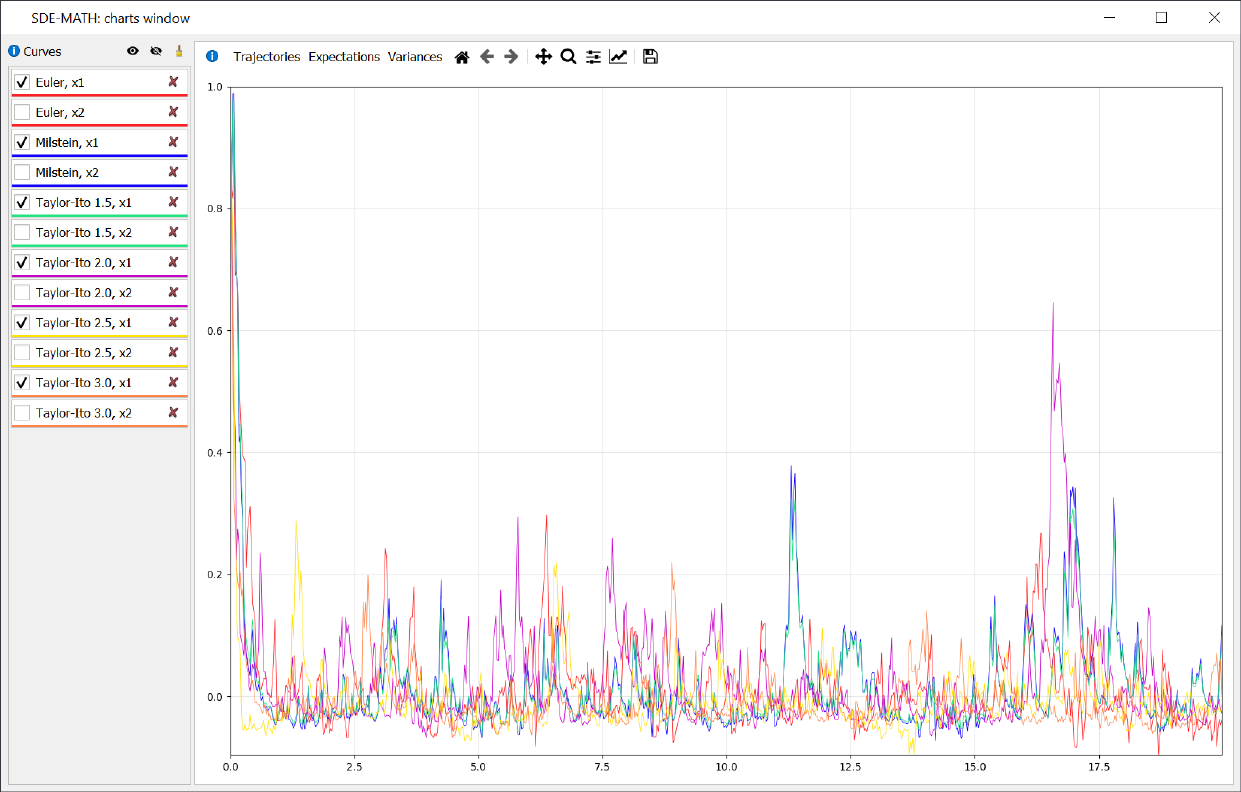}
    \caption{Strong Taylor--It\^o schemes of orders 0.5, 1.0, 1.5, 2.0, 2.5, and 3.0 (${\bf x}_t^{(1)}$ component, $C = 14000,$ $dt = 0.025$)\label{fig:ito_3p0_small_6}}
\end{figure}

\begin{figure}[H]
    \vspace{10mm}
    \centering
    \includegraphics[width=.9\textwidth]{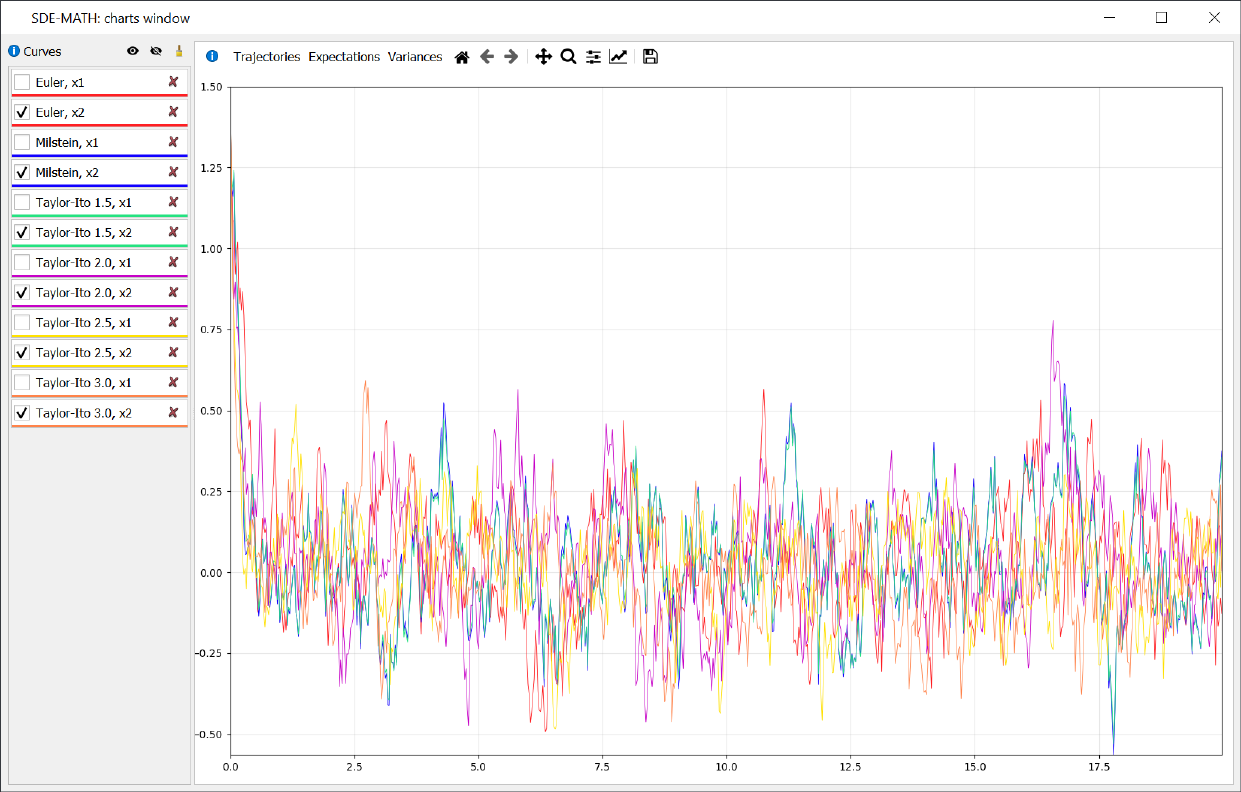}
    \caption{Strong Taylor--It\^o schemes of orders 0.5, 1.0, 1.5, 2.0, 2.5, and 3.0 (${\bf x}_t^{(2)}$ component, $C = 14000,$ $dt = 0.025$)\label{fig:ito_3p0_small_6}}
\end{figure}

%
%

%
%

%
%

\begin{figure}[H]
    \vspace{7mm}
    \centering
    \includegraphics[width=.9\textwidth]{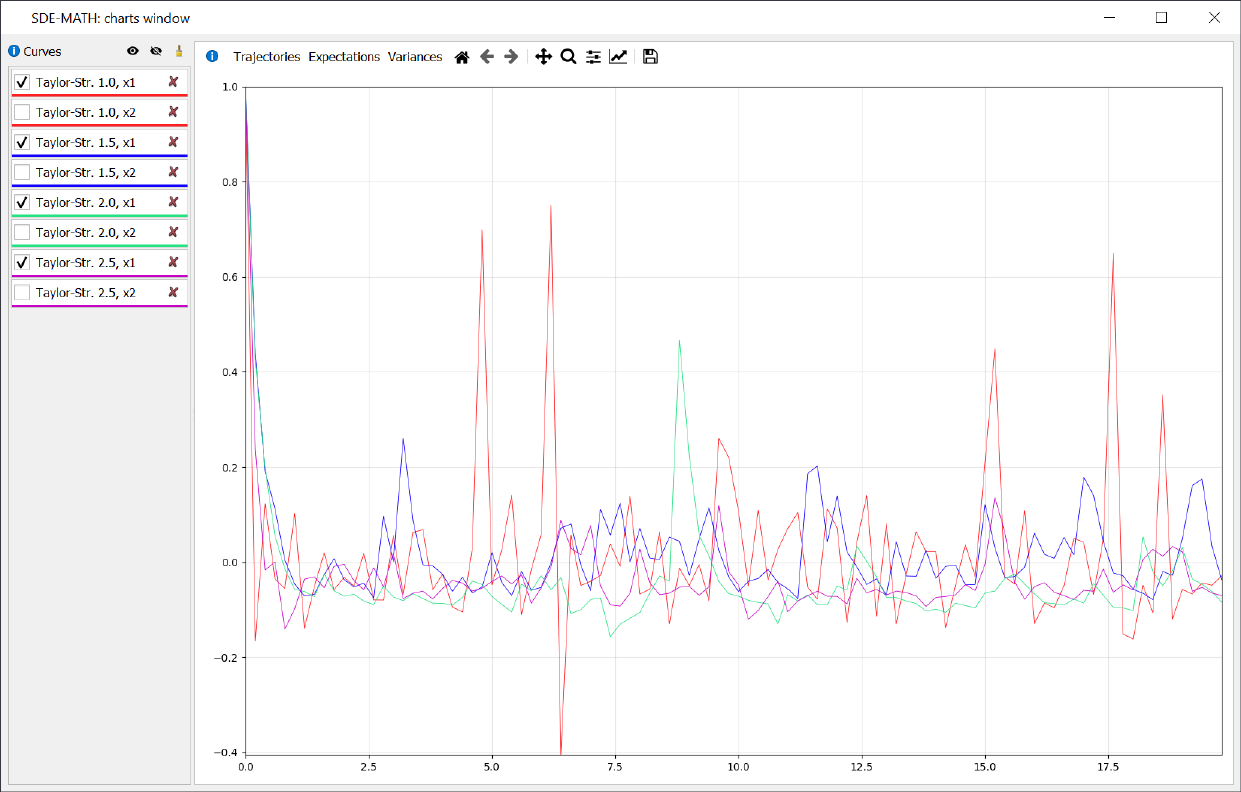}
    \caption{Strong Taylor--Stratonovich schemes of orders 1.0, 1.5, 2.0, and 2.5 (${\bf x}_t^{(1)}$ component, $C = 0.8,$ $dt = 0.2$)\label{fig:straton_2p5_big_5}}
\end{figure}

\begin{figure}[H]
    \centering
    \includegraphics[width=.9\textwidth]{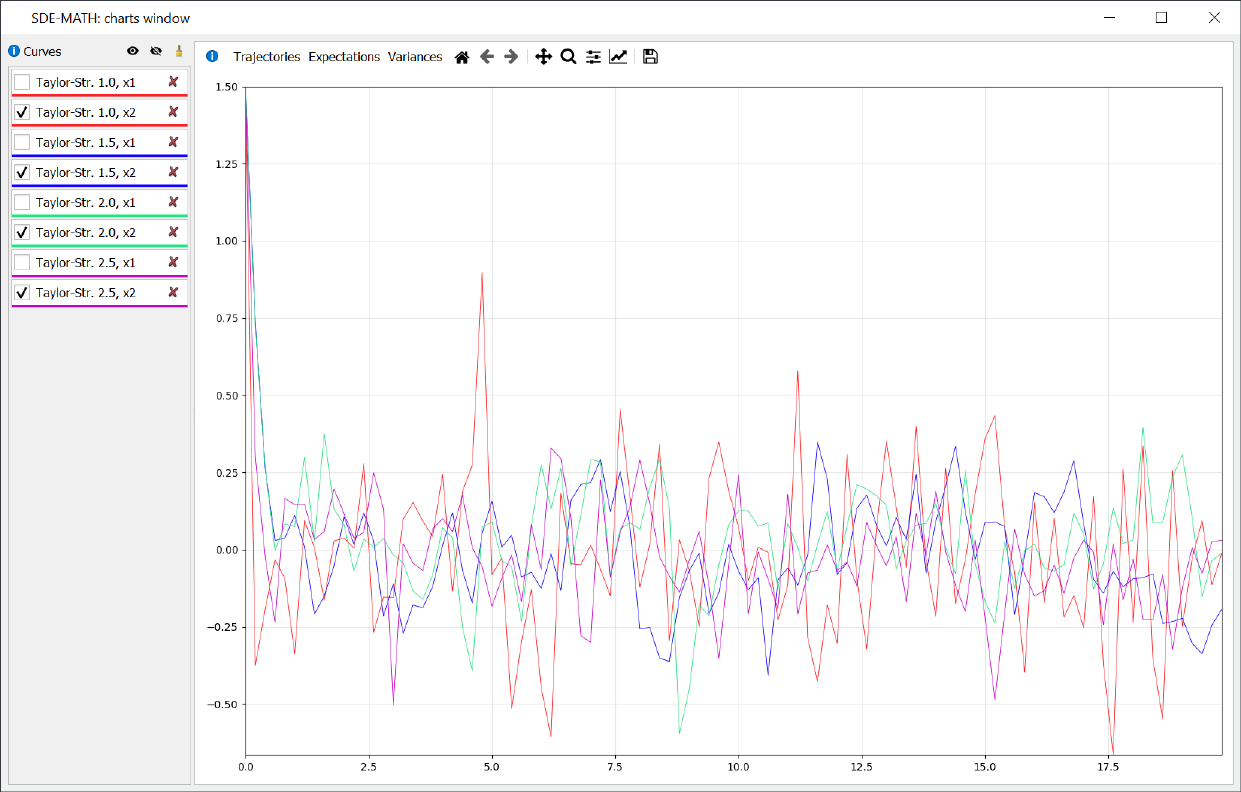}
    \caption{Strong Taylor--Stratonovich schemes of orders 1.0, 1.5, 2.0, and 2.5 (${\bf x}_t^{(2)}$ component, $C = 0.8,$ $dt = 0.2$)\label{fig:straton_2p5_big_6}}
\end{figure}

\begin{figure}[H]
    \centering
    \includegraphics[width=.9\textwidth]{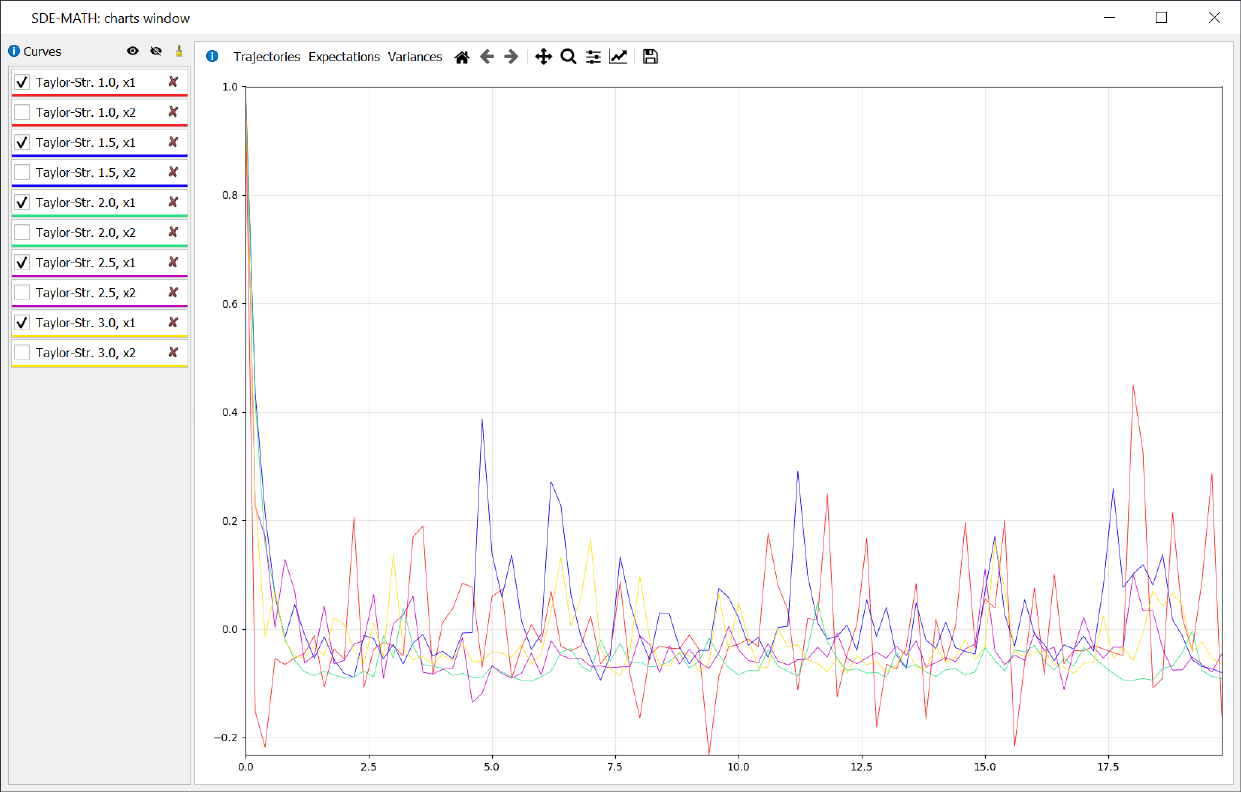}
    \caption{Strong Taylor--Stratonovich schemes of orders 1.0, 1.5, 2.0, 2.5, and 3.0 (${\bf x}_t^{(1)}$ component, $C = 4,$ $dt = 0.2$)\label{fig:straton_3p0_big_6}}
\end{figure}

\subsection{Example of Linear System of It\^o SDEs (Solar Activity)}


Consider a mathematical model of the solar activity without its average value
 in a form of the system of linear It\^{o} SDEs (\ref{lin1}) \cite{4}. In (\ref{lin1}) 
we choose \cite{4} $n=2,$ $m=1,$ $k=2,$ ${\bf x}_0^{(1)}=7,$ ${\bf x}_0^{(2)}=-0.25,$

\begin{equation}
    \label{lol1}
    A = \left(\begin{matrix}
        0 & 1\\\\
        -0.3205 & -0.14
    \end{matrix}\right),\ \ 
    B = \left(\begin{matrix}
        ~0 & 0~\\\\
        ~0 & 0~
    \end{matrix}\right),
\end{equation}

\begin{equation}
    \label{lol2}
    {\bf u}(t) \equiv \left(\begin{matrix}
        ~0~\\\\
        ~0~
    \end{matrix}\right),\ \ 
    F = \left(\begin{matrix}
        0\\\\
        5.08
    \end{matrix}\right).
\end{equation}

\subsection{Visualization and Numerical Results for Solar Activity Model}

This subsection is devoted to the visualization and numerical results 
for the model of solar activity (\ref{lin1}), (\ref{lol1}), (\ref{lol2}).

\begin{figure}[H]
    \vspace{13mm}
    \centering
    \includegraphics[width=.45\textwidth]{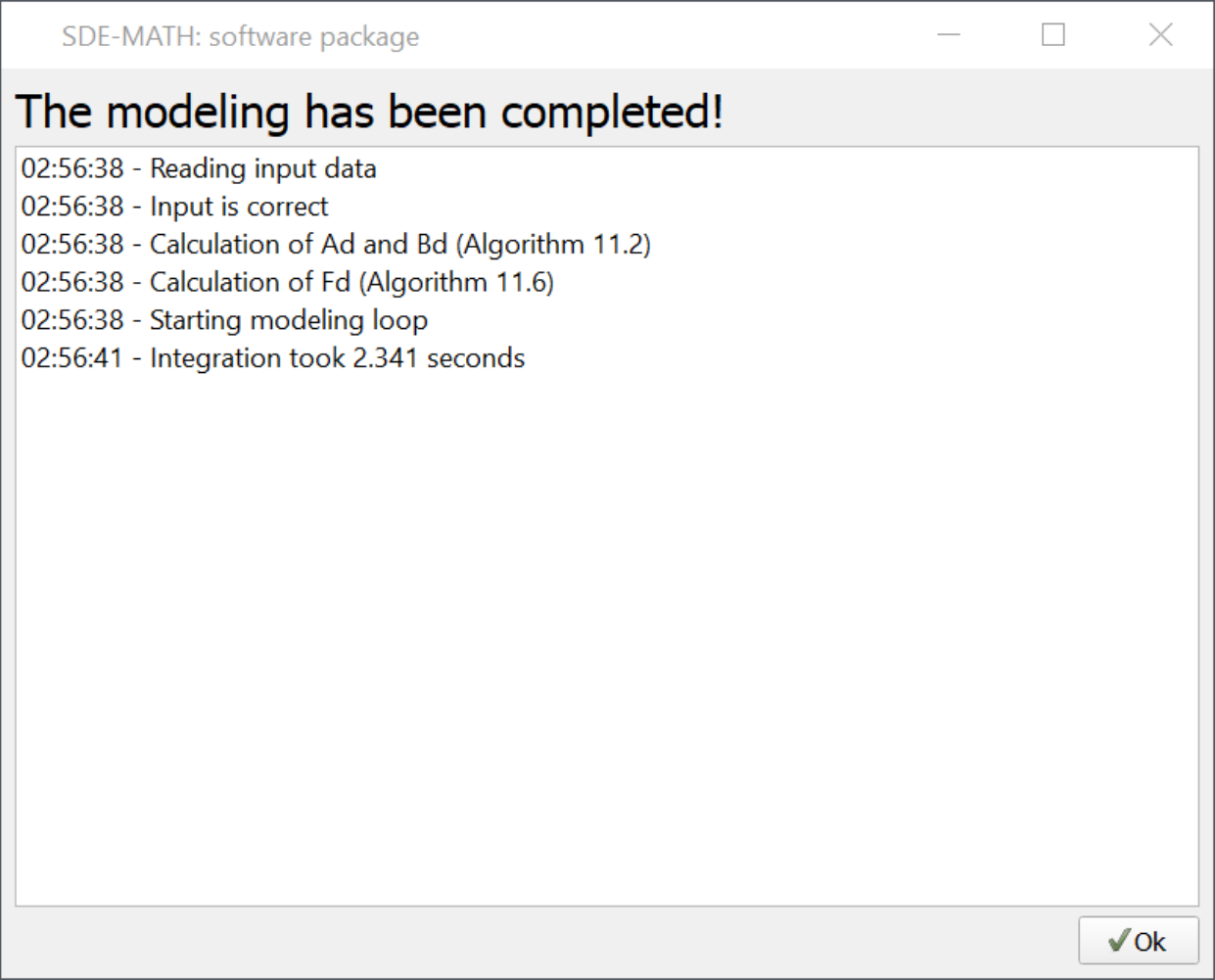}
    \caption{Modeling logs for solar activity model\label{fig:solar_1}}
\end{figure}

\begin{figure}[H]
    \vspace{10mm}
    \centering
    \includegraphics[width=.9\textwidth]{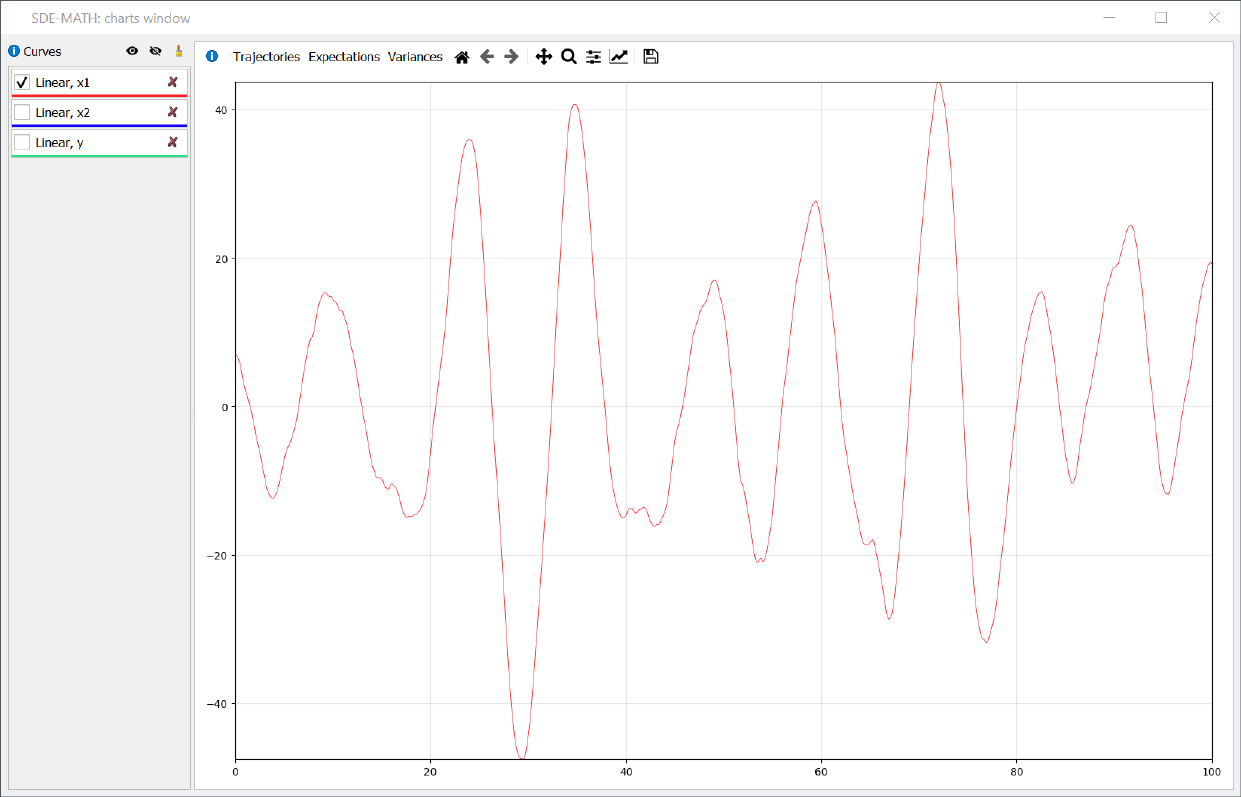}
    \caption{Solar activity model ($x^{(1)}_t$ component)\label{fig:solar_2}}
\end{figure}

\begin{figure}[H]
    \vspace{5mm}
    \centering
    \includegraphics[width=.9\textwidth]{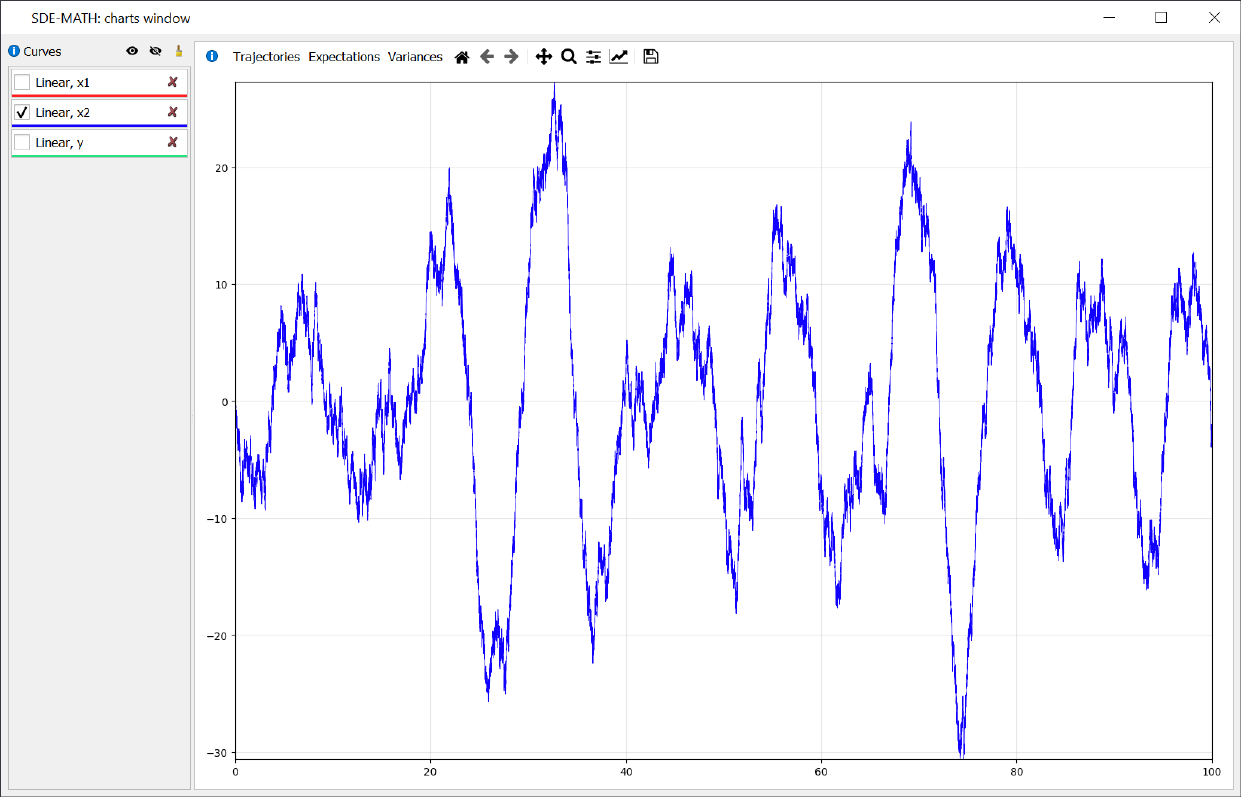}
    \caption{Solar activity model ($x^{(2)}_t$ component)\label{fig:solar_3}}
\end{figure}

\begin{figure}[H]
    \centering
    \includegraphics[width=.9\textwidth]{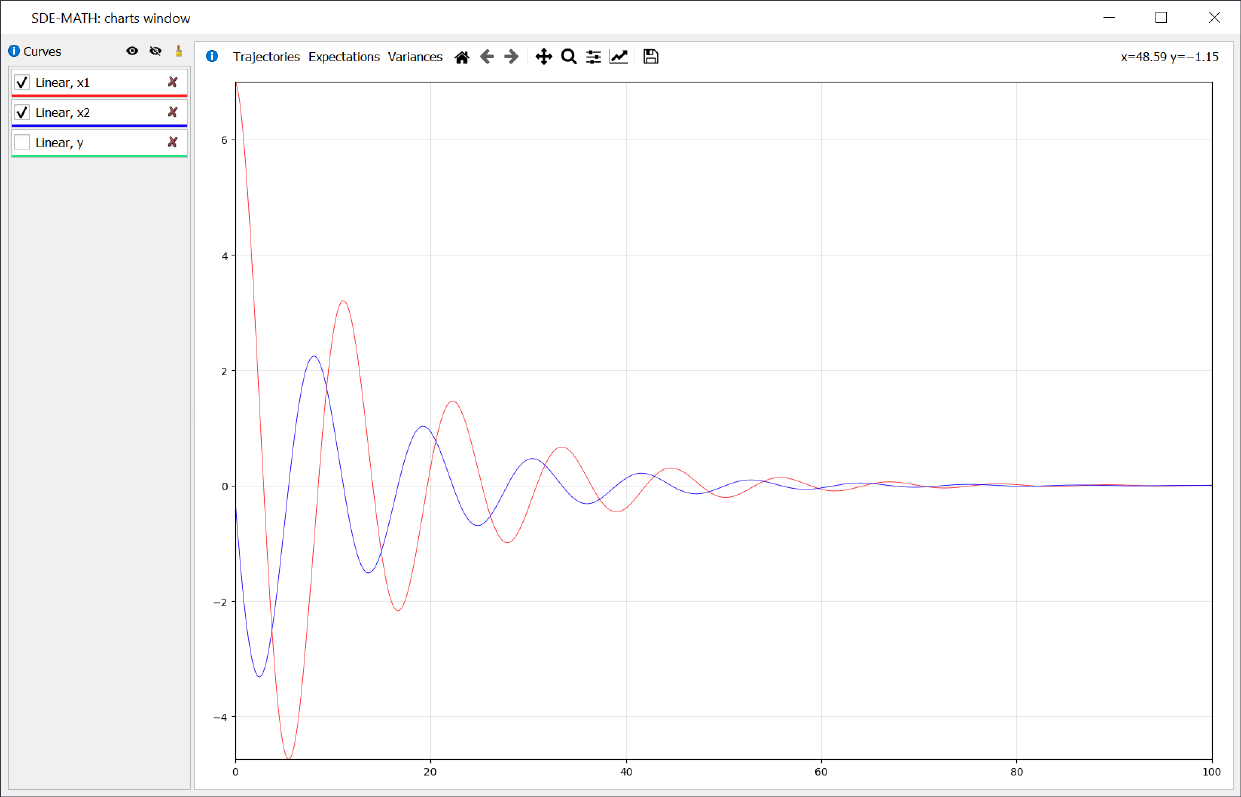}
    \caption{Solar activity model (expectations)\label{fig:solar_4}}
\end{figure}

\begin{figure}[H]
    \centering
    \includegraphics[width=.9\textwidth]{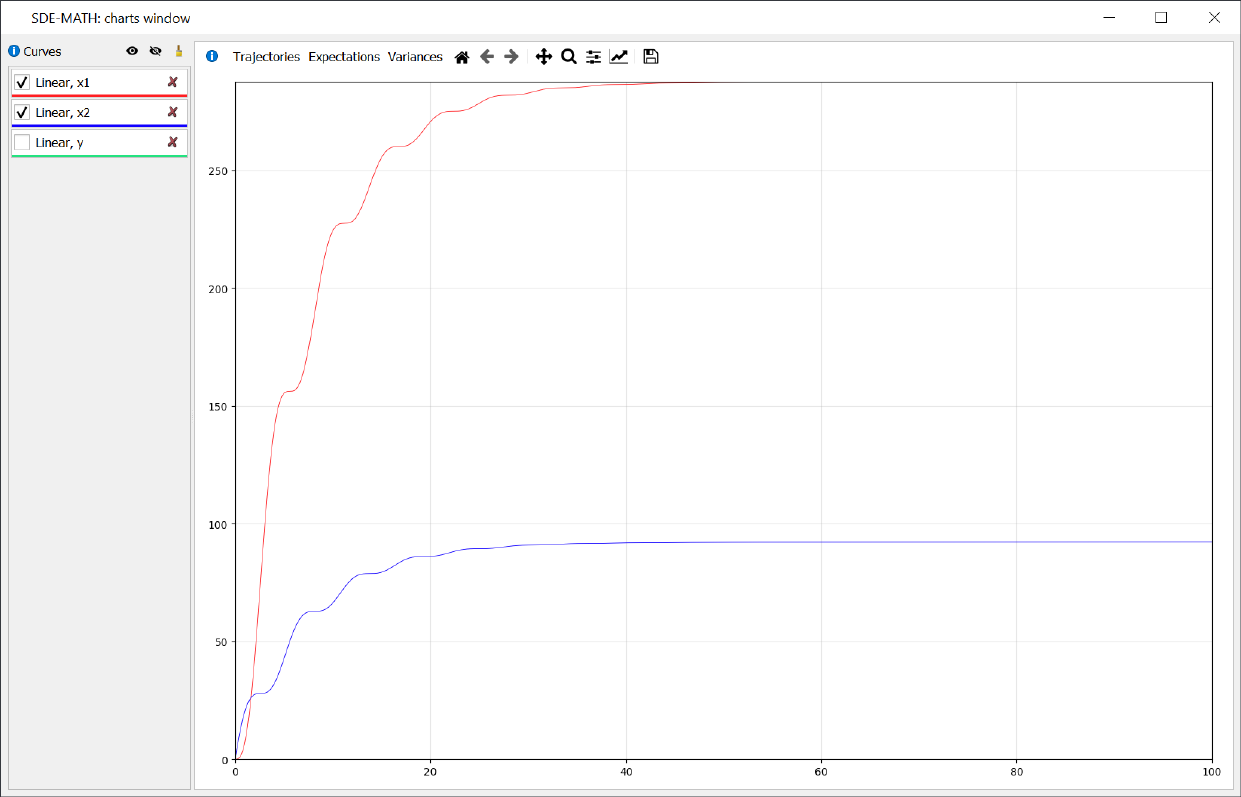}
    \caption{Solar activity model (variances)\label{fig:solar_5}}
\end{figure}

\subsection{Example of Abstract Linear System of It\^o SDEs}

Now consider the system of linear It\^{o} SDEs (\ref{lin1}) with the following data

\vspace{-3mm}
\begin{equation}
    \label{kek}
    n=4,\ \  m=5,\ \  k=3,\ \  {\bf x}_0^{(1)}=1,\ \ 
    {\bf x}_0^{(2)}=2,\ \ {\bf x}_0^{(3)}=-1,\ \  {\bf x}_0^{(4)}=-2,
\end{equation}

\vspace{-3mm}
\begin{equation}
    \label{kek1}
    A = \left(\begin{matrix}
        -1 & 0 & 0 & 0\\\\
        0 & -2 & 0 & 0\\\\
        0 & 0 & -3 & 0\\\\
        0 & 0 & 0 & -4
    \end{matrix}\right)\hspace{-1mm},\ \ 
    B = \left(\begin{matrix}
        ~1 & 1 & 1~\\\\
        ~1 & 1 & 1~\\\\
        ~1 & 1 & 1~\\\\
        ~1 & 1 & 1~
    \end{matrix}\right)\hspace{-1mm},\ \ 
    F = \left(\begin{matrix}
        0.2 & 0.1 & 0.1 & 0.1 & 0.1\\\\
        0.1 & 0.2 & 0.1 & 0.1 & 0.1\\\\
        0.1 & 0.1 & 0.2 & 0.1 & 0.1\\\\
        0.1 & 0.1 & 0.1 & 0.2 & 0.1
    \end{matrix}\right)\hspace{-1mm},
\end{equation}

\begin{equation}
    \label{kek2}
    {\bf u}(t) \equiv \left(\begin{matrix}
        0 & 0 & 0 
    \end{matrix}\right)^{{\sf T}},\ \ 
    H = \left(\begin{matrix}
        0.1 & 0.1 & 0.1 & 0.1
    \end{matrix}\right).
\end{equation}

\subsection{Visualization and Numerical Results for Abstract Linear System of It\^o SDEs Obtained via the SDE-MATH Software Package}

This subsection is devoted to the visualization and numerical results 
for the model (\ref{lin1}), (\ref{kek})--(\ref{kek2}).

\vspace{4mm}
\begin{figure}[H]
    \centering
    \includegraphics[width=.45\textwidth]{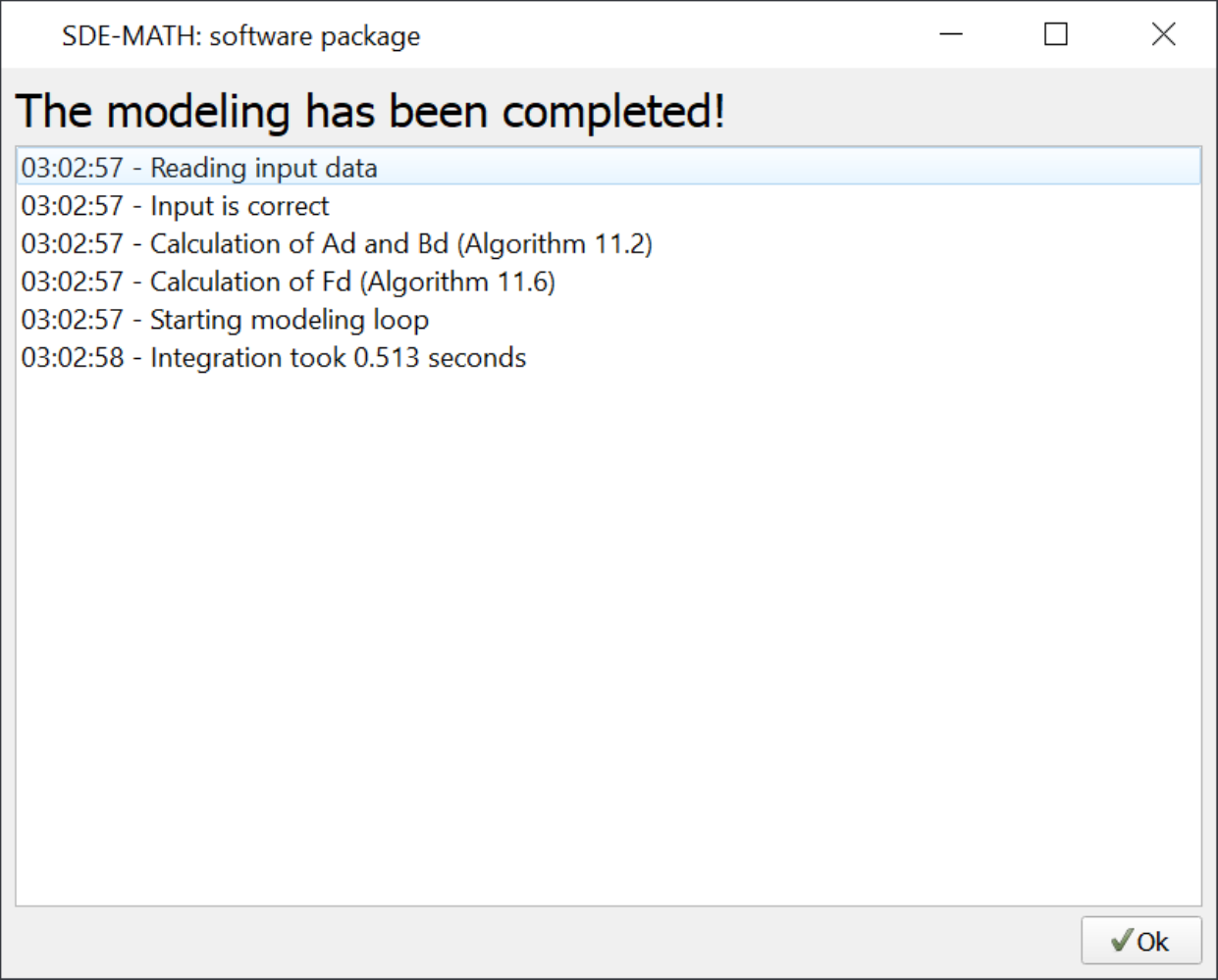}
    \caption{Modeling logs (linear system of It\^o SDEs (\ref{lin1}), (\ref{kek})--(\ref{kek2}))\label{fig:abstract_1}}
\end{figure}

\begin{figure}[H]
    \centering
    \includegraphics[width=.9\textwidth]{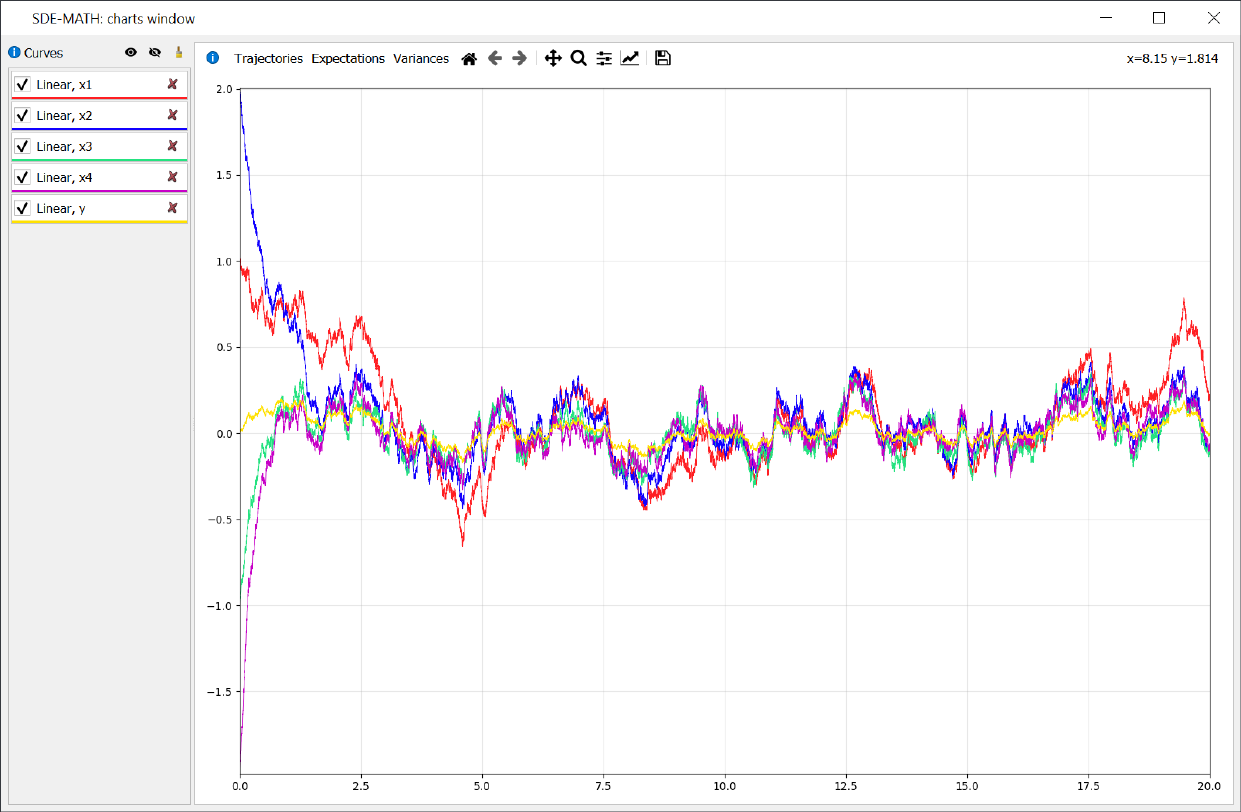}
    \caption{Linear system of It\^o SDEs (\ref{lin1}), (\ref{kek})--(\ref{kek2}) (components of solution)\label{fig:abstract_2}}
\end{figure}

\begin{figure}[H]
    \centering
    \includegraphics[width=.9\textwidth]{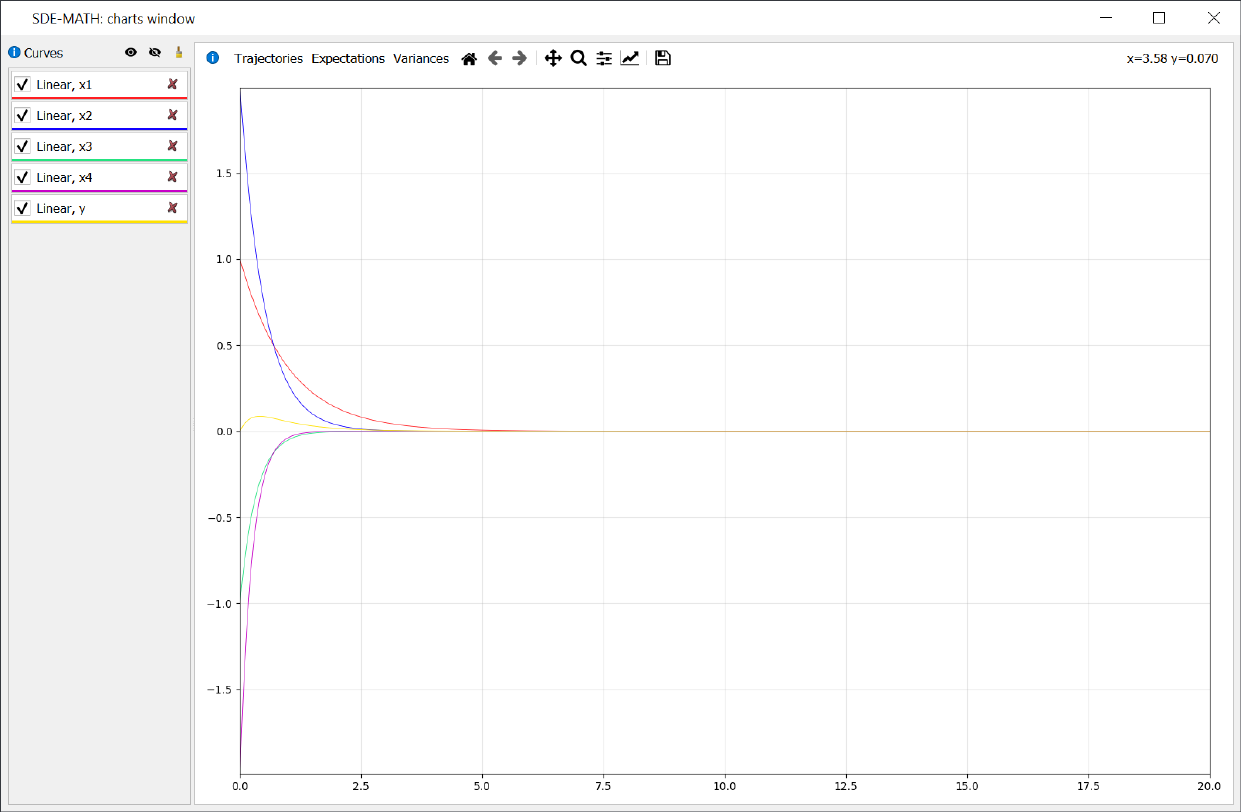}
    \caption{Linear system of It\^o SDEs (\ref{lin1}), (\ref{kek})--(\ref{kek2}) (expectations)\label{fig:abstract_3}}
\end{figure}

\begin{figure}[H]
    \centering
    \includegraphics[width=.9\textwidth]{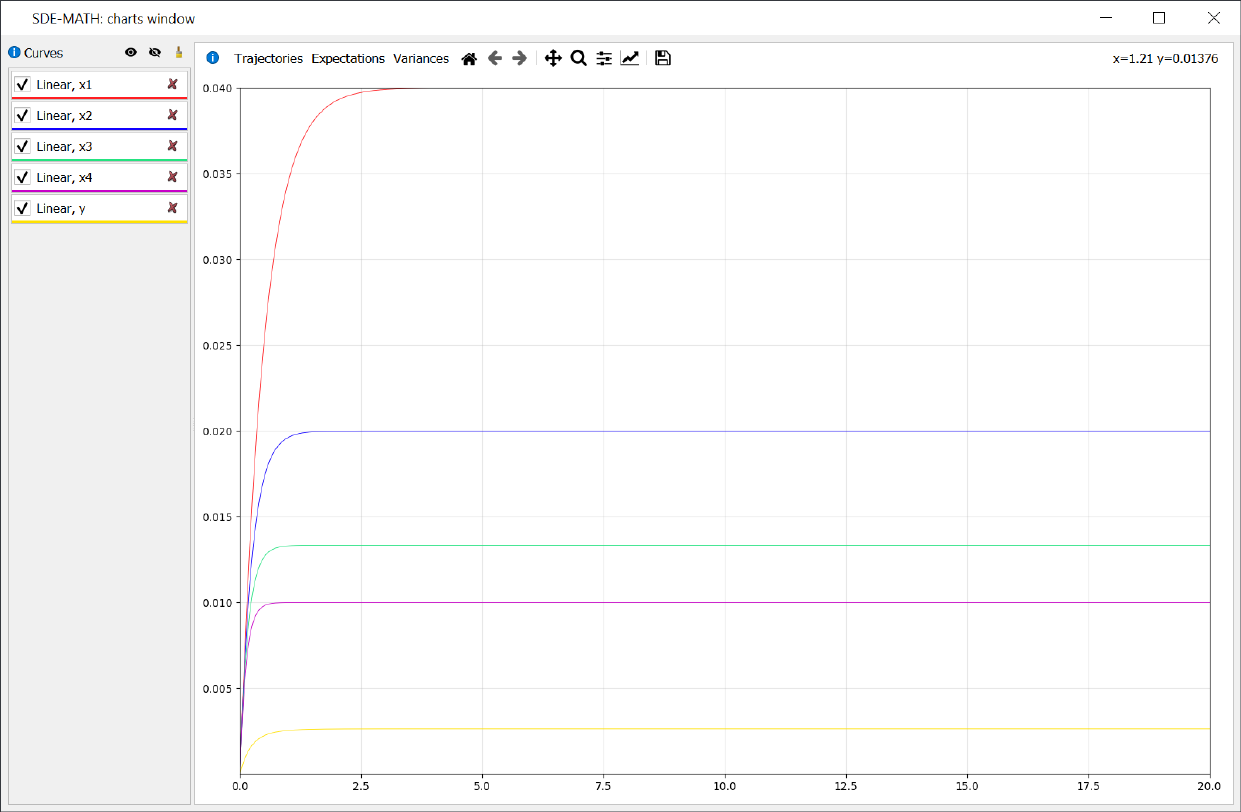}
    \caption{Linear system of It\^o SDEs (\ref{lin1}), (\ref{kek})--(\ref{kek2}) (variances)\label{fig:abstract_4}}
\end{figure}


\section{Source Codes of the SDE-MATH Software Package in the Python Programming Language}

\subsection{Source Codes of Graphical User Interface}

\subsubsection{Source Codes of Main Menu}

\begin{lstlisting}[language=Python, label={lst:gui_config}, caption=\bfseries{Configuration file example}]
"""
Configuration file example, change this paths to yours
"""

# Paths to resources
resources = "./resources/"

# Path to database
database = "./resources/database.db"

# Size of read buffer for the Fourier-Legendre coefficients
read_buffer_size = 8192

# Recursion limit for the Fourier-Legendre calculations
recursion_limit = 10 ** 8
\end{lstlisting}
\begin{lstlisting}[language=Python, label={lst:gui_main_gui}, caption=\bfseries{Program entry}]
#!/usr/bin/env python
import logging
import os
import sys

from PyQt5 import QtGui, QtWidgets
from PyQt5.QtWidgets import QApplication
from PyQt5.QtWinExtras import QWinTaskbarButton

from config import database, images
from mathematics.sde.nonlinear.symbolic.coefficients.c import C
from tools.database import connect, disconnect
from ui.main.main_window import MainWindow

def main():

    logging.basicConfig(
        level=logging.INFO,
        format="%(asctime)s - %(levelname)s - %(message)s",
        datefmt="%H:%M:%S"
    )

    app = QApplication(sys.argv)
    app.setWindowIcon(QtGui.QIcon(os.path.join(images, "function.png")))
    app.setStyle(QtWidgets.QStyleFactory.create('Fusion'))

    main_window = MainWindow()

    main_window.taskbar_button = QWinTaskbarButton()
    main_window.taskbar_button.setOverlayIcon(QtGui.QIcon("resources/function.svg"))

    exit(app.exec())

if __name__ == "__main__":
    main()
\end{lstlisting}
\begin{lstlisting}[language=Python, label={lst:gui_main_window}, caption=\bfseries{Main window}]
from PyQt5.QtCore import QThreadPool, pyqtSignal
from PyQt5.QtWidgets import QStackedWidget, QMainWindow
from sympy.physics.mechanics.tests.test_system import lam

from init.initialization import initialization
from tools.fsys import is_locked
from ui.async_calls.worker import Worker
from ui.charts.charts_window import PlotWindow
from ui.main.greetings import GreetingsWidget
from ui.main.menu.base import MainMenuWidget
from ui.main.modeling.linear.base import LinearModelingWidget
from ui.main.modeling.nonliear.base import NonlinearModelingWidget
from ui.main.progress.complex_progress import ComplexProgressWidget
from ui.main.progress.simple_progress import SimpleProgressWidget

class MainWindow(QMainWindow):
    
    main_window_close = pyqtSignal()
    start_simple_progress = pyqtSignal(str)
    stop_simple_progress = pyqtSignal(str)

    def __init__(self):
        super(QMainWindow, self).__init__()

        self.plot_window = PlotWindow()

        self.stack_widget = QStackedWidget(self)

        self.main_menu = MainMenuWidget(self.stack_widget)
        self.complex_progress = ComplexProgressWidget(self.stack_widget)
        self.simple_progress = SimpleProgressWidget(self.stack_widget)
        self.greetings = GreetingsWidget(self.stack_widget)
        self.linear_modeling = LinearModelingWidget(self.stack_widget)
        self.nonlinear_modeling = NonlinearModelingWidget(self.stack_widget)

        self.stack_widget.addWidget(self.main_menu)
        self.stack_widget.addWidget(self.nonlinear_modeling)
        self.stack_widget.addWidget(self.linear_modeling)
        self.stack_widget.addWidget(self.greetings)
        self.stack_widget.addWidget(self.simple_progress)
        self.stack_widget.addWidget(self.complex_progress)

        self.setCentralWidget(self.stack_widget)

        self.exec_init()

        self.setWindowTitle("SDE-MATH: software package")
        self.setMinimumSize(640, 480)
        self.resize(800, 600)
        self.show()

        self.main_menu.group1.show_nonlinear_dialog.connect(self.show_nonlinear)
        self.main_menu.group2.show_nonlinear_dialog.connect(self.show_nonlinear)
        self.main_menu.group3.show_linear_dialog.connect(
            lambda: self.stack_widget.setCurrentWidget(self.linear_modeling))

        self.nonlinear_modeling.show_main_menu.connect(
            lambda: self.stack_widget.setCurrentWidget(self.main_menu))
        self.nonlinear_modeling.start_progress.connect(
            lambda: self.stack_widget.setCurrentWidget(self.complex_progress))

        self.linear_modeling.show_main_menu.connect(
            lambda: self.stack_widget.setCurrentWidget(self.main_menu))
        self.linear_modeling.start_progress.connect(
            lambda: self.stack_widget.setCurrentWidget(self.complex_progress))

        self.greetings.show_main_menu.connect(
            lambda: self.stack_widget.setCurrentWidget(self.main_menu))

        self.greetings.show_main_menu.connect(
            lambda: self.stack_widget.setCurrentWidget(self.main_menu))

        self.greetings.show_main_menu.connect(
            lambda: self.stack_widget.setCurrentWidget(self.main_menu))

        self.complex_progress.back_btn.clicked.connect(
            lambda: self.stack_widget.setCurrentWidget(self.main_menu))

        # plot events

        self.main_window_close.connect(self.plot_window.close)

        self.main_menu.charts_check.clicked.connect(self.plot_window.setVisible)
        self.plot_window.charts_show.connect(
            lambda: self.main_menu.charts_check.setChecked(True))
        self.plot_window.charts_hide.connect(
            lambda: self.main_menu.charts_check.setChecked(False))

        self.nonlinear_modeling.draw_chart.connect(self.plot_window.charts_list.new_items)
        self.nonlinear_modeling.draw_chart.connect(self.plot_window.plot_widget.new_items)
        self.nonlinear_modeling.draw_chart.connect(self.plot_window.show)

        self.nonlinear_modeling.charts_check.stateChanged.connect(self.plot_window.setVisible)
        self.plot_window.charts_show.connect(
            lambda: self.nonlinear_modeling.charts_check.setChecked(True))
        self.plot_window.charts_hide.connect(
            lambda: self.nonlinear_modeling.charts_check.setChecked(False))

        self.linear_modeling.draw_chart.connect(self.plot_window.charts_list.new_items)
        self.linear_modeling.draw_chart.connect(self.plot_window.plot_widget.new_items)
        self.linear_modeling.draw_chart.connect(self.plot_window.show)

        self.linear_modeling.charts_check.clicked.connect(self.plot_window.setVisible)
        self.plot_window.charts_show.connect(
            lambda: self.linear_modeling.charts_check.setChecked(True))
        self.plot_window.charts_hide.connect(
            lambda: self.linear_modeling.charts_check.setChecked(False))

        self.linear_modeling.start_progress.connect(self.complex_progress.spin)
        self.nonlinear_modeling.start_progress.connect(self.complex_progress.spin)
        self.linear_modeling.stop_progress.connect(self.complex_progress.stop)
        self.nonlinear_modeling.stop_progress.connect(self.complex_progress.stop)

    def closeEvent(self, event):
        self.main_window_close.emit()

    def exec_init(self):
        self.simple_progress.spin("Preparing the database...")
        self.stack_widget.setCurrentWidget(self.simple_progress)

        worker = Worker(initialization)
        worker.signals.finished.connect(self.init_done)

        QThreadPool.globalInstance().start(worker)

    def init_done(self):
        if not is_locked(".welcome.lock"):
            self.stack_widget.setCurrentWidget(self.greetings)
        else:
            self.stack_widget.setCurrentWidget(self.main_menu)
        self.simple_progress.stop()

    def show_nonlinear(self, scheme_id):
        self.nonlinear_modeling.set_scheme(scheme_id)
        self.stack_widget.setCurrentWidget(self.nonlinear_modeling)
\end{lstlisting}
\begin{lstlisting}[language=Python, label={lst:gui_greetings}, caption=\bfseries{Greetings window}]
from PyQt5.QtCore import pyqtSignal, Qt
from PyQt5.QtWidgets import QPushButton, QVBoxLayout, QWidget, QSizePolicy, \
    QSpacerItem, QHBoxLayout, QLabel, QCheckBox, QApplication, QStyle

from tools.fsys import lock, unlock
from ui.main.svg import SVG

class GreetingsWidget(QWidget):
    
    show_main_menu = pyqtSignal()

    def __init__(self, parent=None):
        super(QWidget, self).__init__(parent)

        header = QLabel("Welcome to SDE-MATH Software Package for "
                        "the Numerical Solution of Systems of Ito SDEs")
        font = header.font()
        font.setPointSize(15)
        header.setAlignment(Qt.AlignJustify)
        header.setWordWrap(True)
        header.setFont(font)

        welcome = QLabel(
            "Exact solutions of Ito SDEs are known in rare cases. For this "
            "reason, it becomes necessary to construct numerical methods for "
            "Ito SDEs. Moreover, the problem of numerical solution of Ito SDEs "
            "often occurs even in cases when the exact solution of Ito SDE is known. "
            "This means that in some cases, knowing the exact solution to the Ito "
            "SDE does not allow us to simulate it numerically in a simple way.", self
        )

        font = welcome.font()
        welcome.setFont(font)

        welcome.setAlignment(Qt.AlignJustify)
        welcome.setWordWrap(True)
        welcome.setSizePolicy(QSizePolicy(QSizePolicy.Expanding, QSizePolicy.Minimum))

        check_again = QCheckBox("Do not show again", self)
        next_btn = QPushButton("Ok", self)
        next_btn.setIcon(QApplication.style().standardIcon(QStyle.SP_DialogApplyButton))

        check_again.clicked.connect(self.check_lock)
        next_btn.clicked.connect(lambda: self.show_main_menu.emit())

        controls = QHBoxLayout()
        controls.addItem(QSpacerItem(0, 0, QSizePolicy.Expanding, QSizePolicy.Minimum))
        controls.addWidget(check_again)
        controls.addWidget(next_btn)
        controls.addItem(QSpacerItem(0, 0, QSizePolicy.Expanding, QSizePolicy.Minimum))

        eq1 = QHBoxLayout()
        eq1.addItem(QSpacerItem(0, 0, QSizePolicy.Expanding, QSizePolicy.Minimum))
        eq1.addWidget(SVG("equation1.svg", scale_factor=1.))
        eq1.addItem(QSpacerItem(0, 0, QSizePolicy.Expanding, QSizePolicy.Minimum))

        eq2 = QHBoxLayout()
        eq2.addItem(QSpacerItem(0, 0, QSizePolicy.Expanding, QSizePolicy.Minimum))
        eq2.addWidget(SVG("equation2.svg", scale_factor=1.))
        eq2.addItem(QSpacerItem(0, 0, QSizePolicy.Expanding, QSizePolicy.Minimum))

        column = QVBoxLayout()
        column.addItem(QSpacerItem(0, 0, QSizePolicy.Expanding, QSizePolicy.Expanding))
        column.addWidget(header)
        column.addItem(QSpacerItem(0, 15, QSizePolicy.Expanding, QSizePolicy.Minimum))
        column.addLayout(eq1)
        column.addItem(QSpacerItem(0, 15, QSizePolicy.Expanding, QSizePolicy.Minimum))
        column.addLayout(eq2)
        column.addItem(QSpacerItem(0, 15, QSizePolicy.Expanding, QSizePolicy.Minimum))
        column.addWidget(welcome)
        column.addItem(QSpacerItem(0, 30, QSizePolicy.Expanding, QSizePolicy.Minimum))
        column.addLayout(controls)
        column.addItem(QSpacerItem(0, 0, QSizePolicy.Expanding, QSizePolicy.Expanding))

        layout = QHBoxLayout()
        layout.addItem(QSpacerItem(50, 0, QSizePolicy.Expanding, QSizePolicy.Minimum))
        layout.addLayout(column)
        layout.addItem(QSpacerItem(50, 0, QSizePolicy.Expanding, QSizePolicy.Minimum))

        self.setLayout(layout)

    def check_lock(self):
        if self.sender().isChecked():
            lock(".welcome.lock")
        else:
            unlock(".welcome.lock")
\end{lstlisting}
\begin{lstlisting}[language=Python, label={lst:gui_info}, caption=\bfseries{Info icon}]
from PyQt5.QtCore import QSize
from PyQt5.QtWidgets import QWidget, QApplication, QStyle, QLabel

class InfoIcon(QLabel):

    def __init__(self, text: str, parent=None):
        super(QWidget, self).__init__(parent)

        self.setToolTip(text)
        self.setStyleSheet("QToolTip {background: white;}")
        self.setPixmap(QApplication.style().standardIcon(
            QStyle.SP_MessageBoxInformation).pixmap(QSize(16, 16)))
\end{lstlisting}
\begin{lstlisting}[language=Python, label={lst:gui_error}, caption=\bfseries{Error widget}]
from PyQt5.QtCore import QSize
from PyQt5.QtWidgets import QWidget, QApplication, QStyle, QLabel, QHBoxLayout, QSpacerItem, QSizePolicy

class ErrorWidget(QWidget):

    def __init__(self, text: str, parent=None):
        super(QWidget, self).__init__(parent)

        msg_m = QLabel(text)
        msg_m.setStyleSheet("QLabel { color: rgb(230, 0, 0); }")

        msg_i = QLabel()
        msg_i.setStyleSheet("QToolTip { background: white; }")
        msg_i.setPixmap(QApplication.style().standardIcon(
            QStyle.SP_MessageBoxCritical).pixmap(QSize(16, 16)))

        layout = QHBoxLayout()
        layout.setContentsMargins(0, 0, 0, 0)
        layout.addWidget(msg_i)
        layout.addWidget(msg_m)
        layout.addItem(QSpacerItem(0, 0, QSizePolicy.Expanding, QSizePolicy.Minimum))

        self.setLayout(layout)
\end{lstlisting}
\begin{lstlisting}[language=Python, label={lst:gui_svg}, caption=\bfseries{Svg picture}]
import os

from PyQt5.QtCore import QSize
from PyQt5.QtSvg import QSvgWidget
from PyQt5.QtWidgets import QSizePolicy

from config import images

class SVG(QSvgWidget):
    
    def __init__(self, name: str, scale_factor=1.):
        super(QSvgWidget, self).__init__()

        self.load(os.path.join(images, name))

        self.scale_factor = scale_factor
        self.setSizePolicy(QSizePolicy.Fixed, QSizePolicy.Fixed)

    def sizeHint(self):
        size = self.renderer().defaultSize()
        return QSize(size.width() * self.scale_factor,
                     size.height() * self.scale_factor)
\end{lstlisting}
\begin{lstlisting}[language=Python, label={lst:gui_menu_base}, caption=\bfseries{Main menu (base part)}]
from PyQt5.QtWidgets import QWidget, QSizePolicy, QSpacerItem, QHBoxLayout, QVBoxLayout, QCheckBox, QLabel

from ui.main.info import InfoIcon
from ui.main.menu.linear import LinearGroupWidget
from ui.main.menu.taylor_ito import ItoGroupWidget
from ui.main.menu.taylor_stratonovich import StratonovichGroupWidget

class MainMenuWidget(QWidget):

    def __init__(self, parent=None):
        super(QWidget, self).__init__(parent)

        self.charts_check = QCheckBox("Charts window", self)

        icon = InfoIcon("This is charts window checkbox, it will\n"
                        "follow you on every application dialog, so\n"
                        "you can easily open or close window with available charts")

        bar_layout = QHBoxLayout()
        bar_layout.addItem(QSpacerItem(0, 35, QSizePolicy.Expanding, QSizePolicy.Minimum))
        bar_layout.addWidget(icon)
        bar_layout.addWidget(self.charts_check)

        header = QLabel("Strong Numerical Schemes for Ito SDEs", parent=self)
        font = header.font()
        font.setPointSize(15)
        header.setFont(font)

        icon = InfoIcon("You are now in main menu, you can choose\n"
                        "any scheme to perform modeling")

        header_layout = QHBoxLayout()
        header_layout.addItem(QSpacerItem(0, 0, QSizePolicy.Expanding, QSizePolicy.Expanding))
        header_layout.addWidget(icon)
        header_layout.addWidget(header)
        header_layout.addItem(QSpacerItem(0, 0, QSizePolicy.Expanding, QSizePolicy.Expanding))

        self.group3 = LinearGroupWidget(self)
        self.group1 = ItoGroupWidget(self)
        self.group2 = StratonovichGroupWidget(self)

        menu_layout = QHBoxLayout()
        menu_layout.addWidget(self.group1)
        menu_layout.addWidget(self.group2)
        menu_layout.addWidget(self.group3)

        layout = QVBoxLayout()
        layout.addLayout(bar_layout)
        layout.addItem(QSpacerItem(0, 0, QSizePolicy.Expanding, QSizePolicy.Expanding))
        layout.addLayout(header_layout)
        layout.addItem(QSpacerItem(0, 0, QSizePolicy.Expanding, QSizePolicy.Expanding))
        layout.addLayout(menu_layout)
        layout.addItem(QSpacerItem(0, 0, QSizePolicy.Expanding, QSizePolicy.Expanding))

        self.setLayout(layout)
\end{lstlisting}
\begin{lstlisting}[language=Python, label={lst:gui_linear}, caption=\bfseries{Main menu (linear part)}]
from PyQt5.QtCore import pyqtSignal
from PyQt5.QtWidgets import QPushButton, QVBoxLayout, QSizePolicy, QSpacerItem, QGroupBox

class LinearGroupWidget(QGroupBox):
    
    show_linear_dialog = pyqtSignal()

    def __init__(self, parent=None):
        super(QGroupBox, self).__init__(parent)

        linear_btn = QPushButton("Dispersion Spectral Decomposition")

        layout = QVBoxLayout()
        layout.addWidget(linear_btn)
        layout.addItem(QSpacerItem(0, 0, QSizePolicy.Expanding, QSizePolicy.Expanding))

        self.setLayout(layout)

        self.setTitle("Linear Ito SDEs Systems Modeling")
        self.setSizePolicy(QSizePolicy(QSizePolicy.Expanding, QSizePolicy.Expanding))

        linear_btn.clicked.connect(lambda: self.show_linear_dialog.emit())
\end{lstlisting}
\begin{lstlisting}[language=Python, label={lst:gui_taylor_ito}, caption=\bfseries{Main menu (Taylor-It\^{o} part)}]
from PyQt5.QtCore import pyqtSignal
from PyQt5.QtWidgets import QPushButton, QVBoxLayout, QSizePolicy, QSpacerItem, QGroupBox

class ItoGroupWidget(QGroupBox):
    
    show_nonlinear_dialog = pyqtSignal(int)

    def __init__(self, parent=None):
        super(QGroupBox, self).__init__(parent)

        btn1 = QPushButton("Euler")
        btn2 = QPushButton("Milstein")
        btn3 = QPushButton("Convergence Order 1.5")
        btn4 = QPushButton("Convergence Order 2.0")
        btn5 = QPushButton("Convergence Order 2.5")
        btn6 = QPushButton("Convergence Order 3.0")

        layout = QVBoxLayout()
        layout.addWidget(btn1)
        layout.addWidget(btn2)
        layout.addWidget(btn3)
        layout.addWidget(btn4)
        layout.addWidget(btn5)
        layout.addWidget(btn6)
        layout.addItem(QSpacerItem(0, 0, QSizePolicy.Expanding, QSizePolicy.Expanding))

        self.setLayout(layout)

        self.setTitle("Taylor-Ito Schemes")
        self.setSizePolicy(QSizePolicy(QSizePolicy.Expanding, QSizePolicy.Expanding))

        btn1.clicked.connect(lambda: self.show_nonlinear_dialog.emit(0))
        btn2.clicked.connect(lambda: self.show_nonlinear_dialog.emit(1))
        btn3.clicked.connect(lambda: self.show_nonlinear_dialog.emit(2))
        btn4.clicked.connect(lambda: self.show_nonlinear_dialog.emit(3))
        btn5.clicked.connect(lambda: self.show_nonlinear_dialog.emit(4))
        btn6.clicked.connect(lambda: self.show_nonlinear_dialog.emit(5))
\end{lstlisting}
\begin{lstlisting}[language=Python, label={lst:gui_taylor_stratonovich}, caption=\bfseries{Main menu (Taylor--Stratonovich part)}]
from PyQt5.QtCore import pyqtSignal
from PyQt5.QtWidgets import QPushButton, QVBoxLayout, QSizePolicy, QSpacerItem, QGroupBox

class StratonovichGroupWidget(QGroupBox):
    
    show_nonlinear_dialog = pyqtSignal(int)

    def __init__(self, parent=None):
        super(QGroupBox, self).__init__(parent)

        btn1 = QPushButton("Convergence Order 1.0")
        btn2 = QPushButton("Convergence Order 1.5")
        btn3 = QPushButton("Convergence Order 2.0")
        btn4 = QPushButton("Convergence Order 2.5")
        btn5 = QPushButton("Convergence Order 3.0")

        layout = QVBoxLayout()
        layout.addWidget(btn1)
        layout.addWidget(btn2)
        layout.addWidget(btn3)
        layout.addWidget(btn4)
        layout.addWidget(btn5)
        layout.addItem(QSpacerItem(0, 0, QSizePolicy.Expanding, QSizePolicy.Expanding))

        self.setLayout(layout)

        self.setTitle("Taylor-Stratonovich Schemes")
        self.setSizePolicy(QSizePolicy(QSizePolicy.Expanding, QSizePolicy.Expanding))

        btn1.clicked.connect(lambda: self.show_nonlinear_dialog.emit(6))
        btn2.clicked.connect(lambda: self.show_nonlinear_dialog.emit(7))
        btn3.clicked.connect(lambda: self.show_nonlinear_dialog.emit(8))
        btn4.clicked.connect(lambda: self.show_nonlinear_dialog.emit(9))
        btn5.clicked.connect(lambda: self.show_nonlinear_dialog.emit(10))
\end{lstlisting}
\begin{lstlisting}[language=Python, label={lst:gui_complex_progress}, caption=\bfseries{Complex progress view}]
import logging

from PyQt5.QtWidgets import QWidget, QHBoxLayout, QVBoxLayout, QSpacerItem, QSizePolicy, \
    QListWidget, QLabel, QPushButton, QApplication, QStyle
from pyqtspinner.spinner import WaitingSpinner

from ui.main.progress.log_handler import LogHandler

class ComplexProgressWidget(QWidget):

    def __init__(self, parent=None):
        super(QWidget, self).__init__(parent)

        self.spinner = WaitingSpinner(self, 
                                      radius=5.0, 
                                      lines=10, 
                                      line_length=5.0, 
                                      centerOnParent=False)
        self.list_widget = QListWidget(self)
        self.handler = LogHandler(self.handle_message)

        self.label = QLabel(self)
        font = self.label.font()
        font.setPointSize(15)
        self.label.setFont(font)

        self.back_btn = QPushButton("Ok")
        self.back_btn.setIcon(QApplication.style().standardIcon(QStyle.SP_DialogApplyButton))
        self.back_btn.hide()

        spinner_layout = QHBoxLayout()
        spinner_layout.addWidget(self.spinner)
        spinner_layout.addWidget(self.label)
        spinner_layout.addSpacerItem(QSpacerItem(0, 0, QSizePolicy.Expanding, 
                                                       QSizePolicy.Minimum))

        bottom_bar = QHBoxLayout()
        bottom_bar.addItem(QSpacerItem(0, 0, QSizePolicy.Expanding, QSizePolicy.Minimum))
        bottom_bar.addWidget(self.back_btn)

        layout = QVBoxLayout()
        layout.addLayout(spinner_layout)
        layout.addWidget(self.list_widget)
        layout.addLayout(bottom_bar)

        self.setLayout(layout)

    def handle_message(self, text):
        self.list_widget.addItem(text)
        self.list_widget.scrollToBottom()

    def spin(self, text):
        self.list_widget.clear()
        logging.getLogger().addHandler(self.handler)
        self.back_btn.hide()
        self.spinner.start()
        self.label.setText(text)

    def stop(self, text):
        self.back_btn.show()
        self.list_widget.scrollToBottom()
        self.label.setText(text)
        self.spinner.stop()
        logging.getLogger().removeHandler(self.handler)
\end{lstlisting}
\begin{lstlisting}[language=Python, label={lst:gui_simple_progress}, caption=\bfseries{Simple progress view}]
from PyQt5.QtCore import Qt
from PyQt5.QtWidgets import QVBoxLayout, QWidget, QLabel, QSpacerItem, QSizePolicy
from pyqtspinner.spinner import WaitingSpinner

class SimpleProgressWidget(QWidget):

    def __init__(self, parent=None):
        super(QWidget, self).__init__(parent)

        self.spinner = WaitingSpinner(self, radius=15.0, lines=10, line_length=15.0)

        self.label = QLabel()
        self.label.setAlignment(Qt.AlignCenter)
        font = self.label.font()
        font.setPointSize(15)
        self.label.setFont(font)

        layout = QVBoxLayout()
        layout.addItem(QSpacerItem(0, 0, QSizePolicy.Expanding, QSizePolicy.Expanding))
        layout.addWidget(self.spinner)
        layout.addItem(QSpacerItem(0, 50, QSizePolicy.Minimum, QSizePolicy.Minimum))
        layout.addWidget(self.label)
        layout.addItem(QSpacerItem(0, 0, QSizePolicy.Expanding, QSizePolicy.Expanding))

        self.setLayout(layout)

    def spin(self, text):
        self.spinner.start()
        self.label.setText(text)

    def stop(self):
        self.spinner.stop()
\end{lstlisting}
\begin{lstlisting}[language=Python, label={lst:gui_log_handler}, caption=\bfseries{Log handler for application}]
import logging

class LogHandler(logging.Handler):
    
    def __init__(self, callback):
        super().__init__()
        self.callback = callback
        self.setFormatter(logging.Formatter("%(asctime)s - %(message)s", "%H:%M:%S"))

    def handle(self, record):
        self.callback(self.format(record))
\end{lstlisting}
\begin{lstlisting}[language=Python, label={lst:gui_matrix_widget}, caption=\bfseries{Matrix widget}]
from PyQt5.QtWidgets import QTableWidgetItem, QTableWidget, QSizePolicy

class MatrixWidget(QTableWidget):

    def __init__(self, parent=None):
        super(QTableWidget, self).__init__(parent)

        self.itemChanged.connect(self.item_changed)
        self.m = [["0"]]

        self.setRowCount(1)
        self.setColumnCount(1)
        self.setItem(0, 0, CustomItem("0"))

        self.setSizePolicy(QSizePolicy(QSizePolicy.Expanding, QSizePolicy.Expanding))

    def resize_w(self, w: int):

        self.blockSignals(True)

        old_w = self.columnCount()
        self.setColumnCount(w)
        h = self.rowCount()

        self.m = [[self.m[i][j] if i < h and j < old_w else "0"
                   for j in range(w)]
                  for i in range(h)]

        for i in range(h):
            for j in range(w):
                item = self.item(i, j)
                if item is not None:
                    item.setText(self.m[i][j])
                else:
                    self.setItem(i, j, CustomItem(self.m[i][j]))

        self.blockSignals(False)

    def resize_h(self, h: int):

        self.blockSignals(True)

        old_h = self.rowCount()
        w = self.columnCount()
        self.setRowCount(h)

        self.m = [[self.m[i][j] if i < old_h and j < w else "0"
                   for j in range(w)]
                  for i in range(h)]

        for i in range(h):
            for j in range(w):
                item = self.item(i, j)
                if item is not None:
                    item.setText(self.m[i][j])
                else:
                    self.setItem(i, j, CustomItem(self.m[i][j]))

        self.blockSignals(False)

    def item_changed(self, item):
        self.m[item.row()][item.column()] = item.text()

class CustomItem(QTableWidgetItem):

    def __init__(self, value: str):
        super(QTableWidgetItem, self).__init__(value)

        self.valid = True
\end{lstlisting}

\subsubsection{Source Codes of Charts Window}

\begin{lstlisting}[language=Python, label={lst:gui_charts_window}, caption=\bfseries{Charts window}]
from PyQt5.QtCore import pyqtSignal
from PyQt5.QtWidgets import QWidget, QHBoxLayout, QMainWindow, QSplitter

from ui.charts.side.available_charts_widget import AvailableChartsWidget
from ui.charts.visuals.charts_widget import ChartsWidget

class PlotWindow(QMainWindow):
    
    charts_show = pyqtSignal()
    charts_hide = pyqtSignal()

    def __init__(self):
        super(QMainWindow, self).__init__()

        self.plot_widget = ChartsWidget(self)
        self.charts_list = AvailableChartsWidget(self)

        splitter = QSplitter()
        splitter.addWidget(self.charts_list)
        splitter.addWidget(self.plot_widget)
        splitter.setSizes([splitter.width() / 0.85,
                           splitter.width() / 0.15])

        layout = QHBoxLayout(self)
        layout.addWidget(splitter)

        central_widget = QWidget()
        central_widget.setLayout(layout)
        self.setCentralWidget(central_widget)

        self.setWindowTitle("SDE-MATH: charts window")
        self.resize(1200, 800)

        self.charts_list.on_show_all.connect(self.plot_widget.show_all)
        self.charts_list.on_hide_all.connect(self.plot_widget.hide_all)
        self.charts_list.on_remove_all.connect(self.plot_widget.delete_all)

    def showEvent(self, event):
        self.charts_show.emit()

    def closeEvent(self, event):
        self.charts_hide.emit()
\end{lstlisting}
\begin{lstlisting}[language=Python, label={lst:gui_available_charts_widget}, caption=\bfseries{Curves list}]
import os

from PyQt5 import QtGui
from PyQt5.QtCore import pyqtSignal
from PyQt5.QtWidgets import QVBoxLayout, QWidget, QSizePolicy, QSpacerItem, \
     QScrollArea, QLabel, QHBoxLayout, QPushButton, QApplication, QStyle

from config import images
from ui.charts.side.item_widget import ItemWidget
from ui.main.info import InfoIcon

class AvailableChartsWidget(QWidget):
    
    on_hide_all = pyqtSignal()
    on_show_all = pyqtSignal()
    on_remove_all = pyqtSignal()

    def __init__(self, parent=None):
        super(QWidget, self).__init__(parent)

        self.items = dict()

        self.spacer = QSpacerItem(0, 0, QSizePolicy.Minimum, QSizePolicy.Expanding)
        self.plot_widget = self.parent().plot_widget

        remove_all = QPushButton()
        remove_all.setFlat(True)
        remove_all.setIcon(
            QApplication.style().standardIcon(QStyle.SP_DialogResetButton))

        hide_all = QPushButton()
        hide_all.setFlat(True)
        hide_all.setIcon(QtGui.QIcon(os.path.join(images, "crossed.png")))

        show_all = QPushButton()
        show_all.setFlat(True)
        show_all.setIcon(QtGui.QIcon(os.path.join(images, "eye.png")))

        header_layout = QHBoxLayout()
        header_layout.addWidget(
            InfoIcon("Here You will see all modeling series\n"
                     "You can hide them or delete, if you need to"))
        header_layout.addItem(
            QSpacerItem(5, 0, QSizePolicy.Minimum, QSizePolicy.Minimum))
        header_layout.addWidget(QLabel("Curves"))
        header_layout.addItem(
            QSpacerItem(5, 0, QSizePolicy.Minimum, QSizePolicy.Minimum))
        header_layout.setContentsMargins(0, 0, 0, 0)
        header_layout.setSpacing(0)
        header_layout.addItem(
            QSpacerItem(0, 0, QSizePolicy.Expanding, QSizePolicy.Minimum))
        header_layout.addWidget(show_all)
        header_layout.addWidget(hide_all)
        header_layout.addWidget(remove_all)

        self.layout = QVBoxLayout()
        self.layout.setContentsMargins(3, 3, 3, 3)
        self.layout.setSpacing(2)

        scroll_widget = QWidget(self)
        scroll_widget.setLayout(self.layout)

        scroll_area = QScrollArea(self)
        scroll_area.setWidgetResizable(True)
        scroll_area.setWidget(scroll_widget)

        self.layout.addItem(self.spacer)

        layout = QVBoxLayout()
        layout.setContentsMargins(0, 0, 0, 0)
        layout.addLayout(header_layout)
        layout.addWidget(scroll_area)

        self.setLayout(layout)

        self.setSizePolicy(
            QSizePolicy(QSizePolicy.MinimumExpanding, 
                        QSizePolicy.MinimumExpanding))

        show_all.clicked.connect(self.show_all)
        hide_all.clicked.connect(self.hide_all)
        remove_all.clicked.connect(self.delete_all)

    def new_items(self, lines: list):
        self.layout.removeItem(self.spacer)

        for i in range(len(lines)):
            item_widget = ItemWidget(lines[i].name, lines[i].color, parent=self)
            item_widget.uid = lines[i].uid
            item_widget.on_show.connect(self.plot_widget.show_item)
            item_widget.on_hide.connect(self.plot_widget.hide_item)
            item_widget.on_delete.connect(self.plot_widget.delete_item)
            item_widget.on_delete.connect(self.delete_item)
            self.items[lines[i].uid] = item_widget

            self.plot_widget.hide_label.connect(lambda uid: self.items[uid].hide())
            self.plot_widget.show_label.connect(lambda uid: self.items[uid].show())

            self.layout.addWidget(item_widget)

        self.layout.addItem(self.spacer)

    def delete_item(self):
        s = self.sender()
        s.setParent(None)
        self.items.pop(s.uid)
        self.layout.removeWidget(s)

    def delete_all(self):
        for item in self.items.values():
            item.setParent(None)
            self.layout.removeWidget(item)
        self.items.clear()
        self.on_remove_all.emit()

    def hide_all(self):
        for item in self.items.values():
            item.checkbox.blockSignals(True)
            item.checkbox.setChecked(False)
            item.checkbox.blockSignals(False)
        self.on_hide_all.emit()

    def show_all(self):
        for item in self.items.values():
            item.checkbox.blockSignals(True)
            item.checkbox.setChecked(True)
            item.checkbox.blockSignals(False)
        self.on_show_all.emit()
\end{lstlisting}
\begin{lstlisting}[language=Python, label={lst:gui_item_widget}, caption=\bfseries{Curves list item}]
from PyQt5.QtCore import pyqtSignal
from PyQt5.QtWidgets import QPushButton, QSizePolicy, QHBoxLayout, QStyle, \
    QApplication, QCheckBox, QVBoxLayout, QLabel, QSpacerItem, QFrame

class ItemWidget(QFrame):
    
    on_delete = pyqtSignal(object)
    on_hide = pyqtSignal(int)
    on_show = pyqtSignal(int)

    def __init__(self, name, color, parent=None):
        super(QFrame, self).__init__(parent)

        self.uid = 0

        self.setFrameShape(QFrame.StyledPanel)
        self.setStyleSheet("QFrame { background: white; }")

        self.checkbox = QCheckBox(name)
        self.checkbox.setChecked(True)

        btn = QPushButton()
        btn.setFlat(True)
        btn.setIcon(QApplication.style().standardIcon(QStyle.SP_DialogCancelButton))

        underline = QLabel()
        underline.setSizePolicy(QSizePolicy(QSizePolicy.Expanding, QSizePolicy.Maximum))
        underline.setMaximumHeight(3)
        underline.setStyleSheet(f"QLabel {{ background: {color}; }}")

        layout = QHBoxLayout()
        layout.setContentsMargins(3, 3, 3, 3)
        layout.setSpacing(0)
        layout.addWidget(self.checkbox)
        layout.addItem(QSpacerItem(0, 0, QSizePolicy.Expanding, QSizePolicy.Minimum))
        layout.addWidget(btn)

        layout_underlined = QVBoxLayout()
        layout_underlined.setContentsMargins(0, 0, 0, 0)
        layout_underlined.setSpacing(0)
        layout_underlined.addLayout(layout)
        layout_underlined.addWidget(underline)

        self.setLayout(layout_underlined)

        self.checkbox.stateChanged.connect(self.checkbox_changed)
        btn.clicked.connect(lambda: self.on_delete.emit(self.uid))

    def checkbox_changed(self, v):
        if v > 0:
            self.on_show.emit(self.uid)
        else:
            self.on_hide.emit(self.uid)
\end{lstlisting}
\begin{lstlisting}[language=Python, label={lst:gui_charts_widget}, caption=\bfseries{Charts area}]
import matplotlib.pyplot as plt
from PyQt5.QtCore import pyqtSignal
from PyQt5.QtWidgets import QSizePolicy, QFrame, QHBoxLayout, QPushButton, \
    QSpacerItem
from PyQt5.QtWidgets import QVBoxLayout
from matplotlib.backends.backend_qt5agg import FigureCanvasQTAgg as FigureCanvas

from ui.charts.visuals.color import Color
from ui.charts.visuals.toolbar import ToolBar
from ui.main.info import InfoIcon

class ChartsWidget(QFrame):
    
    hide_label = pyqtSignal(int)
    show_label = pyqtSignal(int)

    def __init__(self, parent=None):
        super(QFrame, self).__init__(parent)

        self.plots = dict()
        self.mode = 0

        self.setFrameStyle(QFrame.StyledPanel)
        self.setStyleSheet("QFrame { background: white; }")

        self.figure = plt.figure()
        self.canvas = FigureCanvas(self.figure)
        self.canvas.mpl_connect('resize_event', self.on_resize)
        self.toolbar = ToolBar(self.canvas, self)

        self.ax = self.figure.add_subplot(111)
        self.ax.margins(0)
        self.ax.grid(axis='both', alpha=.3)
        self.ax.relim(visible_only=True)
        self.ax.autoscale()

        self.figure.tight_layout()
        self.canvas.draw()

        self.btn_to_fn = QPushButton("Trajectories")
        self.btn_to_fn.setFlat(True)
        self.btn_to_mx = QPushButton("Expectations")
        self.btn_to_mx.setFlat(True)
        self.btn_to_dx = QPushButton("Variances")
        self.btn_to_dx.setFlat(True)

        toolbar_layout = QHBoxLayout()
        toolbar_layout.setContentsMargins(0, 0, 0, 0)
        toolbar_layout.setSpacing(0)
        toolbar_layout.addItem(QSpacerItem(15, 0, QSizePolicy.Minimum, QSizePolicy.Minimum))
        toolbar_layout.addWidget(InfoIcon("Click this buttons to switch plot modes\n"
                                          "between trajectories, expectations and variances"))
        toolbar_layout.addItem(QSpacerItem(15, 0, QSizePolicy.Minimum, QSizePolicy.Minimum))
        toolbar_layout.addWidget(self.btn_to_fn)
        toolbar_layout.addWidget(self.btn_to_mx)
        toolbar_layout.addWidget(self.btn_to_dx)
        toolbar_layout.addWidget(self.toolbar)
        toolbar_layout.addItem(QSpacerItem(15, 0, QSizePolicy.Minimum, QSizePolicy.Minimum))

        layout = QVBoxLayout()
        layout.setContentsMargins(0, 0, 0, 0)
        layout.setSpacing(0)
        layout.addLayout(toolbar_layout)
        layout.addWidget(self.canvas)

        self.setLayout(layout)

        self.setSizePolicy(QSizePolicy(QSizePolicy.Expanding, QSizePolicy.Expanding))

        self.btn_to_fn.pressed.connect(self.fn_mode)
        self.btn_to_mx.pressed.connect(self.mx_mode)
        self.btn_to_dx.pressed.connect(self.dx_mode)

    def rescale(self):
        self.ax.relim(visible_only=True)
        self.ax.autoscale()
        self.canvas.draw()

    def clear(self):
        for f in self.plots.values():
            self.hide_label.emit(f.uid)
            if f.line_fn is not None:
                f.line_fn.set_visible(False)
            if f.line_mx is not None:
                f.line_mx.set_visible(False)
            if f.line_dx is not None:
                f.line_dx.set_visible(False)

    def fn_mode(self):

        self.mode = 0

        self.clear()

        for f in self.plots.values():
            if f.line_fn is not None:
                self.show_label.emit(f.uid)
                if f.visible:
                    f.line_fn.set_visible(True)

        self.rescale()

    def mx_mode(self):

        self.mode = 1

        self.clear()

        for f in self.plots.values():
            if f.line_mx is not None:
                self.show_label.emit(f.uid)
                if f.visible:
                    f.line_mx.set_visible(True)

        self.rescale()

    def dx_mode(self):

        self.mode = 2

        self.clear()

        for f in self.plots.values():
            if f.line_dx is not None:
                self.show_label.emit(f.uid)
                if f.visible:
                    f.line_dx.set_visible(True)

        self.rescale()

    def new_items(self, lines: list):

        for line in lines:
            self.plots[line.uid] = line

        for f in self.plots.values():
            if f.line_fn is None and f.fn is not None:
                f.line_fn = self.ax.plot(f.t, f.fn, linewidth=1, color=f.color)[0]
            if f.line_mx is None and f.mx is not None:
                f.line_mx = self.ax.plot(f.t, f.mx, linewidth=1, color=f.color)[0]
            if f.line_dx is None and f.dx is not None:
                f.line_dx = self.ax.plot(f.t, f.dx, linewidth=1, color=f.color)[0]

        if self.mode == 0:
            self.fn_mode()

        if self.mode == 1:
            self.mx_mode()

        if self.mode == 2:
            self.dx_mode()

    def delete_item(self, uid: int):
        item = self.plots.pop(uid)
        Color.free(item.color)
        if item.line_fn is not None:
            item.line_fn.remove()
        if item.line_mx is not None:
            item.line_mx.remove()
        if item.line_dx is not None:
            item.line_dx.remove()

        self.rescale()

    def hide_item(self, uid: int):
        item = self.plots[uid]
        if item.line_fn is not None and self.mode == 0:
            item.line_fn.set_visible(False)
        if item.line_mx is not None and self.mode == 1:
            item.line_mx.set_visible(False)
        if item.line_dx is not None and self.mode == 2:
            item.line_dx.set_visible(False)
        item.visible = False

        self.rescale()

    def show_item(self, uid: int):
        item = self.plots[uid]
        if item.line_fn is not None and self.mode == 0:
            item.line_fn.set_visible(True)
        if item.line_mx is not None and self.mode == 1:
            item.line_mx.set_visible(True)
        if item.line_dx is not None and self.mode == 2:
            item.line_dx.set_visible(True)
        item.visible = True

        self.rescale()

    def show_all(self):
        for item in self.plots.values():
            if item.line_fn is not None and self.mode == 0:
                item.line_fn.set_visible(True)
            if item.line_mx is not None and self.mode == 1:
                item.line_mx.set_visible(True)
            if item.line_dx is not None and self.mode == 2:
                item.line_dx.set_visible(True)
            item.visible = True

        self.rescale()

    def hide_all(self):
        for item in self.plots.values():
            if item.line_fn is not None and self.mode == 0:
                item.line_fn.set_visible(False)
            if item.line_mx is not None and self.mode == 1:
                item.line_mx.set_visible(False)
            if item.line_dx is not None and self.mode == 2:
                item.line_dx.set_visible(False)
            item.visible = False

        self.rescale()

    def delete_all(self):
        for item in reversed(self.plots.values()):
            Color.free(item.color)
            if item.line_fn is not None:
                item.line_fn.remove()
            if item.line_mx is not None:
                item.line_mx.remove()
            if item.line_dx is not None:
                item.line_dx.remove()

        self.plots.clear()

        self.rescale()

    def on_resize(self, event):
        self.figure.tight_layout()
        self.canvas.draw()
\end{lstlisting}
\begin{lstlisting}[language=Python, label={lst:gui_color}, caption=\bfseries{Curves color}]
from random import choice

class Color:

    reserved_colors = [
        "#ff834a",
        "#ffe100",
        "#c700c7",
        "#24e280",
        "#1100ff",
        "#ff1e22",
    ]
    available_colors = [
        "#ff834a",
        "#ffe100",
        "#c700c7",
        "#24e280",
        "#1100ff",
        "#ff1e22",
    ]

    def __new__(cls, *args, **kwargs):

        try:
            return cls.available_colors.pop()

        except IndexError:
            return f"#{''.join([choice('0123456789ABCDEF') for j in range(6)])}"

    @classmethod
    def free(cls, code: str):

        if code in cls.reserved_colors:
            cls.available_colors.append(code)
\end{lstlisting}
\begin{lstlisting}[language=Python, label={lst:gui_line}, caption=\bfseries{Curve}]
from ui.charts.visuals.color import Color

class Line:
    count = 0

    def __init__(self, name, t, fn, mx=None, dx=None):
        self.name = name
        self.t = t
        self.fn = fn
        self.mx = mx
        self.dx = dx
        self.line_fn = None
        self.line_mx = None
        self.line_dx = None
        self.visible = True
        self.color = Color()

        self.uid = Line.count
        Line.count += 1
\end{lstlisting}

\subsubsection{Source Codes of Input for Nonlinear Systems of It\^o SDEs}

\lstinputlisting[language=Python, label={lst:gui_nonlinear_base}, caption=\bfseries{Base part of data input for nonlinear systems}]{listings/ui/main/modeling/nonliear/base.py}
\begin{lstlisting}[language=Python, label={lst:gui_nonlinear_step1}, caption=\bfseries{Step 1 of data input for nonlinear systems}]
from PyQt5.QtCore import pyqtSignal
from PyQt5.QtWidgets import QWidget, QHBoxLayout, QGridLayout, QLineEdit, QLabel, \
    QVBoxLayout, QSpacerItem, QSizePolicy, QPushButton, QApplication, QStyle

from ui.main.error import ErrorWidget
from ui.main.info import InfoIcon
from ui.main.svg import SVG

class Step1(QWidget):
    
    n_valid = pyqtSignal(int)
    m_valid = pyqtSignal(int)

    def __init__(self, parent=None):
        super(QWidget, self).__init__(parent)

        self.n_is_valid = False
        self.m_is_valid = False

        self.input_stack = self.parent()

        info_n = InfoIcon("Dimension of linear system of Ito SDEs")
        info_m = InfoIcon("Dimension of vector Wiener process")

        label_n = QLabel("n")
        label_m = QLabel("m")

        self.lineedit_n = QLineEdit()
        self.lineedit_m = QLineEdit()

        self.msg_n = ErrorWidget("Wrong value!")
        self.msg_n.hide()

        self.msg_m = ErrorWidget("Wrong value!")
        self.msg_m.hide()

        grid_layout = QGridLayout()
        grid_layout.addWidget(self.msg_n, 0, 2)
        grid_layout.addWidget(self.msg_m, 2, 2)
        grid_layout.addWidget(info_n, 1, 0)
        grid_layout.addWidget(info_m, 3, 0)
        grid_layout.addWidget(label_n, 1, 1)
        grid_layout.addWidget(label_m, 3, 1)
        grid_layout.addWidget(self.lineedit_n, 1, 2)
        grid_layout.addWidget(self.lineedit_m, 3, 2)

        header = QLabel("Dimensions settings", parent=self)
        font = header.font()
        font.setPointSize(15)
        header.setFont(font)

        self.next_btn = QPushButton("Next", self)
        self.next_btn.setIcon(QApplication.style().standardIcon(QStyle.SP_ArrowForward))
        self.next_btn.setEnabled(False)

        header_layout = QHBoxLayout()
        header_layout.addItem(QSpacerItem(0, 0, QSizePolicy.Expanding, QSizePolicy.Minimum))
        header_layout.addWidget(header)
        header_layout.addItem(QSpacerItem(0, 0, QSizePolicy.Expanding, QSizePolicy.Minimum))

        bottom_bar = QHBoxLayout()
        bottom_bar.addItem(QSpacerItem(0, 0, QSizePolicy.Expanding, QSizePolicy.Minimum))
        bottom_bar.addWidget(self.next_btn)

        eq1 = QHBoxLayout()
        eq1.addItem(QSpacerItem(0, 0, QSizePolicy.Expanding, QSizePolicy.Minimum))
        eq1.addWidget(SVG("equation1.svg", scale_factor=1.))
        eq1.addItem(QSpacerItem(0, 0, QSizePolicy.Expanding, QSizePolicy.Minimum))

        eq2 = QHBoxLayout()
        eq2.addItem(QSpacerItem(0, 0, QSizePolicy.Expanding, QSizePolicy.Minimum))
        eq2.addWidget(SVG("equation2.svg", scale_factor=1.))
        eq2.addItem(QSpacerItem(0, 0, QSizePolicy.Expanding, QSizePolicy.Minimum))

        equalities_layout = QVBoxLayout()
        equalities_layout.addLayout(eq1)
        equalities_layout.addItem(QSpacerItem(0, 20, QSizePolicy.Minimum, QSizePolicy.Minimum))
        equalities_layout.addLayout(eq2)

        equalities_wrap = QHBoxLayout()
        equalities_wrap.addItem(QSpacerItem(0, 0, QSizePolicy.Expanding, QSizePolicy.Minimum))
        equalities_wrap.addLayout(equalities_layout)
        equalities_wrap.addItem(QSpacerItem(0, 0, QSizePolicy.Expanding, QSizePolicy.Minimum))

        grid_wrap = QHBoxLayout()
        grid_wrap.addItem(QSpacerItem(0, 0, QSizePolicy.Expanding, QSizePolicy.Minimum))
        grid_wrap.addLayout(grid_layout)
        grid_wrap.addItem(QSpacerItem(0, 0, QSizePolicy.Expanding, QSizePolicy.Minimum))

        layout = QVBoxLayout()
        layout.addItem(QSpacerItem(0, 0, QSizePolicy.Minimum, QSizePolicy.Expanding))
        layout.addLayout(equalities_wrap)
        layout.addItem(QSpacerItem(0, 25, QSizePolicy.Minimum, QSizePolicy.Minimum))
        layout.addLayout(header_layout)
        layout.addItem(QSpacerItem(0, 5, QSizePolicy.Minimum, QSizePolicy.Minimum))
        layout.addLayout(grid_wrap)
        layout.addItem(QSpacerItem(0, 0, QSizePolicy.Minimum, QSizePolicy.Expanding))
        layout.addItem(QSpacerItem(0, 25, QSizePolicy.Minimum, QSizePolicy.Minimum))
        layout.addLayout(bottom_bar)

        self.setLayout(layout)

        self.lineedit_n.textChanged.connect(self.validate_n)
        self.lineedit_m.textChanged.connect(self.validate_m)

    def validate_form(self):
        if self.n_is_valid and self.m_is_valid:
            self.next_btn.setEnabled(True)
        else:
            self.next_btn.setEnabled(False)

    def validate_n(self, value):
        try:
            typed_value = int(value)
            if typed_value <= 0:
                raise ValueError()

            self.n_is_valid = True
            self.n_valid.emit(typed_value)
            self.msg_n.hide()

        except ValueError:
            self.n_is_valid = False
            self.msg_n.show()

        finally:
            self.validate_form()

    def validate_m(self, value):
        try:
            typed_value = int(value)
            if typed_value <= 0:
                raise ValueError()

            self.m_is_valid = True
            self.m_valid.emit(typed_value)
            self.msg_m.hide()

        except ValueError:
            self.m_is_valid = False
            self.msg_m.show()

        finally:
            self.validate_form()
\end{lstlisting}
\lstinputlisting[language=Python, label={lst:gui_nonlinear_step2}, caption=\bfseries{Step 2 of data input for nonlinear systems}]{listings/ui/main/modeling/nonliear/step2.py}
\lstinputlisting[language=Python, label={lst:gui_nonlinear_step3}, caption=\bfseries{Step 3 of data input for nonlinear systems}]{listings/ui/main/modeling/nonliear/step3.py}
\lstinputlisting[language=Python, label={lst:gui_nonlinear_step4}, caption=\bfseries{Step 4 of data input for nonlinear systems}]{listings/ui/main/modeling/nonliear/step4.py}
\lstinputlisting[language=Python, label={lst:gui_nonlinear_step5}, caption=\bfseries{Step 5 of data input for nonlinear systems}]{listings/ui/main/modeling/nonliear/step5.py}

\subsubsection{Source Codes of Input for Linear Systems of It\^o SDEs}

\lstinputlisting[language=Python, label={lst:gui_linear_base}, caption=\bfseries{Base part of data input for linear systems}]{listings/ui/main/modeling/linear/base.py}
\lstinputlisting[language=Python, label={lst:gui_linear_step1}, caption=\bfseries{Step 1 of data input for linear systems}]{listings/ui/main/modeling/linear/step1.py}
\lstinputlisting[language=Python, label={lst:gui_linear_step2}, caption=\bfseries{Step 2 of data input for linear systems}]{listings/ui/main/modeling/linear/step2.py}
\lstinputlisting[language=Python, label={lst:gui_linear_step3}, caption=\bfseries{Step 3 of data input for linear systems}]{listings/ui/main/modeling/linear/step3.py}
\lstinputlisting[language=Python, label={lst:gui_linear_step4}, caption=\bfseries{Step 4 of data input for linear systems}]{listings/ui/main/modeling/linear/step4.py}
\lstinputlisting[language=Python, label={lst:gui_linear_step5}, caption=\bfseries{Step 5 of data input for linear systems}]{listings/ui/main/modeling/linear/step5.py}
\lstinputlisting[language=Python, label={lst:gui_linear_step6}, caption=\bfseries{Step 6 of data input for linear systems}]{listings/ui/main/modeling/linear/step6.py}
\lstinputlisting[language=Python, label={lst:gui_linear_step7}, caption=\bfseries{Step 7 of data input for linear systems}]{listings/ui/main/modeling/linear/step7.py}
\lstinputlisting[language=Python, label={lst:gui_linear_step8}, caption=\bfseries{Step 8 of data input for linear systems}]{listings/ui/main/modeling/linear/step8.py}

\subsection{Source Codes for Nonlinear Systems of It\^o SDEs}

\subsubsection{Source Codes for Calculation of the Fourier--Legendre Coefficients}

\begin{lstlisting}[language=Python, label={lst:polinomial}, caption=\bfseries{Symbolic function of the Legendre polinomial}]
from sympy import Rational, factorial, diff

def polynomial(n: int):
    """
    Returns the Legendre polynomial in symbolic format
    Parameters
    ==========
    n : int
        degree of the Legendre polynomial
    Returns
    =======
    sympy.Expr
    """
    from sympy.abc import x
    return Rational(1, 2) ** n / factorial(n) * diff((x ** 2 - 1) ** n, x, n)
\end{lstlisting}
\lstinputlisting[language=Python, label={lst:get_c}, caption=\bfseries{Symbolic function of the Fourier--Legendre coefficient calculation}]{listings/nonlinear/c.py}
\vspace{-2mm}
\lstinputlisting[language=Python, label={lst:c}, caption=\bfseries{Symbolic function of the Fourier--Legendre coefficient}]{listings/nonlinear/coefficients/c.py}
\lstinputlisting[language=Python, label={lst:c000}, caption=\bfseries{Calculation of the Fourier--Legendre coefficients $C^{000}_{j_3 j_2 j_1}$}]{listings/nonlinear/coefficients/c000.py}
\lstinputlisting[language=Python, label={lst:c10}, caption=\bfseries{Calculation of the Fourier--Legendre coefficients $C^{10}_{j_2 j_1}$}]{listings/nonlinear/coefficients/c10.py}
\lstinputlisting[language=Python, label={lst:c01}, caption=\bfseries{Calculation of the Fourier--Legendre coefficients $C^{01}_{j_2 j_1}$}]{listings/nonlinear/coefficients/c01.py}
\lstinputlisting[language=Python, label={lst:c0000}, caption=\bfseries{Calculation of the Fourier--Legendre coefficients $C^{0000}_{j_4 j_3 j_2 j_1}$}]{listings/nonlinear/coefficients/c0000.py}
\lstinputlisting[language=Python, label={lst:c100}, caption=\bfseries{Calculation of the Fourier--Legendre coefficients $C^{100}_{j_3 j_2 j_1}$}]{listings/nonlinear/coefficients/c100.py}
\lstinputlisting[language=Python, label={lst:c010}, caption=\bfseries{Calculation of the Fourier--Legendre coefficients $C^{010}_{j_3 j_2 j_1}$}]{listings/nonlinear/coefficients/c010.py}
\lstinputlisting[language=Python, label={lst:c001}, caption=\bfseries{Calculation of the Fourier--Legendre coefficients $C^{001}_{j_3 j_2 j_1}$}]{listings/nonlinear/coefficients/c001.py}
\lstinputlisting[language=Python, label={lst:c00000}, caption=\bfseries{Calculation of the Fourier--Legendre coefficients $C^{00000}_{j_5 j_4 j_3 j_2 j_1}$}]{listings/nonlinear/coefficients/c00000.py}
\lstinputlisting[language=Python, label={lst:c20}, caption=\bfseries{Calculation of the Fourier--Legendre coefficients $C^{20}_{j_2 j_1}$}]{listings/nonlinear/coefficients/c20.py}
\lstinputlisting[language=Python, label={lst:c11}, caption=\bfseries{Calculation of the Fourier--Legendre coefficients $C^{11}_{j_2 j_1}$}]{listings/nonlinear/coefficients/c11.py}
\lstinputlisting[language=Python, label={lst:c02}, caption=\bfseries{Calculation of the Fourier--Legendre coefficients $C^{02}_{j_2 j_1}$}]{listings/nonlinear/coefficients/c02.py}
\lstinputlisting[language=Python, label={lst:c1000}, caption=\bfseries{Calculation of the Fourier--Legendre coefficients $C^{1000}_{j_4 j_3 j_2 j_1}$}]{listings/nonlinear/coefficients/c1000.py}
\lstinputlisting[language=Python, label={lst:c0100}, caption=\bfseries{Calculation of the Fourier--Legendre coefficients $C^{0100}_{j_4 j_3 j_2 j_1}$}]{listings/nonlinear/coefficients/c0100.py}
\lstinputlisting[language=Python, label={lst:c0010}, caption=\bfseries{Calculation of the Fourier--Legendre coefficients $C^{0010}_{j_4 j_3 j_2 j_1}$}]{listings/nonlinear/coefficients/c0010.py}
\lstinputlisting[language=Python, label={lst:c0001}, caption=\bfseries{Calculation of the Fourier--Legendre coefficients $C^{0001}_{j_4 j_3 j_2 j_1}$}]{listings/nonlinear/coefficients/c0001.py}
\lstinputlisting[language=Python, label={lst:c000000}, caption=\bfseries{Calculation of the Fourier--Legendre coefficients $C^{000000}_{j_6 j_5 j_4 j_3 j_2 j_1}$}]{listings/nonlinear/coefficients/c000000.py}
\lstinputlisting[language=Python, label={lst:main_new_c}, caption=\bfseries{Program entry for generation of new Fourier--Legendre coefficients}]{listings/main_new_c.py}
\lstinputlisting[language=Python, label={lst:new_c}, caption=\bfseries{Module for the Fourier--Legendre coefficients generation}]{listings/nonlinear/new_c.py}

\subsubsection{Source Codes for Supplementary Differential Operators and Functions}

\lstinputlisting[language=Python, label={lst:l}, caption=\bfseries{Implementation of the differential operator $L$}]{listings/nonlinear/l.py}
\lstinputlisting[language=Python, label={lst:g}, caption=\bfseries{Implementation of the differential operator $G_0^{(i)}$}]{listings/nonlinear/g.py}
\vspace{-3mm}
\begin{lstlisting}[language=Python, label={lst:aj}, caption=\bfseries{Implementation of the function $\bar {\bf a}({\bf x}, t)$}]
from sympy import sympify, Matrix, Add

from mathematics.sde.nonlinear.symbolic.g import G
from mathematics.sde.nonlinear.symbolic.operator import Operator

class Aj(Operator):
    nargs = 4

    def __new__(cls, *args, **kwargs):
        """
        Creates new Aj object with given args
        Parameters
        ==========
        args
            bunch of necessary arguments
        Returns
        =======
        sympy.Expr
            formula to simplify and substitute
        """
        i, a, b, dxs = sympify(args)
        n = b.shape[0]
        m = b.shape[1]

        return Matrix([a[i, 0] -
                       (Add(*[0.5 * G(b[:, j], b[i, j], dxs)
                              for j in range(m)]))
                       for i in range(n)])
\end{lstlisting}
\vspace{-3mm}
\lstinputlisting[language=Python, label={lst:lj}, caption=\bfseries{Implementation of the differential operator $\bar L$}]{listings/nonlinear/lj.py}
\lstinputlisting[language=Python, label={lst:ind}, caption=\bfseries{Implementation of the indicator function ${\bf 1}_{\{i_1 = i_2\}}$}]{listings/nonlinear/ind.py}

\subsubsection{Source Codes for Iterated It\^{o} Stochastic Integrals Approximations Subprograms}

\lstinputlisting[language=Python, label={lst:i0}, caption=\bfseries{Approximation of It\^{o} stochastic integral $I_{(0)\tau_{p+1},\tau_p}^{(i_1)}$}]
{listings/nonlinear/ito/i0.py}
\lstinputlisting[language=Python, label={lst:i00}, caption=\bfseries{Approximation of iterated It\^{o} stochastic integral $I_{(00)\tau_{p+1},\tau_p}^{(i_1 i_2)}$}]
{listings/nonlinear/ito/i00.py}
\lstinputlisting[language=Python, label={lst:i1}, caption=\bfseries{Approximation of It\^{o} stochastic integral $I_{(1)\tau_{p+1},\tau_p}^{(i_1)}$}]
{listings/nonlinear/ito/i1.py}
\lstinputlisting[language=Python, label={lst:i000}, caption=\bfseries{Approximation of iterated It\^{o} stochastic integral $I_{(000)\tau_{p+1},\tau_p}^{(i_1 i_2 i_3)}$}]
{listings/nonlinear/ito/i000.py}
\lstinputlisting[language=Python, label={lst:i10}, caption=\bfseries{Approximation of iterated It\^{o} stochastic integral $I_{(10)\tau_{p+1},\tau_p}^{(i_1 i_2)}$}]
{listings/nonlinear/ito/i10.py}
\lstinputlisting[language=Python, label={lst:i01}, caption=\bfseries{Approximation of iterated It\^{o} stochastic integral $I_{(01)\tau_{p+1},\tau_p}^{(i_1 i_2)}$}]
{listings/nonlinear/ito/i01.py}
\lstinputlisting[language=Python, label={lst:i0000}, caption=\bfseries{Approximation of iterated It\^{o} stochastic integral $I_{(0000)\tau_{p+1},\tau_p}^{(i_1 i_2 i_3 i_4)}$}]{listings/nonlinear/ito/i0000.py}
\lstinputlisting[language=Python, label={lst:i2}, caption=\bfseries{Approximation of It\^{o} stochastic integral $I_{(2)\tau_{p+1},\tau_p}^{(i_1)}$}]
{listings/nonlinear/ito/i2.py}
\lstinputlisting[language=Python, label={lst:i100}, caption=\bfseries{Approximation of iterated It\^{o} stochastic integral $I_{(100)\tau_{p+1},\tau_p}^{(i_1 i_2 i_3)}$}]
{listings/nonlinear/ito/i100.py}
\lstinputlisting[language=Python, label={lst:i010}, caption=\bfseries{Approximation of iterated It\^{o} stochastic integral $I_{(010)\tau_{p+1},\tau_p}^{(i_1 i_2 i_3)}$}]
{listings/nonlinear/ito/i010.py}
\lstinputlisting[language=Python, label={lst:i001}, caption=\bfseries{Approximation of iterated It\^{o} stochastic integral $I_{(001)\tau_{p+1},\tau_p}^{(i_1 i_2 i_3)}$}]
{listings/nonlinear/ito/i001.py}
\lstinputlisting[language=Python, label={lst:i00000}, caption=\bfseries{Approximation of iterated It\^{o} stochastic integral $I_{(00000)\tau_{p+1},\tau_p}^{(i_1 i_2 i_3 i_4 i_5)}$}]{listings/nonlinear/ito/i00000.py}
\lstinputlisting[language=Python, label={lst:i20}, caption=\bfseries{Approximation of iterated It\^{o} stochastic integral $I_{(20)\tau_{p+1},\tau_p}^{(i_1 i_2)}$}]
{listings/nonlinear/ito/i20.py}
\lstinputlisting[language=Python, label={lst:i02}, caption=\bfseries{Approximation of iterated It\^{o} stochastic integral $I_{(02)\tau_{p+1},\tau_p}^{(i_1 i_2)}$}]
{listings/nonlinear/ito/i02.py}
\lstinputlisting[language=Python, label={lst:i11}, caption=\bfseries{Approximation of iterated It\^{o} stochastic integral $I_{(11)\tau_{p+1},\tau_p}^{(i_1 i_2)}$}]
{listings/nonlinear/ito/i11.py}
\lstinputlisting[language=Python, label={lst:i1000}, caption=\bfseries{Approximation of iterated It\^{o} stochastic integral $I_{(1000)\tau_{p+1},\tau_p}^{(i_1 i_2 i_3 i_4)}$}]
{listings/nonlinear/ito/i1000.py}
\lstinputlisting[language=Python, label={lst:i0100}, caption=\bfseries{Approximation of iterated It\^{o} stochastic integral $I_{(0100)\tau_{p+1},\tau_p}^{(i_1 i_2 i_3 i_4)}$}]
{listings/nonlinear/ito/i0100.py}
\lstinputlisting[language=Python, label={lst:i0010}, caption=\bfseries{Approximation of iterated It\^{o} stochastic integral $I_{(0010)\tau_{p+1},\tau_p}^{(i_1 i_2 i_3 i_4)}$}]
{listings/nonlinear/ito/i0010.py}
\lstinputlisting[language=Python, label={lst:i0001}, caption=\bfseries{Approximation of iterated It\^{o} stochastic integral $I_{(0001)\tau_{p+1},\tau_p}^{(i_1 i_2 i_3 i_4)}$}]
{listings/nonlinear/ito/i0001.py}
\lstinputlisting[language=Python, label={lst:i000000}, caption=\bfseries{Approximation of iterated It\^{o} stochastic integral $I_{(000000)\tau_{p+1},\tau_p}^{(i_1 i_2 i_3 i_4 i_5 i_6)}$}]
{listings/nonlinear/ito/i000000.py}

\subsubsection{Source Codes for Iterated Stratonovich Stochastic Integrals Approximations Subprograms}

\lstinputlisting[language=Python, label={lst:j0}, caption=\bfseries{Approximation of Stratonovich stochastic integral $I_{(0)\tau_{p+1},\tau_p}^{*(i_1)}$}]
{listings/nonlinear/stratonovich/j0.py}
\lstinputlisting[language=Python, label={lst:j00}, caption=\bfseries{Approximation of iterated Stratonovich stochastic integral $I_{(00)\tau_{p+1},\tau_p}^{*(i_1 i_2)}$}]
{listings/nonlinear/stratonovich/j00.py}
\vspace{-3mm}
\lstinputlisting[language=Python, label={lst:j1}, caption=\bfseries{Approximation of Stratonovich stochastic integral $I_{(1)T,t}^{*(i_1)}$ }]
{listings/nonlinear/stratonovich/j1.py}
\vspace{-3mm}
\lstinputlisting[language=Python, label={lst:j000}, caption=\bfseries{Approximation of iterated Stratonovich stochastic integral $I_{(000)\tau_{p+1},\tau_p}^{*(i_1 i_2 i_3)}$}]
{listings/nonlinear/stratonovich/j000.py}
\vspace{-3mm}
\lstinputlisting[language=Python, label={lst:j10}, caption=\bfseries{Approximation of iterated Stratonovich stochastic integral $I_{(10)\tau_{p+1},\tau_p}^{*(i_1 i_2)}$}]
{listings/nonlinear/stratonovich/j10.py}
\vspace{3mm}
\lstinputlisting[language=Python, label={lst:j01}, caption=\bfseries{Approximation of iterated Stratonovich stochastic integral $I_{(01)\tau_{p+1},\tau_p}^{*(i_1 i_2)}$ }]
{listings/nonlinear/stratonovich/j01.py}
\lstinputlisting[language=Python, label={lst:j0000}, caption=\bfseries{Approximation of iterated Stratonovich stochastic integral $I_{(0000)\tau_{p+1},\tau_p}^{*(i_1 i_2 i_3 i_4)}$}]{listings/nonlinear/stratonovich/j0000.py}
\lstinputlisting[language=Python, label={lst:j2}, caption=\bfseries{Approximation of Stratonovich stochastic integral $I_{(2)\tau_{p+1},\tau_p}^{*(i_1)}$}]
{listings/nonlinear/stratonovich/j2.py}
\lstinputlisting[language=Python, label={lst:j100}, caption=\bfseries{Approximation of iterated Stratonovich stochastic integral $I_{(100)\tau_{p+1},\tau_p}^{*(i_1 i_2 i_3)}$}]
{listings/nonlinear/stratonovich/j100.py}
\lstinputlisting[language=Python, label={lst:j010}, caption=\bfseries{Approximation of iterated Stratonovich stochastic integral $I_{(010)\tau_{p+1},\tau_p}^{*(i_1 i_2 i_3)}$}]
{listings/nonlinear/stratonovich/j010.py}
\lstinputlisting[language=Python, label={lst:j001}, caption=\bfseries{Approximation of iterated Stratonovich stochastic integral $I_{(001)\tau_{p+1},\tau_p}^{*(i_1 i_2 i_3)}$}]
{listings/nonlinear/stratonovich/j001.py}
\lstinputlisting[language=Python, label={lst:j00000}, caption=\bfseries{Approximation of iterated Stratonovich stochastic integral $I_{(00000)\tau_{p+1},\tau_p}^{*(i_1 i_2 i_3 i_4 i_5)}$}]{listings/nonlinear/stratonovich/j00000.py}
\lstinputlisting[language=Python, label={lst:j20}, caption=\bfseries{Approximation of iterated Stratonovich stochastic integral $I_{(20)\tau_{p+1},\tau_p}^{*(i_1 i_2)}$}]
{listings/nonlinear/stratonovich/j20.py}
\lstinputlisting[language=Python, label={lst:j02}, caption=\bfseries{Approximation of iterated Stratonovich stochastic integral $I_{(02)\tau_{p+1},\tau_p}^{*(i_1 i_2)}$}]
{listings/nonlinear/stratonovich/j02.py}
\lstinputlisting[language=Python, label={lst:j11}, caption=\bfseries{Approximation of iterated Stratonovich stochastic integral $I_{(11)\tau_{p+1},\tau_p}^{*(i_1 i_2)}$}]
{listings/nonlinear/stratonovich/j11.py}
\lstinputlisting[language=Python, label={lst:j1000}, caption=\bfseries{Approximation of iterated Stratonovich stochastic integral $I_{(1000)\tau_{p+1},\tau_p}^{*(i_1 i_2 i_3 i_4)}$}]{listings/nonlinear/stratonovich/j1000.py}
\lstinputlisting[language=Python, label={lst:j0100}, caption=\bfseries{Approximation of iterated Stratonovich stochastic integral $I_{(0100)\tau_{p+1},\tau_p}^{*(i_1 i_2 i_3 i_4)}$}]{listings/nonlinear/stratonovich/j0100.py}
\lstinputlisting[language=Python, label={lst:j0010}, caption=\bfseries{Approximation of iterated Stratonovich stochastic integral $I_{(0010)\tau_{p+1},\tau_p}^{*(i_1 i_2 i_3 i_4)}$}]{listings/nonlinear/stratonovich/j0010.py}
\lstinputlisting[language=Python, label={lst:j0001}, caption=\bfseries{Approximation of iterated Stratonovich stochastic integral $I_{(0001)\tau_{p+1},\tau_p}^{*(i_1 i_2 i_3 i_4)}$}]{listings/nonlinear/stratonovich/j0001.py}
\lstinputlisting[language=Python, label={lst:j000000}, caption=\bfseries{Approximation of iterated Stratonovich stochastic integral $I_{(000000)\tau_{p+1},\tau_p}^{*(i_1 i_2 i_3 i_4 i_5 i_6)}$}]{listings/nonlinear/stratonovich/j000000.py}

\subsubsection{Source Codes for Calculation of the Numbers $q,$ $q_1$,$\ldots$, $q_{15}$}

\lstinputlisting[language=Python, label={lst:q}, caption=\bfseries{Calculation of the numbers $q,$ $q_1$,$\ldots$, $q_{15}$}]{listings/nonlinear/q.py}

\subsubsection{Source Codes for Strong Taylor--It\^{o} Numerical Schemes with
Convergence Orders $0.5$, $1.0$, $1.5$, $2.0$, $2.5$, and $3.0$ for It\^{o} SDEs}

\vspace{-3mm}
\lstinputlisting[language=Python, label={lst:euler_dri}, caption=\bfseries{Euler scheme modeling subprogram}]{listings/nonlinear/drivers/euler.py}
\lstinputlisting[language=Python, label={lst:euler}, caption=\bfseries{Euler scheme}]{listings/nonlinear/schemes/euler.py}
\vspace{3mm}
\lstinputlisting[language=Python, label={lst:milstein_dri}, caption=\bfseries{Milstein scheme modeling subprogram}]{listings/nonlinear/drivers/milstein.py}
\lstinputlisting[language=Python, label={lst:milstein}, caption=\bfseries{Milstein scheme}]{listings/nonlinear/schemes/milstein.py}
\vspace{-2mm}
\lstinputlisting[language=Python, label={lst:strong_taylor_ito_1p5_dri}, caption=\bfseries{Strong Taylor--It\^{o} scheme with convergence order $1.5$ modeling subprogram}]{listings/nonlinear/drivers/strong_taylor_ito_1p5.py}
\vspace{-2mm}
\lstinputlisting[language=Python, label={lst:strong_taylor_ito_1p5}, caption=\bfseries{Strong Taylor--It\^{o} scheme with convergence order $1.5$}]{listings/nonlinear/schemes/strong_taylor_ito_1p5.py}
\lstinputlisting[language=Python, label={lst:strong_taylor_ito_2p0_dri}, caption=\bfseries{Strong Taylor--It\^{o} scheme with convergence order $2.0$ modeling subprogram}]{listings/nonlinear/drivers/strong_taylor_ito_2p0.py}
\lstinputlisting[language=Python, label={lst:strong_taylor_ito_2p0}, caption=\bfseries{Strong Taylor--It\^{o} scheme with convergence order $2.0$}]{listings/nonlinear/schemes/strong_taylor_ito_2p0.py}
\lstinputlisting[language=Python, label={lst:strong_taylor_ito_2p5_dri}, caption=\bfseries{Strong Taylor--It\^{o} scheme with convergence order $2.5$ modeling subprogram}]{listings/nonlinear/drivers/strong_taylor_ito_2p5.py}
\lstinputlisting[language=Python, label={lst:strong_taylor_ito_2p5}, caption=\bfseries{Strong Taylor--It\^{o} scheme with convergence order $2.5$}]{listings/nonlinear/schemes/strong_taylor_ito_2p5.py}
\lstinputlisting[language=Python, label={lst:strong_taylor_ito_3p0_dri}, caption=\bfseries{Strong Taylor--It\^{o} scheme with convergence order $3.0$ modeling subprogram}]{listings/nonlinear/drivers/strong_taylor_ito_3p0.py}
\lstinputlisting[language=Python, label={lst:strong_taylor_ito_3p0}, caption=\bfseries{Strong Taylor--It\^{o} scheme with convergence order $3.0$}]{listings/nonlinear/schemes/strong_taylor_ito_3p0.py}

\subsubsection{Source Codes for Strong Taylor--Stratonovich Numerical Schemes 
with Convergence Orders $1.0$, $1.5$, $2.0$, $2.5$, and $3.0$ for It\^{o} SDEs}

\lstinputlisting[language=Python, label={lst:strong_taylor_stratonovich_1p0_dri}, caption=\bfseries{Strong Taylor--Stratonovich scheme with convergence order $1.0$ modeling subprogram}]{listings/nonlinear/drivers/strong_taylor_stratonovich_1p0.py}
\lstinputlisting[language=Python, label={lst:strong_taylor_stratonovich_1p0}, caption=\bfseries{Strong Taylor--Stratonovich scheme with convergence order $1.0$}]{listings/nonlinear/schemes/strong_taylor_stratonovich_1p0.py}
\lstinputlisting[language=Python, label={lst:strong_taylor_stratonovich_1p5_dri}, caption=\bfseries{Strong Taylor--Stratonovich scheme with convergence order $1.5$ modeling subprogram}]{listings/nonlinear/drivers/strong_taylor_stratonovich_1p5.py}
\lstinputlisting[language=Python, label={lst:strong_taylor_stratonovich_1p5}, caption=\bfseries{Strong Taylor--Stratonovich scheme with convergence order $1.5$}]{listings/nonlinear/schemes/strong_taylor_stratonovich_1p5.py}
\lstinputlisting[language=Python, label={lst:strong_taylor_stratonovich_2p0_dri}, caption=\bfseries{Strong Taylor--Stratonovich scheme with convergence order $2.0$ modeling subprogram}]{listings/nonlinear/drivers/strong_taylor_stratonovich_2p0.py}
\lstinputlisting[language=Python, label={lst:strong_taylor_stratonovich_2p0}, caption=\bfseries{Strong Taylor--Stratonovich scheme with convergence order $2.0$}]{listings/nonlinear/schemes/strong_taylor_stratonovich_2p0.py}
\lstinputlisting[language=Python, label={lst:strong_taylor_stratonovich_2p5_dri}, caption=\bfseries{Strong Taylor--Stratonovich scheme with convergence order $2.5$ modeling subprogram}]{listings/nonlinear/drivers/strong_taylor_stratonovich_2p5.py}
\lstinputlisting[language=Python, label={lst:strong_taylor_stratonovich_2p5}, caption=\bfseries{Strong Taylor--Stratonovich scheme with convergence order $2.5$}]{listings/nonlinear/schemes/strong_taylor_stratonovich_2p5.py}
\lstinputlisting[language=Python, label={lst:strong_taylor_stratonovich_3p0_dri}, caption=\bfseries{Strong Taylor--Stratonovich scheme with convergence order $3.0$ modeling subprogram}]{listings/nonlinear/drivers/strong_taylor_stratonovich_3p0.py}
\lstinputlisting[language=Python, label={lst:strong_taylor_stratonovich_3p0}, caption=\bfseries{Strong Taylor--Stratonovich scheme with convergence order $3.0$}]{listings/nonlinear/schemes/strong_taylor_stratonovich_3p0.py}

\subsection{Source Codes for Linear Stationary Systems of It\^{o} SDEs}

\lstinputlisting[language=Python, label={lst:matrix}, caption=\bfseries{Implementation of supplementary functions}]{listings/linear/matrix.py}
\lstinputlisting[language=Python, label={lst:dindet}, caption=\bfseries{Implementation of Algorithm 11.2 \cite{56}}]{listings/linear/dindet.py}
\lstinputlisting[language=Python, label={lst:stoch}, caption=\bfseries{Implementation of Algorithm 11.6 \cite{56}}]{listings/linear/stoch.py}
\vspace{3mm}
\lstinputlisting[language=Python, label={lst:distortions}, caption=\bfseries{Implementation of the vector function ${\bf u}(t)$}]{listings/linear/distortions.py}
\lstinputlisting[language=Python, label={lst:integration}, caption=\bfseries{Modeling of linear system of It\^{o} SDEs}]{listings/linear/integration.py}

\subsection{Source Codes for Utilities and Initialization}

\lstinputlisting[language=Python, label={lst:initialization}, caption=\bfseries{Initialization module}]{listings/initialization/initialization.py}
\lstinputlisting[language=Python, label={lst:database}, caption=\bfseries{Module for database initialization}]{listings/initialization/database.py}
\lstinputlisting[language=Python, label={lst:database}, caption=\bfseries{Database module}]{listings/tools/database.py}
\lstinputlisting[language=Python, label={lst:fsys}, caption=\bfseries{File system utilities}]{listings/tools/fsys.py}

%% file: chapters/bibliography.tex
\section{Future Work}

Considering the future work, it is important to say that symbolic algebra gives a wide field for 
optimizations of modeling process. Symbolic operations is actually operations with strings.
Such operations has relatively high complexity and slows down modeling process significantly.
One of possible ways to improve modeling performance is to parallelize computations. Since 
strong numerical schemes for It\^o SDEs have massive amount of terms this idea appears justified.

The strong numerical schemes for It\^o SDEs seem to be easily optimizable, on the other hand, 
superpositions of the differential operators (\ref{2.3}), (\ref{2.4}), and (\ref{2.4xxx}) are not.
They are called recursively during calculation process which is more difficult to parallelize 
than strong numerical schemes for It\^o SDEs. Differential operators obviously include differentiating which is high cost and
optimization of them is a dedicated issue.

In the future, it is possible to improve the SDE-MATH software package 
in a number of other directions. In particular, high-order strong numerical 
methods of the Runge-Kutta type \cite{2}, \cite{7}, \cite{40}, \cite{56} (including multistep numerical
methods \cite{2}, \cite{7}, \cite{40}, \cite{56})  for It\^o SDEs can be implemented. In addition,
software for solving the filtering problem and the problem of stochastic optimal
control can also be developed. These improvements will lead to changes of the 
graphical user interface due to new features.

\vspace{3mm}

\renewcommand{\refname}{\rm
{\bf Bibliography}}

\vspace{3mm}

\addcontentsline{toc}{section}{{\normalsize References}}